\documentclass[paper=letter, fontsize=11pt, chapterprefix=on, oneside]{scrbook}
\usepackage{amsfonts}
\usepackage{amssymb}
\usepackage{rotating}
\usepackage{amsmath}
\usepackage{amsrefs}
\usepackage{theorem}
\usepackage{epsfig}
\usepackage{color}
\usepackage{url}
\usepackage[top=2.54cm, bottom=2.54cm, left=3.05cm, right=2.54cm]{geometry}
\usepackage{tikz}
\usepackage{pdflscape}
\usepackage{kbordermatrix}
\usepackage{chngcntr}
\usepackage{tocloft}
\usepackage[font={footnotesize}]{caption}
\usepackage{float}
\usetikzlibrary{decorations.pathreplacing}
\usetikzlibrary{arrows,automata}
\usepackage[latin1]{inputenc}
\usepackage{enumerate}
\usepackage{multirow}
\usepackage{array}
\usepackage{afterpage}
\usepackage{listings}
\usepackage{color}

\DeclareOldFontCommand{\rm}{\normalfont\rmfamily}{\mathrm}
\DeclareOldFontCommand{\sf}{\normalfont\sffamily}{\mathsf}
\DeclareOldFontCommand{\tt}{\normalfont\ttfamily}{\mathtt}
\DeclareOldFontCommand{\bf}{\normalfont\bfseries}{\mathbf}
\DeclareOldFontCommand{\it}{\normalfont\itshape}{\mathit}
\DeclareOldFontCommand{\sl}{\normalfont\slshape}{\@nomath\sl}
\DeclareOldFontCommand{\sc}{\normalfont\scshape}{\@nomath\sc}
\DeclareRobustCommand*\cal{\@fontswitch\relax\mathcal}
\DeclareRobustCommand*\mit{\@fontswitch\relax\mathnormal}

\newcolumntype{L}{>{\centering\arraybackslash}m{3.5cm}}
\usetikzlibrary{patterns,decorations.markings,decorations.pathreplacing}
\makeatletter
\pgfmathdeclarefunction{square}{1}{%
\begingroup
 \pgfmathparse{#1*#1}
 \pgfmath@smuggleone\pgfmathresult%
\endgroup}
\counterwithin*{equation}{chapter}
\usepackage[colorlinks,citecolor=black,linkcolor=black,urlcolor=black]{hyperref} 

\def\previousdegrees{\gdef\@previousdegrees}
\def\degreetitle{\gdef\@degreetitle}

\def\faculty{\gdef\@faculty}
\def\department{\gdef\@department}
\def\institution{\gdef\@institution}
\def\campus{\gdef\@campus}
\def\majorprofessor{\gdef\@majorprofessor}
\def\committeeA{\gdef\@committeeA}
\def\committeeB{\gdef\@committeeB}
\def\committeeC{\gdef\@committeeC}
\def\admin{\gdef\@admin}

\def\month{\gdef\@month}
\def\year{\gdef\@year}

\def\@institution{University of Idaho}
\def\@faculty{College of Graduate Studies}
\def\@campus{Moscow}
\def\@majorprofessor{Alexander Woo, Ph.D.}
\def\@committeeA{Jennifer Johnson-Leung, Ph.D.}
\def\@committeeB{Terence Soule, Ph.D.}
\def\@committeeC{Stefan Tohaneanu, Ph.D.}
\def\@admin{Christopher Williams, Ph.D.}

\def\maketitle{
\thispagestyle{empty}
\begin{center}
\vspace*{\fill}
  {\LARGE\@title}
  \bigskip\bigskip\bigskip

  \textit{A Dissertation}\\[12pt]
  \textit{Presented in Partial Fulfillment of the Requirements for the Degree of}\\[12pt]
  {\large\textrm\@degreetitle\\[12pt]}
  \textit{with a}\\[12pt]
  {Major in \@department\\[12pt]
   \textit{in the}\\[12pt]
    \@faculty\\[12pt]
    \@institution\\}
  \bigskip\bigskip

  \textit{by}\\[12pt]
  {\large\textrm\@author\\[12pt]}
  \bigskip

  \textit{Major Professor}\\
  {\textrm\@majorprofessor}\\[12pt]
  \textit{Committee}\\
  {\textrm\@committeeA}\\
  {\textrm\@committeeB}\\
  {\textrm\@committeeC}\\[12pt]
  \textit{Department Administrator}\\
  {\textrm\@admin}\\[12pt]
  \@month\ \@year\\[12pt]
\vspace*{\fill}
\end{center}
\newpage
}

\title{ENUMERATION OF PERMUTATIONS\\
  INDEXING LOCAL COMPLETE INTERSECTION\vspace{10pt}\\
  SCHUBERT VARIETIES}
\author{Masaki Ikeda}
\month{\textrm{June}} \year{2016}
\previousdegrees{M.Sc., University of Idaho, 2012}
\degreetitle{Doctor of Philosophy}
\department{Mathematics}

\usepackage{gitinfo}
\usepackage{fancyhdr}
\newif\ifversioninfo\versioninfofalse
\ifversioninfo
\fancypagestyle{plain}{
  \fancyhf{}
  \fancyfoot[LE,RO]{\footnotesize\thepage}
  \fancyfoot[LO,RE]{\footnotesize\gitAuthorDate~(\gitAbbrevHash)}}
\else
\fancypagestyle{plain}{\fancyhf{}\fancyhead[LE,RO]{\footnotesize\thepage}}
\fi
\theoremstyle{plain}
\theorembodyfont{\normalfont}

\newtheorem{theorem}{Theorem}[chapter]
\newtheorem{lemma}[theorem]{Lemma}

\newtheorem{proposition}[theorem]{Proposition}

\def \C{{\mathcal C}}

\def \Z{{\mathbb Z}}
\def \A{{\mathcal A}}

\def \S{{\mathcal S}}

\def \a{{\alpha}}
\def \b{{\beta}}
\def \d{{\delta}}
\def \e{\varepsilon}
\def \l{{\ell}}

\def \s{{\sigma}}
\def \cont{\preceq}
\def \avd{\not\preceq}
\newcommand {\PS}[1]{\pi^{-1}({#1})}

\DeclareMathOperator{\LRmax}{LRmax}
\DeclareMathOperator{\RLmin}{RLmin}
\DeclareMathOperator{\nwsum}{\,\!\.\.\raisebox{-0.032cm}{\begin{tikzpicture}[->]
\draw[>=stealth,black] (-0.01,0.01) -- (-0.19,0.19);
\draw (-0.1,0.1) circle (0.123);
\end{tikzpicture}\,\!\.\.}}
\DeclareMathOperator{\sesum}{\raisebox{0.195cm}{\rotatebox{180}{$\nwsum$}}}
\DeclareMathOperator{\nwssum}{\,\!\.\.\raisebox{-0.032cm}{\begin{tikzpicture}[scale=0.8,->]
\draw[>=stealth,black] (-0.01,0.01) -- (-0.19,0.19);
\draw (-0.1,0.1) circle (0.123);
\end{tikzpicture}\,\!\.\.}}
\DeclareMathOperator{\sessum}{\raisebox{0.195cm}{\rotatebox{180}{$\nwssum$}}}

\def \Av{\textrm{Av}}
\def \nn{{\nonumber}}

\def \ra{{\,\,\rightarrow\,\,}}

\def \vhh{\vspace{.25in}}
\def \vhhh{\vspace{.125in}}
\def \.{\:\!}
\def \_{\text{\quad\.}}

\def \sbst{\subseteq}

\DeclareFontFamily{U}{MnSymbolF}{}
\DeclareSymbolFont{mnsymbols}{U}{MnSymbolF}{m}{n}
\DeclareFontShape{U}{MnSymbolF}{m}{n}{
<-6> MnSymbolF5
<6-7> MnSymbolF6
<7-8> MnSymbolF7
<8-9> MnSymbolF8
<9-10> MnSymbolF9
<10-12> MnSymbolF10
<12-> MnSymbolF12}{}
\DeclareMathSymbol{\bigominus}{\mathop}{mnsymbols}{55}

\DeclareFontFamily{U}{mathx}{\hyphenchar\font45}
\DeclareFontShape{U}{mathx}{m}{n}{
      <5> <6> <7> <8> <9> <10>
      <10.95> <12> <14.4> <17.28> <20.74> <24.88>
      mathx10
      }{}
\DeclareSymbolFont{mathx}{U}{mathx}{m}{n}
\DeclareMathSymbol{\bigboxplus}      {1}{mathx}{"D0}

\def \ds{\displaystyle}

\newcommand{\plotPerm}[1]
{
  \foreach \y [count=\len] in {#1} ;
  \foreach \x in {1,...,\len} \draw[black,thick] (\x,1)--(\x,\len) (1,\x)--(\len,\x);
  \foreach \y [count=\x] in {#1} \draw [color=black,fill=black] (\x,\y) circle [radius=0.1];
}

\setkomafont{disposition}{\normalfont\bfseries}
\addtokomafont{chapter}{\Large}
\addtokomafont{section}{\large}

\usepackage{titletoc}
\titlecontents{chapter}[0pt]%
{\addvspace{.7pc}}{\hspace*{1.2em}\rm\bf\makebox[0pt][r]{\thecontentslabel\enspace}}%
{\rm\bf}{\titlerule*[0.75pc]{.} \contentspage}
\usepackage{etoolbox}
\makeatletter
\patchcmd{\@caption}{\csname the#1\endcsname}{\hspace{-0.23in}\csname fnum@#1\endcsname:}{}{}
\renewcommand*\l@figure{\@dottedtocline{1}{1.5em}{4.5em}} 
\let\l@table\l@figure 
\makeatother
\makeatletter
\g@addto@macro\bfseries{\boldmath}
\makeatother

\newcolumntype{C}[1]{>{\centering\arraybackslash$}p{#1}<{$}}
\newlength{\mycolwd}                                         
\renewcommand*{\thesection}{\arabic{chapter}.\arabic{section}}

\DeclareFixedFont{\ttb}{T1}{txtt}{bx}{n}{12} 
\DeclareFixedFont{\ttm}{T1}{txtt}{m}{n}{12}  
\definecolor{deepblue}{rgb}{0,0,0.5}
\definecolor{deepred}{rgb}{0.6,0,0}
\definecolor{deepgreen}{rgb}{0,0.5,0}

\newcommand\pythonstyle{\lstset{
    language=Python,
    basicstyle=\footnotesize,
    otherkeywords={self},
    keywordstyle=\footnotesize\color{deepblue},
    emph={__init__},
    emphstyle=\footnotesize\color{deepred},
    stringstyle=\color{deepgreen},
    frame=single,
    showstringspaces=false  ,
    breaklines=true,
    numbers=left,
    numberstyle=\footnotesize,
    tabsize=3,
    breakatwhitespace=false
}}

\lstnewenvironment{python}[1][]
{
    \pythonstyle
    \lstset{#1}
}{}

\newcommand\pythoninline[1]{{\pythonstyle\lstinline!#1!}}

\begin{document}
\counterwithin{figure}{chapter}
\counterwithin{table}{chapter}
\counterwithin{theorem}{chapter}
\counterwithout{footnote}{chapter}
\maketitle
\phantomsection
\addcontentsline{toc}{chapter}{Authorization to Submit Dissertation}
\chapter*{Authorization to Submit Dissertation}
\pagenumbering{roman}
\setcounter{page}{2}
This dissertation of Masaki Ikeda, submitted for the degree of Doctor of Philosophy with a Major in Mathematics and titled ``Enumeration of permutations  indexing local complete intersection Schubert varieties,'' has been reviewed in final form. Permission, as indicated by the signatures and dates below, is now granted to submit final copies to the College of Graduate Studies for approval.

\vspace{1in}

\begin{center}
\begin{tabular}{ p{.3\textwidth}  p{.4\textwidth}  p{.18\textwidth} }
Major Professor: & \hrulefill & \hrulefill\\[-8pt]
 & Alexander Woo, Ph.D. & Date\\[18pt]
Committee Members: & \hrulefill & \hrulefill\\[-8pt]
 & Jennifer Johnson-Leung, Ph.D. & Date\\[18pt]
 & \hrulefill & \hrulefill\\[-8pt]
 & Terence Soule, Ph.D. & Date\\[18pt]
& \hrulefill & \hrulefill\\[-8pt]
 & Stefan Tohaneanu, Ph.D. & Date\\[18pt]
Department \\
Administrator: & \hrulefill & \hrulefill\\[-8pt]
 & Christopher Williams, Ph.D. & Date\\[18pt]
\end{tabular}
\end{center}
\newpage
\phantomsection
\addcontentsline{toc}{chapter}{Abstract}
\chapter*{Abstract}
We find the generating function for the permutation class $\mathcal{A}'=\text{Av}(52341,53241,52431,35142,\allowbreak 42513,351624)$. Partial motivation for this work comes from algebraic geometry. In particular, certain classes of Schubert varieties are indexed by permutations in some permutation classes. For example, as shown by Lakshmibai and Sandhya \cite{LS1990}, a smooth Schubert variety is indexed by a permutation in the class $\textrm{Av}(3412,4231)$, and a Schubert variety defined by inclusions (dbi for short) has an index in $\mathcal{A}=\textrm{Av}(4231,35142,42513,351624)$, as shown by Gasharov and Reiner \cite{GR2002}. In addition to those results, \'{U}lfarsson and Woo \cite{UW2013} showed Schubert varieties which are local complete intersections (lci for short) are indexed by permutations in $\A'$.\\

The enumeration of the permutations indexing smooth Schubert varieties was initially found by Haiman \cite{H1992}, and then also discussed by Bousquet-M\'{e}lou and Butler \cite{BB2007}. Furthermore, Albert and Brignall \cite{AB2013} discovered the enumeration of Schubert varieties defined by inclusions. This dissertation completes the enumeration of the class $\mathcal{A}'$ by extending the method used by Albert and Brignall to enumerate $\mathcal{A}$ in \cite{AB2013}.
\newpage
\phantomsection
\addcontentsline{toc}{chapter}{Acknowledgements}
\chapter*{Acknowledgements}
I couldn't have come this far academically without support provided by my family, friends, and former/current professors. The biggest appreciation goes to my advisor, Alexander Woo. He not only guided me to the right direction with his intelligence and knowledge, but also became a great connection to the community of permutation patterns. I could tell that he always cared about me and my career as a mathematician. One day, I hope to become a great inspiration to my own students.\\

I am truly thankful that I chose the University of Idaho for my graduate career. Teaching assistantship they provided me each year helped me highly, not only financially, but also for me to gain remarkable experiences as an instructor. Thanks to my committee members, Jennifer, Stefan and Terry, and to wonderful staff members, Jana, Melissa, Jaclyn and Stacey, and all other great professors who helped me to become a better scholar every day. Special thanks goes to Monte Boisen, the former department chair who supported me even after his retirement. My stay in Moscow was a truly amazing experience.\\

I would like to thank Western Oregon University for offering a great position to start my professional career in academia as well as my wonderful undergraduate experience. It is a shame that I had to have so many people wait for such a long time as I worked to finish my thesis. I appreciate all the support from Sue, Debbie and everyone in mathematics department here at WOU. Especially, Mike Ward, who first taught me the beauty of mathematics. Without his courses, I would have never stepped into the path I took.\\

Permutation patterns community was truly supportive. I met wonderful scholars who were helpful and motivating to complete my thesis. Thanks to Michael Albert, who gave me such kind words to push me through the toughest time as well as his extremely helpful program, PermLab to visualize my study. I am very happy to join such friendly community.\\

Lastly, I want to thank my host family in Oregon and friends. I am grateful to have so many caring people surrounding me in my life. Jonny Olson, my best friend, was a true motivation and distraction throughout my graduate career. Congratulations to him for achieving Ph.D. degree at Louisiana State University. Thanks to Veronica, Frank, Emily, Josh, Amanda and Caitlin, my second family since 2003, and my best friends, Jeff, Sierra, Ben, Jesse, Jenna, Hannah, Trevor, Morgan, Josh, Don, Mallory, Veronica, Diadra, Karissa and the ones that I shamefully forgot.
\newpage
\phantomsection
\addcontentsline{toc}{chapter}{Dedication}
\begin{center}
\vspace*{\fill}
\Large\textbf{Dedication}\\
\normalsize\textit{To my mother, Yuriko, and my father, Katsunobu.}
\vspace*{\fill}
\end{center}
\newpage
\phantomsection
\addcontentsline{toc}{chapter}{Table of Contents}
\textrm{}\\
\vspace{-1.1in}\tableofcontents
\newpage
\phantomsection
\addcontentsline{toc}{chapter}{List of Figures}
\textrm{}\\
\vspace{-1.1in}
\listoffigures
\newpage
\phantomsection
\addcontentsline{toc}{chapter}{List of Tables}
\textrm{}\\
\vspace{-1.1in}
\listoftables
\newpage
\pagenumbering{arabic}
\setcounter{page}{1}
\chapter{Introduction}
\section{Result}
Denote by $\Av(B)$ the set of all permutations avoiding every permutation in $B$. One essential study of permutation classes is to find the generating function for $\Av(B)$ with a particular set of permutations $B$. Our final goal is to prove the following theorem.\\
\\
\textbf{Theorem.} \textit{The generating function for the class $\A'$ is defined by}\\
\[
f_{\A'}=\frac{\sum_{i=0}^{5} a_i \bar{G}^i}{\sum_{i=0}^{6} b_i \bar{G}^i}
\]
\\
\textit{where $\bar{G}=G-1$ and $G$ is the generating function satisfying the equation}\\
\[
\bar{G}=1+\frac{x\bar{G}}{1-x\bar{G}^2},
\]
\textit{and}\\
\[
\begin{array}{l}
a_0=-1 + 14 x - 39 x^2 + 28 x^3 + 9 x^4 - 11 x^5 + x^6,\\
a_1=-12 + 81 x - 100 x^2 + 15 x^3 + 46 x^4 - 19 x^5,\\
a_2=-8 + 35 x - 20 x^2 - 25 x^3 + 31 x^4 - 6 x^5 - x^6,\\
a_3=7,\qquad a_4=1,\qquad a_5=-2.
\end{array}\]

\[
\begin{array}{l}
b_0=-1 + 57 x - 125 x^2 + 143 x^3 - 48 x^4 - 64 x^5 + 51 x^6 - 2 x^8,\\
b_1=-54 + 260 x - 386 x^2 + 250 x^3 + 81 x^4 - 226 x^5 + 74 x^6 + 15 x^7 - 3 x^8,\\
b_2=-18 + 114 x - 104 x^2 - 22 x^3 + 148 x^4 - 123 x^5 + 11 x^6 + 14 x^7 - x^8,\\
b_3=24,\qquad b_4=-2,\qquad b_5=-5,\qquad b_6=1.
\end{array}
\]
\\
In this chapter, we introduce history of permutation class study as well as the place of dissertation in the literature. Detailed definitions are given in the next chapter.\\
\section{History}
The concept of permutation avoidance first appeared in the literature in 1915. In \cite{M1915}, MacMahon proved that the number of permutations of length $n$ which can be partitioned into two decreasing subsequences is the Catalan numbers. Also in 1935, Erd\H{o}s and Szekeres \cite{ES1935} showed that, given two positive integers $a$, $b$ and a sequence of $n$ real numbers $x = x_1,\ldots,x_n$ with $n=ab+1$, $x$ either contains a strictly increasing subsequence of length $a+1$ or a strictly decreasing subsequence of length $b+1$. Despite these early results, we consider the study of permutation classes to have begun in 1968 with Knuth's \textit{The Art of Computer Programming} \cite{K1968}. Knuth proved a permutation $\pi$ is stack-sortable if and only if $\pi$ avoids 231, and stack-sortable permutations are also counted by the Catalan numbers. The result of Knuth brought up the notion of permutation avoidance and the study of permutation classes. Within two decades of Knuth's contribution, many enumeration results were discovered by various researchers. In particular, Simion and Schmidt \cite{SS1985} summarized permutation classes avoiding permutations of length 3.\\

Based on some earlier results, Stanley and Wilf separately conjectured in the late 1980s that, for any permutation $\pi$, there exists a constant $C_\pi$ whose $n$-th power is an upper bound for the number of permutations of length $n$ avoiding $\pi$. In other words, while the growth rate of the set of length $n$ permutations is factorial, having a restriction of avoiding an arbitrary single permutation reduces the growth rate to be exponential. This is known as the Stanley-Wilf conjecture. Some partial results were proved by B\'{o}na \cite{B1999} and Alon together with Friedgut \cite{AF2000}. A breakthrough was made by Klazar in 2000 \cite{K2000}, when he showed that the {F}\"uredi-{H}ajnal conjecture implies the {S}tanley-{W}ilf conjecture. The {S}tanley-{W}ilf conjecture remained unproven for almost two decades until Marcus and Tardos proved it in 2004 by proving the {F}\"uredi-{H}ajnal conjecture \cite{MT2004}.\\

Although the proof of Marcus and Tardos gives an upper bound for the growth rate constant $C_\pi$ depending on the length of $\pi$, the precise growth rates for most permutation classes are still unknown. Hence, one of the main problems in present research is to characterize the growth rates of permutation classes. Pratt \cite{P1973} as well as Spielman together with B\'{o}na \cite{BS2000} showed the existence of permutation classes containing infinite antichains. Since such a class contains uncountably many distinct subclasses, their result implied there exist uncountably many distinct growth rates. Then in \cite{KK2003}, Kaiser and Klazar showed the only possible growth rates of any permutation class that is less than 2 are positive solutions to $1-2x^k+x^{k+1}$ for some $k\geq 0$. Later, Klazar also showed there are only countably many permutation classes with growth rate less than 2. As an extension of this result, Vatter showed two noteworthy results in \cite{V2011}: The smallest possible growth rate greater than 2 for a permutation class is the unique positive root of $1+2x+x^2+x^3-x^4$, which is approximately 2.06599, and the smallest growth rate for which there are uncountably many permutation classes is the unique positive root of $1+2x^2-x^3$, which is approximately 2.20557.\\

Another notable contribution to the study of growth rates was made by Fox in 2014 \cite{F2014}. It was believed that the growth rates of permutation classes avoiding a single permutation of length $n$ grows at most quadratically. In \cite{F2014}, Fox showed this is false, but the function $g(k)=\max\{\C_\pi:|\pi|=k\}$ grows mildly exponentially.\\

While understanding growth rate constants is the primary interest of many researchers, classifying necessary and sufficient conditions for a permutation class to have a specific type of generating function is also an important question. For example, in \cite{AABRV2013}, authors defined the notion of geometric griddable class, and showed every geometrically griddable class has a rational generating function. In 1996, Noonan and Zeilberger conjectured that the generating function of a finitely based permutation class is $D$-finite \cite{NZ1996}, \textit{i.e.,} it is the solution to some differential equation with polynomial coefficient. However, Zeilberger himself later conjectured to the contrary that the generating function counting the permutations avoiding 1324 is not $D$-finite \cite{EV2005}. This, the Noonan-Zeilberger conjecture, was disproved by Garrabrant and Pak in 2015, who give a general method for generating counterexamples \cite{GP2015}.\\

As we briefly mentioned, the generating function for the permutation class avoiding 1324 remains unknown. In fact, this is the only class (up to symmetry) avoiding a single permutation of length 4 that has not been enumerated. To illustrate the difficulty of finding the generating function for this class, Zeilberger stated ``Not even God knows $|\Av_{1000}(1324)|$," where $|\Av_{1000}(1324)|$ is the number of permutations of length 1000 avoiding 1324. On the other hand, Steingr\'{\i}msson disagrees with Zeilberger by saying ``I'm not sure how good Zeilberger's God is at math, but I believe that some humans will
find this number in the not so distant future." Perhaps, our ultimate goal is to prove Steingr\'{\i}msson's claim is correct. The growth rate of this class is also unknown. Currently, the best known lower bound of the growth rate is approximately 9.81, as shown by Bevan \cite{Be2015}, and the best known upper bound is approximately 13.74, as shown by B\'{o}na \cite{Bo2015}.\\

In order to find the generating functions for more permutation classes as well as for the one of permutations avoiding 1324, researchers have worked on establishing many enumeration techniques which can be applied to specific classes. In particular, with the notion of simple permutation introduced in \cite{AA2003,AA2005}, we can find the generating functions for certain classes by enumerating simple permutations in these classes first. Also, authors of \cite{AABRV2013} discovered concrete enumeration methods for geometric griddable classes. In \cite{AB2013}, Albert and Brignall combine these two ideas to find the generating function for the class avoiding permutations 4231, 35142, 42513 and 351624.\\
\section{Place of dissertation in the literature}
There are objects known as Schubert varieties in algebraic geometry, and certain types of Schubert varieties are indexed by permutations in some permutation classes. For instance, as shown by Gasharov and Reiner \cite{GR2002}, permutations in the class Albert and Brignall enumerated in \cite{AB2013} index Schubert varieties defined by inclusions. Also, Lakshmibai and Sandhya \cite{LS1990} showed that smooth Schubert varieties are indexed by permutations avoiding 4231 and 3412, and \'{U}lfarsson and Woo \cite{UW2013} proved local complete intersection ones are indexed by permutations avoiding 52341, 53241, 52431, 35142, 42513 and 351624. The enumeration of the permutations indexing smooth Schubert varieties was initially found by Haiman \cite{H1992}, and then also discussed by Bousquet-M\'{e}lou and Butler \cite{BB2007}. In this dissertation, we enumerate permutations indexing Schubert varieties that are local complete intersections.\\

As discussed in the previous section, three aspects of the study of permutation classes are growth rates, relations between types of generating functions and classes, and enumeration techniques. With this dissertation, we contribute to all of these aspects. Since we prove our final result by extending the methods Albert and Brignall used in \cite{AB2013}, we primarily contribute to the study of enumeration techniques. In the future, we plan to use our result to characterize the permutation classes to which these methods can be applied. In addition, we can see that the generating function we obtain, like the one in \cite{AB2013}, is algebraic. Although we don't have enough evidence to establish any conjecture about the types of generating functions related to classes that can be enumerated by methods we use, there may be some possible connections.\\

In Chapter \ref{chap:2}, we define necessary terminology and provide examples to understand the study of permutation classes. Chapter \ref{chap:3} specifically covers two examples of enumeration. These results are stated as lemmas and will be referred in Chapter \ref{chap:6}. We spend the entirety of Chapter \ref{chap:4} to describe the methods which Albert and Brignall used in \cite{AB2013} in detail. By extending this idea in Chapter \ref{chap:5} and \ref{chap:6}, we complete our final result.\\

For a more detailed history of the study of permutation classes, we refer the reader to the excellent survey of Vatter \cite{V2015}.
\chapter{Definitions and prerequisites}\label{chap:2}
\section{Permutations and permutation classes}
\subsection{Permutations}
A \textit{permutation $\pi$} is a bijective function from $\{1,2,\ldots,n\}$ to itself for some positive integer $n$, which is called the \textit{length of $\pi$}, denoted by $|\pi|$. We call an integer in the domain of $\pi$ a \textit{position}, and an integer in the image of $\pi$ a \textit{value}. In this dissertation, we will write permutations in one-line notation; a permutation $\pi$ will be written as a sequence of $n$ positive integers $\pi_1\ldots\pi_n$, indicating for each $i\in\{1,2,\ldots,n\}$ that $\pi(i)=\pi_i$. (We usually write $\pi(i)$ instead of $\pi_i$ to refer to the value in the $i$-th position of $\pi$). For example, $\pi=316254$ is a permutation of length 6, and $\pi(4)=2$. Let $\mathcal{S}$ and $\mathcal{S}_n$ be the set of all permutations and the set of all permutations of length $n$, respectively. (For $n=0$, $\mathcal{S}_0$ is the set containing only the empty permutation, which is denoted by $\e$). One can easily verify that $|\mathcal{S}_n|$, the number of length $n$ permutations, is $n!$.\\ 

\begin{figure}[H]
\[\begin{tikzpicture}[scale=0.8]
\plotPerm{3,1,6,2,5,4}
\end{tikzpicture}\]
\caption{The graph of the permutation $\pi=316254$.}
\label{fig:2.1}
\end{figure}
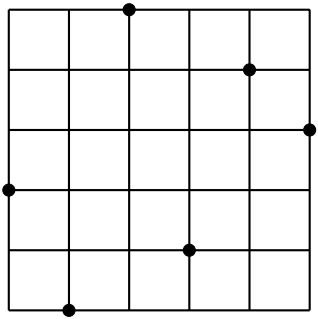

Here, we also introduce the \textit{graph of a permutation}, as this idea will help to illustrate some important concepts. Given a permutation of length $n$, we draw a point at $(i,\pi(i))$ on the Cartesian plane for each $i$ with $1\leq i\leq n$. Note that we always have a unique point on each of $x=1,\ldots,x=n$ and $y=1,\ldots,y=n$ lines. We then draw an $n\times n$ grid containing these points. Figure \ref{fig:2.1} shows the graph of $\pi=316254$.\\
We frequently abuse notation and think of values as points in the graph. For instance, by saying the value $x$ of $\pi$ is located below and to the left of $y$, we mean $x<y$ and $\pi^{-1}(x)<\pi^{-1}(y)$.\\
\subsection{Constructions of new permutations}
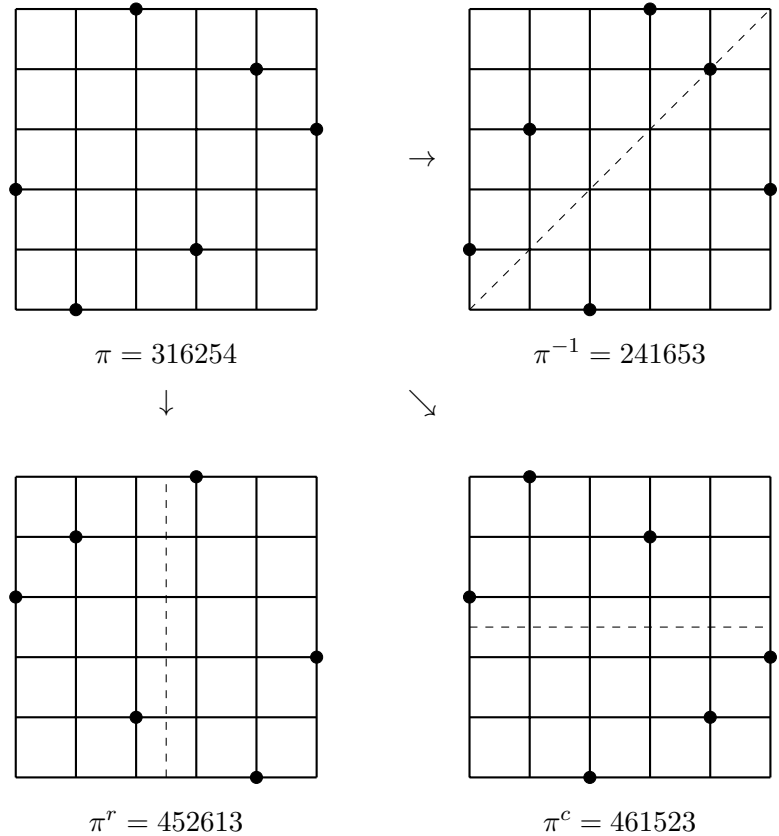
\begin{figure}[H]
\[\begin{array}{ccc}\begin{tikzpicture}[scale=0.8]
\plotPerm{3,1,6,2,5,4}
\end{tikzpicture} &\qquad\raisebox{2cm}{$\rightarrow$}\qquad& \begin{tikzpicture}[scale=0.8]
\plotPerm{2,4,1,6,5,3}
\draw [dashed] (1.0,1.0) -- (6.0,6.0);
\end{tikzpicture}\\
\pi=316254& &\pi^{-1}=241653 \\
\downarrow &\qquad\searrow\qquad& \\
& & \\
\begin{tikzpicture}[scale=0.8]
\plotPerm{4,5,2,6,1,3}
\draw [dashed] (3.5,1.0) -- (3.5,6.0);
\end{tikzpicture} & & \begin{tikzpicture}[scale=0.8]
\plotPerm{4,6,1,5,2,3}
\draw [dashed] (1.0,3.5) -- (6.0,3.5);
\end{tikzpicture}\\
\pi^r=452613& &\pi^c=461523 \\
\end{array}\]
\caption{The graphs of $\pi^{-1}$, $\pi^{r}$ and $\pi^{c}$.}
\label{fig:2.2}
\end{figure}

In this section, we introduce several methods to construct other permutations from a given permutation. Let $\pi\in\S_n$. For every $i$ ($1\leq i\leq n$), the \textit{inverse of $\pi$} is the permutation $\pi^{-1}\in\S_n$ such that $\pi(\pi^{-1}(i))=\pi^{-1}(\pi(i))=i$, the \textit{reverse of $\pi$} is the permutation $\pi^r\in\S_n$ given by $\pi^r(i)=\pi(n+1-i)$, and the \textit{complement of $\pi$} is the permutation $\pi^c\in\S_n$ defined by $\pi^c(i)=n+1-\pi(i)$. Examples of each operation are shown in Figure \ref{fig:2.2} with $\pi=316254$.\\

As we can observe, the inverse, the reverse, and the complement of a permutation $\pi$ can be obtained by reflecting $\pi$ respectively over SW-NE diagonal, vertical, and horizontal lines on the graph of $\pi$. We denote by $\textrm{Sym}(\pi)$ the set of permutations obtained by applying compositions of these operations and call this set the \textit{symmetry class of $\pi$}. In other words, $\textrm{Sym}(\pi)$ is the orbit of the dihedral group of order 8 acting on the graph of $\pi$. For $\pi= 316254$, one can easily verify that $\textrm{Sym}(316254)=\{316254,241653,452613,461523,421635,356142,325164,536124\}$.\\

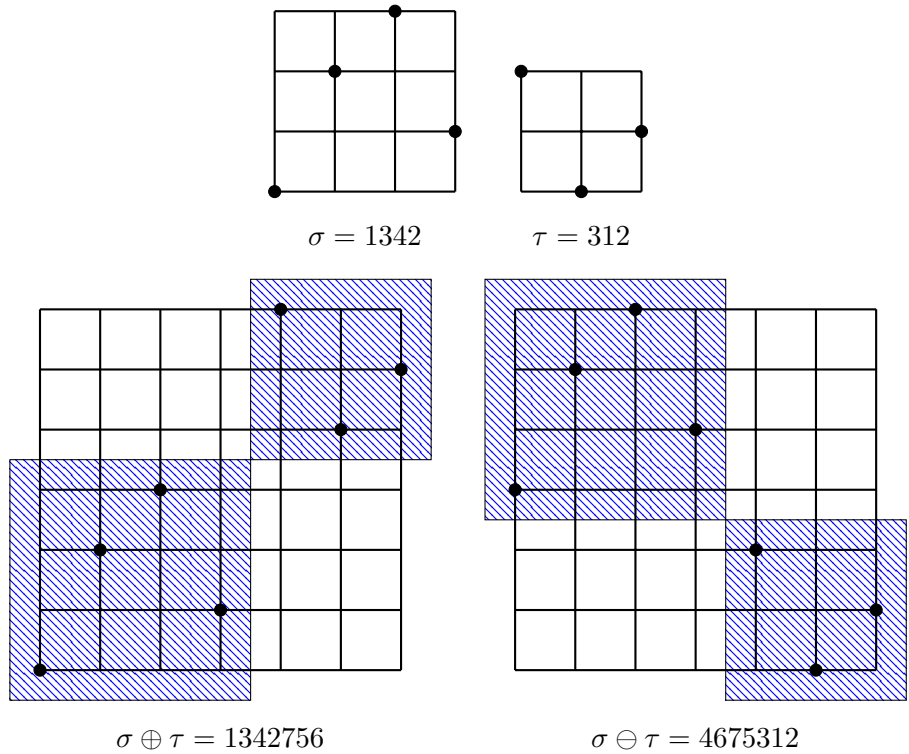
\begin{figure}[H]
\[\begin{array}{ccc}\begin{tikzpicture}[scale=0.8]
\plotPerm{1,3,4,2}
\end{tikzpicture} &\qquad & \begin{tikzpicture}[scale=0.8]
\plotPerm{3,1,2}
\end{tikzpicture}\\
\s=1342 & &\tau=312
\end{array}\]
\[\begin{array}{ccc}
\begin{tikzpicture}[scale=0.8]
\draw[pattern=north west lines, pattern color=blue](0.5,0.5)--(0.5,4.5)--(4.5,4.5)--(4.5,0.5)--cycle;
\draw[pattern=north west lines, pattern color=blue](4.5,4.5)--(4.5,7.5)--(7.5,7.5)--(7.5,4.5)--cycle;
\plotPerm{1,3,4,2,7,5,6}
\end{tikzpicture} &\qquad & \begin{tikzpicture}[scale=0.8]
\draw[pattern=north west lines, pattern color=blue](0.5,3.5)--(0.5,7.5)--(4.5,7.5)--(4.5,3.5)--cycle;
\draw[pattern=north west lines, pattern color=blue](4.5,0.5)--(4.5,3.5)--(7.5,3.5)--(7.5,0.5)--cycle;
\plotPerm{4,6,7,5,3,1,2}
\end{tikzpicture}\\
\s\oplus\tau=1342756 & &\s\ominus\tau=4675312
\end{array}\]
\caption{The graphs of $\s\oplus\tau$ and $\s\ominus\tau$.}
\label{fig:2.3}
\end{figure}

The previous three operations construct a new permutation of the same length from a given permutation. Next, we introduce two ways to ``glue" two permutations together to construct a lager permutation. Given $\s\in\S_m$ and $\tau\in\S_n$, the \textit{sum of $\s$ and $\tau$} is the permutation defined by\[
\s\oplus\tau=\s(1)\s(2)\cdots\s(m)\tau'(1)\tau'(2)\cdots\tau'(n)
\]
where $\tau'(i)=\tau(i)+m$ for each $i$ with $1\leq i\leq n$. Similarly, we define the \textit{skew-sum of $\s$ and $\tau$} to be the permutation\[
\s\ominus\tau=\s'(1)'\s(2)\cdots\s'(m)\tau(1)\tau(2)\cdots\tau(n)
\]
where $\s'(i)=\s(i)+n$ for each $i$ with $1\leq i\leq n$. Figure \ref{fig:2.3} shows the sum and skew-sum of $\s=1342$ and $\tau=312$. Note that neither of these sum operations is commutative, but they are associative. For instance, while $\s\oplus\tau=1342756$, we have $\tau\oplus\s=3124675$.\\
If $\pi$ can be constructed by a sum $\s_1\oplus\s_2$ for some nonempty permutations $\s_1$ and $\s_2$, then we say $\pi$ is \textit{sum-decomposable}. Similarly, if $\pi=\s_1\ominus\s_2$ for some nonempty permutations $\s_1$ and $\s_2$, then we say $\pi$ is \textit{skew-decomposable}. More importantly, if $\pi$ is not sum-decomposable (respectively skew-decomposable), \textit{i.e.}, there does not exist two permutations $\s_1$ and $\s_2$ such that $\pi=\s_1\oplus\s_2$ (respectively $\pi=\s_1\ominus\s_2$), then we say $\pi$ is \textit{sum-indecomposable} (respectively \textit{skew-indecomposable}).\\

Lastly, there is another method, called inflation, to construct a new permutation that will be important in this dissertation, but it makes more sense to introduce this construction when we discuss simple permutations, so we postpone discussing it until Section 1.3.\\
\subsection{Permutation avoidance and permutation classes}
We now introduce the concept of permutation avoidance. The \textit{flattening} of a sequence of $n$ distinct positive real numbers is the unique permutation $\s$ of length $n$ where $\s$ has the same relative order as the sequence. For example, the flattening of $364$ is $\s=132$. If a permutation $\pi$ has a subsequence (not necessarily consecutive) whose flattening is $\s$, we say the \textit{pattern $\s$ is contained in $\pi$} and write $\s\cont\pi$, since pattern containment is a partial order. Notice that the permutation $\pi=316254$ has a subsequence $364$, so $\s=132$ is contained in $\pi$. More importantly, if $\s$ is not contained in $\pi$, we say \textit{$\pi$ avoids $\s$ pattern} and write $\s\avd\pi$. With the same permutation $\pi$ for example, $\pi$ avoids $\tau=4231$ because $\pi$ does not contain a subsequence of four positive integers whose flattening is $\tau$.\\

We can visualize the concepts of containment and avoidance with graphs of permutations. In the graph of $\pi=316254$, if we take the points corresponding to the subsequence $364$ and disregard all the others, we can construct the graph of $\s=132$ by removing the lines with no points and squeezing together the lines that remain into a $3\times 3$ grid, as shown in Figure \ref{fig:2.4}. On the other hand, there is no subset of points that we can choose to construct $\tau=4231$ in this way, and hence, $\tau\avd\pi$.\\
\begin{figure}[b!]
\[\begin{array}{ccc}\begin{tikzpicture}[scale=0.8]
\plotPerm{3,1,6,2,5,4}
\draw[color=black] (1,3) circle [radius=.25];
\draw[color=black] (3,6) circle [radius=.25];
\draw[color=black] (6,4) circle [radius=.25];
\draw[black] (1.75,0.75) -- (2.25,1.25);
\draw[black] (1.75,1.25) -- (2.25,0.75);
\draw[black] (3.75,1.75) -- (4.25,2.25);
\draw[black] (3.75,2.25) -- (4.25,1.75);
\draw[black] (4.75,4.75) -- (5.25,5.25);
\draw[black] (4.75,5.25) -- (5.25,4.75);
\end{tikzpicture} &\qquad\raisebox{2cm}{$\rightarrow$}\qquad& \raisebox{1.25cm}{\begin{tikzpicture}[scale=0.8]
\plotPerm{1,3,2}
\end{tikzpicture}}
\end{array}\]
\caption{The graph describing $132\cont 316254$.}
\label{fig:2.4}
\end{figure}

Define a \textit{permutation class} to be a set $\C$ of permutations with the property that, if $\pi\in\mathcal{C}$ and $\s\cont\pi$, then $\s\in\mathcal{C}$. The set of permutations $\b$ that are minimal (with respect to containment) among those not in $\C$ is called the \textit{basis of $\C$}. Let $B$ be a set of permutations and $\textrm{Av}(B)$ be the set of permutations avoiding every permutation in $B$. We call this set the \textit{permutation class of $B$}. Note that if $B$ is the basis of $\mathcal{C}$, then we have\[
\mathcal{C}=\Av(B)=\{\pi:\b\avd\pi\textrm{ for all }\b\in B\}.
\]

A trivial example is $\mathcal{C}=\textrm{Av}(21)$ (we usually omit superfluous braces). It is a straightforward exercise to see that for each positive integer $n$, $12\ldots n$ is the only permutation in $\mathcal{S}_n$ that avoids $21$. Thus, $\textrm{Av}(21)=\{\e,1,12,123,1234,\ldots\}$.\\

We conclude this section by introducing some basic properties of permutation classes. The definition of a permutation class immediately implies the following proposition.
\begin{proposition}\label{prop:2.1}
\textit{Let $B_1$ and $B_2$ be distinct sets bases of permutation classes. If, for every $\b_2\in B_2$, there exists $\b_1\in B_1$ such that $\b_1\cont\b_2$, then $\Av(B_1)\subsetneq\Av(B_2)$.}
\end{proposition}
\vhhh
For example, $\Av(132,213)\subsetneq\Av(132,3142)$ since $132\cont 132$ and $213\cont 3142$.\\

The following result is less obvious but not too difficult to prove.
\begin{proposition}\label{prop:2.2}
\textit{Let $\textrm{op}$ be a fixed composition of the inverse, reverse and complement operations. Given a basis $B$, define $op(B)=\{op(\b):\b\in B\}$. Then $\pi\in\Av(B)$ if and only if $op(\pi)\in\Av(op(B))$. Hence, $|\Av(B)|=|\Av(\textrm{op}(B))|$}.\end{proposition}
\section{Generating functions}
Given an arbitrary class $\C$, we are interested in discovering a formula called a generating function which expresses the number of length $n$ permutations in $\C$, if possible. In this section, we define the notion of generating function for an integer sequence and discuss a few enumerative results which have been found in recent years.\\

Suppose we have a sequence of positive integers $\{a_n\}_{n=0}^\infty$. The \textit{generating function for $\{a_n\}_{n=0}^\infty$} is a formal power series of the form\[
f=\sum_{n=0}^\infty a_nx^n=a_0+a_1x+a_2x^2+a_3x^3+\cdots.
\]
In particular, a generating function for a permutation class $\mathcal{C}$ is a formal power series $f_\C$ defined by\[
f_\C=\sum_{n=0}^\infty s_n(\C) x^n,
\]
where $s_n(\C)=|\S_n\cap\C|$ (or simply $s_n$, if it is clear from context). In other words, for each $i\geq 0$, the coefficient of $x^i$ is the number of length $i$ permutations in $\C$. For instance, the generating function for the class $\Av(21)$ is\\
\[
f_{\Av(21)}=\frac{1}{1-x}=1+x+x^2+x^3+x^4+\cdots.
\]

It is sometimes convenient (and necessary) to exclude the constant term $1$, which corresponds to the empty permutation $\e$. In this dissertation, we let $\bar{f}=f-1$.\\

We will also need the concept of multivariate generating functions. The \textit{multivariate generating function} with $i$ variables for a sequence of positive integers $\{a_{n_1,n_2,\ldots,n_i}\}_{n_1,n_2,\ldots,n_i\geq 0}$ is a formal power series of the form\[
f=\sum_{n_1,n_2,\ldots,n_i\geq 0}a_{n_1,n_2,\ldots,n_i}x_1^{n_1}x_2^{n_2}\cdots x_{i}^{n_i}.
\]

Once we discover the generating function for a permutation class, we are often able to find the growth rate of the class. As the Stanley-Wilf conjecture states, any permutation class has exponential growth. Hence, by applying the Ratio Test and Taylor's Theorem from complex analysis, the growth rate is the reciprocal of the distance from 0 to the closest pole (\textit{i.e.} the radius of convergence). In particular, $\textrm{gr}(\C)$, the growth rate of a permutation class $\C$, is obtained by\\
\[
\textrm{gr}(\C)=\lim_{n\to\infty}\frac{|a_{n+1}|}{|a_n|}=\frac{1}{R},
\]
\\
where $R$ is the radius of convergence of the generating function for $\C$.\\

There are many different techniques for finding generating function of permutation classes. We postpone examples of finding generating functions until Chapter \ref{chap:3}, where we show how to enumerate two classes, $\Av(123,213,132)$ and $\Av(4123,4213,4132)$. In the remainder of this section, we will discuss some historical enumerative results.\\

Since the enumerative study of permutation classes blossomed in 1980s, generating functions for various classes have been discovered. The first well known result is the following.\vhhh
\begin{theorem}\label{thm:2.1}
\textit{For every permutation $\b$ in $\S_3$, the number of permutations of length $n$ in $\Av(\b)$ is the $n$-th Catalan number. Hence,}\\
\[
f_{\Av(\b)}=\frac{1-\sqrt{1-4x}}{2x}=1+x+2x^2+5x^3+14x^4+42x^5+\cdots
\]
\end{theorem}\vhh
For $\b=123$ and $\b=231$, the above result was classically known \cite{M1915,K1973}. In \cite{SS1985}, Simion and Schmidt summarize this result as well as various results for bases having more than one permutation of length 3, including $\Av(123,213,132)$, which we introduce in Chapter \ref{chap:3}.\\

In the 1990s, researchers discovered numerous results with bases containing permutations of length 4 \cite{G1990,B1997}. In \cite{B1997}, B\'{o}na not only finds the generating function for $\Av(1342)$, but furthermore finds the exact formula for $s_n(1342)$ for every $n$ by giving a bijection between permutations in $\Av(1342)$ and labeled plane trees of a certain type on $n$ vertices. The class $\textrm{Av}(3412,4231)$, the set of permutations indexing smooth Schubert varieties, was also enumerated in the early 1990s by Haiman \cite{H1992}; this result was also discussed in \cite{BB2007}. Other various results are also discussed in \cite{AAV2009,AAB2011}, and there are many more.\\

In addition to these results, researchers have recently established some concrete enumeration techniques which can be applied to certain types of classes. In \cite{B2010}, Brignall summarizes a variety of techniques coming from simple permutations (introduced in the next section). Also, \cite{AABRV2013} defines so-called geometric grid classes and introduces some techniques that can be used to enumerate this type of class. In Chapter 3, we investigate the method Albert and Brignall use to enumerate the class $\A$ in \cite{AB2013}.\\
\section{Simple permutations}
\subsection{Definition}
We now move onto the discussion of simple permutations. First, we define two kinds of intervals. Given a permutation $\pi$ of length $n$, a \textit{segment} of $\pi$ is a set of consecutive positions $i,i+1,\ldots, j$ in $\pi$, and a \textit{range} of $\pi$ is a set of consecutive values $a,a+1,\ldots, b$ of $\pi$. We use the standard notations of intervals for them. Although permutations only contain positive integer values, we may sometimes use open intervals to exclude boundary positions or values. For these intervals, we carry the notion of length from a permutation to denote the number of elements in $[i,j]$, which is simply $j-i+1$. We call a segment and a range that are not $[1,n]$ \textit{proper}.\\

Now, we define simple permutations. Let $[i,j]$ be a segment of $\pi$ whose length is $m$. If the set of values $\{\pi(i),\ldots,\pi(j)\}$ is an interval $[a,b]$, then $[i,j]$ is called a \textit{block of $\pi$}. In this case, we may denote by $\pi([i,j])$ the corresponding range to the segment $[i,j]$ forming a block. Every permutation of length $n$ has $n$ singleton blocks as well as the block $[1,n]$. If a permutation $\pi$ only contains these blocks, $\pi$ is called \textit{simple}. By convention, $\pi=1$ is not considered simple. As an example and a non-example of simple permutations, suppose we have $\pi=25314$ and $\s=4127563$. Notice that $\s$ contains segments $[2,3]$ and $[4,6]$ which are mapped to ranges $[1,2]$ and $[5,7]$ by $\pi$ respectively. Thus, they are blocks of length 2 and 3, so $\s$ is not simple. On the contrary, since the only blocks $\pi$ contains are singletons and the block $[1,5]$, $\pi$ is a simple permutation. Given a permutation class $\C$, we denote by $\textrm{Si}(\C)$ the set of simple permutations in $\C$.\\
\subsection{Inflation}
An alternative definition of simple permutations can be given by introducing the following method to construct a new permutation.\\

Let $\pi$ be a permutation of length $n$ and $\s_1,\s_2,\ldots,\s_n$ be $n$ non-empty permutations of various lengths. Denote by $i_j$ the length of $\s_j$ for every $j$ with $1\leq j\leq n$, and let $i_0=0$. Finally, define intervals $I_j=[i_0+i_1+\cdots+i_{j-1}+1,i_0+i_1+\cdots+i_j]$ for each $j$ with $1\leq j\leq n$. The \textit{inflation of $\pi$ by $\s_1,\s_2,\ldots,\s_n$} is the permutation $\a$ of length $i_1+\cdots+i_n$, denoted by $\pi[\s_1,\s_2,\ldots,\s_n]$, where\begin{enumerate}
\item For every $j$ ($1\leq j\leq n$), there exists a constant $t$ such that\[
(\a(i_0+\cdots+i_{j-1}+1)-t)\cdots(\a(i_0+\cdots+i_j)-t)=\s_j.
\]
\item For distinct $j$ and $k$ ($1\leq j,k\leq n$), if $\l\in I_j$, $m\in I_k$ and $\pi(\l)\leq\pi(m)$, then $\a(\l)\leq\a(m)$.
\end{enumerate}
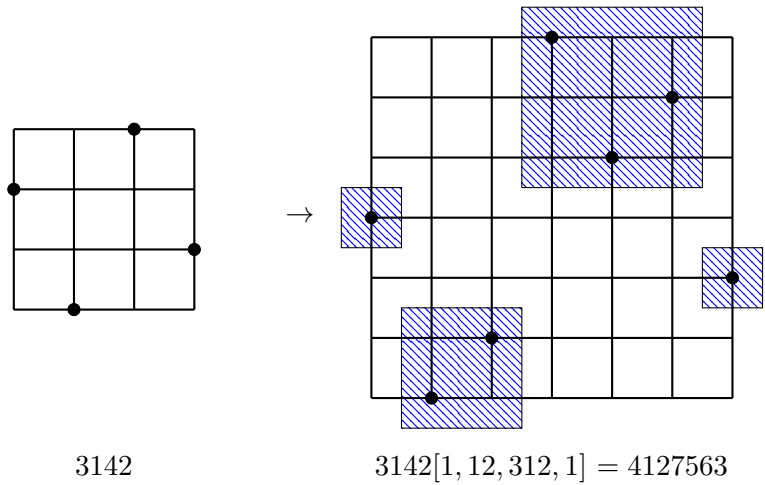
\begin{figure}[b!]
\[\begin{array}{ccc}\raisebox{1.5cm}{\begin{tikzpicture}[scale=0.8]
\plotPerm{3,1,4,2}
\end{tikzpicture}} &\qquad\raisebox{2.75cm}{$\rightarrow$}\qquad& \begin{tikzpicture}[scale=0.8]
\draw[pattern=north west lines, pattern color=blue](0.5,3.5)--(1.5,3.5)--(1.5,4.5)--(0.5,4.5)--cycle;
\draw[pattern=north west lines, pattern color=blue](1.5,0.5)--(3.5,0.5)--(3.5,2.5)--(1.5,2.5)--cycle;
\draw[pattern=north west lines, pattern color=blue](3.5,4.5)--(6.5,4.5)--(6.5,7.5)--(3.5,7.5)--cycle;
\draw[pattern=north west lines, pattern color=blue](6.5,2.5)--(7.5,2.5)--(7.5,3.5)--(6.5,3.5)--cycle;
\plotPerm{4,1,2,7,5,6,3}
\end{tikzpicture}\\
3142 & & 3142[1,12,312,1]=4127563
\end{array}
\]
\caption{The inflation of 3142 by $1$, $12$, $312$ and $1$.}
\label{fig:2.5}
\end{figure}
In other words, $\a=\pi[\s_1,\s_2,\ldots,\s_n]$ is the permutation obtained by replacing each $\pi(j)$ with a block whose flattening is $\s_j$, so that for every distinct $j$ and $k$, the relative ordering of $\a(\l)$ and $\a(m)$ is the same as the relative ordering of $\pi(\l)$ and $\pi(m)$. An example of inflation is shown in Figure \ref{fig:2.5}, as this concept is best illustrated using graphs.\\

If the only ways to obtain a permutation $\pi$ by inflation are $\pi[1,1,\ldots,1]$ and $1[\pi]$, then $\pi$ is simple. Hence, as discussed previously and shown in Figure \ref{fig:2.5}, $\s=4127563$ is not simple, whereas $\pi=25314$ is simple.\\
\subsection{The importance of simple permutations}
The notion of simple permutations was introduced in \cite{AA2003}, where the authors describe how this notion is useful in the study of various permutation classes. The most essential proposition in \cite{AA2003} for our purpose is the following.
\begin{proposition}\label{prop:2.3}
\textit{(Albert and Atkinson, 2005) For every permutation $\a$, there exists a simple permutation $\pi$ of some length $n$ and permutations $\s_1,\s_2,\ldots,\s_n$ such that}\[
\a=\pi[\s_1,\s_2,\ldots,\s_n].
\]
\textit{Furthermore, if $\pi\neq 12,21$, then the permutation $\pi$ is uniquely determined by $\a$. If $\pi=12$ or 21, then $\s_1,\s_2$ are uniquely determined so long as we require that $\s_1$ is sum-indecomposable or skew-indecomposable respectively.}
\end{proposition}

In other words, any permutation $\a$ can be decomposed as the inflation of a unique simple permutation $\pi$ by $\s_1,\ldots,\s_n$. In this case, $\pi$ is called the \textit{skeleton of $\a$}. As stated, we must ensure that $|\pi|\neq 2$ (so $|\pi|\geq 4$, since there are no simple permutations of length 3) or, equivalently, that $\a$ is sum/skew-indecomposable, to say $\s_1,\ldots,\s_n$ are also uniquely determined. For example, $1234=12[1,123]=12[12,12]=12[123,1]$. We can enforce uniqueness for sum/skew-decomposable permutations by insisting $\s_1$ be sum-indecomposable, so $1234=12[1,123]$. Hence with Proposition \ref{prop:2.3} being analogous to the Fundamental Theorem of Arithmetic, simple permutations are to permutations as prime numbers are to integers. This essential idea often helps in discovering more general methods for finding generating functions for certain classes. In particular, we have the following propositions.
\begin{proposition}\label{prop:2.4}
\textit{Let $\C$ be a permutation class, $\pi$ be a simple permutation of length $n\geq 4$ in $\C$ and $\textrm{ifl}_\C(\pi)$ be the set of permutations in $\C$ which can be inflated from $\pi$. Suppose there exist subclasses $\C_i$ of $\C$ (for each $i$, $1\leq i\leq n$) such that $\pi[\s_1,\ldots,\s_n]\in\C$ if and only if $\s_i\in\C_i$. Then the generating function for $\textrm{ifl}(\pi)$ is}
\[
f_{\textrm{ifl}_\C(\pi)}=\prod_{i=1}^n \bar{f}_{\C_i}.
\]
\textit{If $\pi=12$ (respectively $\pi=21$) satisfies the above hypothesis, then}
\[
f_{\textrm{ifl}_\C(12)}=\bar{f}_{C_1}^{\oplus}\cdot\bar{f}_{C_2}\qquad(\textrm{respectively }f_{\textrm{ifl}_\C(21)}=\bar{f}_{C_1}^{\ominus}\cdot\bar{f}_{C_2})
\]
\textit{where $\bar{f}_{C_1}^{\oplus}$ (respectively $\bar{f}_{C_1}^{\ominus}$) is the generating function for sum-indecomposable (respectively skew-indecomposable) permutations in $C_1$, excluding the empty permutation.}
\end{proposition}
\begin{proposition}\label{prop:2.5}
\textit{Let $\C$ be a permutation class. If Proposition \ref{prop:2.4} is applicable for all simple permutations in $\C$, then the generating function for $\C$ is}\[
f_\C=\sum_{\pi\in\textrm{Si}(\C)}f_{\textrm{ifl}_\C(\pi)}.
\]
\end{proposition}

Proposition \ref{prop:2.4} states that if every set of choices for each $\s_i$ forms a subclass of $\C$ independently, then by the combinatorial meaning of multiplication together with the definition of generating functions, the generating function for $\textrm{ifl}_\C(\pi)$ can be obtained by multiplying each $\bar{f}_{\C_i}$. If the basis of a class $\C$ contains only simple permutations, then $\C_i=\C$ for all $i$, so the hypothesis for Proposition \ref{prop:2.4} is automatically satisfied.\\

If Proposition \ref{prop:2.4} is applicable to every simple permutation in $\C$, then by the combinatorial meaning of addition and the definition of generating functions, Proposition \ref{prop:2.5} is an immediate consequence.\\

In \cite{AA2005}, authors also discuss the asymptotic result stated in Theorem \ref{thm:2.3} as well as Theorem \ref{thm:2.4}, which is a strong result about generating functions for permutation classes containing finitely many simple permutations.\vhhh
\begin{theorem}\label{thm:2.3}
\textit{(Albert and Atkinson, 2003 \cite{AA2003}) Let $p_n$ be the number of simple permutations of length $n$. Then}
\[
p_n=\frac{n!}{e^2}\left(1-\frac{4}{n}+\frac{2}{n(n-1)}+O(n^{-3})\right).
\]
\end{theorem}
\vhh
\begin{theorem}\label{thm:2.4}
\textit{(Albert and Atkinson, 2005 \cite{AA2005}) If a permutation class contains only finitely many simple permutations, then it has a finite basis and an algebraic generating function. Furthermore, if such a class does not contain the permutation $n(n-1)\cdots 321$ for some $n$, then it has a rational generating function.}
\end{theorem}

Theorem \ref{thm:2.4} is not applicable to $\A$ and $\A'$, the classes we enumerate in this dissertation, since they have infinitely many simple permutations. However, the method we use for them heavily depends on the structure of simple permutations in these classes.\\
\section{Automata and the transfer matrix method}
We conclude this chapter with the discussion of elementary automata theory and the so-called transfer matrix method. Given a digraph with finitely many vertices and edges, the transfer matrix method allows us to find the generating function (according to some weight function on the edges) for paths from a specified vertex to another. We can apply this to the state diagram of an automaton to find the generating function for the language it accepts. This will be the key to finding $f_\A$ and $f_{\A'}$.\\
\subsection{Definition and example}
We start with the definition of a deterministic finite-state automaton. An \textit{alphabet} $\Sigma$ is a finite set, and we call the elements of the alphabet \textit{letters}. For example, $\Sigma=\{a,b,c\}$ is an alphabet. A string of letters $\a_1,\a_2,\ldots,\a_n$ where $\a_i\in\Sigma$ ($1\leq i\leq n$) is called a \textit{word}. The string with no letters is called the \textit{empty word}, denoted by $\lambda$. The set of all words associated with $\Sigma$, including the empty word, is denoted by $\Sigma^*$. In particular,\[
\Sigma^*=\bigcup_{i=0}^\infty\Sigma^i,\quad\textrm{ where }\Sigma^i=\underbrace{\Sigma\times\cdots\times\Sigma}_{i\textrm{ times}}.
\]
With $\Sigma=\{a,b,c\}$, we have $\Sigma^*=\{\lambda,a,b,c,aa,ab,ac,ba,bb,bc,ca,cb,cc,aaa,\ldots\}$. Any subset of $\Sigma^*$ is called a \textit{language}.\\

A \textit{deterministic finite-state automaton} is a 5-tuple $M=(Q,\Sigma,\d,q_0,F)$ where:\begin{itemize}
\item $Q$ is a finite set. Elements in $Q$ are called \textit{states} and denoted by $q$ with some subscript.
\item $\Sigma$ is an alphabet.
\item $\d:Q\times\Sigma\to Q$ is a function called the \textit{transition function}.
\item $q_0\in Q$ is called the \textit{initial state}.
\item $F$ is a subset of $Q$. Any state in $F$ is called a \textit{accept state}.
\end{itemize}

\begin{figure}[b!]
\[\begin{array}{ccc}\begin{tikzpicture}[->,>=stealth',shorten >=1pt,auto,node distance=2.8cm,semithick,initial/.style={}]
  \tikzstyle{every state}=[circle,fill=black!25,minimum size=17pt,inner sep=0pt]
  \node[state,initial]   (A)                    {$A$};
  \node[state]           (B) [above right of=A] {$B$};
  \node[state,accepting] (C) [below right of=B] {$C$};
  \node[state]           (J) [below right of=A] {$J$};
  \path (A) edge [bend left]  node {$b$} (B)
            edge [bend right] node[below left] {$c$} (C)
            edge              node[below left] {$a$} (J)
        (B) edge              node {$a$} (A)
            edge              node[below left] {$c$} (C)
            edge              node {$b$} (J)
        (C) edge [bend right] node[above right] {$b$} (B)
            edge [loop right] node {$c$} (C)
            edge              node {$a$} (J)
        (J) edge [loop below] node {$a,b,c$} (J);
  \draw[<-] (A) -- node[above] {} ++(-1cm,0);
\end{tikzpicture} & \qquad\qquad &
\raisebox{1.65cm}{\begin{tikzpicture}[->,>=stealth',shorten >=1pt,auto,node distance=2.8cm,semithick,initial/.style={}]
  \tikzstyle{every state}=[circle,fill=black!25,minimum size=17pt,inner sep=0pt]
  \node[state,initial]   (A)                    {$A$};
  \node[state]           (B) [above right of=A] {$B$};
  \node[state,accepting] (C) [below right of=B] {$C$};
  \path (A) edge [bend left]  node {$b$} (B)
            edge              node {$c$} (C)
        (B) edge              node {$a$} (A)
            edge              node {$c$} (C)
        (C) edge [bend right] node[above right] {$b$} (B)
            edge [loop right] node {$c$} (C);
  \draw[<-] (A) -- node[above] {} ++(-1cm,0);
\end{tikzpicture}}\\
(a) & & (b)\end{array}\]
\caption{The state diagram of the example automaton.}
\label{fig:2.6}
\end{figure}
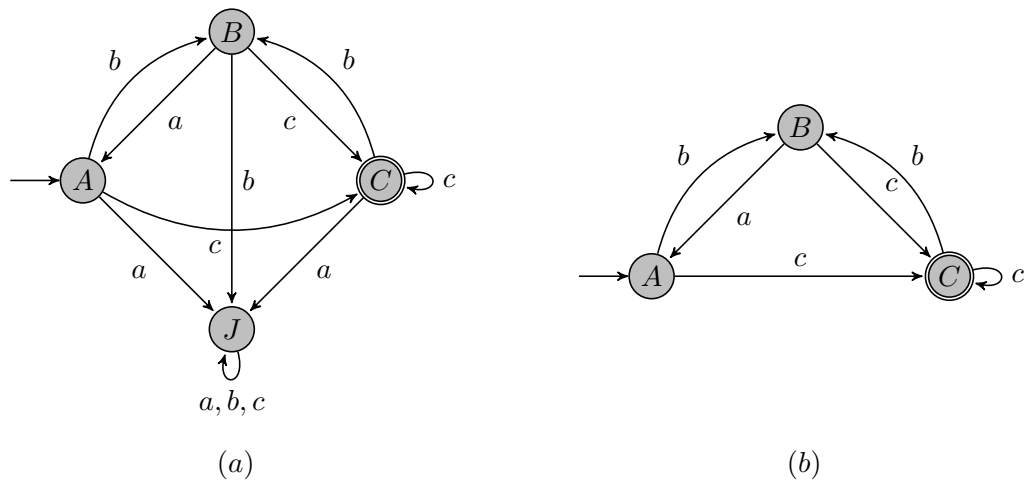

As an example, let $Q=\{A,B,C,J\}$ with $A$ being the initial state, $F=\{C\}$, $\Sigma=\{a,b,c\}$ and $\d$ defined by\\
\[
\d(q,\a)=\left\{\begin{array}{cl}
A & \textrm{ if }(q,\a)=(B,a)\\
B & \textrm{ if }(q,\a)=(A,b),(C,b)\\
C & \textrm{ if }(q,\a)=(A,c),(B,c),(C,c)\\
J & \textrm{ otherwise}
\end{array}\right..
\]\\
Together, we have a deterministic finite-state automaton $M=(Q,\Sigma,\d,A,F)$. We can provide a graphical representation of an automaton called the \textit{state diagram}, which is an edge-labelled directed graph. States are represented by vertices. If $\d(q_i,\a)=q_j$, we draw a directed edge from $q_i$ to $q_j$ which is labeled with $\a$. Such a directed edge is called a \textit{transition} from $q_i$ to $q_j$. The initial state is indicated by the arrow with no label, and accept states are shown by double circles. For the example above, we have the state diagram shown in Figure \ref{fig:2.6}$(a)$.\\

Associated to each automaton is a language $\mathcal{L}(M)$. Let $w=\a_0,\a_1,\ldots,\a_n$ be a word in $\Sigma^*$. The \textit{run} of an automaton on $w$ is a sequence of states $q_0,q_1,\ldots,q_n$ (not necessarily distinct) where $q_0$ is the initial state, and $q_i=\d(q_{i-1},\a_i)$ for every $i$ with $1\leq i\leq n$. A word $w$ is said to be \textit{accepted} by an automaton $(Q,\Sigma,\d,q_0,F)$ if $q_n\in F$. In other words, an automaton reads letters in $w$ in sequence. As it reads the letter $\a_i$, it moves to the state $q_i$. If the final state $q_n$ is in $F$, then automaton accepts $w$. The language $\mathcal{L}(M)$ is the set of all words accepted by the automaton.\\

In our example automaton, once we arrive at the state $J$, the transition function always gives us $J$ thereon and we are no longer able to go to any other states and in particular any accept states. We call such a state a \textit{jail state} and often omit all transitions to jail states as well as jail states themselves from the state diagram. Thus, we may represent the diagram of the above example as Figure \ref{fig:2.6}$(b)$ instead.\\

A language $K$ is said to be \textit{regular} if there exists a deterministic finite-state automaton $M$ such that $\mathcal{L}(M)=K$. Note that in most textbooks, a regular language is defined as an element of the set of languages $\mathcal{R}$ over $\Sigma$ which is defined recursively as the following.\begin{itemize}
\item $\emptyset,\{\lambda\}\in\mathcal{R}$ and for every $\a\in\Sigma$, $\{\a\}\in\mathcal{R}$.
\item If $K,L\in\mathcal{R}$, then $K\cup L,KL=\{\underbrace{\a_1,\ldots,\a_m}_{w_K},\underbrace{\b_1,\ldots,b_n}_{w_L}:w_K\in K,w_L\in L\},K^*\in\mathcal{R}$.
\item $\mathcal{R}$ is the minimal set satisfying above two conditions.
\end{itemize}
A classical theorem states that the definition we provide and the definition above are equivalent.\\
\subsection{Transfer matrix method}
Here, we introduce a useful application of an automaton. With the method we discuss in this section, we are able to find the generating function giving the number of $n$-letter words that are accepted by a given automaton. To start, we define a weight function on a digraph.\\

Let $D$ be a finite digraph with a vertex set $V$ and an edge set $E$. A \textit{weight function on $E$} is a mapping from $E$ to some commutative ring $R$. For an arbitrary walk $\Gamma=e_1 e_2\cdots e_n$, the \textit{weight of $\Gamma$ with respect to $w$} is defined by $w(e_1)w(e_2)\cdots w(e_n)$. Intuitively, we may think that every time we pass through an edge, we ``count" it by  multiplying by the assigned weight.\\

Given a digraph $D=(V,E)$ with $|V|=m$ and a weight function on $D$, an \textit{adjacency matrix $P$ of $D$ with respect to $w$} is the $m\times m$ matrix where row $i$ ($1\leq i\leq m$) and column $j$ ($1\leq i\leq m$) are labeled by vertices $v_i$ and $v_j$ respectively, and each entry is\[
P_{ij}=\sum_e w(e)
\]
where the sum is over all edges from $v_i$ to $v_j$.\\

As an example, we take the state diagram in Figure \ref{fig:2.6}$(a)$. Let $w$ be the weight function defined by $w(t)=x$ for every transition $t$. Then we obtain the following adjacency matrix.\[
P=\kbordermatrix{ & A & B & C & J\\
A & 0 & x & x & x\\
B & x & 0 & x & x\\
C & 0 & x & x & x\\
J & 0 & 0 & 0 & 3x}
\]
\\
Again, we are not interested in jail states for our purpose, so we omit the row and column designated for $J$ and express the adjacency matrix of this example as the following.\begin{equation}
P=\kbordermatrix{ & A & B & C\\
A & 0 & x & x\\
B & x & 0 & x\\
C & 0 & x & x}
\end{equation}

With this definition of adjacency matrix, we derive the following theorem.\vhh
\begin{theorem}\label{thm:2.5}
\textit{Let $P$ be the adjacency matrix of some digraph $D=(V,E)$, where $|V|=m$, with respect to a weight function $w$. For any positive integer $n$ and any $i$ and $j$ with $1\leq i,j\leq m$, the $(i,j)$-entry of $P^n$ is}\[
(P^n)_{ij}=\sum_{\Gamma}w(\Gamma)
\]
\textit{where the sum is over all walks $\Gamma$ in $D$ of length $n$ from $v_i$ to $v_j$. (By convention, we define $P^0=I$ even if $P$ is not invertible.)}
\end{theorem}
\vhhh
\textit{Proof.} The proof is by induction on $n$. The base case with $n=1$ is obvious by the definition of adjacency matrix. Assume the statement is true for some positive integer $k$. For any $i,j$ with $1\leq i,j\leq m$, the $(i,j)$-entry of $P^{k+1}$ is obtained by\[
(P^{k+1})_{ij}=(P^k)_{i1}P_{1j}+(P^k)_{i2}P_{2j}+\cdots+(P^k)_{im}P_{mj}=\sum_{1\leq\l\leq m}(P^k)_{i\l}P_{\l j}
\]

By the inductive hypothesis, for each $\l$ with $1\leq\l\leq m$, $(P^k)_{i\l}=\sum_{\Delta}w(\Delta)$ where the sum is over all walks $\Delta$ in $D$ of length $k$ from $v_i$ to $v_\l$. By multiplying $P_{\l j}=\sum_e w(e)$ where the sum is over all edges from $v_\l$ to $v_j$, we have\[
(P^k)_{i\l}P_{\l j}=\sum_{\Delta}w(\Delta)\sum_e w(e)= \sum_e\sum_{\Delta}w(\Delta)w(e)=\sum_{\Gamma_\l}w(\Gamma_\l)
\]
where the sum is over all walks $\Gamma\l$ in $D$ of length $k+1$ from $v_i$ to $v_j$, such that the $k$-th vertex is $v_\l$. Hence,\[
(P^{k+1})_{ij}=\sum_{1\leq\l\leq m}(P^k)_{i\l}P_{\l j}=\sum_{1\leq\l\leq m}\sum_{\Gamma_\l}w(\Gamma_\l)=\sum_{\Gamma}w(\Gamma)
\]
where the sum is over all walks $\Gamma$ in $D$ of length $k+1$ from $v_i$ to $v_j$.\hfill$\blacksquare$\\

We provide an example to illustrate the use of this theorem. For the adjacency matrix $P$ $(1)$ for Figure \ref{fig:2.6}$(b)$, we compute $P^3$. We get\[
P^3=\kbordermatrix{ & A & B & C\\
A & x^3 & 3x^3 & 4x^3\\
B & 2x^3 & 2x^3 & 4x^3\\
C & x^3 & 3x^3 & 4x^3}.
\]\\
We are interested in the entry of row for $A$ and the column $C$, since $A$ is the initial state and $C$ is the only accept state. We have $4x^3$ in the $(A,C)$-entry, so we look at the details of where this term comes from. $(P^3)_{AC}$ is the following sum of the nonzero terms:\begin{equation}
(P^3)_{AC}=P_{AB}P_{BC}P_{CC}+P_{AB}P_{BA}P_{AC}+P_{AC}P_{CB}P_{BC}+P_{AC}P_{CC}P_{CC}.
\end{equation}
Every time we pass through a transition, we multiply by $x$. For a walk with three transitions from $A$ to $C$, we get $x\cdot x\cdot x=x^3$. With four distinct walks from $A$ to $C$ with exactly three transitions, we obtain $x^3$ four times, and therefore the $(A,C)$-entry of $P^3$ is $4x^3$.\\

Hence, by looking at the $(A,C)$-entry of $P^n$, we can find the number of distinct walks from $A$ to $C$ with $n$ transitions. Specifically, this number is the coefficient given our weight function. This is equivalent to the number of distinct $n$-letter words that are accepted by the automaton of the diagram in Figure \ref{fig:2.6}$(b)$. In the case where we have multiple accept states, we consider all the entries in the row corresponding to the initial state and a column corresponding to an accept state and sum them up.\\

Our next goal is to find the generating function $\sum s_n x^n$, where $s_n$ is the number of distinct $n$ length walks from one vertex to another by applying this technique. Suppose we have a digraph $D=(V,E)$ and the weight function $w$ defined by $w(e)=x$ for every edge $e\in E$. Because the coefficient of $x^n$ in the entry $(P^n)_{ij}$ describes the number of distinct walks of length $n$ from $v_i$ to $v_j$, the generating function for this is simply given by\[
f=\sum_{n=0}^\infty (P^n)_{ij},\quad\textrm{or equivalently,}\quad\left(\sum_{n=0}^\infty P^n\right)_{ij}.
\]

Notice that\[\begin{array}{rcl}
I & = & I+(P-P)+(P^2-P^2)+\cdots=(I+P+P^2+\cdots)-(P+P^2+P^3+\cdots)\\
  & = & (I+P+P^2+\cdots)-P(I+P+P^2+\cdots)=(I+P+P^2+\cdots)(I-P)
\end{array}
\]
This gives us $(I-P)^{-1}=\sum_{n=0}^\infty P^n$, so the $(i,j)$-entry of the matrix $(I-P)^{-1}$ is the desired function.\\

With the previous example, we obtain\[
(I-P)^{-1}=\kbordermatrix{ & A & B & C\\
 A & \ds\frac{-x^2-x+1}{-2 x^2-x+1} & \ds\frac{x}{-2 x^2-x+1} & \ds\frac{x^2+x}{-2 x^2-x+1}\vhhh\\
 B & \ds\frac{x-x^2}{-2 x^2-x+1} & \ds\frac{1-x}{-2 x^2-x+1} & \ds\frac{x^2+x}{-2 x^2-x+1}\vhhh\\
 C & \ds\frac{x^2}{-2 x^2-x+1} & \ds\frac{x}{-2 x^2-x+1} & \ds\frac{1-x^2}{-2 x^2-x+1}
}
\]
The $(A,C)$-entry is\[
\frac{x^2+x}{-2 x^2-x+1}=x+2x^2+4x^3+8x^4+\cdots,
\]
which is the generating function for the number of distinct $n$-letter words that are accepted by the automaton described in Figure \ref{fig:2.6}$(b)$.\\

Finally, we provide another example of the technique described above, which allows us to obtain more information from $P$ than just the number of distinct $n$-letter words accepted by the associated automaton. In the previous example, let us replace the weight function with the following one.
\[
w:T\to\Z[a,b,c]\quad\textrm{defined by}\quad w(t)=\left\{\begin{array}{cl} a &\textrm{ if }t\textrm{ is from }B\textrm{ to }A\\
 b &\textrm{ if }t\textrm{ is from }A\textrm{ to }B \textrm{ or from }C\textrm{ to }B\\
 c &\textrm{ if }t\textrm{ is from }A\textrm{ to }C\textrm{, from }B\textrm{ to }C \textrm{ or from }C\textrm{ to }C
 \end{array}\right.
\]
where $T$ is the set of all transitions and $\Z[a,b,c]$ is the commutative ring generated by elements $a$, $b$ and $c$ with integer coefficients. With this weight function, the adjacency matrix is\[
P=\kbordermatrix{ & A & B & C\\
A & 0 & b & c\\
B & a & 0 & c\\
C & 0 & b & c}.
\]\\
Again, we look at the $(A,C)$-entry of $P^3$ as an example. Since the computation of matrix multiplication does not change, we still get the equation $(2)$. Hence, it is\[
(P^3)_{AC}=b\cdot c\cdot c+b\cdot a\cdot c+c\cdot b\cdot c+c\cdot c\cdot c=abc+2bc^2+c^3.
\]

Here, we need to keep in mind that $a$, $b$ and $c$ are elements of the commutative ring, not letters for a regular language associated with the automaton. However, the defined weight function ``counts" each transition as a variable according to its label, so the word $bacbc\in\Sigma^*$ for example is weighted as $ab^2c^2\in\Z[a,b,c]$ by the weight function.\\

Thus, the equation above of three variables tells us that not only are there four $3$-letter words accepted by the automaton of Figure \ref{fig:2.6}$(b)$, but also that one of them consists of one $a$, one $b$ and one $c$, two of them consist of two $b$'s and one $c$, and the other one contains three $c$'s. Just as before, we can find the multivariate generating function describing the number of $a^{n_1}b^{n_2}c^{n_3}$ words as its coefficients by looking at the $(A,C)$-entry of $(I-P)^{-1}$. This time, we obtain\[
(I-P)^{-1}=\kbordermatrix{ & A & B & C\\
A & \ds\frac{-b c-c+1}{-a b-c b-c+1} & \ds\frac{b}{-a b-c b-c+1} & \ds\frac{b c+c}{-a b-c b-c+1} \vhhh\\
B & \ds\frac{a-a c}{-a b-c b-c+1} & \ds\frac{1-c}{-a b-c b-c+1} & \ds\frac{a c+c}{-a b-c b-c+1} \vhhh\\
C & \ds\frac{a b}{-a b-c b-c+1} & \ds\frac{b}{-a b-c b-c+1} & \ds\frac{1-a b}{-a b-c b-c+1}
}\]
The $(A,C)$-entry is\[
\frac{b c+c}{-a b-c b-c+1}=c+(bc+c^2)+(abc+2bc^2+c^3)+(a b^2 c + 2 a b c^2 + b^2 c^2 + 3 b c^3 + c^4)+\cdots
\]

In Chapter \ref{chap:4} and \ref{chap:6}, we apply this method to find generating functions for simple permutations of length greater than or equal to 4 in $\A$ and $\A'$ respectively.
\newpage
\chapter{Examples of finding generating functions}\label{chap:3}
In this chapter, we prove enumerative results for two classes, $\Av(123,213,132)$ and $\Av(4123,4213,\allowbreak4132)$, as lemmas. Both are proved by showing that a generating function for the class satisfies a specific functional equation.
\\
\section{Enumeration of the class $\Av(123,213,132)$}
We first state the result.\vhhh
\begin{lemma}\label{lem:3.1}
\textit{(Simion and Schmidt, 1985 \cite{SS1985}) The numbers of permutations of length $n$ in $\Av(123,\allowbreak213,132)$ form the Fibonacci sequence. Thus, the generating function for $\Av(123,213,132)$ is}
\[
f_{\Av(123,213,132)}=\frac{1}{1-x-x^2}.
\]
\end{lemma}\vhh
\textit{Proof.} For convenience, let $\C=\Av(123,213,132)$ and $F=f_{\Av(123,213,132)}$. Since clearly, $s_0=s_1=1$, we show this result by proving the number of length $n$ permutations in $\C$ is equal to the sum of the number of length $n-1$ permutations in $\C$ and the number of length $n-2$ permutations in $\C$ for all $n\geq 2$. Let $\pi\in\S_n\cap\C$ $(n\geq 2)$ be arbitrary. Suppose the biggest value $n$ appears after the second position of $\pi$, meaning $\pi(i)=n$ for some $i\geq 3$. This is an immediate contradiction, because either $\pi(1)<\pi(2)$ or $\pi(1)>\pi(2)$, and hence the flattening of positions 1, 2 and $i$ is either 123 or 213.\\

Now suppose $\pi(2)=n$. If $n-1$ shows up after the second position, then the values $\pi(1)$, $n$ and $n-1$ together create a subsequence whose flattening is 132, so we must have $\pi(1)=n-1$. Thus for the case $\pi(2)=n$, we have
\[
\pi=(n-1)n\pi(3)\pi(4)\cdots\pi(n).
\]

Notice that $\pi=12\ominus\s=21[12,\s]$ where $\s\in\S_{n-2}\cap\C$. Since 12 is skew-indecomposable, Proposition \ref{prop:2.3} guarantees that for each distinct $\s$, we obtain a unique $\pi$. Thus, we have an obvious bijection between $\S_{n-2}\cap\C$ and the set $\{\pi\in\C_n\cap\C:\pi(2)=n\}$, namely $\phi(\s)=12\ominus\s$, so the number of length $n$ permutations in $\C$ such that $\pi(2)=n$ is $s_{n-2}$.\\

Similarly, if $\pi(1)=n$, then we have $\pi=21[1,\s]=21[1,\s]$ where $\s\in\S_{n-1}\cap\C$. With the same argument above, we have $s_{n-1}$ distinct $\pi$ of length $n$ in $\C$ with $\pi(1)=n$.\\

Consequently, combining these observations together, we have the relation
\[
s_n=s_{n-1}+s_{n-2}\textrm{ for }n\geq 2,\qquad s_0=s_1=1,
\]
showing that $s_n$ forms the Fibonacci sequence.\\

We now translate this into a functional equation. We have\\
\[\begin{array}{rcl}
F & = & \ds\sum_{n=0}^\infty s_nx^n=1+x+\sum_{n=2}^\infty s_nx^n=1+x+\sum_{n=2}^\infty(s_{n-1}+s_{n-2})x^n\vhhh\\
  & = & \ds 1+x+\sum_{n=2}^\infty s_{n-1}x^n+\sum_{n=2}^\infty s_{n-2}x^n=1+x\left(\sum_{n=0}^\infty s_nx^n\right)+x^2\left(\sum_{n=0}^\infty s_nx^n\right)\vhhh\\
  & = & 1+xF+xF^2.
\end{array}
\]
\\
Thus, $F=1+xF+x^2F$. Solving this for $F$, we obtain the desired result.$\hfill\blacksquare$\\

The proof we presented suggests that any permutation in $\Av(123,213,132)$ can be written as\\
\[
\bigominus_{\s\in\{\e,1,12\}}\s.
\]\\
As it was previously mentioned, Lemma \ref{lem:3.1} was first proved in \cite{SS1985}. Note that with $\textrm{op}(\b)=\b^r$ where $\b$ is a permutation, $\textrm{op}(\{123,213,132\})=\{321,312,231\}$, so by Proposition \ref{prop:2.2}, we have the same enumeration result for $\Av(321,312,231)$. Permutations in this class are called \textit{free permutations} in \cite{PT2012}, and authors describe another way to enumerate this class.\\
\section{Enumeration of the class $\Av(4123,4213,4132)$}
\subsection{Number of permutations in $\Av(4123,4213,4132)$}
Next, we derive the generating function for the class $\Av(4123,4213,4132)$. Note that an alternate derivation can be found in \cite{AHPSV2015}, particularly, in Section 4.2.\\
\begin{lemma}\label{lem:3.2}
\textit{Let $G=f_{\Av(4123,4213,4132)}$. Then $G$ satisfies the equation}\\
\begin{equation}\label{eqn:3.1}
G=1+\frac{xG}{1-xG^2}.
\end{equation}
\end{lemma}\vhh
\textit{Proof.} For convenience, let $\C=\Av(4123,4213,4132)$. Given that $s_0=1$, Equation \ref{eqn:3.1} can be written as\\
\[\begin{array}{rcl}
G & = & 1+xG-xG^2+xG^3=1+xG+G\cdot G\cdot(G-1)x\\
 & = & 1+(s_0x+s_1x^2+\cdots)+(s_0+s_1x+\cdots)(s_0+s_1x+\cdots)(s_1x+s_2x^2+\cdots)x
\end{array}
\]
\\
Thus, equation \ref{eqn:3.1} claims that, for each $n\geq 1$,\\
\begin{eqnarray}
s_n & = s_{n-1} &\ds+\sum_{0\leq i,j\leq n-2,\,\,1\leq k\leq n-1\,\atop i+j+k=n-1}s_i s_j s_k\label{eqn:3.2}\vhhh\\
 & = s_{n-1} & +(s_0s_{n-2}+s_1s_{n-3}+\cdots+s_{n-2}s_0)s_1\nn\\
 & & +(s_0s_{n-3}+s_1s_{n-4}+\cdots+s_{n-3}s_0)s_2\nn\\
 & & \hspace{0.7in}\vdots\nn\\
 & & +(s_0s_1+s_1s_0)s_{n-2}\nn\\
 & & +s_0 s_0 s_{n-1}\nn
\end{eqnarray}
\\
We prove the lemma by showing Equation \ref{eqn:3.2} is true for all $n\geq 1$.\\

Let $n\geq 1$ be arbitrary. We claim that for each $k$ with $1\leq k\leq n-1$, $(s_0s_{n-(k+1)}+s_1s_{n-(k+2)}+\cdots+s_{n-(k+1)}s_0)s_k$ is the number of permutations of length $n$ in $\C$ whose last element is $k+1$. For example, if $n=5$ and $k=2$, $(s_0s_2+s_1s_1+s_2s_0)s_2$ is the total number of permutations $\pi$ such that $|\pi|=5$ and $\pi(5)=3$.\\

\begin{figure}
\[\hspace{-.35cm}\begin{tikzpicture}
\tikzstyle{vertex}=[circle,fill=black,inner sep=0pt,font=\tiny]
\coordinate (O) at (0,0);
\draw[gray!80,fill=gray!80](-3.2,0) rectangle (0,3.2);
\draw[gray!80,fill=gray!80](0,1.6) rectangle (1.6,3.2);
\draw[gray!80,fill=gray!80](1.6,-2.4) rectangle (3.2,1.6);
\draw[gray!80,fill=gray!80](0,-2.4) rectangle (1.6,0);
\path[draw,thick] (-3.2,-2.4) rectangle (3.2,3.2);
\draw [thick,dashed] (0,-2.4) -- (0,3.2);
\draw [thick,dashed] (1.6,-2.4) -- (1.6,3.2);
\draw [thick,dashed] (-3.2,0) -- (3.2,0);
\draw [thick,dashed] (0,1.6) -- (3.2,1.6);
\draw[color=black,fill=black] (3.2,0) circle [radius=.08] node [black,right] {$k+1$};
\draw[color=black,fill=black] (0,-0.8) circle [radius=.08] node [black,right] {$x$};
\draw[color=black,fill=black] (1.6,-1.6) circle [radius=.08] node [black,right] {$y$};
\draw [decorate,decoration={brace,amplitude=10pt},xshift=-4pt,yshift=0pt]
(-3.2,-2.4) -- (-3.2,0) node [black,midway,xshift=-.5cm,text width=1cm]
{$K$};
\node [black] at (0.8,0.8) {$I$};
\node [black] at (2.4,2.4) {$J$};
\end{tikzpicture}\]
\caption{The graphs of a permutation in $\C_{ik}$.}
\label{fig:3.1}
\end{figure}
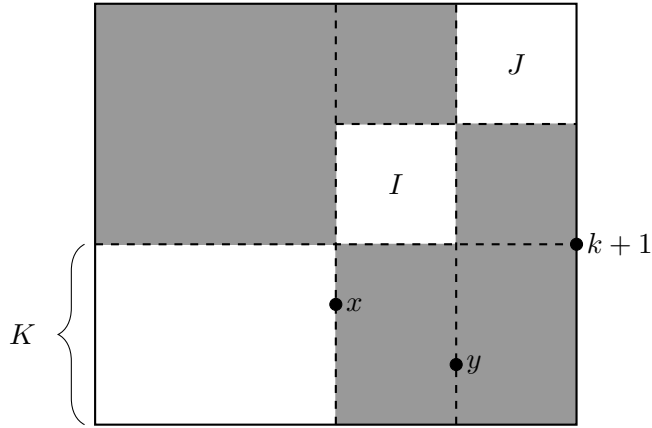

To prove this claim, let $k$ be arbitrary with $1\leq k\leq n-1$, and denote by $\C_k$ the set of permutations of length $n$ in $\C$ having $k+1$ in the last position. We also define particular subsets of $\C_k$. First, let $K$ be the range $[1,k]$. Let $x,y\in K$ be values that have the second right-most position and the last position respectively, so $x=\pi(a)$ and $y=\pi(b)$ where $a<b$ are the two largest integers with $\pi(a),\pi(b)\in K$. If $k=1$, we let $y=1$ and do not define $x$. For $k\neq 1$, let $I=(\pi^{-1}(x),\pi^{-1}(y))$, the segment between the positions of $x$ and $y$, but excluding $\pi^{-1}(x)$ and $\pi^{-1}(y)$ themselves. If $k=1$, let $I=[1,\pi^{-1}(y))$, the segment between the first position and the one immediately to the left of the position of $y=1$. Similarly, let $J=(\pi^{-1}(y),\pi^{-1}(k+1))$. Notice that for all $a\in I$ and $b\in J$, $\pi(a)<\pi(b)$ because, otherwise, $\pi$ would contain 4132 with $y$ and $k+1$ corresponding to $1$ and $2$ respectively. Furthermore, positions of every value greater than $k+1$ must belong to either $I$ or $J$, because if this is not the case, there exists a value $z\geq k+1$ such that $\pi^{-1}(z)<\pi^{-1}(x)$, and $zxy(k+1)$ forms either a 4123 or a 4213 pattern. These observations imply both segments $I$ and $J$ form blocks. The length of these blocks are some nonnegative integers $i$ and $j$ such that $i+j+k=n-1$. For all $i$ with $0\leq i\leq n-2$, we define $\C_{ik}$ to be the set of permutations in $\C_k$ such that $|I|=i$. Figure \ref{fig:3.1} shows what a permutation in $\C_{ik}$ looks like.\\

Let $i$ and $j$ be positive integers with $0\leq i\leq n-2$ and $j=(n-1)-k-i$. We now define a bijection between $\C_{ik}$ and $(\S_i\cap\C)\times(\S_j\cap\C)\times(\S_k\cap\C)$. Let $\phi:\C_{ik}\to (\S_i\cap\C)\times(\S_j\cap\C)\times(\S_k\cap\C)$ be a function defined by\[
\phi(\pi)=(\s_1,\s_2,\s_3)
\]
for all $\pi\in\C_{ik}$ where\[\begin{array}{c}
\s_1\textrm{ is the flattening of }\pi(I),\qquad\s_2\textrm{ is the flattening of }\pi(J),\\
\s_3\textrm{ is the flattening of }\pi(1)\pi(2)\cdots xy.
\end{array}
\]

Needless to say, $\s_1$, $\s_2$ and $\s_3$ are in $\C$ simply by the definition of a permutation class. Also, how we defined each segment and range clearly indicates $|\s_1|=i$, $|\s_2|=j$ and $|\s_3|=k$, so $(\s_1,\s_2,\s_3)\in(\S_i\cap\C)\times(\S_j\cap\C)\times(\S_k\cap\C)$.\\

We prove that $\phi$ is a bijection by constructing $\phi^{-1}$, the inverse of $\phi$, and show it is a function. The inverse map $\phi^{-1}:(\S_i\cap\C)\times(\S_j\cap\C)\times(\S_k\cap\C)\to\C_{ik}$ is defined by $\phi^{-1}[(\s_1,\s_2,\s_3)]=\pi$ where\\
\[
\pi(\l)=\left\{\begin{array}{cl}
\s_3(\l) & \textrm{ if }1\leq\l\leq k-1\\
\s_1(\l-(k-1))+(k+1) & \textrm{ if }k\leq\l\leq k-1+i\\
\s_3(\l-i) & \textrm{ if }\l=k+i\\
\s_2(\l-(k+i))+(k+1+i) & \textrm{ if }k+i+1\leq\l\leq k+i+j\\
k+1 & \textrm{ if }\l=n
\end{array}\right..
\]\\
In one-line notation, this is equivalent to
\[
\phi^{-1}[(\s_1,\s_2,\s_3)]=\s_3(1)\cdots\underbrace{\s_3(k-1)}_{x}\underbrace{\s_1'(1)\cdots\s_1'(i)}_{I}\underbrace{\s_3(k)}_{y}\underbrace{\s_2'(1)\cdots\s_2'(j)}_{J}(k+1),
\]
for all $(\s_1,\s_2,\s_3)\in(\S_i\cap\C)\times(\S_j\cap\C)\times(\S_k\cap\C)$ where $\s_1'(\l)=\s_1(\l)+(k+1)$ and $\s_2'(m)=\s_2(m)+(k+1+i)$ for all $\l$ and $m$ ($1\leq\l\leq i$, $1\leq m \leq j$).\\

We show $\phi^{-1}$ maps into $\C_{ik}$ by contradiction. Suppose there exists $(\s_1,\s_2,\s_3)\in(\S_i\cap\C)\times(\S_j\cap\C)\times(\S_k\cap\C)$ such that $\pi=\phi^{-1}[(\s_1,\s_2,\s_3)]$ is not in $\C_{ik}$. Since $\pi$ has $k+1$ in $n$-th position, and the way $\pi$ is constructed forces the segment $I$ to have length $i$, $\pi\notin\C_{ik}$ implies $\b\cont\pi$ for some $\b\in\{4123,4213,4132\}$. Because every permutation in the basis has the value 4 in the $1^{\textrm{st}}$ position, $\pi$ must have the value $v_4$ corresponding to the 4 in $\b$ before the positions of the other values $v_1$, $v_2$ and $v_3$ corresponding to 1, 2 and 3 respectively. Let $m_1$, $m_2$, $m_3$ and $m_4$ be the positions of $v_1$, $v_2$, $v_3$ and $v_4$ respectively. Suppose $1\leq m_4\leq k-1$. Since this implies that the value $v_4$ is at most $k$, any position $m$ with $m_4+1\leq m\leq n$ such that $\pi(m)<v_4$ must satisfy $m_4+1\leq m\leq k-1$ or $m=k+1$. Hence, all of $v_1$, $v_2$, $v_3$ and $v_4$ are determined by $\s_3$, but this means $\b\in\s_3$. We can similarly show that we cannot have $k\leq m_4\leq k-1+i$ or $k+i+1\leq m_4\leq k+i+j$. In addition, $m_4\neq k+1$ and $m_4\neq n$ because there is no position greater than $m_4$ whose value is less than $v_4$. Consequently, $\pi=\phi^{-1}[(\s_1,\s_2,\s_3)]\in\C_{ik}$ for every $(\s_1,\s_2,\s_3)\in(\S_i\cap\C)\times(\S_j\cap\C)\times(\S_k\cap\C)$.\\

The way $\phi^{-1}$ is constructed clearly shows that it is the inverse of $\phi$, so $\phi$ is a bijection. Hence, $|\C_{ik}|=|(\S_i\cap\C)\times(\S_j\cap\C)\times(\S_k\cap\C)|=s_is_js_k$, and\[
\C_k=\bigcup_{0\leq i\leq n-2}\C_{ik},
\]
which implies
\[
|\C_k|=\sum_{i=0}^{n-(k+1)}|\C_{ik}|=(s_0s_{n-(k+1)}+\cdots+s_{n-(k+1)}s_0)s_k,
\]\\
so this completes the proof of the claim.\\

Finally, we consider permutations in $\C$ of length $n$ having a 1 in the last position. Notice that every permutation $\pi\in\C$ where $|\pi|=n$ and $\pi(n)=1$ can be written as $\tau\ominus 1$ where $\tau$ is some permutation in $\C$ whose length is $n-1$ (including the case of $\pi=1$), and for every permutation $\tau\in\S_{n-1}\cap\C$, $\tau\ominus 1\in\C$. Thus, the number of such permutations is the same as the number of permutations of length $n-1$ in $\C$, which is simply $s_{n-1}$. Summing up every possible case, we obtain\\
\[s_n=s_{n-1}+\sum_{k=1}^{n-1}|\C_k|,\]\\
for all $n\geq 1$, and this is equation \ref{eqn:3.2}. This proves the generating function $G$ for $\C$ satisfies equation \ref{eqn:3.1} and completes the proof.\hfill$\blacksquare$\\

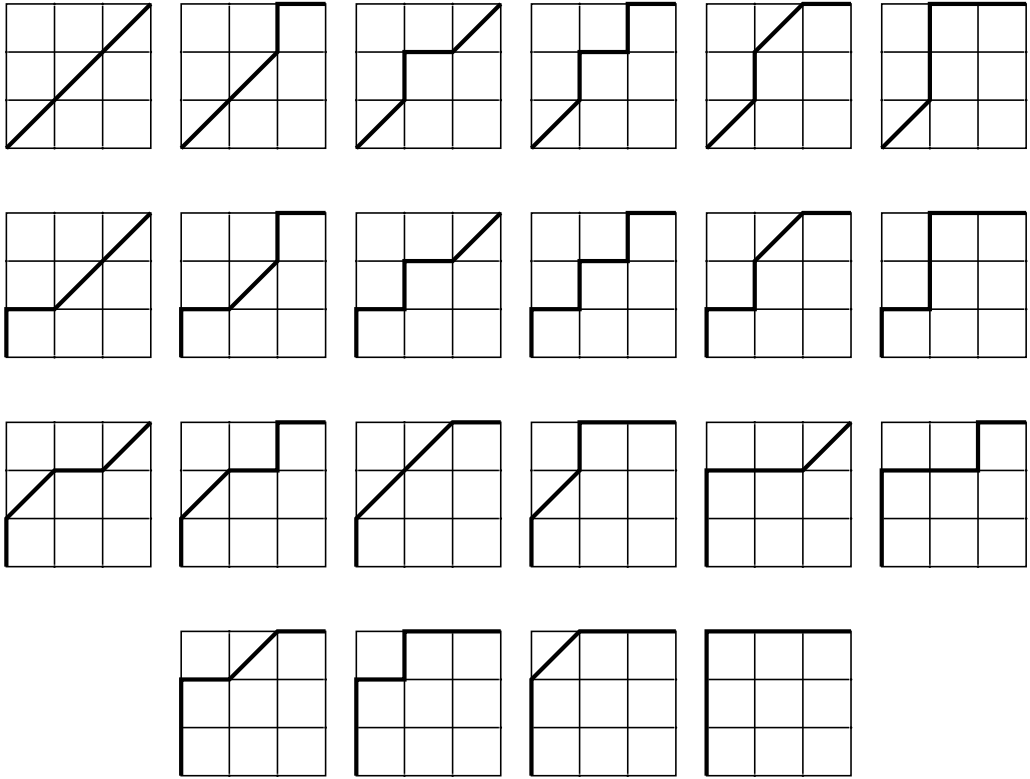
\begin{figure}[b!]
\[\begin{array}{cccccc}\begin{tikzpicture}[xscale=0.8,yscale=0.8]
\tikzstyle{vertex}=[circle,fill=black,inner sep=0pt,font=\tiny]
\coordinate (O) at (0,0);
\path[draw,semithick] (-1.2,-1.2) rectangle (1.2,1.2);
 \foreach \name/\x/\y in {a/-0.4/1.2, b/0.4/1.2, A/-0.4/-1.2, B/0.4/-1.2, 1/1.2/-0.4, 2/1.2/0.4, 3/-1.2/-0.4, 4/-1.2/0.4}
     \node[vertex] (\name) at (\x,\y) {};
 \foreach \from/\to in {a/A, b/B, 1/3, 2/4}
 \draw[black, semithick] (\from) -- (\to);
\draw[black, ultra thick]  (-1.2,-1.2) -- (1.2,1.2);
\end{tikzpicture}\quad & \begin{tikzpicture}[xscale=0.8,yscale=0.8]
\tikzstyle{vertex}=[circle,fill=black,inner sep=0pt,font=\tiny]
\coordinate (O) at (0,0);
\path[draw,semithick] (-1.2,-1.2) rectangle (1.2,1.2);
 \foreach \name/\x/\y in {a/-0.4/1.2, b/0.4/1.2, A/-0.4/-1.2, B/0.4/-1.2, 1/1.2/-0.4, 2/1.2/0.4, 3/-1.2/-0.4, 4/-1.2/0.4}
     \node[vertex] (\name) at (\x,\y) {};
 \foreach \from/\to in {a/A, b/B, 1/3, 2/4}
 \draw[black, semithick] (\from) -- (\to);
\draw[black, ultra thick]  (-1.2,-1.2) -- (0.4,0.4) -- (0.4,1.2) -- (1.2,1.2);
\end{tikzpicture}\quad & \begin{tikzpicture}[xscale=0.8,yscale=0.8]
\tikzstyle{vertex}=[circle,fill=black,inner sep=0pt,font=\tiny]
\coordinate (O) at (0,0);
\path[draw,semithick] (-1.2,-1.2) rectangle (1.2,1.2);
 \foreach \name/\x/\y in {a/-0.4/1.2, b/0.4/1.2, A/-0.4/-1.2, B/0.4/-1.2, 1/1.2/-0.4, 2/1.2/0.4, 3/-1.2/-0.4, 4/-1.2/0.4}
     \node[vertex] (\name) at (\x,\y) {};
 \foreach \from/\to in {a/A, b/B, 1/3, 2/4}
 \draw[black, semithick] (\from) -- (\to);
\draw[black, ultra thick]  (-1.2,-1.2) -- (-0.4,-0.4) -- (-0.4,0.4) -- (0.4,0.4) -- (1.2,1.2);
\end{tikzpicture}\quad & \begin{tikzpicture}[xscale=0.8,yscale=0.8]
\tikzstyle{vertex}=[circle,fill=black,inner sep=0pt,font=\tiny]
\coordinate (O) at (0,0);
\path[draw,semithick] (-1.2,-1.2) rectangle (1.2,1.2);
 \foreach \name/\x/\y in {a/-0.4/1.2, b/0.4/1.2, A/-0.4/-1.2, B/0.4/-1.2, 1/1.2/-0.4, 2/1.2/0.4, 3/-1.2/-0.4, 4/-1.2/0.4}
     \node[vertex] (\name) at (\x,\y) {};
 \foreach \from/\to in {a/A, b/B, 1/3, 2/4}
 \draw[black, semithick] (\from) -- (\to);
\draw[black, ultra thick]  (-1.2,-1.2) -- (-0.4,-0.4) -- (-0.4,0.4) -- (0.4,0.4) -- (0.4,1.2) -- (1.2,1.2);
\end{tikzpicture}\quad & \begin{tikzpicture}[xscale=0.8,yscale=0.8]
\tikzstyle{vertex}=[circle,fill=black,inner sep=0pt,font=\tiny]
\coordinate (O) at (0,0);
\path[draw,semithick] (-1.2,-1.2) rectangle (1.2,1.2);
 \foreach \name/\x/\y in {a/-0.4/1.2, b/0.4/1.2, A/-0.4/-1.2, B/0.4/-1.2, 1/1.2/-0.4, 2/1.2/0.4, 3/-1.2/-0.4, 4/-1.2/0.4}
     \node[vertex] (\name) at (\x,\y) {};
 \foreach \from/\to in {a/A, b/B, 1/3, 2/4}
 \draw[black, semithick] (\from) -- (\to);
\draw[black, ultra thick]  (-1.2,-1.2) -- (-0.4,-0.4) -- (-0.4,0.4) -- (0.4,1.2) -- (1.2,1.2);
\end{tikzpicture}\quad & \begin{tikzpicture}[xscale=0.8,yscale=0.8]
\tikzstyle{vertex}=[circle,fill=black,inner sep=0pt,font=\tiny]
\coordinate (O) at (0,0);
\path[draw,semithick] (-1.2,-1.2) rectangle (1.2,1.2);
 \foreach \name/\x/\y in {a/-0.4/1.2, b/0.4/1.2, A/-0.4/-1.2, B/0.4/-1.2, 1/1.2/-0.4, 2/1.2/0.4, 3/-1.2/-0.4, 4/-1.2/0.4}
     \node[vertex] (\name) at (\x,\y) {};
 \foreach \from/\to in {a/A, b/B, 1/3, 2/4}
 \draw[black, semithick] (\from) -- (\to);
\draw[black, ultra thick]  (-1.2,-1.2) -- (-0.4,-0.4) -- (-0.4,1.2) -- (1.2,1.2);
\end{tikzpicture}\\ 
\\
\begin{tikzpicture}[xscale=0.8,yscale=0.8]
\tikzstyle{vertex}=[circle,fill=black,inner sep=0pt,font=\tiny]
\coordinate (O) at (0,0);
\path[draw,semithick] (-1.2,-1.2) rectangle (1.2,1.2);
 \foreach \name/\x/\y in {a/-0.4/1.2, b/0.4/1.2, A/-0.4/-1.2, B/0.4/-1.2, 1/1.2/-0.4, 2/1.2/0.4, 3/-1.2/-0.4, 4/-1.2/0.4}
     \node[vertex] (\name) at (\x,\y) {};
 \foreach \from/\to in {a/A, b/B, 1/3, 2/4}
 \draw[black, semithick] (\from) -- (\to);
\draw[black, ultra thick]  (-1.2,-1.2) -- (-1.2,-0.4) -- (-0.4,-0.4) -- (1.2,1.2);
\end{tikzpicture}\quad & \begin{tikzpicture}[xscale=0.8,yscale=0.8]
\tikzstyle{vertex}=[circle,fill=black,inner sep=0pt,font=\tiny]
\coordinate (O) at (0,0);
\path[draw,semithick] (-1.2,-1.2) rectangle (1.2,1.2);
 \foreach \name/\x/\y in {a/-0.4/1.2, b/0.4/1.2, A/-0.4/-1.2, B/0.4/-1.2, 1/1.2/-0.4, 2/1.2/0.4, 3/-1.2/-0.4, 4/-1.2/0.4}
     \node[vertex] (\name) at (\x,\y) {};
 \foreach \from/\to in {a/A, b/B, 1/3, 2/4}
 \draw[black, semithick] (\from) -- (\to);
\draw[black, ultra thick]  (-1.2,-1.2) -- (-1.2,-0.4) -- (-0.4,-0.4) -- (0.4,0.4) -- (0.4,1.2) -- (1.2,1.2);
\end{tikzpicture}\quad & \begin{tikzpicture}[xscale=0.8,yscale=0.8]
\tikzstyle{vertex}=[circle,fill=black,inner sep=0pt,font=\tiny]
\coordinate (O) at (0,0);
\path[draw,semithick] (-1.2,-1.2) rectangle (1.2,1.2);
 \foreach \name/\x/\y in {a/-0.4/1.2, b/0.4/1.2, A/-0.4/-1.2, B/0.4/-1.2, 1/1.2/-0.4, 2/1.2/0.4, 3/-1.2/-0.4, 4/-1.2/0.4}
     \node[vertex] (\name) at (\x,\y) {};
 \foreach \from/\to in {a/A, b/B, 1/3, 2/4}
 \draw[black, semithick] (\from) -- (\to);
\draw[black, ultra thick]  (-1.2,-1.2) -- (-1.2,-0.4) -- (-0.4,-0.4) -- (-0.4,-0.4) -- (-0.4,0.4) -- (0.4,0.4) -- (1.2,1.2);
\end{tikzpicture}\quad & \begin{tikzpicture}[xscale=0.8,yscale=0.8]
\tikzstyle{vertex}=[circle,fill=black,inner sep=0pt,font=\tiny]
\coordinate (O) at (0,0);
\path[draw,semithick] (-1.2,-1.2) rectangle (1.2,1.2);
 \foreach \name/\x/\y in {a/-0.4/1.2, b/0.4/1.2, A/-0.4/-1.2, B/0.4/-1.2, 1/1.2/-0.4, 2/1.2/0.4, 3/-1.2/-0.4, 4/-1.2/0.4}
     \node[vertex] (\name) at (\x,\y) {};
 \foreach \from/\to in {a/A, b/B, 1/3, 2/4}
 \draw[black, semithick] (\from) -- (\to);
\draw[black, ultra thick]  (-1.2,-1.2) -- (-1.2,-0.4) -- (-0.4,-0.4) -- (-0.4,-0.4) -- (-0.4,0.4) -- (0.4,0.4) -- (0.4,1.2) -- (1.2,1.2);
\end{tikzpicture}\quad & \begin{tikzpicture}[xscale=0.8,yscale=0.8]
\tikzstyle{vertex}=[circle,fill=black,inner sep=0pt,font=\tiny]
\coordinate (O) at (0,0);
\path[draw,semithick] (-1.2,-1.2) rectangle (1.2,1.2);
 \foreach \name/\x/\y in {a/-0.4/1.2, b/0.4/1.2, A/-0.4/-1.2, B/0.4/-1.2, 1/1.2/-0.4, 2/1.2/0.4, 3/-1.2/-0.4, 4/-1.2/0.4}
     \node[vertex] (\name) at (\x,\y) {};
 \foreach \from/\to in {a/A, b/B, 1/3, 2/4}
 \draw[black, semithick] (\from) -- (\to);
\draw[black, ultra thick]  (-1.2,-1.2) -- (-1.2,-0.4) -- (-0.4,-0.4) -- (-0.4,-0.4) -- (-0.4,0.4) -- (0.4,1.2) -- (1.2,1.2);
\end{tikzpicture}\quad & \begin{tikzpicture}[xscale=0.8,yscale=0.8]
\tikzstyle{vertex}=[circle,fill=black,inner sep=0pt,font=\tiny]
\coordinate (O) at (0,0);
\path[draw,semithick] (-1.2,-1.2) rectangle (1.2,1.2);
 \foreach \name/\x/\y in {a/-0.4/1.2, b/0.4/1.2, A/-0.4/-1.2, B/0.4/-1.2, 1/1.2/-0.4, 2/1.2/0.4, 3/-1.2/-0.4, 4/-1.2/0.4}
     \node[vertex] (\name) at (\x,\y) {};
 \foreach \from/\to in {a/A, b/B, 1/3, 2/4}
 \draw[black, semithick] (\from) -- (\to);
\draw[black, ultra thick]  (-1.2,-1.2) -- (-1.2,-0.4) -- (-0.4,-0.4) -- (-0.4,-0.4) -- (-0.4,1.2) -- (1.2,1.2);
\end{tikzpicture}\\ 
\\
\begin{tikzpicture}[xscale=0.8,yscale=0.8]
\tikzstyle{vertex}=[circle,fill=black,inner sep=0pt,font=\tiny]
\coordinate (O) at (0,0);
\path[draw,semithick] (-1.2,-1.2) rectangle (1.2,1.2);
 \foreach \name/\x/\y in {a/-0.4/1.2, b/0.4/1.2, A/-0.4/-1.2, B/0.4/-1.2, 1/1.2/-0.4, 2/1.2/0.4, 3/-1.2/-0.4, 4/-1.2/0.4}
     \node[vertex] (\name) at (\x,\y) {};
 \foreach \from/\to in {a/A, b/B, 1/3, 2/4}
 \draw[black, semithick] (\from) -- (\to);
\draw[black, ultra thick]  (-1.2,-1.2) -- (-1.2,-0.4) -- (-0.4,0.4) -- (0.4,0.4) -- (1.2,1.2);
\end{tikzpicture}\quad & \begin{tikzpicture}[xscale=0.8,yscale=0.8]
\tikzstyle{vertex}=[circle,fill=black,inner sep=0pt,font=\tiny]
\coordinate (O) at (0,0);
\path[draw,semithick] (-1.2,-1.2) rectangle (1.2,1.2);
 \foreach \name/\x/\y in {a/-0.4/1.2, b/0.4/1.2, A/-0.4/-1.2, B/0.4/-1.2, 1/1.2/-0.4, 2/1.2/0.4, 3/-1.2/-0.4, 4/-1.2/0.4}
     \node[vertex] (\name) at (\x,\y) {};
 \foreach \from/\to in {a/A, b/B, 1/3, 2/4}
 \draw[black, semithick] (\from) -- (\to);
\draw[black, ultra thick]  (-1.2,-1.2) -- (-1.2,-0.4) -- (-0.4,0.4) -- (0.4,0.4) -- (0.4,1.2) -- (1.2,1.2);
\end{tikzpicture}\quad & \begin{tikzpicture}[xscale=0.8,yscale=0.8]
\tikzstyle{vertex}=[circle,fill=black,inner sep=0pt,font=\tiny]
\coordinate (O) at (0,0);
\path[draw,semithick] (-1.2,-1.2) rectangle (1.2,1.2);
 \foreach \name/\x/\y in {a/-0.4/1.2, b/0.4/1.2, A/-0.4/-1.2, B/0.4/-1.2, 1/1.2/-0.4, 2/1.2/0.4, 3/-1.2/-0.4, 4/-1.2/0.4}
     \node[vertex] (\name) at (\x,\y) {};
 \foreach \from/\to in {a/A, b/B, 1/3, 2/4}
 \draw[black, semithick] (\from) -- (\to);
\draw[black, ultra thick]  (-1.2,-1.2) -- (-1.2,-0.4) -- (0.4,1.2) -- (1.2,1.2);
\end{tikzpicture}\quad & \begin{tikzpicture}[xscale=0.8,yscale=0.8]
\tikzstyle{vertex}=[circle,fill=black,inner sep=0pt,font=\tiny]
\coordinate (O) at (0,0);
\path[draw,semithick] (-1.2,-1.2) rectangle (1.2,1.2);
 \foreach \name/\x/\y in {a/-0.4/1.2, b/0.4/1.2, A/-0.4/-1.2, B/0.4/-1.2, 1/1.2/-0.4, 2/1.2/0.4, 3/-1.2/-0.4, 4/-1.2/0.4}
     \node[vertex] (\name) at (\x,\y) {};
 \foreach \from/\to in {a/A, b/B, 1/3, 2/4}
 \draw[black, semithick] (\from) -- (\to);
\draw[black, ultra thick]  (-1.2,-1.2) -- (-1.2,-0.4) -- (-0.4,0.4) -- (-0.4,1.2) -- (1.2,1.2);
\end{tikzpicture}\quad & \begin{tikzpicture}[xscale=0.8,yscale=0.8]
\tikzstyle{vertex}=[circle,fill=black,inner sep=0pt,font=\tiny]
\coordinate (O) at (0,0);
\path[draw,semithick] (-1.2,-1.2) rectangle (1.2,1.2);
 \foreach \name/\x/\y in {a/-0.4/1.2, b/0.4/1.2, A/-0.4/-1.2, B/0.4/-1.2, 1/1.2/-0.4, 2/1.2/0.4, 3/-1.2/-0.4, 4/-1.2/0.4}
     \node[vertex] (\name) at (\x,\y) {};
 \foreach \from/\to in {a/A, b/B, 1/3, 2/4}
 \draw[black, semithick] (\from) -- (\to);
\draw[black, ultra thick]  (-1.2,-1.2) -- (-1.2,0.4) -- (0.4,0.4) -- (1.2,1.2);
\end{tikzpicture}\quad & \begin{tikzpicture}[xscale=0.8,yscale=0.8]
\tikzstyle{vertex}=[circle,fill=black,inner sep=0pt,font=\tiny]
\coordinate (O) at (0,0);
\path[draw,semithick] (-1.2,-1.2) rectangle (1.2,1.2);
 \foreach \name/\x/\y in {a/-0.4/1.2, b/0.4/1.2, A/-0.4/-1.2, B/0.4/-1.2, 1/1.2/-0.4, 2/1.2/0.4, 3/-1.2/-0.4, 4/-1.2/0.4}
     \node[vertex] (\name) at (\x,\y) {};
 \foreach \from/\to in {a/A, b/B, 1/3, 2/4}
 \draw[black, semithick] (\from) -- (\to);
\draw[black, ultra thick]  (-1.2,-1.2) -- (-1.2,0.4) -- (0.4,0.4) -- (0.4,1.2) -- (1.2,1.2);
\end{tikzpicture}\\ 
\\
 & \begin{tikzpicture}[xscale=0.8,yscale=0.8]
\tikzstyle{vertex}=[circle,fill=black,inner sep=0pt,font=\tiny]
\coordinate (O) at (0,0);
\path[draw,semithick] (-1.2,-1.2) rectangle (1.2,1.2);
 \foreach \name/\x/\y in {a/-0.4/1.2, b/0.4/1.2, A/-0.4/-1.2, B/0.4/-1.2, 1/1.2/-0.4, 2/1.2/0.4, 3/-1.2/-0.4, 4/-1.2/0.4}
     \node[vertex] (\name) at (\x,\y) {};
 \foreach \from/\to in {a/A, b/B, 1/3, 2/4}
 \draw[black, semithick] (\from) -- (\to);
\draw[black, ultra thick]  (-1.2,-1.2) -- (-1.2,0.4) -- (-0.4,0.4) -- (0.4,1.2) -- (1.2,1.2);
\end{tikzpicture}\quad & \begin{tikzpicture}[xscale=0.8,yscale=0.8]
\tikzstyle{vertex}=[circle,fill=black,inner sep=0pt,font=\tiny]
\coordinate (O) at (0,0);
\path[draw,semithick] (-1.2,-1.2) rectangle (1.2,1.2);
 \foreach \name/\x/\y in {a/-0.4/1.2, b/0.4/1.2, A/-0.4/-1.2, B/0.4/-1.2, 1/1.2/-0.4, 2/1.2/0.4, 3/-1.2/-0.4, 4/-1.2/0.4}
     \node[vertex] (\name) at (\x,\y) {};
 \foreach \from/\to in {a/A, b/B, 1/3, 2/4}
 \draw[black, semithick] (\from) -- (\to);
\draw[black, ultra thick]  (-1.2,-1.2) -- (-1.2,0.4) -- (-0.4,0.4) -- (-0.4,1.2) -- (1.2,1.2);
\end{tikzpicture}\quad & \begin{tikzpicture}[xscale=0.8,yscale=0.8]
\tikzstyle{vertex}=[circle,fill=black,inner sep=0pt,font=\tiny]
\coordinate (O) at (0,0);
\path[draw,semithick] (-1.2,-1.2) rectangle (1.2,1.2);
 \foreach \name/\x/\y in {a/-0.4/1.2, b/0.4/1.2, A/-0.4/-1.2, B/0.4/-1.2, 1/1.2/-0.4, 2/1.2/0.4, 3/-1.2/-0.4, 4/-1.2/0.4}
     \node[vertex] (\name) at (\x,\y) {};
 \foreach \from/\to in {a/A, b/B, 1/3, 2/4}
 \draw[black, semithick] (\from) -- (\to);
\draw[black, ultra thick]  (-1.2,-1.2) -- (-1.2,0.4) -- (-0.4,1.2)-- (1.2,1.2);
\end{tikzpicture}\quad & \begin{tikzpicture}[xscale=0.8,yscale=0.8]
\tikzstyle{vertex}=[circle,fill=black,inner sep=0pt,font=\tiny]
\coordinate (O) at (0,0);
\path[draw,semithick] (-1.2,-1.2) rectangle (1.2,1.2);
 \foreach \name/\x/\y in {a/-0.4/1.2, b/0.4/1.2, A/-0.4/-1.2, B/0.4/-1.2, 1/1.2/-0.4, 2/1.2/0.4, 3/-1.2/-0.4, 4/-1.2/0.4}
     \node[vertex] (\name) at (\x,\y) {};
 \foreach \from/\to in {a/A, b/B, 1/3, 2/4}
 \draw[black, semithick] (\from) -- (\to);
\draw[black, ultra thick]  (-1.2,-1.2) -- (-1.2,1.2) -- (1.2,1.2);
\end{tikzpicture}
\end{array}\]
\caption{The list of Schr\"{o}der 3-paths.}
\label{fig:3.2}
\end{figure}

There is another famous combinatorial object whose enumeration is related to the class $\Av(4123,4213,4132)$. \textit{Schr\"oder $n$-paths} are lattice paths in the Cartesian plane using $\langle 0,1\rangle$, $\langle 1,0\rangle$ and $\langle 1,1\rangle$ steps that start at $(0,0)$, end at $(n,n)$ and stay on or above the $x=y$ line. The number of distinct Schr\"oder $n$-paths is called the \textit{$n$-th Schr\"oder number.} Figure \ref{fig:3.2} is the list of all Schr\"oder $3$-paths. As we can see, there are 22 distinct Schr\"oder $3$-paths.\\

It turns out the generating function for Schr\"{o}der $(n-1)$-paths with no three consecutive up-steps is known to satisfy the equation $(1)$ as well (as listed as A106228 in \cite{Sloane}). Hence, the number of length $n$ permutations in $\Av(4123,4213,4132)$ is\\
\[
s_n=\left\{\begin{array}{cl}
1 & \textit{if }n=0\\
\textit{the number of Schr\"{o}der }(n-1)\textit{-paths}\atop\textit{with no three consecutive up-steps} & \textit{if }n\geq 1
\end{array}\right..
\]

The last path shown in Figure \ref{fig:3.2} is the only one with three consecutive up-steps, so there are 21 desired Schr\"{o}der 3-paths. In fact, there are 21 permutations of length 4 in the class $\Av(4123,4213,4132)$, namely every permutation of length 4 except the ones in the basis.\\

Proving Lemma \ref{lem:3.2} by finding a bijection between $\S_n\cap\Av(4123,4213,4132)$ and the set of Schr\"{o}der $(n-1)$-paths having no triple up-steps for each $n\geq 1$ would be ideal. However, this problem remains unsolved. It is worth noting that if we take any two permutations $\b_1$, $\b_2$ from $\{4123,4213,4132\}$ and form a subbasis $\{\b_1,\b_2\}$, then $s_n(\b_1,\b_2)$ is the $(n-1)$-th Schr\"{o}der number. Bijective proofs are given in \cite{Kr2000,K2003}, although the bijections given are to Schr\"{o}der generating trees rather than Schr\"{o}der paths themselves. Thus, one could possibly prove this result bijectively by showing a permutation of length $n$ containing the other permutation $\b_3$ corresponds to a tree that represents a Schr\"{o}der path having triple up-steps.\\
\subsection{Skew-indecomposable permutations in $\Av(4123,4213,4132)$}
For the remainder of this chapter, we determine the generating function describing the number of skew-indecomposable permutations of length $n$ in $\Av(4123,4213,4132)$ as this will be necessary for the main result in Chapter \ref{chap:6}.\\

As it is discussed at the end of the proof for Lemma \ref{lem:3.2}, every permutation $\pi\in\Av(4123,\allowbreak 4213,4132)$ such that $|\pi|=n\geq 1$ and $\pi(n)=1$ can be written as $\tau\ominus 1$ for some length $n-1$ permutations $\tau$ in the same class. Except when $\pi=1$, $\tau$ is nonempty, so whether $\tau$ is skew-indecomposable or not, such a permutation $\pi$ of length $n\geq 2$ is skew-decomposable.\\

In addition to the above case, there are some permutations in $\Av(4123,4213,4132)$ ending with 12 (\textit{i.e.} $\pi(n-1)=1$ and $\pi(n)=2$). These are the permutations for which $\pi(n)=2$ and the segment $J$, as defined in the proof of Lemma \ref{lem:3.2}, is empty. Such permutations can be written as $\rho\ominus 12$ where $\rho$ is a nonempty permutation in $\Av(4123,4213,4132)$ of length $n-2$. Therefore, permutations of length $n\geq 3$ having $\pi(n-1)=1$ and $\pi(n)=2$ are also skew-decomposable.\\

We claim that these two are the only ways for a permutation in $\Av(4123,4213,4132)$ to be skew-decomposable. Suppose to the contrary that there exists a skew-decomposable permutation $\pi$ in $\Av(4123,4213,4132)$ which can be written as $\pi_1\ominus\pi_2$ where $\pi_2\neq 1$ and $\pi_2\neq 12$. If $\pi_2=21$, then this is just the case $\pi_2=1$, so the length of $\pi_2$ must be strictly greater than $2$. If $\pi_2$ has the value 2 in the last position, then $\pi_2$ must contain 132, but with $\pi_1$, this results in $\pi$ containing 4132, so we achieve a contradiction. Likewise, if $\pi_2$ has some value greater than or equal to 3 in the last position, then either $123\cont\pi_2$ or $213\cont\pi_2$, implying $4123\cont\pi$ or $4213\cont\pi$, which is again a contradiction. Consequently, a permutation $\pi$ in $\Av(4123,4213,4132)$ is skew-decomposable if and only if $\pi=\tau\ominus 1$ or $\pi=\rho\ominus 12$, where $\tau$ and $\rho$ are nonempty permutations in $\Av(4123,4213,4132)$, and $|\tau|=n-1$, $|\rho|=n-2$.\\

By excluding these two kinds of permutations, we can obtain the number of skew-indecomp-\allowbreak osable permutations in $\Av(4123,4213,4132)$ as the following.\vhhh
\begin{lemma}\label{lem:3.3}
\textit{Let $\bar{G}=\bar{f}_{\Av(4123,4213,4132)}$. The generating function for the number of skew-indecomposable permutations of length $n$ in $\Av(4123,4213,4132)$, excluding the empty permutation, is}\\
\[
(1-x-x^2)\bar{G}.
\]
\end{lemma}\vhhh
\textit{Proof.} Let $t_n$ be the number of skew-indecomposable permutations of length $n$ ($n\geq 1$) in $\Av(4123,4213,4132)$. As we have already discussed, all we have to do is take out the permutations having forms of $\pi=\tau\ominus 1$ or $\pi=\rho\ominus 12$, where $\tau$ and $\rho$ are nonempty permutations in $\Av(4123,4213,4132)$. The numbers of each are $s_{n-1}$ for $n\geq 2$ and $s_{n-2}$ for $n\geq 3$ respectively, so we have\\
\[
t_n=\left\{\begin{array}{cl}
s_n &\textrm{ if }n=1\\
s_n-s_{n-1} &\textrm{ if }n=2\\
s_n-s_{n-1}-s_{n-2} &\textrm{ if }n\geq 3
\end{array}\right..
\]\\
Therefore,\\
\[\begin{array}{rcl}
t_1 x+t_2 x^2+t_3 x^3+\cdots & = & s_1 x+(s_2-s_{1})x^2+(s_3-s_{2}-s_{1})x^3+\cdots\\
 & = & (s_1x+s_2x^2+\cdots)-(s_1x+s_2x^2+\cdots)x-(s_1x+s_2x^2+\cdots)x^2\\
 & = & \bar{G}-\bar{G}x-\bar{G}x^2=\bar{G}(1-x-x^2).
\end{array}
\]
\text{ }\hfill$\blacksquare$
\\
\newpage
\chapter{Enumeration of the class $\mathcal{A}$}\label{chap:4}
\section{Overview}
This chapter is entirely based on the paper \textit{Enumerating indices of Schubert varieties defined by inclusions} by Albert and Brignall \cite{AB2013}. We repeat and expand upon the details of that paper here since the methods we use to enumerate $\A'$ are an extension of theirs.\\

To begin, we give a short overview of their method of enumeration. They first enumerate simple permutations in $\A$. To do this, they characterize the structure of the simple permutations in $\A$ of length greater than or equal to 4. Using this characterization, they define an encoding of simple permutations in $\A$ into words. They then construct an automaton that accepts precisely these words to show that the set of encoded words form a regular language and apply the transfer matrix method to complete the enumeration of simple permutations. The whole class is enumerated by applying Proposition \ref{prop:2.4} and adding in the case of sum and skew decomposable permutations.\\

Recall $\A=\Av\{4231,35142,42513,351624\}$. Before we start, we note two symmetry properties of the class $\A$, introduce new terminology called the extreme pattern of a permutation, and give the general idea of the enumeration method. Let $\textrm{op}_1$ and $\textrm{op}_2$ be the inverse operation and the reverse complement operation, \textit{i.e.} $\textrm{op}_1(\pi)=\pi^{-1}$ and $\textrm{op}_2(\pi)=(\pi^r)^c$ for every permutation $\pi$. Notice that $\textrm{op}_1(\{4231,35142,42513,351624\})=\textrm{op}_2(\{4231,35142,42513,351624\})=\{4231,35142,42513,\allowbreak 351624\}$. Thus, by Proposition \ref{prop:2.2}, if $\pi$ is in the class $\A$, then $\pi^{-1}$ and $(\pi^{r})^{c}$ are also in $\A$.\\

Now, we define the extreme pattern of a permutation. The \textit{extreme pattern} of a permutation is the flattening of the first, the last, the greatest and the least values of a permutation. For instance, the extreme pattern of $\pi=47128365$ is $2143$ due to the subsequence $4185$ where 4, 1, 8 and 5 correspond to the first, the least, the greatest and the last value respectively. It is not always the case that the extreme pattern is length 4, as $\s=52413$ has the greatest value 5 in the first position, so its extreme pattern is $312$ due to the subsequence $513$. However, if $\pi$ is simple and $|\pi|\geq 4$, then its extreme pattern must be one of 2143, 2413, 3142 and 3412, since a simple permutation cannot begin or end with its greatest or least value.\\
\section{Extreme patterns 2413, 3142 and 3412}
In the next section, we will establish the structure of the simple permutations $\pi$ in $\A$ with $|\pi|\geq 4$ and $\pi(2)\neq 1$. Before we do so, we first study some special cases defined by their extreme patterns. We start with the simple permutations having extreme pattern 2413.\vhhh
\begin{proposition}\label{prop:4.1}
\textit{Let $\pi$ be a simple permutation in $\A$ with extreme pattern 2413. Let $b$, $d$, $a$ and $c$ be the first, the greatest, the least and the last values of $\pi$ respectively. Then the graph of $\pi$ is $N$-shaped, that is, values corresponding to positions in $[\pi^{-1}(b),\pi^{-1}(d)]$ are increasing, values corresponding to positions in  $[\pi^{-1}(d),\pi^{-1}(a)]$ are decreasing and values corresponding to positions in  $[\pi^{-1}(a),\pi^{-1}(c)]$ are increasing.}
\end{proposition}

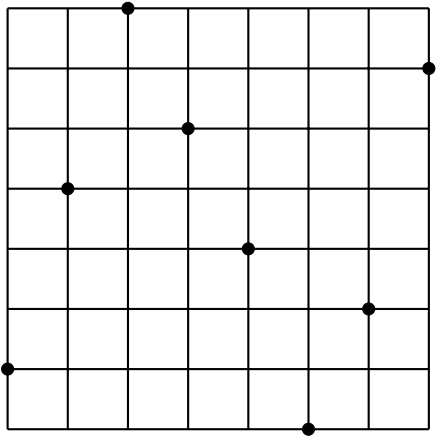
\begin{figure}[H]
\[
  \begin{tikzpicture}[scale=0.8]
    \plotPerm{2,5,8,6,4,1,3,7}
  \end{tikzpicture}
\]
\caption{Graph of 25864137.}
\label{fig:4.1}
\end{figure}

\begin{figure}
\[
  \begin{tikzpicture}[scale=0.8]
    \plotPerm{2,4,1,3}
    \node[left] at (1,2) {$b$};
    \node[above] at (2,4) {$d$};
    \node[below] at (3,1) {$a$};
    \node[right] at (4,3) {$c$};
    \node at (1.5,1.5) {$B_{31}$};
    \node at (3.5,3.5) {$B_{13}$};
  \end{tikzpicture}
\]
\caption{Partial graph of $\pi$ of extreme pattern 2413.}
\label{fig:4.2}
\end{figure}
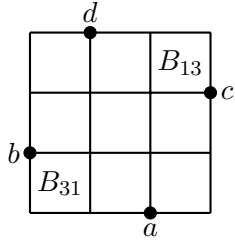

For instance, Figure \ref{fig:4.1} shows the graph of $\pi=25864137$, a length 8 simple permutation of extreme pattern 2413 in $\A$. As we can see, drawn points form an $N$-shape.\\
\noindent
\textit{Proof.} Let $\pi\in\A$ be a simple permutation of the extreme pattern 2413. We provide the graph of extreme pattern 2413 in Figure \ref{fig:4.2}. A permutation $\pi$ of extreme pattern 2413 has a graph that can be drawn by filling in more points in the interior regions of this graph. We first claim that, for $\pi\in\textrm{Si}(\A)$, there cannot be any points in the region denoted by $B_{31}$. In other words, $\pi$ cannot have a value less than $b$ whose position is in the segment $(\pi^{-1}(b),\pi^{-1}(d))$.\\

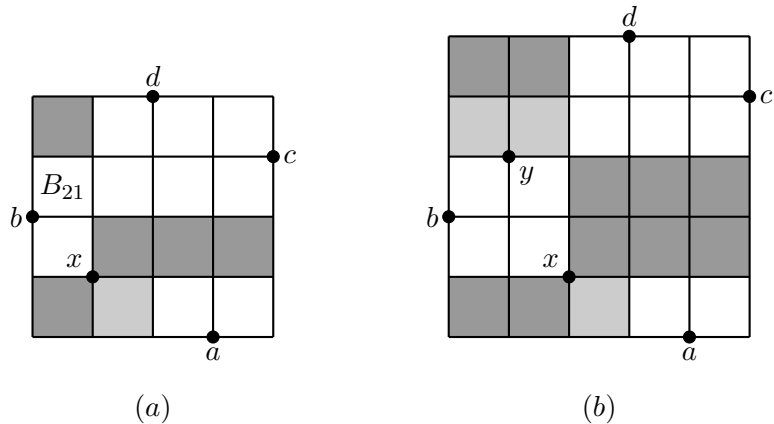
\begin{figure}[H]
\[\begin{array}{ccc}
  \begin{tikzpicture}[scale=0.8]
    \draw[white, fill=gray!80](1,1)--(2,1)--(2,2)--(1,2)--cycle;
    \draw[white, fill=gray!80](2,2)--(5,2)--(5,3)--(2,3)--cycle;
    \draw[white, fill=gray!80](1,4)--(2,4)--(2,5)--(1,5)--cycle;
    \draw[white, fill=gray!40](2,1)--(2,2)--(3,2)--(3,1)--cycle;
    \plotPerm{3,2,5,1,4}
    \node[left] at (1,3) {$b$};
    \node[above left] at (2,2) {$x$};
    \node[above] at (3,5) {$d$};
    \node[below] at (4,1) {$a$};
    \node[right] at (5,4) {$c$};
    \node at (1.5,3.5) {$B_{21}$};
  \end{tikzpicture} & \qquad\qquad & \begin{tikzpicture}[scale=0.8]
    \draw[white, fill=gray!80](1,1)--(3,1)--(3,2)--(1,2)--cycle;
    \draw[white, fill=gray!80](3,2)--(6,2)--(6,4)--(3,4)--cycle;
    \draw[white, fill=gray!80](1,5)--(3,5)--(3,6)--(1,6)--cycle;
    \draw[white, fill=gray!40](3,1)--(3,2)--(4,2)--(4,1)--cycle;
    \draw[white, fill=gray!40](1,4)--(3,4)--(3,5)--(1,5)--cycle;
    \plotPerm{3,4,2,6,1,5}
    \node[left] at (1,3) {$b$};
    \node[below right] at (2,4) {$y$};
    \node[above left] at (3,2) {$x$};
    \node[above] at (4,6) {$d$};
    \node[below] at (5,1) {$a$};
    \node[right] at (6,5) {$c$};
  \end{tikzpicture}\\
(a) & & (b)\end{array}\]
\caption{Partial graphs of $\pi$ with the assumption of having a value in $B_{31}$.}
\label{fig:4.3}
\end{figure}

Suppose to the contrary that there is a point in the region $B_{31}$. Let $x$ be the point in $B_{31}$ of minimum value as shown in Figure \ref{fig:4.3}$(a)$. We further examine where other values of $\pi$ can be located. Notice that if $\pi$ has a point in any of the dark grey regions, $\pi$ will contain some permutation in $\{4231,35142,42513,351624\}$, so these regions cannot contain any points. For example, if there exists a point $z$ in the dark grey region $(\pi^{-1}(x),\pi^{-1}(d))\times(x,b)$, then $bxza$ would form a 4231 pattern. In addition, there are no points in the light grey region because we chose $x$ to be the point of minimum value with position in the segment $(\pi^{-1}(b),\pi^{-1}(d))$. We use this two-color coding to differentiate forbidden regions in future proofs also.\\
Since $\pi$ is a simple permutation, the segment defined by the positions of $b$ and $x$ cannot be a block. Hence, there exists a point in the region denoted by $B_{21}$. (We will say such a point \textit{splits} the potential block $[\pi^{-1}(b),\pi^{-1}(x)]$.) Choose the point of greatest value in $B_{21}$ and call this value $y$. Then we obtain the graph shown in Figure \ref{fig:4.3}$(b)$. Notice now that $[\pi^{-1}(b),\pi^{-1}(x)]$ must be a block, since there cannot be any points directly above, below, to the right or to the left of the region $[\pi^{-1}(b),\pi^{-1}(x)]\times[x,y]$. Therefore, $\pi$ is not simple, a contradiction.

Referring back to Figure \ref{fig:4.2}, since $\A$ is preserved by reverse complement, we can rotate our previous argument by $180^\circ$ to show there are no points in the region $B_{13}$.\\

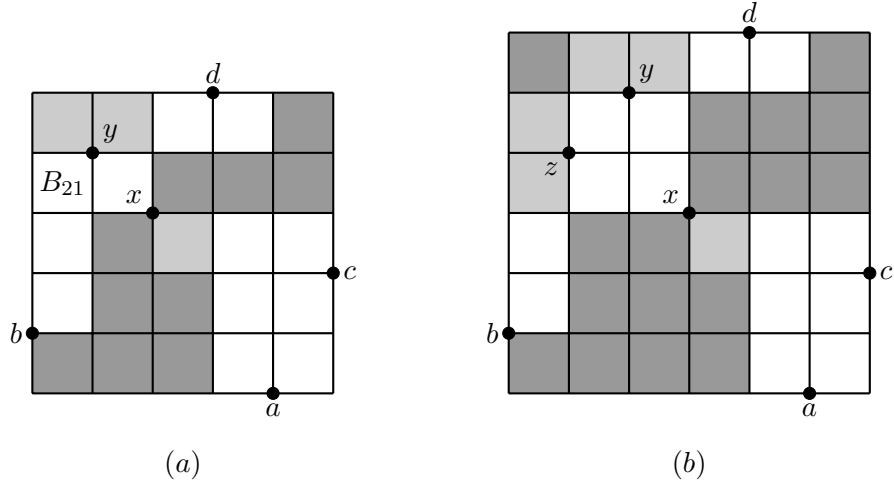
\begin{figure}
\[\begin{array}{ccc}
  \begin{tikzpicture}[scale=0.8]
    \draw[white, fill=gray!80](2,1)--(4,1)--(4,4)--(2,4)--cycle;
    \draw[white, fill=gray!80](3,4)--(6,4)--(6,5)--(3,5)--cycle;
    \draw[white, fill=gray!80](5,5)--(6,5)--(6,6)--(5,6)--cycle;
    \draw[white, fill=gray!80](1,1)--(2,1)--(2,2)--(1,2)--cycle;
    \draw[white, fill=gray!40](1,5)--(3,5)--(3,6)--(1,6)--cycle;
    \draw[white, fill=gray!40](3,3)--(4,3)--(4,4)--(3,4)--cycle;
    \plotPerm{2,5,4,6,1,3}
    \node[left] at (1,2) {$b$};
    \node[above right] at (2,5) {$y$};
    \node[above left] at (3,4) {$x$};
    \node[above] at (4,6) {$d$};
    \node[below] at (5,1) {$a$};
    \node[right] at (6,3) {$c$};
    \node at (1.5,4.5) {$B_{21}$};
  \end{tikzpicture} & \qquad\qquad &
  \begin{tikzpicture}[scale=0.8]
    \draw[white, fill=gray!80](2,1)--(5,1)--(5,4)--(2,4)--cycle;
    \draw[white, fill=gray!80](4,4)--(7,4)--(7,6)--(4,6)--cycle;
    \draw[white, fill=gray!80](6,6)--(7,6)--(7,7)--(6,7)--cycle;
    \draw[white, fill=gray!80](1,1)--(2,1)--(2,2)--(1,2)--cycle;
    \draw[white, fill=gray!80](1,6)--(2,6)--(2,7)--(1,7)--cycle;
    \draw[white, fill=gray!40](2,6)--(4,6)--(4,7)--(2,7)--cycle;
    \draw[white, fill=gray!40](4,3)--(5,3)--(5,4)--(4,4)--cycle;
    \draw[white, fill=gray!40](1,4)--(2,4)--(2,6)--(1,6)--cycle;
    \plotPerm{2,5,6,4,7,1,3}
    \node[left] at (1,2) {$b$};
    \node[below left] at (2,5) {$z$};
    \node[above right] at (3,6) {$y$};
    \node[above left] at (4,4) {$x$};
    \node[above] at (5,7) {$d$};
    \node[below] at (6,1) {$a$};
    \node[right] at (7,3) {$c$};
  \end{tikzpicture}\\
(a) & & (b)\end{array}\]
\caption{Partial graphs of $\pi$ with the assumption of having a decreasing sub-segment in $[\pi^{-1}(b),\pi^{-1}(d)]$.}
\label{fig:4.4}
\end{figure}

Now, we show that values corresponding to positions in $[\pi^{-1}(b),\pi^{-1}(d)]$ are increasing. We again show this by contradiction using the graph of $\pi$. Suppose the values corresponding to the segment $[\pi^{-1}(b),\pi^{-1}(d)]$ are not strictly increasing. This means that there is at least one sub-segment in $[\pi^{-1}(b),\pi^{-1}(d)]$ whose values are decreasing. Let $[\pi^{-1}(y),\pi^{-1}(x)]$ be the decreasing sub-segment with the greatest possible $y$ and the least possible $x$ given $y$. This provides the graph shown in Figure \ref{fig:4.4}$(a)$. Since the segment $[\pi^{-1}(y),\pi^{-1}(x)]$ cannot be a block, we have a point in the region $B_{21}$ of Figure \ref{fig:4.4}$(a)$. By choosing the left-most such point and denoting by $z$ the value of this point, we obtain the graph in Figure \ref{fig:4.4}$(b)$. Since we have a block $[\pi^{-1}(z),\pi^{-1}(x)]$ which cannot be split, we achieve a contradiction.\\

The reverse complement of this argument shows values corresponding to positions in $[\pi^{-1}(a),\pi^{-1}(c)]$ are strictly increasing.\\

Finally, values corresponding to the segment $[\pi^{-1}(d),\pi^{-1}(a)]$ must be decreasing because, otherwise, $\pi$ would contain $4231$ with $d$ and $a$ corresponding to the $4$ and $1$ respectively. This completes the proof.\hfill$\blacksquare$\\

If we apply the inverse symmetry to the previous proof, we obtain the following result for a simple permutation of extreme pattern 3142.\vhhh
\begin{proposition}\label{prop:4.2}
\textit{
Given a simple permutation $\pi$ in $\A$ with extreme pattern 3142. Let $c$, $a$, $d$ and $b$ be the first, the least, the greatest and the last values of $\pi$ respectively. Then the graph of $\pi$ is $S$-shaped, that is, values from the range $[a,b]$ are increasing, values from the range $[b,c]$ are decreasing and values from the range $[c,d]$ are increasing.}
\end{proposition}

\begin{figure}[H]
\[\begin{tikzpicture}[scale=0.8]
    \draw[white, fill=gray!80](1,1)--(2,1)--(2,2)--(1,2)--cycle;
    \draw[white, fill=gray!80](3,3)--(3,4)--(4,4)--(4,3)--cycle;
    \plotPerm{3,4,1,2}
    \node[left] at (1,3) {$c$};
    \node[above] at (2,4) {$d$};
    \node[below] at (3,1) {$a$};
    \node[right] at (4,2) {$b$};
    \node at (1.5,2.5) {$B_{21}$};
    \node at (2.5,3.5) {$B_{12}$};
  \end{tikzpicture}
\]
\caption{Partial graph of $\pi$ of extreme pattern 3412.}
\label{fig:4.5}
\end{figure}
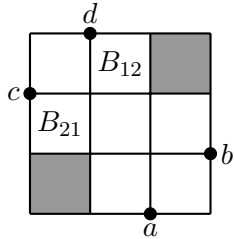

Lastly, we prove the following proposition.\vhhh
\begin{proposition}\label{prop:4.3}
\textit{No simple permutation in $\A$ has extreme pattern 3412.}
\end{proposition}
\textit{Proof.} Suppose the statement is false, and let $\pi$ be a simple permutation in $\A$ whose extreme pattern is 3412. Just as in the proof of Proposition \ref{prop:4.1}, we start with the graph of extreme pattern 3412 with $c$, $d$, $a$ and $b$ representing the first, the greatest, the least and the last values respectively. As shown in Figure \ref{fig:4.5}, the segment $[\pi^{-1}(c),\pi^{-1}(d)]$ would form a block without the presence of a point in either $B_{21}$ or $B_{12}$.\\

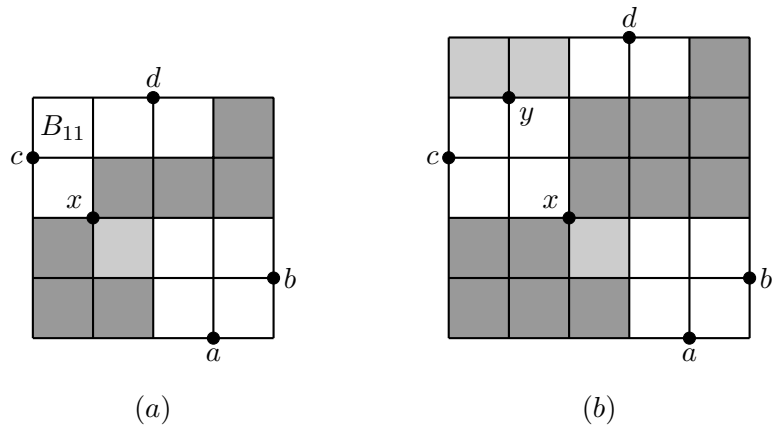
\begin{figure}[H]
\[\begin{array}{ccc}
  \begin{tikzpicture}[scale=0.8]
    \draw[white, fill=gray!80](1,1)--(3,1)--(3,2)--(2,2)--(2,3)--(1,3)--cycle;
    \draw[white, fill=gray!80](2,3)--(5,3)--(5,5)--(4,5)--(4,4)--(2,4)--cycle;
    \draw[white, fill=gray!40](2,2)--(2,3)--(3,3)--(3,2)--cycle;
    \plotPerm{4,3,5,1,2}
    \node[left] at (1,4) {$c$};
    \node[above left] at (2,3) {$x$};
    \node[above] at (3,5) {$d$};
    \node[below] at (4,1) {$a$};
    \node[right] at (5,2) {$b$};
    \node at (1.5,4.5) {$B_{11}$};
  \end{tikzpicture} & \qquad\qquad & \begin{tikzpicture}[scale=0.8]
    \draw[white, fill=gray!80](1,1)--(4,1)--(4,2)--(3,2)--(3,3)--(1,3)--cycle;
    \draw[white, fill=gray!80](3,3)--(6,3)--(6,6)--(5,6)--(5,5)--(3,5)--cycle;
    \draw[white, fill=gray!40](3,2)--(4,2)--(4,3)--(3,3)--cycle;
    \draw[white, fill=gray!40](1,5)--(3,5)--(3,6)--(1,6)--cycle;
    \plotPerm{4,5,3,6,1,2}
    \node[left] at (1,4) {$c$};
    \node[below right] at (2,5) {$y$};
    \node[above left] at (3,3) {$x$};
    \node[above] at (4,6) {$d$};
    \node[below] at (5,1) {$a$};
    \node[right] at (6,2) {$b$};
  \end{tikzpicture}\\
(a) & & (b)\end{array}\]
\caption{Partial graphs of $\pi$ with the assumption of having a value in $B_{21}$.}
\label{fig:4.6}
\end{figure}

Assume there is a point in $B_{21}$, and let $x$ denote the minimum value of any point in $B_{21}$. This is shown in Figure \ref{fig:4.6}$(a)$. Now, we must have a point in $B_{11}$ in Figure \ref{fig:4.6}$(a)$ to prevent the segment $[\pi^{-1}(c),\pi^{-1}(x)]$ from being a block. Let $y$ be the greatest value of any point in $B_{11}$, we obtain the graph in Figure \ref{fig:4.6}$(b)$. Now $[\pi^{-1}(c),\pi^{-1}(x)]$ must be a block, which is a contradiction.\\

If we assume the existence of a point in $B_{12}$ in Figure \ref{fig:4.5}, we end up with the same result, since $\A$ and the pattern 3412 are preserved by the inverse operation followed by the reverse complement operation. Hence, there is no simple permutation of extreme pattern 3412 in $\A$.\hfill$\blacksquare$\\
\section{General simple permutations in $\A$.}
\subsection{Structure theorem}
We have characterized the structure of simple permutations having extreme patterns $2413$, $3142$ and $3412$. The remaining extreme pattern is $2143$. Instead of immediately examining this case, we first discuss the structure of simple permutations $\pi$ where $|\pi|\geq 4$ and $\pi(2)\neq 1$. Later, we show that the condition $\pi(2)\neq 1$ is equivalent to the condition $\pi(1)=2$. Thus, every simple permutation in $\A$ has either the value 2 in the first position or the value $1$ in the second position, but not both, since a permutation with both is not simple. Observing that any simple permutation of extreme pattern 2413 has 2 in the first position and any permutation of extreme pattern 3142 has 1 in the second position, we can provide the classification of simple permutations as shown in Table \ref{tab:4.1}. For the remainder of this chapter, let $H=\{\pi\in\textrm{Si}(\A):|\pi|\geq 4\textrm{ and }\pi(2)\neq 1\}$.\\
\begin{table}[H]
\begin{center}
\begin{tabular}{cccc}
\phantom{Column 1} & \phantom{Column 1} & \phantom{Column 1} & \phantom{Column 1}\\
\multicolumn{4}{c}{Simple permutation $\pi$ in $\A$ with $|\pi|\geq 4$}\\
\hline
\multicolumn{2}{c|}{$\pi(1)=2$} & \multicolumn{2}{c}{$\pi(2)=1$}\\
\hline
Extreme pattern 2413 & \multicolumn{2}{|c|}{Extreme pattern 2143} & Extreme pattern 3142\\
\end{tabular}
\end{center}
\caption{Classification of simple permutations in $\A$.}
\label{tab:4.1}
\end{table}

In order to describe the structure of $\pi$ in $H$, we need to define four special sum-like operations which we call glue sums that combine two permutations satisfying certain conditions into a longer permutation. They are called the type 1-0 NW glue sum, the type 1-1 NW glue sum, the type 1-0 SE glue sum and the type 1-1 SE glue sum. In Chapter \ref{chap:5}, we will define more operations that are similar to these four, and type numbering will be explained there.\\

We first define the type 1-0 NW glue sum. Let $\s$ and $\tau$ be simple permutations in $\A$ of length $m$ and $n$ respectively. Let $i=\s^{-1}(m)$ and $j=\tau(1)$. Furthermore, suppose $i\leq m-2$, $\s(m)=m-1$, $j\geq 3$ and $\tau(2)=1$. For $\s$ and $\tau$ satisfying these conditions, we define the \textit{type 1-0 NW glue sum}, denoted by $\s\nwsum_1^0\tau$, as the following permutation.
\[
\s\nwsum_1^0\tau=\s'(1)\s'(2)\cdots\s'(m-1)\tau'(3)\tau'(4)\cdots\tau'(n),
\]
where $\s'(i)=m+(j-3)$ and $\s'(k)=\s(k)$ for $k\neq i$, and $\tau'(k)=\tau(k)+(m-3)$ for $k$ with $3\leq k\leq n$.\\

The type 1-0 NW glue sum identifies the greatest value $m$ in $\s$ and the first value $j$ in $\tau$. These two points are combined into one with the position $\s^{-1}(m)$ and the value $j+m-3$, which is the value $j$ shifted up by $m-3$ as all other points of $\tau$ are shifted up by $m+3$. The type 1-0 NW glue sum also eliminates $\s(m)=m-1$ and $\tau(2)=1$. The remaining values in $\tau$ are attached just as in the usual sum $\oplus$ except they are shifted up by $m-3$. Notice that $|\s\nwsum_1^0\tau|=m+n-3$, since one pair of values is combined and two values are eliminated.\\

Any simple permutations in $\A$ of extreme pattern 2413 and simple permutations in $\A$ of extreme pattern 3142 respectively satisfy the conditions for $\s$ and $\tau$ required in the definition of $\nwsum_1^0$. Figure \ref{fig:4.15} illustrates the type 1-0 NW glue sum with $\s_1=2753146$ and $\tau_1=5162473$.\\
\begin{figure}[b!]
\[
\begin{array}{ccc}
\begin{tikzpicture}[scale=0.5]
\draw[color=black,fill=black] (1,2) circle [radius=0.08];
\draw[color=black,fill=black] (2,7) circle [radius=0.08];
\draw[color=black,fill=black] (3,5) circle [radius=0.08];
\draw[color=black,fill=black] (4,3) circle [radius=0.08];
\draw[color=black,fill=black] (5,1) circle [radius=0.08];
\draw[color=black,fill=black] (6,4) circle [radius=0.08];
\draw[color=black,fill=black] (7,6) circle [radius=0.08];
\draw[thick,dashed] (8,0) -- (8,12);
\node at (4,-.5) {$\s_1=2753146$};
\node at (12,12.5) {$\tau_1=5162473$};
\draw[color=black,fill=black] (9,9) circle [radius=0.08];
\draw[color=black,fill=black] (10,5) circle [radius=0.08];
\draw[color=black,fill=black] (11,10) circle [radius=0.08];
\draw[color=black,fill=black] (12,6) circle [radius=0.08];
\draw[color=black,fill=black] (13,8) circle [radius=0.08];
\draw[color=black,fill=black] (14,11) circle [radius=0.08];
\draw[color=black,fill=black] (15,7) circle [radius=0.08];
\draw[thick,black] (6.7,6.3) -- (7.3,5.7);
\draw[thick,black] (9.7,5.3) -- (10.3,4.7);
\draw[arrows=->,thick,black] (2,7.2) -- (2,8.8);
\draw[arrows=->,thick,black] (8.8,9) -- (2.2,9);
\end{tikzpicture}
& \quad\raisebox{3cm}{$\rightarrow$}\quad &
\raisebox{-2.8pt}{\begin{tikzpicture}[scale=0.5]
\draw[color=black,fill=black] (1,2) circle [radius=0.08];
\draw[color=black,fill=black] (2,9) circle [radius=0.08];
\draw[color=black,fill=black] (3,5) circle [radius=0.08];
\draw[color=black,fill=black] (4,3) circle [radius=0.08];
\draw[color=black,fill=black] (5,1) circle [radius=0.08];
\draw[color=black,fill=black] (6,4) circle [radius=0.08];
\draw[color=black,fill=black] (7,10) circle [radius=0.08];
\draw[color=black,fill=black] (8,6) circle [radius=0.08];
\draw[color=black,fill=black] (9,8) circle [radius=0.08];
\draw[color=black,fill=black] (10,11) circle [radius=0.08];
\draw[color=black,fill=black] (11,7) circle [radius=0.08];
\node at (6,-.5) {$\s_1\nwsum_1^0\tau_1=2$ 9 5 3 1 4 10 6 8 11 7};
\end{tikzpicture}}
\end{array}
\]
\caption{Illustration of $\s_1\nwsum_1^0\tau_1$.}
\label{fig:4.15}
\end{figure}

Next, we define the type 1-1 NW glue sum. Let $\s$ and $\tau$ be simple permutations in $\A$ satisfying the same conditions as in the definition of the type 1-0 NW glue sum. The \textit{type 1-1 NW glue sum} denoted by $\s\nwsum_1^1\tau$ is defined as the following.\\
\[
\s\nwsum_1^1\tau=\s'(1)\s'(2)\cdots\s'(m)\tau'(3)\tau'(4)\cdots\tau'(n)
\]
where $\s'(i)=m+(j-2)$ and $\s'(k)=\s(k)$ for $k\neq i$, and $\tau'(k)=\tau(k)+(m-2)$ for $k$ with $3\leq k\leq n$. The only differences between $\s\nwsum_1^1\tau$ and $\s\nwsum_1^0\tau$ are that $\s\nwsum_1^1\tau$ has an extra value $\s(m)=m-1$ between the subsequences coming from $\s$ and $\tau$ (as $m-1$ is not deleted), so $\s(i)$ and $\tau(k)$ ($3\leq k\leq n$) are shifted up by $m-2$ instead of $m-3$. Therefore, $|\s\nwsum_1^1\tau|=m+n-2$. For $\s=2753146$ and $\tau=5162473$, we obtain\[
\s\nwsum_1^1\tau=2\textrm{ 10 5 3 1 4 6 11 7 9 12 8.}
\]

Note that, if we restrict to sum simple permutations of extreme pattern 2413 and 3142 together, both NW glue sums are injective operations. To show this for $\nwsum_1^0$, let $\s_1$ and $\s_2$ be simple permutations of extreme pattern 2413 and $\tau_1$ and $\tau_2$ be simple permutations of extreme pattern 3142. Let $\pi_1=\s_1\nwsum_0^1\tau_1$ and $\pi_2=\s_2\nwsum_0^1\tau_2$. Suppose $\pi_1=\pi_2$. If $|\s_1|=|\s_2|$, then it is obvious that $\s_1=\s_2$ and $\tau_1=\tau_2$ by definition. So suppose $|\s_1|<|\s_2|$. Let $|\s_1|=m$ and $|\s_2|=m+1$. Then $\pi_1(m)=\tau'_1(3)=\s'_2(m)=\pi_2(m)$. Because $\tau_1(2)=1$ by Proposition \ref{prop:4.2}, $\tau_1(3)$ must be greater than $\tau_1(n)$ ($n=|\tau_1|$), because, otherwise, due to Proposition \ref{prop:4.2}, we must have $\tau_1(3)=2$, which is a contradiction to $\tau_1$ being simple. Hence, we obtain $\s'_2(m)>\pi_2(m+n-3)$. This is impossible, while $\s'_2(m)=\s_2(m)<m-1$ and $\pi_2(m+n-3)$ is at least $m-1$ by definition of $\nwsum_1^0$. We have the same contradiction for the case $|\s_2|>m+1$, so $|\s_1|=|\s_2|$, implying the type 1-0 NW glue sum is injective. The same argument can be applied for the type 1-1 NW glue sum as well.\\

The type 1-0 and type 1-1 NW glue sums both combine the greatest value in $\s$ and the first value in $\tau$. Next, we define the inverse notions of these two sums, called the type 1-0 and type 1-1 SE glue sums, denoted by $\sesum_1^0$ and $\sesum_1^1$, respectively, with the property that\[
(\s\sesum_1^0\tau)^{-1}=\s^{-1}\nwsum_1^0\tau^{-1}\quad\textrm{and}\quad(\s\sesum_1^1\tau)^{-1}=\s^{-1}\nwsum_1^0\tau^{-1}.
\]
Specifically, let $\s$ and $\tau$ be simple permutations in $\A$ of length $m$ and $n$ respectively. Let $i=\s(m)$ and $j=\tau^{-1}(1)$. We require that $i\leq m-2$, $m=\s(m-1)$, $j\geq 3$ and $2=\tau(1)$. For $\s$ and $\tau$ satisfying these conditions, we define the \textit{type 1-0 SE glue sum} as the following.
\[
\s\sesum_1^0\tau=\s(1)\s(2)\cdots\s(m-2)\tau'(2)\tau'(3)\cdots\tau'(n)
\]
where $\tau'(j)=i$ and $\tau'(k)=\tau(k)+(m-3)$ for $k\neq j$. Similarly, the \textit{type 1-1 SE glue sum} is defined as
\[
\s\sesum_1^1\tau=\s(1)\s(2)\cdots\s(m-1)\tau'(2)\tau'(3)\cdots\tau'(n)
\]
where $\tau'(j)=i$ and $\tau'(k)=\tau(k)+(m-2)$ for $k\neq j$. Both SE glue sums are injective when we sum simple permutations of extreme pattern 3142 and 2413.\\

Figure \ref{fig:4.16} shows $\s_2\sesum_1^0\tau_2$ where $\s_2=5146372$ and $\tau_2=2475136$. Notice that $\s_2^{-1}=\s_1$ and $\tau_2^{-1}=\tau_1$ from the previous example. Indeed, $\s_2\sesum_1^0\tau_2=5$ 1 4 6 3 8 11 9 2 7 10, which is the inverse of $\s_1\nwsum_1^0\tau_1=2$ 9 5 3 1 4 10 6 8 11 7.\\

Note that glue sums we have defined so far are associative operations only if the lengths of all summands are at least 4. For our purpose, we set the convention that when we sum permutations with multiple glue sums, we always ensure to operate from left to right.\\

We are now ready to state the theorem for the simple permutations in $H$.\begin{figure}
\[
\begin{array}{ccc}
\begin{tikzpicture}[scale=0.5]
\draw[color=black,fill=black] (1,5) circle [radius=0.08];
\draw[color=black,fill=black] (2,1) circle [radius=0.08];
\draw[color=black,fill=black] (3,4) circle [radius=0.08];
\draw[color=black,fill=black] (4,6) circle [radius=0.08];
\draw[color=black,fill=black] (5,3) circle [radius=0.08];
\draw[color=black,fill=black] (6,7) circle [radius=0.08];
\draw[color=black,fill=black] (7,2) circle [radius=0.08];
\draw[thick,dashed] (0,8) -- (12,8);
\node at (4,-.5) {$\s_2=5146372$};
\node at (8,16.5) {$\tau_2=2475136$};
\draw[color=black,fill=black] (5,10) circle [radius=0.08];
\draw[color=black,fill=black] (6,12) circle [radius=0.08];
\draw[color=black,fill=black] (7,15) circle [radius=0.08];
\draw[color=black,fill=black] (8,13) circle [radius=0.08];
\draw[color=black,fill=black] (9,9) circle [radius=0.08];
\draw[color=black,fill=black] (10,11) circle [radius=0.08];
\draw[color=black,fill=black] (11,14) circle [radius=0.08];
\draw[thick,black] (5.7,7.3) -- (6.3,6.7);
\draw[thick,black] (4.7,10.3) -- (5.3,9.7);
\draw[arrows=->,thick,black] (7.2,2) -- (8.8,2);
\draw[arrows=->,thick,black] (9,8.8) -- (9,2.2);
\end{tikzpicture}
& \quad\raisebox{4.4cm}{$\rightarrow$}\quad &
\raisebox{1.2cm}{\begin{tikzpicture}[scale=0.5]
\draw[color=black,fill=black] (1,5) circle [radius=0.08];
\draw[color=black,fill=black] (2,1) circle [radius=0.08];
\draw[color=black,fill=black] (3,4) circle [radius=0.08];
\draw[color=black,fill=black] (4,6) circle [radius=0.08];
\draw[color=black,fill=black] (5,3) circle [radius=0.08];
\draw[color=black,fill=black] (6,8) circle [radius=0.08];
\draw[color=black,fill=black] (7,11) circle [radius=0.08];
\draw[color=black,fill=black] (8,9) circle [radius=0.08];
\draw[color=black,fill=black] (9,2) circle [radius=0.08];
\draw[color=black,fill=black] (10,7) circle [radius=0.08];
\draw[color=black,fill=black] (11,10) circle [radius=0.08];
\node at (6,-.5) {$\s_2\sesum_1^0\tau_2=5$ 1 4 6 3 8 11 9 2 7 10};
\end{tikzpicture}}
\end{array}
\]
\caption{Illustration of $\s_2\sesum_1^0\tau_2$.}
\label{fig:4.16}
\end{figure}
\begin{theorem}\label{thm:4.1}
\textit{Let $\pi$ be a permutation in $H$. Then there exist simple permutations in $\A$ of extreme pattern 2413 $\s_i$ ($i$ odd) and simple permutations in $\A$ of extreme pattern 3142 $\tau_i$ ($i$ even) such that\begin{equation}\label{eqn:4.1}
\pi=\left\{\begin{array}{clc}
\s_1\nwsum_1^{k_1}\tau_2\sesum_1^{k_2}\s_3\nwsum_1^{k_3}\tau_4\sesum_1^{k_4}\cdots\sesum_1^{k_{m-1}}\s_m & \textrm{ if }m\textrm{ is odd} &\qquad (a)\vhhh\\
\s_1\nwsum_1^{k_1}\tau_2\sesum_1^{k_2}\s_3\nwsum_1^{k_3}\tau_4\sesum_1^{k_4}\cdots\nwsum_1^{k_{m-1}}\tau_m & \textrm{ if }m\textrm{ is even} &\qquad (b)
\end{array}\right.
\end{equation}\\
where $m$ is a positive integer and $k_\l\in\{0,1\}$ ($1\leq\l\leq m-1$). Hence, $\pi$ has one of the structures illustrated in Figure \ref{fig:4.7}. Moreover, every simple permutation of these forms is in $H$.}
\end{theorem}
\vhhh

In both Equation \ref{eqn:4.1}$(a)$ and \ref{eqn:4.1}$(b)$, we always make sure to sum permutations from left to right. As glue sums are injective operations when we sum simple permutations of 2413 and 3142, every glue sum in Equation \ref{eqn:4.1}$(a)$ and \ref{eqn:4.1}$(b)$ are also injective.\\

We give the description of Figure \ref{fig:4.7}. Each point of $\pi$ is either one of the isolated points denoted by $d_i$ or located on one of the sequences of lines. These isolated points are the ones identified by NW and SE glue sums except the first, the least, the greatest and the last ones. The precise identification of each $d_i$ is given in the proof. We call the sequence of diagonal lines shown in Figure \ref{fig:4.7} a \textit{crenellation} according to its suggested shape. This structure can be arbitrarily long depending on the length of $\pi$. Each line segment of the crenellation can contain any number of points or be empty, but points must be placed so that $\pi$ avoids having blocks. Unlike the lines of the crenellation, the isolated points must be present so long as the structure continues.\\
\begin{figure}
\[
\begin{array}{ccc}
\begin{tikzpicture}[scale=0.8]
\path[draw,thick] (-1.5,-0.5) -- (-1.05,-0.05);
\path[draw,thick] (-1,-0.05) -- (-0.55,-0.5);
\path[draw,thick] (-0.5,-0.5) -- (-0.05,-0.05);
\path[draw,thick] (0,0) -- (0.45,0.45);
\path[draw,thick] (0.45,0.5) -- (0,0.95);
\path[draw,thick] (0,1) -- (0.45,1.45);
\path[draw,thick] (0.5,1.5) -- (0.95,1.95);
\path[draw,thick] (1,1.95) -- (1.45,1.5);
\path[draw,thick] (1.5,1.5) -- (1.95,1.95);
\draw[dotted,ultra thick] (2.5,3) -- (2.9,3.4);
\draw[color=black,fill=black] (-2,-1) circle [radius=.05] node [left]{$d_1$};
\draw[color=black,fill=black] (-0.525,-1.525) circle [radius=.05] node [right]{$d_2$};
\draw[color=black,fill=black] (-1.025,0.975) circle [radius=.05] node [left]{$d_3$};
\draw[color=black,fill=black] (1.475,0.475) circle [radius=.05] node [right]{$d_4$};
\draw[color=black,fill=black] (0.975,2.975) circle [radius=.05] node [left]{$d_5$};
\path[draw,thick] (-1.5+4.95,-0.5+4.95) -- (-1.05+4.95,-0.05+4.95);
\path[draw,thick] (-1+4.95,-0.05+4.95) -- (-0.55+4.95,-0.5+4.95);
\path[draw,thick] (-0.5+4.95,-0.5+4.95) -- (-0.05+4.95,-0.05+4.95);
\draw[color=black,fill=black] (-0.525+4.95,-1.525+4.95) circle [radius=.05] node [right]{$d_{m+1}$};
\draw[color=black,fill=black] (-1.025+4.95,0.975+4.95) circle [radius=.05] node [above]{$d_{m+2}$};
\draw[color=black,fill=black] (0.45+4.95,0.45+4.95) circle [radius=.05] node [above]{$d_{m+3}$};
\end{tikzpicture} & \quad\quad &\begin{tikzpicture}[scale=0.8]
\path[draw,thick] (-1.5,-0.5) -- (-1.05,-0.05);
\path[draw,thick] (-1,-0.05) -- (-0.55,-0.5);
\path[draw,thick] (-0.5,-0.5) -- (-0.05,-0.05);
\path[draw,thick] (0,0) -- (0.45,0.45);
\path[draw,thick] (0.45,0.5) -- (0,0.95);
\path[draw,thick] (0,1) -- (0.45,1.45);
\path[draw,thick] (0.5,1.5) -- (0.95,1.95);
\path[draw,thick] (1,1.95) -- (1.45,1.5);
\path[draw,thick] (1.5,1.5) -- (1.95,1.95);
\draw[dotted,ultra thick] (2.5,3) -- (2.9,3.4);
\draw[color=black,fill=black] (-2,-1) circle [radius=.05] node [left]{$d_1$};
\draw[color=black,fill=black] (-0.525,-1.525) circle [radius=.05] node [right]{$d_2$};
\draw[color=black,fill=black] (-1.025,0.975) circle [radius=.05] node [left]{$d_3$};
\draw[color=black,fill=black] (1.475,0.475) circle [radius=.05] node [right]{$d_4$};
\draw[color=black,fill=black] (0.975,2.975) circle [radius=.05] node [left]{$d_5$};
\path[draw,thick] (0+4.05,0+4.05) -- (0.45+4.05,0.45+4.05);
\path[draw,thick] (0.45+4.05,0.5+4.05) -- (0+4.05,0.95+4.05);
\path[draw,thick] (0+4.05,1+4.05) -- (0.45+4.05,1.45+4.05);
\draw[color=black,fill=black] (-1.025+4.05,0.975+4.05) circle [radius=.05] node [left]{$d_{m+2}$};
\draw[color=black,fill=black] (1.475+4.05,0.475+4.05) circle [radius=.05] node [above]{$d_{m+3}$};
\draw[color=black,fill=black] (0.95+4.05,1.95+4.05) circle [radius=.05] node [above]{$d_{m+4}$};
\end{tikzpicture}\\
(a) & & (b)
\end{array}
\]
\caption{Structure of $\pi$ with $|\pi|\geq 4$ and $\pi(2)\neq 1$.}
\label{fig:4.7}
\end{figure}
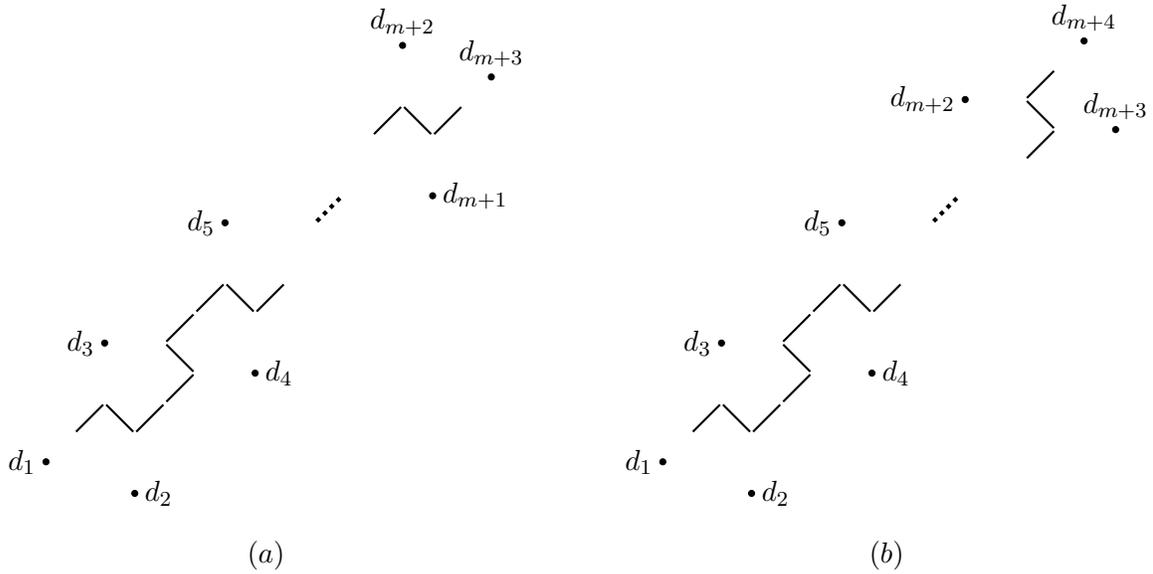
\begin{figure}[b!]
\[
\begin{tikzpicture}[scale=0.9]
\path[draw,thick] (-3.5,-2.5) -- (-3.05,-2.05);
\draw[color=black,fill=black] (-3.15,-2.15) circle [radius=.05];
\path[draw,thick] (-3,-2.05) -- (-2.55,-2.5);
\draw[color=black,fill=black] (-2.65,-2.4) circle [radius=.05];
\path[draw,thick] (-2.5,-2.5) -- (-2.05,-2.05);
\draw[color=black,fill=black] (-2.275,-2.275) circle [radius=.05];
\path[draw,thick] (-2,-2) -- (-1.55,-1.55);
\draw[color=black,fill=black] (-1.775,-1.775) circle [radius=.05];
\path[draw,thick] (-1.55,-1.5) -- (-2,-1.05);
\draw[color=black,fill=black] (-1.65,-0.65) circle [radius=.05];
\path[draw,thick] (-2,-1) -- (-1.55,-0.55);
\draw[color=black,fill=black] (-1.9,-1.15) circle [radius=.05];
\path[draw,thick] (-1.5,-0.5) -- (-1.05,-0.05);
\draw[color=black,fill=black] (-1.15,-0.15) circle [radius=.05];
\path[draw,thick] (-1,-0.05) -- (-0.55,-0.5);
\path[draw,thick] (-0.5,-0.5) -- (-0.05,-0.05);
\draw[color=black,fill=black] (-0.4,-0.4) circle [radius=.05];
\path[draw,thick] (0,0) -- (0.45,0.45);
\draw[color=black,fill=black] (0.225,0.225) circle [radius=.05];
\path[draw,thick] (0.45,0.5) -- (0,0.95);
\draw[color=black,fill=black] (0.35,0.6) circle [radius=.05];
\draw[color=black,fill=black] (0.1,0.85) circle [radius=.05];
\path[draw,thick] (0,1) -- (0.45,1.45);
\path[draw,thick] (0.5,1.5) -- (0.95,1.95);
\draw[color=black,fill=black] (0.725,1.725) circle [radius=.05];
\path[draw,thick] (1,1.95) -- (1.45,1.5);
\draw[color=black,fill=black] (1.35,1.6) circle [radius=.05];
\draw[color=black,fill=black] (1.1,1.85) circle [radius=.05];
\path[draw,thick] (1.5,1.5) -- (1.95,1.95);
\draw[color=black,fill=black] (-3.8,-2.8) circle [radius=.05];
\draw[color=black,fill=black] (-2.525,-3.3) circle [radius=.05];
\draw[color=black,fill=black] (-3.025,-1.025) circle [radius=.05];
\draw[color=black,fill=black] (-0.525,-1.525) circle [radius=.05];
\draw[color=black,fill=black] (-1.025,0.975) circle [radius=.05];
\draw[color=black,fill=black] (1.475,0.475) circle [radius=.05];
\draw[color=black,fill=black] (0.975,2.75) circle [radius=.05];
\draw[color=black,fill=black] (2.25,2.25) circle [radius=.05];
\end{tikzpicture}
\]
\caption{Graph of $\pi=2$ 5 9 3 1 4 8 6 10 12 17 7 11 16 13 15 19 22 20 18 14 21.}
\label{fig:4.8}
\end{figure}
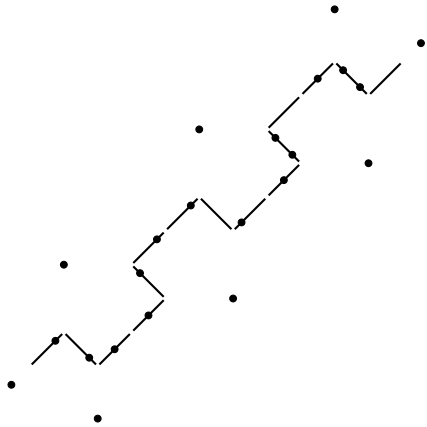

Figure \ref{fig:4.8} shows the graph of a simple permutation $\pi$ of length 21 in $\A$. In particular,\[
\pi=2573146\nwsum_1^0 514263\sesum_1^1 246135 \nwsum_1^0 6152473 \sesum_1^0 2475316.
\]
In order to actually place each point on the crenellation, we usually have to adjust the spacing between points. Thus, the points may not be located at $(i,\pi(i))$ any longer, but this does not change which permutation the graph represents.\\

Note that every permutation of extreme pattern 2413 is in $H$. Indeed, an $N$-shaped structure is a special case of a crenellation, namely one with only 4 isolated points. Thus, simple permutations of extreme pattern 2413 follow the structure described in Figure \ref{fig:4.7}$(a)$ with $m=1$.
\subsection{Proof of Theorem \ref{thm:4.1} (Part 1)}
The proof of Theorem \ref{thm:4.1} is much longer than those of the propositions we discussed earlier, so we break this proof into two propositions.
\begin{proposition}\label{prop:4.4}
If $\pi$ is a simple permutation in $H$, then $\pi$ has one of the structures illustrated in Figure \ref{fig:4.7}.
\end{proposition}
\textit{Proof.} Suppose $\pi$ is a simple permutation in $\A$ with $|\pi|=n\geq 4$ and $\pi(2)\neq 1$. We define a sequence $d_1,\ldots,d_{m+3}$ of values of $\pi$. Let $d_1=\pi(1)$ and\[
d_i=\left\{\begin{array}{cl}
\pi(\max\{s:\pi(s)<\pi(d_{i-1})\}) &\textrm{ if }i\textrm{ is even}\\
\max\{t:\pi^{-1}(t)<\pi^{-1}(d_{i-1})\} &\textrm{ if }i\textrm{ is odd}
\end{array}\right.
\]\\
for $i$ with $1\leq i\leq m+3$. In other words, $d_i$ is the right-most value that is less than $d_{i-1}$ if $i$ is even or the greatest value located to the left of $d_{i-1}$ if $i$ is odd. Note that by definition, $d_2=1$ for any $\pi\in H$. We will soon show that $d_i\neq d_j$ for $i\neq j$. Therefore, since $\pi$ has a finite length, we eventually obtain either $d_{m+2}=\pi(n)$ and $d_{m+3}=n$ for some even integer $m$ or $d_{m+2}=n$ and $d_{m+3}=\pi(n)$ for some odd integer $m$. We let $m$ denote this integer in the remainder of this section.\\

First, we show $d_i\neq d_j$ for $i\neq j$. By the definition of $d_i$, it is clear that $d_i\leq d_{i+2}\leq d_{i+4}\leq\cdots$ for odd $i$ and $\pi^{-1}(d_i) \leq \pi^{-1}(d_{i+2})\leq \pi^{-1}(d_{i+4})\leq\cdots$ for even $i$. We claim that, in addition, $d_i\leq d_{i+2}\leq d_{i+4}\leq\cdots$ for even $i$ and $\pi^{-1}(d_i) \leq \pi^{-1}(d_{i+2})\leq \pi^{-1}(d_{i+4})\leq\cdots$ for odd $i$. Suppose $d_i>d_j$ for positive even integers $i$ and $j$ where $i<j$. Then $d_j<d_{i-1}$, but by definition, $d_i$ is the right-most value that is less than $d_{i-1}$, so we achieve a contradiction. Similarly, if $\pi^{-1}(d_i)>\pi^{-1}(d_j)$ for odd $i,j$ where $i<j$, then $\pi^{-1}(d_j)<\pi^{-1}(d_{i-1})$, which contradicts the fact that $d_i$ is the greatest value located to the left of $d_{i-1}$. Therefore, $d_i\leq d_{i+2}\leq d_{i+4}\leq\cdots$ and $\pi^{-1}(d_i) \leq \pi^{-1}(d_{i+2})\leq \pi^{-1}(d_{i+4})\leq\cdots$ for any positive integer $i$.\\

Next, we show $d_i\neq d_j$ where $i$ is odd and $j$ is even. It is immediate by definition that $d_i\neq d_{i+1}$ for any $i$. For the case $|i-j|>1$, assume for contradiction that we have $d_i=d_j$ and $i<j$. Then by definition, $d_{j-1}$ must be located to the right of $d_j=d_i$, which cannot be true since $\pi^{-1}(d_i)\leq\pi^{-1}(d_j)$. Likewise, if $i>j$, then $d_{i-1}<d_i=d_j$, which is also a contradiction.\\

We now claim that, if $d_i\neq n$ for odd $i$, then $d_i\neq d_{i+2}$. Suppose to the contrary that $d_i=d_{i+2}$. Note that $d_{i+1}\neq\pi(n)$ since $d_{i+1}=\pi(n)$ implies $d_{i+2}=d_i=n$. By definition, $d_{i+1}$ is the right-most value which is less than $d_i$, and $d_i$ is the greatest whose position is in $[1,\pi^{-1}(d_{i+1})]$. These imply that any position corresponding to a value less than $d_i$ is in $[1,\pi^{-1}(d_{i+1})]$ and any value corresponding to a position less than $d_{i+1}$ is in $[i,d_i]$. Hence, $[1,\pi^{-1}(d_{i+1})]$ is a block with $d_{i+1}\neq\pi(n)$, which is a contradiction.\\

To show $d_i\neq d_j$ for two distinct odd $i$ and $j$ where $|i-j|>2$, we assume $d_i=d_j$. Without loss of generality, say $i<j$. With $d_i\leq d_{i+2}\leq d_{i+4}\leq\cdots\leq d_j$, we have $d_i=d_{i+2}$, which cannot be true. Thus, $d_i\neq d_j$ for any two distinct odd integers $i$ and $j$, so long as $d_i\neq n$. By the inverse argument, $d_i\neq d_j$ for even $i$ and $j$ so long as $d_i\neq\pi(n)$.\\

Hence, we have shown that $d_i\neq d_j$ for $i\neq j$.\\

Note that $d_2\neq\pi(n)$ since $d_2=1$ and $\pi$ is simple. Hence, every simple permutation in $H$ has at least 4 values denoted by $d_i$ for some $i\geq 1$. For a given $\pi$, suppose $\pi$ has $m+3$ values denoted by $d_i$ $(1\leq i\leq m+3)$. We show that $\pi$ satisfies one of Equations \ref{eqn:4.1} by induction on $m\geq 1$.\\

For the base case, suppose $m=1$. Then $d_4$ is the last value of $\pi$. Notice that each $d_i$ ($1\leq i\leq 4$) is an extreme point of $\pi$, and since the flattening of $d_1d_3d_2d_4$ is 2413, $\pi$ has extreme pattern 2413, so we are done.\\

Now, suppose that every $\pi\in H$ with $m+3$ values denoted by $d_i$ $(1\leq i\leq m+3)$ satisfies Equation \ref{eqn:4.1}$(a)$ for some positive odd integer $m$. We first need to show that an arbitrary $\pi\in H$ of length $n$ with $m+4$ values denoted by $d_i$ $(1\leq i\leq m+4)$ satisfies Equation \ref{eqn:4.1}$(b)$. Let $\pi$ be a permutation in $H$ of length $n$ with $m+4$ values denoted by $d_i$. In this case, $d_{m+4}=n$. Let $p_m$ be the following.\[
p_m=\pi(\min\{\pi^{-1}(s):s>d_{m+3},\pi^{-1}(s)>\pi^{-1}(d_{m+1})\}).
\]
In other words, $p_m$ is the left-most value greater than $d_{m+3}$ and located to the right of $d_{m+1}$. Since $d_{m+4}>d_{m+3}$ and $d_{m+4}$ is located to the right of $d_{m+1}$, $\pi^{-1}(p_m)\leq \pi^{-1}(d_{m+4})$.\\

Next, we define $q_m$ and $r_m$ by
\[
q_m=\pi(\pi^{-1}(p_m)-1)\qquad\textrm{and}\qquad r_m=\max\{\pi(s)<d_{m+3}:s\in[1,\pi^{-1}(q_m)]\}
\]
\textit{i.e.} $q_m$ is the value immediately to the left of $p_m$, and $r_m$ is greatest value less than $d_{m+3}$ whose position is in the segment $[1,\pi^{-1}(q_m)]$. It is possible that $q_m=d_{m+1}$ or $q_m=r_m$, but not both, since $r_m\geq d_m>d_{m+1}$. Note that $q_m<d_{m+3}$ because, $q_m\geq d_{m+3}$ contradicts the definition of $p_m$ as $q_m$ is located to the left of $p_m$ and greater than $d_{m+3}$.\\
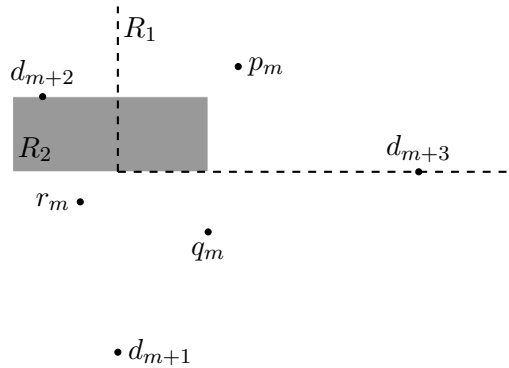
\begin{figure}[b!]
\[\begin{tikzpicture}
\draw[white,fill=gray!80](-2.2,-0.8) rectangle (0.4,0.2);
\draw[color=black,fill=black] (-1.8,0.2)node[above]{$d_{m+2}$} circle [radius=.04];
\draw[color=black,fill=black] (-0.8,-3.2)node[right]{$d_{m+1}$} circle [radius=.04];
\draw[color=black,fill=black] (0.4,-1.6)node[below]{$q_m$} circle [radius=.04];
\draw[color=black,fill=black] (0.8,0.6)node[right]{$p_m$} circle [radius=.04];
\draw[color=black,fill=black] (-1.3,-1.2)node[left]{$r_m$} circle [radius=.04];
\draw[color=black,fill=black] (3.2,-0.8)node[above]{$d_{m+3}$} circle [radius=.04];
\draw[thick,dashed] (-0.8,-0.8) -- (-.8,1.4);
\draw[thick,dashed] (-0.8,-0.8) -- (4.4,-0.8);
\node at (-0.5,1.1) {$R_1$};
\node at (-1.9,-0.5) {$R_2$};
\end{tikzpicture}
\]
\caption{Illustration of relations among $p_m$, $q_m$, $r_m$ and $d_i$.}
\label{fig:4.11}
\end{figure}

We provide Figure \ref{fig:4.11} to show the relations among $p_m$, $q_m$, $r_m$ and $d_i$. The value $p_m$ is the left-most value in the region denoted by $R_1$. It is possible to have $p_m<d_{m+2}$. The value $q_m$ is immediately to the left of $q_m$, and as noted previously, it has to be less than $d_{m+3}$. The position of $r_m$ can also be to the left of $d_{m+2}$.\\
\begin{figure}
\[\begin{array}{ccc}
\begin{tikzpicture}[scale=0.8]
    \draw[white, fill=gray!40](1,4)--(3,4)--(3,5)--(1,5)--cycle;
    \draw[white, fill=gray!80](0.5,3)--(1,3)--(1,5)--(0.5,5)--cycle;
    \draw[white, fill=gray!80](4,0.5)--(5,0.5)--(5,2)--(4,2)--cycle;
    \draw[white, fill=gray!40](3,0.5)--(4,0.5)--(4,2)--(3,2)--cycle;
    \plotPerm{2,4,1,5,3}
    \foreach \x in {1,...,5} \draw[black,thick] (\x,1)--(\x,0.5) (1,\x)--(0.5,\x);
    \node[above left] at (1,2) {$d_m$};
    \node[above] at (2,4) {$d_{m+2}$};
    \node[above] at (3,1) {$d_{m+1}$};
    \node[above] at (4,5) {$d_{m+4}$};
    \node[right] at (5,3) {$d_{m+3}$};
    \node at (1.5,3.5) {$B_{21}$};
    \node at (2.5,3.5) {$B_{22}$};
  \end{tikzpicture} & \quad\quad & \begin{tikzpicture}[scale=0.8]
    \draw[white, fill=gray!40](1,5)--(4,5)--(4,6)--(1,6)--cycle;
    \draw[white, fill=gray!80](0.5,3)--(1,3)--(1,6)--(0.5,6)--cycle;
    \draw[white, fill=gray!80](5,0.5)--(6,0.5)--(6,2)--(5,2)--cycle;
    \draw[white, fill=gray!80](2,0.5)--(3,0.5)--(3,4)--(2,4)--cycle;
    \draw[white, fill=gray!80](3,4)--(6,4)--(6,5)--(3,5)--cycle;
    \draw[white, fill=gray!40](4,0.5)--(5,0.5)--(5,2)--(4,2)--cycle;
    \draw[white, fill=gray!40](1,4)--(3,4)--(3,5)--(1,5)--cycle;
    \plotPerm{2,5,4,1,6,3}
    \foreach \x in {1,...,6} \draw[black,thick] (\x,1)--(\x,0.5) (1,\x)--(0.5,\x);
    \node[above left] at (1,2) {$d_m$};
    \node[above] at (2,5) {$d_{m+2}$};
    \node[above right] at (3,4) {$z$};
    \node[above] at (4,1) {$d_{m+1}$};
    \node[above] at (5,6) {$d_{m+4}$};
    \node[right] at (6,3) {$d_{m+3}$};
  \end{tikzpicture}\\
  (a) & & (b)
  \end{array}
\]
\[\begin{array}{c}
\begin{tikzpicture}[scale=0.8]
    \draw[white, fill=gray!80](0.5,3)--(1,3)--(1,5)--(2,5)--(2,6)--(0.5,6)--cycle;
    \draw[white, fill=gray!80](2,0.5)--(3,0.5)--(3,3)--(2,3)--cycle;
    \draw[white, fill=gray!80](4,4)--(6,4)--(6,5)--(4,5)--cycle;
    \draw[white, fill=gray!80](5,0.5)--(6,0.5)--(6,2)--(5,2)--cycle;
    \draw[white, fill=gray!40](1,4)--(4,4)--(4,6)--(2,6)--(2,5)--(1,5)--cycle;
    \draw[white, fill=gray!40](4,0.5)--(5,0.5)--(5,2)--(4,2)--cycle;
    \plotPerm{2,4,5,1,6,3}
    \foreach \x in {1,...,6} \draw[black,thick] (\x,1)--(\x,0.5) (1,\x)--(0.5,\x);
    \node[above left] at (1,2) {$d_m$};
    \node[above left] at (2,4) {$z$};
    \node[above] at (3,5) {$d_{m+2}$};
    \node[above] at (4,1) {$d_{m+1}$};
    \node[above] at (5,6) {$d_{m+4}$};
    \node[right] at (6,3) {$d_{m+3}$};
    \node at (2.5,3.5) {$B_{32}$};
  \end{tikzpicture}\\
  (c)
  \end{array}
\]
\[\begin{array}{ccc}
\begin{tikzpicture}[scale=0.8]
    \draw[white, fill=gray!80](0.5,3)--(1,3)--(1,6)--(3,6)--(3,7)--(0.5,7)--cycle;
    \draw[white, fill=gray!80](2,0.5)--(4,0.5)--(4,3)--(3,3)--(3,4)--(2,4)--cycle;
    \draw[white, fill=gray!80](3,4)--(7,4)--(7,6)--(5,6)--(5,5)--(3,5)--cycle;
    \draw[white, fill=gray!80](6,0.5)--(7,0.5)--(7,2)--(6,2)--cycle;
    \draw[white, fill=gray!40](1,5)--(5,5)--(5,7)--(3,7)--(3,6)--(1,6)--cycle;
    \draw[white, fill=gray!40](5,0.5)--(6,0.5)--(6,2)--(5,2)--cycle;
    \draw[white, fill=gray!40](3,3)--(4,3)--(4,4)--(3,4)--cycle;
    \plotPerm{2,5,4,6,1,7,3}
    \foreach \x in {1,...,7} \draw[black,thick] (\x,1)--(\x,0.5) (1,\x)--(0.5,\x);
    \node[above left] at (1,2) {$d_m$};
    \node[above left] at (2,5) {$z$};
    \node[below right] at (3,4) {$x$};
    \node[above] at (4,6) {$d_{m+2}$};
    \node[above] at (5,1) {$d_{m+1}$};
    \node[above] at (6,7) {$d_{m+4}$};
    \node[right] at (7,3) {$d_{m+3}$};
    \node at (1.5,4.5) {$B_{31}$};
  \end{tikzpicture} & \quad\quad & \begin{tikzpicture}[scale=0.8]
    \draw[white, fill=gray!80](0.5,3)--(1,3)--(1,6)--(2,6)--(2,7)--(4,7)--(4,8)--(0.5,8)--cycle;
    \draw[white, fill=gray!80](2,0.5)--(5,0.5)--(5,3)--(4,3)--(4,4)--(2,4)--cycle;
    \draw[white, fill=gray!80](4,4)--(8,4)--(8,7)--(6,7)--(6,6)--(4,6)--cycle;
    \draw[white, fill=gray!80](7,0.5)--(8,0.5)--(8,2)--(7,2)--cycle;
    \draw[white, fill=gray!40](1,4)--(2,4)--(2,6)--(1,6)--cycle;
    \draw[white, fill=gray!40](2,6)--(6,6)--(6,8)--(4,8)--(4,7)--(2,7)--cycle;
    \draw[white, fill=gray!40](6,0.5)--(7,0.5)--(7,2)--(6,2)--cycle;
    \draw[white, fill=gray!40](4,3)--(5,3)--(5,4)--(4,4)--cycle;
    \plotPerm{2,5,6,4,7,1,8,3}
    \foreach \x in {1,...,8} \draw[black,thick] (\x,1)--(\x,0.5) (1,\x)--(0.5,\x);
    \node[above left] at (1,2) {$d_m$};
    \node[above left] at (2,5) {$y$};
    \node[above left] at (3,6) {$z$};
    \node[below right] at (4,4) {$x$};
    \node[above] at (5,7) {$d_{m+2}$};
    \node[above] at (6,1) {$d_{m+1}$};
    \node[above] at (7,8) {$d_{m+4}$};
    \node[right] at (8,3) {$d_{m+3}$};
  \end{tikzpicture} \\
  (d) & & (e)
  \end{array}
\]
\caption{Partial graphs of $\pi$ to show that there exists no value in $R_2$.}
\label{fig:4.17}
\end{figure}

We claim that there is no value $z$ with $d_{m+3}<z<d_{m+2}$ whose position is in the segment $[1,\pi^{-1}(q_m)]$. In other words, there is no point in the shaded region denoted by $R_2$ in \ref{fig:4.11}. First, there is no value in the intersected region of $R_1$ and $R_2$, because $p_m$ is the left-most value in the region $R_1$. Suppose we have a value in the region $R_2\setminus R_1$. That is, there exists a value $z$ in $B_{21}$ or $B_{22}$ of the graph in Figure \ref{fig:4.17}$(a)$. If $z$ is in $B_{22}$, then we have the graph shown in Figure \ref{fig:4.17}$(b)$. The segment $[\pi^{-1}(d_{m+2}),\pi^{-1}(x)]$ of Figure \ref{fig:4.17}$(b)$ is an unsplittable block, so we achieve an immediate contradiction.\\

Next, assume $z$ is in the region $B_{21}$ of Figure \ref{fig:4.17}$(a)$. Consequently, we have the graph shown in Figure \ref{fig:4.17}$(c)$ with the block $[\pi^{-1}(z),\pi^{-1}(d_{m+2})]$ that needs to be split. The only way to do this is by assuming the existence of a point in the region $B_{32}$. Choose a point with the least value, say $x$, so we have the graph shown in Figure \ref{fig:4.17}$(d)$. To split the block $[\pi^{-1}(z),\pi^{-1}(x)]$ in Figure \ref{fig:4.17}$(d)$, choose the left-most value $y$ in the region $B_{31}$, so we finally have the graph as in Figure \ref{fig:4.17}$(e)$. The segment $[\pi^{-1}(y),\pi^{-1}(x)]$ is a block that cannot be split, so we have a contradiction. Therefore, there is no value $z$ with $d_{m+3}<z<d_{m+2}$ and $\pi^{-1}(z)\in[1,\pi^{-1}(q_m)]$.\\

Denote by $\pi_2$ the flattening of the subsequence of $\pi$ obtained by removing every value except $d_{m+2}$ and $d_{m+1}$ corresponding to a position in $[1,\pi^{-1}(q_m)]$. We claim that $\pi_2$ is a simple permutation in $\A$ with $|\pi_2|\geq 4$ and $\pi_2(2)=1$. From here, we use the hat notation on a variable, such as $\hat{x}$, to refer to the value of $\pi_2$ corresponding to the value of $\pi$ denoted by $x$. It is obvious that $\pi_2$ is still in $\A$ since $\pi_2\cont\pi$ and $\pi\in\A$. Also, since $\hat{d}_{m+1}$, $\hat{d}_{m+2}$, $\hat{d}_{m+3}$ and $\hat{d}_{m+4}$ are distinct values of $\pi_2$ of $\pi_2$, $|\pi_2|\geq 4$. Notice that $\pi_2(1)=\hat{d}_{m+2}$ and $\pi_2(2)=\hat{d}_{m+1}=1$. Hence, the only thing we need to show is that $\pi_2$ is simple.\\

Suppose $\pi_2$ is not simple. This implies that we have a proper non-singleton segment $I$ of $\pi_2$ which is a block. There are three cases to consider. First, assume $I=[\pi_2^{-1}(\hat{x}),\pi_2^{-1}(\hat{y})]$ for some values $\hat{x}$ and $\hat{y}$ where $\pi_2^{-1}(\hat{x})>\pi_2^{-1}(\hat{d}_{m+1})$. Because every point with a position less than $\pi^{-1}(q_m)$ was removed, we must have $\pi^{-1}(q_m)<\pi^{-1}(x)$, so $[\pi^{-1}(x),\pi^{-1}(y)]$ is also a block in $\pi$. Since $\pi$ is simple, $\pi$ having a block is a contradiction. Next, suppose $I=[\pi_2^{-1}(\hat{d}_{m+1}),\pi_2^{-1}(\hat{y})]$ for some $\hat{y}$. In this case, the value immediately to the right of $\hat{d}_{m+1}$ must be in $\pi(I)$. By construction, this value is $\hat{p}_m$. Because $\hat{d}_{m+1}<\hat{d}_{m+3}<\hat{p}_m$, $\hat{d}_{m+3}\in\pi(I)$, implying $\pi_2^{-1}(\hat{d}_{m+3})\in I$. However, since $\pi_2^{-1}(\hat{d}_{m+1})<\pi_2^{-1}(\hat{d}_{m+4})<\pi_2^{-1}(\hat{d}_{m+3})$, we have $\pi_2^{-1}(\hat{d}_{m+4})\in I$, which implies $\hat{d}_{m+4}\in\pi(I)$. Since $\hat{d}_{m+4}$ is the greatest value in $\pi_2$, $\pi(I)=[\hat{d}_{m+1},\hat{d}_{m+4}]$, so $I$ is all of $\pi_2$, contradiction. Finally, assume $I=[\pi_2^{-1}(\hat{d}_{m+2}),\pi_2^{-1}(\hat{y})]$ for some value $\hat{y}$. Since $\hat{d}_{m+1}\in\pi_2(I)$ and $\hat{d}_{m+3}<\hat{d}_{m+2}$, $\hat{d}_{m+3}$ also has to be in $\pi_2(I)$, so $\pi_2^{-1}(\hat{d}_{m+3})\in I$. However, since $\hat{d}_{m+3}$ is the last value of $\pi_2$, this implies $\hat{y}=\hat{d}_{m+3}$, which means $I$ is all of $\pi_2$. Consequently, we conclude that $\pi_2$ is simple. Moreover, $\pi_2$ has extreme pattern 3142.\\

Next, let $\pi_1$ be the flattening of the subsequence $\pi(1)\pi(2)\cdots \pi(\pi^{-1}(q_m)-1) d_{m+3}$ of $\pi$ if $q_m=r_m$. In other words, if $q_m$ is the second greatest value whose position is in the segment $[1,\pi^{-1}(q_m)]$, then $\pi_1$ is the flattening of the sequence containing all values up to the one immediately to the left of $q_m$, and $d_{m+3}$. Otherwise, let $\pi_1$ be the flattening of the subsequence $\pi(1)\pi(2)\cdots q_m d_{m+3}$ of $\pi$. We again use the hat notation to refer to the value of $\pi_1$ corresponding to the value of $\pi$. We claim that $\pi_1$ is simple.\\

First, we show that every value we removed to construct $\pi_1$ is greater than or equal to $r_m$. Suppose to the contrary that there exists a value $z<r_m$ that was removed. Since $p_m>d_{m+3}$ by definition and $d_{m+3}>r_m$, we have $p_m>r_m$, so $z$ cannot be $p_m$. Also, if $z=q_m$, then $q_m=r_m$ as well, so we have a contradiction. Therefore, we must have $\pi^{-1}(z)>\pi^{-1}(p_m)$. There are four subcases to consider. If $\pi^{-1}(r_m)>\pi^{-1}(d_{m+2})$ and $p_m<d_{m+2}$, then $\pi$ contains 4231 pattern with $d_{m+2}r_m p_m z$. If $\pi^{-1}(r_m)>\pi^{-1}(d_{m+2})$ and $p_m>d_{m+2}$, then $\pi$ contains 42513 pattern with $d_{m+2}r_m p_m zd_{m+3}$. If $\pi^{-1}(r_m)<\pi^{-1}(d_{m+2})$ and $p_m<d_{m+2}$, then $\pi$ contains 35142 pattern with $r_m d_{m+2} d_{m+1} p_m z$, and finally, if $\pi^{-1}(r_m)<\pi^{-1}(d_{m+2})$ and $p_m>d_{m+2}$, then $\pi$ contains 351624 pattern with $r_m d_{m+2}d_{m+1} p_m z d_{m+3}$. Hence, every value we removed to construct $\pi_1$ must be greater than or equal to $r_m$.\\

Now, assume $\pi_1$ is not simple. Let $I$ be a proper non-singleton block of $\pi_1$. Suppose $\pi_1(I)=[\hat{x},\hat{y}]$ for some values $\hat{x}$ and $\hat{y}$ where $\hat{y}<\hat{d}_{m+3}$. Note that $\hat{y}<\hat{d}_{m+3}$ implies $\hat{y}\leq\hat{q}_m-1$ if $q_m=r_m$ and $\hat{y}\leq \hat{q}_m$ if $q_m\neq r_m$. In either case, for all $\hat{z}\in[\hat{x},\hat{y}]$, $\hat{z}=z$ since every point we removed had a value greater than or equal to $r_m$. Hence, the block $I$ in $\pi_1$ is also a block in $\pi$, so we achieve a contradiction. Next, assume $\pi_1(I)=[\hat{x},\hat{d}_{m+3}]$. Since $I$ is a block, $\hat{d}_{m+3}-1$ must be in $[\hat{x},\hat{d}_{m+3}]$. By the way we constructed $\pi_1$, the position of $\hat{d}_{m+3}-1$ must be either to the left of $\hat{d}_{m+2}$ or in between of $\hat{d}_{m+2}$ and $\hat{d}_{m+1}$. In either case, it implies that $\pi^{-1}_1(\hat{d}_{m+1})\in I$, so $\hat{d}_{m}\in\pi_1(I)$ as $\hat{d}_{m+1}<\hat{d}_m<\hat{d}_{m+3}$. Hence, $\pi^{-1}(\hat(d)_m)\in I$, so $\pi^{-1}(\hat(d)_{m+2})\in I$, which means $\hat{d}_{m+2}\in\pi_1(I)$. However, $\hat{d}_{m+2}>\hat{d}_{m+3}$, so we have a contradiction. Finally, suppose $\pi_1(I)=[\hat{x},\hat{d}_{m+2}]$. Because $r_m<d_{m+3}<d_{m+2}$, we have $\hat{d}_{m+3}=\hat{d}_{m+2}-1$. Thus, $\hat{d}_{m+3}\in[\hat{x},\hat{d}_{m+2}]$. Because the position of $\hat{d}_{m+1}$ is in between $\hat{d}_{m+2}$ and $\hat{d}_{m+3}$, we must have $\hat{d}_{m+1}$ in $[\hat{x},\hat{d}_{m+2}]$. Now, since $\hat{d}_{m+1}<\hat{d}_m<\hat{d}_{m+2}$, $\hat{d}_m$ must be also in $[\hat{x},\hat{d}_{m+2}]$, but this again implies $\hat{d}_{m-1}\in[\hat{x},\hat{d}_{m+2}]$, because $\pi_1^{-1}(\hat{d}_m)<\pi_1^{-1}(\hat{d}_{m-1})<\pi_1^{-1}(\hat{d}_{m+1})$. Continuing in this way, we must have $\hat{d}_1=1$ in $[\hat{x},\hat{d}_{m+2}]$, but then $I$ is all of $\pi_1$. Since $I$ is a proper block, this is a contradiction. Consequently, $\pi_1$ is simple.\\

Since $\pi_1$ is a simple permutation of length 4 or more with $\pi(2)\neq 1$, $\pi_1$ is in $H$ with $m+3$ values denoted by $d_i$ ($1\leq i\leq m+3$), $\pi_1$ has the form expressed in Equation \ref{eqn:4.1}$(a)$ by the induction hypothesis. Let $n_1=|\pi_1|$ and $n_2=|\pi_2|$. If $q_m=r_m$, then we apply $\nwsum_1^1$. The greatest value $n_1$ of $\pi_1$ is shifted upward by $\pi_2(1)-2$, which is how much $d_{m+2}$ was shifted down to construct $\pi_1$ by flattening. The value $\pi_1(n_1)$ is $n_1-1$ in $\pi_1$, and it is the second greatest value up to the position of itself in $\pi$. By definition, this is $r_m=q_m$ in $\pi$. Finally, to construct $\pi_2$, we removed $n_1-2$ values from $\pi$, so shifting up each value of $\pi_2$, except $\pi_2(1)$ and $\pi_2(2)$, by $n_1-2$ will recover the values of $\pi$ corresponding to $\pi_2$. Hence, $\pi_1\nwsum_1^1\pi_2=\pi_1(1)\cdots d_{m+2}\cdots q_m p_m\cdots d_{m+3}$. Similarly, if $q_m\neq r_m$, then $\pi=\pi_1(1)\cdots d_{m+2}\cdots q_m p_m\cdots d_{m+3}$. In either case, $\pi\in H$ with $m+4$ values denoted by $d_i$ has the form expressed in Equation \ref{eqn:4.1}$(b)$, so we are done.\\

For our purpose, in the process of combining $\pi_1$ and $\pi_2$ with $\nwsum_1^0$, it is more appropriate to think that $\pi_1(n_1)$ and $\pi_2(1)$ are combined into $d_{m+3}$ and $d_{m+1}$ respectively, as these are the values corresponding to them in $\pi_1$ and $\pi_2$ after the flattening. For $\nwsum_1^1$, even though it appears as if $\pi_1(n_1)$ stays where it is, it is still proper to think that $\pi_1(n_1)$ is merged into the right-most value of $\pi_2$ to become $d_{m+3}$, and leaving a copy of itself at where it used to be as the value $q_m$.\\

Next, assume $m$ is even, that is, $\pi$ in $H$ with $m+3$ values denoted by $d_i$ ($1\leq i\leq m+3$) satisfies Equation \ref{eqn:4.1}$(b)$ for some positive even integer $m$. Then again, we need to show that for $\pi\in H$ of length $n$ with $m+4$ values denoted by $d_i$ $(1\leq i\leq m+4)$ satisfies Equation \ref{eqn:4.1}$(a)$. This time, let $p_m$, $q_m$ and $r_m$ denote the following.
\[\begin{array}{c}
p_m=\min\{s:\pi^{-1}(s)>\pi^{-1}(d_{m+3}),s>d_{m+1}\},\qquad q_i=p_m-1,\vhhh\\
r_m=\pi(\max\{\pi^{-1}(t)<\pi^{-1}(d_{m+3}):t\in[1,q_m]\}).
\end{array}
\]

The rest of the proof is the the same argument applied to the inverses of all permutations involved. At the end, we acquire Equation \ref{eqn:4.1}$(a)$ for $\pi$, and this completes the proof.\hfill$\blacksquare$\\

Because, for any $\pi\in H$, $\pi$ is simple and the values corresponding to the segment $[1,\pi^{-1}(d_3)]$ are increasing, we conclude $d_1=2$.\\

The crenellation can be viewed as a repetition of $N$ and $S$ structures. Each $N$-shape corresponds to a simple permutation $\s_i$ ($i$ odd) of extreme pattern 2413 in Equations \ref{eqn:4.1}, whereas each $S$-shape corresponds to a simple permutation $\tau_i$ ($i$ even) of extreme pattern 3142. Notice that $m$ in the previous proof indicates the total number of $N$ and $S$ structures. With the condition $\pi(2)\neq 1$, the structure always starts with an $N$-shape. Therefore, $\pi$ with $m+3$ values denoted by $d_i$ $(1\leq i\leq m+3)$ has $(m+1)/2$ $N$-shapes and $(m-1)/2$ $S$-shapes if $m$ is odd or $m/2$ $N$-shapes and $m/2$ $S$-shapes if $m$ is even. We call $m$ the number of components of $\pi$.\\
\subsection{Proof of Theorem \ref{thm:4.1} (Part 2)}
We now prove that any simple permutation of the structure described in Figure \ref{fig:4.7} is in $H$. Before we start, we introduce some new terminology. A value of a permutation $\pi(i)$ is called a \textit{left-to-right maximum} if for all $j$ with $1\leq j\leq i$, $\pi(j)\leq\pi(i)$. In other words, $\pi(i)$ is a left-to-right maximum if it is greater than every value on its left. We define a \textit{right-to-left minimum} analogously. For example, with the permutation\begin{center}
$\pi=2$ 5 9 3 1 4 8 6 10 12 17 7 11 16 13 15 19 22 20 18 14 21
\end{center}
which was provided in Figure \ref{fig:4.8}, the subsequence of left-to-right maxima is 2 5 9 10 12 17 19 22 and the subsequence of right-to-left minima is 1 4 6 7 11 13 14 21. We denote by $\LRmax(\pi)$ and $\RLmin(\pi)$ the set of left-to-right maxima values of $\pi$ and the set of right-to-minima values of $\pi$ respectively.\\
\begin{proposition}\label{prop:4.5}
Let $\pi$ be a simple permutation whose structure is described in Figure \ref{fig:4.7}. Then $\pi$ is in $H$.
\end{proposition}
\textit{Proof.} Suppose $\pi$ is a simple permutation of length $n$ whose structure is as described in Figure \ref{fig:4.7} with $m$ components that was explained previously. We show $\pi$ avoids each permutation in the basis $\{4231,35142,42513,351624\}$. First, we break down the crenellation and set notation for the sets of values in different parts of the crenellation.\\
\begin{figure}[b!]
\[
\begin{tikzpicture}[scale=3]
\draw[dotted,thick] (-1.525,-0.525) rectangle (-0.025,-0.025);
\node[above] at (-0.775,-0.025) {$N_i$};
\draw[dotted,thick] (-0.025,-0.025) rectangle (0.475,1.475);
\node[right] at (0.475,0.725) {$S_{i+1}$};
\path[draw,thick] (-1.5,-0.5) -- (-1.05,-0.05);
\node[above left] at (-1.275,-0.275) {$A_i$};
\path[draw,thick] (-1,-0.05) -- (-0.55,-0.5);
\node[below left] at (-0.775,-0.275) {$B_i$};
\path[draw,thick] (-0.5,-0.5) -- (-0.05,-0.05);
\node[above left] at (-0.275,-0.275) {$C_i$};
\draw[dotted,ultra thick] (-1.45,-1.05) -- (-1.65,-1.25);
\path[draw,thick] (0,0) -- (0.45,0.45);
\node[above left] at (0.225,0.225) {$A_i$};
\path[draw,thick] (0.45,0.5) -- (0,0.95);
\node[below left] at (0.225,0.775) {$B_i$};
\path[draw,thick] (0,1) -- (0.45,1.45);
\node[above left] at (0.225,1.225) {$C_i$};
\draw[dotted,ultra thick] (1,1.4) -- (1.2,1.6);
\draw[color=black,fill=black] (-3.025,-1.025)node[above]{$d_{i}$} circle [radius=.0166];
\draw[color=black,fill=black] (-0.525,-1.525)node[above]{$d_{i+1}$} circle [radius=.0166];
\draw[color=black,fill=black] (-1.025,0.975)node[above]{$d_{i+2}$} circle [radius=.0166];
\draw[color=black,fill=black] (1.475,0.475)node[above]{$d_{i+3}$} circle [radius=.0166];
\end{tikzpicture}
\]
\caption{Notations for sets of values in the crenellation.}
\label{fig:4.12}
\end{figure}
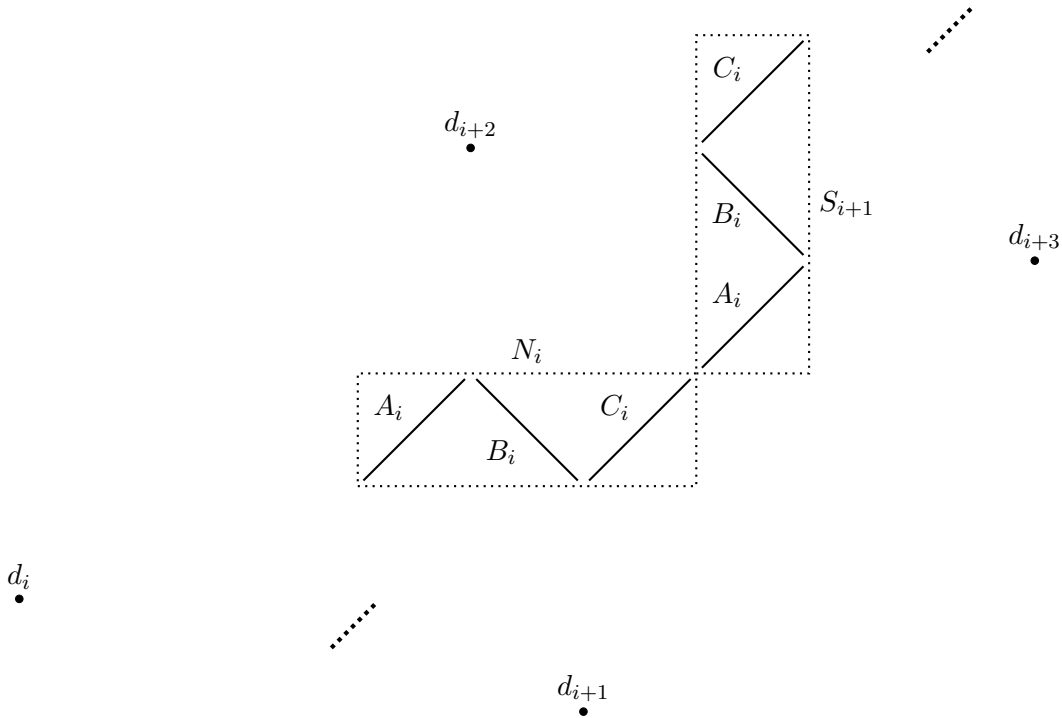

Define each $d_i$ ($1\leq i\leq m+3$) in the same way as in the previous section. Let $N_i$ ($i$ odd) be the set of values corresponding to the permutation $\s_i$ of extreme pattern 2413 in$\pi$. Similarly, let $S_i$ ($i$ even) denote the set of values corresponding to the permutation $\tau_i$ of extreme pattern 31422 in $\pi$. Let $u_i$ be the first value of $N_i$ if $i$ is odd or the least value of $S_i$ if $i$ is even. If $i$ is odd, let $A_i$, $B_i$ and $C_i$ be the set of values corresponding to positions in $[\pi^{-1}(u_i),\pi^{-1}(d_{i+2}))$, the set of values corresponding to positions in $(\pi^{-1}(d_{i+2}),\pi^{-1}(d_{i+1}))$ and the set of values corresponding to positions in $(\pi^{-1}(d_{i+1}),\pi^{-1}(q_i)]$ respectively. If $i$ is even, let $A_i=[u_i,d_{i+2})$, $B_i=(d_{i+2},d_{i+1})$ and $C_i=(d_{i+2},q_i]$ respectively. Figure \ref{fig:4.12} illustrates the definitions of each set of values.\\

Furthermore, let
\[D=\bigcup_{i=1}^{m+3}\{d_i\},\quad A_{\textrm{odd}}=\bigcup_{i\textrm{ odd}} A_i\quad\textrm{ and }\quad A_{\textrm{even}}=\bigcup_{i\textrm{ even}} A_i,
\]
and define $B_{\textrm{odd}}$, $B_{\textrm{even}}$, $C_{\textrm{odd}}$, $C_{\textrm{even}}$ analogously. Hence, for $i$ with $1\leq i\leq m$, $A_i\cup B_i\cup C_i=N_i$ if $i$ is odd, $A_i\cup B_i\cup C_i=S_i$ if $i$ is even and\[
D\cup\left(\bigcup_{i\textrm{ odd}}N_i\right)\cup\left(\bigcup_{i\textrm{ even}}S_i\right)=\{1,\ldots, n\}.
\]

Now we are ready to show that $\pi$ avoids 4231. Note that\[
\LRmax(\pi)=\left(\bigcup_{i\textrm{ odd}}^{m+3}\{d_i\}\right)\cup A_{\textrm{odd}}\cup C_{\textrm{even}}\qquad\textrm{and}\qquad
\RLmin(\pi)=\left(\bigcup_{i\textrm{ even}}^{m+3}\{d_i\}\right)\cup A_{\textrm{even}}\cup C_{\textrm{odd}}
\]
Hence, if any two values from these two sets play the role of 2 or 3, we cannot find a value corresponding to either 4 or 1. This implies that role of 2 and 3 must be played by two values from $B_{\textrm{odd}}\cup B_{\textrm{even}}$. Moreover, since values in each $B_i$ are decreasing, a value playing the role of 2 and a value for 3 must come from distinct $B_i$ and $B_j$ where $i<j$. If $j=i+1$, we are able to assign $d_{i+2}$ to play the role of either 4 or 1 (depending on whether $i$ is even or odd), but then we cannot find another value for the other. Consequently, we are unable to find a subsequence whose flattening is 4231, so $4231\avd\pi$.\\

For 35142, a value corresponding to 4 must come from the set $B$ for the same reason as for the 2 and 3 in 4231. If it is from $B_{\textrm{odd}}$, say $B_i$ for some odd $i$, then a value playing the role of 5 must either also come from $A_i$, from $B_i$ or be $d_{i+2}$, but there is no value that can be assigned to 1 in between. Similarly, if a value for 4 of 35142 is from $B_i$ for some even $i$, a value for 2 must be from $A_i$, from $B_i$ or be $d_{i+2}$, but then we cannot choose a value for 3. Thus, $35142\avd\pi$. We can apply the reverse complement argument to show that $\pi$ avoids 42513 as well.\\

Lastly, suppose $\pi$ contains 351624. Suppose a value for 3 comes from $B_i$ for some odd $i$. Then a value for 1 must be in $B_i\cup C_i\cup\{d_{m+1}\}$. Whichever it is, a value for 2 must be from $C_i$, but this prevents us from assigning a value to 6. Next, assume a value from $B_i$ for some even $i$ plays the role of 3. If $d_{i+2}$ is for 1, then there is no value for 2, so a value for 1 must come from $A_i$ or $B_i$. This forces a value corresponding to 5 to come from $C_i$, but then there is no possible value for 4. Hence, a value for 3 must be from $\LRmax(\pi)$ or $\RLmin(\pi)$. Carrying out similar arguments, we can show that every value comes from $\LRmax(\pi)$ or $\RLmin(\pi)$. Since values in $\RLmin(\pi)$ are right-to-left minima, a value for 3 must be from $\LRmax(\pi)$. So suppose it is from $A_i$ where $i$ is odd. Then a value playing the role of 1 must belong to $B_i$ or it is $d_{i+1}$. Either way, we are able to assign a value from $C_i$ to the role of 2, but then not for 6. If a value for 3 is from $C_i$ where $i$ is even, then a value for 1 cannot be $d_{i+2}$ as this leaves no choice for 2. Hence, a value playing the role of 1 must be from either $A_i$ or $B_i$, forcing a value for 5 to be from $C_i$. However, we now don't have a value for 4, which is a contradiction. Thus, $d_{i}$ for some odd $i$ must play the role of 3. A value for 5 may come from either $C_{i-1}$ or $A_i$, or it is $d_{i+2}$. If it's $d_{i+2}$ or from $A_i$, a value for 1 must be $d_{i+1}$, but as before, we cannot find a value for 2. Thus, a value for 5 must be from $C_{i-1}$, but again, we cannot assign a value for 4, which is another contradiction. Therefore, we achieve a contradiction in every case, implying $\pi$ avoids 351624.\\

Since $\pi$ avoids every permutation in the basis, we have the desired result.\hfill$\blacksquare$\\

Consequently, with Proposition \ref{prop:4.4} and \ref{prop:4.5}, we have Theorem \ref{thm:4.1}.\\
\section{Enumeration}
We are now ready to enumerate the class $\A$. As briefly mentioned in Section 4.1, we first find the generating function for the set $H=\{\pi\in\textrm{Si}(\A):|\pi|\geq 4\textrm{ and }\pi(2)\neq 1\}$. Since the inverses of these permutations give us all the permutations whose second value is 1, doubling the generating function for $H$ gives us the generating function for all simple permutations of length greater than or equal to 4.\\
\subsection{Enumeration of simple permutations in $\A$}
Let us first state the result.\vhhh
\begin{theorem}\label{thm:4.2}
\textit{Let $f_{\textrm{Si}(\A)\setminus\S_2}$ be the generating function for the set of simple permutations in $\A$ excluding $\S_2=\{12,21\}$. Then}\\
\[
f_{\textrm{Si}(\A)\setminus\S_2}=\frac{2x^4}{(1-3x)(1+x)}.
\]
\end{theorem}\vhhh

Once we prove Theorem \ref{thm:4.2}, it is trivial to find $f_{\textrm{Si}(\A)}$, the generating function for all simple permutations in $\A$; simply, we add $2x^2$ to the result to include the permutations 12 and 21.\\
\\
Let $\Sigma=\{a,b,c,d,d_\l\}$. We take the following six steps to accomplish the proof of Theorem \ref{thm:4.2}.\begin{enumerate}
\item Define an encoding function $\phi$ from $H$ to $\Sigma^*$.
\item Define a language $L\sbst\Sigma^*$.
\item Prove $\phi$ is a bijection between $H$ and $L$.
\item Define another language $\overline{L}\sbst\Sigma^*$ which is related to $L$.
\item Define an automaton $M$ such that $\mathcal{L}(M)=\overline{L}$.
\item Apply the transfer matrix method to $M$ to enumerate $|\mathcal{L}(M)|=|\overline{L}|=|L|=|H|$.
\end{enumerate}

We start by defining $\phi$ which maps an arbitrary simple permutation $\pi$ in $H$ to $\Sigma^*$. Let $\pi\in H$ and suppose the number of components of $\pi$ is $m$. By looking at the structure of $\pi$, we can find which of $A_i$, $B_i$, $C_i$ $(1\leq i\leq m)$ or $D$ each value of $\pi$ belongs to. The encoding function $\phi$ reads each value of $\pi$ in certain order and writes out a unique word $w$ consisting of letters in $\Sigma$. Lay out each set of values in the following order.\[
\{d_1,d_2\}\to N_1\to\{d_3\}\to S_2\to\{d_4\}\to N_3\to\{d_5\}\to S_4\to\cdots
\]
This ends as\[
\cdots\to\{d_{m+1}\}\to N_m \to\{d_{m+2},d_{m+3}\}\textrm{ if }m\textrm{ is odd}
\]
or
\[
\cdots\to\{d_{m+1}\}\to S_m \to\{d_{m+2},d_{m+3}\}\textrm{ if }m\textrm{ is even.}
\]

Simply, $\phi$ writes $a$, $b$, $c$ and $d$ for any value in $A_i$, $B_i$, $C_i$ and $D$ respectively, except that it writes $d_\l$ for $d_{m+2}$. For $N_i$ ($1\leq i\leq m$ and $i$ odd), $\phi$ encodes values from bottom to top, \textit{i.e.} $\phi$ first writes the smallest value, then the second smallest value, and so on. On the other hand, $\phi$ encodes values of $S_i$ ($1\leq i\leq m$ and $i$ even) from left to right.\\

The order of encoding for $\{d_1,d_2\}$ and $\{d_{m+2},d_{m+3}\}$ actually does not matter, however, when we establish a similar encoding function for the case of $\A'$ in Chapter \ref{chap:6}, it becomes important to set a convention. For this reason, $\phi$ shall decode $d_2$ first and then $d_1$ next, both into the letter $d$. Similarly, $\phi$ encodes $d_{m+3}$ with $d$ and finally $d_{m+2}$ with $d_\l$.\\

As an example, we encode $\pi$ described in Figure \ref{fig:4.8} using the encoding function $\phi$. First, $\phi$ recognizes $d_2=1$ and $d_1=2$, and encodes both as $d$'s. Moving onto $N_1=\{3,4,5\}$, $\phi$ encodes 3 as $b$, 4 as $c$ and 5 as $a$ in this order, since $3\in B_1$, $4\in C_1$ and $5\in A_1$. Next, $d_3=9$ is encoded as $d$, and $\phi$ continues to $S_2=\{6,8,10\}$. Since $\phi$ encodes values of $S_2$ from left to right, it reads values in the order 10, 6, 8. Again, based on where they are placed, 10, 6 and 8 are encoded as $c$, $a$ and $b$ respectively. Continuing this process up to $d_7=22$, the encoded word $w=\phi(\pi)$ comes out as\[
w=d\,\,d\,\,b\,\,c\,\,a\,\,d\,\,b\,\,a\,\,c\,\,d\,\,c\,\,a\,\,d\,\,b\,\,a\,\,b\,\,d\,\,b\,\,a\,\,b\,\,d\,\,d_\l
\]
\begin{figure}[b!]
\[
\begin{tikzpicture}[scale=1.2]
\draw[arrows=->,thick,black] (-3.65,-3.35) -- (-3.65,-0.975);
\draw[dotted,thick] (-3.525,-2.525) rectangle (-2.025,-2.025);
\node[above] at (-2.775,-2.025) {$N_1$};
\draw[arrows=->,thick,black] (-2.025,-2.15) -- (-0.475,-2.15);
\draw[dotted,thick] (-2.025,-2.025) rectangle (-1.525,-0.525);
\node[right] at (-1.525,-1.275) {$S_2$};
\path[draw,thick] (-3.5,-2.5) -- (-3.05,-2.05);
\draw[arrows=->,thick,black] (-1.65,-0.525) -- (-1.65,1.025);
\draw[dotted,thick] (-1.525,-0.525) rectangle (-0.025,-0.025);
\node[above] at (-0.775,-0.025) {$N_3$};
\draw[arrows=->,thick,black] (-0.025,-0.15) -- (1.525,-0.15);
\draw[dotted,thick] (-0.025,-0.025) rectangle (0.475,1.475);
\node[right] at (0.475,0.725) {$S_4$};
\draw[arrows=->,thick,black] (0.35,1.475) -- (0.35,2.8);
\draw[dotted,thick] (0.475,1.475) rectangle (1.975,1.975);
\node[above] at (1.225,1.975) {$N_5$};
\path[draw,thick] (-3.5,-2.5) -- (-3.05,-2.05);
\draw[color=black,fill=black] (-3.15,-2.15) circle [radius=.05];
\path[draw,thick] (-3,-2.05) -- (-2.55,-2.5);
\draw[color=black,fill=black] (-2.65,-2.4) circle [radius=.05];
\path[draw,thick] (-2.5,-2.5) -- (-2.05,-2.05);
\draw[color=black,fill=black] (-2.275,-2.275) circle [radius=.05];
\path[draw,thick] (-2,-2) -- (-1.55,-1.55);
\draw[color=black,fill=black] (-1.775,-1.775) circle [radius=.05];
\path[draw,thick] (-1.55,-1.5) -- (-2,-1.05);
\draw[color=black,fill=black] (-1.65,-0.65) circle [radius=.05];
\path[draw,thick] (-2,-1) -- (-1.55,-0.55);
\draw[color=black,fill=black] (-1.9,-1.15) circle [radius=.05];
\path[draw,thick] (-1.5,-0.5) -- (-1.05,-0.05);
\draw[color=black,fill=black] (-1.15,-0.15) circle [radius=.05];
\path[draw,thick] (-1,-0.05) -- (-0.55,-0.5);
\path[draw,thick] (-0.5,-0.5) -- (-0.05,-0.05);
\draw[color=black,fill=black] (-0.4,-0.4) circle [radius=.05];
\path[draw,thick] (0,0) -- (0.45,0.45);
\draw[color=black,fill=black] (0.225,0.225) circle [radius=.05];
\path[draw,thick] (0.45,0.5) -- (0,0.95);
\draw[color=black,fill=black] (0.35,0.6) circle [radius=.05];
\draw[color=black,fill=black] (0.1,0.85) circle [radius=.05];
\path[draw,thick] (0,1) -- (0.45,1.45);
\path[draw,thick] (0.5,1.5) -- (0.95,1.95);
\draw[color=black,fill=black] (0.725,1.725) circle [radius=.05];
\path[draw,thick] (1,1.95) -- (1.45,1.5);
\draw[color=black,fill=black] (1.35,1.6) circle [radius=.05];
\draw[color=black,fill=black] (1.1,1.85) circle [radius=.05];
\path[draw,thick] (1.5,1.5) -- (1.95,1.95);
\draw[color=black,fill=black] (-3.8,-2.8)node[left]{$d_1$} circle [radius=.05];
\draw[color=black,fill=black] (-2.525,-3.3)node[right]{$d_2$} circle [radius=.05];
\draw[color=black,fill=black] (-3.025,-1.025)node[left]{$d_3$} circle [radius=.05];
\draw[color=black,fill=black] (-0.525,-1.525)node[right]{$d_4$} circle [radius=.05];
\draw[color=black,fill=black] (-1.025,0.975)node[left]{$d_5$} circle [radius=.05];
\draw[color=black,fill=black] (1.475,0.475)node[right]{$d_6$} circle [radius=.05];
\draw[color=black,fill=black] (0.975,2.75)node[left]{$d_7$} circle [radius=.05];
\draw[color=black,fill=black] (2.25,2.25)node[right]{$d_8$} circle [radius=.05];
\end{tikzpicture}
\]
\caption{Encoding of $\pi=2$ 5 9 3 1 4 8 6 10 12 17 7 11 16 13 15 19 22 20 18 14 21.}
\label{fig:4.13}
\end{figure}
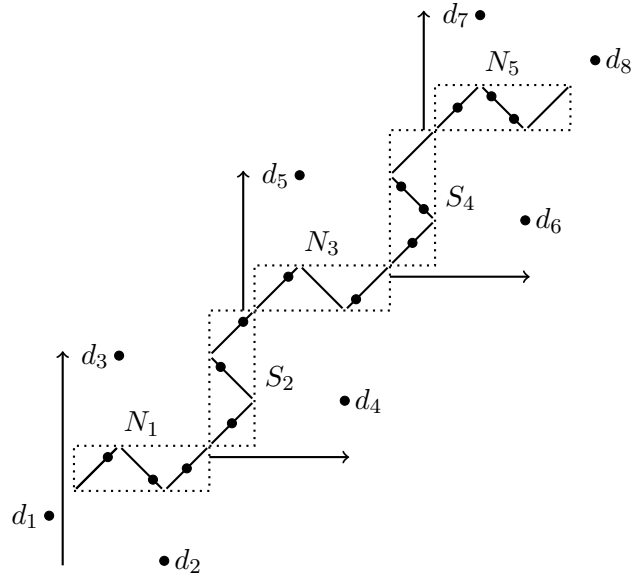
\vspace{-0.1in}

Arrows in Figure \ref{fig:4.13} shows the order of encoding by $\phi$. First arrow reads $d_2$, $d_1$, then values in $N_1$ from bottom to top, and ends with $d_3$. Once we hit intermediate $d$'s, we move to the next arrow. The last one encodes $N_5$, then finally $d_8$ and $d_7$.\\

Next, we define a language $L\sbst\Sigma^*$ with the following conditions for every $w\in L$.\begin{itemize}
\item $w$ must begin with $dd$ and end with $dd_\l$.
\item $w$ must not contain $aa$, $bb$ or $cc$.
\item $d_\l$ is only allowed at the very end.
\item $w$ cannot begin with $dda$ or end with $cdd_\l$.
\item $w$ must not contain $da$.
\end{itemize}
\vhhh

We will show $\phi$ is a bijection between $H$ and $L$, but before we do so, we define a decoding function $\psi$ from $L$ to $\A$. Suppose $w$ is a word in $L$. Let $w_i$ ($1\leq i \leq m$) be the sub-word of $w$ defined as the following.\[
w=\cdots \underbrace{d}_{(i+1)\textrm{-th }d}\underbrace{\cdots\quad\cdots}_{w_i}\underbrace{d}_{(i+2)\textrm{-th }d}\cdots
\]
In particular, $w_i$ is the consecutive letters of $w$ from the letter immediately after the $(i+1)$-th $d$ up to the letter immediately before the $(i+2)$-th $d$. Note that $w_i$ can be empty. Let us divide $w$ as below.\[
w=dd\,\,\ra\,\, w_1\,\,\ra\,\,d\,\,\ra\,\,w_2\,\,\ra\,\,\cdots\,\,\ra\,\,d\,\,\ra\,\, w_m\,\,\ra\,\,dd_\l
\]
The decoding function $\psi$ reads $w$ from left to right. As $\psi$ reads the $k$-th letter, it draws a point at the fixed location $(x_k,y_k)$ $(1\leq k\leq n)$. At the end, $\psi$ constructs the graph of a permutation by drawing points iteratively. We define the decoding function $\psi$ by the following algorithm called \textsf{DECODE}.\\
\\
\textbf{Algorithm} \textsf{DECODE}\\
\textsf{INPUT}: A word $w$ in $L$.\\
\textsf{OUTPUT}: A permutation $\pi$ in $H$.\\
\begin{enumerate}[\hspace{7pt}\sffamily Case 1:]
\item[\sffamily Initialize:] \textsf{Draw first two points and initialize variables.}\\
      Draw points at $(2,1)$ and $(1,2)$. Let $P_a=(2,1)$ and $P_b=P_c=(1,2)$. Let $t=2$. Let $\a$ be the third letter in $w$.
\item \textsf{Draw points for $w_i$ ($i$ odd) which corresponds to the set $N_i$.}\\
      If $\a$ is in $w_i$ ($i$ odd), then \textsf{BEGIN}\begin{enumerate}[\sffamily a.]
          \item If $\a=a$, then draw a point at $(x,y)$ where $P_a^{(x)}<x<P_b^{(x)}$ and $t<y$. Set $P_a$ to be this new point and $t=P_a^{(y)}$. \textsf{GOTO STEP 1} with setting $\a$ to be the next letter.
          \item If $\a=b$, then draw a point at $(x,y)$ where $P_a^{(x)}<x<P_b^{(x)}$ and $t<y$. Set $P_b$ to be this new point and $t=P_b^{(y)}$. \textsf{GOTO STEP 1} with setting $\a$ to be the next letter.
          \item If $\a=c$, then draw a point at $(x,y)$ where $P_c^{(x)}<x$ and $t<y$. Set $P_c$ to be this new point and $t=P_c^{(y)}$. \textsf{GOTO STEP 1} with setting $\a$ to be the next letter.
          \end{enumerate}
      Otherwise, \textsf{GOTO STEP 2}.
\item \textsf{Draw points for $w_i$ ($i$ even) which corresponds to the set $S_i$.}\\
      If $\a$ is in $w_i$ ($i$ even), then \textsf{BEGIN}\begin{enumerate}[\sffamily a.]
          \item If $\a=a$, then draw a point at $(x,y)$ where $t<x$ and $P_a^{(y)}<y<P_b^{(y)}$. Set $P_a$ to be this new point and $t=P_a^{(x)}$. \textsf{GOTO STEP 2} with setting $\a$ to be the next letter.
          \item If $\a=b$, then draw a point at $(x,y)$ where $t<x$ and $P_a^{(y)}<y<P_b^{(y)}$. Set $P_b$ to be this new point and $t=P_b^{(x)}$. \textsf{GOTO STEP 2} with setting $\a$ to be the next letter.
          \item If $\a=c$, then draw a point at $(x,y)$ where $t<x$ and $P_c^{(y)}<y$. Set $P_c$ to be this new point and $t=P_c^{(x)}$. \textsf{GOTO STEP 2} with setting $\a$ to be the next letter.
          \end{enumerate}
      Otherwise, \textsf{GOTO STEP 3}.
\item \textsf{Draw points for $d$'s which correspond to points $d_i$ ($1\leq i\leq m+3$).}\\
      If $\a=d$, then \textsf{BEGIN}\begin{enumerate}[\sffamily a.]
          \item If it is the last $d$ (\textit{i.e.} second letter from the last in $w$), then \textsf{BEGIN}\begin{enumerate}[\sffamily i.]
                \item If it is immediately after $w_i$ with $i$ odd (possibly empty), then draw a point at $(x,y)$ where $P_c^{(x)}<x$ and $t<y$. Set $P_c$ to be this new point and $t=P_c^{(y)}$. \textsf{GOTO STEP 4} with setting $\a$ to be the next letter.
                \item If it is immediately after $w_i$ with $i$ even (possibly empty), then draw a point at $(x,y)$ where $t<x$ and $P_c^{(y)}<y$. Set $P_c$ to be this new point and $t=P_c^{(x)}$. \textsf{GOTO STEP 4} with setting $\a$ to be the next letter.
                \end{enumerate}
          \item Otherwise, \textsf{BEGIN}\begin{enumerate}[\sffamily i.]
                \item If it is immediately after $w_i$ with $i$ odd (possibly empty), then draw a point at $(x,y)$ where $P_a^{(x)}<x<P_b^{(x)}$ and $t<y$. Set $P_a=(P_c^{(x)},t)$, $P_b$ and $P_c$ to be this new point and $t=P_a^{(x)}$. \textsf{GOTO STEP 2} with setting $\a$ to be the next letter.
                \item If it is immediately after $w_i$ with $i$ even (possibly empty), then draw a point at $(x,y)$ where $t<x$ and $P_y^{(x)}<y<P_b^{(y)}$. Set $P_a=(t,P_c^{(y)})$, $P_b$ and $P_c$ to be this new point and $t=P_a^{(y)}$. \textsf{GOTO STEP 1} with setting $\a$ to be the next letter.
                \end{enumerate}
          \end{enumerate}
\item \textsf{Draw a point for $d_\l$ which correspond to points $d_{m+3}$.}\\
      If $\a=d_\l$, then \textsf{BEGIN}\begin{enumerate}[\sffamily a.]
          \item If $m$ is odd (\textit{i.e.} the last sub-word $w_m$ corresponds to $N_m$), then draw a point at $(x,y)$ where $P_a^{(x)}<x<P_b^{(x)}$ and $t<y$. \textsf{GOTO STEP 5}.
          \item If $m$ is even (\textit{i.e.} the last sub-word $w_m$ corresponds to $S_m$), then draw a point at $(x,y)$ where $t<x$ and $P_a^{(y)}<y<P_b^{(y)}$. \textsf{GOTO STEP 5}.
          \end{enumerate}
\item Let $\pi$ be a permutation obtained by flattening the constructed graph. \textsf{OUTPUT} $\pi$.
\end{enumerate}
\vhhh

We visualize how $\psi$ decodes a word $w$ into a permutation. After $\psi$ draws points on $(2,1)$ and $(1,2)$ for the first and second $d$'s, it draws points from bottom to top in intervals $I_1=(1,2)$ and $I_2=(2,\infty)$ for each letter up to the next $d$. For $a$'s and $c$'s, $\psi$ draws points from the left side of $I_1$ and $I_2$ respectively towards right and from the right side of $I_1$ towards left for $b$'s. Additionally, the point for the proceeding $d$ is drawn as if it is for another $b$, \textit{i.e.} immediately left of the last $b$. Figure \ref{fig:4.14}$(a)$ describes how points are drawn until the third $d$.\\

Drawing points for $w_2$ and the fourth $d$ is shown in Figure \ref{fig:4.14}$(b)$. The point with the greatest $x$-coordinate value (\textit{i.e.} the left-most point) in the section denoted by $\overline{N}_1$ in Figure \ref{fig:4.14}$(a)$ determines the $x$-coordinate of the point corresponding to the first letter in $w_2$. Also, the point with the greatest $y$-coordinate in $\overline{N}_1$ and the point corresponding to the third $d$ together determines the next two intervals $I_3$ and $I_4$. Points are drawn from left to right. For $a$'s and $c$'s, $\psi$ draws points upward from the bottom of $I_3$ and $I_4$ respectively, whereas points for $b$'s are drawn downward from the top of $I_3$. Again, $\psi$ draws a point for the fourth $d$ as if it is for $b$. Thus, points for letters in $w_2$ and proceeding $d$ are drawn as shown in Figure \ref{fig:4.14}$(b)$. We continue decoding by repeating this process. The last $d$ and $d_\l$ are drawn as if they are $cb$, both belonging to the last $\overline{N}_i$ or $\overline{S}_i$ section. At the end, we obtain the structure shown in Figure \ref{fig:4.7}. Hence, $m$, the number of sub-words in $w$ is the number of total $N$ and $S$ sets described in the second part of the proof of Theorem \ref{thm:4.1}. Each $d$ in $w$ corresponds to $d_i$ ($1\leq i\leq m+3$) in the proof of Theorem \ref{thm:4.1}, and $a$, $b$ and $c$ are placed to be values in the set $A_i$, $B_i$ and $C_i$ $(1\leq i\leq m$) respectively.\\
\begin{figure}
\[\begin{array}{ccc}\raisebox{-0.57in}{
\begin{tikzpicture}[scale=2]
\draw[dotted,thick] (-1.6,-0.95) rectangle (-0.025,-0.025);
\node[left] at (-1.6,-0.4375) {$\overline{N}_1$};
\path[draw,thick] (-1.5,-0.5) -- (-1.05,-0.05);
\path[draw,thick] (-1,-0.05) -- (-0.55,-0.5);
\path[draw,thick] (-0.5,-0.5) -- (-0.05,-0.05);
\draw[->,>=stealth',semithick] (-1.525,-1.05) -- (-1.03,-1.05);
\node[below] at (-1.2775,-1.1) {$a$};
\draw[->,>=stealth',semithick] (-0.55,-1.05) -- (-1.02,-1.05);
\node[below] at (-0.785,-1.05) {$b$};
\draw[->,>=stealth',semithick] (-0.5,-1.05) -- (-0.03,-1.05);
\node[below] at (-0.265,-1.09) {$c$};
\draw[thick] (-0.525,-1.1) -- (-0.525,-1);
\draw[thick] (-1.55,-1.1) -- (-1.55,-1);
\draw[color=black,fill=black] (-1.55,-0.6)node[right]{$(1,2)$} circle [radius=.0166];
\draw[color=black,fill=black] (-0.525,-0.8)node[right]{$(2,1)$} circle [radius=.0166];
\draw[color=black,fill=black] (-1.025,0.975)circle [radius=.0166];
\draw [decorate,decoration={brace,amplitude=3pt},xshift=0pt,yshift=0pt]
(-0.545,-1.3) -- (-1.53,-1.3) node [below,midway,xshift=10,yshift=-5,text width=1cm]
{$I_1$};
\draw [decorate,decoration={brace,amplitude=3pt},xshift=0pt,yshift=0pt]
(-0.025,-1.3) -- (-0.505,-1.3) node [below,midway,xshift=10,yshift=-5,text width=1cm]
{$I_2$};
\end{tikzpicture}} & \quad\raisebox{.9in}{$\ra$}\quad &\begin{tikzpicture}[scale=2]
\draw[dotted,thick] (-1.6,-0.95) rectangle (0.475,1.475);
\node[left] at (-1.6,0.2625) {$\overline{S}_2$};
\path[draw,thick] (-1.5,-0.5) -- (-1.05,-0.05);
\path[draw,thick] (-1,-0.05) -- (-0.55,-0.5);
\path[draw,thick] (-0.5,-0.5) -- (-0.05,-0.05);
\path[draw,thick] (0,0) -- (0.45,0.45);
\path[draw,thick] (0.45,0.5) -- (0,0.95);
\path[draw,thick] (0,1) -- (0.45,1.45);
\draw[->,>=stealth',semithick] (-0.15,0) -- (-0.15,0.47);
\node[left] at (-0.2,0.235) {$a$};
\draw[->,>=stealth',semithick] (-0.15,0.95) -- (-0.15,0.48);
\node[left] at (-0.2,0.715) {$b$};
\draw[->,>=stealth',semithick] (-0.15,1) -- (-0.15,1.47);
\node[left] at (-0.2,1.235) {$c$};
\draw[thick] (-0.2,-0.025) -- (-0.1,-0.025);
\draw[thick] (-0.2,0.975) -- (-0.1,0.975);
\draw[color=black,fill=black] (-1.55,-0.6) circle [radius=.0166];
\draw[color=black,fill=black] (-0.525,-0.8) circle [radius=.0166];
\draw[color=black,fill=black] (-1.025,0.975) circle [radius=.0166];
\draw[color=black,fill=black] (1.475,0.475) circle [radius=.0166];
\draw [decorate,decoration={brace,amplitude=3pt},xshift=0pt,yshift=0pt]
(-0.4,-0.005) -- (-0.4,0.955) node [below,midway,xshift=-3.5,yshift=7,text width=1cm]
{$I_3$};
\draw [decorate,decoration={brace,amplitude=3pt},xshift=0pt,yshift=0pt]
(-0.4,0.995) -- (-0.4,1.475) node [below,midway,xshift=-3.5,yshift=7,text width=1cm]
{$I_4$};
\end{tikzpicture}\\
(a) & & (b)
\end{array}
\]
\caption{Illustration of the decoding function $\psi$.}
\label{fig:4.14}
\end{figure}
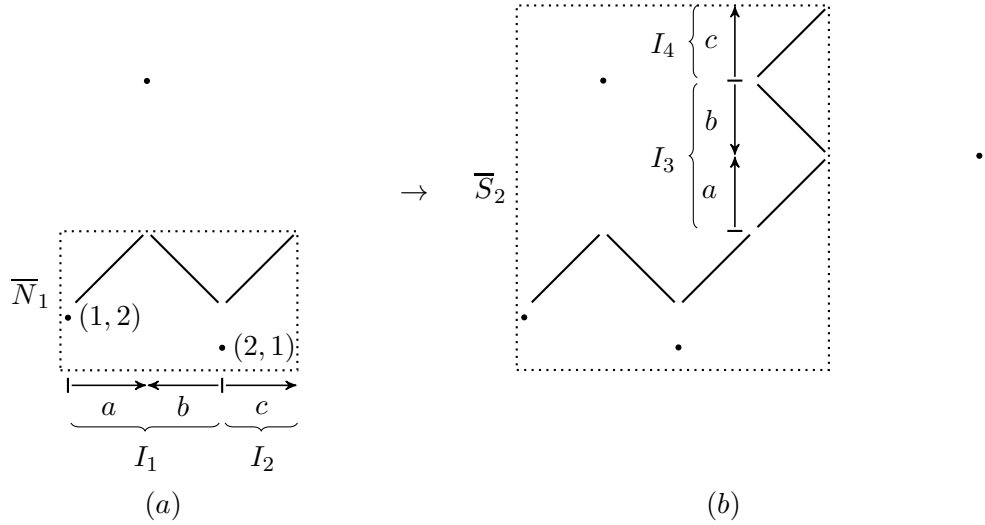
\begin{lemma}\label{lem:4.1}
\textit{The encoding function $\phi$ is a bijection between $H$ and $L$.}
\end{lemma}
\textit{Proof.} We first show the image of $\phi$ is in $L$. Let $\pi$ be in $H$. Our goal is to show that $\phi(\pi)=w$ satisfies the conditions of $L$. It is clear that $w$ begins with $dd$ and end with $dd_\l$, since the first two values and last two values that are encoded by $\phi$ are $d_2$, $d_1$, $d_{m+3}$ and $d_{m+2}$.\\

Now, suppose $w$ contains $aa$. Two values corresponding to these two $a$'s must come from the same $A_i$, because, otherwise, they would be separated by at least one $d$. For $i$ odd, the positions of all values of $A_i$ are consecutive. Therefore, since $\pi$ is simple, the values of two points in $A_i$ cannot be consecutive, as otherwise, those two points would form a block. Hence, as $\phi$ encodes $N_i$ from bottom to top, it cannot write $aa$. On the other hand, for $i$ even, the values of $A_i$ are consecutive, so positions cannot be consecutive. Thus, $\phi$ also cannot write $aa$ while encoding $S_i$. Applying the same argument, we see that $w$ does not contain $bb$ and $cc$ either.\\

If $w$ begins with $dda$, then the first value encoded in $N_1$ belongs to $A_1$. Since $\phi$ encodes $N_1$ from bottom to top, this value is 3. However, this forms a block $[1,2]$ with $\pi([1,2])=[2,3]$. We can show $w$ does not end with $cdd_\l$ in the same way.\\

Finally, suppose $w$ contains $da$. Further, suppose $a$ comes from $A_{i+1}$ for odd $i$ (so $A_{i+1}$ belongs to $S_i$). By the definition of $p_i$ from the proof of Proposition \ref{prop:4.4}, the first value being encoded by $\phi$ is $p_i$ so long as $S_{i+1}$ is nonempty. Thus, the letter $a$ of $da$ corresponds to $p_i$. However, $p_i$ cannot be less than $d_{i+3}$ by definition, so $w$ cannot contain $da$. If $i$ is odd, we use the same argument with the inverse symmetry applied.\\

Next, we need to show that the image of $\psi$ is in $H$. Theorem \ref{thm:4.1} states every simple permutation of the form described in Figure \ref{fig:4.7} is in $\A$. Since the graph of $\pi$ follows one of the structures in \ref{fig:4.7} by construction, all we need to show is for any $w\in L$, $\psi(w)$ is simple. Suppose $\pi=\psi(w)$ is not simple for some $w\in L$. Hence, there exists a proper non-singleton segment $I$ of $\pi$ which is a block.\\

Suppose our non-singleton proper block $I$ contains at least two positions $s$ and $t$ corresponding to $d_i$ and $d_j$ ($1\leq i,j\leq m+3$). Notice that $I$ cannot contain positions corresponding to both $d_1$ and $d_{m+3}$ because, in that case, $I=[1,n]$. Suppose $I$ does not contain $1$, which is the position corresponding to $d_1$. Choose the least $s\in I$ such that $\pi(s)=d_i$ for some $i$ with $1\leq i\leq m+3$. First, assume $i$ is even. If $i=m+1$, then $I$ must contain $t$ ($t=n$ in particular) such that $\pi(t)=d_{i+2}=d_{m+3}$ as well, because $d_{m+3}$ is the only one located to the left of $d_{m+1}$ decoded from a letter $d$ in $w$. However, this implies that $I$ contains a position $u$ such that $\pi(u)=d_{i-1}$, since $d_{i}<d_{i-1}<d_{i+2}$. Since $u<s$, this contradicts our assumption. On the other hand, if $i\leq m$, we know there exists $t\in I$ such that $\pi(t)=d_{i+3}$, but then this implies that $I$ contains $u$ such that $\pi(u)=d_{i-1}$ again, and we achieve a contradiction. Now, suppose $i$ is odd. Then there exists $t$ in $I$ such that $\pi(t)=d_{i-1}$, and this means that there exists a position $u$ in $I$ corresponding to $d_{i-2}$. Since $u<s$, we again have a contradiction. We can apply a similar argument for the case where I does not contain the position corresponding to $d_{m+3}$, so $I$ cannot contain two or more positions corresponding to isolated points.\\

Now, assume $I$ contains only one position $s$ such that $\pi(s)=d_i$ for some $i$ with $1\leq i\leq m+3$. However, this is only possible if $i=1$ or $i=m+3$. Otherwise, if $i$ is odd, then $I$ contains a value $r$ which is either immediately to the left of $d_i$ or immediately to the right of $d_i$. Since the algorithm \textsf{DECODE} forces $r<d_{i+1}<d_i$, $I$ must contain $t$ such that $\pi(t)=d_{i+1}$. We achieve a similar result for the case of $i$ being even. So suppose $I$ contains $s=1$ so that $\pi(1)=d_1$. Since $I$ does not contain a position $t$ such that $\pi(t)=d_3$, $I$ must only contain $1$ and any positions whose corresponding points are placed in $A_1$. However, because $dda$ is not allowed in $w$, this cannot be the case. Similarly, if $I$ contains $s$ such that $\pi(s)=d_{m+3}$, then other positions in $I$ had to be decoded from $c$, but since $cdd_\l$ is not allowed, we again achieve a contradiction.\\

Finally, suppose $I$ contains no points in $D$. This implies that $\pi(I)$ is a subset of $B_i$ or $C_i\cup A_{i+1}$ for some $i$ ($1\leq i\leq m)$. Because $aa$, $bb$ and $cc$ are not in $w$, the only case $I$ can form a block is $I=[s,s+1]$ where the points $(s,\pi(s))$ and $(s+1,\pi(s+1))$ correspond to the last $c$ in $w_i$ and the first $a$ in $w_{i+1}$ respectively. Moreover, they are the last and the first letters in $w_i$ and $w_{i+1}$. However, this would mean $da$ occurs in $w$, which is not allowed, so we achieve a contradiction.\\

Consequently, $\psi(w)=\pi$ cannot have a proper non-singleton segment $I$ forming a block, so $\pi$ is simple. Therefore, the image of $\psi$ is in $H$. Due to how we construct $\phi$ and $\psi$, it is obvious that $\psi(\phi(\pi))=\pi$ and $\phi(\psi(w))=w$ for any $\pi\in H$ and $w\in L$, so $\phi$ is a bijection.\hfill$\blacksquare$\\

Next, we define another language $\overline{L}\sbst\Sigma^*$. We let $\overline{L}$ be the set of words that can be constructed by removing the first two $d$'s of an arbitrary $w$ in $L$. Hence, $\overline{L}=\{w\in\Sigma^*:ddw\in L\}$. There is an obvious bijection from $L$ to $\overline{L}$, namely the one erases the first two $d$'s of $w$ in $L$. The conditions of $\overline{L}\sbst\Sigma^*$ are:\begin{itemize}
\item $w$ must end with $dd_\l$.
\item $w$ must contain no $aa$, $bb$ or $cc$.
\item $d_\l$ is only allowed at the very end.
\item $w$ cannot begin with $a$ or end with $cdd_\l$.
\item $w$ must contain no $da$.
\end{itemize}

Now, we define an automaton $M=(Q,\Sigma,\d,A,\{D_\l\})$ where $Q=\{A,B,C,CD,D,D_\l\}$, $\Sigma=\{a,b,c,d,d_\l\}$, and $\d$ is described in Table \ref{tab:4.2}. Jail states and transitions to them are omitted. Our last task to complete before the enumeration is to show $\mathcal{L}(M)$, the set of words accepted by $M$ is equal to $\overline{L}$.\\
\begin{table}[H]
\[
\begin{array}{c|ccccc}
 & a & b & c & d & d_\l\\
\hline
A & & B & C & D & \\
B & A & & C & D & \\
C & A & B & & CD & \\
CD & & B & C & D & \\
D & & B & C & D & D_\l\\
D_\l & & & & &
\end{array}\]
\caption{Transitions of $M$.}
\label{tab:4.2}
\end{table}
\begin{lemma}\label{lem:4.2}
\textit{$\mathcal{L}(M)=\overline{L}$.}
\end{lemma}
\vhh
\textit{Proof.} We first show $\mathcal{L}(M)\sbst \overline{L}$. Let $w\in\mathcal{L}(M)$. We need to show that $w$ does not violate any condition of $\overline{L}$. From the initial state $A$, transitions that are allowed are $b$, $c$ and $d$. Similarly, in order to reach the only accept state $D_\l$, we have to pass through the state $D$. Since the transition to arrive at $D$ is $d$, we have
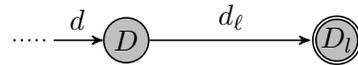
\[\begin{tikzpicture}[shorten >=1pt,auto,node distance=2.8cm,semithick,initial/.style={}]
  \tikzstyle{every state}=[circle,fill=black!25,minimum size=17pt,inner sep=0pt]
  \node[state] (D)               {$D$};
  \node[state,accepting] (Dl) [right of=D] {$D_l$};
  \path[->,>=stealth'] (D) edge node {$d_\l$} (Dl);
  \draw[dotted,thick] (-1.5,0) -- (-1,0);
  \draw[<-,>=stealth'] (D) -- node[above] {$d$} ++(-1cm,0);
\end{tikzpicture}
\]
This implies $w$ must begin with $b$, $c$ or $d$ and end with $dd_\l$. Furthermore, if $w$ ends with $cdd_\l$, then the third transition from the last had to be $c$. This takes to the state $C$, and the next transition $d$ takes to the state $CD$. However, there is no transition to the accept state $D_\l$ from $CD$, so ending with $cdd_\l$ is impossible.\\

Any instance of the letter $a$ in $w$ sends us into state $A$. Since there is no transition using the letter $a$ from $A$, $w$ cannot contain $aa$. By applying similar argument, we can easily show that $w$ does not contain $bb$, $cc$ or $da$. Hence, $w\in\overline{L}$.\\

Next, we show $\overline{L}\sbst\mathcal{L}(M)$. Suppose $w$ is not in $\mathcal{L}(M)$, that is, $w$ is not accepted by the automaton $M$. The only ways $w$ cannot be accepted by $M$ are either the run of $M$ on $w$ contains the jail state or the last state is not $D_\l$. The latter implies that the last letter of $w$ is not $d_\l$, so this violates the first condition of $\overline{L}$, and hence, $w\notin \overline{L}$. For the case the run of $M$ on $w$ contains the jail state, we show that every transition to the jail state is due to a failure of $w$ to meet one of the conditions of $\overline{L}$.\\
\\
\textit{Cases $(A,a)$, $(B,b)$ and $(C,c)$:} To get to states $A$, $B$ and $C$, the previous transitions must be $a$, $b$ and $c$ respectively. Hence, having these transitions implies that $w$ contains $aa$, $bb$ and $cc$ respectively, so $w$ violates the second condition.\\
\\
\textit{Cases $(CD,a)$ and $(D,a)$:} Since the previous transition is $d$ for both cases, having $a$ next means that $w$ contains $da$, so $w$ fails to meet the fourth condition.\\
\\
\textit{Cases $(A,d_\l)$, $(B,d_\l)$, $(C,d_\l)$ and $(CD,d_\l)$:} These transitions imply that $w$ ends with $ad_\l$, $bd_\l$, $cd_\l$ and $cdd_\l$ respectively. Since $w$ must end with $dd_\l$, but not with $cdd_\l$, none of them are allowed.\\
\\
\textit{Cases $(D_\l,a)$, $(D_\l,b)$, $(D_\l,c)$, $(D_\l,d)$ and $(D_\l,d_\l)$:} This causes $d_\l$ to appear in the middle of $w$, so the third condition is not met.\\

Therefore, any run of $M$ containing the jail state implies $w$ that violates at least one condition of $\overline{L}$. We now have proved $\overline{L}\sbst\mathcal{L}(M)$, and this completes the proof of \textit{$\mathcal{L}(M)=\overline{L}$}.\hfill$\blacksquare$\\

Finally, with Lemma \ref{lem:4.1} and \ref{lem:4.2} together, we are ready to prove Theorem \ref{thm:4.2}.\\
\\
\textit{Proof of Theorem \ref{thm:4.2}.} We apply the transfer matrix method to $\mathcal{L}(M)$. With the weight function giving $x$ for all transitions to non-jail states, the adjacency matrix is\[
P=\kbordermatrix{
    & A & B & C & CD & D & D_\l\\
 A  & 0 & x & x & 0  & x &  0  \\
 B  & x & 0 & x & 0  & x &  0  \\
 C  & x & x & 0 & x  & 0 &  0  \\
 CD & 0 & x & x & 0  & x &  0  \\
 D  & 0 & x & x & 0  & x &  x  \\
D_\l& 0 & 0 & 0 & 0  & 0 &  0
}
\]
\\
By observing the $(A,D_\l)$-entry of $(I-P)^{-1}$, we obtain\footnote{The computation of $(I-P)^{-1}$ is done by using Mathematica.}\\
\[
\frac{x^2}{(1-3x)(1+x)}.
\]\\
This is the generating function for $n$-letter words in $\mathcal{L}(M)$. By Lemma \ref{lem:4.2}, this is also the generating function for $n$-letter words in $\overline{L}$. Since the number of $(n-2)$-letter words in $L$ is equal to the number of $n$-letter words in $\overline{L}$, the generating function for $n$-letter words in $L$ is\\
\[
x^2\cdot\frac{x^2}{(1-3x)(1+x)}=\frac{x^4}{(1-3x)(1+x)},
\]\\
which is, by Lemma \ref{lem:4.1}, also the generating function for permutations of length $n$ in the set $H$. As it was explained previously, doubling the generating function gives us\[
f_{\textrm{Si}(\A)\setminus\S_2}=\frac{2x^4}{(1-3x)(1+x)}.
\]
\text{ }\hfill$\blacksquare$\\
\subsection{Enumeration of the whole class of $\A$}
We are close to completing the enumeration of the whole class. From here, we show that every simple permutation in $\A$ satisfies the hypothesis of Proposition \ref{prop:2.4}, then apply both Proposition \ref{prop:2.4} and \ref{prop:2.5}. Let $K=\textrm{Si}(\A)-H-\{12,21\}$. If we show that for all $\pi\in\textrm{Si}(\A)$, $\pi$ satisfies the hypothesis of Proposition \ref{prop:2.4}, we have the equation of generating functions as the following.\\
\[
f_\A=1+x+\sum_{\pi\in\textrm{Si}(\A)} f_{\textrm{ifl}(\pi)}=1+x+f_{\textrm{ifl}(12)}+f_{\textrm{ifl}(21)}+\sum_{\pi\in H} f_{\textrm{ifl}(\pi)}+\sum_{\pi\in K} f_{\textrm{ifl}(\pi)}.
\]

First, using what we have established, we want to acquire the generating function for length $n$ permutations in $\A$ which can be obtained by inflation of simple permutations in $H$. By the inverse argument, we get the exact same result for inflation of simple permutations in $K$, so we will simply multiply 2 as before.\\

Let $\pi\in H$. As it is explained in the second part of the proof of Theorem \ref{thm:4.1}, any value of $\pi$ belongs to either the set $\LRmax(\pi)$, the set $\RLmin(\pi)$ or the set $B:=B_{\textrm{odd}}\cup B_{\textrm{even}}$. In order to apply Proposition \ref{prop:2.4}, we prove the following lemma.\\
\begin{lemma}\label{lem:4.3}
\textit{The condition $\a=\pi[\s_1,\ldots,\s_n]\in\A$ is equivalent to the condition stating that for all $i$ with $1\leq i\leq n$,}\begin{itemize}
\item \textit{if $\pi(i)\in\LRmax(\pi)$, then $\s_i\in\Av(312)$,}
\item \textit{if $\pi(i)\in\RLmin(\pi)$, then $\s_i\in\Av(231)$, and}
\item \textit{if $\pi(i)\in B$, then $\s_i\in\Av(12)$.}
\end{itemize}
\end{lemma}
\vhhh
\textit{Proof.} Suppose the latter condition is false. That is, at least one of the above three conditions is not met. Assume it is the first one. Then there exists $i$ ($1\leq i\leq n$) such that $\pi(i)\in\LRmax(\pi)$ and $\s_i$ contains 312. Now, whichever point $(i,\pi(i))$ is, there is at least one point $(j,\pi(j))$ where $i<j$ and $\pi(i)>\pi(j)$. Namely, this point has value $d_k$ if $\s_i\in C_{k-2}\cup A_{k-1}\cup\{d_{k-1},d_{k+1}\}$. Any point within the block corresponding to $\s_j$ together with the subsequence 312 of $\a$ within the block corresponding to $\s_i$, $\a$ contains 4231, so $\a\notin\A$. If the second condition is not met, we can apply the reverse complement argument of the previous one to show $\a$ contains 4231. If the third condition is false, say $12\cont\s_i$ for some $i$ and the value $\pi(i)$ is in $B_k$ for some $k$, then $d_{k+1}$ and $d_{k+2}$ along with $\s_i$ cause $\a$ to contain 4231, so again, $\a\notin\A$. By the contrapositive argument, $\a=\pi[\s_1,\ldots,\s_n]\in\A$ implies the latter condition.\\

Next, assume a permutation $\a=\pi[\s_1,\ldots,\s_n]$ where $\pi\in H$ is not in $\A$. Thus, $\a$ contains at least one permutation $\b$ in the basis. Since $\pi$ avoids every permutation in the basis, it means there exist $\s_{i_1},\ldots,\s_{i_k}$ (for all $j\in\{1,\ldots,k\}$, $1\leq i_j\leq n$) such that $\a$ contains $\b$ within the union of subintervals corresponding to $\s_{i_1},\ldots,\s_{i_k}$.\\

Since every permutation in $\{35142,42513,351624\}$ is simple, so we only need to consider the case $\b=4231$. If there exist $\s_{i_1},\s_{i_2},\s_{i_3},\s_{i_4}$ such that each point of $\b$ is contained in intervals corresponding to $\s_{i_1},\s_{i_2},\s_{i_3},\s_{i_4}$ respectively, then again, $\b\cont\pi$ which cannot be true. Also, with the same reason as the previous one, if there exists a single $\s_i$ such that $\b\cont\s_i$, then the latter conditions are false. Hence, there exist two or three subintervals of $\a$ such that the containment of $\b$ is involved in. So suppose there exist $\s_{i_1}$ and $\s_{i_2}$ which involve $\b$ together in $\a$. Then it has to be either $312\cont\s_{i_1}$ or $231\cont\s_{i_2}$. Assume $312\cont\s_{i_1}$ is true. Then notice that $\pi(i_1)$ cannot be in $\RLmin(\pi)$, because there is no point which can play the role of 1. This implies $\pi(i_1)\in\LRmax(\pi)$ and $312\cont\s_{i_1}$ or $\pi(i_1)\in B$ and $12\cont 312\cont\s_{i_1}$, so the second condition is not met. With the reverse complement argument, we conclude the same for the case of $231\in\s_{i_2}$.\\

Finally, suppose the containment of $\b$ is shared by three subintervals of $\a$, say the ones of $\s_{i_1}$, $\s_{i_2}$ and $\s_{i_3}$. Then this must imply $12\in\s_{i_2}$. $\pi(i_2)$ can be in neither $\LRmax(\pi)$ nor $\RLmin(\pi)$. Hence, $\pi(i_2)\in B$, so again, the latter condition is false. With every observation we made and contrapositive argument, the second condition implies $\a=\pi[\s_1,\ldots,\s_n]\in\A$, so those two statements are equivalent.\hfill$\blacksquare$\\

Consequently, every simple permutation in $H$ satisfies the hypothesis of Proposition \ref{prop:2.4}. Both $\Av(312)$ and $\Av(231)$ are counted by the Catalan numbers, as explained in Theorem \ref{thm:2.1}. Hence, excluding the empty permutation, we have\\
\[
\bar{f}_{\Av(312)}=\bar{f}_{\Av(231)}=\bar{f}_{\textrm{cat}}=\frac{1-\sqrt{1-4x}}{2x}-1=\frac{1-2x-\sqrt{1-4x}}{2x}.
\]
\\
For $\Av(12)\setminus\{\e\}=\{1,21,321,4321,\ldots\}$, the generating function is the geometric power series minus 1, so\[
\bar{f}_{\Av(12)}=\bar{f}_{\textrm{geom}}=\frac{1}{1-x}-1=\frac{x}{1-x}.
\]

Now, we go back to $\mathcal{L}(M)$. Using the transfer matrix method again, we want to find $f_{\textrm{ifl}(H)}=\sum_{\pi\in H}f_{\textrm{ifl}(\pi)}$ where $\textrm{ifl}(H)=\bigcup_{\pi\in H}\textrm{ifl}(\pi)$. Because each transition embeds the information that which set of $A=A_{\textrm{odd}}\cup A_{\textrm{even}}$, $B$, $C=C_{\textrm{odd}}\cup C_{\textrm{even}}$, or $D$ the point corresponding to the letter belongs to, we can effectively define the weight function so that we can obtain the desired generating function $f_{\textrm{ifl}(H)}$. Define $w(t)=\bar{f}_{\textrm{cat}}$ if $t=a,c,d$ or $d_\l$, and $w(t)=\bar{f}_{\textrm{geom}}$ if $t=b$. The adjacency matrix $\hat{P}$ with this weight function $w$ is\\
\[
\hat{P}=\kbordermatrix{
    & A & B & C & CD & D & D_\l\\
 A  & 0 & \bar{f}_{\textrm{geom}} & \bar{f}_{\textrm{cat}} & 0  & \bar{f}_{\textrm{cat}} &  0  \\
 B  & \bar{f}_{\textrm{cat}} & 0 & \bar{f}_{\textrm{cat}} & 0  & \bar{f}_{\textrm{cat}} &  0  \\
 C  & \bar{f}_{\textrm{cat}} & \bar{f}_{\textrm{geom}} & 0 & \bar{f}_{\textrm{cat}}  & 0 &  0  \\
 CD & 0 & \bar{f}_{\textrm{geom}} & \bar{f}_{\textrm{cat}} & 0  & \bar{f}_{\textrm{cat}} &  0  \\
 D  & 0 & \bar{f}_{\textrm{geom}} & \bar{f}_{\textrm{cat}} & 0  & \bar{f}_{\textrm{cat}} &  \bar{f}_{\textrm{cat}}  \\
D_\l& 0 & 0 & 0 & 0  & 0 &  0
}
\]\\

By computing the $(A,D_\l)$-entry of $(I-\hat{P})^{-1}$ and multiplying by $(\bar{f}_{\textrm{cat}})^2$ to include the initial two $d$'s, we obtain the generating function for $\textrm{ifl}(H)$ by Proposition \ref{prop:2.4} and \ref{prop:2.5}. Since $\textrm{ifl}(K)$ can be obtained by inverting every permutation in $\textrm{ifl}(H)$, we can multiply 2 to the generating function $f_{\textrm{ifl}(H)}$ to include this result. By doing so, we arrive at\\
\[
f_{\textrm{ifl}(\textrm{Si}(\A)\setminus\mathcal{S}_2)}=\frac{\left(-2 x-\sqrt{1-4 x}+1\right)^4}{8 x^2 \left(x^2-\sqrt{1-4
   x} x+3 x+\sqrt{1-4 x}-1\right)}.
\]\\

We now move onto the case where the skeleton is $\pi=21$. Recall that, to ensure $\s_1$ and $\s_2$ are uniquely determined by $\a$ when $\a=21[\s_1,\s_2]$, we must require that $\s_1$ to be skew-indecomposable. We claim that $\a=21[\s_1,\s_2]\in\A$ with skew-indecomposable $\s_1$ is equivalent to the condition that $\s_1\in\Av(312)$ and $\s_1$ is skew-indecomposable, and $\s_2\in\Av(231)$. The condition of $\s_1$ being skew-indecomposable cannot be dropped to enforce the uniqueness of inflation. It is clear that either of $312\cont\s_1$ or $231\cont\s_2$ implies $\a\notin\A$. So assume that $\a\notin\A$. Then $\b\in\a$ for some $\b\in\{4231,35142,42513,351624\}$. As before, since 35142, 42513 and 351624 are simple, if $\b\in\{35142,42513,351624\}$, then $\b\cont\s_1$ or $\b\cont\s_2$. In either case, $312\cont\s_1$ or $231\cont\s_2$, so the second condition is not met. If $\b=4231$, then it is immediate that $312\cont\s_1$ or $231\cont\s_2$, so the condition $\a=21[\s_1,\s_2]\in\A$ with skew-indecomposable $\s_1$ and the condition $\s_1\in\Av(312)$ where $\s_1$ is skew-indecomposable are equivalent.\\

By Proposition \ref{prop:2.4}, $\bar{f}_{\textrm{ifl}(21)}=\bar{f}_{\Av(312)}^\ominus\cdot \bar{f}_{\Av(231)}$. We need to derive $\bar{f}_{\Av(312)}^\ominus$, the generating function for skew-indecomposable permutation in $\Av(312)$. Notice that every skew-decomposable permutation $\pi$ in $\Av(312)$ can be written as $\s\ominus 1$ where $\s$ is a nonempty permutation avoiding 312. To show this, suppose it is not true. Then $\pi=\s\ominus\tau$ where $|\tau|\geq 2$. If $12\cont\tau$, then $\pi$ contains 312. So assume $\tau$ avoids 12. Then $\tau$ must be strictly decreasing, but then, it is possible to write $\pi$ as $\s'\ominus 1$ for some $\s'$, so if $\pi$ is skew-decomposable in $\Av(312)$, it can be always written as $\s\ominus 1$. To find $\bar{f}_{\Av(312)}^\ominus$, we need to exclude skew-decomposable permutations. Since the generating function for skew-decomposable permutations is $x\bar{f}_{\Av(312)}$, we obtain $\bar{f}_{\Av(312)}^\ominus=\bar{f}_{\Av(312)}-x\bar{f}_{\Av(312)}$. Consequently,\[
f_{\textrm{ifl}(21)}=f_{\Av(312)}^\ominus\cdot f_{\Av(231)}=(\bar{f}_{\textrm{cat}}-x\bar{f}_{\textrm{cat}})\cdot\bar{f}_{\textrm{cat}}=\frac{\left(1-2 x-\sqrt{1-4 x}\right)^2 (1-x)}{4 x^2}.
\]\vhhh

Lastly, for the case $\pi=12$, it is possible to inflate both 1 and 2 by any permutations $\s_1$ and $\s_2$ of $\A$ itself, provided that $\s_1$ is a sum-indecomposable permutation in $\A$. Three cases for $\s_1$ being sum-indecomposable are $\s_1=1$, $\s_1$ is skew-decomposable, or $\s_1$ is an inflated permutation of $\pi$ in $H\cup K$ (possibly with $1,\ldots,1$). Generating functions for each case are $x$, $f_{\textrm{ifl}(21)}$ and $f_{\textrm{ifl}(\textrm{Si}(\A)\setminus\mathcal{S}_2)}$ respectively. Thus,
\[
f_{\textrm{ifl}(12)}=(x+f_{\textrm{ifl}(21)}+f_{\textrm{ifl}(\textrm{Si}(\A)\setminus\mathcal{S}_2)})\cdot\bar{f}_\A.
\]

Consequently, the generating function for $\A$ satisfies the functional equation\\
\[\begin{array}{rcl}
f_\A & = & 1+x+f_{\textrm{ifl}(12)}+f_{\textrm{ifl}(21)}+f_{\textrm{ifl}(\textrm{Si}(\A)\setminus\mathcal{S}_2)}\\
 & = & 1+x+(x+f_{\textrm{ifl}(21)}+f_{\textrm{ifl}(\textrm{Si}(\A)\setminus\mathcal{S}_2)})\cdot \bar{f}_\A+f_{\textrm{ifl}(21)}+f_{\textrm{ifl}(\textrm{Si}(\A)\setminus\mathcal{S}_2)}.
 \end{array}
\]
\\
With $\bar{f}_\A=f_\A-1$, we solve for $f_\A$. Then, we obtain\\
\[\begin{array}{rcl}
f_\A & = & \ds\frac{1}{1-x-f_{\textrm{ifl}(21)}-f_{\textrm{ifl}(\textrm{Si}(\A)\setminus\mathcal{S}_2)}}=\frac{2(1-x^2-3x-(1-x)\sqrt{1-4 x})}{1-3x-\sqrt{1-4 x} \left(2 x^2-x+1\right)}\vhhh\\
 & = & 1+x+2x^2+6x^3+23x^4+101x^5+477x^6+2343x^7+11762x^8+\cdots
\end{array}
\]
\newpage
\chapter{Structure of general simple permutations in $\mathcal{A}'$}\label{chap:5}
Recall $\A'=\Av(52341,53241,52431,35142,42513,351624)$. In this chapter, we establish the theorem analogous to Theorem \ref{thm:4.1}. Note that the class $\A'$ is also preserved by the inverse operation and the reverse complement operation.
\section{Extreme patterns 2413, 3142 and 3412}
\subsection{Structural propositions}
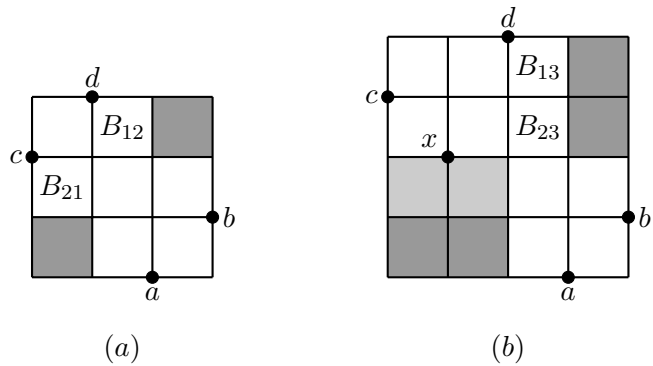
\begin{figure}[H]
\[
\begin{array}{ccc}
  \begin{tikzpicture}[scale=0.8]
    \draw[white, fill=gray!80](1,1)--(2,1)--(2,2)--(1,2)--cycle;
    \draw[white, fill=gray!80](3,3)--(3,4)--(4,4)--(4,3)--cycle;
    \plotPerm{3,4,1,2}
    \node[left] at (1,3) {$c$};
    \node[above] at (2,4) {$d$};
    \node[below] at (3,1) {$a$};
    \node[right] at (4,2) {$b$};
    \node at (1.5,2.5) {$B_{21}$};
    \node at (2.5,3.5) {$B_{12}$};
  \end{tikzpicture} & \qquad\qquad & \begin{tikzpicture}[scale=0.8]
    \draw[white, fill=gray!80](1,1)--(3,1)--(3,2)--(1,2)--cycle;
    \draw[white, fill=gray!80](4,3)--(5,3)--(5,5)--(4,5)--cycle;
    \draw[white, fill=gray!40](1,2)--(1,3)--(3,3)--(3,2)--cycle;
    \plotPerm{4,3,5,1,2}
    \node[left] at (1,4) {$c$};
    \node[above left] at (2,3) {$x$};
    \node[above] at (3,5) {$d$};
    \node[below] at (4,1) {$a$};
    \node[right] at (5,2) {$b$};
    \node at (3.5,4.5) {$B_{13}$};
    \node at (3.5,3.5) {$B_{23}$};
  \end{tikzpicture}\\
  (a) & & (b)
  \end{array}
\]
\caption{Partial graphs of $\pi$ of extreme pattern 3412.}
\label{fig:5.1}
\end{figure}

We first prove the following.
\begin{proposition}\label{prop:5.1}
\textit{No simple permutation in $\A'$ has extreme pattern 3412.}
\end{proposition}
\textit{Proof.} Suppose the statement is false. Let $\pi$ be a simple permutation in $\A'$ of extreme pattern 3412. We start with the graph of extreme pattern of $\pi$, which is shown in Figure \ref{fig:5.1}$(a)$. As before, dark grey indicates a point in the region would create a forbidden pattern, and light grey indicates we have made specific assumptions that there does not exist a point in the region.\\

There must exist a point in $B_{21}$ or $B_{12}$ to avoid $[\pi^{-1}(c),\pi^{-1}(d)]$ being a block. Just like the proof for Proposition \ref{prop:4.3}, one of these two cases can be obtained by the inverse operation followed by the reverse complement operation, so we only give a proof for the case of $B_{21}$. Letting $x$ be the least value of all possible points in $B_{21}$, we obtain the graph shown in Figure \ref{fig:5.1}$(b)$.\\

In order to split the block $[\pi^{-1}(c),\pi^{-1}(d)]$, we need to have a point in either $B_{13}$ or $B_{23}$. Suppose there is a point with some value $y$ in $B_{23}$. Then we obtain the graph in Figure \ref{fig:5.2}$(a)$, forcing it to have a point in $B_{14}$ due to the block $[\pi^{-1}(c),\pi^{-1}(y)]$. However, having a point with a value $z$ in $B_{14}$ results in the graph in Figure \ref{fig:5.2}$(b)$, which contains the unsplittable block $[\pi^{-1}(c),\pi^{-1}(z)]$.\\

\begin{figure}
\[
\begin{array}{ccc}
  \begin{tikzpicture}[scale=0.8]
    \draw[white, fill=gray!80](1,1)--(4,1)--(4,3)--(6,3)--(6,6)--(5,6)--(5,5)--(2,5)--(2,4)--(1,4)--cycle;
    \plotPerm{5,3,6,4,1,2}
    \node[left] at (1,5) {$c$};
    \node[above left] at (2,3) {$x$};
    \node[above] at (3,6) {$d$};
    \node[above left] at (4,4) {$y$};
    \node[below] at (5,1) {$a$};
    \node[right] at (6,2) {$b$};
    \node at (4.5,5.5) {$B_{14}$};
  \end{tikzpicture} & \qquad\qquad & \begin{tikzpicture}[scale=0.8]
    \draw[white, fill=gray!80](1,1)--(5,1)--(5,3)--(7,3)--(7,7)--(4,7)--(4,6)--(3,6)--(3,5)--(2,5)--(2,4)--(1,4)--cycle;
    \draw[white, fill=gray!80](1,6)--(2,6)--(2,7)--(1,7)--cycle;
    \plotPerm{5,3,7,4,6,1,2}
    \node[left] at (1,5) {$c$};
    \node[above left] at (2,3) {$x$};
    \node[above] at (3,7) {$d$};
    \node[above left] at (4,4) {$y$};
    \node[above left] at (5,6) {$z$};
    \node[below] at (6,1) {$a$};
    \node[right] at (7,2) {$b$};
  \end{tikzpicture}\\
  (a) & & (b)
  \end{array}
\]
\caption{Partial graphs of $\pi$ with the assumption of having a value in $B_{23}$ in Figure \ref{fig:5.1}.}
\label{fig:5.2}
\end{figure}
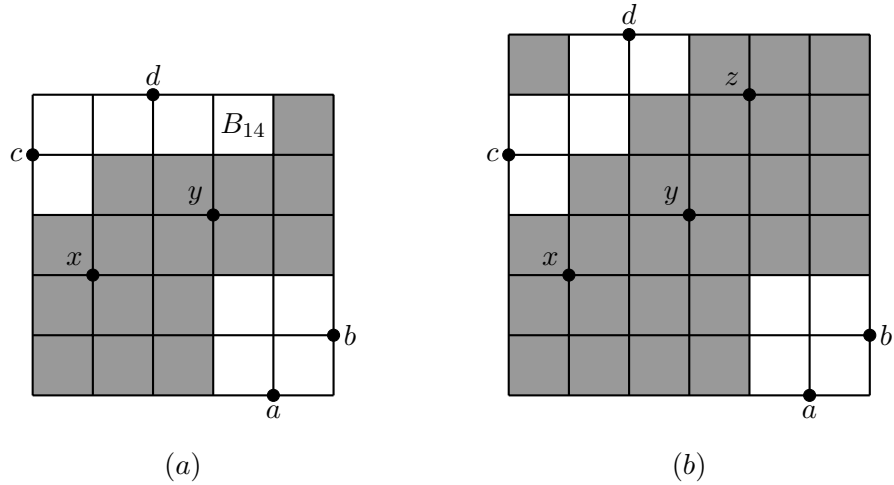
\begin{figure}[b!]
\[
\begin{array}{ccc}
  \begin{tikzpicture}[scale=0.8]
    \draw[white, fill=gray!80](1,1)--(4,1)--(4,2)--(1,2)--cycle;
    \draw[white, fill=gray!80](5,3)--(6,3)--(6,6)--(5,6)--cycle;
    \draw[white, fill=gray!40](1,2)--(3,2)--(3,3)--(1,3)--cycle;
    \draw[white, fill=gray!40](3,3)--(5,3)--(5,6)--(4,6)--(4,4)--(3,4)--cycle;
    \plotPerm{4,3,6,5,1,2}
    \node[left] at (1,4) {$c$};
    \node[above left] at (2,3) {$x$};
    \node[above] at (3,6) {$d$};
    \node[above left] at (4,5) {$y$};
    \node[below] at (5,1) {$a$};
    \node[right] at (6,2) {$b$};
    \node at (3.5,3.5) {$B_{33}$};
    \node at (4.5,3.5) {$B_{34}$};
    \node at (3.5,2.5) {$B_{43}$};
  \end{tikzpicture} & \qquad\qquad & \begin{tikzpicture}[scale=0.8]
    \draw[white, fill=gray!80](1,1)--(5,1)--(5,3)--(7,3)--(7,7)--(4,7)--(4,6)--(3,6)--(3,3)--(1,3)--cycle;
    \draw[white, fill=gray!80](1,6)--(2,6)--(2,7)--(1,7)--cycle;
    \draw[white, fill=gray!40](1,3)--(2,3)--(2,2)--(3,2)--(3,4)--(1,4)--cycle;
    \plotPerm{5,4,7,3,6,1,2}
    \node[left] at (1,5) {$c$};
    \node[above left] at (2,4) {$x$};
    \node[above] at (3,7) {$d$};
    \node[above left] at (4,3) {$z$};
    \node[above left] at (5,6) {$y$};
    \node[below] at (6,1) {$a$};
    \node[right] at (7,2) {$b$};
  \end{tikzpicture}\\
  (a) & & (b)
  \end{array}
\]
\caption{Partial graphs of $\pi$ with the assumption of having a value in $B_{13}$ in Figure \ref{fig:5.1}.}
\label{fig:5.3}
\end{figure}
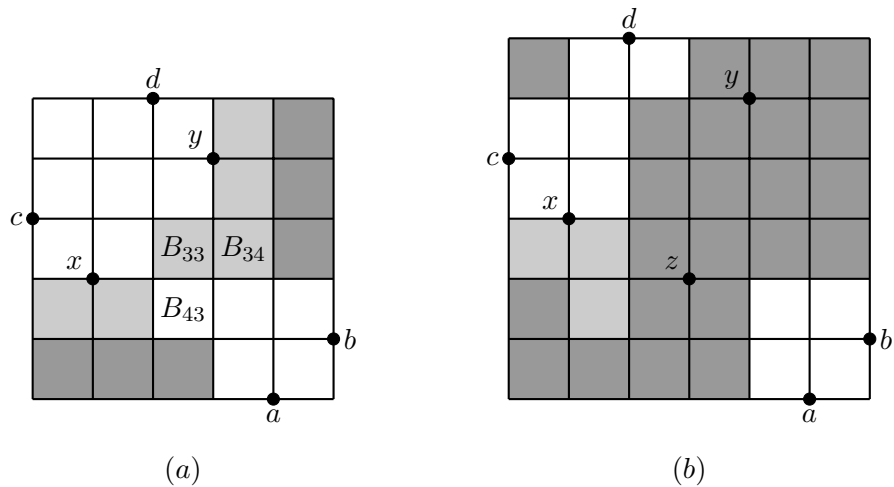

Now, assume there exists a point in $B_{13}$ of Figure \ref{fig:5.1}$(b)$, and none in $B_{23}$. Calling the value of left-most point $y$, we have the graph shown in Figure \ref{fig:5.3}$(a)$. The regions indicated by $B_{33}$ and $B_{34}$ are shaded in light grey since having a point in either of these regions is equivalent to having a point in $B_{23}$ of Figure \ref{fig:5.1}$(b)$. We must have a point in $B_{43}$ of Figure \ref{fig:5.3}$(a)$ to split the block $[\pi^{-1}(c),\pi^{-1}(y)]$, but this results in the block $[\pi^{-1}(c),\pi^{-1}(y)]$ shown in Figure \ref{fig:5.3}$(b)$ which we cannot split any longer.\\

Consequently, it is impossible to have a simple permutation in $\A'$ of extreme pattern 3412.\\
\text{}\hfill$\blacksquare$\\

Next, we describe the structure of simple permutations in $\A'$ of extreme pattern 2413. Let $b$, $d$, $a$ and $c$ be the first, the greatest, the least and the last values of $\pi$ respectively. We denote by $A$, $B$ and $C$ each segment $[\pi^{-1}(b),\pi^{-1}(d))$, $[\pi^{-1}(d),\pi^{-1}(a)]$ and $(\pi^{-1}(a),\pi^{-1}(c)]$ respectively. First, we prove the following lemma for values corresponding to positions in $B$.\\
\begin{lemma}\label{lem:5.1}
\textit{Let $\pi$ be a simple permutation in $\A'$ of extreme pattern 2413. Then values corresponding to positions in $B$ have a pattern of the form}\\
\begin{equation}\label{eqn:5.1}
\bigominus_{i=1}^k\s_i,
\end{equation}
\textit{where each $\s_i$ is either 1 or 12.}
\end{lemma}

For Equation \ref{eqn:5.1}, Note $k=0$ is possible, in which case $B$ consists of only $d$ and $a$.\\

\begin{figure}[H]
\[\begin{tikzpicture}[scale=.5]
\draw[color=black,fill=black] (1,10) circle [radius=.08];
\draw[color=black,fill=black] (2,8) circle [radius=.08];
\draw[color=black,fill=black] (3,9) circle [radius=.08];
\draw[color=black,fill=black] (4,7) circle [radius=.08];
\draw[color=black,fill=black] (5,5) circle [radius=.08];
\draw[color=black,fill=black] (6,6) circle [radius=.08];
\draw [dotted,ultra thick] (7,4) -- (8,3);
\draw[color=black,fill=black] (9,2) circle [radius=.08];
\draw[color=black,fill=black] (10,0) circle [radius=.08];
\draw[color=black,fill=black] (11,1) circle [radius=.08];
\draw[color=black,fill=black] (12,-1) circle [radius=.08];
\node[left] at (1,10) {$d$};
\node[right] at (12,-1) {$a$};
\end{tikzpicture}
\]
\caption{Structure of values corresponding positions in $B$.}
\label{fig:5.4}
\end{figure}
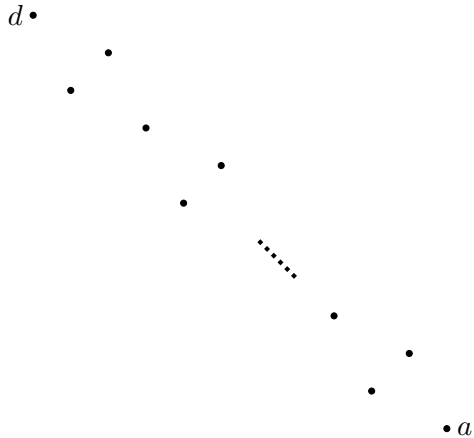

The proof is an immediate consequence of 52341, 53241 and 52431 avoidance conditions with the observation we made with Lemma \ref{lem:3.1}. Therefore, values corresponding to positions in $B$ form the structure described in Figure \ref{fig:5.4}. Except for $d$ and $a$, every point corresponding to 1 and pair of points corresponding to 12 in Equation \ref{eqn:5.1} can be empty.\\

Next, we establish the lemma describing the structure of points whose positions are in the segment $A$. In order to do this, we need to define a special sum of two permutations. Given $\s\in\S_m$ and $\tau\in\S_n$, the \textit{value-interchange sum of $\s$ and $\tau$ with $1$ shift} is the permutation defined by\[
\s\oplus_1\tau=\s'(1)\s'(2)\cdots\s'(m)\tau'(1)\tau'(2)\cdots\tau'(n)
\]
where\[
\s'(i)=\left\{\begin{array}{cl}
\s(i) & \textrm{ if }\s(i)\leq m-1\\
m+1 & \textrm{ if }\s(i)=m
\end{array}\right.\qquad\textrm{for each }i\,\,(1\leq i\leq m)
\]
and\[\tau'(j)=\left\{\begin{array}{cl}
m & \textrm{ if }\tau(j)=1\\
\tau(j)+(m+1) & \textrm{ if }\tau(j)\geq 2
\end{array}\right.\qquad\textrm{for each }j\,\,(1\leq j\leq n)
\]\\
In other words, $\s\oplus_1\tau$ is just like $\s\oplus\tau$ except we interchange the positions of the values corresponding the greatest value of $\s$ and the least value of $\tau$. So, algebraically, $\s\oplus_1\tau=(m\,\,m+1)(\s\oplus\tau)$, the product of $(\s\oplus\tau)$ with the adjacent transposition $(m\,\,m+1)$. For example, with $\s=1342$ and $\tau=312$, $\s\oplus_1\tau=1352746$ whereas $\s\oplus\tau=1342756$. Note that $\oplus_1$ is not associative if one of the summands is length $1$.\\
\begin{figure}[b!]
\[\begin{tikzpicture}[scale=.5]
\draw[color=black,fill=black] (1,1) circle [radius=.08];
\draw[color=black,fill=black] (2,2) circle [radius=.08];
\draw[color=black,fill=black] (3,3) circle [radius=.08];
\draw[color=black,fill=black] (4,4) circle [radius=.08];
\draw[color=black,fill=black] (2.5,-0.5) circle [radius=.08];
\draw[color=black,fill=black] (4.5,1.5) circle [radius=.08];
\draw [dotted,ultra thick] (5.5,4) -- (6.5,5);
\draw[color=black,fill=black] (7,7) circle [radius=.08];
\draw[color=black,fill=black] (8,8) circle [radius=.08];
\draw[color=black,fill=black] (8.5,5.5) circle [radius=.08];
\draw[black,thick] (0.8,0.8)--(1.2,0.8) (2.8,2.8)--(3.2,2.8) (6.8,6.8)--(7.2,6.8);
\end{tikzpicture}
\]
\caption{Structure of a 231-value chain.}
\label{fig:5.5}
\end{figure}
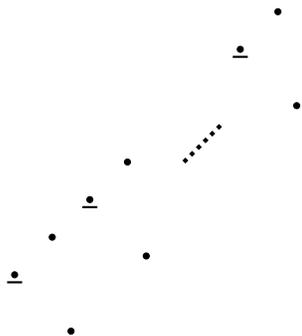

\begin{figure}[b!]
\[\begin{array}{ccc}
  \begin{tikzpicture}[scale=0.8]
    \draw[white, fill=gray!80](1,1)--(2,1)--(2,2)--(1,2)--cycle;
    \draw[white, fill=gray!80](3,3)--(5,3)--(5,2)--(6,2)--(6,4)--(3,4)--cycle;
    \draw[white, fill=gray!80](1,5)--(3,5)--(3,6)--(1,6)--cycle;
    \draw[white, fill=gray!40](2,1)--(4,1)--(4,2)--(3,2)--(3,3)--(1,3)--(1,2)--(2,2)--cycle;
    \plotPerm{4,3,2,6,1,5}
    \node[left] at (1,4) {$b$};
    \node[above left] at (2,3) {$y_0$};
    \node[above left] at (3,2) {$z$};
    \node[above] at (4,6) {$d$};
    \node[below] at (5,1) {$a$};
    \node[right] at (6,5) {$c$};
    \node at (1.5,4.5) {$B_{21}$};
    \node at (2.5,3.5) {$B_{32}$};
  \end{tikzpicture} & \quad\quad & \begin{tikzpicture}[scale=0.8]
    \draw[white, fill=gray!80](1,1)--(3,1)--(3,3)--(6,3)--(6,2)--(7,2)--(7,5)--(2,5)--(2,4)--(1,4)--cycle;
    \draw[white, fill=gray!80](1,6)--(4,6)--(4,7)--(1,7)--cycle;
    \draw[white, fill=gray!40](3,1)--(5,1)--(5,2)--(4,2)--(4,3)--(3,3)--cycle;
    \plotPerm{5,3,4,2,7,1,6}
    \node[left] at (1,5) {$b$};
    \node[above left] at (2,3) {$y_0$};
    \node[above left] at (3,4) {$y_1$};
    \node[above left] at (4,2) {$z$};
    \node[above] at (5,7) {$d$};
    \node[below] at (6,1) {$a$};
    \node[right] at (7,6) {$c$};
    \node at (2.5,5.5) {$B_{22}$};
  \end{tikzpicture}\\
(a) & & (b)\\
\\
\begin{tikzpicture}[scale=0.8]
    \draw[white, fill=gray!80](1,1)--(4,1)--(4,3)--(7,3)--(7,2)--(8,2)--(8,6)--(5,6)--(5,5)--(2,5)--(2,4)--(1,4)--cycle;
    \draw[white, fill=gray!80](1,6)--(2,6)--(2,7)--(5,7)--(5,8)--(1,8)--cycle;
    \draw[white, fill=gray!40](2,6)--(4,6)--(4,7)--(2,7)--cycle;
    \draw[white, fill=gray!40](4,1)--(6,1)--(6,2)--(5,2)--(5,3)--(4,3)--cycle;
    \plotPerm{5,3,6,4,2,8,1,7}
    \node[left] at (1,5) {$b$};
    \node[above left] at (2,3) {$y_0$};
    \node[above left] at (3,6) {$x_1$};
    \node[above left] at (4,4) {$y_1$};
    \node[above left] at (5,2) {$z$};
    \node[above] at (6,8) {$d$};
    \node[below] at (7,1) {$a$};
    \node[right] at (8,7) {$c$};
    \node at (4.5,5.5) {$B_{34}$};
  \end{tikzpicture} & \quad\quad & \begin{tikzpicture}[scale=0.8]
    \draw[white, fill=gray!80](1,1)--(5,1)--(5,3)--(8,3)--(8,2)--(9,2)--(9,7)--(4,7)--(4,6)--(3,6)--(3,5)--(2,5)--(2,4)--(1,4)--cycle;
    \draw[white, fill=gray!80](1,6)--(2,6)--(2,8)--(6,8)--(6,9)--(1,9)--cycle;
    \draw[white, fill=gray!40](2,7)--(4,7)--(4,8)--(2,8)--cycle;
    \draw[white, fill=gray!40](5,1)--(7,1)--(7,2)--(6,2)--(6,3)--(5,3)--cycle;
    \plotPerm{5,3,7,4,6,2,9,1,8}
    \node[left] at (1,5) {$b$};
    \node[above left] at (2,3) {$y_0$};
    \node[above left] at (3,7) {$x_1$};
    \node[above left] at (4,4) {$y_1$};
    \node[above left] at (5,6) {$y_2$};
    \node[above left] at (6,2) {$z$};
    \node[above] at (7,9) {$d$};
    \node[below] at (8,1) {$a$};
    \node[right] at (9,8) {$c$};
    \node at (4.5,7.5) {$B_{24}$};
  \end{tikzpicture}\\
(c) & & (d)
\end{array}\]
\caption{Partial graphs of $\pi$ with the assumption of $b$ corresponding to the 3 in 321 and having a value in $B_{32}$.}
\label{fig:5.6}
\end{figure}
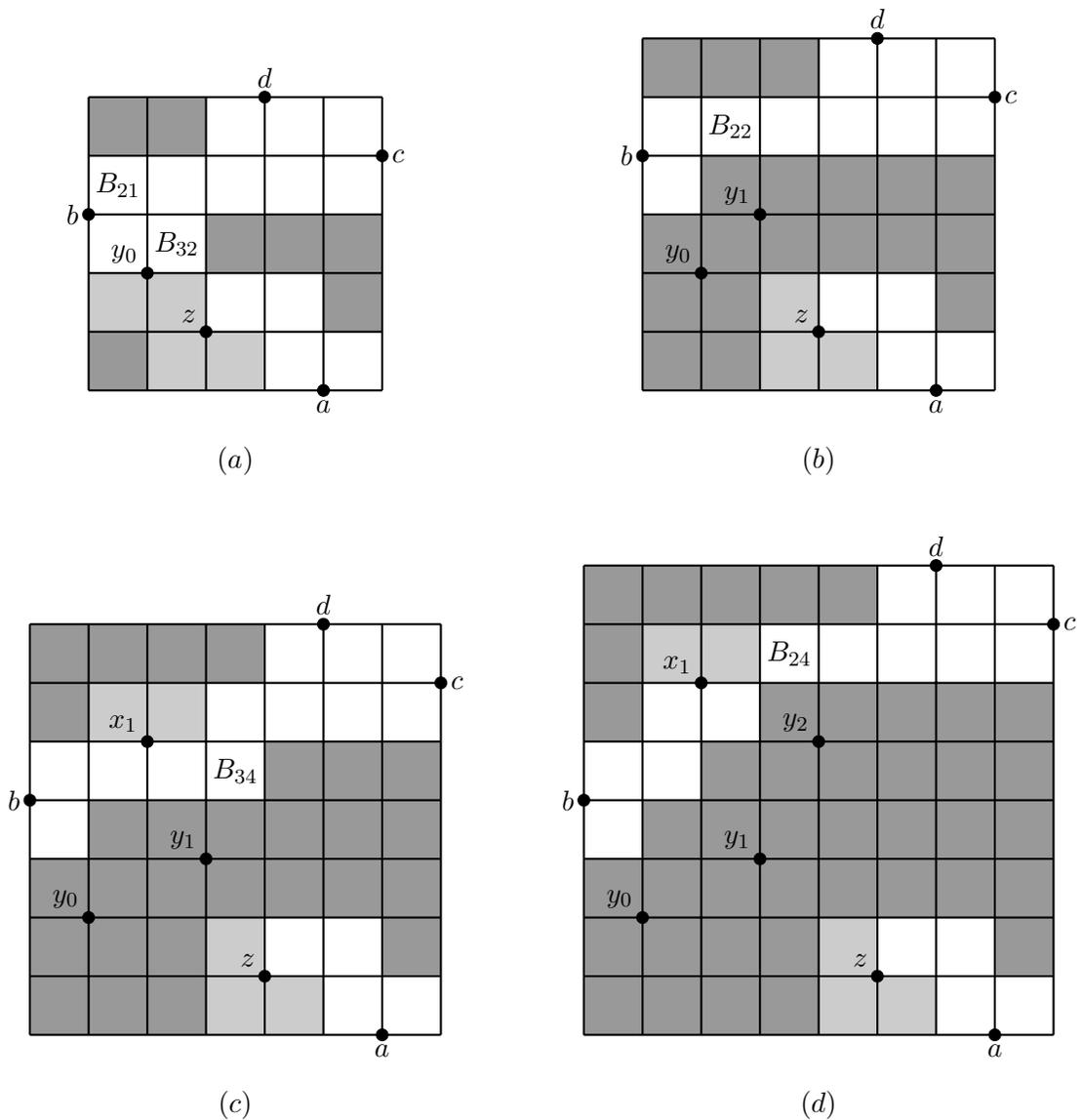

Let us now describe a 231-value chain using the sum we discussed above. A \textit{231-value chain} is a sequence of values of the form\[
\bigoplus_{i=1}^k\hspace{-1.5pt}\raisebox{-3pt}{${}_1$}\,\s_i
\]
where each $\s_i$ is in $\{21,231\}$. The structure of a 231-value chain is shown in Figure \ref{fig:5.5}. The underlined points are optional; if they do not exist, we have added 21 rather than 231. All other points in the figure are required. We will state the lemma describing the structure of points whose positions are in $A$, but in order to prove it, we first establish a few lemmas.\\
\begin{lemma}\label{lem:5.2}
\textit{Let $\pi$ be a simple permutation in $\A'$ of extreme pattern 2413. Then values corresponding to positions in $A$ avoid 321 pattern and 23451 pattern.}
\end{lemma}
\textit{Proof.} We first prove that values of $\pi$ corresponding to positions in $A$ avoid 321 pattern. Suppose to the contrary that values corresponding to positions in $A$ contain 321. Furthermore, assume $b$ corresponds to the 3 of 321. Let $z$ be the least possible value that can play the role of 1 for $b$ and $y_0$ be the least possible value that can play the role of 2 for given $b$ and $z$. Hence, the graph of $\pi$ is as shown in Figure \ref{fig:5.6}$(a)$.\\

Since $\pi$ is simple, we cannot have a block in $\pi$. Hence, there exists a point in either $B_{21}$ or $B_{32}$ of Figure \ref{fig:5.6}$(a)$ to prevent $[\pi^{-1}(b),\pi^{-1}(y_0)]$ from being a block. Suppose there is a point in the region $B_{32}$ and call its value $y_1$. Now, we have a block $[\pi^{-1}(y_0),\pi^{-1}(y_1)]$ as shown in Figure \ref{fig:5.6}$(b)$, so there must exist a point in $B_{22}$ of \ref{fig:5.6}$(b)$. Choose the point of greatest value in $B_{22}$ of \ref{fig:5.6}$(b)$ and let $x_1$ be the value of this point. Then we have the graph shown in Figure \ref{fig:5.6}$(c)$. The block defined by the positions of $b$ and $y_1$ still needs to be split. The only way to do so is by assuming the existence of a point in the region $B_{34}$ of Figure \ref{fig:5.6}$(c)$, say $y_2$. Having $y_2$ results in the graph shown in Figure \ref{fig:5.6}$(d)$.\\

From here on, attempting to split the block $[\pi^{-1}(b),\pi^{-1}(y_i)]$ can only be done by assuming the existence of points with values $x_i$ and $y_{i+1}$ alternatively. In particular, $x_i$ is the greatest possible value such that $\pi^{-1}(y_{i-1})<\pi^{-1}(x_i)<\pi^{-1}(y_i)$ and $x_{i-1}<x_i<c$, and $y_i$ is the value such that $\pi^{-1}(y_{i})<\pi^{-1}(y_{i+1})<\pi^{-1}(z)$ and $x_{i-1}<y_{i+1}<x_i$. However, splitting $[\pi^{-1}(b),\pi^{-1}(y_i)]$ by assuming there exists $x_i$ still results in $[\pi^{-1}(b),\pi^{-1}(y_i)]$ being a block. Similarly, splitting it by assuming there exists $y_{i+1}$ constructs another block $[\pi^{-1}(b),\pi^{-1}(y_{i+1})]$. Hence, attempting to split blocks constructs an infinite strict chain of permutations contained in $\pi$, but this is impossible since the length of $\pi$ is finite.\\

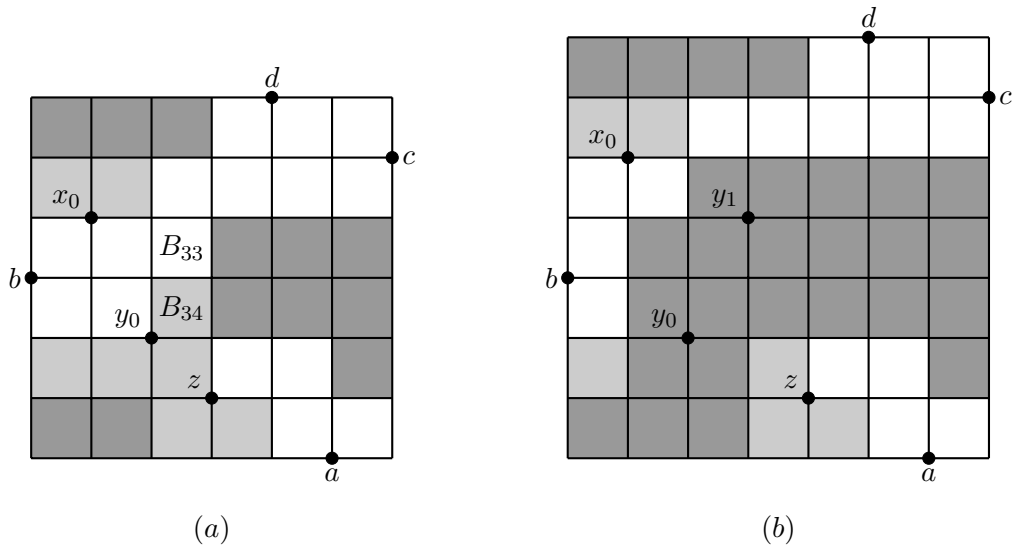
\begin{figure}[b!]
\[\begin{array}{ccc}
  \begin{tikzpicture}[scale=0.8]
    \draw[white, fill=gray!80](1,1)--(3,1)--(3,2)--(1,2)--cycle;
    \draw[white, fill=gray!80](4,3)--(6,3)--(6,2)--(7,2)--(7,5)--(4,5)--cycle;
    \draw[white, fill=gray!80](1,6)--(4,6)--(4,7)--(1,7)--cycle;
    \draw[white, fill=gray!40](1,5)--(3,5)--(3,6)--(1,6)--cycle;
    \draw[white, fill=gray!40](1,2)--(3,2)--(3,1)--(5,1)--(5,2)--(4,2)--(4,4)--(3,4)--(3,3)--(1,3)--cycle;
    \plotPerm{4,5,3,2,7,1,6}
    \node[left] at (1,4) {$b$};
    \node[above left] at (2,5) {$x_0$};
    \node[above left] at (3,3) {$y_0$};
    \node[above left] at (4,2) {$z$};
    \node[above] at (5,7) {$d$};
    \node[below] at (6,1) {$a$};
    \node[right] at (7,6) {$c$};
    \node at (3.5,4.5) {$B_{33}$};
    \node at (3.5,3.5) {$B_{34}$};
  \end{tikzpicture} & \qquad\qquad & \begin{tikzpicture}[scale=0.8]
    \draw[white, fill=gray!80](1,1)--(4,1)--(4,3)--(7,3)--(7,2)--(8,2)--(8,6)--(3,6)--(3,5)--(2,5)--(2,2)--(1,2)--cycle;
    \draw[white, fill=gray!80](1,7)--(5,7)--(5,8)--(1,8)--cycle;
    \draw[white, fill=gray!40](1,6)--(3,6)--(3,7)--(1,7)--cycle;
    \draw[white, fill=gray!40](1,2)--(2,2)--(2,3)--(1,3)--cycle;
    \draw[white, fill=gray!40](4,1)--(6,1)--(6,2)--(5,2)--(5,3)--(4,3)--cycle;
    \plotPerm{4,6,3,5,2,8,1,7}
    \node[left] at (1,4) {$b$};
    \node[above left] at (2,6) {$x_0$};
    \node[above left] at (3,3) {$y_0$};
    \node[above left] at (4,5) {$y_1$};
    \node[above left] at (5,2) {$z$};
    \node[above] at (6,8) {$d$};
    \node[below] at (7,1) {$a$};
    \node[right] at (8,7) {$c$};
  \end{tikzpicture}\\
(a) & & (b)
\end{array}\]
\caption{Partial graphs of $\pi$ with the assumption of $b$ corresponding to the 3 in 321 and having a value in $B_{21}$.}
\label{fig:5.7}
\end{figure}

Referring back to Figure \ref{fig:5.6}$(a)$, assume there is a point in the region $B_{21}$ this time. Let $x_0$ be the greatest value possible of all points in $B_{21}$. We then have a graph shown in Figure \ref{fig:5.7}$(a)$. In order to avoid the segment $[\pi^{-1}(b),\pi^{-1}(y_0)]$ being a block, we must have a point in $B_{33}$ shown in Figure \ref{fig:5.7}$(a)$. Note that $B_{34}$ is shaded in light grey because having a point here is equivalent to having a point in $B_{32}$ of Figure \ref{fig:5.6}$(a)$, which we have discussed previously. If we have a point $y_1$ in $B_{33}$, we have a block $[\pi^{-1}(b),\pi^{-1}(y_1)]$ shown in \ref{fig:5.7}$(b)$. The rest of the argument is identical to the one for the previous case. Again, trying to split blocks constructs an infinite strict chain of permutations contained in $\pi$, so we achieve a contradiction. Therefore, $b$ cannot correspond to the 3 in 321.\\

\begin{figure}[h!]
\[\begin{array}{ccc}
  \begin{tikzpicture}[scale=0.8]
    \draw[white, fill=gray!80](2,1)--(3,1)--(3,3)--(2,3)--cycle;
    \draw[white, fill=gray!80](4,4)--(6,4)--(6,3)--(7,3)--(7,5)--(4,5)--cycle;
    \draw[white, fill=gray!80](1,6)--(4,6)--(4,7)--(1,7)--cycle;
    \draw[white, fill=gray!40](1,4)--(2,4)--(2,6)--(1,6)--cycle;
    \draw[white, fill=gray!40](2,3)--(3,3)--(3,1)--(5,1)--(5,3)--(4,3)--(4,4)--(2,4)--cycle;
    \plotPerm{2,5,4,3,7,1,6}
    \node[left] at (1,2) {$b$};
    \node[above left] at (2,5) {$\hat{b}$};
    \node[above left] at (3,4) {$y_0$};
    \node[above left] at (4,3) {$z$};
    \node[above] at (5,7) {$d$};
    \node[below] at (6,1) {$a$};
    \node[right] at (7,6) {$c$};
  \end{tikzpicture} & \qquad\qquad & \begin{tikzpicture}[scale=0.8]
    \draw[white, fill=gray!80](4,5)--(6,5)--(6,4)--(7,4)--(7,7)--(6,7)--(6,6)--(4,6)--cycle;
    \draw[white, fill=gray!80](2,1)--(5,1)--(5,3)--(3,3)--(3,4)--(2,4)--cycle;
    \draw[white, fill=gray!40](1,5)--(2,5)--(2,7)--(1,7)--cycle;
    \draw[white, fill=gray!40](2,4)--(3,4)--(3,3)--(5,3)--(5,4)--(4,4)--(4,5)--(2,5)--cycle;
    \plotPerm{2,6,5,4,7,1,3}
    \node[left] at (1,2) {$b$};
    \node[above left] at (2,6) {$\hat{b}$};
    \node[above left] at (3,5) {$y_0$};
    \node[above left] at (4,4) {$z$};
    \node[above] at (5,7) {$d$};
    \node[below] at (6,1) {$a$};
    \node[right] at (7,3) {$c$};
  \end{tikzpicture}\\
(a) & & (b)
\end{array}\]
\caption{Partial graphs of $\pi$ with $\hat{b}$ corresponding to the 3 in 321.}
\label{fig:5.8}
\end{figure}
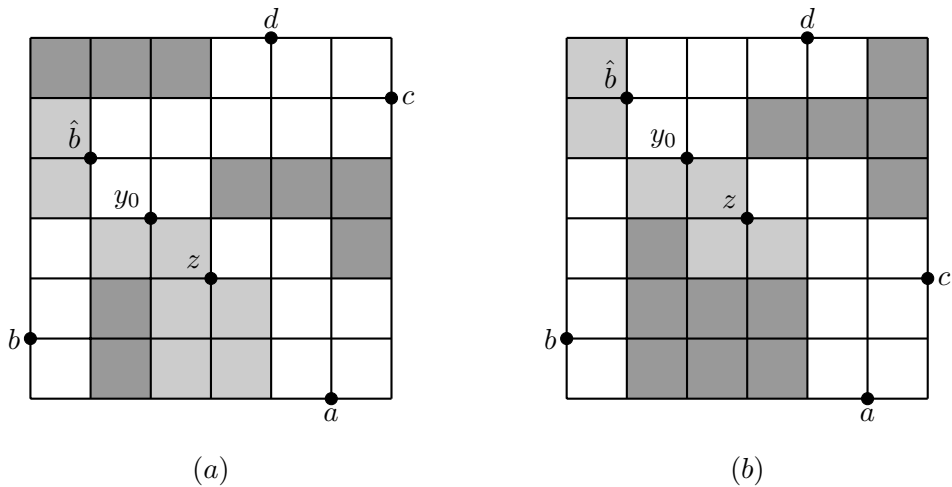

Next, let $\hat{b}$ be the left-most value which can play the role of 3 in 321. Furthermore, let $z$ be the least possible value that can play the role of 1 for $\hat{b}$ and $y_0$ be the least possible value that can play the role of 2 for given $\hat{b}$ and $z$. Depending on weather $\hat{b}<c$ or $c<\hat{b}$, we have one of the graphs shown in Figure \ref{fig:5.8}. For both cases, achieving a contradiction from here is identical to the case of $b$ playing the role of 3 in 321.\\

Consequently, values corresponding to positions in $A$ avoid 321.\\

\begin{figure}[t!]
\[  \begin{tikzpicture}[scale=0.8]
    \draw[white, fill=gray!80](1,6)--(2,6)--(2,7)--(3,7)--(3,8)--(6,8)--(6,9)--(1,9)--cycle;
    \draw[white, fill=gray!80](8,3)--(9,3)--(9,7)--(8,7)--cycle;
    \draw[white, fill=gray!40](1,3)--(6,3)--(6,8)--(3,8)--(3,7)--(2,7)--(2,6)--(1,6)--cycle;
    \draw[white, fill=white](1,3)--(2,3)--(2,4)--(1,4)--cycle;
    \draw[white, fill=white](2,4)--(3,4)--(3,5)--(2,5)--cycle;
    \draw[white, fill=white](3,5)--(3,6)--(4,6)--(4,5)--cycle;
    \draw[white, fill=white](4,6)--(4,7)--(5,7)--(5,6)--cycle;
    \draw[white, fill=white](5,7)--(5,8)--(6,8)--(6,7)--cycle;
    \draw[white, fill=gray!40](6,1)--(7,1)--(7,3)--(6,3)--cycle;
    \plotPerm{2,4,5,6,7,3,9,1,8}
    \node[left] at (1,2) {$b$};
    \node[above left] at (2,4) {$x_1$};
    \node[above left] at (3,5) {$x_2$};
    \node[above left] at (4,6) {$x_3$};
    \node[above left] at (5,7) {$x_4$};
    \node[above left] at (6,3) {$y$};
    \node[above] at (7,9) {$d$};
    \node[below] at (8,1) {$a$};
    \node[right] at (9,8) {$c$};
    \node at (3.5,1.5) {$B_{73}$};
    \node at (3.5,2.5) {$B_{83}$};
    \node at (6.5,5.5) {$B_{46}$};
    \node at (7.5,5.5) {$B_{47}$};
  \end{tikzpicture}\]
\caption{Partial graphs of $\pi$ with the assumption of values corresponding to positions in $A$ contain 23451.}
\label{fig:5.9}
\end{figure}
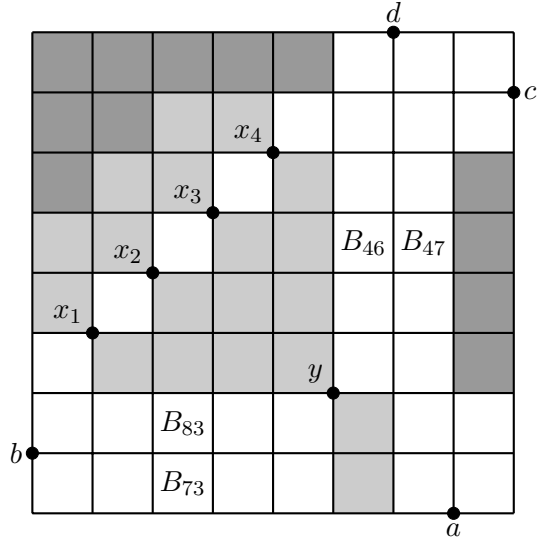

Now, we show that values corresponding to positions in $A$ avoid 23451 pattern as well. Suppose to the contrary that values whose positions are in $A$ contain 23451. Since these values must avoid 321, we obtain Figure \ref{fig:5.9}. Since $[\pi^{-1}(x_2),\pi^{-1}(x_3)]$ is a block, we must assume an existence of a point in either $B_{73}$, $B_{83}$, $B_{46}$ or $B_{47}$. However, if there is a point in either $B_{73}$ or $B_{83}$, then $[\pi^{-1}(x_1),\pi^{-1}(x_2)]$ becomes an unsplittable block. Similarly, assuming a point in either $B_{46}$ or $B_{47}$ implies $[\pi^{-1}(x_3),\pi^{-1}(x_4)]$ is an unsplittable block, so values whose positions are in $A$ cannot form 23451 pattern. We obtain the exact same result even if $a=\pi(1)$ is a part of 23451 pattern.\\

Hence, we have shown that values corresponding to positions in $A$ of $\pi$ cannot contain both 321 pattern and 23451 pattern, so we are done.\hfill$\blacksquare$\\
\\

Before we proceed to the next lemma, let us introduce some new terminologies. First, we call a pair of values forming a 21 pattern with consecutive positions a \textit{descent}. Next, for a permutation $\pi$, the expression\[
\pi=\bigoplus_{i=1}^k \s_i
\]
where $k$ is the greatest value possible is called the \textit{sum-decomposition of $\pi$}. Note that if $k\geq 2$, then $\pi$ is sum-decomposable. Furthermore, let $\pi'$ be the flattening of values of $\pi$ whose positions are in a segment $[a,b]$. Then the above expression for $\pi'$ is called the \textit{sum-decomposition within $[a,b]$.} In both cases, we call the consecutive positions of values corresponding to each $\s_i$ a \textit{sum block}. Notice that a sum block within the whole segment $[1,n]$ is equivalent to a block, but this is in general not true for a sum block within a proper segment $[a,b]$.\\

\begin{figure}[h!]
\[\begin{tikzpicture}[scale=0.8]
    \draw[white, fill=gray!40](1,3)--(3,3)--(3,1)--(5,1)--(5,7)--(2,7)--(2,6)--(1,6)--cycle;
    \draw[white, fill=gray!80](1,6)--(2,6)--(2,7)--(5,7)--(5,8)--(1,8)--cycle;
    \draw[white, fill=gray!80](7,3)--(8,3)--(8,6)--(7,6)--cycle;
    \plotPerm{2,4,5,6,3,8,1,7}
    \node[left] at (1,2) {$b$};
    \node[above left] at (2,4) {$\pi(s)$};
    \node[above left] at (3,5) {$$};
    \node[above left] at (4,6) {$x$};
    \node[above left] at (5,3) {$y$};
    \node[above] at (6,8) {$d$};
    \node[below] at (7,1) {$a$};
    \node[right] at (8,7) {$c$};
    \node at (6.5,5.5) {$B_{36}$};
    \node at (6.5,4.5) {$B_{46}$};
    \node at (6.5,3.5) {$B_{56}$};
    \end{tikzpicture}\]
\caption{Partial graphs of $\pi$ with a 2341 pattern and one descent.}
\label{fig:5.14}
\end{figure}
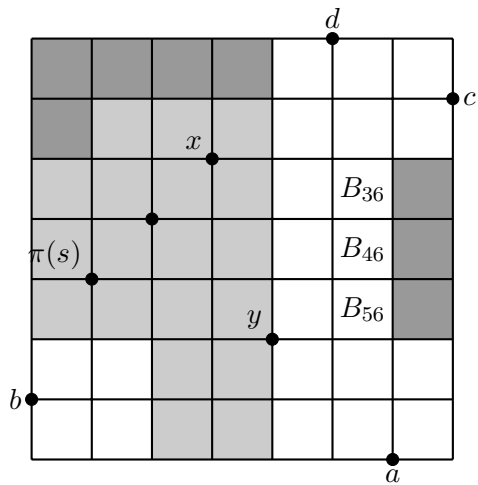

We now state and prove the following lemmas.
\begin{lemma}\label{lem:5.13}
\textit{Let $\pi$ be a simple permutation in $\A'$ of extreme pattern 2413. Then a sum block within the segment A with a descent must end on a descent.}
\end{lemma}
\textit{Proof.} Consider a sum block $I$ within the segment $A$ with at least one descent. Let $x$ and $y$ ($x>y$) be the last descent in $I$. Suppose to the contrary that there is a value located to the right of $y$ within $I$. For any such value $z$, $z<x$ because, otherwise, $\pi^{-1}(z)\notin I$. Also, $y<z$ because 321 pattern is forbidden. Thus, we obtain the graph shown in Figure \ref{fig:5.11}. Attempting to split the block $[\pi^{-1}(y),\pi^{-1}(z)]$ by having a value in $B_{41}$ will cause a construction of an infinite strict chain as observed in the proof of Lemma \ref{lem:5.2}. Therefore, the only way to split the block $[\pi^{-1}(y),\pi^{-1}(z)]$ is by having a point in $B_{23}$. However, this implies $z$ is a part of a descent located to the right of $x$ and $y$, which is a contradiction. Thus, $I$ must end on a descent.\hfill$\blacksquare$
\begin{figure}[t!]
\[\begin{tikzpicture}[scale=0.8]
    \draw[white, fill=gray!40](1,4)--(3,4)--(3,6)--(1,6)--cycle;
    \draw[white, fill=gray!80](1,6)--(4,6)--(4,7)--(1,7)--cycle;
    \draw[white, fill=gray!80](2,1)--(4,1)--(4,3)--(7,3)--(7,5)--(3,5)--(3,4)--(2,4)--cycle;
    \plotPerm{2,5,3,4,7,1,6}
    \node[left] at (1,2) {$b$};
    \node[above left] at (2,5) {$x$};
    \node[above left] at (3,3) {$y$};
    \node[above left] at (4,4) {$z$};
    \node[above] at (5,7) {$d$};
    \node[below] at (6,1) {$a$};
    \node[right] at (7,6) {$c$};
    \node at (3.5,5.5) {$B_{23}$};
    \node at (1.5,3.5) {$B_{41}$};
    \end{tikzpicture}\]
\caption{Partial graph of $\pi$ for Lemma \ref{lem:5.13}.}
\label{fig:5.11}
\end{figure}
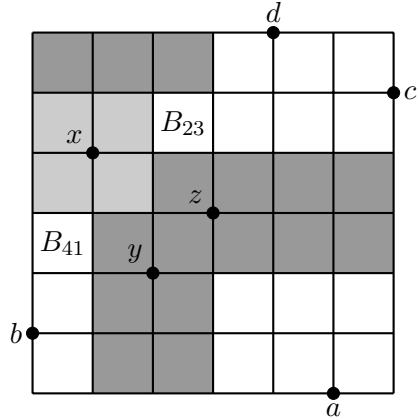
\begin{lemma}\label{lem:5.14}
\textit{Let $\pi$ be a simple permutation in $\A'$ of extreme pattern 2413. Let $[s,t]$ be a sum block within the segment $A$. If the values of positions in the segment $[s,t']$ ($t'\leq t$) ends with a descent, then these values form a 231-value chain.}
\end{lemma}
\textit{Proof.} We prove the statement by induction on the number of descents in $[s,t']$. First, assume that values of positions in the segment $[s,t']$ has only one descent $x$ and $y$ ($x>y$), \textit{i.e.} $\pi(t')=y$. By Lemma \ref{lem:5.2}, values with consecutive positions from $s$ up to $\pi^{-1}(y)$ must form a 21, 231 or 2341 pattern. Suppose to the contrary that we have a 2341 pattern. The graph shown in Figure \ref{fig:5.14} shows this case. As observed, the block formed by the four values for the pattern can only be split by having a value to the right. However, attempting to split this way will immediately cause another block to appear, which is unsplittable. Hence, only a 21 or 231 pattern is permitted in the case which we have a single descent. In either case, we have a 231-value chain.\\

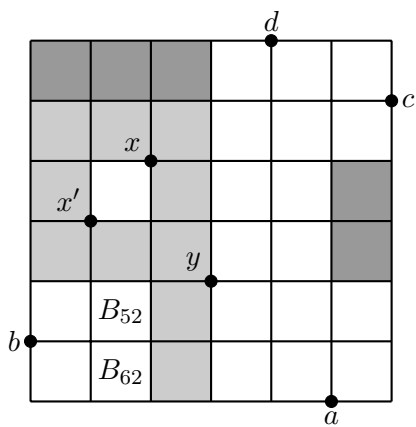
\begin{figure}[H]
\[\begin{tikzpicture}[scale=0.8]
    \draw[white, fill=gray!40](1,3)--(3,3)--(3,1)--(4,1)--(4,6)--(1,6)--cycle;
    \draw[white, fill=gray!80](1,6)--(4,6)--(4,7)--(1,7)--cycle;
    \draw[white, fill=gray!80](6,3)--(7,3)--(7,5)--(6,5)--cycle;
    \draw[white, fill=white](2,4)--(3,4)--(3,5)--(2,5)--cycle;
    \plotPerm{2,4,5,3,7,1,6}
    \node[left] at (1,2) {$b$};
    \node[above left] at (2,4) {$x'$};
    \node[above left] at (3,5) {$x$};
    \node[above left] at (4,3) {$y$};
    \node[above] at (5,7) {$d$};
    \node[below] at (6,1) {$a$};
    \node[right] at (7,6) {$c$};
    \node at (2.5,2.5) {$B_{52}$};
    \node at (2.5,1.5) {$B_{62}$};
    \end{tikzpicture}\]
\caption{Partial graphs of $\pi$ for the inductive case.}
\label{fig:5.12}
\end{figure}

This time, suppose $[s,t']$ has $k$ descents. Let $x$ and $y$ be the ending descent. Also, let $x'$ be the smallest value with $x'>y$ located to the left of $x$. At this point, we have a graph as shown in Figure \ref{fig:5.12}. Notice that it is possible to have one value in $B_{32}$, but no more due to 321 and 23451 avoiding conditions. If we don't have a value in $B_{52}$ or $B_{62}$, then we only have one descent in $[s,t']$, and this is our base case. So suppose we have a value $y'$ in $B_{52}$ or $B_{62}$. Then by inductive hypothesis, values of positions in the segment $[s,\pi^{-1}(y')]$ form a 231-value chain. With $x$, $y$ and potentially another value in $B_{32}$ of Figure \ref{fig:5.12}, we have a 231-value chain with values whose positions are in $[s,t']$, so we are done.
\text{ }\hfill$\blacksquare$\\

By applying Lemma \ref{lem:5.13} and \ref{lem:5.14}, we have the following lemma, which is the description of the structure of points whose positions are in $A$.\\
\begin{lemma}\label{lem:5.3}
\textit{Let $\pi$ be a simple permutation in $\A'$ of extreme pattern 2413. Then values corresponding to positions in $A$ can be expressed as}\begin{equation}\label{eqn:5.2}
\bigoplus_{i=1}^k\s_i,
\end{equation}
\textit{where each $\s_i$ is 1 or a 231-value chain.}
\end{lemma}
\textit{Proof.} By Lemma \ref{lem:5.13}, each sum block in $A$ is either a singleton segment or a segment ending on a descent. Also, by Lemma \ref{lem:5.14}, every segment ending on a descent is a 231-value chain, so we have the desired result.\hfill$\blacksquare$\\

\begin{figure}[h!]
\[\begin{tikzpicture}[scale=.5,rotate=180]
\draw[color=black,fill=black] (1,1) circle [radius=.08];
\draw[color=black,fill=black] (2,2) circle [radius=.08];
\draw[color=black,fill=black] (5,5) circle [radius=.08];
\draw[color=black,fill=black] (6,6) circle [radius=.08];
\draw[color=black,fill=black] (2.5,-0.5) circle [radius=.08];
\draw[color=black,fill=black] (6.5,3.5) circle [radius=.08];
\draw [dotted,ultra thick] (3.5,2) -- (4.5,3);
\draw[color=black,fill=black] (7,7) circle [radius=.08];
\draw[color=black,fill=black] (8,8) circle [radius=.08];
\draw[color=black,fill=black] (8.5,5.5) circle [radius=.08];
\draw[black,thick] (0.8,1.2)--(1.2,1.2) (4.8,5.2)--(5.2,5.2) (6.8,7.2)--(7.2,7.2);
\end{tikzpicture}\]
\caption{Structure of a 312-value chain.}
\label{fig:5.16}
\end{figure}
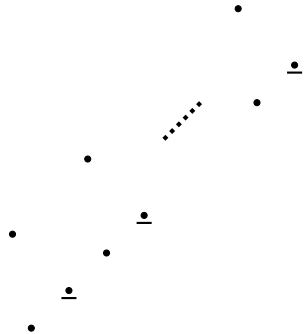

Next, define a \textit{312-value chain} to be the reverse complement of a 231-value chain, \textit{i.e.} it is a sequence of values of the form\[
\bigoplus_{i=1}^k\hspace{-1.5pt}\raisebox{-3pt}{${}_1$}\,\s_i
\]
where each $\s_i$ is in $\{21,312\}$. Hence, the structure of a 312-value chain is as shown in Figure \ref{fig:5.16}. As before, underlined points can be empty, but no others can be so long as the chain continues.\\

With the reverse complement property, we obtain the following lemma.
\begin{lemma}\label{lem:5.4}
\textit{Let $\pi$ be a simple permutation in $\A'$ of extreme pattern 2413. Then values corresponding to positions in $C$ can be expressed as}\begin{equation}\label{eqn:5.3}
\bigoplus_{i=1}^k\s_i,
\end{equation}
\textit{where each $\s_i$ is 1 or a 312-value chain.}
\end{lemma}

Finally, we summarize the structure of a simple permutation $\pi$ in $\A'$ of extreme pattern 2413. By Lemma \ref{lem:5.1}, \ref{lem:5.3} and \ref{lem:5.4}, we have the following proposition.
\begin{proposition}\label{prop:5.2}
\textit{Let $\pi$ be a simple permutation in $\A'$ of extreme pattern 2413. Let $b$, $d$, $a$ and $c$ be the first, the greatest, the least and the last values of $\pi$ respectively. Then values corresponding to positions in $[\pi^{-1}(b),\pi^{-1}(d))$ can be expressed as Equation \ref{eqn:5.2}, values corresponding to positions in $[\pi^{-1}(d),\pi^{-1}(a)]$ can be expressed as Equation \ref{eqn:5.1}, and values corresponding to positions in $(\pi^{-1}(a),\pi^{-1}(c)]$ can be expressed as Equation \ref{eqn:5.3}.}
\end{proposition}

We can state the structure of a simple permutation $\pi$ in $A'$ of extreme pattern 3142 with the inverse property, but in order to do so, we need to define the position-interchange sum, which is the inverted value-interchange sum. Given $\s\in\S_m$ and $\tau\in\S_n$, The \textit{position-interchange sum of $\s$ and $\pi$ with 1 shift} is the permutation defined by\[
\s\oplus^1\tau=\s(1)\s(2)\cdots\s(m-1)\tau'(1)\s(m)\tau'(2)\cdots\tau'(n)
\]
where $\tau'(i)=\tau(i)+m$ for each $i$ with $1\leq i\leq n$. We then define the \textit{231-position chain} and the \textit{312-position chain} as
\[
\bigoplus_{i=1}^k\hspace{-1.5pt}\raisebox{1pt}{${}^1$}\,\s_i\qquad\textrm{where each }\s_i\textrm{ is in \{21,231\}}
\]
and
\[
\bigoplus_{i=1}^k\hspace{-1.5pt}\raisebox{1pt}{${}^1$}\,\s_i\qquad\textrm{where each }\s_i\textrm{ is in \{21,312\}}
\]\\
respectively. By using these, the following proposition is the structure of a simple permutation $\pi$ in $\A'$ of extreme pattern 3142.
\begin{proposition}\label{prop:5.3}
\textit{Let $\pi$ be a simple permutation in $\A'$ of extreme pattern 3142. Let $c$, $a$, $d$ and $b$ be the first, the least, the greatest and the last values of $\pi$ respectively. Then values from the range $[b,c]$ can be expressed as Equation \ref{eqn:5.1}, and values from the range $[a,b)$ and values from the range $(c,d]$ can be expressed as}\begin{equation}\label{eqn:5.4}
\bigoplus_{i=1}^k\s_i\qquad\textit{where each }\s_i\textit{ is 1 or a 312-position chain}
\end{equation}
\textit{and}\begin{equation}\label{eqn:5.5}
\bigoplus_{i=1}^k\s_i\qquad\textit{where each }\s_i\textit{ is 1 or a 231-position chain}
\end{equation}
\textit{respectively.}
\end{proposition}
\subsection{Detailed structures}
Here, we discuss more details of how values of a simple permutation in $\A'$ of extreme 2413 can be placed. We state Proposition \ref{prop:5.4} for simple permutations of extreme pattern 2413, and divide into four lemmas to prove the statement. As before, we apply the inverse symmetry to state the analogous proposition for a permutation of extreme pattern 3142. Afterwards, we state propositions that are converses of Proposition \ref{prop:5.4} and \ref{prop:5.5}, which will be important in Chapter \ref{chap:6}. At the end, we observe how first and last few values of a simple permutation of extreme pattern 2413 can be placed, and what first and last few values of a simple permutation of extreme pattern 3142 can be, so we can establish all the glue sum operations we need to prove the main theorem in the next section.\\

Consider a simple permutation $\pi$ in $\A'$ of extreme pattern 2413. For the following proposition and four lemmas, let $b$, $d$, $a$ and $c$ be the first, the greatest, the least and the last values of $\pi$ respectively, and denote by $A$, $B$ and $C$ each segment $[\pi^{-1}(b),\pi^{-1}(d))$, $[\pi^{-1}(d),\pi^{-1}(a)]$ and $(\pi^{-1}(a),\pi^{-1}(c)]$ respectively.
\begin{proposition}\label{prop:5.4}
\textit{Let $\pi$ be a simple permutation in $\A'$ of extreme pattern 2413. Then:}\begin{enumerate}[\itshape 1.]
\item \textit{in between two values for $1\oplus 1$ in Equation \ref{eqn:5.2}, there exists a value $x$ such that $\pi^{-1}(x)\in B$ or $\pi^{-1}(x)\in C$.}
\item \textit{in between two values for $1\ominus 1$ in Equation \ref{eqn:5.1}, there exists a value $x$ such that $\pi^{-1}(x)\in A$ or $\pi^{-1}(x)\in C$.}
\item \textit{in between two values for $1\oplus 1$ in Equation \ref{eqn:5.3}, there exists a value $x$ such that $\pi^{-1}(x)\in A$ or $\pi^{-1}(x)\in B$.}
\item \textit{for a 231-value chain $\a$ in the segment $A$, let $m$ and $M$ be the minimum and maximum values of $\a$ respectively. Then so long as the chain continues, $\pi^{-1}(x)\in A$ for all $x\in[m,M]\setminus\{M-1\}$, and $\pi^{-1}(M-1)\in B$.}
\item \textit{for a 312-value chain $\b$ in the segment $C$, let $m$ and $M$ be the minimum and maximum values of $\b$ respectively. Then so long as the chain continues, $\pi^{-1}(x)\in C$ for all $x\in[m,M]\setminus\{m+1\}$, and $\pi^{-1}(m+1)\in B$.}
\item \textit{in between two values $s$ and $t$ playing roles of 1 and 2 of 12 in Equation \ref{eqn:5.1}, there exists a value whose position is either in $A$ or $C$. Further, there can be at most four such values $x_1$, $x_2$, $y_1$ and $y_2$ with $x_1<x_2<y_1<y_2$ where $\pi^{-1}(x_1),\pi^{-1}(x_2)\in A$ and $\pi^{-1}(y_1),\pi^{-1}(y_2)\in C$.}
\item \textit{the positions of the value $a+1$ must be in $A$ or $B$ and the position of the value $d-1$ must be in $B$ or $C$.}
\end{enumerate}
\end{proposition}

\vhhh
Note that, unlike $1\oplus 1$, in between $\s_1\oplus 1$, $1\oplus\s_1$ and $\s_1\oplus\s_2$ in Equation \ref{eqn:5.2} with 231-value chains $\s_1$ and $\s_2$, there can be a value $x$ such that $\pi^{-1}(x)\in B$ or $\pi^{-1}(x)\in C$, but there does not have to be such a value. Similarly, in between $12\ominus 1$, $1\ominus 12$ and $12\ominus 12$ in Equation \ref{eqn:5.1}, there can be an optional value $x$ such that $\pi^{-1}(x)\in A$ or $\pi^{-1}(x)\in C$, and in between $\s_1\oplus 1$, $1\oplus\s_1$ and $\s_1\oplus\s_2$ in Equation \ref{eqn:5.3} with 312-value chains $\s_1$ and $\s_2$, there can be an optional value $x$ such that $\pi^{-1}(x)\in A$ or $\pi^{-1}(x)\in B$.\\

We now prove the following lemma, which is for the first three statements of Proposition \ref{prop:5.4}.
\begin{lemma}\label{lem:5.5}
\textit{Let $\pi$ be a simple permutation in $\A'$ of extreme pattern 2413. Then:}\begin{itemize}
\item \textit{in between two values for $1\oplus 1$ in Equation \ref{eqn:5.2}, there exists a value $x$ such that $\pi^{-1}(x)\in B$ or $\pi^{-1}(x)\in C$.}
\item \textit{in between two values for $1\ominus 1$ in Equation \ref{eqn:5.1}, there exists a value $x$ such that $\pi^{-1}(x)\in A$ or $\pi^{-1}(x)\in C$.}
\item \textit{in between two values for $1\oplus 1$ in Equation \ref{eqn:5.3}, there exists a value $x$ such that $\pi^{-1}(x)\in A$ or $\pi^{-1}(x)\in B$.}
\end{itemize}
\end{lemma}
\textit{Proof.} Since the proof for each statement is similar, we only prove the first statement. Let $s$ and $t$ be the values for the first 1 and the second 1 of a $1\oplus 1$ respectively. Since their positions are both in $A$ and consecutive, we cannot have $t=s+1$. Let $x=s+1$. Suppose to the contrary that $\pi^{-1}(x)\in A$. Now, due to the structure described in Lemma \ref{lem:5.3}, if $\pi^{-1}(x)<\pi^{-1}(s)$, then $s$ and $x$ would be parts of a $231$-value chain, so $s$ would not be playing the role of the first 1 of a $1\oplus 1$. On the other hand, if $\pi^{-1}(x)>\pi^{-1}(t)$, then $t$ and $x$ would be a part of a 231-value chain, so we have a contradiction again. Thus, the position of $x$ must be either in $B$ or $C$.\hfill$\blacksquare$\\

Next, we prove the following.
\begin{lemma}\label{lem:5.6}
\textit{Let $\pi$ be a simple permutation in $\A'$ of extreme pattern 2413. Then for a 231-value chain $\a$ in the segment $A$, let $m$ and $M$ be the minimum and maximum values of $\a$ respectively. Then so long as the chain continues, $\pi^{-1}(x)\in A$ for all $x\in[m,M]\setminus\{M-1\}$, and $\pi^{-1}(M-1)\in B$.}
\end{lemma}
\vhhh
\textit{Proof.} We first show $\pi^{-1}(M-1)\in B$. If the position of $M-1$ is in $A$, then it must be a part of $\a$ due to the structure of $\pi$. Further, if $\pi^{-1}(M-1)>\pi^{-1}(M)$, then as discussed in the proof of Lemma \ref{lem:5.3}, $M-1$ must have its descent pair, but the value of this point would be greater than $M$, so the position of $M-1$ must be to the left of $M$. However, this implies the segment corresponding to $\a$ forms a block, so $\pi^{-1}(M-1)\notin A$. If $\pi^{-1}(M-1)\in C$, then $M(M-2)da(M-1)$ or $M(M-3)da(M-1)$ gives us 42513 pattern, so the position of $M-1$ must be in $B$.\\

\begin{figure}[t]
\[\begin{tikzpicture}[scale=.5]
\draw[white, fill=gray!80](0.5,1.5)--(17.5,1.5)--(17.5,8)--(14.5,8)--(14.5,7)--(11.5,7)--(11.5,8)--(0.5,8)--cycle;
\draw[color=black,fill=black] (0.5,-1) circle [radius=.08];
\node[left] at (0.5,-1) {$b$};
\draw [dotted,ultra thick] (1.5,0) -- (2.5,1);
\draw[color=black,fill=black] (11.5,10) circle [radius=.08];
\node[above] at (11.5,10) {$d$};
\draw [dotted,ultra thick] (9.5,8) -- (10.5,9);
\draw[color=black,fill=black] (3,3) circle [radius=.08];
\draw[color=black,fill=black] (4,4) circle [radius=.08];
\draw[color=black,fill=black] (4.5,1.5) circle [radius=.08];
\draw [dotted,ultra thick] (5.5,4) -- (6.5,5);
\draw[color=black,fill=black] (7,7) circle [radius=.08];
\draw[color=black,fill=black] (8,8) circle [radius=.08];
\draw[color=black,fill=black] (8.5,5.5) circle [radius=.08];
\draw[black,thick] (2.8,2.8)--(3.2,2.8) (6.8,6.8)--(7.2,6.8);
\draw[color=black,fill=black] (14.5,-3) circle [radius=.08];
\draw[color=black,fill=black] (13,7.5) circle [radius=.08];
\node[above] at (13,7.5) {$M-1$};
\node[above] at (8,8) {$M$};
\node[above] at (4.5,1.5) {$m$};
\node[below] at (14.5,-3) {$a$};
\node[above left] at (11.5,1.5) {$R_1$};
\node[above] at (13,1.5) {$R_2$};
\node[above right] at (14.5,1.5) {$R_3$};
\draw [thick,dashed] (14.5,-3) -- (14.5,10);
\draw [thick,dashed] (11.5,-3) -- (11.5,10);
\draw [decorate,decoration={brace,amplitude=10pt},xshift=0pt,yshift=0pt]
(8.5,1.3) -- (3,1.3) node [black,midway,,xshift=0.2cm,yshift=-1cm,text width=2.8cm]
{an arbitrary 231-value chain};
\end{tikzpicture}
\]
\caption{Forbidden regions with a 231-value chain.}
\label{fig:5.17}
\end{figure}
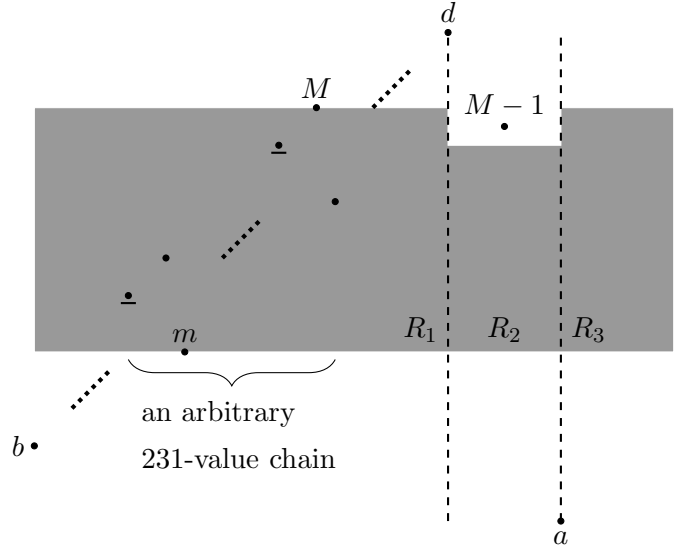

Next, we explain why other values in $[m,M]$ must be involved in the chain by using Figure \ref{fig:5.17}. We cannot have any point in the region $R_1$ by \ref{lem:5.3}. If we have a point in $R_2$, then its value, say $w$, must be between values corresponding to 2 and 1 of a descent 21 or 3 and 1 of 231 in 231-value chain. So suppose $w$ is between values $x$ and $y$ which correspond to 2 and 1 respectively of $\a_i=21$ for some positive integer $i$. Along with the value corresponding to 1 of $\a_{i+1}$, say $z$, $\pi$ contains 52341 or 52431 due to the subsequence $xyzwa$. Thus, we cannot have a point in $R_2$. We can apply the same argument if $w$ is between values 3 and 1 of 231 in 231-value chain. Finally, assume we have a point $w$ in $R_3$. Again, $w$ must be between values corresponding to 2 and 1 of 21 or 3 and 1 of 231 in 231-value chain. Letting $x$ and $y$ denote the same as before, $xydaw$ gives 42513 pattern, so we cannot have a point in $R_3$.\hfill$\blacksquare$\\

With the reverse complement symmetry, the following lemma is an immediate consequence of Lemma \ref{lem:5.6}.
\begin{lemma}\label{lem:5.7}
\textit{Let $\pi$ be a simple permutation in $\A'$ of extreme pattern 2413. Then for a 312-value chain $\b$ in the segment $C$, let $m$ and $M$ be the minimum and maximum values of $\b$ respectively. Then so long as the chain continues, $\pi^{-1}(x)\in C$ for all $x\in[m,M]\setminus\{m+1\}$, and $\pi^{-1}(m+1)\in B$.}
\end{lemma}

Thus, for a 231-value chain with values of positions in $A$ and a 312-value chain with values of positions in $C$ have one value with its position in $B$. We refer to this value multiple times in later discussion, so it is worth naming this value. For a 231-value chain $\a$ in the segment $A$ with the minimum value $m$ and the maximum value $M$, we call the value $M-1$ the \textit{scissor of a 231-value chain}. Similarly, for a 312-value chain $\b$ in the segment $C$ with the minimum value $m$ and the maximum value $M$, we call the value $m+1$ the \textit{scissor of a 312-value chain}.\\

We move onto the next lemma describing the positions of values that are in between 1 and 2 of 12 in Equation \ref{eqn:5.1}.
\begin{lemma}\label{lem:5.8}
\textit{Let $\pi$ be a simple permutation in $\A'$ of extreme pattern 2413. Then in between two values $s$ and $t$ playing roles of 1 and 2 of 12 in Equation \ref{eqn:5.1}, there exists a value whose position is either in $A$ or $C$. Further, there can be at most four such values $x_1$, $x_2$, $y_1$ and $y_2$ with $x_1<x_2<y_1<y_2$ where $\pi^{-1}(x_1),\pi^{-1}(x_2)\in A$ and $\pi^{-1}(y_1),\pi^{-1}(y_2)\in C$.}
\end{lemma}
\vhhh
\textit{Proof.} Given an arbitrary values corresponding to 12 in Equation \ref{eqn:5.1}, call values playing the roles of 1 and 2 of 12 $s$ and $t$ respectively. It is necessary to have at least one value in $[s,t]$ whose position is in $A$ or $C$ due to the simplicity of $\pi$. If we have one value, say $x$ in $A$ and another value, say $y$ in $C$ that are in $[s,t]$, then $y>x$, because, otherwise, 35142 is contained due to $xdsty$. Because each 231-value chain needs the position of its scissor to be in $B$, it is not possible to have all values of an entire 231-value chain to be in between 1 and 2 of 12 as it requires for $B$ to involve 123, 213 or 132 with the scissor, $s$ and $t$, resulting in 52341, 53241 or 52431 pattern containment.\\
\begin{figure}[b!]
\[  \begin{tikzpicture}[scale=0.8]
    \draw[white, fill=gray!80](1,6)--(3,6)--(3,8)--(1,8)--cycle;
    \draw[white, fill=gray!80](4,1)--(6,1)--(6,3)--(8,3)--(8,5)--(7,5)--(7,8)--(5,8)--(5,6)--(4,6)--cycle;
    \draw[white, fill=gray!40](1,3)--(4,3)--(4,6)--(1,6)--cycle;
    \plotPerm{2,4,5,8,3,6,1,7}
    \node[left] at (1,2) {$b$};
    \node[above] at (4,8) {$d$};
    \node[above left] at (2,4) {$x_1$};
    \node[above left] at (3,5) {$x_2$};
    \node[above left] at (5,3) {$s$};
    \node[above left] at (6,6) {$t$};
    \node[below] at (7,1) {$a$};
    \node[right] at (8,7) {$c$};
  \end{tikzpicture}\]
\caption{Partial graph of $\pi$ with the assumption of having two values in $[s,t]$ whose positions are in $A$.}
\label{fig:5.18}
\end{figure}
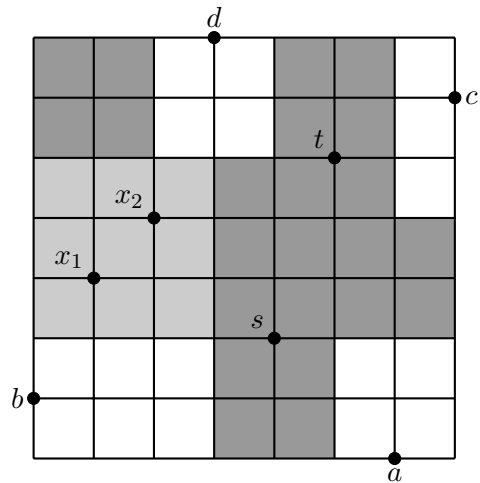

Now, suppose we have two values $x_1$, $x_2$ in $[s,t]$ with $x_1<x_2$ whose positions are in $A$. As we can see in Figure \ref{fig:5.18}, splitting the block $[\pi^{-1}(x_1),\pi^{-1}(x_2)]$ is only possible by having a point with the value less than $s$ and the position in $[\pi^{-1}(x_1),\pi^{-1}(x_2)]$. This makes $x_1$ a part of a 231-value chain with $M=x_1$. Since this is the only way to split $[\pi^{-1}(x_1),\pi^{-1}(x_2)]$, it is impossible to have a $1\oplus 1$ in Equation \ref{eqn:5.2} in between $s$ and $t$. Hence, on the side of $A$, we can have at most two values $x_1$ and $x_2$ with $x_1$ being the maximum of a 231-value chain. Applying the reverse complement symmetry, we also have at most two values $y_1$ and $y_2$ $(y_1<y_2)$ on the side of $C$ with $y_2$ begin the minimum of a 312-value chain in $C$.\hfill$\blacksquare$\\

Finally, we prove the lemma for the last statement of Proposition \ref{prop:5.4}.
\begin{lemma}\label{lem:5.9}
\textit{Let $\pi$ be a simple permutation in $\A'$ of extreme pattern 2413. Then the positions of the value $a+1$ must be in $A$ or $B$ and the position of the value $d-1$ must be in $B$ or $C$.}
\end{lemma}
\vhhh
\textit{Proof.} Suppose to the contrary that $\pi^{-1}(a+1)\in C$. Then $a+1$ cannot play the role of 1 in Equation \ref{eqn:5.3} because, otherwise, $[\pi^{-1}(a),\pi^{-1}(a+1)]$ would be an unsplittable block. So $a+1$ must be playing the role of the minimum value of a 312-value chain. In either case, the position of $a+2$ must be in $B$ due to Lemma \ref{lem:5.7}. However, now we have an unsplittable block $[\pi^{-1}(a+2),\pi^{-1}(\l)]$ where $\l$ is the last value of the 312-value chain that $a+1$ is a part of. Hence, the position of the value $a+1$ cannot be in $C$. Proving $\pi^{-1}(d-1)\notin A$ can be done by the reverse complement argument.\hfill$\blacksquare$\\

With Lemma \ref{lem:5.5}, \ref{lem:5.6}, \ref{lem:5.7}, \ref{lem:5.8} and \ref{lem:5.9}, we have Proposition \ref{prop:5.4}.\\

In Proposition \ref{prop:5.4}, we observed the necessary conditions of simple permutations in $\A'$ of extreme pattern 2413. Proposition \ref{prop:5.5} summarize the inverse statement of these propositions for simple permutations in $\A'$ of extreme pattern 3142. Consider a simple permutation $\pi$ in $\A'$ of extreme pattern 3142. This time, let $c$, $a$, $d$ and $b$ be the first, the least, the greatest and the last values of $\pi$ respectively.
\begin{proposition}\label{prop:5.5}
\textit{Let $\pi$ be a simple permutation in $\A'$ of extreme pattern 3142. Then:}
\begin{enumerate}[\itshape 1.]
\item \textit{in between two positions of values for $1\oplus 1$ in Equation \ref{eqn:5.4}, there exists a value $x$ in $(b,c)$ or in $(c,d]$ such that $\pi^{-1}(s)\leq\pi^{-1}(x)\leq\pi^{-1}(t)$.}
\item \textit{in between two positions of values for $1\ominus 1$ in Equation \ref{eqn:5.1}, there exists a value $x$ in $[a,b)$ or in $(c,d]$ such that $\pi^{-1}(s)\leq\pi^{-1}(x)\leq\pi^{-1}(t)$.}
\item \textit{in between two positions of values for $1\oplus 1$ in Equation \ref{eqn:5.5}, there exists a value $x$ in $[a,b)$ or in $(b,c)$ such that $\pi^{-1}(s)\leq\pi^{-1}(x)\leq\pi^{-1}(t)$.}
\item \textit{for a 312-value chain $\a$ in the range $[a,b)$, let $m$ and $M$ be the values of $\a$ whose positions are the first and the last in $\a$ respectively. Then so long as the chain continues, $x\in[a,b)$ for all $x$ with $\pi^{-1}(x)\in[\pi^{-1}(m),\pi^{-1}(M)]\setminus\{\pi^{-1}(M)-1\}$ and $\pi(\pi^{-1}(M)-1)\in (b,c)$.}
\item \textit{for a 231-value chain $\b$ in the range $(c,d]$, let $m$ and $M$ be the values of $\b$ whose positions are the first and the last in $\b$ respectively. Then so long as the chain continues, $x\in(c,d]$ for all $x$ with $\pi^{-1}(x)\in[\pi^{-1}(m),\pi^{-1}(M)]\setminus\{\pi^{-1}(m)+1\}$ and $\pi(\pi^{-1}(m)+1)\in (b,c)$.}
\item \textit{there exists a value in $[a,b)$ or $(c,d]$ whose position is in between $\pi^{-1}(s)$ and $\pi^{-1}(t)$ where $s$ and $t$ are values playing roles of 1 and 2 of 12 in Equation \ref{eqn:5.1}. Further, there could be at most four such values $x_1$, $x_2$, $y_1$ and $y_2$ with $\pi^{-1}(x_1)<\pi^{-1}(x_2)<\pi^{-1}(y_1)<\pi^{-1}(y_2)$ where $x_1,x_2\in [a,b)$ and $y_1,y_2\in (c,d]$.}
\item \textit{the value of the position $\pi^{-1}(c)+1$ must be in $[a,b)$ or $[b,c]$ and the value of the position $\pi^{-1}(b)-1$ must be in $[b,c]$ or $(c,d]$.}
\end{enumerate}
\end{proposition}

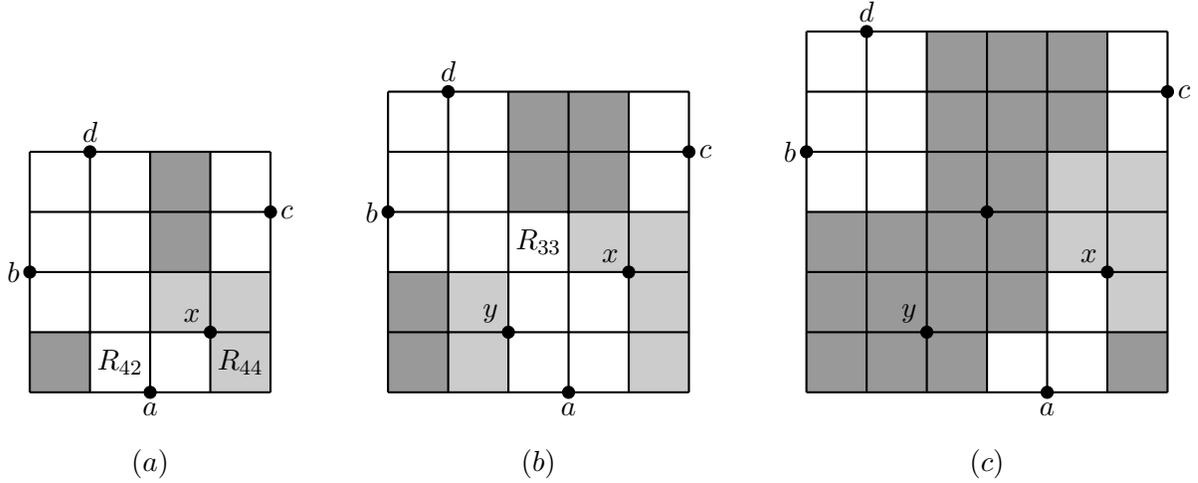
\begin{figure}[h!]
\[\begin{array}{ccccc}
  \begin{tikzpicture}[scale=0.8]
    \draw[white, fill=gray!80](1,1)--(2,1)--(2,2)--(1,2)--cycle;
    \draw[white, fill=gray!80](3,3)--(4,3)--(4,5)--(3,5)--cycle;
    \draw[white, fill=gray!40](3,2)--(4,2)--(4,1)--(5,1)--(5,3)--(3,3)--cycle;
    \plotPerm{3,5,1,2,4}
    \node[left] at (1,3) {$b$};
    \node[above] at (2,5) {$d$};
    \node[below] at (3,1) {$a$};
    \node[above left] at (4,2) {$x$};
    \node[right] at (5,4) {$c$};
    \node at (2.5,1.5) {$R_{42}$};
    \node at (4.5,1.5) {$R_{44}$};
  \end{tikzpicture} & & \begin{tikzpicture}[scale=0.8]
    \draw[white, fill=gray!80](1,1)--(2,1)--(2,3)--(1,3)--cycle;
    \draw[white, fill=gray!40](2,1)--(3,1)--(3,3)--(2,3)--cycle;
    \draw[white, fill=gray!80](3,4)--(5,4)--(5,6)--(3,6)--cycle;
    \draw[white, fill=gray!40](4,3)--(5,3)--(5,1)--(6,1)--(6,4)--(4,4)--cycle;
    \plotPerm{4,6,2,1,3,5}
    \node[left] at (1,4) {$b$};
    \node[above] at (2,6) {$d$};
    \node[above left] at (3,2) {$y$};
    \node[below] at (4,1) {$a$};
    \node[above left] at (5,3) {$x$};
    \node[right] at (6,5) {$c$};
    \node at (3.5,3.5) {$R_{33}$};
  \end{tikzpicture} & & \begin{tikzpicture}[scale=0.8]
    \draw[white, fill=gray!80](1,1)--(4,1)--(4,2)--(5,2)--(5,5)--(6,5)--(6,7)--(3,7)--(3,4)--(1,4)--cycle;
    \draw[white, fill=gray!40](5,3)--(6,3)--(6,2)--(7,2)--(7,5)--(5,5)--cycle;
    \draw[white, fill=gray!80](6,1)--(7,1)--(7,2)--(6,2)--cycle;
    \plotPerm{5,7,2,4,1,3,6}
    \node[left] at (1,5) {$b$};
    \node[above] at (2,7) {$d$};
    \node[above left] at (3,2) {$y$};
    \node[below] at (5,1) {$a$};
    \node[above left] at (6,3) {$x$};
    \node[right] at (7,6) {$c$};
  \end{tikzpicture}\\
(a) & & (b) & & (c)
\end{array}\]
\caption{Partial graphs of $\pi$ with the assumption of having a value $x$ that is not a part of a 312-value chain.}
\label{fig:5.19}
\end{figure}

\vhhh
So far, we have proven that if a permutation $\pi$ of extreme pattern 2413 is in $\textrm{Si}(\A')$, then $\pi$ satisfies every condition listed in Proposition \ref{prop:5.4} as well as the structural conditions described in Proposition \ref{prop:5.2}. Proposition \ref{prop:5.6} states that if an arbitrary permutation $\pi$ of extreme pattern 2413 having a structure described in Proposition \ref{prop:5.2} and \ref{prop:5.4}, then $\pi\in\textrm{Si}(\A')$. In order to prove this statement, we need the following lemma. As usual, let $b$, $d$, $a$ and $c$ be the first, the greatest, the least and the last values of $\pi$.
\begin{lemma}\label{lem:5.10}
There is no value $x$ less than $b$ such that $\pi^{-1}(x)\in C$.
\end{lemma}
\vhhh
\textit{Proof.} Suppose to the contrary that there exists a value less than $b$ whose position is in $C$. Let $x$ be the greatest such value. First, assume $x$ is not a part of a 312-value chain, so we have the graph shown in Figure \ref{fig:5.19}$(a)$. The only way to split the block $[\pi^{-1}(a),\pi^{-1}(x)]$ is by having a point in the region denoted by $R_{42}$. Denoting by $y$ the value of the right-most point in $R_{42}$ of \ref{fig:5.19}$(a)$, we obtain the graph in \ref{fig:5.19}$(b)$. We attempt to split the block $[\pi^{-1}(y),\pi^{-1}(x)]$ by having a point in $R_{33}$, but the consequential graph has an unsplittable block $[\pi^{-1}(y),\pi^{-1}(x)]$ as shown in \ref{fig:5.19}$(c)$.\\
\begin{figure}[p]
\[\begin{array}{ccc}
  \begin{tikzpicture}[scale=0.8]
    \draw[white, fill=gray!80](1,1)--(2,1)--(2,3)--(1,3)--cycle;
    \draw[white, fill=gray!80](3,4)--(5,4)--(5,6)--(3,6)--cycle;
    \draw[white, fill=gray!40](3,3)--(5,3)--(5,1)--(6,1)--(6,4)--(3,4)--cycle;
    \plotPerm{4,6,1,3,2,5}
    \node[left] at (1,4) {$b$};
    \node[above] at (2,6) {$d$};
    \node[below] at (3,1) {$a$};
    \node[above left] at (4,3) {$x$};
    \node[above left] at (5,2) {$z$};
    \node[right] at (6,5) {$c$};
    \node at (5.5,2.5) {$R_{45}$};
    \node at (3.5,1.5) {$R_{53}$};
  \end{tikzpicture} & \qquad & \begin{tikzpicture}[scale=0.8]
    \draw[white, fill=gray!80](1,1)--(2,1)--(2,2)--(5,2)--(5,3)--(6,3)--(6,7)--(3,7)--(3,4)--(1,4)--cycle;
    \draw[white, fill=gray!40](6,1)--(7,1)--(7,5)--(6,5)--cycle;
    \plotPerm{5,7,3,1,4,2,6}
    \node[left] at (1,5) {$b$};
    \node[above] at (2,7) {$d$};
    \node[above left] at (3,3) {$z+1$};
    \node[below] at (4,1) {$a$};
    \node[above left] at (5,4) {$x$};
    \node[above left] at (6,2) {$z$};
    \node[right] at (7,6) {$c$};
    \node at (2.5,1.5) {$R_{62}$};
  \end{tikzpicture}\\
(a) & & (b)
\end{array}\]
\[\begin{array}{c}
\begin{tikzpicture}[scale=0.8]
    \draw[white, fill=gray!80](1,1)--(4,1)--(4,2)--(5,2)--(5,3)--(6,3)--(6,4)--(7,4)--(7,8)--(3,8)--(3,5)--(1,5)--cycle;
    \draw[white, fill=gray!80](6,1)--(8,1)--(8,2)--(6,2)--cycle;
    \draw[white, fill=gray!40](7,2)--(8,2)--(8,6)--(7,6)--cycle;
    \plotPerm{6,8,2,4,1,5,3,7}
    \node[left] at (1,6) {$b$};
    \node[above] at (2,8) {$d$};
    \node[above left] at (3,2) {$y$};
    \node[above left] at (4,4) {$z+1$};
    \node[below] at (5,1) {$a$};
    \node[above left] at (6,5) {$x$};
    \node[above left] at (7,3) {$z$};
    \node[right] at (8,7) {$c$};
  \end{tikzpicture}\\
(c)\end{array}\]
\[\begin{array}{ccc}
\begin{tikzpicture}[scale=0.8]
    \draw[white, fill=gray!80](1,1)--(2,1)--(2,5)--(1,5)--cycle;
    \draw[white, fill=gray!80](6,1)--(8,1)--(8,3)--(7,3)--(7,2)--(6,2)--cycle;
    \draw[white, fill=gray!80](3,6)--(7,6)--(7,8)--(3,8)--cycle;
    \draw[white, fill=gray!40](2,2)--(3,2)--(3,5)--(2,5)--cycle;
    \draw[white, fill=gray!40](3,5)--(7,5)--(7,3)--(8,3)--(8,6)--(3,6)--cycle;
    \plotPerm{6,8,3,1,2,5,4,7}
    \node[left] at (1,6) {$b$};
    \node[above] at (2,8) {$d$};
    \node[above left] at (3,3) {$m+1$};
    \node[below] at (4,1) {$a$};
    \node[above left] at (5,2) {$m$};
    \node[above left] at (6,5) {$x$};
    \node[above left] at (7,4) {$z$};
    \node[right] at (8,7) {$c$};
    \node at (2.5,1.5) {$R_{82}$};
  \end{tikzpicture} & \qquad & \begin{tikzpicture}[scale=0.8]
    \draw[white, fill=gray!80](1,1)--(4,1)--(4,2)--(5,2)--(5,7)--(8,7)--(8,9)--(3,9)--(3,4)--(2,4)--(2,6)--(1,6)--cycle;
    \draw[white, fill=gray!80](6,1)--(9,1)--(9,4)--(8,4)--(8,3)--(7,3)--(7,2)--(6,2)--cycle;
    \draw[white, fill=gray!40](2,4)--(3,4)--(3,6)--(2,6)--cycle;
    \draw[white, fill=gray!40](5,6)--(8,6)--(8,4)--(9,4)--(9,7)--(5,7)--cycle;
    \plotPerm{7,9,2,4,1,3,6,5,8}
    \node[left] at (1,7) {$b$};
    \node[above] at (2,9) {$d$};
    \node[above left] at (3,2) {$y$};
    \node[above left] at (4,4) {$m+1$};
    \node[below] at (5,1) {$a$};
    \node[above left] at (6,3) {$m$};
    \node[above left] at (7,6) {$x$};
    \node[above right] at (8,5) {$z$};
    \node[right] at (9,8) {$c$};
  \end{tikzpicture}\\
  (d) & & (e)\end{array}
\]
\caption{Partial graphs of $\pi$ with the assumption of having a value $x$ that is a part of a 312-value chain.}
\label{fig:5.20}
\end{figure}

Next, suppose $x$ is a part of a 312-value chain. If $x$ plays the role of 2 of 21 in the chain, we have a point with the value $z$ playing the role of 1 of 21 in the region $R_{44}$ of Figure \ref{fig:5.19}, resulting in the graph shown in Figure \ref{fig:5.20}$(a)$. Note that if $x$ played the role of $3$ of 312 in the chain, then $z$ shown in Figure \ref{fig:5.20}$(a)$ would be for 1 of 312, and we would have another value for 2 of 312 in the region $R_{45}$. The rest of the proof is the same for the case of $x$ playing the role of 2 of 21, so we only show the case of $x$ for 2 of 21.\\

Suppose $z$ is the minimum value in the 312-value chain. Then as discussed before, the position of $z+1$ must be in $B$. Thus, we have the graph shown in Figure \ref{fig:5.20}$(b)$. Splitting the block $[\pi^{-1}(y_1),\pi^{-1}(z)]$ can be only done by having another value $y$ in the region $R_{62}$, but then we have an unsplittable block $[\pi^{-1}(y),\pi^{-1}(z)]$ as shown in Figure \ref{fig:5.20}$(c)$, so we have a contradiction. Now, suppose $z$ is not the minimum value in the 312-value chain. That means we have the minimum value $m$ of the 312-value chain in $R_{53}$ in Figure \ref{fig:5.20}$(a)$. With the value $m+1$ whose position is in $B$, we have the graph as in Figure \ref{fig:5.20}$(d)$. The segment $[\pi^{-1}(m+1),\pi^{-1}(z)]$ is a block, and the only way to split it is by having a point with the value $y$ in $R_{82}$. Consequently, we obtain the graph shown in Figure \ref{fig:5.20}$(e)$, but $[\pi^{-1}(y),\pi^{-1}(z)]$ is an unsplittable box. This completes the proof of the lemma.\hfill$\blacksquare$\\

With the reverse complement property of Lemma \ref{lem:5.9}, we have the following lemma. With these lemmas, we are ready to prove Proposition \ref{prop:5.6}.
\begin{lemma}\label{lem:5.11}
\textit{There is no value $x$ greater than $a$ such that $\pi^{-1}(x)\in A$.}
\end{lemma}
\begin{proposition}\label{prop:5.6}
\textit{Let $\pi$ be a permutation of extreme pattern 2413 where $|\pi|=n$. If $\pi$ obeys all the conditions of Proposition \ref{prop:5.2} and \ref{prop:5.4}, then $\pi\in\textrm{Si}(\A')$.}
\end{proposition}
\vhhh
\textit{Proof.} First, we prove $\pi$ is simple. Suppose $\pi$ is not simple, and let $I$ be a non-singleton proper block of $\pi$. Let $b$, $d$, $a$ and $c$ be the first, the greatest, the least and the last values of $\pi$ respectively. Additionally, let $A$, $B$ and $C$ be segments $[\pi^{-1}(b),\pi^{-1}(d))$, $[\pi^{-1}(d),\pi^{-1}(a)]$ and $(\pi^{-1}(a),\pi^{-1}(c)]$ respectively. Obviously, $I$ cannot contain both the positions of $b$ and $c$ because this implies $I=[\pi^{-1}(b),\pi^{-1}(c)]=[1,n]$. The same result can be quickly obtained if $I$ contains any two positions of $b$, $d$, $a$ and $c$. For example, assume $I$ contains the positions of $b$ and $d$. Since $b<c<d$, we must have $\pi^{-1}(c)\in I$, so again, $I=[1,n]$. Same can be done for any other two of $b$, $d$, $a$ and $c$.\\

Assume $I$ contains only one position of $b$, $d$, $a$ or $c$. We first assume $\pi^{-1}(d)\in I$. Then $\pi^{-1}(d-1)$ also has to be in $I$. By the seventh condition listed in Proposition \ref{prop:5.4}, we have the assumption $\pi^{-1}(d-1)\notin A$. If $\pi^{-1}(d-1)\in C$, then $\pi^{-1}(a)$ has to be in $I$, which is a contradiction. So suppose $\pi^{-1}(d-1)\in B$. Note that $d$ and $d-1$ cannot be forming $1\ominus 1$ due to the second condition listed in Proposition \ref{prop:5.4}. Hence, $d$ and $d-1$ are 1 of the initial 1 and 2 of the 12 in the expression $1\ominus 12\ominus\cdots$. However, by the sixth condition of Proposition \ref{prop:5.4}, we must have a value in between 1 and 2 of 12 whose position is either in $A$ or $C$, \textit{i.e.} $\pi^{-1}(d-2)$ must be in $A$ or $C$. Moreover, by Lemma \ref{lem:5.10}, every value whose position is in $A$ must be less than $c$, so $\pi^{-1}(d-2)\in C$. This implies $\pi^{-1}(a)\in I$, which is a contradiction. Hence, $I$ cannot contain the position of $d$. We have the same result for the position of $a$ by the reverse complement argument.\\

Thus, suppose $I$ contains $1=\pi^{-1}(b)$. Because $I$ does not contain the position of $d$, we have $I\sbst A$. Now, if $b$ is a part of a 231-value chain, then $[m,M]\sbst \pi(I)$ where $m$ and $M$ are the minimum and maximum values of the 231-value chain that $b$ is a part of. However, since we know $\pi^{-1}(M-1)\in B$ by the fourth condition of Proposition \ref{prop:5.4}, this implies $\pi^{-1}(d)\in I$. So assume $b$ is playing the role of the first 1 in Equation \ref{eqn:5.2}. If $\pi(2)$ is also playing the role of another 1 in Equation \ref{eqn:5.2}, then there must be a value $x$ such that $b<x<\pi(2)$ where $\pi^{-1}(x)$ is in $B$ or in $C$, so again, we have a contradiction of $\pi^{-1}(d)\in I$. Hence, $\pi(2)$ is a part of a 231-value chain, but this implies $[m,M]\sbst \pi(I)$ where $m$ and $M$ are the minimum and maximum values of the 231-value chain that $\pi(2)$ is a part of, resulting in $\pi^{-1}(d)\in I$ again. Similar result can be obtained for the case of $I$ containing $n=\pi^{-1}(c)$. Consequently, $I$ cannot contain any of $\pi^{-1}(b)$, $\pi^{-1}(d)$, $\pi^{-1}(a)$ or $\pi^{-1}(c)$.\\

So suppose $I$ does not contain any position of $b$, $d$, $a$ or $c$. This implies $I\sbst A$, $I\sbst B$ or $I\sbst C$. Assume $I\sbst A$. Then $\pi(I)$ must contain at least two values $s$ and $t$ whose positions are consecutive in $A$ ($\pi^{-1}(s)<\pi^{-1}(t)$). If $s$ and $t$ are playing roles of a $1\oplus 1$ in Equation \ref{eqn:5.2}, then there exists a value $x$ such that $s<x<t$ where $\pi^{-1}(x)$ is in $B$ or in $C$. As before, this implies $\pi^{-1}(d)\in I$, which is a contradiction. If $s$ and $t$ are both parts of the same 231-value chain, then $[m,M]\sbst \pi(I)$ where $m$ and $M$ are the minimum and maximum values of the 231-value chain that $s$ and $t$ are parts of. Since $\pi^{-1}(M-1)\in B$, this implies $\pi^{-1}(d)\in I$ again, so we have a contradiction. Similar result can be shown if $s$ is playing a role of 1 in Equation \ref{eqn:5.2} and $t$ is a part of a 231-value chain, or vice versa. Thus, we cannot have $I\sbst A$. By the reverse complement argument, $I\sbst C$ also cannot be true.\\

Finally, suppose $I\sbst B$. Then $\pi(I)$ must contain at least two values $s$ and $t$ whose positions are consecutive in $B$ ($\pi^{-1}(s)<\pi^{-1}(t)$). If $s$ and $t$ are playing roles of a $1\ominus 1$ in Equation \ref{eqn:5.1}, then there exists a value $x$ whose position is in $A$ or $C$. Hence, $x\in\pi(I)$, so $\pi^{-1}(x)\in I$. If $\pi^{-1}(x)\in A$, then the position of $d$ would be in $I$. On the other hand, if $\pi^{-1}(x)\in C$, then the position of $a$ would be in $I$. Thus, in either way, we achieve a contradiction. Similarly, if $s$ and $t$ are playing roles of 1 and 2 of a 12 in Equation \ref{eqn:5.1}, then there exists a value $x$ whose position is in $A$ or $C$. Again either the position of $d$ or $a$ would be in $I$, which is a contradiction. Similar result can be shown for the case if $s$ is playing a role of 2 of a 12 and $t$ is playing a role of 1 of a 1 in Equation \ref{eqn:5.1}, and for the case if $s$ is playing a role of 1 of a 1 and $t$ is playing a role of 1 of a 12 in Equation \ref{eqn:5.1}.\\

Consequently, we have considered every case of where a non-singleton proper block $I$ can belong to, and we achieved a contradiction in each case. Hence, $\pi$ is simple.\\

Next, we prove $\pi$ avoids every $\b$ in $\{52341,53241,52431,35142,42513,351624\}$. Suppose to the contrary that $\b\cont\pi$ for some $\b\in\{52341,53241,52431,35142,42513,351624\}$. We first let $\b\in\{52341,53241,52431\}$. Notice that $\LRmax(\pi)=\{x:\pi^{-1}(x)\in A\textrm{ and }x\textrm{ plays}\linebreak\textrm{the role of 2 or 3 of a 231-value chain, or 1 of 1 in Equation \ref{eqn:5.2}}\}\cup\{d\}$ and $\RLmin(\pi)=\{x:\pi^{-1}(x)\in C\textrm{ and }x\textrm{ plays the role of 1 or 2 of a 312-value chain, or 1 of 1 in Equation \ref{eqn:5.3}}\}\cup\{a\}$. If any two values from these two sets play the role of 2, 3 or 4, we cannot find a value corresponding to either 5 or 1. Thus, each value playing the role of 2, 3 and 4 must be in one of the following three sets.\begin{itemize}
\item $S_1=\{x:\pi^{-1}(x)\in A,\,x\textrm{ corresponds to 1 of 21 or 231 in a 231-value chain}\}$
\item $S_2=\{x:\pi^{-1}(x)\in B,\,x\neq 1,n\}$
\item $S_3=\{x:\pi^{-1}(x)\in C,\,x\textrm{ corresponds to 2 of 21 or 3 of 312 in a 312 value chain}\}$
\end{itemize}

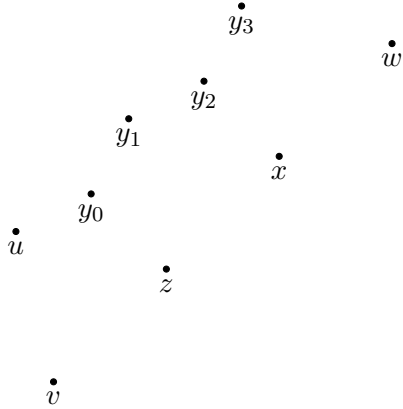
\begin{figure}[t!]
\[\begin{tikzpicture}[scale=.5]
\draw[color=black,fill=black] (-1,1) circle [radius=.08] node [below]{$u$};
\draw[color=black,fill=black] (0,-3) circle [radius=.08] node [below]{$v$};
\draw[color=black,fill=black] (1,2) circle [radius=.08] node [below]{$y_0$};
\draw[color=black,fill=black] (2,4) circle [radius=.08] node [below]{$y_1$};
\draw[color=black,fill=black] (3,0) circle [radius=.08] node [below]{$z$};
\draw[color=black,fill=black] (4,5) circle [radius=.08] node [below]{$y_2$};
\draw[color=black,fill=black] (5,7) circle [radius=.08] node [below]{$y_3$};
\draw[color=black,fill=black] (6,3) circle [radius=.08] node [below]{$x$};
\draw[color=black,fill=black] (9,6) circle [radius=.08] node [below]{$w$};
\end{tikzpicture}\]
\caption{Partial graph of $\pi$ to show that $52341,53241,52431\avd\pi$.}
\label{fig:5.29}
\end{figure}

Let $x\in\S_1$. Furthermore, suppose $x$ is a value corresponding to 1 of 231 in some 231-value chain. We show that $x$ cannot play the role of any of 2, 3 and 4. Suppose to the contrary that $x$ plays the role of one of the values of $\b(2)$, $\b(3)$ or $\b(4)$. There are three values $y_1$, $y_2$ and $y_3$ that are greater than $x$ located to the left of $x$ in increasing order, as shown in Figure \ref{fig:5.29}, so one of them has to play the role of $\b(1)=5$. As we choose one of them to be playing the role of $5$, other two cannot be assigned to any other values of $\b$, as these two values are either located to the left of, or greater than the value for $\b(1)=5$.\\

Suppose $x$ plays the role of $\b(4)$. If $\b=52341$ or $53241$, then $\b(4)=4$, but then we need two values less than $x$ whose positions are in between $\pi^{-1}(y_i)$ and $\pi^{-1}(x)$ which we don't have for each of $y_1$, $y_2$ and $y_3$. If $\b=52431$, then we assign $y_1$ and its descent pair $z$ to play the roles of $5$ and $2$ respectively, but we do not have any value for $4$. Hence, $x$ cannot correspond to $\b(4)$.\\

This time, assume $x$ plays the role of $\b(3)$. For $\b=53241$, we don't have any value corresponding to 3. For $\b=52341$, we assign $y_1$ and $z$ to 5 and 2 respectively as this is the only way to have a value for 2, but since $y_1=x+1$, we do not have a value corresponding to 4. Lastly, for $\b=52431$, we again assign $y_1$ and $z$ to 5 and 2 respectively. However, if $y_1$ and $z$ correspond to 3 and 1 of 231 in the 231-value chain respectively as in Figure \ref{fig:5.29}, then the only values greater than $z$ and less than $x$ is $y_0$ and $u$ shown in Figure \ref{fig:5.29}, which are located to the left of $x$, so we don't have a value for 3. On the other hand, if $y_1$ and $z$ correspond to 2 and 1 of 21 in the 231-value chain, then the only value in between $z$ and $x$ is $u$. It is also possible that $u$ and $v$ do not exist, as $y_1$ and $z$ correspond to the first summand of the 231-value chain expression. In this case, there is no value greater than $z$ and less than $x$, because $x=z+1$, so we do not have a value for 3 in either way. Therefore, $x$ cannot play the role of $\b(3)$.\\

The last case is $x$ playing the role of $\b(2)$. If $y_3$ plays the role of $\b(1)=5$, then there is one value $w$ that is greater than $x$ and less than $y_3$ whose position is to the right of $x$. Note that $w$ can be a value in $S_1$, \textit{i.e.} a part of the next summand in the 231-value chain expression, or $\pi^{-1}(w)\in B$, if the 231 that $y_3$ and $x$ are parts of is the last summand in the expression. If $\b=52341$ or $52431$, we do not have enough values to assign to 3 and 4. If $\b=53241$, then $w$ must play the role of 4, but since there is no value less than $x$ whose position is in between $\pi^{-1}(x)$ and $\pi^{-1}(w)$, we do not have a value corresponding to 2 of $\b$.\\

Consequently, we achieve a contradiction in every case. If $x$ is a value corresponding to 1 of 21 in a 231-value chain, then we obtain the same result, since the only difference is the absence of $y_2$ in Figure \ref{fig:5.29}. So $x$ cannot play the role of any value of $\b(2)$, $\b(3)$ and $\b(4)$. Since these three values correspond to 2, 3 and 4 for every $\b$, $x$ cannot play the role of any of 2, 3 and 4.\\

We can apply the reverse complement argument of the above to show that any value in $S_3$ cannot play the role of any of 2, 3 and 4. Therefore, values playing roles of 2, 3 and 4 must be in $S_2$, implying their positions are in $B$. However, since values whose positions are in $B$ must be expressed by Equation \ref{eqn:5.1}, these values avoid 123, 213 and 231 patterns. If the positions of the values corresponding to 2, 3 and 4 of $\b$ are in $B$, then these values must form either 123, 213 or 231 pattern. Hence, $\pi$ avoids every $\b\in\{52341,53241,52431\}$.\\
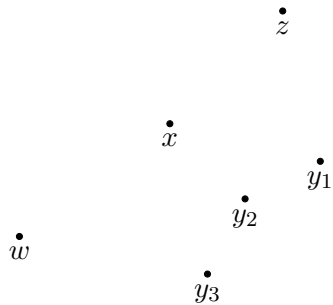
\begin{figure}[b!]
\[\begin{tikzpicture}[scale=.5]
\draw[color=black,fill=black] (-3,2) circle [radius=.08] node [below]{$w$};
\draw[color=black,fill=black] (1,5) circle [radius=.08] node [below]{$x$};
\draw[color=black,fill=black] (2,1) circle [radius=.08] node [below]{$y_3$};
\draw[color=black,fill=black] (3,3) circle [radius=.08] node [below]{$y_2$};
\draw[color=black,fill=black] (4,8) circle [radius=.08] node [below]{$z$};
\draw[color=black,fill=black] (5,4) circle [radius=.08] node [below]{$y_1$};
\end{tikzpicture}\]
\caption{Partial graph of $\pi$ to show that $35142\avd\pi$.}
\label{fig:5.30}
\end{figure}

Now, let $\b=35142$. For the same reason that the roles of 2, 3 and 4 of 52341 cannot be played by any value in $\LRmax(\pi)$ and $\RLmin(\pi)$, the role of 4 cannot be taken by any value in $\LRmax(\pi)$ and $\RLmin(\pi)$. Let $x$ be the value playing the role of 4. We first let $x\in S_1$. As before, we may assume $x$ corresponds to 1 of 231 in a 231-value chain, since the case of $x$ corresponding to 1 of 21 in a chain gives the same result. Let us refer to Figure \ref{fig:5.29} again. The only value we can assign to 5 is $y_1$, as this is the only way to have a value for 1, which is $z$. Either $u$ or $y_0$ must play the role of 3. If $u$ plays the role of 3, then we do not have a value for $2$ because $u=z+1$. On the other hand, if $y_0$ plays the role of 3, then $u$ is the only value in between $z$ and $y_0$, if $u$ and $v$ exist. Since it is located to the left of $y_0$, it cannot take the role of 2. Therefore, a value in $S_1$ cannot be assigned to 4 of $\b$.\\

Next, assume $x\in S_3$. Thus, we have Figure \ref{fig:5.30}. The only way to assign values to both 3 and 2 are by setting $w$ and $y_3$ to take roles of them respectively, but we now don't have a value corresponding to 5. Consequently, a value in $S_3$ cannot play the role of 4 of $\b$.\\

Finally, suppose $x\in S_2$. In this case, the position of a value $y$ playing the role of 5 must be in $A$ or $B$ with $\pi^{-1}(y)<\pi^{-1}(x)$. In order to assign a value to 1, $x$ must correspond to 2 of 12 in Equation \ref{eqn:5.1} or the scissor of some 231-value chain. In the former case, the paired value $z$ corresponding to 1 of 12 being assigned to 1 of $\b$. Due to the structure of values whose positions are in $B$, the position of the value $v$ corresponding to 3 must come from $C$ and the position of the value $w$ corresponding to 2 must come from $A$. However, with the sixth condition of Proposition \ref{prop:5.4}, we ensure that $v<w$, so these assignments are not possible. In case of $x$ corresponding to the scissor of some 231-value chain with the minimum value $m$ and maximum value $M$, $M$ and its descent paired value must be assigned to 3 and 1 of $\b$ respectively. Then we do not have a value corresponding to 2, because except $x$, the positions of the values in the range $[m,M]$ must be in $A$. Hence, no value in $S_3$ can play the role of $4$ of $\b$.\\

Consequently, we have $\b=35142\avd\pi$. We can show that $\pi$ does not contain $\b=42513$ by applying the reverse complement of the argument for $35142$. Thus, the remained case is $\b=351624$.\\

Let $x$ be the value playing the role of 6. Note that $x$ cannot be in $\RLmin(\pi)$, since this leaves no values to be assigned to 2 and 4. So assume $x\in S_3$. Referring to Figure \ref{fig:5.30}, the value playing the role of 2 must be $y_1$ or $y_2$ while the value for 4 must be $y_2$ or $y_3$. Therefore, Suppose $y_1$ plays the role of 2. Then the values for 3 and 5 must come from $\LRmax(\pi)$, $S_1$ or $S_2$. However, since $x$ and $y_1$ are parts of a 312-value chain, the only value that is allowed to have a position outside of $C$ is the splitter, so we don't have enough values to assign for both 3 and 5. We obtain the same result if $y_2$ plays the roles of 2. Note that, by the reverse complement argument, the value playing the role of 1 cannot be in $\LRmax(\pi)$ or $S_1$.\\

Now, suppose $x\in S_2$. Because we cannot assign any value in $\LRmax(\pi)$ or $S_1$ to 1 of $\b$, we also have to have a value $z$ in $S_2$ which plays the role of 1. Due to the structure of values with positions in $B$, the only way to assign values to both 1 and 6 is by setting $z$ and $x$ to correspond to 1 and 2 of the same 12 in Equation \ref{eqn:5.1} respectively. Then the values for 3 and 5 must come from $\LRmax(\pi)$ or $S_1$ and the values for 2 and 4 must come from $S_2$ or $\RLmin(\pi)$. This is, however, a contradiction, because four values in between of values corresponding to 1 and 2 of 12 in Equation \ref{eqn:5.1} must be in increasing order from left to right due to the sixth condition of Proposition \ref{prop:5.4}. So this concludes that $x$ cannot be in $S_2$.\\

The remained two cases are $x\in S_1$ and $x\in\LRmax(\pi)$. In either case, the value $z$ playing the role of 1 must be in $S_1$ or $\LRmax(\pi)$ as well, but as discussed previously, the role of 1 cannot be played by a value in $S_1$ or $\LRmax(\pi)$. Hence, $\b=351624\avd\pi$.\\

Consequently, $\pi$ avoids every permutation in $\{52341,53241,52431,35142,42513,351624\}$, so $\pi$ is in $\A'$. With the simplicity we proved earlier, this completes the proof.\hfill$\blacksquare$\\

We state the analogous proposition for permutations of extreme pattern 3142. Proposition \ref{prop:5.6} and \ref{prop:5.7} will be referred in Chapter \ref{chap:6}.
\begin{proposition}\label{prop:5.7}
\textit{Let $\pi$ be a permutation of extreme pattern 3142 where $|\pi|=n$. If $\pi$ has a structure described in Proposition \ref{prop:5.3} and \ref{prop:5.5}, then $\pi\in\textrm{Si}(\A')$.}
\end{proposition}

\vhhh
For the reminder of this section, we discuss how first few values of a simple permutation $\pi$ of extreme pattern 2413 can be placed. In particular, we divide the arrangement of the values in $\{x:1\leq x\leq\pi(1)\}$ into six distinct cases. We will use these six cases in the encoding rule, which we define in Chapter \ref{chap:6}. Then we apply the reverse complement property and the inverse property to establish analogous six distinct cases for each of:\begin{itemize}
\item where the values $x$ with $\pi(n)\leq x \leq n$ are located for a simple permutation $\pi$ of extreme pattern 2413,
\item what values $x$ with $1\leq \pi^{-1}(x)\leq \pi^{-1}(1)$ can be for a simple permutation $\pi$ of extreme pattern 3142 and
\item what values $x$ with $\pi^{-1}(n)\leq \pi^{-1}(x)\leq n$ can be for a simple permutation $\pi$ of extreme pattern 3142.
\end{itemize}

\begin{figure}
\[
\begin{array}{c|c|c}
& &\\
\begin{tikzpicture}[scale=.5]
\draw[color=black,fill=black] (3,2) circle [radius=.08] node [below]{$b=2$};
\draw[color=black,fill=black] (9,1) circle [radius=.08] node [below]{$a=1$};
\draw [thick,dashed] (6,0.5) -- (6,4.5);
\end{tikzpicture}\quad&\quad\begin{tikzpicture}[scale=.5]
\draw[color=black,fill=black] (2.4,3) circle [radius=.08] node [below]{$b=3$};
\draw[color=black,fill=black] (7.2,2) circle [radius=.08] node [below]{$2$};
\draw[color=black,fill=black] (9.6,1) circle [radius=.08] node [below]{$1$};
\draw [thick,dashed] (4.8,0.5) -- (4.8,4.5);
\end{tikzpicture}\quad&\quad\begin{tikzpicture}[scale=.5]
\draw[color=black,fill=black] (2.4,3) circle [radius=.08] node [below]{$b=3$};
\draw[color=black,fill=black] (4.8,2) circle [radius=.08] node [below]{$2$};
\draw[color=black,fill=black] (9.6,1) circle [radius=.08] node [below]{$1$};
\draw [thick,dashed] (7.2,0.5) -- (7.2,4.5);
\end{tikzpicture}
\\
\textrm{Case 1} & \textrm{Case 2} & \textrm{Case 3}\\
 & &\\
 \hline
\end{array}\]
\vspace{-0.465cm}
\[
\begin{array}{c|c}
& \\
\begin{tikzpicture}[scale=.5]
\draw[color=black,fill=black] (2,4) circle [radius=.08] node [below]{$b=4$};
\draw[color=black,fill=black] (4,2) circle [radius=.08] node [below]{$2$};
\draw[color=black,fill=black] (8,3) circle [radius=.08] node [below]{$3$};
\draw[color=black,fill=black] (10,1) circle [radius=.08] node [below]{$1$};
\draw [thick,dashed] (6,0.5) -- (6,4.5);
\end{tikzpicture}\qquad\qquad&\qquad\qquad\begin{tikzpicture}[scale=.5]
\draw[color=black,fill=black] (2,4) circle [radius=.08] node [below]{$b=4$};
\draw[color=black,fill=black] (4,2) circle [radius=.08] node [below]{$2$};
\draw[color=black,fill=black] (6,3) circle [radius=.08] node [below]{$3$};
\draw[color=black,fill=black] (10,1) circle [radius=.08] node [below]{$1$};
\draw [thick,dashed] (8,0.5) -- (8,4.5);
\end{tikzpicture}\\
\textrm{Case 4 and 5}\qquad\qquad&\qquad\qquad\textrm{Case 6}
\end{array}\]
\caption{Positions of values $x$ with $1\leq x\leq \pi(1)$ for $\pi$ with extreme pattern 2413.}
\label{fig:5.21}
\end{figure}
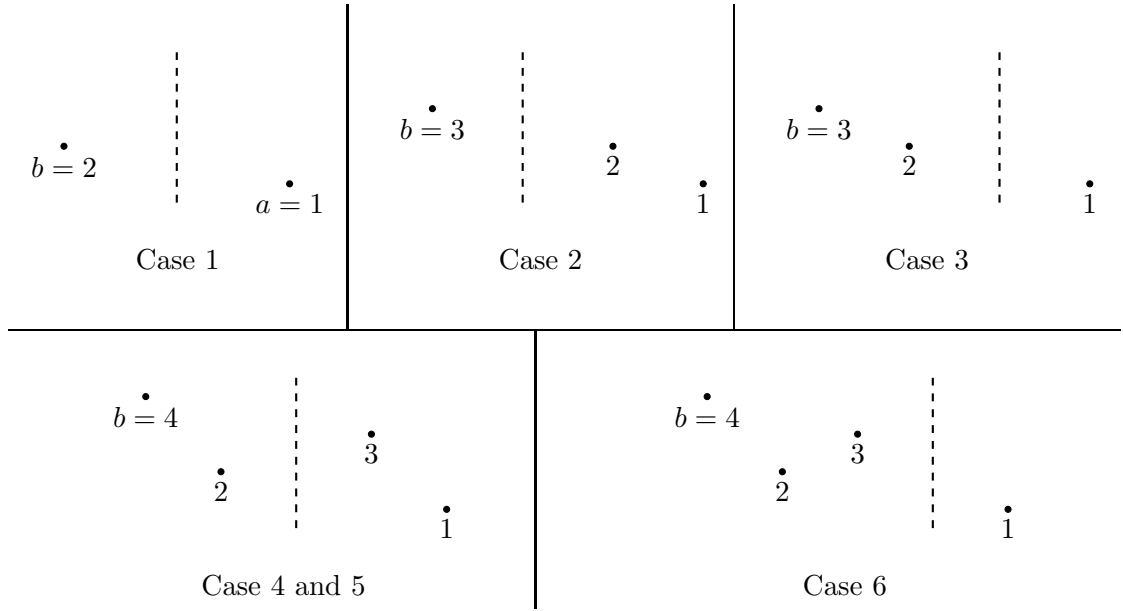

\vhhh
As before, let $b=\pi(1)$, $d=n$, $a=1$ and $c=\pi(n)$ where $\pi$ is a simple permutation in $\A'$ of extreme pattern 2413. Six distinct starting cases are shown in Figure \ref{fig:5.21}. The dashed line in each case indicates the position of $d$.\\

In every case, $a=1$ plays the role of the last 1 in Equation \ref{eqn:5.1}. The difference is what each value $x$ with $2\leq x\leq b$ is playing the role of. In Case 1, simply, $b=2$ and it is playing the role of the first 1 in segment $A$, when the values in segment $A$ are expressed as $1\oplus\cdots$. Also, in this case, the values of positions in the segment $B$ are written as $\cdots\ominus 1\ominus 1$. In this case,\[
\pi=2\ldots d\ldots 1\ldots c.\]

In Case 2, $b$ again plays the role of the first 1 when the values of positions in the segment $A$ are formed as $1\oplus\cdots$. However, the values of positions in the segment $B$ are expressed as $\cdots\ominus 12\ominus 1$, and $b$ is in between the two values taking the roles of 1 and 2 of the last 12. Thus, the value 2 is playing the role of 1 of the 12 here, which makes $b=3$. By Proposition \ref{prop:5.4}, there are possibly two more values with positions in segment $C$ which are between the values playing the roles of 1 and 2 of 12, so $\pi$ should be written as\[
\pi=3\ldots d\ldots 2\,z\,1\ldots c\]
where $z=4$, 5 or 6.\\

Case 3 may appear similar to Case 2, but the role $b$ is taking is completely different. In this case, segment $A$ looks like $\a_1\oplus\cdots$ where $\a_1$ is the first 231-value chain which starts with 231 pattern, \textit{i.e.} $\a_1=231\oplus_1\cdots$. The values $b=3$ and 2 are playing the roles of 2 and 1 of this 231 respectively. Moreover, the value 5 is playing the role of 3 of this 231, so we have\[
\pi=3\,5\,2\,\ldots d\ldots 1\ldots c.\]
If we simply have $\a_1=231$, then the position of the value 4 must be in the segment $B$. If the 231-value chain $\a_1$ continues, then $\pi^{-1}(4)$ must be in $A$.\\

In Case 4, 5 and 6, segment $A$ starts as $\a_1\oplus\cdots$ where $\a_1$ is the first 231-value chain where $\a_1=21\oplus_1\cdots$. The values $b=4$ and 2 are playing the roles of 2 and 1 of the first 21 in $\a_1$ respectively. In Case 4 and Case 5, $\a_1=21$, so the position of the value 3 is in the segment $B$. If the values of positions in $B$ is expressed as $\cdots\ominus 1\ominus 1$, then the value 3 is playing the role of the second 1 from the last. In contrast, if the values of positions in $B$ is written as $\cdots\ominus 12\ominus 1$, the value 3 is playing the role of 1 of the last 12. Hence, in Case 4 and Case 5, respectively, we have\[
\pi=4\,2\ldots d\ldots 3\,1\ldots c\qquad\textrm{and}\qquad\pi=4\,2\ldots d\ldots 3\,z\,1\ldots c\]
with $5\leq z\leq 8$.
On the other hand, in Case 6, the 231-value chain $\a_1$ continues after the initial 21, resulting in $\pi^{-1}(3)\in A$. Thus,\[
\pi=4\,2\ldots 3\ldots d\ldots 1\ldots c.\]
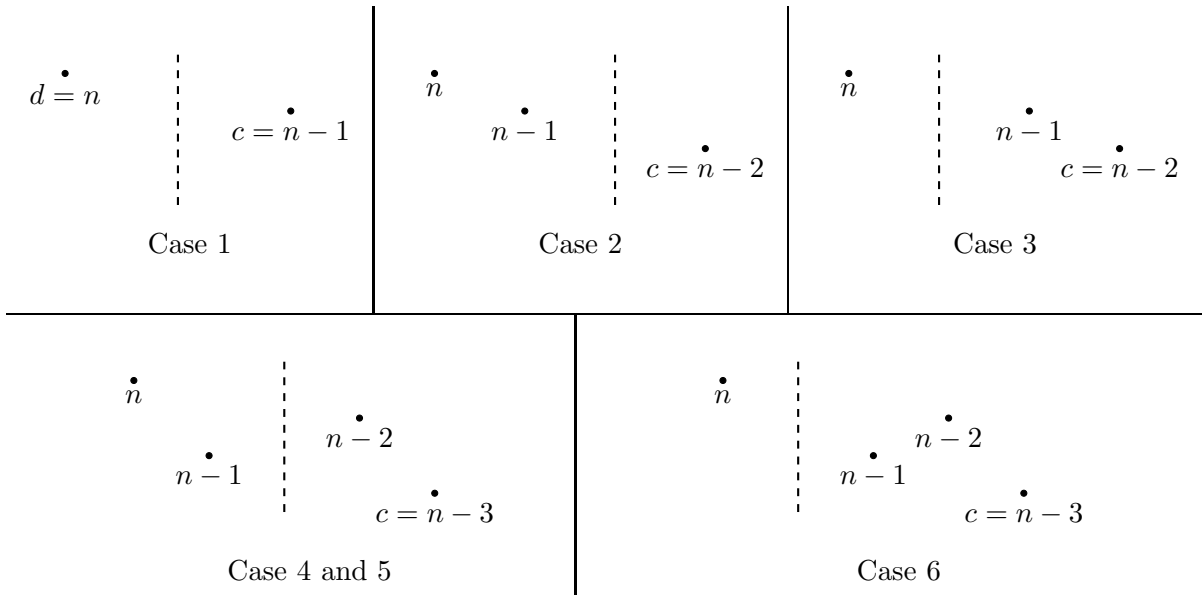
\begin{figure}[b!]
\[
\begin{array}{c|c|c}
& &\\
\begin{tikzpicture}[scale=.5]
\draw[color=black,fill=black] (3,4) circle [radius=.08] node [below]{$d=n$};
\draw[color=black,fill=black] (9,3) circle [radius=.08] node [below]{$c=n-1$};
\draw [thick,dashed] (6,0.5) -- (6,4.5);
\end{tikzpicture}&\quad\begin{tikzpicture}[scale=.5]
\draw[color=black,fill=black] (2.4,4) circle [radius=.08] node [below]{$n$};
\draw[color=black,fill=black] (4.8,3) circle [radius=.08] node [below]{$n-1$};
\draw[color=black,fill=black] (9.6,2) circle [radius=.08] node [below]{$c=n-2$};
\draw [thick,dashed] (7.2,0.5) -- (7.2,4.5);
\end{tikzpicture}&\quad\begin{tikzpicture}[scale=.5]
\draw[color=black,fill=black] (2.4,4) circle [radius=.08] node [below]{$n$};
\draw[color=black,fill=black] (7.2,3) circle [radius=.08] node [below]{$n-1$};
\draw[color=black,fill=black] (9.6,2) circle [radius=.08] node [below]{$c=n-2$};
\draw [thick,dashed] (4.8,0.5) -- (4.8,4.5);
\end{tikzpicture}
\\
\textrm{Case 1} & \textrm{Case 2} & \textrm{Case 3}\\
 & &\\
 \hline
\end{array}\]
\vspace{-0.465cm}
\[
\begin{array}{c|c}
& \\
\begin{tikzpicture}[scale=.5]
\draw[color=black,fill=black] (2,4) circle [radius=.08] node [below]{$n$};
\draw[color=black,fill=black] (4,2) circle [radius=.08] node [below]{$n-1$};
\draw[color=black,fill=black] (8,3) circle [radius=.08] node [below]{$n-2$};
\draw[color=black,fill=black] (10,1) circle [radius=.08] node [below]{$c=n-3$};
\draw [thick,dashed] (6,0.5) -- (6,4.5);
\end{tikzpicture}\qquad\qquad&\qquad\qquad\begin{tikzpicture}[scale=.5]
\draw[color=black,fill=black] (2,4) circle [radius=.08] node [below]{$n$};
\draw[color=black,fill=black] (6,2) circle [radius=.08] node [below]{$n-1$};
\draw[color=black,fill=black] (8,3) circle [radius=.08] node [below]{$n-2$};
\draw[color=black,fill=black] (10,1) circle [radius=.08] node [below]{$c=n-3$};
\draw [thick,dashed] (4,0.5) -- (4,4.5);
\end{tikzpicture}\\
\textrm{Case 4 and 5}\qquad\qquad&\qquad\qquad\textrm{Case 6}
\end{array}\]
\caption{Positions of values $x$ with $\pi(n)\leq x \leq n$ for $\pi$ with extreme pattern 2413.}
\label{fig:5.22}
\end{figure}

Based on the structure we discussed in Proposition \ref{prop:5.2} and \ref{prop:5.4}, these six cases the only possible ways that first several values of simple permutations of extreme pattern 2413 can be placed. By applying the reverse complement, we obtain the possible behaviors of points with the values $x$ with $\pi(n)\leq x \leq n$, where $\pi$ is a simple permutation of extreme pattern 2413. Analogous cases are shown in Figure \ref{fig:5.22}. Also, for a simple permutation $\pi$ of extreme pattern 3142, the first and the last few values of $\pi$ are obtained by the inverse property, shown in Figure \ref{fig:5.23} and \ref{fig:5.24} respectively. In Figure \ref{fig:5.22}, \ref{fig:5.23} and \ref{fig:5.24}, dashed lines indicate $\pi^{-1}(a)$, $c$ and $b$ respectively.\\
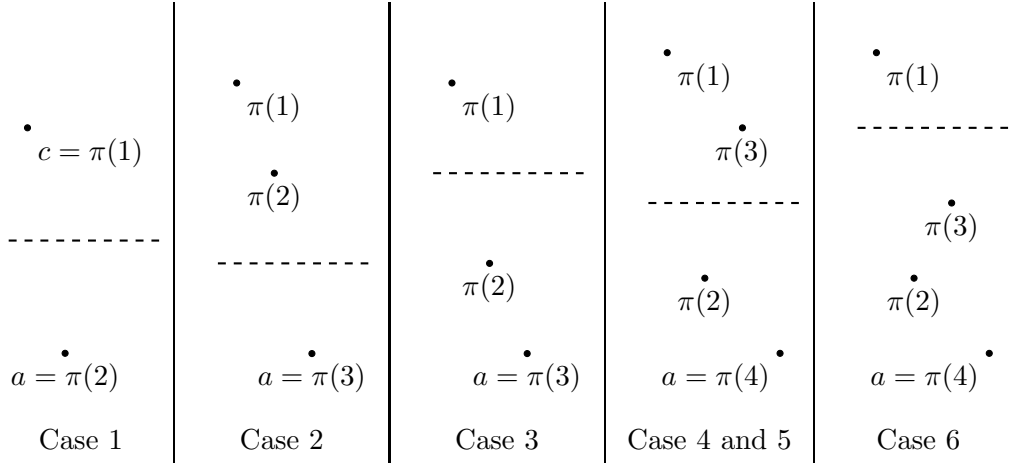
\begin{figure}[t!]
\[
\begin{array}{c|c|c|c|c}
& & & &\\
\begin{tikzpicture}[scale=.5]
\draw[color=black,fill=black] (1,9) circle [radius=.08] node [below right]{$c=\pi(1)$};
\draw[color=black,fill=black] (2,3) circle [radius=.08] node [below]{$a=\pi(2)$};
\draw [thick,dashed] (0.5,6) -- (4.5,6);
\end{tikzpicture}\quad&\quad\begin{tikzpicture}[scale=.5]
\draw[color=black,fill=black] (1,9.6) circle [radius=.08] node [below right]{$\pi(1)$};
\draw[color=black,fill=black] (2,7.2) circle [radius=.08] node [below]{$\pi(2)$};
\draw[color=black,fill=black] (3,2.4) circle [radius=.08] node [below]{$a=\pi(3)$};
\draw [thick,dashed] (0.5,4.8) -- (4.5,4.8);
\end{tikzpicture}\quad&\quad\begin{tikzpicture}[scale=.5]
\draw[color=black,fill=black] (1,9.6) circle [radius=.08] node [below right]{$\pi(1)$};
\draw[color=black,fill=black] (2,4.8) circle [radius=.08] node [below]{$\pi(2)$};
\draw[color=black,fill=black] (3,2.4) circle [radius=.08] node [below]{$a=\pi(3)$};
\draw [thick,dashed] (0.5,7.2) -- (4.5,7.2);
\end{tikzpicture}\quad&\quad\begin{tikzpicture}[scale=.5]
\draw[color=black,fill=black] (1,10) circle [radius=.08] node [below right]{$\pi(1)$};
\draw[color=black,fill=black] (2,4) circle [radius=.08] node [below]{$\pi(2)$};
\draw[color=black,fill=black] (3,8) circle [radius=.08] node [below]{$\pi(3)$};
\draw[color=black,fill=black] (4,2) circle [radius=.08] node [below left]{$a=\pi(4)$};
\draw [thick,dashed] (0.5,6) -- (4.5,6);
\end{tikzpicture}\quad&\quad\begin{tikzpicture}[scale=.5]
\draw[color=black,fill=black] (1,10) circle [radius=.08] node [below right]{$\pi(1)$};
\draw[color=black,fill=black] (2,4) circle [radius=.08] node [below]{$\pi(2)$};
\draw[color=black,fill=black] (3,6) circle [radius=.08] node [below]{$\pi(3)$};
\draw[color=black,fill=black] (4,2) circle [radius=.08] node [below left]{$a=\pi(4)$};
\draw [thick,dashed] (0.5,8) -- (4.5,8);
\end{tikzpicture}\quad\\
\textrm{Case 1} & \textrm{Case 2} & \textrm{Case 3} & \textrm{Case 4 and 5} & \textrm{Case 6}
\end{array}\]
\caption{Values of positions $s$ with $1\leq s \leq \pi^{-1}(1)$ for $\pi$ with extreme pattern 3142.}
\label{fig:5.23}
\end{figure}
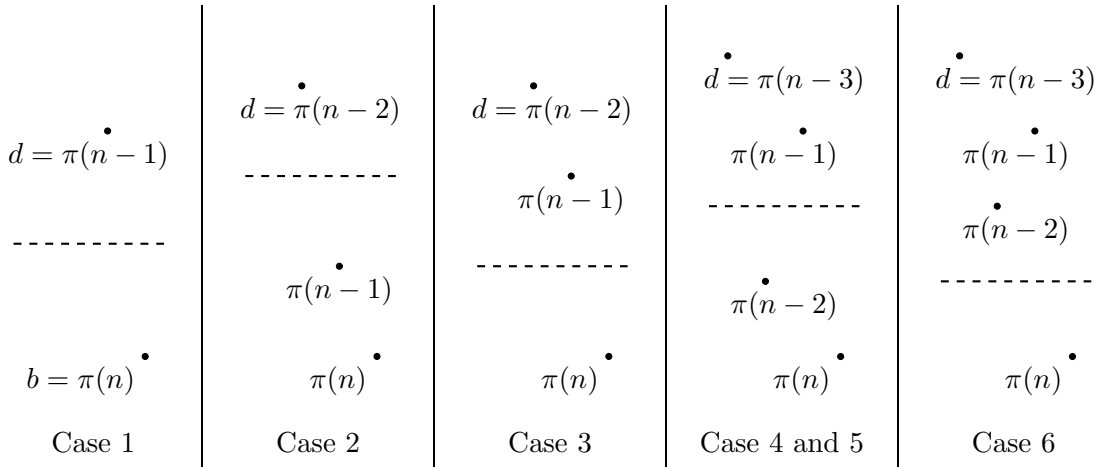
\begin{figure}[t!]
\[
\begin{array}{c|c|c|c|c}
& & & &\\
\begin{tikzpicture}[scale=.5]
\draw[color=black,fill=black] (3,9) circle [radius=.08];
\node[below] at (2.5,9) {$d=\pi(n-1)$};
\draw[color=black,fill=black] (4,3) circle [radius=.08] node [below left]{$b=\pi(n)$};
\draw [thick,dashed] (0.5,6) -- (4.5,6);
\end{tikzpicture}\,\,\,&\,\,\,\begin{tikzpicture}[scale=.5]
\draw[color=black,fill=black] (2,9.6) circle [radius=.08];
\node[below] at (2.5,9.6) {$d=\pi(n-2)$};
\draw[color=black,fill=black] (3,4.8) circle [radius=.08] node [below]{$\pi(n-1)$};
\draw[color=black,fill=black] (4,2.4) circle [radius=.08] node [below left]{$\pi(n)$};
\draw [thick,dashed] (0.5,7.2) -- (4.5,7.2);
\end{tikzpicture}\,\,\,&\,\,\,\begin{tikzpicture}[scale=.5]
\draw[color=black,fill=black] (2,9.6) circle [radius=.08];
\node[below] at (2.5,9.6) {$d=\pi(n-2)$};
\draw[color=black,fill=black] (3,7.2) circle [radius=.08] node [below]{$\pi(n-1)$};
\draw[color=black,fill=black] (4,2.4) circle [radius=.08] node [below left]{$\pi(n)$};
\draw [thick,dashed] (0.5,4.8) -- (4.5,4.8);
\end{tikzpicture}\,\,\,&\,\,\,\begin{tikzpicture}[scale=.5]
\draw[color=black,fill=black] (1,10) circle [radius=.08];
\node[below] at (2.5,10) {$d=\pi(n-3)$};
\draw[color=black,fill=black] (2,4) circle [radius=.08];
\node[below] at (2.5,4) {$\pi(n-2)$};
\draw[color=black,fill=black] (3,8) circle [radius=.08];
\node[below] at (2.5,8) {$\pi(n-1)$};
\draw[color=black,fill=black] (4,2) circle [radius=.08] node [below left]{$\pi(n)$};
\draw [thick,dashed] (0.5,6) -- (4.5,6);
\end{tikzpicture}\,\,\,&\,\,\,\begin{tikzpicture}[scale=.5]
\draw[color=black,fill=black] (1,10) circle [radius=.08];
\node[below] at (2.5,10) {$d=\pi(n-3)$};
\draw[color=black,fill=black] (2,6) circle [radius=.08];
\node[below] at (2.5,6) {$\pi(n-2)$};
\draw[color=black,fill=black] (3,8) circle [radius=.08];
\node[below] at (2.5,8) {$\pi(n-1)$};
\draw[color=black,fill=black] (4,2) circle [radius=.08] node [below left]{$\pi(n)$};
\draw [thick,dashed] (0.5,4) -- (4.5,4);
\end{tikzpicture}\,\,\,\\
\textrm{Case 1} & \textrm{Case 2} & \textrm{Case 3} & \textrm{Case 4 and 5} & \textrm{Case 6}
\end{array}\]
\caption{Values of positions $s$ with $1\leq \pi^{-1}(n)\leq n$ for $\pi$ with extreme pattern 3142.}
\label{fig:5.24}
\end{figure}
\section{General simple permutations in $\A'$}
In the previous section, we looked at the structures of simple permutations of extreme pattern 2413 and the structures of simple permutations of extreme pattern 3142 in great detail. We are ready to discuss the structure of general simple permutations in $\A'$. Recall that in $\A$, we established Theorem \ref{thm:4.1}, which states that the half of simple permutations can be written as alternated NW and SE glue sums of simple permutations of extreme pattern 2413 and 3142, starting with one of extreme pattern 2413. In this section, we extend Theorem \ref{thm:4.1} to describe the structure of the half of simple permutations in $\A'$. In order to state the theorem extending Theorem \ref{thm:4.1}, we need to define several more glue sums.
\subsection{Glue sums and the structure theorem}
We carry the same definitions for type 1-0 NW glue sum, type 1-0 SE glue sum, type 1-1 NW glue sum and type 1-1 SE glue sum for simple permutations in $\A'$. Recall that, in the definition of $\s\nwsum_1^0\tau$ and $\s\nwsum_1^1\tau$, $\s$ ($|\s|=m$) and $\tau$ ($|\tau|=n$) must satisfy the conditions that $\s^{-1}(m)\leq m-2$, $\s(m)=m-1$, $\tau(1)\geq 3$ and $\tau(2)=1$. In the language of the previous section, $\s$ must have the positions of the last two values $m-1$ and $m$ as Case 1 of Figure \ref{fig:5.22} and $\tau$ must have the values of the first two positions $\tau(1)$ and $\tau(2)$ as Case 1 of Figure \ref{fig:5.23}.\\

Similarly, for $\s\sesum_1^0\tau$ and $\s\sesum_1^1\tau$, two simple permutations $\s$ ($|\s|=m$) and $\tau$ ($|\tau|=n$) must satisfy the conditions of $\s(m)\leq m-2$, $m=\s(m-1)$, $\tau^{-1}(1)\geq 3$ and $2=\tau(1)$. Thus, $\s$ must have the values of the last two positions $\s(m-1)$ and $\s(m)$ as Case 1 of Figure \ref{fig:5.24} and $\tau$ must have the positions of the first two values $1$ and $2$ as Case 1 of Figure \ref{fig:5.21}.\\

For every glue sums that we define to combine simple permutations $\s$ and $\tau$ in this section, $\s$ and $\tau$ have to satisfy specific conditions similar to the ones for type 1-0 and 1-1 glue sums. In particular, $\s$, the left permutation of the sum, must follow a particular structure described in Figure \ref{fig:5.22} for NW glue sums and \ref{fig:5.24} for SE glue sums. As noted, Type 1-0 and 1-1 are for Case 1 structure. Type 2-0 and 2-1 are for Case 2, type 3-0 is for Case 3, and type 4-0 is for Case 4 and Case 5. According to which structure $\s$ has, $\tau$ also has to have a certain structure. The second number indicates whether the last value of $\s$ is left as a copy or not in the process of glue summing. If the second number is 1, the glue sum leaves a copy of the last point of $\s$ when it is merged into the last point of $\tau$, just as it was explained with type 1-1 at the end of the proof of Proposition \ref{prop:4.4}.\\

We now introduce all types of NW glue sums in Table \ref{tab:5.1}. Note that every simple permutation of extreme pattern 2413 satisfies the conditions of $\s$, which is the left side of each sum, and similarly, every simple permutation of extreme pattern 3142 satisfies the conditions of $\tau$, the right side of each sum. Let $\s$ and $\tau$ be simple permutations in $\A'$ of length $m$ and $n$ respectively, $i=\s^{-1}(m)$ and $j=\tau(1)$. For each NW glue sum, it is still proper to think the last value of $\s$ and the least value of $\tau$ are merged into the last value of $\tau$ and the least value of $\s$ respectively. Just like type 1-0 and 1-1, each sum combines the point of $\s$ whose value is $m$ and the point of $\tau$ whose value is $j$ into one point. This point keeps the position $i$, and its value is $m$ shifted upward by certain amount. In addition, for type 2-0, 2-1 and 3-0, the point whose value is $m-1$ is also shifted up by the same amount $m$ is shifted. Likewise, for type 4-0, the point whose value is $m-2$ is shifted up together with $m$ by the same amount.\\

The only differences between type 1-0 and type 2-0 are the conditions of $\s$ and whether $m-1$ is shifted up or not. Same can be said for type 1-1 and 2-1.\\
\afterpage{
\begin{landscape}
\begin{table}[p]
\begin{footnotesize}
\begin{center}
\begin{tabular}{ |c|L|L|c| }
\hline
Type & Conditions of $\s$ & Conditions of $\tau$ & Notation and Definition \\
\hline
\multirow{2}{*}{1-0} & \multirow{4}{*}{\parbox{3.5cm}{\centering $i\leq m-2$, $\s(m)=m-1$.\\(Case 1 of Figure \ref{fig:5.22})}} & \multirow{10}{*}{\parbox{3.5cm}{\centering $j\geq 3$, $\tau(2)=1$.\\(Case 1 of Figure \ref{fig:5.23})}} & \multirow{2}{*}{\parbox{12cm}{\centering $\s\nwssum_1^0\tau=\s'(1)\s'(2)\cdots\s'(m-1)\tau'(3)\tau'(4)\cdots\tau'(n)$ where $\s'(i)=m+(j-3)$ and $\s'(k)=\s(k)$ for $k\neq i$, and $\tau'(k)=\tau(k)+(m-3)$ for $k$ with $3\leq k\leq n$.}}\\
 & & & \\
\cline{1-1} \cline{4-1}
\multirow{2}{*}{1-1} &  &  & \multirow{2}{*}{\parbox{12cm}{\centering $\s\nwssum_1^1\tau=\s'(1)\s'(2)\cdots\s'(m)\tau'(3)\tau'(4)\cdots\tau'(n)$ where $\s'(i)=m+(j-2)$ and $\s'(k)=\s(k)$ for $k\neq i$, and $\tau'(k)=\tau(k)+(m-2)$ for $k$ with $3\leq k\leq n$.}}\\
 & & & \\
\cline{1-2} \cline{4-1}
\multirow{3}{*}{2-0} & \multirow{6}{*}{\parbox{3.5cm}{\centering $i\leq m-4$, $\s(i+2)=m-1$, $\s(m)=m-2$.\\(Case 2 of Figure \ref{fig:5.22})}} & & \multirow{3}{*}{\parbox{12cm}{\centering $\s\nwssum_2^0\tau=\s'(1)\s'(2)\cdots\s'(m-1)\tau'(3)\tau'(4)\cdots\tau'(n)$ where $\s'(i)=m+(j-3)$, $\s'(i+2)=(m-1)+(j-3)$ and $\s'(k)=\s(k)$ for $k\neq i,i+2$, and\\$\tau'(k)=\tau(k)+(m-3)$ for $k$ with $3\leq k\leq n$.}}\\
 & & & \\
 & & & \\
\cline{1-1} \cline{4-1}
\multirow{3}{*}{2-1} & & & \multirow{3}{*}{\parbox{13cm}{\centering $\s\nwssum_2^1\tau=\s'(1)\s'(2)\cdots\s'(m)\tau'(3)\tau'(4)\cdots\tau'(n)$ where $\s'(i)=m+(j-2)$, $\s'(i+2)=(m-1)+(j-2)$ and $\s'(k)=\s(k)$ for $k\neq i,i+2$, and\\$\tau'(k)=\tau(k)+(m-2)$ for $k$ with $3\leq k\leq n$.}}\\
 & & & \\
 & & & \\
\hline
\multirow{7}{*}{3-0} & \multirow{6}{*}{\parbox{3.5cm}{\centering $i\leq m-5$, $\s(m-2)=m-1$, $\s(m)=m-2$ and $\s(m)$ is a part of the last 312-value chain $\a$.\\(Case 3 of Figure \ref{fig:5.22})}} & \multirow{6}{*}{\parbox{3.5cm}{\centering $j\geq 6$, $\tau(2)=3$, $\tau(3)=1$ and $1$ is a part of the least 312-position chain $\b$.\\(Case 3 of Figure \ref{fig:5.23})}} & \multirow{7}{*}{\parbox{13cm}{\centering $\s\nwssum_3^0\tau=\s'(1)\s'(2)\cdots\s'(m-1)\tau'(\l+3)\tau'(\l+4)\cdots\tau'(n)$\\where $\l$ is the number of values in the 312-position chain $\b$,\\ $\s'(i)=m+(j-\l-3)$, $\s'(m-2)=(m-1)+(h-\l-2)$ where $h=\tau(\l+1)$ and\\$\s'(k)=\s(k)$ for $k\neq i,m-2$, and $\tau'(k)=\tau(k)+(m-\l-3)$ for $k$ with $\l+3\leq k\leq n$.}}\\
 & & & \\
 & & & \\
 & & & \\
 & & & \\
 & & & \\
\cline{2-3}
 & \multicolumn{2}{c|}{\centering Additionally, $\a$ must be similar to $\b$.} & \\ \hline
\multirow{6}{*}{4-0} & \multirow{6}{*}{\parbox{3.5cm}{\centering $i\leq m-4$, $\s(m-1)=m-1$, $\s(m)=m-3$.\\(Case 4 or 5 of\\Figure \ref{fig:5.22})}} & \multirow{6}{*}{\parbox{3.5cm}{\centering $j\geq 5$, $\tau(2)=j-2$, $\tau(3)=1$.\\(Case 2 of Figure \ref{fig:5.23})}} & \multirow{6}{*}{\parbox{12cm}{\centering $\s\nwssum_4^0\tau=\s'(1)\s'(2)\cdots\s'(m-2)\tau'(4)\tau'(5)\cdots\tau'(n)$ where $\s'(i)=m+(j-3)$, $\s'(s)=(m-2)+(j-3)$ ($s$ position of the value $m-2$) and $\s'(k)=\s(k)$ for $k\neq i,s$, and $\tau'(k)=\tau(k)+(m-3)$ for $k$ with $4\leq k\leq n$.}}\\
 & & & \\
 & & & \\
 & & & \\
 & & & \\
 & & & \\
\hline
\end{tabular}
\end{center}
\end{footnotesize}
\caption{Definitions of NW glue sums ($i=\s^{-1}(m)$ and $j=\tau(1)$).}
\label{tab:5.1}
\end{table}
\end{landscape}
}

Before we discuss type 3-0 NW glue sum, we define another terminology. Let $\a$ and $\b$ be a 312-value chain and a 312-position chain respectively. We say $\a$ and $\b$ are \textit{similar} if the flattening of $21\oplus_1\a$ and the flattening of $\b\oplus^1 21$ are equal to each other. For instance, $\a=21\oplus_1 312$ is similar to $\b=312\oplus^1 21$ since $21\oplus_1\a=\b\oplus^1 21=3152746$. Likewise, we say a 231-value chain $\a$ and a 231-position chain $\b$ are \textit{similar} if the flattening $21\oplus_1\a$ and the flattening of $\b\oplus^1 21$ are equal to each other.\\

As indicated in Table \ref{tab:5.1}, for type 3-0 NW glue sum, $\s$ and $\tau$ must have a 312-value chain with $\l$ values involving $\s(m)$ and a 312-position chain with $\l$ values involving $1$ respectively that are similar to each other. For now, let $\s$ and $\tau$ be simple permutations of extreme pattern 2413 and 3142 respectively. We visualize $\s\oplus_3^0\tau$ in Figure \ref{fig:5.25}. In $\s$, we must have a scissor $s$ as indicated in Figure \ref{fig:5.25} because of the fifth condition of Proposition \ref{prop:5.4}. Similarly, we have a scissor $h$ as shown in Figure \ref{fig:5.25}. Note that the position of $h$ is $\l+1$, so $h=\tau(\l+1)$. The 312-value chain of $\s$ together with 1 and $s$ has the same pattern as the pattern of 312-position chain of $\tau$ with $h$ and $\tau(n)$. Type 3-0 NW glue sum identifies each boxed point of $\s$ shown in Figure \ref{fig:5.25}, and combines them into boxed points of $\tau$ from left to right. In this process, $\l+2$ values of $\tau$ that are less than $j$ are combined, so the amount we shift $\s(i)=m$ upward is $(j-1)-(\l+2)=j-\l-3$. Similarly, $\l+1$ values of $\tau$ below $h$ are combined so the amount we shift $\s(m-2)=m-1$ upward is $(h-1)-(\l+1)=h-\l-2$. Since $\tau(k)$ $(1\leq k\leq \l+2)$ are combined into corresponding values of $\s$, the glue sum starts referring to values of $\tau$ from $\tau(\l+3)$. We need to shift $\tau(k)$ $(\l+3\leq k\leq n)$ upward by $m-(\l+3)$ because this is the number of values in $\s$ that are not merged with values of $\tau$. As a result, type 3-0 NW glue sum creates a chain containing both a 312-value chain and a 231-position chain, as indicated in the box in Figure \ref{fig:5.27}. Note that this chain has both a scissor $s$ with respect to its 312-value chain aspect and a scissor $(m-2)+(h-\l-2)$ with respect to its 231-value chain aspect.\\
\begin{figure}[p]
\[
\begin{array}{c}\begin{tikzpicture}[scale=.4,rotate=180]
\draw[color=black,fill=white] (10.2,5.2) rectangle (10.8,5.8);\node[below] at (10.5,5.8) {$s$};
\draw[color=black,fill=white] (8.2,10.2) rectangle (8.8,10.8);\node[below] at (8.5,10.8) {$1$};
\draw[color=black,fill=white] (2.2,-0.8) rectangle (2.8,-0.2);\node[left] at (2.8,-0.5) {$m-1$};
\draw[color=black,fill=white] (1.3,1.3) rectangle (0.7,0.7);\node[right] at (0.7,1) {$m-2$};
\draw[color=black,fill=white] (2.3,2.3) rectangle (1.7,1.7);\draw[color=black,fill=white] (5.3,5.3) rectangle (4.7,4.7);\draw[color=black,fill=white] (6.3,6.3) rectangle (5.7,5.7);\draw[color=black,fill=white] (6.8,3.8) rectangle (6.2,3.2);
\draw[color=black,fill=black] (1,1) circle [radius=.1];
\draw[color=black,fill=black] (2,2) circle [radius=.1];
\draw[color=black,fill=black] (5,5) circle [radius=.1];
\draw[color=black,fill=black] (6,6) circle [radius=.1];
\draw[color=black,fill=black] (2.5,-0.5) circle [radius=.1];
\draw[color=black,fill=black] (6.5,3.5) circle [radius=.1];
\draw [dotted,ultra thick] (3.5,2) -- (4.5,3);
\draw[color=black,fill=black] (10.5,5.5) circle [radius=.1];
\draw[color=black,fill=black] (8.5,10.5) circle [radius=.1];
\draw[color=black,fill=black] (13.5,-2) circle [radius=.1] node [below]{$m$};
\draw [dotted,ultra thick] (15,5.5) -- (16.4,5.5);
\draw[arrows=->,thick,black] (13.5,-2.2) -- (13.5,-5.2);
\draw[arrows=->,thick,black] (2.5,-0.8) -- (2.5,-2.2);
\draw[thick,dashed] (-2.4,12) -- (-2.4,-8);
\end{tikzpicture}\hspace{-6.15cm}\raisebox{1.4cm}{\begin{tikzpicture}[scale=.4]
\draw[color=black,fill=white] (10.8,8.8) rectangle (10.2,8.2);\node[below] at (10.5,8.2) {$\tau(n)$};
\draw[color=black,fill=white] (5.8,10.8) rectangle (5.2,10.2);\node[below] at (5.5,10.2) {$h$};
\draw[color=black,fill=white] (-0.2,2.8) rectangle (-0.8,2.2);\node[below] at (-0.5,2.2) {$3$};
\draw[color=black,fill=white] (1.3,1.3) rectangle (0.7,0.7);\node[below] at (1,0.7) {$1$};
\draw[color=black,fill=white] (2.3,2.3) rectangle (1.7,1.7);\draw[color=black,fill=white] (5.3,5.3) rectangle (4.7,4.7);\draw[color=black,fill=white] (6.3,6.3) rectangle (5.7,5.7);\draw[color=black,fill=white] (3.8,6.8) rectangle (3.2,6.2);
\draw[color=black,fill=black] (1,1) circle [radius=.1];
\draw[color=black,fill=black] (2,2) circle [radius=.1];
\draw[color=black,fill=black] (5,5) circle [radius=.1];
\draw[color=black,fill=black] (6,6) circle [radius=.1];
\draw[color=black,fill=black] (-0.5,2.5) circle [radius=.1];
\draw[color=black,fill=black] (3.5,6.5) circle [radius=.1];
\draw [dotted,ultra thick] (2,3.5) -- (3,4.5);
\draw[color=black,fill=black] (5.5,10.5) circle [radius=.1];
\draw[color=black,fill=black] (10.5,8.5) circle [radius=.1];
\draw[color=black,fill=black] (-2,13.5) circle [radius=.1] node [below]{$j$};
\draw [dotted,ultra thick] (5.5,15) -- (5.5,16.4);
\draw[arrows=->,thick,black] (-2.,13.5) -- (-19,13.5);
\draw[arrows=->,thick,black] (5.2,10.5) -- (-8,10.5);
\end{tikzpicture}}\\
\s\hspace{3.2in}\tau\\
\\
\downarrow\\
\\
\begin{tikzpicture}[scale=.4,rotate=180]
\draw[color=black,fill=white] (1.5,1) rectangle (7,6.5);
\draw[color=black,fill=black] (2,2) circle [radius=.1];
\draw[color=black,fill=black] (5,5) circle [radius=.1];
\draw[color=black,fill=black] (6,6) circle [radius=.1];
\draw[color=black,fill=black] (4.5,1.5) circle [radius=.1];
\draw[color=black,fill=black] (2.5,-2.6) circle [radius=.1] node [left]{$(m-1)+(h-\l-2)$};
\draw[color=black,fill=black] (6.5,3.5) circle [radius=.1];
\draw [dotted,ultra thick] (4.5,3) -- (5.1,3.6);
\draw[color=black,fill=black] (10.5,5.5) circle [radius=.1] node [below]{$s$};
\draw[color=black,fill=black] (8.5,10.5) circle [radius=.1] node [below]{$1$};
\draw[color=black,fill=black] (13.5,-5.4) circle [radius=.1] node [below]{$m+(j-\l-3)$};
\draw [dotted,ultra thick] (15,5.5) -- (16.4,5.5);
\draw[color=black,fill=black] (-3.5,-0.6) circle [radius=.1] node [below]{$\tau(n)+(m-(\l+3))$};
\draw [dotted,ultra thick] (2.5,-7.1) -- (2.5,-8.5);
\end{tikzpicture}
\end{array}
\]
\caption{Illustration of $\s\nwsum_3^0\tau$.}
\label{fig:5.25}
\end{figure}
\afterpage{
\begin{landscape}
\begin{table}[p]
\begin{footnotesize}
\begin{center}
\begin{tabular}{ |c|L|L|c| }
\hline
Type & Conditions of $\s$ & Conditions of $\tau$ & Notation and Definition \\
\hline
\multirow{2}{*}{1-0} & \multirow{4}{*}{\parbox{3.5cm}{\centering $i\leq m-2$, $m=\s(m-1)$.\\(Case 1 of Figure \ref{fig:5.24})}} & \multirow{10}{*}{\parbox{3.5cm}{\centering $j\geq 3$, $2=\tau(1)$.\\(Case 1 of Figure \ref{fig:5.21})}} & \multirow{2}{*}{\parbox{12cm}{\centering $\s\sessum_1^0\tau=\s(1)\s(2)\cdots\s(m-2)\tau'(2)\tau'(3)\cdots\tau'(n)$ where $\tau'(j)=i$ and $\tau'(k)=\tau(k)+(m-3)$ for $k\neq j$.}}\\
 & & & \\
\cline{1-1} \cline{4-1}
\multirow{2}{*}{1-1} &  &  & \multirow{2}{*}{\parbox{12cm}{\centering $\s\sessum_1^1\tau=\s(1)\s(2)\cdots\s(m-1)\tau'(2)\tau'(3)\cdots\tau'(n)$ where $\tau'(j)=i$ and $\tau'(k)=\tau(k)+(m-2)$ for $k\neq j$.}}\\
 & & & \\
\cline{1-2} \cline{4-1}
\multirow{3}{*}{2-0} & \multirow{6}{*}{\parbox{3.5cm}{\centering $i\leq m-4$, $i+2=\s(m-1)$, $m=\s(m-2)$.\\(Case 2 of Figure \ref{fig:5.24})}} & & \multirow{3}{*}{\parbox{12cm}{\centering $\s\sessum_2^0\tau=\s(1)\s(2)\cdots\s(m-3)\tau'(2)\tau'(3)\cdots\tau'(n+1)$ where $\tau'(j)=i+2$, $\tau'(j+1)=i$, and $\tau'(k)=\tau(k)+(m-3)$ for $k\leq j-1$ and $\tau'(k)=\tau(k-1)+(m-3)$ for $j+2\leq k\leq n$.}}\\
 & & & \\
 & & & \\
\cline{1-1} \cline{4-1}
\multirow{3}{*}{2-1} & & & \multirow{3}{*}{\parbox{13cm}{\centering $\s\sessum_2^1\tau=\s(1)\s(2)\cdots\s(m-3)\s(m-2)\tau'(2)\tau'(3)\cdots\tau'(n+1)$ where $\tau'(j)=i+2$, $\tau'(j+1)=i$, and $\tau'(k)=\tau(k)+(m-3)$ for $k\leq j-1$\\and $\tau'(k)=\tau(k-1)+(m-3)$ for $j+2\leq k\leq n$.}}\\
 & & & \\
 & & & \\
\hline
\multirow{7}{*}{3-0} & \multirow{6}{*}{\parbox{3.5cm}{\centering $i\leq m-5$, $m-2=\s(m-1)$, $m=\s(m-2)$ and $m$ is a part of the greatest 231-position chain $\a$.\\(Case 3 of Figure \ref{fig:5.24})}} & \multirow{6}{*}{\parbox{3.5cm}{\centering $j\geq 6$, $2=\tau(3)$, $3=\tau(1)$ and $\tau(1)$ is a part of the first 231-value chain $\b$.\\(Case 3 of Figure \ref{fig:5.21})}} & \multirow{7}{*}{\parbox{13cm}{\centering $\s\sessum_3^0\tau=\s(1)\s(2)\cdots\s(m-3)\tau'(\l+1)\tau'(\l+2)\cdots\tau'(n)$\\where $\l$ is the number of values in the 231-value chain $\b$,\\and $\tau(j-1)=\s(m-1)$, $\tau(j)=i$, and $\tau'(k)=\tau(k)+(m-\l-3)$\\for $k$ with $\l+1\leq k\leq n$ and $k\neq j-1,j$.}}\\
 & & & \\
 & & & \\
 & & & \\
 & & & \\
 & & & \\
\cline{2-3}
 & \multicolumn{2}{c|}{\centering Additionally, $\a$ must be similar to $\b$.} & \\ \hline
\multirow{6}{*}{4-0} & \multirow{6}{*}{\parbox{3.5cm}{\centering $i\leq m-4$, $m-1=\s(m-1)$, $m=\s(m-3)$.\\(Case 4 or 5 of\\Figure \ref{fig:5.24})}} & \multirow{6}{*}{\parbox{3.5cm}{\centering $j\geq 5$, $2=\tau(j-2)$, $3=\tau(1)$.\\(Case 2 of Figure \ref{fig:5.21})}} & \multirow{6}{*}{\parbox{12cm}{\centering $\s\sessum_4^0\tau=\s(1)\s(2)\cdots\s(m-4)\tau'(2)\tau'(3)\cdots\tau'(n)$ where $\tau'(j-2)=\s(m-2)$, $\tau'(j)=i$, and $\tau'(k)=\tau(k)+(m-3)$ for $k$ with $2\leq k\leq n$ and $k\neq j-2,j$.}}\\
 & & & \\
 & & & \\
 & & & \\
 & & & \\
 & & & \\
\hline
\end{tabular}
\end{center}
\end{footnotesize}
\caption{Definitions of SE glue sums ($i=\s(m)$ and $j=\tau^{-1}(1)$).}
\label{tab:5.2}
\end{table}
\end{landscape}
}

Type 4-0 is rather similar to type 1-0, 1-1, 2-0 and 2-1. It is appropriate to visualize that $\s(m-1)$ is merged into either $(4,\tau(4))$ or $(5,\tau(5))$, depending on whether $\tau(4)<j$ or not. If $\tau(4)<j$, then $(m-1,m-1)$ is merged into $(4,\tau(4))$, otherwise, $(m-1,m-1)$ is merged into $(5,\tau(5))$. In addition to the combined point of $(i,m)$ and $(1,j)$, type 4-0 also combines $(s,m-2)$ and $(2,\tau(2))$. The resulting point keeps the position $s$, and the value is $m-2+(j-3)$. Since three starting values of $\tau$ are merged into certain points of $\s$, values of $\tau$ are reserved from the fourth position instead of the third one. Note that $s=i+1$ if $\s$ has the structure of Case 4 in Figure \ref{fig:5.22}, and $s=i+2$ if $\s$ has the structure of Case 5 in Figure \ref{fig:5.23}.\\

As before, we define every SE glue sum precisely, so that\[
(\s\sesum_x^y\tau)^{-1}=\s^{-1}\nwsum_x^y\tau^{-1}
\]
for each of $(x,y)\in\{(1,0),(1,1),(2,0),(2,1),(3,0),(4,0)\}$. Thus, each type of SW glue sum is defined as in Table \ref{tab:5.2}. Here, $\s$ and $\tau$ are simple permutations in $\A'$ of length $m$ and $n$ respectively, $i=\s(m)$ and $j=\tau^{-1}$.\\

As we discussed in Chapter 4, type 1-0 and 1-1 glue sums are injective operations. It is not difficult to show that all other glue sums are also injective. Also, we note that all glues sums are associative only under certain conditions. In particular, both the left and the right summands must have certain lengths for each glue sum to be associative. We set the convention that when we sum multiple permutations with sequence of glue sums, we always sum from left to right.\\

By using all of the glue sums we defined, we describe the structure of half of the simple permutations in $\A'$, just as we did in Chapter \ref{chap:4}. In particular, we find the generating function for the set $H'=\{\pi\in\textrm{Si}(\A'):|\pi|\geq 4, 2\leq\pi(1)\leq 4\textrm{ and }\pi(2)\neq 1\}$. Notice that this set is the complement of the set of simple permutation $\pi$ in $\A'$ of length 4 or more such that $\pi(1)\geq 5$ or $\pi(2)=2$. Later, we will show that the condition of $H'$ is equivalent to $\pi^{-1}(1)\geq 5$ or $\pi(1)=2$. With what we observed for simple permutations of extreme pattern 2413, we can verify the set of extreme pattern 2413 simple permutations are a subset of $H'$. Hence, in summary, simple permutations in $\A'$ of length 4 or more can be classified as shown in Table \ref{tab:5.3}.
\begin{table}[H]
\begin{center}
\begin{tabular}{cccc}
\phantom{Column 1} & \phantom{Column 1} & \phantom{Column 1} & \phantom{Column 1}\\
\multicolumn{4}{c}{Simple permutation $\pi$ in $\A'$ with $|\pi|\geq 4$}\\
\hline
\multicolumn{2}{c|}{$\pi^{-1}(1)\geq 5$ or $\pi(1)=2$} & \multicolumn{2}{c}{$\pi(1)\geq 5$ or $\pi(2)=1$}\\
\hline
Extreme pattern 2413 & \multicolumn{2}{|c|}{Extreme pattern 2143} & Extreme pattern 3142\\
\end{tabular}
\end{center}
\caption{Classification of simple permutations in $\A'$.}
\label{tab:5.3}
\end{table}

The theorem describing the structure of simple permutations in $H'$ is the following.
\begin{theorem}\label{thm:5.1}
\textit{Let $\pi$ be a permutation in $H'$. Then there uniquely exist simple permutations $\s_i$ ($i$ odd) in $\A'$ of extreme pattern 2413 and simple permutations $\tau_i$ ($i$ even) in $\A'$ of extreme pattern 3142 such that\begin{equation}\label{eqn:5.6}
\pi=\left\{\begin{array}{clc}
\s_1\nwsum_{x_1}^{y_1}\tau_2\sesum_{x_2}^{y_2}\s_3\nwsum_{x_3}^{y_3}\tau_4\sesum_{x_4}^{y_4}\cdots\sesum_{x_{m-1}}^{y_{m-1}}\s_m & \textrm{ if }m\textrm{ is odd} &\qquad (a)\vhhh\\
\s_1\nwsum_{x_1}^{y_1}\tau_2\sesum_{x_2}^{y_2}\s_3\nwsum_{x_3}^{y_3}\tau_4\sesum_{x_4}^{y_4}\cdots\nwsum_{x_{m-1}}^{y_{m-1}}\tau_m & \textrm{ if }m\textrm{ is even} &\qquad (b)
\end{array}\right.
\end{equation}\\
where $m$ is a positive integer and $(x_\l,y_\l)\in\{(1,0),(1,1),(2,0),(2,1),(3,0),(4,0)\}$ ($1\leq\l\leq m-1$). Moreover, every permutation written as Equation \ref{eqn:5.1}$(a)$ or \ref{eqn:5.1}$(b)$ is in $H'$.}
\end{theorem}
\subsection{Proof of Theorem \ref{thm:5.1} (Part 1)}
As we did in Chapter \ref{chap:4}, we break the proof of Theorem \ref{thm:5.1} into two propositions and prove each one separately.
\begin{proposition}\label{prop:5.8}
\textit{If $\pi$ is a simple permutation in $H'$, then there uniquely exist simple permutations in $\A'$ of extreme pattern 2413 $\s_i$ ($i$ odd) and simple permutations in $\A'$ of extreme pattern 3142 $\tau_i$ ($i$ even) such that\begin{equation}\nn
\pi=\left\{\begin{array}{cl}
\s_1\nwsum_{x_1}^{y_1}\tau_2\sesum_{x_2}^{y_2}\s_3\nwsum_{x_3}^{y_3}\tau_4\sesum_{x_4}^{y_4}\cdots\sesum_{x_{m-1}}^{y_{m-1}}\s_m & \textrm{ if }m\textrm{ is odd}\vhhh\\
\s_1\nwsum_{x_1}^{y_1}\tau_2\sesum_{x_2}^{y_2}\s_3\nwsum_{x_3}^{y_3}\tau_4\sesum_{x_4}^{y_4}\cdots\nwsum_{x_{m-1}}^{y_{m-1}}\tau_m & \textrm{ if }m\textrm{ is even}
\end{array}\right.
\end{equation}\\
where $m$ is a positive integer and $(x_\l,y_\l)\in\{(1,0),(1,1),(2,0),(2,1),(3,0),(4,0)\}$ ($1\leq\l\leq m-1$).}
\end{proposition}
\textit{Proof.} Note that uniqueness of this decomposition follows from injectivity of the glue sums.\\

Suppose $\pi$ is a simple permutation in $\A'$ with $|\pi|\geq 4$ which satisfies either $\pi^{-1}(1)\geq 5$ or $\pi(1)=2$. We define a sequence of values of $\pi$, $d_1,\ldots, d_{m+3}$ in the exact same way as we did in Chapter \ref{chap:4}. Therefore, $d_1=\pi(1)$ and\[
d_i=\left\{\begin{array}{cl}
\pi(\max\{s:\pi(s)<\pi(d_{i-1})\}) &\textrm{ if }i\textrm{ is even}\\
\max\{t:\pi^{-1}(t)<\pi^{-1}(d_{i-1})\} &\textrm{ if }i\textrm{ is odd}
\end{array}\right.
\]\\
for $i$ with $1\leq i\leq m+3$. Thus, $d_i$ is the right-most value that is less than $d_{i-1}$ if $i$ is even or the greatest value located to the left of $d_{i-1}$ if $i$ is odd. As before, we have $d_{m+2}=\pi(n)$ and $d_{m+3}=n$ for some even integer $m$ or $d_{m+2}=n$ and $d_{m+3}=\pi(n)$ for some odd integer $m$.\\

We omit the proof of $d_i\neq d_j$ for $i\neq j$, since it is exactly the same as we did in Chapter \ref{chap:4}. Since $d_2=1$, we know $d_2\neq\pi(n)$, so every permutation in $H$ has at least four values denoted by $d_i$ for some $i\geq 1$. As before, we use mathematical induction on $m\geq 1$ to show that $\pi$ with $m+1$ values denoted by $d_i$ satisfies either Equation \ref{eqn:5.6}$(a)$ or \ref{eqn:5.6}$(b)$.\\

The base case is identical to the one in Chapter \ref{chap:4}. That is, if $m=1$, then there are four values $d_1$, $d_2$, $d_3$ and $d_4$ where $d_4=\pi(n)$. This implies $\pi$ is a simple permutation of extreme pattern 2413, so we are done.\\

Suppose every $\pi\in H'$ with $m+3$ values denoted by $d_i$ $(1\leq i\leq m+3)$ satisfies Equation \ref{eqn:5.6}$(a)$ for some positive odd integer $m$. We show that a permutation $\pi\in H'$ of length $n$ with $m+4$ values $d_i$ $(1\leq i\leq m+4)$ satisfies Equation \ref{eqn:5.6}$(b)$. Let $\pi$ be a permutation in $H'$ of length $n$ with $m+4$ values denoted by $d_i$. Thus, $d_{m+4}=n$.\\
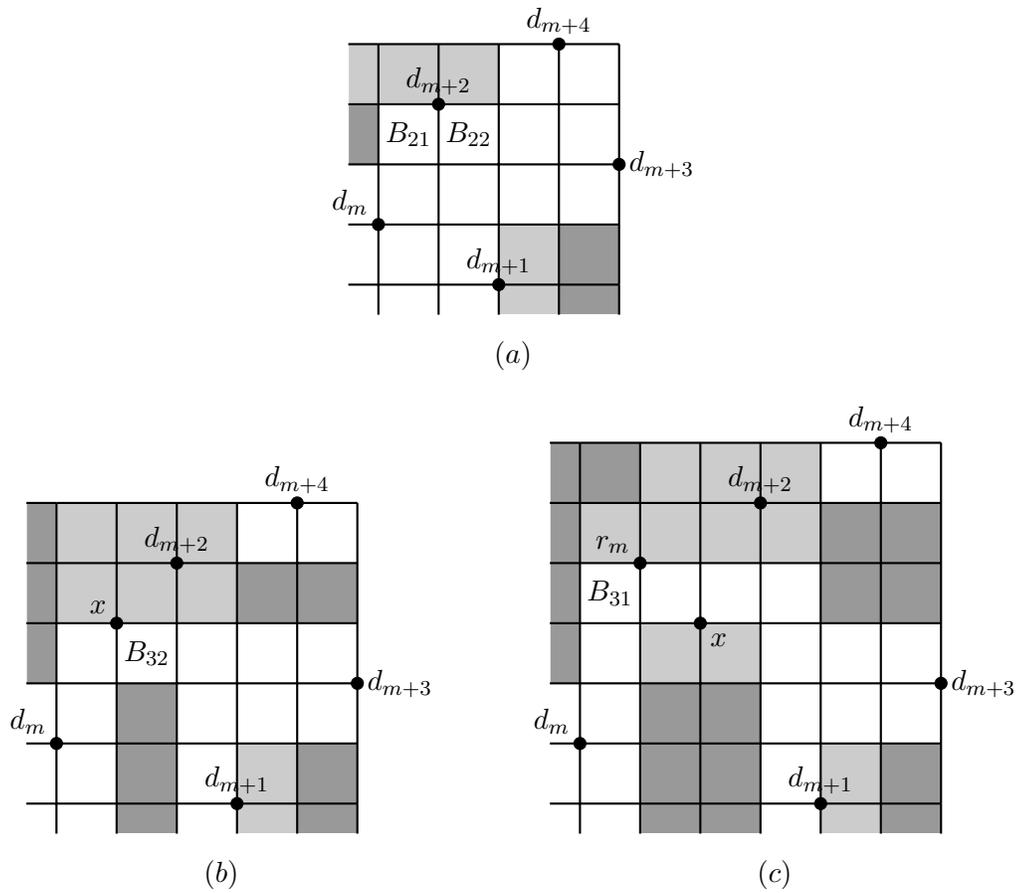
\begin{figure}[b!]
\[\begin{array}{c}
\begin{tikzpicture}[scale=0.8]
    \draw[white, fill=gray!40](0.5,4)--(3,4)--(3,5)--(0.5,5)--cycle;
    \draw[white, fill=gray!80](0.5,3)--(1,3)--(1,4)--(0.5,4)--cycle;
    \draw[white, fill=gray!80](4,0.5)--(5,0.5)--(5,2)--(4,2)--cycle;
    \draw[white, fill=gray!40](3,0.5)--(4,0.5)--(4,2)--(3,2)--cycle;
    \plotPerm{2,4,1,5,3}
    \foreach \x in {1,...,5} \draw[black,thick] (\x,1)--(\x,0.5) (1,\x)--(0.5,\x);
    \node[above left] at (1,2) {$d_m$};
    \node[above] at (2,4) {$d_{m+2}$};
    \node[above] at (3,1) {$d_{m+1}$};
    \node[above] at (4,5) {$d_{m+4}$};
    \node[right] at (5,3) {$d_{m+3}$};
    \node at (1.5,3.5) {$B_{21}$};
    \node at (2.5,3.5) {$B_{22}$};
\end{tikzpicture}\\
(a)
\end{array}\]
\[\begin{array}{ccc}
\begin{tikzpicture}[scale=0.8]
    \draw[white, fill=gray!80](0.5,3)--(1,3)--(1,6)--(0.5,6)--cycle;
    \draw[white, fill=gray!80](2,0.5)--(3,0.5)--(3,3)--(2,3)--cycle;
    \draw[white, fill=gray!80](4,4)--(6,4)--(6,5)--(4,5)--cycle;
    \draw[white, fill=gray!80](5,0.5)--(6,0.5)--(6,2)--(5,2)--cycle;
    \draw[white, fill=gray!40](1,4)--(4,4)--(4,6)--(1,6)--cycle;
    \draw[white, fill=gray!40](4,0.5)--(5,0.5)--(5,2)--(4,2)--cycle;
    \plotPerm{2,4,5,1,6,3}
    \foreach \x in {1,...,6} \draw[black,thick] (\x,1)--(\x,0.5) (1,\x)--(0.5,\x);
    \node[above left] at (1,2) {$d_m$};
    \node[above left] at (2,4) {$x$};
    \node[above] at (3,5) {$d_{m+2}$};
    \node[above] at (4,1) {$d_{m+1}$};
    \node[above] at (5,6) {$d_{m+4}$};
    \node[right] at (6,3) {$d_{m+3}$};
    \node at (2.5,3.5) {$B_{32}$};
  \end{tikzpicture}& \quad\quad & \begin{tikzpicture}[scale=0.8]
    \draw[white, fill=gray!80](0.5,3)--(1,3)--(1,6)--(2,6)--(2,7)--(0.5,7)--cycle;
    \draw[white, fill=gray!80](2,0.5)--(4,0.5)--(4,3)--(2,3)--cycle;
    \draw[white, fill=gray!80](5,4)--(7,4)--(7,6)--(5,6)--cycle;
    \draw[white, fill=gray!80](6,0.5)--(7,0.5)--(7,2)--(6,2)--cycle;
    \draw[white, fill=gray!40](1,5)--(5,5)--(5,7)--(2,7)--(2,6)--(1,6)--cycle;
    \draw[white, fill=gray!40](5,0.5)--(6,0.5)--(6,2)--(5,2)--cycle;
    \draw[white, fill=gray!40](2,3)--(4,3)--(4,4)--(2,4)--cycle;
    \plotPerm{2,5,4,6,1,7,3}
    \foreach \x in {1,...,7} \draw[black,thick] (\x,1)--(\x,0.5) (1,\x)--(0.5,\x);
    \node[above left] at (1,2) {$d_m$};
    \node[above left] at (2,5) {$r_{m}$};
    \node[below right] at (3,4) {$x$};
    \node[above] at (4,6) {$d_{m+2}$};
    \node[above] at (5,1) {$d_{m+1}$};
    \node[above] at (6,7) {$d_{m+4}$};
    \node[right] at (7,3) {$d_{m+3}$};
    \node at (1.5,4.5) {$B_{31}$};
  \end{tikzpicture}\\
  (b) & & (c)
\end{array}\]
\caption{Partial graphs of $\pi$ to show that there exists no value in $B_{21}$.}
\label{fig:5.26}
\end{figure}

Before we define the values $p_m$, $q_m$ and $r_m$, we claim that there is at most one value $r'_m$ besides $d_{m+2}$ such that $r'_m>d_{m+3}$ and $\pi^{-1}(r'_m)<\pi^{-1}(d_{m+1})$. Moreover, if such a point $r'_m$ exists, then the value $r'_m$ must satisfy $r'_m=d_{m+2}-1$ or $r'_m=d_{m+2}-2$, and the position $\pi^{-1}(r'_m)$ must satisfy $\pi^{-1}(r'_m)=\pi^{-1}(d_{m+2})+1$ or $\pi^{-1}(r'_m)=\pi^{-1}(d_{m+2})+2$. So let us start from the graph shown in Figure \ref{fig:5.26}$(a)$. We first show that we cannot have a point in the region $B_{21}$.\\

\begin{figure}[b!]
\[\begin{array}{ccc}
\begin{tikzpicture}[scale=0.8]
    \draw[white, fill=gray!40](1,3)--(2,3)--(2,4)--(3,4)--(3,5)--(1,5)--cycle;
    \draw[white, fill=gray!80](0.5,3)--(1,3)--(1,5)--(0.5,5)--cycle;
    \draw[white, fill=gray!80](4,0.5)--(5,0.5)--(5,2)--(4,2)--cycle;
    \draw[white, fill=gray!40](3,0.5)--(4,0.5)--(4,2)--(3,2)--cycle;
    \plotPerm{2,4,1,5,3}
    \foreach \x in {1,...,5} \draw[black,thick] (\x,1)--(\x,0.5) (1,\x)--(0.5,\x);
    \node[above left] at (1,2) {$d_m$};
    \node[above] at (2,4) {$d_{m+2}$};
    \node[above] at (3,1) {$d_{m+1}$};
    \node[above] at (4,5) {$d_{m+4}$};
    \node[right] at (5,3) {$d_{m+3}$};
    \node at (2.5,3.5) {$B_{22}$};
  \end{tikzpicture} & \quad\quad & \begin{tikzpicture}[scale=0.8]
    \draw[white, fill=gray!40](1,3)--(2,3)--(2,4)--(4,4)--(4,6)--(1,6)--cycle;
    \draw[white, fill=gray!80](0.5,3)--(1,3)--(1,6)--(0.5,6)--cycle;
    \draw[white, fill=gray!80](5,0.5)--(6,0.5)--(6,2)--(5,2)--cycle;
    \draw[white, fill=gray!80](2,0.5)--(3,0.5)--(3,1)--(2,1)--cycle;
    \draw[white, fill=gray!40](4,0.5)--(5,0.5)--(5,2)--(4,2)--cycle;
    \plotPerm{2,5,4,1,6,3}
    \foreach \x in {1,...,6} \draw[black,thick] (\x,1)--(\x,0.5) (1,\x)--(0.5,\x);
    \node[above left] at (1,2) {$d_m$};
    \node[above] at (2,5) {$d_{m+2}$};
    \node[above right] at (3,4) {$r'_m$};
    \node[above] at (4,1) {$d_{m+1}$};
    \node[above] at (5,6) {$d_{m+4}$};
    \node[right] at (6,3) {$d_{m+3}$};
    \node at (2.5,3.5) {$B_{32}$};
    \node at (3.5,3.5) {$B_{33}$};
    \node at (2.5,2.5) {$B_{42}$};
    \node at (2.5,1.5) {$B_{52}$};
    \node at (4.5,4.5) {$B_{24}$};
    \node at (5.5,4.5) {$B_{25}$};
  \end{tikzpicture}\\
  (a) & & (b)
  \end{array}
\]
\[\begin{array}{ccc}
\begin{tikzpicture}[scale=0.8]
    \draw[white, fill=gray!80](0.5,3)--(1,3)--(1,7)--(0.5,7)--cycle;
    \draw[white, fill=gray!80](2,0.5)--(4,0.5)--(4,4)--(7,4)--(7,6)--(3,6)--(3,5)--(2,5)--cycle;
    \draw[white, fill=gray!80](6,0.5)--(7,0.5)--(7,2)--(6,2)--cycle;
    \draw[white, fill=gray!40](1,3)--(2,3)--(2,5)--(3,5)--(3,6)--(5,6)--(5,7)--(1,7)--cycle;
    \draw[white, fill=gray!40](5,0.5)--(6,0.5)--(6,2)--(5,2)--cycle;
    \plotPerm{2,6,4,5,1,7,3}
    \foreach \x in {1,...,7} \draw[black,thick] (\x,1)--(\x,0.5) (1,\x)--(0.5,\x);
    \node[above left] at (1,2) {$d_m$};
    \node[above] at (2,6) {$d_{m+2}$};
    \node[below right] at (3,4) {$x$};
    \node[above right] at (4,5) {$r'_m$};
    \node[above] at (5,1) {$d_{m+1}$};
    \node[above] at (6,7) {$d_{m+4}$};
    \node[right] at (7,3) {$d_{m+3}$};
  \end{tikzpicture} & \quad\quad & \begin{tikzpicture}[scale=0.8]
    \draw[white, fill=gray!80](0.5,3)--(1,3)--(1,7)--(0.5,7)--cycle;
    \draw[white, fill=gray!80](2,0.5)--(4,0.5)--(4,1)--(3,1)--(3,4)--(2,4)--cycle;
    \draw[white, fill=gray!80](4,5)--(7,5)--(7,6)--(4,6)--cycle;
    \draw[white, fill=gray!80](6,0.5)--(7,0.5)--(7,2)--(6,2)--cycle;
    \draw[white, fill=gray!40](5,0.5)--(6,0.5)--(6,2)--(5,2)--cycle;
    \draw[white, fill=gray!40](1,3)--(2,3)--(2,4)--(3,4)--(3,5)--(4,5)--(4,6)--(5,6)--(5,7)--(1,7)--cycle;
    \plotPerm{2,6,5,4,1,7,3}
    \foreach \x in {1,...,7} \draw[black,thick] (\x,1)--(\x,0.5) (1,\x)--(0.5,\x);
    \node[above left] at (1,2) {$d_m$};
    \node[above] at (2,6) {$d_{m+2}$};
    \node[above right] at (3,5) {$r'_m$};
    \node[above right] at (4,4) {$x$};
    \node[above] at (5,1) {$d_{m+1}$};
    \node[above] at (6,7) {$d_{m+4}$};
    \node[right] at (7,3) {$d_{m+3}$};
  \end{tikzpicture} \\
  (c) & & (d)
  \end{array}
\]
\caption{Partial graphs of $\pi$ to show the possible existence of $r'_m$.}
\label{fig:5.27}
\end{figure}
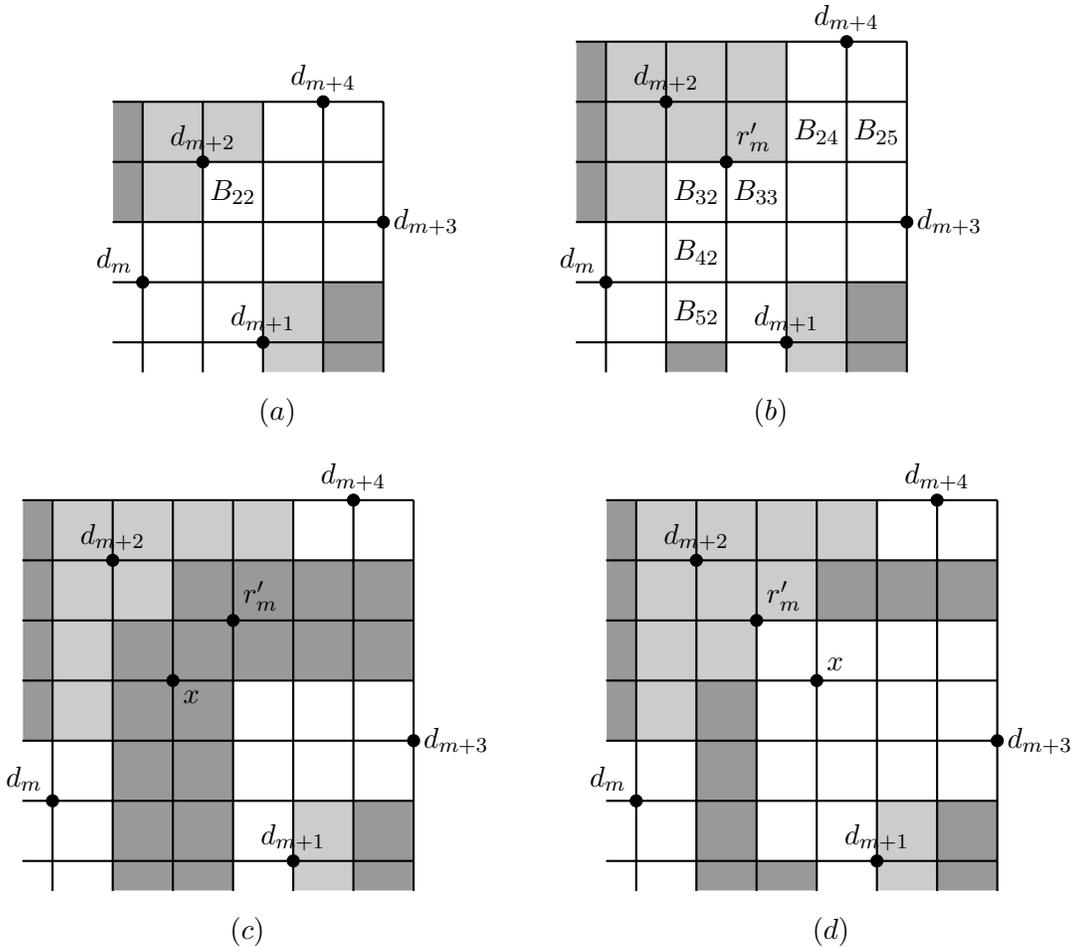

As usual, proceed by assuming there is a value in $B_{21}$, and let $x$ be the greatest such value, so we obtain the graph in Figure \ref{fig:5.26}$(b)$. In order to split the block $[\pi^{-1}(x),\pi^{-1}(d_{m+2})]$, we must have a point in the region $B_{32}$, so let $y$ be the least value there. We achieve the graph shown in Figure \ref{fig:5.26}$(c)$, requiring to having a point in either $B_{31}$ or $B_{34}$. However, having a point in either region implies the infinite chain contradiction that we have observed in previous proofs, such as the one for Lemma \ref{lem:5.2}. Therefore, having a point in the region $B_{21}$ of Figure \ref{fig:5.26}$(a)$ is prohibited, and we must have the structure of Figure \ref{fig:5.27}$(a)$.\\

Now, let $r'_m$ be the point with the greatest value in the region $B_{22}$ of Figure \ref{fig:5.27}$(a)$, so we have the graph as in Figure \ref{fig:5.27}$(b)$. Suppose to the contrary there is another point with the value $x$ such that $x>d_{m+3}$ and $\pi^{-1}(x)<\pi^{-1}(d_{m+1})$. If $x$ is in $B_{32}$ of Figure \ref{fig:5.27}$(b)$, then we have the graph shown in Figure \ref{fig:5.27}$(c)$ as a result, which has an unsplittable block $[\pi^{-1}(d_{m+2},\pi^{-1}(r'_m)]$. So we prohibit the region $B_{32}$ of Figure \ref{fig:5.27}$(b)$, and assume $x$ is in $B_{33}$ instead. Then we obtain the graph as in Figure \ref{fig:5.27}$(d)$, which also has a block $[\pi^{-1}(d_{m+2},\pi^{-1}(r'_m)]$ that cannot be split. Hence, we can only have one point, besides $d_{m+2}$ that has a value greater than $d_{m+3}$ and is located to the left of $d_{m+1}$.\\

When we have $r'_m$ as shown in Figure \ref{fig:5.27}$(b)$, having at least one point in either $B_{24}$, $B_{25}$, $B_{42}$ or $B_{52}$ is necessary to ensure $\pi$ is simple. However, due to 52341, 53241 and 52431 patterns avoidance, we can have only one point whose position is in $[\pi^{-1}(d_{m+2},\pi^{-1}(r'_m)]$ and only one point whose value is in $[r'_m,d_{m+2}]$. It is possible to have both of these two points. Thus, there are three cases, which are\begin{enumerate}
\item $\pi^{-1}(r'_m)=\pi^{-1}(d_{m+2})+1$ and $r'_m=d_{m+2}-2$,
\item $\pi^{-1}(r'_m)=\pi^{-1}(d_{m+2})+2$ and $r'_m=d_{m+2}-1$ and
\item $\pi^{-1}(r'_m)=\pi^{-1}(d_{m+2})+2$ and $r'_m=d_{m+2}-2$.
\end{enumerate}\vhhh

We are now ready to define the values $p_m$, $q_m$ and $r_m$. We define $p_m$ in the same way as before, that is,\[
p_m=\pi(\min\{\pi^{-1}(s):s>d_{m+3},\pi^{-1}(s)>\pi^{-1}(d_{m+1})\}).
\]

The definition of $q_m$ is much more complex than the one of Chapter \ref{chap:4}. We define $q_m$ by the following algorithm called \textsf{TRACE}.\\

\textbf{Algorithm} \textsf{TRACE}

\textsf{INPUT}: A simple permutation $\pi$ in $H'$.

\textsf{OUTPUT}: Value $q_m$. In addition, value $q'_m$ for \textsf{Step 3a}.\\
\begin{enumerate}[\sffamily Step 1:]
\item \textsf{Define the value $T_1=u_{m}$ and $T_2$.}\\
Let $T_1=u_{m}=\pi(\pi^{-1}(p_m)+1)$, \textit{i.e.} the value immediately to the right of $p_m$. Let $T_2=u_{m}+1$.
\item \textsf{Determine if the algorithm continues.}\\
If $T_1<p_m<d_{m+2}$ and $\pi^{-1}(T_2)<\pi^{-1}(T_1)$, then \textsf{GOTO STEP 3}. Otherwise, \textsf{OUTPUT} $q_m=\pi(\pi^{-1}(p_m)-1)$.
\item \textsf{Find the beginning of the 312-position chain.}\\
\textsf{WHILE} $\pi(\pi^{-1}(T_2)+1)< T_1-1$ \textsf{AND} $\pi^{-1}(T_2)>\pi^{-1}(d_{m+1})$:\\
Newly define $T_1=\pi(\pi^{-1}(T_2)+1)$ and $T_2=T_1+1$.
\item \textsf{Determine which condition of Loop was violated.}\\
\textsf{BEGIN}\begin{enumerate}[\sffamily a.]
    \item If $\pi(\pi^{-1}(T_2)+1)\geq T_1-1$, then \textsf{OUTPUT} $q_m=\pi(\pi^{-1}(T_2)-1)$.
    \item If $\pi^{-1}(T_2) < \pi^{-1}(d_{m+1})$, then \textsf{OUTPUT} $q_m=u_m$ and $q'_m=T_2$.
    \end{enumerate}
\end{enumerate}

Finally, we let $r_m=\max\{\pi(s)<d_{m+3}:s\in[1,\pi^{-1}(q_m)]\}$. Hence, $r_m$ is the greatest value less than $d_{m+3}$, whose position is in $[1,\pi^{-1}(q_m)]$.\\

Before we move on, we explain what the roles of $p_m$, $q_m$, $q'_m$ and $r_m$ are. For each NW glue sum, the value $q_m$ is the last value of the left hand side of $(\s_1\nwsum_{x_1}^{y_1}\tau_2\sesum_{x_2}^{y_2}\cdots\sesum_{x_{m-1}}^{y_{m-1}}\s_m)\nwsum_{x_m}^{y_m}\tau_{m+1}$. Hence, in Table \ref{tab:5.1}, $q_m$ is $\s'(m-1)$ for $\nwsum_{1}^{0}$, $\nwsum_{2}^{0}$ and $\nwsum_{3}^{0}$, $\s'(m)$ for $\nwsum_{1}^{1}$ and $\nwsum_{2}^{1}$, and $\s'(m-2)$ for $\nwsum_{4}^{0}$. Informally, we may consider that $q_m$ is where the left permutation and the right permutation are connected. The value $p_m$ is merely used to locate $q_m$. Except for the case of type 3-0 and some special cases of type 1-0 and type 1-1, $q_m$ is immediately to the left of $p_m$, so $p_m$ can be viewed as the first value to be glued. The other exceptional cases are explained with more details next. The value $r_m$ is to distinguish the type 1-1 and 2-1 from the type 1-0 and 2-0 respectively.\\

For a permutation in $H'$, the right-most summand of its glue sum decomposition can either begin with a value greater $d_{m+3}$ or a 312-position chain. The loop in \textsf{TRACE} finds the beginning of this chain. The two exit conditions of this loop correspond to the cases where the 312-position chain ends on its own, or with a scissor on its left, as in the following examples. Suppose we have $\mu_1=264135$, $\mu_2=2\,\,10\,\,5\,\,1\,\,3\,\,7\,\,4\,\,9\,\,6\,\,8$, $\nu_1=71426385$ and $\nu_2=831527496$. Notice that these are simple permutations in $\A'$, and furthermore, $\mu_1$ and $\mu_2$ are of extreme pattern 2413 and $\nu_1$ and $\nu_2$ are of extreme pattern 3142. We obtain\[
\mu_1\nwsum_1^0\nu_1=2\,\,10\,\,4\,\,1\,\,3\,\,7\,\,5\,\,9\,\,6\,\,11\,\,8\qquad\textrm{and}\qquad\mu_2\nwsum_3^0\nu_2=2\,\,10\,\,5\,\,1\,\,3\,\,7\,\,4\,\,9\,\,6\,\,11\,\,8\]
\begin{figure}[t!]
\[\begin{array}{ccc}\begin{tikzpicture}[scale=.4]
\draw[color=black,fill=black] (1,2) circle [radius=.1];
\draw[color=black,fill=black] (2,10) circle [radius=.1];
\draw[color=black,fill=black] (3,4) circle [radius=.1];
\draw[color=black,fill=black] (4,1) circle [radius=.1];
\draw[color=black,fill=black] (5,3) circle [radius=.1] node[below] {$q_1$};
\draw[color=black,fill=black] (6,7) circle [radius=.1];
\draw[color=black,fill=black] (7,5) circle [radius=.1] node[below] {$T_1$};
\draw[color=black,fill=black] (8,9) circle [radius=.1] node[above] {$p_1$};
\draw[color=black,fill=black] (9,6) circle [radius=.1] node[below] {$T_2$};
\draw[color=black,fill=black] (10,11) circle [radius=.1];
\draw[color=black,fill=black] (11,8) circle [radius=.1];
\end{tikzpicture}& \qquad\qquad\qquad\qquad & \begin{tikzpicture}[scale=0.4]
\draw[color=black,fill=black] (1,2) circle [radius=.1];
\draw[color=black,fill=black] (2,10) circle [radius=.1];
\draw[color=black,fill=black] (3,5) circle [radius=.1] node[above] {$q'_1=T_2$};
\draw[color=black,fill=black] (4,1) circle [radius=.1];
\draw[color=black,fill=black] (5,3) circle [radius=.1];
\draw[color=black,fill=black] (6,7) circle [radius=.1];
\draw[color=black,fill=black] (7,4) circle [radius=.1] node[below] {$T_1$};
\draw[color=black,fill=black] (8,9) circle [radius=.1] node[above] {$p_1$};
\draw[color=black,fill=black] (9,6) circle [radius=.1] node[right] {$q_1$};
\draw[color=black,fill=black] (10,11) circle [radius=.1];
\draw[color=black,fill=black] (11,8) circle [radius=.1];
\end{tikzpicture}\\
\\
\mu_1\nwsum_1^0\nu_1 & & \mu_2\nwsum_3^0\nu_2
\end{array}\]
\caption{Graphs of $\mu_1\nwsum_1^0\nu_1$ and $\mu_2\nwsum_3^0\nu_2$.}
\label{fig:5.31}
\end{figure}

The graphs of $\mu_1\nwsum_1^0\nu_1$ and $\mu_2\nwsum_3^0\nu_2$ are shown in Figure \ref{fig:5.31}. What we have in common in these two sums is that the value corresponding to $p_1$ is a scissor of a 312-position chain before the sum. In case of $\mu_1\nwsum_1^0\nu_1$, the values in $[1,5)$ of $\nu_1$ are of the form $1\oplus 312$. Since 312 alone is a 312-position chain, we have the value 6 as the scissor, and this value corresponds to $p_1$ in $\mu_1\nwsum_1^0\nu_1$. When the right summand does not begin with a 312-position chain, the value immediately to the left of $p_m$ is $q_m$, but when it does as in this example, trace follows this chain to find $q_1=m_1'(5)=3$ in \textsf{Step 4a.}\\

On the other hand, the values in $[1,6)$ of $\nu_2$ are of the form $\a=312\oplus^1 21$, a 312-position chain. The value 7 corresponding to $p_1$ in $\mu_2\nwsum_3^0\nu_2$ is still a scissor, so we again find the beginning of the 312-position chain using \textsf{TRACE.} However, in this case, we end with $T_2$ to the left of $1=d_{2}=d_{m+1}$, which means the 312-position chain also has a scissor on its left. We define this scissor to be $q'_1$, and the value immediately to the right of $p_1$ to be $q_1$, according to \textsf{Step 4b}. Hence, the value 4 corresponding to $q'_1$ is the scissor of 312-value chain of $\mu_2$ that is similar to the 312-position chain of $\nu_2$. Thus, this whole procedure is done to determine whether two permutations are glued by type 1-0 (or type 1-1) with the values of the right summand starting with $1\oplus\a$ where $\a$ is a 312-position chain or they are glued by type 3-0.\\

In summary, $p_m$ is used to located $q_m$, $q_m$ is the last value of the left summand, $q'_m$ is the scissor of the last 312-value chain in the left summand, and $r_m$ is used to distinguish the type 1-1 and 2-1 from the type 1-0 and 2-0 respectively.\\

From here, we divide into four cases based on the existence of $r'_m$, the value greater than $d_{m+3}$ that is located to the left of $d_{m+1}$, and how $q_m$ is determined by \textsf{TRACE}. These four cases are namely the following, and they correspond to the type of the glue sum used.\begin{enumerate}[\textit{Case} \itshape A:]
\item There does not exist $r'_m$ and $q_m$ is output by \textsf{STEP 2} or \textsf{STEP 4a}.
\item There exists $r'_m=d_{m+2}-1$.
\item There does not exist $r'_m$ and $q_m$ is output by \textsf{STEP 4b}.
\item There exists $r'_m=d_{m+2}-2$.
\end{enumerate}

In each case, we slice $\pi$ into $\pi_1$ and $\pi_2$ according to the position of $q_m$. Then we show $\pi_1$ is a permutation in $H'$ with $m+3$ values denoted by $d_i$, and $\pi_2$ is a simple permutation in $\A'$ of extreme pattern 3142. Then by using the appropriate NW glue sum, we show $\pi=\pi_1\nwsum_x^y\pi_2$ ($(x,y)\in\{(1,0),(1,1),(2,0),(2,1),(3,0),(4,0)\}$) to complete the proof.\\
\\
\textit{Case A:} Denote by $\pi_2$ the flattening of the subsequence of $\pi$ obtained by removing every value whose position is in $[1,\pi^{-1}(q_m)]$, expect $d_{m+2}$ and $d_{m+1}$. Showing $\pi_2$ is a simple permutation in $\A'$ of extreme pattern 3142 is exactly the same as we did in Chapter \ref{chap:4}, so we omit it.\\

Next, let $\pi_1$ be the flattening of the subsequence $\pi(1)\pi(2)\cdots\pi(\pi^{-1}(q_m)-1)d_{m+3}$ of $\pi$ if $q_m=r_m$, and $\pi(1)\pi(2)\cdots q_m d_{m+3}$ of $\pi$ if $q_m\neq r_m$. We claim that $\pi_1$ is simple.\\

As we did in Chapter \ref{chap:4}, we first show that every value we removed to construct $\pi_1$ is greater than or equal to $r_m$. Suppose $q_m=r_m$. Then the positions of values that are removed are in $[\pi^{-1}(q_m),\pi^{-1}(d_{m+3})-1]$. Assume to the contrary that there exists $z<r_m$ such that $\pi^{-1}(z)\in [\pi^{-1}(q_m),\pi^{-1}(d_{m+3})-1]$. Note that $z\neq p_m$ since $p_m>r_m$. If $p_m>d_{m+2}$, then $d_{m+2} r_m p_m z d_{m+3}$ forms a 42513 pattern. So suppose $p_m<d_{m+2}$. Let $w$ be the value whose position is immediately to the right of the position of $p_m$. For now, assume $w>r_m$. Now, if $w>d_{m+2}$, then $\pi$ contains a 42513 pattern with $d_{m+2}r_m w z d_{m+3}$. On the other hand, if $w<d_{m+2}$, then $\pi$ contains either a 52341 pattern or a 52431 pattern with $d_{m+2} q_m p_m w z$. Thus, we have $w<r_m$, which is impossible since $d_{m+2} q_m (w+1) p_m w$ forms a 53241 pattern.\\

Now, assume $q_m\neq r_m$. Then the positions of values that are removed are in $[\pi^{-1}(p_m),\allowbreak\pi^{-1}(d_{m+3})-1]$. Suppose there exists $z<r_m$ such that $\pi^{-1}(z)\in [\pi^{-1}(p_m),\pi^{-1}(d_{m+3})-1]$. Again, $z\neq p_m$, so we must have $\pi^{-1}(z)>\pi^{-1}(p_m)$. If we assume $\pi^{-1}(r_m)>\pi^{-1}(d_{m+2})$ and $p_m<d_{m+3}$, we achieve a similar contradiction as before. Also, the rest of the proof is identical to the argument in Chapter \ref{chap:4}. That is, if $\pi^{-1}(r_m)>\pi^{-1}(d_{m+2})$ and $p_m>d_{m+2}$, then $\pi$ contains 42513 pattern with $d_{m+2}r_m p_m zd_{m+3}$. If $\pi^{-1}(r_m)<\pi^{-1}(d_{m+2})$ and $p_m<d_{m+2}$, then $\pi$ contains 35142 pattern with $r_m d_{m+2} d_{m+1} p_m z$, and finally, if $\pi^{-1}(r_m)<\pi^{-1}(d_{m+2})$ and $p_m>d_{m+2}$, then $\pi$ contains 351624 pattern with $r_m d_{m+2}d_{m+1} p_m z d_{m+3}$.\\

Hence, every value we removed to construct $\pi_1$ must be greater than or equal to $r_m$. With this, we can show $\pi_1$ is simple in the exact same way as in Chapter \ref{chap:4}. Therefore, $\pi_1\in H'$ with $m+3$ values denoted by $d_i$. Moreover, by the way we constructed $\pi_1$ and $\pi_2$, $\pi_1$ and $\pi_2$ respectively have structures of Case 1 in Figure \ref{fig:5.22} and Case 1 in Figure \ref{fig:5.23}. If $q_m=r_m$, then $\pi=\pi_1\nwsum_1^1\pi_2$. Otherwise, $\pi=\pi_1\nwsum_1^0\pi_2$, so we are done with Case $A$.\\
\\
\textit{Case B:} As before, denote by $\pi_2$ the flattening of the subsequence of $\pi$ obtained by removing every value whose position is in $[1,\pi^{-1}(q_m)]$, except $d_{m+2}$ and $d_{m+1}$. We can show $\pi_2$ is a simple permutation in $\A'$ of extreme pattern 3142 in the same way as Chapter \ref{chap:4}. We also define $\pi_1$ in the same way as Case $A$, so it is the flattening of the subsequence $\pi(1)\pi(2)\cdots\pi(\pi^{-1}(q_m)-1)d_{m+3}$ of $\pi$ if $q_m=r_m$, and $\pi(1)\pi(2)\cdots q_m d_{m+3}$ of $\pi$ if $q_m\neq r_m$. Showing every value we removed to construct $\pi_1$ is greater than or equal to $r_m$ is identical to Case $A$.\\

Suppose to the contrary $\pi_1$ is not simple. Let us use the hat notation as we did in Chapter \ref{chap:4} to refer to the value of $\pi_1$ corresponding to a value of $\pi$. Let $I$ be a proper non-singleton block of $\pi_1$. We have three cases to consider. First, if $\pi_1(I)=[\hat{x},\hat{y}]$ for some values $\hat{x}$ and $\hat{y}$ where $\hat{y}\leq\hat{d}_{m+3}$, we use the same argument as Chapter \ref{chap:4}. Second, suppose $\pi_1(I)=[\hat{x},\hat{r}'_m]$. Because $r_m<d_{m+3}<r'_{m}$, we have $\hat{d}_{m+3}=\hat{r}'_{m}-1$. Thus, $\hat{d}_{m+3}\in[\hat{x},\hat{r}'_{m}]$. Since $\pi_1{-1}(\hat{r}'_m)<\pi_1^{-1}(\hat{d}_{m+1})<\pi_1^{-1}(\hat{d}_{m+3})$, we have $\hat{d}_{m+1}\in[\hat{x},\hat{r}'_m]$. We can then argue that $\hat{d}_{m}$, and hence $\hat{d}_{m+2}$, are also in $[\hat{x},\hat{r}'_m]$, which is a contradiction, since $\hat{d}_{m+2}>\hat{r}'_m$. Finally, if $\pi_1(I)=[\hat{x},\hat{d}_{m+2}]$, then we use the same argument as Chapter \ref{chap:4}. (It is impossible to have $\hat{x}=\hat{r}'_m$ because the point in the region $B_{42}$ or $B_{52}$ as shown in Figure \ref{fig:5.30}$(b)$ is included in $\pi_1$.) Consequently, $\pi_1$ is simple.\\

With $r'_m=\hat{d}_{m+2}-1$ in $\pi_1$, the structure of $\pi_1$ is Case 2 in Figure \ref{fig:5.23}, whereas $\pi_2$ still has the structure of Case 2 in Figure \ref{fig:5.24}. By definitions of $\nwsum_2^0$ and $\nwsum_2^1$, $\pi=\pi_1\nwsum_2^0\pi_2$ if $q_m=r_m$, and $\pi=\pi_1\nwsum_2^1\pi_2$ otherwise.\\
\\
\textit{Case C:} This time, we let $\pi_2$ be the flattening of the subsequence of $\pi$ obtained by removing every value whose position is in $[1,\pi^{-1}(q'_m-1)-2]$, except $d_{m+2}$, $q'_m$ and $d_{m+1}$. We show $\pi_2$ is a simple permutation in $\A'$ of extreme pattern 3142.\\

Suppose $\pi_2$ is not simple, so there exists a proper non-singleton block $I$. If $I=[\pi_2^{-1}(\hat{x}),\allowbreak\pi_2^{-1}(\hat{y})]$ for some values $\hat{x}$ and $\hat{y}$ where $\pi_2^{-1}(\hat{x})>\pi_2^{-1}(\hat{d}_{m+1})$, then $[\pi^{-1}(\hat{x}),\pi^{-1}(\hat{y})]$ is also a block in $\pi$. Next, suppose $I=[\pi_2^{-1}(\hat{d}_{m+1}),\pi_2^{-1}(\hat{y})]$ for some $\hat{y}$. By the way $q'_m$ is defined, the value immediately to the right of the value $\hat{d}_{m+1}$ must be greater than $\hat{q}'_m$. Hence, $\pi_2^{-1}(\hat{q}'_m)\in I$, so we have a contradiction. So suppose $I=[\pi_2^{-1}(\hat{q}'_m),\pi_2^{-1}(\hat{y})]$ for some $\hat{y}$. Let $T'_1=q'_m-1$ and $T'_2=\pi_2(\pi_2^{-1}(T'_1)-1)$. Then $\pi_2^{-1}(T'_1),\pi_2^{-1}(T'_2)\in I$. Next, if $T'_2\neq\hat{p}_m$, we let $T'_3=T'_2-1$ and $T'_4=\pi_2(\pi_2^{-1}(T'_3)-1)$, so we have $\pi_2^{-1}(T'_3),\pi_2^{-1}(T'_4)\in I$. So long as $T'_i$ ($i$ even) is not equal to $\hat{p}_m$, we continue definition $T'_i$ in the same way, so\[
T'_i=\left\{\begin{array}{cl}
T'_{i-1}-1 &\textrm{ if }i\textrm{ is odd}\\
\pi_2(\pi_2^{-1}(T'_{i-1})-1) &\textrm{ if }i\textrm{ is even}
\end{array}\right.
\]
The way $q'_m$ is defined guarantees that there exists an even integer $j$ such that $T'_j=\hat{p}_m$ (and thus, $T'_{j-1}=\hat{q}_m$), so $\pi_2^{-1}(\hat{p}_m)\in[\pi_2^{-1}(\hat{q}'_m),\pi_2^{-1}(\hat{y})]$. Since $\hat{p}_m>\hat{d}_{m+3}$, $\pi_2^{-1}(\hat{d}_{m+3})\in I$, however, this implies $\pi_2^{-1}(\hat{d}_{m+4})\in I$, and therefore, $\pi_2^{-1}(\hat{d}_{m+2})\in I$ as well. Since $\pi_2^{-1}(\hat{d}_{m+2})<\pi_2^{-1}(\hat{q}'_m)$, we have a contradiction. Lastly, assume $I=[\pi_2^{-1}(\hat{d}_{m+2}),\pi_2^{-1}(\hat{y})]$ for some $\hat{y}$. Because $\pi_2^{-1}(\hat{q}'_m)=\pi_2^{-1}(\hat{d}_{m+2})+1$, we must have $\pi_2^{-1}(\hat{q}'_m)\in I$. Since, $\hat{q}'_m<\hat{d}_{m+3}$, we obtain $I=[\pi_2^{-1}(\hat{q}'_m),\pi_2^{-1}(\hat{d}_{m+3})]$, but then $I$ is all of $\pi$.\\

In every case, we achieve a contradiction, so $\pi_2$ is simple. In particular, $\pi_2$ has extreme pattern 3142 and has the structure of Case 3 in Figure \ref{fig:5.23}. By Proposition \ref{prop:5.3}, $\hat{q}'_m$ and $\hat{d}_{m+1}$ are involved in the same 312-position chain. In fact, by Proposition \ref{prop:5.5}, every value whose position is in $[\pi_2^{-1}(\hat{q}'_m),\pi_2^{-1}(\hat{q}_m)]$, except $\hat{p}_m$ together forms this 312-position chain.\\

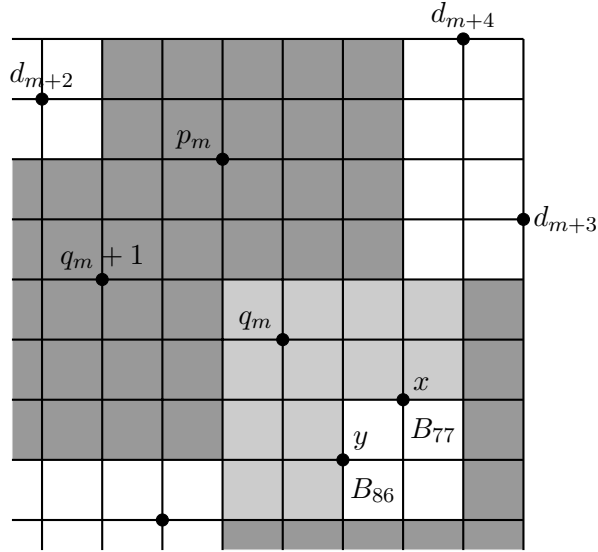
\begin{figure}[t!]
\[\begin{tikzpicture}[scale=0.8]
    \draw[white, fill=gray!40](4,1)--(6,1)--(6,3)--(8,3)--(8,5)--(4,5)--cycle;
    \draw[white, fill=gray!80](4,0.5)--(9,0.5)--(9,5)--(8,5)--(8,1)--(4,1)--cycle;
    \draw[white, fill=gray!80](0.5,2)--(4,2)--(4,5)--(7,5)--(7,9)--(2,9)--(2,7)--(0.5,7)--cycle;
    \plotPerm{8,5,1,7,4,2,3,9,6}
    \foreach \x in {1,...,9} \draw[black,thick] (\x,1)--(\x,0.5) (1,\x)--(0.5,\x);
    \node[above] at (1,8) {$d_{m+2}$};
    \node[above] at (2,5) {$q_m+1$};
    \node[above left] at (4,7) {$p_m$};
    \node[above left] at (5,4) {$q_m$};
    \node[above right] at (6,2) {$y$};
    \node[above right] at (7,3) {$x$};
    \node[above] at (8,9) {$d_{m+4}$};
    \node[right] at (9,6) {$d_{m+3}$};
    \node at (6.5,1.5) {$B_{86}$};
    \node at (7.5,2.5) {$B_{77}$};
  \end{tikzpicture}\]
\caption{Partial graph of $\pi$ with the assumption of $\pi^{-1}(x)>\pi^{-1}(q_m)$.}
\label{fig:5.28}
\end{figure}

Next, let $\pi_1$ be the flattening of $\pi(1)\pi(2)\cdots q_m d_{m+3}$. We first show that every value $z$ removed to construct $\pi_1$ is greater than $q_m$. In order to show this, we first claim that the value $x=q_m-1$ is located to the left of $q_m$. Suppose to the contrary that $\PS{x}>\PS{q_m}$. We divide into two cases. Let $y$ be the value immediately to the right of $q_m$. Suppose $y>q_m$. Then $\pi$ contains 42513 due to $d_{m+2} q_m y x d_{m+3}$ if $d_{m+2}<y$, or 53241 due to $d_{m+2} (q_m+1) q_m y x$ if $y<d_{m+2}$. So assume $y<q_m$. Note that $x\neq y$ because $x=y$ implies that $\pi$ contains an unsplittable block $[\PS{q_m},\PS{x}]$. Hence, we obtain the graph shown in Figure \ref{fig:5.28}. Splitting the block $[\PS{q_m},\PS{y}]$ can be only done by having a point in $B_{77}$ or $B_{86}$. In either case, we obtain an unsplittable block, so we achieve a contradiction. Thus, we must have $\PS{x}<\PS{q_m}$.\\

We now show that every value located to the right of $q_m$ is greater than $q_m$. Suppose to the contrary that there exists a value $z<q_m$ such that $\PS{z}>\PS{q_m}$. If $\PS{x}<\PS{d_{m+2}}$, then $\pi$ contains 35142 due to $x d_{m+2} d_{m+1} p_m z$. If $\PS{d_{m+2}}<\PS{x}<\PS{q_m+1}$, then $\pi$ contains 52341 due to $d_{m+2} x (q_m+1) p_m z$. Finally, if $\PS{q_m+1}<\PS{x}$, then $\pi$ contains 53241 due to $d_{m+2} (q_m+1) x p_m z$. Consequently, every value we removed to construct $\pi_1$ is greater than $q_m$.\\

Showing $\pi_1$ is simple is essentially the argument of inverse followed by reverse complement symmetry of $\pi_2$. Suppose $\pi_1$ is not simple, so there exists a proper non-singleton block $I$ of $\pi_1$. If $\pi_1(I)=[\hat{x},\hat{y}]$ where $\hat{x}<\hat{y}\leq\hat{d}_{m+3}-1$, then $I$ would be also a block in $\pi$. If $\pi_1(I)=[\hat{x},\hat{d}_{m+3}]$, then $\hat{d}_{m+3}-1=\hat{q}_m+1$ must be in $\pi_1(I)$. However, since $\pi_1^{-1}(\hat{q}_m+1)<\pi_1^{-1}(\hat{p}_m)<\pi_1^{-1}(\hat{d}_{m+3})$, $\hat{p}_m$ also has to be in $\pi_i(I)$. Because $\hat{p}_m>\hat{d}_{m+3}$, we achieve a contradiction. So assume $\pi_1(I)=[\hat{x},\hat{p}_{m}]$. Due to the way $q'_m$ in $\pi$ was assigned, $\hat{q}'_m$ must be in $\pi_1(I)$. However, this implies $\hat{d}_{m+1}$ is also in $\pi_1(I)$, causing $\hat{d}_m\in\pi_1(I)$, and therefore, $\hat{d}_{m+2}\in\pi_1(I)$. Since $\hat{d}_{m+2}>\hat{p}_m$, this is a contradiction. Finally, suppose $\pi_1(I)=[\hat{x},\hat{d}_{m+2}]$. Then $\hat{p}_m\in\pi_1(I)$, which results in $\hat{x}=\hat{d}_{m+1}$. Since $I$ is all of $\pi_1$, we again have a contradiction.\\

Therefore, $\pi_1$ is a simple permutation in $H'$ with $m+3$ values denoted by $d_i$. Note that the points $\pi_1$ and $\pi_2$ both contain as results of removal are $d_{m+2}$, $q'_m$, $d_{m+1}$, every value whose position is in $[\pi^{-1}(q'_m-1)-1,\pi^{-1}(q_m)]$ (including $p_m$), and $d_{m+3}$. Hence, $\pi_1$ ends with a 312-value chain $\a$ involving $q_m$, $p_m$ and $d_{m+3}$. Also, $\pi_2$ begins with a 312-position chain $\b$ involving $q'_m$, $d_{m+1}$ and $q'_m-1$. Then $21\oplus_1\a=\b\oplus^1 21$, so $\a$ and $\b$ are similar. Since $\pi$ has the structure of Case 3 in Figure \ref{fig:5.22}, we have $\pi=\pi_1\nwsum_3^0\pi_2$, so we are done with Case $C$.\\
\\
\textit{Case D:} Finally, let $\pi_2$ be the flattening of the subsequence of $\pi$ obtained by removing every value whose position is in $[1,\pi^{-1}(q_m)]$, expect $d_{m+2}$, $r'_m$ and $d_{m+1}$. We start by showing $\pi_2$ is simple. Assume it is not simple. Let $I$ be a proper non-singleton block of $\pi_2$. As usual, if $I=[\pi_2^{-1}(x),\pi_2^{-1}(y)]$ for some values $\hat{x}$ and $\hat{y}$ where $\pi_2^{-1}(x)>\pi_2^{-1}(\hat{d}_{m+1})$, then $I$ would be a block in $\pi$ as well. Next, suppose $I=[\pi_2^{-1}(\hat{d}_{m+1}),\pi_2^{-1}(\hat{y})]$. Then $\pi_2^{-1}(\hat{p}_m)\in I$, resulting in $\pi_2^{-1}(\hat{d}_{m+3})\in I$. Hence, $\pi_2^{-1}(\hat{d}_{m+4})$ is in $I$, and therefore, so is $\pi_2^{-1}(\hat{d}_{m+2})$. Since $\pi_2^{-1}(\hat{d}_{m+2})<\pi_2^{-1}(\hat{d}_{m+1})$, we have a contradiction. This time, let $I=[\pi_2^{-1}(\hat{r}'_m),\pi_2^{-1}(\hat{y})]$. Since $\pi_2^{-1}(\hat{r}'_m)+1=\pi_2^{-1}(\hat{d}_{m+1})$, we know $\pi_2^{-1}(\hat{d}_{m+1})\in I$. The rest is the same as the previous case, $I=[\pi_2^{-1}(\hat{d}_{m+1}),\pi_2^{-1}(\hat{y})]$. We achieve a contradiction of $\pi_2^{-1}(\hat{d}_{m+2})\in I$. Lastly, suppose $I=[\pi_2^{-1}(\hat{d}_{m+2}),\pi_2^{-1}(\hat{y})]$. Then $\pi_2^{-1}(\hat{r}'_m)$ is in $I$. Since $\hat{r}'_m<\hat{p}_m<\hat{d}_{m+2}$, we know $\pi_2^{-1}(\hat{p}_m)\in I$, implying $\pi_2^{-1}(\hat{d}_{m+1})\in I$. Hence, $\pi_2^{-1}(\hat{d}_{m+3})$ is also in $I$, but then now $I$ is all of $\pi_2$. Altogether, $\pi_2$ is simple.\\

Next, let $\pi_1$ be the flattening of $\pi(1)\pi(2)\cdots q_m p'_m d_{m+3}$ where $p'_m=d_{m+2}-1$. Suppose to the contrary that there exists a value $z<r_m$ located to the right of $q_m$. If $\pi^{-1}(z)<\pi^{-1}(p'_m)$, then $p'_m\neq p_m$, and $r'_m p_m z p'_m d_{m+3}$ forms a 35142 pattern. So assume $\pi^{-1}(z)>\pi^{-1}(p'_m)$. Then we have a 53241, 52341 or 35142 pattern depending the position of $r_m$. So every value we removed to construct $\pi_1$ must be greater than $r_m$. Now, we prove $\pi_1$ is simple. Assume it is not simple, so we have a proper non-singleton block $I$. If $\pi_1(I)=[\hat{x},\hat{y}]$ for some $\hat{x},\hat{y}$ with $\hat{x}<\hat{y}<\hat{d}_{m+3}$, then $I$ would be a block in $\pi$. Suppose $\pi_1(I)=[\hat{x},\hat{d}_{m+3}]$. Since $\pi_1^{-1}(\hat{p}'_m)=\pi_1^{-1}(\hat{d}_{m+3})-1$, we must have $\hat{p}'_m\in\pi_1(I)$. Since $\hat{p}'_m>\hat{d}_{m+3}$, we have a contradiction. So assume $\pi_1(I)=[\hat{x},\hat{r}'_{m}]$. Since $\hat{d}_{m+3}=\hat{r}'_m-1$, we have $\hat{d}_{m+3}\in\pi_1(I)$. With $\pi_1^{-1}(\hat{r}'_m)<\pi_1^{-1}(\hat{p}'_m)<\pi_1^{-1}(\hat{d}_{m+3})$, we obtain $\hat{p}'_m\in\pi_1(I)$, but since $\hat{p}'_m>\hat{r}'_m$, we again have a contradiction. Suppose $\pi_1(I)=[\hat{x},\hat{p}'_m]$. Because $\hat{r}'_m=\hat{p}'_m-1$, $\hat{r}'_m\in\pi_1(I)$. With $\pi_1^{-1}(\hat{r}'_m)<\pi_1^{-1}(\hat{d}_{m+1})<\pi_1^{-1}(\hat{p}'_m)$, we know $\pi_1^{-1}(\hat{d}_{m+1})\in I$, and this implies $\hat{d}_m\in\pi_1(I)$. Because $\pi_1^{-1}(\hat{d}_m)<\pi_1^{-1}(\hat{d}_{m+2})<\pi_1^{-1}(\hat{p}'_m)$, $\hat{d}_{m+2}$ must be in $\pi_1(I)$, but since $\hat{d}_{m+2}>\hat{p}'_m$, we achieve a contradiction. Finally, let $\pi_1(I)=[\hat{x},\hat{d}_{m+2}]$. Since $\hat{p}'_m=\hat{d}_{m+2}-1$, we must have $\hat{p}'_m\in\pi_1(I)$, but this implies $\hat{d}_{m+1}\in\pi_1(I)$ due to the fact $\pi_1^{-1}(\hat{d}_{m+2})<\pi_1^{-1}(\hat{d}_{m+1})<\pi_1^{-1}(\hat{p}'_m)$. The block $I$ is now all of $\pi_1$, so this is a contradiction. Hence, $\pi_1$ must be simple.\\

With $\hat{d}_{m+2}$, $\hat{r}'_m$, $\hat{p}'_m$ and $\hat{d}_{m+3}$, $\pi_1$ has the structure described in Case 4 or 5 of Figure \ref{fig:5.22}. On the other hand, with $\hat{d}_{m+2}$,  $\hat{r}'_m$ and $\hat{d}_{m+1}$, $\pi_2$ has the structure described in Case 2 of Figure \ref{fig:5.23}. Together, we can recover $\pi$ by $\pi_1\nwsum_4^0\pi_2$.\\

Consequently, in every case, we can express $\pi$ as $\pi_1\nwsum_x^y\pi_2$ for appropriate $\pi_1$ and $\pi_2$ where $(x,y)\in\{(1,0),(1,1),(2,0),(2,1),(3,0),(4,0)\}$, if $m$ is odd. So suppose $m$ is even. Hence, $\pi\in H'$ with $m+3$ values denoted by $d_i$ satisfies Equation \ref{eqn:5.6}$(b)$. We need to show that for $\pi\in H'$ of length $n$ with $m+4$ values denoted by $d_i$ satisfies Equation \ref{eqn:5.6}$(a)$. The proof can be done by the inverse argument of the proof for the case that $m$ is odd. This completes the proof of Proposition \ref{prop:5.8}.\hfill$\blacksquare$\\
\subsection{Proof of Theorem \ref{thm:5.1} (Part 2)}
We now prove the converse direction of Theorem \ref{thm:5.1} as Proposition \ref{prop:5.9}, which states that an arbitrary permutation which has the form of either Equation \ref{eqn:5.6}$(a)$ or \ref{eqn:5.6}$(b)$ must be in $H'$. Before we start, we discuss a few differences between this proposition and Proposition \ref{prop:4.5}, the analogous proposition in Chapter \ref{chap:4}. Then we state and prove Lemma \ref{lem:5.12}, which is necessary for the proof of Proposition \ref{prop:5.9}.\\

In Proposition \ref{prop:4.5}, we assume \textit{a priori} that the permutation $\pi$ with the structure in Figure \ref{fig:4.7} is simple. Here, we will not make that assumption. This is because in Chapter \ref{chap:6}, we will not show that this assumption holds for a permutation decoded from a word in our language. Hence, we need to explicitly prove that $\pi$ of the form in Equation \ref{eqn:5.6}$(a)$ or \ref{eqn:5.6}$(b)$ is simple in addition to showing that $\pi$ avoids every permutation in $\{52341,53241,52431,35142,42513,351624\}$.\\

Furthermore, we proved Proposition \ref{prop:4.5} using a graphical representation of permutations in $H$. In $\A'$, however, providing a graphical representation of $\pi$ in $H'$ is extremely complicated. Hence, we show Proposition \ref{prop:5.9} by induction on the number of glue sums.\\

We now establish the following lemma, which says that if $\pi$ has the form of Equation \ref{eqn:5.6}$(a)$ or \ref{eqn:5.6}$(b)$, then $\pi$ has exactly $m+3$ values denoted by $d_i$ where each $d_i$ is defined in the same way as in the proof of \ref{prop:5.8}.
\begin{lemma}\label{lem:5.12}
\textit{Let $\pi$ be an arbitrary permutation of length $n$. If there exist simple permutations $\s_i$ ($i$ odd) in $\A'$ of extreme pattern 2413 and simple permutations $\tau_i$ ($i$ even) in $\A'$ of extreme pattern 3142 such that\[
\pi=\s_1\nwsum_{x_1}^{y_1}\tau_2\sesum_{x_2}^{y_2}\s_3\nwsum_{x_3}^{y_3}\tau_4\sesum_{x_4}^{y_4}\cdots\sesum_{x_{m-1}}^{y_{m-1}}\s_m\qquad\textrm{ for some odd integer }m
\]
or\[
\pi=\s_1\nwsum_{x_1}^{y_1}\tau_2\sesum_{x_2}^{y_2}\s_3\nwsum_{x_3}^{y_3}\tau_4\sesum_{x_4}^{y_4}\cdots\nwsum_{x_{m-1}}^{y_{m-1}}\tau_m\qquad\textrm{ for some even integer }m
\]
where $(x_\l,y_\l)\in\{(1,0),(1,1),(2,0),(2,1),(3,0),(4,0)\}$ ($1\leq\l\leq m-1$), then $\pi$ has $m+3$ distinct values denoted by $d_i$ where $d_1=\pi(1)$ and\[
d_i=\left\{\begin{array}{cl}
\pi(\max\{s:\pi(s)<\pi(d_{i-1})\}) &\textrm{ if }i\textrm{ is even}\\
\max\{t:\pi^{-1}(t)<\pi^{-1}(d_{i-1})\} &\textrm{ if }i\textrm{ is odd}
\end{array}\right.
\]
for $i$ with $2\leq i\leq m+3$.}
\end{lemma}
\textit{Proof.} We use mathematical induction to prove the statement. As the base case of $m=1$, suppose $\pi=\s$ where $\s$ is a simple permutation in $\A'$ of extreme pattern 2413 whose length is $n$. Then $d_1=\pi(1)$, $d_2=1$, $d_3=n$ and $d_4=\pi(n)$. By definition, we obtain $d_{i}=d_3$ if $i\geq 5$ is odd, and $d_{i}=d_4$ if $i\geq 6$ is even. Thus, we have $4=m+3$ distinct values denoted by $d_i$, so the statement holds for the base case.\\

Suppose the statement holds for some odd integer $m$. Consider \[\pi=\s_1\nwsum_{x_1}^{y_1}\tau_2\sesum_{x_2}^{y_2}\cdots\sesum_{x_{m-1}}^{y_{m-1}}\s_m\] of length $n$ for some $\s_i$ ($1\leq i\leq m$, $i$ odd) and $\tau_i$ ($2\leq i\leq m-1$, $i$ even) where $(x_\l,y_\l)\in\{(1,0),(1,1),(2,0),(2,1),(3,0),(4,0)\}$ ($1\leq\l\leq m-1$). Assume $\s_m$ has the structure of Case 1 of Figure \ref{fig:5.22}. Let $\tau_{m+1}$ be a simple permutation ($|\tau_{m+1}|=k$) in $\A'$ of extreme pattern 3142 which has the structure described in Case 1 of Figure \ref{fig:5.23}. We let $d_i$ ($1\leq i\leq m+3$) for $\pi$ and $d'_i$ for $\pi\nwsum_1^0\tau_{m+1}$ as defined by $d_i$ in the statement of Lemma \ref{lem:5.12}. Note that $d_{m+2}=n$ and $d_{m+3}=\pi(n)$ by definition. We need to show there are $m+4$ distinct values denoted by $d'_i$. The only values of $\pi$ that are modified to construct $\pi\nwsum_1^0\tau_{m+1}$ are $d_{m+2}=\pi(n)$ and $d_{m+3}=\pi(n)$. The value in $\pi\nwsum_1^0\tau_{m+1}$ corresponding to $d_{m+2}$ is adjusted upward by $(\tau_{m+1}(1)-3)$, and the value $d_{m+3}$ is eliminated.  Moreover, every value in $\pi\nwsum_1^0\tau_{m+1}$ corresponding to $\tau_{m+1}$ is greater than or equal to $d_{m+3}$. Hence, $d'_i=d_i$ for $1\leq i\leq m+1$. Also, the value in $\pi\nwsum_1^0\tau_{m+1}$ corresponding to $d_{m+2}$ is still the greatest value whose position is to the left of $d'_{m+1}$, so it is denoted by $d'_{m+2}$, though it is not equal to $d_{m+2}$. Since $\tau_{m+1}(k)<\tau_{m+1}(1)$, $\pi\nwsum_1^0\tau_{m+1}(n+k-3)$, which is the value corresponding to $\tau_{m+1}(k)$, is denoted by $d'_{m+3}$, and the one corresponding to $k$ is denoted by $d'_{m+4}$. We then have $d'_{i}=d_{m+4}$ if $i\geq m+6$ is odd, and $d'_{i}=d_{m+3}$ if $i\geq m+5$ is even. Consequently, we have $m+4$ distinct values denoted by $d'_i$.\\

Proofs for for the cases that $\s_m$ and $\tau_{m+1}$ having different structures, and thus, different types of NW glue sums, are very similar. We can also show that the statement is true for $m+1$ when we assume $m$ is even by applying the inverse argument. Altogether, we have the desired result.\hfill$\blacksquare$\\

We are now ready to state and prove Proposition \ref{prop:5.9}.
\begin{proposition}\label{prop:5.9}
\textit{Let $\pi$ be an arbitrary permutation of length $n$. If there exist simple permutations $\s_i$ ($i$ odd) in $\A'$ of extreme pattern 2413 and simple permutations $\tau_i$ ($i$ even) in $\A'$ of extreme pattern 3142 such that\[
\pi=\s_1\nwsum_{x_1}^{y_1}\tau_2\sesum_{x_2}^{y_2}\s_3\nwsum_{x_3}^{y_3}\tau_4\sesum_{x_4}^{y_4}\cdots\sesum_{x_{m-1}}^{y_{m-1}}\s_m\qquad\textrm{ for some odd integer }m
\]
or\[
\pi=\s_1\nwsum_{x_1}^{y_1}\tau_2\sesum_{x_2}^{y_2}\s_3\nwsum_{x_3}^{y_3}\tau_4\sesum_{x_4}^{y_4}\cdots\nwsum_{x_{m-1}}^{y_{m-1}}\tau_m\qquad\textrm{ for some even integer }m
\]
where $(x_\l,y_\l)\in\{(1,0),(1,1),(2,0),(2,1),(3,0),(4,0)\}$ ($1\leq\l\leq m-1$), then $\pi$ is in $H'$.}
\end{proposition}
\textit{Proof.} If a permutation is expressed as in Equation \ref{eqn:5.6}$(a)$ or \ref{eqn:5.6}$(b)$, then it satisfies the condition of $\pi^{-1}(1)\geq 5$ or $\pi(1)=2$ because the first summand is $\s_1$ which is a simple permutation of extreme pattern 2413. So we have to show that a permutation expressed as in Equation \ref{eqn:5.6}$(a)$ or \ref{eqn:5.6}$(b)$ is simple and avoids every permutation in $\{52341,53241,52431,35142,42513,351624\}$. We use mathematical induction on $m$ to complete the proof. As the base case, suppose $\pi=\s$ where $\s$ is a simple permutation in $\A'$ of extreme pattern 2413. Then we have the desired result.\\

Now, assume that the statement holds for some odd integer $m$. Let\[
\pi=\s_1\nwsum_{x_1}^{y_1}\tau_2\sesum_{x_2}^{y_2}\cdots\sesum_{x_{m-1}}^{y_{m-1}}\s_m\] of length $n$ for some $\s_i$ ($1\leq i\leq m$, $i$ odd) and $\tau_i$ ($2\leq i\leq m-1$, $i$ even) where $(x_\l,y_\l)\in\{(1,0),(1,1),(2,0),(2,1),(3,0),(4,0)\}$ ($1\leq\l\leq m-1$). We need to show that, for every appropriate $\tau_{m+1}$ of length $k$ and $\s_{m}$ with $\nwsum_{x_m}^{y_m}$ ($(x_m,y_m)\in\{(1,0),(1,1),(2,0),(2,1),(3,0),(4,0)\}$), $\pi\nwsum_{x_m}^{y_m}\tau_{m+1}$ is simple, and avoids permutations in $\{52341,53241,52431,35142,42513,351624\}$. Let $\rho=\pi\nwsum_{x_m}^{y_m}\tau_{m+1}$. We first show that $\rho$ is simple.\\

Define values $d_j$ of $\rho$ for positive integers $i$ as before. By Lemma \ref{lem:5.12}, there are $m+4$ values denoted by $d_{i}$ of $\rho$. Suppose to the contrary that $\rho$ is not simple, so we have a proper non-singleton block $I$. Assume $I$ contains at least two positions corresponding to $d_j$ and $d_k$ for some $j$, $k$ with $1\leq j,k\leq m+4$. Note that $I$ cannot contain both positions of $d_1$ and $d_{m+3}$ because this results in $I=[1,n]$. So suppose $I$ does not contain $1=\rho^{-1}(d_1)$. Let $s$ be the least position in $I$ such that $\rho(s)=d_j$ for some $j$ with $2\leq j\leq m+4$. Suppose $j$ is even. If $j=m+3$, then there is no more value to the right of $d_{m+3}$, so $j\leq m+1$. Hence, assume $j\leq m+1$. Then the position of $d_{j+3}$ must be in $I$, but with $d_{j}<d_{j-1}<d_{j+3}$, we must have $\rho^{-1}(d_{j-1})\in I$ while $\rho^{-1}(d_{j-1})<s$, so we have a contradiction. Now, assume $j$ is odd. Then $I$ contains the position of $d_{j-1}$. Because $d_{j-1}<d_{j-2}<d_{j}$, we have $\rho^{-1}(d_{j-2})\in I$, but $\rho^{-1}(d_{j-2})<s$, so we again achieve a contradiction. We can apply a similar argument for the case of $I$ not containing $n=\rho^{-1}(d_{m+3})$, so $I$ cannot contain two or more positions corresponding to values denoted by $d_i$.\\

This time, suppose $I$ does not contain a position corresponding to any of $d_j$ of $\rho$. For every $(x_m,y_m)\in\{(1,0),(1,1),(2,0),(2,1),(3,0),(4,0)\}$, positions of all values of $\pi$ that are modified by $\nwsum_{x_m}^{y_m}$ to construct $\rho$ are greater than or equal to $\rho^{-1}(d_{m+2})$. Thus, if $I$ is a subset of $[1,\rho^{-1}(d_{m+2})-1]$, then $I$ is also a proper non-singleton block of $\pi$. Since $\pi$ is a simple permutation by induction hypothesis, this is a contradiction. Similarly, if $I$ is a subset of $[\rho^{-1}(d_{m+4})+1,\rho^{-1}(d_{m+3})-1]$, then $I$ is a proper non-singleton block of $\tau_{m+1}$. Hence, we have either $I\sbst[\rho^{-1}(d_{m+2})+1,\rho^{-1}(d_{m+1})-1]$ or $I\sbst[\rho^{-1}(d_{m+1})+1,\rho^{-1}(d_{m+4})-1]$. Suppose $I\sbst[\rho^{-1}(d_{m+2})+1,\rho^{-1}(d_{m+1})-1]$. Then $I$ must contain at least one value in $[\rho^{-1}(d_{m+2})+1,\rho^{-1}(d_{m+1})-1]$ corresponding to the value of $\pi$ that was modified by $\nwsum_{x_m}^{y_m}$, because, otherwise, $I$ would also a proper non-singleton block of $\pi$. In cases of $\nwsum_{1}^{0}$, $\nwsum_{1}^{1}$ and $\nwsum_{3}^{0}$, there is no such value. If $\nwsum_{x_m}^{y_m}=\nwsum_{2}^{0}$, $\nwsum_{2}^{1}$ or $\nwsum_{4}^{0}$, either $d_{m+2}-1$ or $d_{m+2}-2$ is the only such value, but in either case, $d_{m+3}$ is less than this value, implying that $I$ contains $\rho^{-1}(d_{m+3})$, so we achieve a contradiction for every type of NW glue sum.\\

So assume $I\sbst[\rho^{-1}(d_{m+1})+1,\rho^{-1}(d_{m+4})-1]$. We first consider the case where $(x_m,y_m)\in\{(1,0),(1,1),(2,0),(2,1)\}$. If $I$ only contains values corresponding to $\tau_{m+1}$, then $I$ would be a proper non-singleton block of $\tau_{m+1}$ as well, because for every $(x_m,y_m)\in\{(1,0),(1,1),(2,0),\allowbreak(2,1)\}$, $\nwsum_{x_m}^{y_m}$ only shifts values of $\tau_{m+1}$ upward by the same amount to construct $\rho$. For a similar reason, if $I$ only contains values corresponding to $\pi$, then $I$ would be a proper non-singleton block of $\pi$. Thus, $I$ must contain values corresponding to both $\pi$ and $\tau_{m+1}$. Because $I$ contains values from both $\pi$ and $\tau_{m+1}$, the position of $\tau'_{m+1}(3)$, the left-most point of $\tau_{m+1}$ in $[\rho^{-1}(d_{m+1})+1,\rho^{-1}(d_{m+4})-1]$, must be in $I$ since every value to its right comes from $\tau$. If $\tau'_{m+1}(3)>d_{m+3}$, then $\rho^{-1}(d_{m+3})$ would be in $I$ since every value of $\rho$ coming from $\pi$ in $[\rho^{-1}(d_{m+1})+1,\rho^{-1}(d_{m+4})-1]$ is less than $d_{m+3}$, so we must have $\tau'_{m+1}(3)<d_{m+3}$. By Proposition \ref{prop:5.5}, this implies that $\tau_{m+1}(3)$ is a part of a 312-position chain. Due to the structure of 312-position chain, every value of $\rho$ corresponding to some value in this 312-position chain must be in $I$. Denote by $M$ the position of the right-most value involved in the 312-position chain containing $\tau_{m+1}(3)$. Then $\tau_{m+1}(M-1)>d_{m+3}$, since $\tau_{m+1}(M-1)$ is the scissor of this 312-position chain. Hence, $\rho^{-1}(d_{m+3})\in I$, which is a contradiction.\\

For the case of $\nwsum_4^0$, $I$ must also contain values from both $\pi$ and $\tau_{m+1}$. The position of $\tau'_{m+1}(4)$, the left-most point of $\tau_{m+1}$ in $[\rho^{-1}(d_{m+1})+1,\rho^{-1}(d_{m+4})-1]$, must be in $I$ since every value to its right comes from $\tau$. However, by proposition \ref{prop:5.5}, $\tau'_{m+1}(4)>d_{m+3}$. Hence, $\rho^{-1}(d_{m+3})$ would be in $I$ since every value of $\rho$ coming from $\pi$ in $[\rho^{-1}(d_{m+1})+1,\rho^{-1}(d_{m+4})-1]$ is less than $d_{m+3}$.\\

Finally, consider the case of $\nwsum_3^0$. The block $I$ cannot contain only values corresponding to $\tau_{m+1}$ due to the same reason as other NW glue sums. Now, if $I$ only contains values corresponding to $\pi$, then it must contain $\pi'(n-2)$, because it is the only value in $[\rho^{-1}(d_{m+1})+1,\rho^{-1}(d_{m+4})-1]$ that is shifted upward by $\nwsum_3^0$ to construct $\rho$. However, since $\pi'(n-2)>d_{m+3}$, we end up with $\rho^{-1}(d_{m+3})\in I$ as usual, so $I$ cannot contain only values corresponding to $\pi$. Thus, $I$ must contain values corresponding to both $\pi$ and $\tau_{m+1}$. Furthermore, $\rho^{-1}(\pi'(n-2))\notin I$, because, otherwise, we again obtain $\rho^{-1}(d_{m+3})\in I$. So the only possibility is $I=[\rho^{-1}(\pi'(n-1)),z]$ for some position $z$ whose value corresponds to a value of $\tau_{m+1}$. Note that, however, $\pi$ has the structure shown in Case 3 of Figure \ref{fig:5.22}. Thus, $\pi(n-2)=n-1$, $\pi(n)=n-2$ and $\pi(n-1)=n-4$ are all parts of the same 312-value chain, which implies that $\pi(n-1)+1$ is located to the left of $\pi(n-1)$. Consequently, $\pi'(n-1)+1\in\rho(I)$ has a position less than $\rho^{-1}(\pi'(n-1))$, so we achieve a contradiction. Altogether, we have a contradiction for the case of $I$ not containing a position corresponding to any of $d_j$ of $\rho$.\\

Lastly, we assume $I$ contains exactly one position corresponding to $d_j$ for some $j$ with $1\leq j\leq m+3$. For every $(x_m,y_m)\in\{(1,0),(1,1),(2,0),(2,1),(3,0),(4,0)\}$, values of $\pi$ that are modified by $\nwsum_{x_m}^{y_m}$ to construct $\rho$ are all located to the right of $d_{m+2}$. Thus, if $j\leq m$, then $I$ would be a proper non-singleton block in $\pi$. Similarly, if $j=m+3$ or $m+4$, then $I$ would be a proper non-singleton block in $\tau_{m+1}$. So $j$ must be either $m+1$ or $m+2$. For the case $j=m+1$, the arguments of previous paragraphs apply. So suppose $j=m+2$. If $I$ contains $\rho^{-1}(d_{m+2})-1$, then $I$ also has to contain $\rho^{-1}(d_{m+3})$ because $\rho(\rho^{-1}(d_{m+2})-1)<d_{m+3}$. Thus, $I=[\rho^{-1}(d_{m+2}),z]$ for some position $z$. The value immediately to the right of $d_{m+2}$, however, is less than $d_{m+3}$, except when it is $d_{m+2}-2$ and $\nwsum_{x_m}^{y_m}=\nwsum_{4}^{0}$. In this case, $\tau_{m+1}'(4)$, the left-most value from $\tau_{m+1}$ is $d_{m+2}-1$, so $\rho^{-1}(d_{m+2}-1)\in I$. Since $\rho^{-1}(d_{m+2}-1)>\rho^{-1}(d_{m+1})$, $\rho^{-1}(d_{m+1})\in I$, which is a contradiction.\\

Consequently, $\rho$ cannot have a proper non-singleton block $I$, so $\rho$ is simple. Next, we prove that $\rho$ avoids every permutation in $\{52341,53241,52431,35142,42513,351624\}$.\\

Suppose $\rho$ contains some permutation $\b$ in $\{52341,53241,52431,35142,42513,351624\}$. Note that $\pi,\tau_{m+1}\in\A'$, and modifications of every value of $\pi$ and $\tau_{m+1}$ by each $\nwsum_{x_m}^{y_m}$ to construct $\rho$ are either upward shift of certain values by the same amount or removal. Thus, whichever permutation $\b$ in the basis that $\rho$ contains, the containment must involve both values corresponding to $\pi$ and values corresponding to $\tau_{m+1}$, because, otherwise, $\b\cont\pi$ or $\b\cont\tau_{m+1}$. Hence, involving both values from $\pi$ and values from $\tau_{m+1}$ for a containment of $\b$ is a necessary condition of $\b\cont\rho$, but it is not a sufficient condition.\\

We first discuss the case $\b\in\{52341,53241,52431\}$. Note that, every value corresponding to $\tau_{m+1}$ is greater than any value corresponding to $\pi$, except for the shifted ones. Thus, for every $\b\in\{52341,53241,52431\}$, we must assign one of the shifted values to $\b(1)=5$ in order to assign some value $\tau'_{m+1}(i)$ to $\b(5)=1$. Hence, the only values which can play the role of 5 are $d_{m+2}$ for every $(x_m,y_m)\in\{(1,0),(1,1),(2,0),(2,1),(3,0),(4,0)\}$, $d_{m+2}-1$ for $(x_m,y_m)=(2,0)$ and $(2,1)$, $\rho(n-2)=\pi'(n-2)$ for $(x_m,y_m)=(3,0)$, and $d_{m+2}-2$ for $(x_m,y_m)=(4,0)$.\\

Consider the cases of $(x_m,y_m)\in\{(1,0),(1,1)\}$. Then $d_{m+2}$ is playing the role of $\b(1)=5$. In addition, the values playing the roles of $\b(2)$, $\b(3)$ and $\b(4)$ are all from $\tau_{m+1}$, since we have to assign a value corresponding to $\tau_{m+1}$ to $\b(5)=1$. However, note that the flattening of $d_{m+2}\tau'_{m+1}(3)\tau'_{m+1}(4)\cdots\tau'_{m+1}(k)$, call it $\gamma$, is contained in $\tau_{m+1}$. Hence, no matter which values of $\rho$ corresponding to $\tau_{m+1}$ plays the role of $\b(2)$, $\b(3)$ and $\b(4)$, $\b\cont\gamma$ implies $\b\cont\tau_{m+1}$. Since $\tau_{m+1}\in\A'$, this is a contradiction.\\

Next, consider $(x_m,y_m)\in\{(2,0),(2,1)\}$. As explained, at least one shifted value must be playing the role of 5, so in this case, it is either $d_{m+2}$ or $d_{m+2}-1$. Suppose only one of them is involved for a containment of $\b$, say $d_{m+2}$. Then again, $\gamma$, the flattening of $d_{m+2}\tau'_{m+1}(3)\tau'_{m+1}(4)\cdots\allowbreak\tau'_{m+1}(k)$ is contained in $\tau_{m+1}$, and all four values playing the roles of $\b(2)$, $\b(3)$, $\b(4)$ and $\b(5)=1$ respectively must correspond to $\tau_{m+1}$, so $\b\cont\tau_{m+1}$, a contradiction. We achieve the exact same result if $d_{m+2}-1$ plays the role of 5. Thus, suppose both $d_{m+2}$ and $d_{m+2}-1$ are involved for a containment of $\b$. Since $d_{m+2}$ is located to the left of $d_{m+2}-1$, $d_{m+2}$ must play the role of $\b(1)=5$. In addition, we cannot assign $d_{m+2}-1$ to either 2 or 3, because there is no integer value between $d_{m+2}-1$ and $d_{m+2}$, so we don't have a value for 4. Since every $\b$ in $\{52341,53241,52431\}$ has either $\b(2)=2$ or $\b(2)=3$, we cannot involve both $d_{m+2}$ and $d_{m+2}-1$ to have $\b\cont\rho$. Consequently, $\rho=\pi\nwsum_2^0\tau_{m+1}$ cannot contain $\b\in\{52341,53241,52431\}$.\\

We move onto the case $(x_m,y_m)=(3,0)$. Let $s$ be the value corresponding to the scissor of the last 312-value chain of $\pi$ which is merged with the first 312-position chain of $\tau$ by the type 3-0 NW glue sum. In addition, let $t$ be the position of the value corresponding to the left-most value of the last 312-value chain of $\pi$. Note that $t=\rho^{-1}(s-1)-1$. Due to the construction of $\rho$ by $\nwsum_3^0$, the flattening of $d_{m+2}s d_{m+1}\rho(t)\rho(t+1)\cdots\rho(n+k-(\l+3))$ is exactly $\tau_{m+1}$. Let $\overline{\tau}_{m+1}$ be this sequence of values. A containment of $\b$ must involve some value other than the values in $\overline{\tau}_{m+1}$, because, otherwise, $\b\cont\tau_{m+1}$. Whether $d_{m+2}$ or $\rho(n-2)$ plays the role of 5, some value $z$ corresponding to $\tau_{m+1}$ must be assigned to 1. Since every value greater than $z$ is in the sequence $\overline{\tau}_{m+1}$, we cannot assign any value to $\b(2)$, $\b(3)$ and $\b(4)$ that are not in $\overline{\tau}_{m+1}$. Hence, $\b\avd\rho$.\\

Lastly, the case $(x_m,y_m)=(4,0)$ is similar to the case of $(x_m,y_m)\in\{(2,0),(2,1)\}$. If $d_{m+2}$ is involved in a containment of $\b$ as 5, but not $d_{m+2}-1$, then again, we are forced to assign some values of $\tau_{m+1}$ to $\b(2)$, $\b(3)$ and $\b(4)$. Because the flattening of $d_{m+2}\tau'_{m+1}(4)\tau'_{m+1}(5)\cdots\tau'_{m+1}(k)$ is contained in $\tau_{m+1}$, we obtain $\b\cont\tau_{m+1}$, which is a contradiction. We achieve the same result for the case of $d_{m+2}-2$ playing the role of 5. So suppose both $d_{m+2}$ and $d_{m+2}-2$ are involved in a containment of $\b$. Hence, $d_{m+2}$ and $d_{m+2}-2$ play the role of $\b(1)=5$ and $\b(2)$ respectively. Note that we only have one value, $d_{m+2}-1$, that is in between $d_{m+2}-2$ and $d_{m+2}$. Thus, if $\b=52341$ or $52431$, we do not have a value for either $\b(2)$ or $\b(3)$. Therefore, suppose $\b=53241$. Since the value assigned to $\b(3)=2$ is less than $d_{m+2}-2$ and its is in between $\rho^{-1}(d_{m+2}-2)$ and $\rho^{-1}(d_{m+2}-1)$, this value must correspond to $\pi$. However, this value is less than the value corresponding to $\tau_{m+1}$ that is assigned to 1, so we cannot form 53241 pattern. Hence, $53241\avd\rho$.\\

Altogether, we conclude that for every $\b\in\{52341,53241,52431\}$, $\rho$ avoids $\b$. We now look at the case of $\b=35142$.\\

First, consider $(x_m,y_m)\in\{(1,0),(1,1)\}$. If we assign a value of $\rho$ from $\pi$ to $\b(4)=4$, then we are forced to choose another value from $\pi$ to play the role of $\b(5)=2$. Since we don't have a value from $\tau_{m+1}$ corresponding to some value of $\b$, we have $\b\cont\pi$, which is a contradiction. Hence, suppose some value corresponding to $\tau_{m+1}$ plays the role of 4. If a value corresponding to $\tau_{m+1}$ plays the role of $\b(2)=5$, then $\b(3)=1$ is automatically played by a value corresponding to $\tau_{m+1}$ as well. We need to assign a value from $\pi$ to $b(1)=3$, but the only value available is $d_{m+2}$. Since $\gamma$, the flattening of $d_{m+2}\tau'_{m+1}(3)\tau'_{m+1}(4)\cdots\tau'_{m+1}(k)$ is contained in $\tau_{m+1}$, we now have $\b\cont\tau_{m+1}$, a contradiction. Thus, suppose some value corresponding to $\pi$ plays the role of 5. Then again, $d_{m+2}$ is the only value which can play the role of 5. We are forced to assign a value $z$ corresponding to $\pi$ and located to the left of $d_{m+2}$ to 3. However, there is no value from $\tau_{m+1}$ that is less than $z$, so we do not have a value for 2. Hence, 35142 cannot be contained in $\rho$ in the case of $(x_m,y_m)\in\{(1,0),(1,1)\}$.\\

For $(x_m,y_m)\in\{(2,0),(2,1),(4,0)\}$, it is possible to assign $d_{m+2}-1$ (when $x_m=2$) or $d_{m+2}-2$ (when $x_m=4$) to 4, and assign some value $z$ from $\tau_{m+1}$ to 2. Since the only value located to the left of and greater than the one assigned to 4 is $d_{m+2}$, it has to play the role of 5, forcing the position of the value playing the role of 3 to be to the left of $d_{m+2}$. However, there is no value located there that is greater than $z$, so $\b$ cannot be contained in $\rho$ this way. The rest of the proof and the case $(x_m,y_m)=(3,0)$ is very similar to the case of $(x_m,y_m)\in\{(1,0),(1,1)\}$. Essentially, 5 must be played by either $d_{m+2}$ or the other shifted value. Then the value playing the role of 3 is located to the left of the value for 5, but this value cannot be greater than any value corresponding to $\tau_{m+1}$, so we don't have a value to assign to 2. Thus, in every case, $35142\avd\rho$.\\

We now let $\b=42513$. Again, since both values corresponding to $\pi$ and values corresponding to $\tau_{m+1}$ must be involved in a containment of $\b$, we know the value for $\b(1)=4$ is from $\pi$ and the value for $\b(5)=3$ is from $\tau_{m+1}$. For each NW glue sum, there are at most two values corresponding to $\pi$ that are greater than some value corresponding to $\tau$, but they are in decreasing order from left to right. Hence, the value playing the role of 4 must be one of the two shifted values, and the value playing the role of $\b(3)=5$ must be from $\tau_{m+1}$ for each case of NW glue sum. This implies the value for $\b(4)=1$ is also corresponding to $\tau_{m+1}$. If the value for $\b(2)=2$ is from $\tau_{m+1}$, then again, the flattening of values assigned to a value of $\b$ is contained in $\tau_{m+1}$. Thus, the value assigned to $2$ must be from $\pi$. Furthermore, this value has to be the other shifted value to assign, and particularly $d_{m+2}-2$. Otherwise, we do not have values from $\tau_{m+1}$ that we can assign to each of 1 and 3. Therefore, $(x_m,y_m)\in\{(3,0),(4,0)\}$, and we assign $d_{m+2}$ to 4, $d_{m+2}-2$ to 2 and three values from $\tau_{m+1}$ to each of 5, 1 and 3. However, because the flattening of values $d_{m+2}$, $d_{m+2}-2$ and all the values corresponding to $\tau_{m+1}$ is contained in $\tau_{m+1}$, we must have $\b\cont\tau_{m+1}$. This is a contradiction, so we conclude that $42513\avd\rho$.\\

The final case is $\b=351624$. For the same reason as the case of $42513$, the value $\b(2)=5$ must be played by one of the shifted values, and the value $\b(4)=6$ must be from $\tau_{m+1}$ for each NW glue sum. Thus, the value playing the role of $\b(1)=3$ is a value corresponding to $\pi$, and this value is not shifted by $\nwsum_{x_m}^{y_m}$. Because the value playing the role of 6 is from $\tau_{m+1}$, so is the value playing the role of $\b(5)=2$, but we do not have a value corresponding to $\tau_{m+1}$ that is less than the value playing the role of 3. Consequently, $351624\avd\rho$.\\

Therefore, when the statement of the proposition holds for odd $m$, it is also true for $m+1$. Showing the case of even $m$ is the same argument applied to the inverses of all permutations involved. Altogether, we have the desired result.\hfill$\blacksquare$
\newpage
\chapter{Enumeration of the class $\mathcal{A}'$}\label{chap:6}
\section{Enumeration of simple permutations in $\A'$}
With the structure we discussed in Chapter \ref{chap:5}, we are finally ready to enumerate the class $\A'$. How we are going to accomplish this is somewhat similar to the method we used in Chapter \ref{chap:4} for $\A$, but we need to break down the whole procedure into small steps. We first give an alphabet $\Sigma'$ and define an encoding function $\phi'$ from $H'$ to $\Sigma'^*$. Then we define a language $L'\subseteq\Sigma'^*$ and show $\phi'$ is a bijection between $H'$ and $L'$. Afterwards, we construct ten languages $\overline{L}'_i$ ($1\leq i\leq 10$) associated with $L'$. Frankly, there are ten distinct kinds of prefixes a word $w$ in $\overline{L}$ can have. Thus, each $\overline{L}'_i$ is prepared to generate the number of words having different prefixes. We then define ten deterministic finite-state automatons $M'_i$ ($1\leq i\leq 10$) and show $\overline{L}'_i=\mathcal{L}(M_i)$ for each $i$ with $1\leq i\leq 10$. Once we obtain ten distinct generating functions, we combine them to obtain the following result.\\

\begin{theorem}\label{thm:6.1}
\textit{Let $f_{\textrm{Si}(\A')\setminus\mathcal{S}_2}$ be the generating function for the set of simple permutations in $\A'$ excluding $\S_2=\{12,21\}$. Then}\\
\[
f_{\textrm{Si}(\A')\setminus\mathcal{S}_2}=\frac{2x^4 \left(x^{10}+7 x^9+18 x^8+23 x^7+16 x^6+10 x^5+12 x^4+9 x^3+2 x^2+1\right)}{(x+1) \left(2 x^9+12 x^8+16 x^7+3 x^6-11 x^5-5 x^4-3 x^2-3
   x+1\right)}.
\]
\end{theorem}\vhhh

In this section, let $N$ and $S$ be the set of simple permutations of extreme pattern 2413 and 3142 respectively. Let $m$ be the total number of simple permutations in $N$ and $S$ together to construct $\pi\in H'$ with glue sums. Furthermore, let $d_i$, $p_i$, $q_i$, $q'_i$ and $r_i$ be defined as in the proof of Proposition \ref{prop:5.8}. Note that each of $p_i$, $q_i$, $q'_i$ and $r_i$ is defined differently based on whether $i$ is odd or even.
\subsection{Defining the encoding function $\phi'$ and the language $L'$}
Let\[
\Sigma'=\{a,a',a'',b,b_s,b',b'',\underline{b},\overline{b},c,c',c'',d,d_a',d_a'',d_c',d_c'',\underline{d},\overline{d},\underline{\overline{d}},d_\l,x,x',x'',\underline{x},y,y',y'',\overline{y},z\}.
\]
Now, we define an encoding function $\phi'$ from $H'$ to $\Sigma'^*$. Comparing to the encoding function $\phi$ for $H$, the description of $\phi'$ is far more complicated, so we divide into several algorithms to define $\phi'$.\\

First, by Theorem \ref{thm:5.1}, given $\pi\in H'$, $\pi$ can be uniquely expressed as Equation \ref{eqn:5.6} where $\s_i$ ($i$ odd) is in $N$ and $\tau_i$ ($i$ even) is in $S$.\\

Next, we define two encoding algorithms. One is to encode a simple permutation in $N$, called \textsf{N-ENCODE}, and the other one is for a simple permutation in $S$, called \textsf{S-ENCODE}. Both algorithms encode simple permutations in $\A'$ into a word in $\{a,a',a'',b,b_s,b',b'',c,c',c''\}^*$. So let $\Sigma_1=\{a,a',a'',b,b_s,b',b'',c,c',c''\}$ which is a subset of $\Sigma'$.\\

We first define \textsf{N-ENCODE}. The algorithm \textsf{N-ENCODE} assigns each value of $\s$ a letter in $\Sigma_1$, then reads the letters assigned to each value from 1 to $n$ in increasing order. Letters are assigned as follows.\\

First, suppose the position of the value $t$ is in $[1,\pi^{-1}(n))$. Then, if $t$ plays the role of 1 of 21 or 1 of 231 of a 231-value chain, it is encoded as $a'$, if $t$ plays the role of 2 of 21 or 3 of 231 of a 231-value chain, it is encoded as $a''$, and if $t$ plays the role of 2 of 231 of a 231-value chain or 1 of 1 in Equation \ref{eqn:5.2}, then it is encoded as $a$. Next, assume the position of the value $t$ is in $[\pi^{-1}(n),\pi^{-1}(1)]$. Then, if $t$ plays the role of 1 of 12 in Equation \ref{eqn:5.1}, it is encoded as $b'$, if $t$ plays the role of 2 of 12 in Equation \ref{eqn:5.1}, it is encoded as $b''$, if $t$ is a scissor of a 231-value chain or a 312-value chain (and it is not already $b'$ or $b''$), it is encoded as $b_s$, and otherwise, it is encoded as $b$. Finally, suppose the position of the value $t$ is in $(\pi^{-1}(1),n]$. Then, if $t$ plays the role of 1 of 21 or 1 of 312 of a 312-value chain, it is encoded as $c'$, if $t$ plays the role of 2 of 21 or 3 of 312 of a 312-value chain, it is encoded as $c''$, and if $t$ plays the role of 2 of 312 of a 312-value chain or 1 of 1 in Equation \ref{eqn:5.2}, then it is encoded as $c$.\\

Now, let us define the language $K_1\sbst\Sigma_1^*$ to have the following conditions, which will turn out to be the image of \textsf{N-ENCODE}.\begin{enumerate}
\item Prefix condition.\\
      A word $w$ must begin with $ba$, $bb'a$, $ba'a$, $ba'b_sa''$, $ba'b'a''$ or $ba'a'a''$.
\item Suffix condition.\\
      A word $w$ must end with $cb$, $cb''b$, $cc''b$, $c'b_sc''b$, $c'b''c''b$ or $c'c''c''b$.
\item Conditions on $a'$ and $a''$.\\
      Every $a'$ and $a''$ in $w$ is a part of a subword in $\{a',a'a\}\{a'a'',a'a''a\}^*\{b_sa'',b'a''\}$. Note the number of $a'$ and the number of $a''$ are equal in each sequence, and thus, in $w$.
\item Conditions on $c'$ and $c''$.\\
      Every $c'$ and $c''$ in $w$ is a part of a subword in $\{c'b_s,c'b''\}\{c'c'',cc'c''\}^*\{c'',cc''\}$. Note the number of $c'$ and the number of $c''$ are equal in each sequence, and thus, in $w$.
\item Conditions on $b_s$.\\
      The letter $b_s$ is only allowed in the third condition and the fourth condition, \textit{i.e.} as a part of a subword of $a'$ and $a''$ or a subword of $c'$ and $c''$.
\item Conditions on $b'$ and $b''$.\\
      Every $b'$ and $b''$ in $w$ is a part of a subword in $\{b'\}\{a'',\lambda\}\{a,\lambda\}\{c,\lambda\}\{c',\lambda\}\{b''\}$ with at least one of $a''$, $a$, $c$ or $c'$ being present. Note the number of $b'$ and the number of $b''$ are equal in each sequence, and thus, in $w$.
\item Repetition restrictions.\\
      $w$ must not contain $aa$, $bb$ or $cc$.
\end{enumerate}
\vhhh

We define \textsf{N-DECODE} on $K_1$, the inverse algorithm of \textsf{N-ENCODE}, and show that $N$ and $K_1$ are bijective due to these two algorithms. The algorithm \textsf{N-DECODE} takes $w\in K_1$ as an input and outputs a permutation as the following.\\
\\
\textbf{Algorithm} \textsf{N-DECODE}\\
\textsf{INPUT}: A word $w$ in $K_1$.\\
\textsf{OUTPUT}: A permutation $\pi$.\\
\begin{enumerate}[\hspace{7pt}\sffamily {Case} 1:]
\item[\sffamily Initialize:] \textsf{Draw a point for the first $b$ and locate coordinates for $P_a$, $P_b$ and $P_c$.}\\
      Draw a point at $(1,1)$. Let $P_a=(0,1)$ and $P_b=P_c=(1,1)$. Let $t=1$. Let $\a$ be the second letter in $w$.
\item[\sffamily Identify:] \textsf{Determine the case due to what letter $\a$ is.}\begin{enumerate}[\sffamily a.]
    \item If $\a=a'$, then \textsf{GOTO Case 1}.
    \item If $\a=b'$, then \textsf{GOTO Case 2}.
    \item If $\a=c'$, then \textsf{GOTO Case 3}.
    \item If $\a$ is none of the above, then \textsf{GOTO Case 4}.
    \end{enumerate}
\item \textsf{Construct and draw a 231-value chain in the following manner.}\\
    Take the subword of $w$ starting from $\a$ up to $a''$ which equalizes the numbers of $a'$ and $a''$. Hence, this subword is in $\{a',a'a\}\{a'a'',a'a''a\}^*\{b_sa'',b'a''\}$. Remove the first $a'$ and the first $a''$. Additionally, if there exists $a$ in between these $a'$ and $a''$, remove it as well. Iteratively remove $a'$, $a''$ and $a$ until $b_s$ or $b'$ is the only letter left. Every time we remove $a'$ and $a''$ together, sum 21 with $\oplus_1$, and every time we remove $a'$, $a''$ and $a$ together, sum 231 with $\oplus_1$. Draw the 231-value chain that we obtain from this procedure so that the chain is horizontally in between $P_a$ and $P_b$, and vertically greater than $t$. Set $t$ to be the greatest $y$ coordinate of all points in the chain and $P_a$ to be the right-most point in the chain.\\
    For $b_s$ or $b'$, draw a point located horizontally in between $P_a$ and $P_b$, and vertically in between the point with the second greatest $y$ coordinate and the point with the greatest $y$ coordinate in the chain. Set $\a$ to be the letter immediately after the subword.\begin{enumerate}[\sffamily a.]
    \item If $b_s$ is a part of the subword, then set $P_b$ to be the point drawn for $b_s$. Then \textsf{GOTO Identify}.
    \item If $b'$ is a part of the subword, then set $P_\l$ to be the point drawn for $b'$ and identify what letter $\a$ is. If $\a=a$, then \textsf{GOTO a} of \textsf{Case 2}. If $\a=c$, then \textsf{GOTO b} of \textsf{Case 2}. If $\a=b''$, then \textsf{GOTO c} of \textsf{Case 2}. If $\a=c'$, then \textsf{GOTO Case 3}.
    \end{enumerate}
\item \textsf{Draw points for $b'$ and letters up to $b''$.}\\
    Draw a point at $(x,y)$ where $P_a^{(x)}<x<P_b^{(x)}$ and $t<y$. Set $P_\l$ to be this new point and $t=P_\l^{(y)}$. Set $\a$ to be the next letter. If $\a=a$, then \textsf{GOTO a}. If $\a=c$, then \textsf{GOTO b}. If $\a=b''$, then \textsf{GOTO c}. If $\a=c'$, then \textsf{GOTO Case 3}.
    \begin{enumerate}[\sffamily a.]
    \item Draw a point at $(x,y)$ where $P_a^{(x)}<x<P_\l^{(x)}$ and $t<y$. Set $P_a$ to be this new point and $t=P_a^{(y)}$. Set $\a$ to be the next letter. If $\a=c$, then \textsf{GOTO b}. If $\a=b''$, then \textsf{GOTO c}. If $\a=c'$, then \textsf{GOTO Case 3}.
    \item Draw a point at $(x,y)$ where $P_c^{(x)}<x$ and $t<y$. Set $P_c$ to be this new point and $t=P_c^{(y)}$. Set $\a$ to be the next letter. If $\a=b''$, then \textsf{GOTO c}. If $\a=c'$, then \textsf{GOTO Case 3}.
    \item Draw a point at $(x,y)$ where $P_\l^{(x)}<x<P_b^{(x)}$ and $t<y$. Set $P_b=P_\l$ and $t$ to be the $y$ coordinate of the point just drawn. Set $\a$ to be the next letter. Then \textsf{GOTO Identify}.
    \end{enumerate}
\item \textsf{Construct and draw a 312-value chain in the following manner.}\\
    Take the subword of $w$ starting from $\a$ up to $c''$ which equalizes the numbers of $c'$ and $c''$. Hence, this subword is in $\{c'b_s,c'b''\}\{c'c'',cc'c''\}^*\{c'',cc''\}$. Remove the first $c'$ and the first $c''$. Additionally, if there exists $c$ in between these $c'$ and $c''$, remove it as well. Iteratively remove $c'$, $c''$ and $c$ until $b_s$ or $b''$ is the only letter left. Every time we remove $c'$ and $c''$ together, sum 21 with $\oplus_1$, and every time we remove $c'$, $c''$ and $c$ together, sum 312 with $\oplus_1$. Draw the 312-value chain that we obtain from this procedure so that the chain is horizontally located to the right of $P_c$, and vertically greater than $t$. Set $t$ to be the greatest $y$ coordinate of all points in the chain and $P_c$ to be the right-most point in the chain.\begin{enumerate}[\sffamily a.]
    \item If $b_s$ is a part of the subword, then draw a point located horizontally in between $P_a$ and $P_b$, and vertically in between the point with the least $y$ coordinate and the point with the second least $y$ coordinate in the chain. Set $P_b$ to be the point drawn for $b_s$. Set $\a$ to be the letter immediately after the subword. Then \textsf{GOTO Identify}.
    \item If $b''$ is a part of the subword, then draw a point located horizontally in between $P_\l$ and $P_b$, and vertically in between the point with the least $y$ coordinate and the point with the second least $y$ coordinate in the chain. Set $P_b=P_\l$. Set $\a$ to be the letter immediately after the subword. Then \textsf{GOTO Identify}.
    \end{enumerate}
\item \textsf{Draw a point due to which letter $\a$ is.}
    \begin{enumerate}[\sffamily a.]
    \item If $\a=a$ or $\a=b$, then draw a point at $(x,y)$ where $P_a^{(x)}<x<P_b^{(x)}$ and $t<y$. Set $P_a$ to be this new point and $t=P_a^{(y)}$.\begin{enumerate}[\sffamily i.]
        \item If $\a=b$ is the last letter in $w$, then \textsf{GOTO Flatten}.
        \item Otherwise, set $\a$ to be the next letter. Then \textsf{GOTO Identify}.
        \end{enumerate}
    \item If $\a=c$, then draw a point at $(x,y)$ where $P_c^{(x)}<x$ and $t<y$. Set $P_c$ to be this new point and $t=P_c^{(y)}$. Set $\a$ to be the next letter. Then \textsf{GOTO Identify}.
    \end{enumerate}
\item[\sffamily Flatten:] \textsf{Finalize the permutation by flattening.}\\
    Let $\pi$ be a permutation obtained by flattening the constructed graph. \textsf{OUTPUT} $\pi$.
\end{enumerate}
\vhhh

We define the function $\phi_N$ on $N$ that $\phi_N(\s)$ is a word obtained by applying \textsf{N-ENCODE} to $\s$, and the function $\psi_N$ on $K_1$ that $\psi_N(w)$ is a permutation obtained by applying \textsf{N-DECODE} to $w$. We then have the following lemma.\begin{lemma}\label{lem:6.1}
\textit{The function $\phi_N$ is a bijection between $N$ and $K_1$.}
\end{lemma}
\textit{Proof.} We first show that the image of $\phi_N$ is in $K_1$. Let $\s$ be in $N$ and $w=\phi_N(\s)$. It is clear that $w$ is in $\Sigma_1^*$, so we need to show that $w$ satisfies all conditions of $K_1$. We can easily verify that each case of Case 1 through Case 6 in Figure \ref{fig:5.21} corresponds to each prefix of $ba$, $bb'a$, $ba'a$, $ba'b_sa''$, $ba'b'a''$ and $ba'a'a''$ respectively. Similarly, each case of Case 1 through Case 6 in Figure \ref{fig:5.22} respectively corresponds to each suffix of $cb$, $cb''b$, $cc''b$, $c'b_sc''b$, $c'b''c''b$ and $c'c''c''b$.\\

For the condition on $a'$ and $a''$, notice that a 231-value chain of the form $21\oplus_1 21\oplus_1\cdots\oplus_1 21$ together with the scissor of the chain is encoded as $a'a'a''a'a''\cdots a'a''b_sa''$ (or $b'$ instead of $b_s$). Suppose there exists a summand 231 instead of 21, so the value playing the role of 2 is encoded as $a$. If the first summand is 231, then $a$ appears immediately after $a'$, corresponding to 1 of the same 231. Otherwise, due to the structure of a 231-value chain, $a$ appears immediately after $a''$, corresponding to 2 of 21 or 3 of 231 of the previous summand. Since 231-value chains are the only ways for $w$ to have $a'$ and $a''$, the third condition is satisfied. Proving the fourth condition is satisfied is almost identical. The fifth condition is also immediate, since \textsf{N-ENCODE} is defined so that $b_s$ is only used for scissors of 231-value chains and 312-value chains.\\

The condition on $b'$ and $b''$ is tied with the sixth condition of Proposition \ref{prop:5.4}, which states there are at most four values in between two values playing the roles of 1 and 2 of 12 in Equation \ref{eqn:5.1}. The positions of two values $x_1$ and $x_2$ ($x_1<x_2$) are in $[1,\s^{-1}(n))$ and the positions of the other two values $y_1$ and $y_2$ ($y_1<y_2$) are in $(\s^{-1}(1),n]$. Moreover, in the proof of Lemma \ref{lem:5.8}, we showed that $x_1$ and $y_2$ respectively must be the maximum value of a 231-value chain and the minimum value of a 312-chain. Altogether, these values are encoded into a subword in $\{b'\}\{a'',\lambda\}\{a,\lambda\}\{c,\lambda\}\{c',\lambda\}\{\b''\}$. Also, since there must be at least one value in between two values playing the roles of 1 and 2 of 12 in Equation \ref{eqn:5.1}, one of $a''$, $a$, $c$ and $c'$ must be in between $b'$ and $b''$. Thus, $w$ satisfies the sixth condition of $K_1$.\\

Finally, for the seventh condition, suppose $w$ contains $aa$. This implies two values for $1\oplus 1$ in Equation \ref{eqn:5.2} are consecutive, so it violates the first condition of Proposition \ref{prop:5.2}. We can apply the same arguments for $bb$ and $cc$ with the second and third conditions of Proposition \ref{prop:5.2} respectively, so $w$ does not contain $aa$, $bb$ and $cc$. Consequently, the image of $\phi_N$ is in $K_1$.\\

To show that the image of $\psi_N$ is in $N$, notice that together with how \textsf{N-DECODE} is defined and conditions of $K_1$, $\psi_N(w)$ is a permutation of extreme pattern 2413 and obeys all conditions of Proposition \ref{prop:5.2} and \ref{prop:5.4}. Therefore, by Proposition \ref{prop:5.6}, $\psi_N(w)$ is in $N$, so the image of $\psi_N$ is in $N$.\\

Due to the constructions of $\phi_N$ and $\psi_N$, they are inverse to each other, so this completes the proof.\hfill$\blacksquare$\\

Next, we define \textsf{S-ENCODE} and \textsf{S-DECODE}, the algorithm to encode simple permutations in $S$, and vice versa. To make things less complicated, we simply define these algorithms as \textsf{N-ENCODE} and \textsf{N-DECODE} applied to the inverse permutations. To be precise, To \textsf{S-ENCODE} $\tau$ in $S$, we \textsf{N-ENCODE} $\tau^{-1}$, and to \textsf{S-DECODE} a word $w$ in $K_1$, we \textsf{N-DECODE} to obtain a permutation $\s$ in $N$, then output $\s^{-1}$.\\

We define $\phi_S$ and $\psi_S$ that $\phi_S(\tau)$ is a word obtained by applying \textsf{S-ENCODE} to a permutation $\tau$ in $S$, and $\psi_S(w)$ is a permutation obtained by applying \textsf{S-DECODE} to a word $w$ in $K_1$. Thus, $\phi_N(\s)$ and $\phi_S(\tau)$ are equal if and only if $\s$ and $\tau$ are inverse to each other. It is straightforward to prove the following lemma by definition of $\phi_S$ and Lemma \ref{lem:6.1}.
\begin{lemma}\label{lem:6.2}
\textit{The function $\phi_S$ is a bijection between $S$ and $K_1$.}
\end{lemma}
\vhhh

Later with $\phi'$, we first decompose a permutation $\pi$ in $H'$ into the sum of $m$ simple permutations in $N$ and $S$, and apply \textsf{N-ENCODE} and \textsf{S-ENCODE} to each of these simple permutations alternatively. We abuse the notations of glue sums, and obtain $\overline{w}=w_1\nwsum_{x_1}^{y_1}w_2\sesum_{x_2}^{y_2}w_3\nwsum_{x_3}^{y_3}\cdots w_m$ to retain which glue sums were used in the expression of $\pi$. Thus, we obtain total $m$ words in $K_1$ connected by glue sums. In particular, $w_i$ for odd $i$ is obtained by using \textsf{N-ENCODE} to each $\s_i$ and $w_i$ for even $i$ is obtained by using \textsf{S-ENCODE} to each $\tau_i$. To make it simpler, we define $\phi_1$ on $H'$ such that\[
\phi_1(\pi)=\left\{\begin{array}{cl}\phi_N(\s_1)\nwsum_{x_1}^{y_1}\phi_S(\tau_2)\sesum_{x_2}^{y_2}\cdots\sesum_{x_{m-1}}^{y_{m-1}}\phi_N(\s_m) & \textrm{ if }m\textrm{ is odd}\\
\phi_N(\s_1)\nwsum_{x_1}^{y_1}\phi_S(\tau_2)\sesum_{x_2}^{y_2}\cdots\nwsum_{x_{m-1}}^{y_{m-1}}\phi_S(\tau_m) & \textrm{ if }m\textrm{ is even}\end{array}\right.
\]
Let $I_j=\{w_1\nwsum_{x_1}^{y_1}w_2\sesum_{x_2}^{y_2}\cdots w_j:w_j\in K_1\textrm{ for all }i\textrm{ with }1\leq i\leq j\}$ and $K_2$ be the image of $\phi_1$. We discuss that $K_2$ is a subset of $\bigcup_{j=1}^m I_j$, but not equal.\\

Recall that each glue sum requires summands to satisfy specific conditions. For example, type 1-0 and 1-1 require the left summand and the right summand to have Case 1 in Figure \ref{fig:5.22} and Case 1 in Figure \ref{fig:5.23} respectively. These restrictions coming from each glue sum force $w_i$ and $w_{i+1}$ to have specific suffix and prefix for all $i$ ($1\leq i\leq m-1$). In particular, for every $w$ in $K_2$, if the $i$-th glue sum is type 1-0 or 1-1, then $w_i$ ends with $cb$ and $w_{i+1}$ begins with $ba$. If it is type 2-0 or 2-1, then $w_i$ ends with $cb''b$ and $w_{i+1}$ begins with $ba$. If it is type 3-0, then $w_i$ ends with $cc''b$ and $w_{i+1}$ begins with $ba'a$, and finally, if it is type 4-0, then $w_i$ ends with $c'b_sc''b$ or $c'b''c''b$ and $w_{i+1}$ begins with $bb'a$.\\

Now, we define $\Sigma_2=\Sigma'\setminus\{d_a',d_a'',d_c',d_c'',d_\l\}$. The next algorithm, called \textsf{W-COMBINE} iteratively connects these words to construct one word $w$ in $\Sigma_2^*$. So let us define \textsf{W-COMBINE} as the following.\\
\\
\textbf{Algorithm} \textsf{W-COMBINE}\\
\textsf{INPUT}: A sequence of words connected by glue sums as $\overline{w}=w_1\nwsum_{x_1}^{y_1}w_2\sesum_{x_2}^{y_2}\cdots w_m$ in $K_2$.\\
\textsf{OUTPUT}: A word $w$ in $\Sigma_2^*$.\\
\begin{enumerate}[\hspace{7pt}\sffamily {Case} 1:]
\item[\sffamily Initialize:] \textsf{Let $w=w_1$ and write $\overline{w}$ as $w\nwsum_{x_1}^{y_1}w_2\sesum_{x_2}^{y_2}\cdots w_m$.}
\item[\sffamily Identify:] \textsf{If $\overline{w}=w$, then OUTPUT $w$. Otherwise, let $i+1$ be the least index of $w_k$ in $\overline{w}$.}
    \begin{enumerate}[\sffamily a.]
    \item If $(x_{i},y_{i})=(1,0)$ or $(x_{i},y_{i})=(1,1)$, then \textsf{GOTO Case 1.}
    \item If $(x_{i},y_{i})=(2,0)$ or $(x_{i},y_{i})=(2,1)$, then \textsf{GOTO Case 2.}
    \item If $(x_{i},y_{i})=(3,0)$, then \textsf{GOTO Case 3.}
    \item If $(x_{i},y_{i})=(4,0)$, then \textsf{GOTO Case 4.}
    \end{enumerate}
\item \textsf{A suffix of $w$ is $cb$ and a prefix of $w_{i+1}$ is $ba$.}
    \begin{enumerate}[\sffamily a.]
    \item If $(x_{i},y_{i})=(1,0)$, replace the suffix $cb$ of $w$ with $d$ and erase the prefix $ba$ of $w_{i+1}$. Then concatenate $w$ with $w_{i+1}$. Redefine $w$ to be the resultant word. Then \textsf{GOTO Identify}.
    \item If $(x_{i},y_{i})=(1,1)$, replace the suffix $b$ of $w$ with $d$ and erase the prefix $ba$ of $w_{i+1}$. Then concatenate $w$ with $w_{i+1}$. Redefine $w$ to be the resultant word. Then \textsf{GOTO Identify}.
    \end{enumerate}
\item \textsf{A suffix of $w$ is $cb''b$, $c\underline{b}b$ or $cy''b$, and a prefix of $w_{i+1}$ is $ba$.}
    \begin{enumerate}[\sffamily a.]
    \item Suppose a suffix of $w$ is $cb''b$.\begin{enumerate}[\sffamily i.]
        \item If $(x_{i},y_{i})=(2,0)$, replace the suffix $cb''b$ of $w$ with $\underline{d}d$ and erase the prefix $ba$ of $w_{i+1}$. Replace the last $b'$ in $w$ with $\overline{b}$. Then concatenate $w$ with $w_{i+1}$. Redefine $w$ to be the resultant word. Then \textsf{GOTO Identify}.
        \item If $(x_{i},y_{i})=(2,1)$, replace the suffix $b''b$ of $w$ with $\underline{d}d$ and erase the prefix $ba$ of $w_{i+1}$. Replace the last $b'$ in $w$ with $\overline{b}$. Then concatenate $w$ with $w_{i+1}$. Redefine $w$ to be the resultant word. Then \textsf{GOTO Identify}.
        \end{enumerate}
    \item Suppose a suffix of $w$ is $c\underline{b}b$.\begin{enumerate}[\sffamily i.]
        \item If $(x_{i},y_{i})=(2,0)$, replace the suffix $c\underline{b}b$ of $w$ with $\underline{d}d$ and erase the prefix $ba$ of $w_{i+1}$. Then concatenate $w$ with $w_{i+1}$. Redefine $w$ to be the resultant word. Then \textsf{GOTO Identify}.
        \item If $(x_{i},y_{i})=(2,1)$, replace the suffix $\underline{b}b$ of $w$ with $\underline{d}d$ and erase the prefix $ba$ of $w_{i+1}$. Then concatenate $w$ with $w_{i+1}$. Redefine $w$ to be the resultant word. Then \textsf{GOTO Identify}.
        \end{enumerate}
    \item Suppose a suffix of $w$ is $cy''b$.\begin{enumerate}[\sffamily i.]
        \item If $(x_{i},y_{i})=(2,0)$, replace the suffix $cy''b$ of $w$ with $\underline{d}d$ and erase the prefix $ba$ of $w_{i+1}$. Replace the last $y'$ in $w$ with $\overline{y}$. Then concatenate $w$ with $w_{i+1}$. Redefine $w$ to be the resultant word. Then \textsf{GOTO Identify}.
        \item If $(x_{i},y_{i})=(2,1)$, replace the suffix $y''b$ of $w$ with $\underline{d}d$ and erase the prefix $ba$ of $w_{i+1}$. Replace the last $y'$ in $w$ with $\overline{y}$. Then concatenate $w$ with $w_{i+1}$. Redefine $w$ to be the resultant word. Then \textsf{GOTO Identify}.
        \end{enumerate}
    \end{enumerate}
\item \textsf{A suffix of $w$ is $c'b_scc''b$, $c'b''cc''b$, $c'\underline{b}cc''b$, $c'c''cc''b$ or $c'y''cc''b$, and a prefix of $w_{i+1}$ is $ba'avb_sa''$ or $ba'avb'a''$ where $v$ is a word in $\{a'a'',a'a''a\}^*$.}\\
    For a suffix of $w$,\begin{enumerate}[\sffamily a.]
    \item if it is $c'b_scc''b$, then replace it with $zxd$.
    \item if it is $c'b''cc''b$, then replace it with $zx''d$ and the last $b'$ in $w$ with $x'$.
    \item if it is $c'\underline{b}cc''b$, then replace it with $z\underline{x}d$.
    \item if it is $c'c''cc''b$, then replace it with $zc''d$.
    \item if it is $c'y''cc''b$, then replace it with $zy''d$.
    \end{enumerate}
    For a prefix of $w_{i+1}$,\begin{enumerate}[\sffamily a.]
    \item if it is $ba'avb_sa''$ ($v\in\{a'a'',a'a''a\}^*$), then replace the whole prefix $ba'avb_sa''$ with $y$.
    \item if it is $ba'avb'a''$, then replace the whole prefix $ba'avb'a''$ with $y'$. Additionally, replace the first $b''$ in $w_{i+1}$ with $y''$.
    \end{enumerate}
    Concatenate $w$ with $w_{i+1}$. Redefine $w$ to be the resultant word. Then \textsf{GOTO Identify}.
\item \textsf{A suffix of $w$ is $c'b_sc''b$, $c'b''c''b$, $c'\underline{b}c''b$ or $c'y''c''b$, and a prefix of $w_{i+1}$ is $bb'a$.}
    \begin{enumerate}[\sffamily a.]
    \item If the suffix of $w$ is $c'b_sc''b$, then replace the suffix $c'b_sc''b$ of $w$ with $\overline{d}d$ and erase the prefix $bb'a$ of $w_{i+1}$. Replace the first $b''$ in $w_{i+1}$ with $\underline{b}$. Then concatenate $w$ with $w_{i+1}$. Redefine $w$ to be the resultant word. Then \textsf{GOTO Identify}.
    \item If the suffix of $w$ is $c'b''c''b$, then replace the suffix $c'b''c''b$ of $w$ with $\underline{\overline{d}}d$ and erase the prefix $bb'a$ of $w_{i+1}$. Replace the first $b''$ in $w_{i+1}$ with $\underline{b}$. Replace the last $b'$ in $w$ with $\overline{b}$. Then concatenate $w$ with $w_{i+1}$. Redefine $w$ to be the resultant word. Then \textsf{GOTO Identify}.
    \item If the suffix of $w$ is $c'\underline{b}c''b$, then replace the suffix $c'\underline{b}c''b$ of $w$ with $\underline{\overline{d}}d$ and erase the prefix $bb'a$ of $w_{i+1}$. Replace the first $b''$ in $w_{i+1}$ with $\underline{b}$. Then concatenate $w$ with $w_{i+1}$. Redefine $w$ to be the resultant word. Then \textsf{GOTO Identify}.
    \item If the suffix of $w$ is $c'y''c''b$, then replace the suffix $c'y''c''b$ of $w$ with $\underline{\overline{d}}d$ and erase the prefix $bb'a$ of $w_{i+1}$. Replace the last $y'$ in $w$ with $\overline{y}$ and the first $b''$ in $w_{i+1}$ with $\underline{b}$. Then concatenate $w$ with $w_{i+1}$. Redefine $w$ to be the resultant word. Then \textsf{GOTO Identify}.
    \end{enumerate}
\end{enumerate}
\vhhh

Let $K_3$ be the language over $\Sigma_2$ with the following conditions. We claim that $K_3$ is in bijection with $K_2$.\begin{enumerate}
\item Prefix condition.\\
      A word $w$ must begin with $ba$, $bb'a$, $ba'a$, $ba'b_sa''$, $ba'b'a''$, $ba'a'a''$, $b\overline{b}a$, $ba'\overline{b}a''$, $bx'a$ or $ba'x'a''$.
\item Suffix condition.\\
      A word $w$ must end with $cb$, $cb''b$, $cc''b$, $c'b_sc''b$, $c'b''c''b$, $c'c''c''b$, $c\underline{b}b$, $c'\underline{b}c''b$, $cy''b$ or $c'y''c''b$.
\item Conditions on $a'$ and $a''$.\\
      Every $a'$ and $a''$ in $w$ is a part of a subword in $\{a',a'a\}\{a'a'',a'a''a\}^*\{b_sa'',b'a'',\overline{b}a'',x'a''\}$. Note the number of $a'$ and the number of $a''$ are equal in each sequence, and thus, in $w$.
\item Conditions on $c'$ and $c''$.\\
      Every $c'$ and $c''$ in $w$ is a part of a subword in $\{c'b_s,c'b'',c'\underline{b},c'y''\}\{c'c'',cc'c''\}^*\{c'',cc'',zc'',\allowbreak czc''\}$. Note the number of $c'$ and the number of $c''$ are equal in each sequence, and thus, in $w$.
\item Conditions on $b_s$.\\
      The letter $b_s$ is only allowed in the third condition and the fourth condition, \textit{i.e.} as a part of a subword of $a'$ and $a''$ or a subword of $c'$ and $c''$.
\item Conditions on $b'$ and $b''$.\\
      Every $b'$ and $b''$ in $w$ is a part of a subword in $\{b'\}\{a'',\lambda\}\{a,\lambda\}\{c,\lambda\}\{c',\lambda\}\{b''\}$ with at least one of $a''$, $a$, $c$ or $c'$ being present. Note the number of $b'$ and the number of $b''$ are equal in each sequence, and thus, in $w$.
\item Conditions on letters with overlines and underlines.\\
      For every letter with an overline, there is a corresponding letter with an underline.\begin{itemize}
      \item Every $\overline{b}$ in $w$ is a part of a subword in $\{\overline{b}\}\{a'',\lambda\}\{a,\lambda\}\{c,\lambda\}\{\underline{d}d,\underline{\overline{d}}d\}$.
      \item Every $\overline{y}$ in $w$ is a part of a subword in $\{\overline{y}\}\{a,\lambda\}\{c,\lambda\}\{\underline{d}d,\underline{\overline{d}}d\}$.
      \item Every $\overline{d}$ and $\underline{\overline{d}}$ in $w$ is a part of a subword in one of the following.\begin{itemize}
          \item $\{\overline{d}d,\underline{\overline{d}}d\}\{c,\lambda\}\{c',\lambda\}\{\underline{b}\}$.
          \item $\{\overline{d}d,\underline{\overline{d}}d\}\{c,\lambda\}\{z\underline{x}\}$.
          \item $\{\overline{d}d,\underline{\overline{d}}d\}\{c,\lambda\}\{\underline{d}d,\underline{\overline{d}}d\}$.
          \end{itemize}
      \end{itemize}
\item Conditions on $x$, $y$, $z$ and other related letters.\\
      Every $x$, $x'$, $x''$, $\underline{x}$, $y$, $y'$, $y''$, $\overline{y}$ and $z$ is a part of a subword $v_1dv_2$ where $v_1$ is in
          \begin{itemize}
          \item $\{zx,z\underline{x},zc'',zy''\}$ or
          \item $\{x'\}\{a'',\lambda\}\{a,\lambda\}\{c,\lambda\}\{zx''\}$, and
          \end{itemize}
      $v_2$ is in
          \begin{itemize}
          \item $\{y,\overline{y}\}$, or
          \item $\{y'\}\{a,\lambda\}\{c,\lambda\}\{y'',c'y'',zy''\}$.
          \end{itemize}
\item Repetition restrictions.\\
      $w$ must not contain $aa$, $bb$, $cc$ or $da$.
\end{enumerate}
\vhhh

Before we move onto the next algorithm \textsf{W-DECOMPOSE}, let us take a closer look at which letter is replaced by each of letters in $\{\underline{b},\overline{b},d,\underline{d},\overline{d},\underline{\overline{d}},x,x',x'',\underline{x},y,y',y'',\overline{y},z\}$ when words are connected by \textsf{W-COMBINE}. For this observation, we examine the simplest case, which is applying \textsf{W-COMBINE} to $w_\s\nwsum_x^y w_\tau$ where $\s$ and $\tau$ be simple permutations of extreme pattern 2413 and 3142 respectively with $|\s|=m$ and $|\tau|=n$, and $w_\s=\phi_N(\s)$, $w_\tau=\phi_S(\tau)$.\\

When $\s$ and $\tau$ are operated with $\nwsum_1^0$ or $\nwsum_1^1$, $w_\s$ ends with $cb$ which corresponds to values $\s(m)$ and $m$, and $w_\tau$ starts with $ba$ which corresponds to $\tau(1)$ and 1. Recall that, for the case of $\nwsum_1^0$, it eliminates the last value of $\s$ and 1 of $\tau$, and combines $m$ and $\tau(1)$ as it shifts up accordingly. Thus, $c$ of $cb$ and $a$ of $ba$ are these eliminated values, and $d$ takes the role of combined $m$ and $\tau(1)$. If $\nwsum_1^1$ is used instead, then $c$ of $cb$ stays there, as $\nwsum_1^1$ does not erase the last value of $\s$. Therefore, the letter $d$ is the only replacement happened in both cases. Technically, $d$ represents $m$ of $\s$ and $\tau(1)$ of $\tau$, but to refer back what the previous letter was, we set the convention to look at $w_\s$, not $w_\tau$. Thus, $d$ was originally encoded as $b$ in $w_\s$.\\

If $\s$ and $\tau$ are operated with $\nwsum_2^0$ or $\nwsum_2^1$, $w_\s$ ends with $cb''b$, corresponding to values $\s(m)=m-2$, $m-1$ and $m$ respectively, and $w_\tau$ starts with $ba$ corresponding to $\tau(1)$ and 1. For $\nwsum_2^0$, there are two things \textsf{W-COMBINE} does in the process of concatenation of $w_\s$ and $w_\tau$. First, where $w_\s$ and $w_\tau$ are concatenated, we would have $cb''bba$, but \textsf{W-COMBINE} replace this with $\underline{d}d$. We may think that $c$ and $a$ are simply deleted, since the last value of $\s$ and the first value of $\tau$ are eliminated. The first $b$ is combined with the second $b$ which results in $d$, and $b''$ is replaced with $\underline{d}$. Second, \textsf{W-COMBINE} replaces the letter $b'$ that is matched up with $b''$ of $cb''b$ with $\overline{b}$. The only difference for $\nwsum_2^1$ is that it does not erase $c$ of $cb''b$ as $\s\nwsum_2^1\tau$ keeps the last value of $\s$. Hence, in $w_\s$, $d$ was originally encoded as $b$, $\underline{d}$ was encoded as $b''$, and $\overline{b}$ was encoded as $b'$.\\

Next, we look at the case of $\nwsum_3^0$. Recall that $\s$ and $\tau$ must have a 312-value chain involving $\s(m)$ and a 312-position chain involving 1 respectively that are similar to each other. Thus, the last three greatest values of $\s$ are placed as shown in Case 3 of \ref{fig:5.22} where $\s(m)=m-2$, $m-1$ and $m$ are encoded by \textsf{N-ENCODE} as $c$, $c''$ and $b$ respectively. For $m-3$, there are three cases. If the 312-value chain is simply $\a=312$, then $m-3$ would be the scissor of this chain. As the first case, if this scissor is 1 of 1 in Equation \ref{eqn:5.1}, then it is encoded as $b_s$. Alternatively, if it is 2 of 12 in Equation \ref{eqn:5.1}, then it is $b''$. The last case is $\a\neq 312$. Then $m-3$ is a part of the chain, and it is encoded as $c''$. In either case, $m-4$ must be 1 of the last 312 in the 312-value chain, so it is $c'$. Altogether, the five greatest values of $\s$ are encoded as $c'b_scc''b$, $c'b''cc''b$ or $c'c''cc''b$. On the other hand, $\tau$ having the similar 312-position chain to start with, encoding of values of $\tau$ by \textsf{S-ENCODE} up to the end of the chain is either $ba'avb_sa''$ or $ba'avb'a''$ where $v\in\{a'a'',a'a''a\}^*$.\\

Now, for each case of $w_\s$ having $c'b_scc''b$, $c'b''cc''b$ and $c'c''cc''b$, \textsf{W-COMBINE} replaces it with $zxd$, $zx''d$ and $zc''d$ respectively. In particular, we consider as follows.\begin{itemize}
\item The first letter $c'$ is replaced with $z$.
\item The second letter is replaced with $x$, $x''$ and $c''$ accordingly (thus, in the case of $c''$, it remains unchanged). In case of replacing $b''$ with $x''$, \textsf{W-COMBINE} additionally replaces the paired $b'$ with $x'$.
\item The third letter $c$ is erased as $\s\nwsum_3^0\tau$ eliminates $\s(m)$.
\item The fourth letter $c''$ is merged into the scissor $b_s$ of $ba'avb_sa''$ or $b'$ of $ba'avb'a''$ in $w_\tau$, which will be encoded latter.
\item The last letter $b$ is replaced with $d$.
\end{itemize}
Since the whole 312-position chain of $\tau$ is merged into the 312-value chain of $\s$ by computing $\s\nwsum_3^0\tau$, every letter in $ba'avb_sa''$ or $ba'avb'a''$ of $w_\tau$ are erased, except the letter $b_s$ or $b'$ which corresponds to the scissor value. As this value gets combined with the value $m-1$ of $\s$, we replace it with the knew letter $y$ or $y'$ respectively. In case of replacing $b'$ with $y'$, we also replaces the paired $b''$ with $y''$. Hence, to summarize this case, previously in $w_\s$, $z$, $x$, $x''$, $d$, $y$ and $y'$ were all encoded as $c'$, $b_s$, $b''$, $b$, $b_s$ and $b'$ respectively.\\

Finally, for $\nwsum_4^0$, for the four greatest values of $\s$, $w_\s$ can either have $c'b_sc''b$ or $c'b''c''b$, whereas $w_\tau$ must start with $bb'a$. Since $\nwsum_4^0$ eliminates $\s(m-1)=m-1$ and $\s(m)=m-3$ of $\s$, $c'$ and $c''$ of both suffixes of $w_\s$ are erased. Also, $\nwsum_4^0$ erases $\tau(1)$, so $a$ of the prefix $bb'a$ of $w_\tau$ is deleted. The combination of $m$ and $\tau(1)$, corresponding to $b$ of $w_\s$ and $b$ of $w_\tau$ respectively, is replaced with $d$, just like in every other glue sum. If $m-2$ of $\s$ is $b_s$, then it is replaced with $\overline{d}$ together with $b'$ of $w_\tau$. The value $\tau(2)$ corresponding to $b'$ of $bb'a$ is 1 of 12 in Equation \ref{eqn:5.1}, and its paired 2 of 12 encoded as $b''$ is also replaced with $\underline{b}$. If $m-2$ of $\s$ is $b''$, then it is replaced with $\underline{\overline{d}}$ instead. The letter $b''$ of $w_\tau$ that is paired with $b'$ in the prefix is again replaced with $\underline{b}$, and in addition, $b'$ that is paired with $b''$ in the suffix of $w_\s$ is also replaced with $\overline{b}$. Therefore, in $w_\s$, $\overline{b}$, $\overline{d}$, $\underline{\overline{d}}$ and $d$ were encoded as $b'$, $b_s$, $b''$ and $b$ respectively, and $\underline{b}$ in $w_\tau$ was encoded as $b''$.\\

The above observations are just for the case of $w_\s\nwsum_x^y w_\tau$. As \textsf{W-COMBINE} recursively concatenates words, two consecutive glue sums sometimes result in the necessities of $\underline{x}$ and $\overline{y}$. For each, if we examine \textsf{W-COMBINE} closely, we find that they were originally $b''$ and $b'$ respectively.\\

Next, let us define the inverse algorithm called \textsf{W-DECOMPOSE} on $K_3$, and we will prove $K_2$ and $K_3$ are in bijection due to these algorithms.\\
\\
\textbf{Algorithm} \textsf{W-DECOMPOSE}\\
\textsf{INPUT}: A word $w$ in $K_3$.\\
\textsf{OUTPUT}: A sequence of words connected by glue sums as $\overline{w}=w_1\nwsum_{x_1}^{y_1}w_2\sesum_{x_2}^{y_2}\cdots w_m$ in $\bigcup_{j=1}^m I_j$.\\
\begin{enumerate}[\hspace{7pt}\sffamily {Case} 1:]
\item[\sffamily Initialize:] \textsf{Let $m$ be the number of $d$ in $w$. Let $w'_m=w$ Also, write $\overline{w}=w'_m$.}
\item[\sffamily Identify:] \textsf{Let $i$ be the index such that}\[
\overline{w}=\left\{\begin{array}{cl}
w'_i\nwsum_{x_i}^{y_i}w_{i+1}\sesum_{x_{i+1}}^{y_{i+1}}\cdots w_m & \textrm{ if }i\textrm{ is odd}\\
w'_i\sesum_{x_i}^{y_i}w_{i+1}\nwsum_{x_{i+1}}^{y_{i+1}}\cdots w_m & \textrm{ if }i\textrm{ is even}\\
\end{array}\right.
\]
    \textsf{If $i=1$, then OUTPUT $\overline{w}$. Otherwise, locate the last $d$ in $w'_i$.}
    \begin{enumerate}[\sffamily a.]
    \item If it is not a part of a subword $\underline{d}d$, $zsd$ ($s\in\{x,x'',\underline{x},c'',y'\}$), $\overline{d}d$ or $\underline{\overline{d}}d$, then \textsf{GOTO Case 1.}
    \item If it is a part of a subword $\underline{d}d$, then \textsf{GOTO Case 2.}
    \item If it is a part of a subword $zsd$ where $s\in\{x,x'',\underline{x},c'',y''\}$, then \textsf{GOTO Case 3.}
    \item If it is a part of a subword $\overline{d}d$ or $\underline{\overline{d}}d$, then \textsf{GOTO Case 4.}
    \end{enumerate}
\item \textsf{Split $w'_i$ into $w'_{i-1}$ and $w_i$ by type 1-0 or 1-1 glue sum.}\\
    Decatenate the suffix of $w'_i$ starting from the letter immediately after the last $d$. Concatenate $ba$ with this suffix of $w'_i$ and call it $w_i$. Let $w'_{i-1}$ be the remaining subword of $w'_i$.
    \begin{enumerate}[\sffamily a.]
    \item If the letter immediately before the last $d$ in $w'_{i-1}$ is not $c$, then replace the last $d$ in $w'_{i-1}$ with $cb$.
        \begin{enumerate}[\sffamily i.]
        \item If $i$ is odd, then replace $w'_i$ in $\overline{w}$ with $w'_{i-1}\sesum_1^0 w_i$. \textsf{GOTO Identify}.
        \item If $i$ is even, then replace $w'_i$ in $\overline{w}$ with $w'_{i-1}\nwsum_1^0 w_i$. \textsf{GOTO Identify}.
        \end{enumerate}
    \item If the letter immediately before the last $d$ in $w'_{i-1}$ is $c$, then replace the last $d$ in $w'_{i-1}$ with $b$.
        \begin{enumerate}[\sffamily i.]
        \item If $i$ is odd, then replace $w'_i$ in $\overline{w}$ with $w'_{i-1}\sesum_1^1 w_i$. \textsf{GOTO Identify}.
        \item If $i$ is even, then replace $w'_i$ in $\overline{w}$ with $w'_{i-1}\nwsum_1^1 w_i$. \textsf{GOTO Identify}.
        \end{enumerate}
    \end{enumerate}
\item \textsf{Split $w'_i$ into $w'_{i-1}$ and $w_i$ by type 2-0 or 2-1 glue sum.}\\
    Decatenate the suffix of $w'_i$ starting from the letter immediately after the last $d$. Concatenate $ba$ with this suffix of $w'_i$ and call it $w_i$. Let $w'_{i-1}$ be the remaining subword of $w'_i$.
    \begin{enumerate}[\sffamily a.]
    \item Suppose the last letter with an overline in $w'_{i-1}$ is $\overline{b}$. Replace this $\overline{b}$ with $b'$.
        \begin{enumerate}[\sffamily i.]
        \item If the letter immediately before the last $\underline{d}$ in $w'_{i-1}$ is not $c$, then replace the suffix $\underline{d}d$ in $w'_{i-1}$ with $cb''b$.
            \begin{enumerate}[\sffamily A.]
            \item If $i$ is odd, then replace $w'_i$ in $\overline{w}$ with $w'_{i-1}\sesum_2^0 w_i$. \textsf{GOTO Identify}.
            \item If $i$ is even, then replace $w'_i$ in $\overline{w}$ with $w'_{i-1}\nwsum_2^0 w_i$. \textsf{GOTO Identify}.
            \end{enumerate}
        \item If the letter immediately before the last $\underline{d}$ in $w'_{i-1}$ is $c$, then replace the suffix $\underline{d}d$ in $w'_{i-1}$ with $b''b$.
            \begin{enumerate}[\sffamily A.]
            \item If $i$ is odd, then replace $w'_i$ in $\overline{w}$ with $w'_{i-1}\sesum_2^1 w_i$. \textsf{GOTO Identify}.
            \item If $i$ is even, then replace $w'_i$ in $\overline{w}$ with $w'_{i-1}\nwsum_2^1 w_i$. \textsf{GOTO Identify}.
            \end{enumerate}
        \end{enumerate}
    \item Suppose the last letter with an overline in $w'_{i-1}$ is $\overline{d}$ or $\underline{\overline{d}}$.
        \begin{enumerate}[\sffamily i.]
        \item If the letter immediately before the last $\underline{d}$ in $w'_{i-1}$ is not $c$, then replace the suffix $\underline{d}d$ in $w'_{i-1}$ with $c\underline{b}b$.
            \begin{enumerate}[\sffamily A.]
            \item If $i$ is odd, then replace $w'_i$ in $\overline{w}$ with $w'_{i-1}\sesum_2^0 w_i$. \textsf{GOTO Identify}.
            \item If $i$ is even, then replace $w'_i$ in $\overline{w}$ with $w'_{i-1}\nwsum_2^0 w_i$. \textsf{GOTO Identify}.
            \end{enumerate}
        \item If the letter immediately before the last $\underline{d}$ in $w'_{i-1}$ is $c$, then replace the suffix $\underline{d}d$ in $w'_{i-1}$ with $\underline{b}b$.
            \begin{enumerate}[\sffamily A.]
            \item If $i$ is odd, then replace $w'_i$ in $\overline{w}$ with $w'_{i-1}\sesum_2^1 w_i$. \textsf{GOTO Identify}.
            \item If $i$ is even, then replace $w'_i$ in $\overline{w}$ with $w'_{i-1}\nwsum_2^1 w_i$. \textsf{GOTO Identify}.
            \end{enumerate}
        \end{enumerate}
    \item Suppose the last letter with an overline in $w'_{i-1}$ is $\overline{y}$. Replace this $\overline{y}$ with $y'$.
        \begin{enumerate}[\sffamily i.]
        \item If the letter immediately before the last $\underline{d}$ in $w'_{i-1}$ is not $c$, then replace the suffix $\underline{d}d$ in $w'_{i-1}$ with $cy''b$.
            \begin{enumerate}[\sffamily A.]
            \item If $i$ is odd, then replace $w'_i$ in $\overline{w}$ with $w'_{i-1}\sesum_2^0 w_i$. \textsf{GOTO Identify}.
            \item If $i$ is even, then replace $w'_i$ in $\overline{w}$ with $w'_{i-1}\nwsum_2^0 w_i$. \textsf{GOTO Identify}.
            \end{enumerate}
        \item If the letter immediately before the last $\underline{d}$ in $w'_{i-1}$ is $c$, then replace the suffix $\underline{d}d$ in $w'_{i-1}$ with $y''b$.
            \begin{enumerate}[\sffamily A.]
            \item If $i$ is odd, then replace $w'_i$ in $\overline{w}$ with $w'_{i-1}\sesum_2^1 w_i$. \textsf{GOTO Identify}.
            \item If $i$ is even, then replace $w'_i$ in $\overline{w}$ with $w'_{i-1}\nwsum_2^1 w_i$. \textsf{GOTO Identify}.
            \end{enumerate}
        \end{enumerate}
    \end{enumerate}
\item \textsf{Split $w'_i$ into $w'_{i-1}$ and $w_i$ by type 3-0 glue sum.}\\
    Decatenate the suffix of $w'_i$ starting from the letter immediately after the last $d$ and call it $v$. Let $w'_{i-1}$ be the remaining subword of $w'_i$.\begin{itemize}
    \item If $v$ starts with $y$, replace this $y$ with $b_sa''$.
    \item If $v$ starts with $y'$, replace this $y'$ with $b'a''$. Additionally, replace $y''$ in $v$ with $b''$.
    \end{itemize}
    To construct $w_i$,
    \begin{enumerate}[\sffamily a.]
    \item if $w'_{i-1}$ has a suffix $zxd$, $zx''d$, $z\underline{x}d$ or $zy''d$, concatenate $ba'a$ with $v$, and call it $w_i$.
    \item if $w'_{i-1}$ has a suffix $zc''d$, take the whole sequence of $c'$ and $c''$ which $zc''$ of $zc''d$ is a part of, and call it $t$. Thus, $t$ is a subword of $w'_{i-1}$ in $\{c'b_s,c'b'',c'\underline{b},c'y''\}\{c'c'',cc'c''\}^*\allowbreak\{zc'',czc''\}$. Let $u=ba'a$. After the prefix $c'b_s$, $c'b''$, $c'\underline{b}$ or $c'y''$ in $t$, for every $c'c''$ and $zc''$, concatenate $u$ with $a'a''$, and for every $cc'c''$ and $czc''$, concatenate $u$ with $a'a''a$. Once this procedure is completed, concatenate $u$ with $v$, and call it $w_i$.
    \end{enumerate}
    To complete $w'_{i-1}$,
    \begin{enumerate}[\sffamily a.]
    \item if $w'_{i-1}$ has a suffix $zxd$, then replace it with $c'b_scc''b$.
    \item if $w'_{i-1}$ has a suffix $zx''d$, then replace it with $c'b''cc''b$ and the last $x'$ in $w'_{i-1}$ with $b'$.
    \item if $w'_{i-1}$ has a suffix $z\underline{x}d$, then replace it with $c'\underline{b}cc''b$.
    \item if $w'_{i-1}$ has a suffix $zc''d$, then replace it with $c'c''cc''b$.
    \item if $w'_{i-1}$ has a suffix $zy''d$, then replace it with $c'y''cc''b$.
    \end{enumerate}
    Then, finally,\begin{itemize}
    \item if $i$ is odd, then replace $w'_i$ in $\overline{w}$ with $w'_{i-1}\sesum_3^0 w_i$. \textsf{GOTO Identify}.
    \item if $i$ is even, then replace $w'_i$ in $\overline{w}$ with $w'_{i-1}\nwsum_3^0 w_i$. \textsf{GOTO Identify}.
    \end{itemize}
\item \textsf{Split $w'_i$ into $w'_{i-1}$ and $w_i$ by type 4-0 glue sum.}\\
    Decatenate the suffix of $w'_i$ starting from the letter immediately after the last $d$. Concatenate $bb'a$ with this suffix of $w'_i$ and call it $w_i$. Replace $\underline{b}$ with $b''$. Let $w'_{i-1}$ be the remaining subword of $w'_i$.
    For $w'_{i-1}$,
    \begin{enumerate}[\sffamily a.]
    \item if its suffix is $\overline{d}d$, then replace it with $c'b_sc''b$.
    \item if its suffix is $\underline{\overline{d}}d$, observe the last letter with an overline in $w'_{i-1}$.
        \begin{enumerate}[\sffamily i.]
        \item If the last letter with an overline in $w'_{i-1}$ is $\overline{b}$, then replace it with $b'$ and the suffix $\underline{\overline{d}}d$ of $w'_{i-1}$ with $c'b''c''b$.
        \item If the last letter with an overline in $w'_{i-1}$ is $\overline{d}$ or $\overline{\underline{d}}$, then replace the suffix $\underline{\overline{d}}d$ of $w'_{i-1}$ with $c'\underline{d}c''b$.
        \item If the last letter with an overline in $w'_{i-1}$ is $\overline{y}$, then replace it with $y'$ and the suffix $\underline{\overline{d}}d$ of $w'_{i-1}$ with $c'y''c''b$.
        \end{enumerate}
    Then,\begin{itemize}
        \item if $i$ is odd, then replace $w'_i$ in $\overline{w}$ with $w'_{i-1}\sesum_4^0 w_i$. \textsf{GOTO Identify}.
        \item if $i$ is even, then replace $w'_i$ in $\overline{w}$ with $w'_{i-1}\nwsum_4^0 w_i$. \textsf{GOTO Identify}.
        \end{itemize}
    \end{enumerate}
\end{enumerate}
\vhhh

We define $\phi_2$ on $K_2$ that $\phi_2(\overline{w})=w$ is obtained by applying \textsf{N-COMBINE} to an arbitrary $\overline{w}\in K_2$, and $\psi_2$ on $K_3$ that $\psi_2(w)=\overline{w}$ is obtained by applying \textsf{W-DECOMPOSE} to an arbitrary $w\in K_3$. By using these two functions, we prove the following.
\begin{lemma}\label{lem:6.3}
\textit{The function $\phi_2$ is a bijection between $K_2$ and $K_3$.}
\end{lemma}
\textit{Proof.} We first claim that the image of $\phi_2$ is in $K_3$. In each $w_i$ of $\overline{w}=w_1\nwsum_{x_1}^{y_1}w_2\sesum_{x_2}^{y_2}\cdots w_m$, there are no $aa$, $bb$ and $cc$ by definition of $K_1$. Thus, if $w$ contains any of $aa$, $bb$, $cc$ or $da$, then this is caused by the process of \textsf{W-COMBINE}.\\

For $w_i$ and $w_{i+1}$ in $\overline{w}$ to be combined, a certain suffix of $w_i$ is modified to a sequence of letters ending with $d$ first in each case of \textsf{W-COMBINE}. After that, a certain prefix of $w_{i+1}$ is either erased (\textsf{Case 1}, \textsf{Case 2} and \textsf{Case 4}) or replaced (\textsf{Case 3}), and then connected to $w_i$. Erased prefixes are either $ba$ or $bb'a$, so the letter after these prefixes cannot be $a$. Also, in the case of replacement, it is always with $y$ or $y'$. Thus, in all cases, $da$ is not contained in $w$.\\

In some cases, such as \textsf{Case 2,a.i}, $b'$, $b''$ or $y'$ is replaced with another letter. However, none of these is replaced with $a$, $b$, $c$ or $d$, so it is impossible for this arrangement to cause $w$ containing $aa$, $bb$, $cc$ or $da$. Hence, $\phi_2(\overline{w})$ satisfies the ninth condition of $K_3$.\\

The other conditions of $K_3$ can be verified by observing each case of \textsf{W-COMBINE}. Thus, the image of $\phi_2$ is in $K_3$.\\

Next, we show that the image of $\psi_2$ is in $K_2$. Frankly, the algorithm \textsf{W-DECOMPOSE} reads a word $w$ from right to left and splits $w$ every time $d$ occurs. Based on the conditions of $K_3$, all possible cases of subwords that $d$ can be a part of are $\underline{d}d$, $zsd$ where $s\in\{x,x'',\underline{x},c'',y''\}$, $\overline{d}d$, $\underline{\overline{d}}d$ and $sd$ where $s\in\{a,a'',b,b'',c,c''\}$. In each case, \textsf{W-DECOMPOSE} decatenates $w_i$ as all letters in $\Sigma_2\setminus\Sigma_1$ are replaced with letters in $\Sigma_1$. We can observe that the word decatenated by \textsf{W-DECOMPOSE} is in $K_1$ for every case, implying the whole $\overline{w}$ is in $K_2$. Thus, the image of $\psi_2$ is in $K_2$.\\

Due to the constructions of $\phi_2$ and $\psi_2$, they are inverse to each other, so $\phi_2$ is a bijection from $K_2$ to $K_1$.\hfill$\blacksquare$\\

We are only one step away to the regular language $L'$ which is bijective to $H'$. Our final algorithm is to modify a prefix and a suffix of $w$ in $K_3$, which is the following.\\
\\
\textbf{Algorithm} \textsf{AFFIX-CONVERT}\\
\textsf{INPUT}: A word $w$ in $K_3$.\\
\textsf{OUTPUT}: A word in $\Sigma'^*$.\\
\\
\textsf{GOTO Prefix}.
\begin{enumerate}[\hspace{30pt}]
\item[\sffamily Prefix:] \textsf{Convert the prefix of $w$.}
    \begin{itemize}
    \item If it is $ba$, then replace it with $dd$. \textsf{GOTO Suffix.}
    \item If it is $bb'a$, then replace it with $db'd$. \textsf{GOTO Suffix.}
    \item If it is $ba'a$, then replace it with $da'd$. \textsf{GOTO Suffix.}
    \item If it is $ba'b_sa''$, then replace it with $dd_a'b_sd_a''$. \textsf{GOTO Suffix.}
    \item If it is $ba'b'a''$, then replace it with $dd_a'b'd_a''$. \textsf{GOTO Suffix.}
    \item If it is $ba'a'a''$, then replace it with $dd_a'a'd_a''$. \textsf{GOTO Suffix.}
    \item If it is $b\overline{b}a$, then replace it with $d\overline{b}d$. \textsf{GOTO Suffix.}
    \item If it is $ba'\overline{b}a''$, then replace it with $dd_a'\overline{b}d_a''$. \textsf{GOTO Suffix.}
    \item If it is $bx'a$, then replace it with $dx'd$. \textsf{GOTO Suffix.}
    \item If it is $ba'x'a''$, then replace it with $dd_a'x'd_a''$. \textsf{GOTO Suffix.}
    \end{itemize}
\item[\sffamily Suffix:] \textsf{Convert the suffix of $w$.}
    \begin{itemize}
    \item If it is $cb$, then replace it with $dd_\l$. \textsf{OUTPUT $w$}.
    \item If it is $cb''b$, then replace it with $db''d_\l$. \textsf{OUTPUT $w$}.
    \item If it is $cc''b$, then replace it with $dc''d_\l$. \textsf{OUTPUT $w$}.
    \item If it is $c'b_sc''b$, then replace it with $d_c'b_sd_c''d_\l$. \textsf{OUTPUT $w$}.
    \item If it is $c'b''c''b$, then replace it with $d_c'b''d_c''d_\l$. \textsf{OUTPUT $w$}.
    \item If it is $c'c''c''b$, then replace it with $d_c'c''d_c''d_\l$. \textsf{OUTPUT $w$}.
    \item If it is $c\underline{b}b$, then replace it with $d\underline{b}d_\l$. \textsf{OUTPUT $w$}.
    \item If it is $c'\underline{b}c''b$, then replace it with $d_c'\underline{b}d_c''d_\l$. \textsf{OUTPUT $w$}.
    \item If it is $cy''b$, then replace it with $dy''d_\l$. \textsf{OUTPUT $w$}.
    \item If it is $c'y''c''b$, then replace it with $d_c'y''d_c''d_\l$. \textsf{OUTPUT $w$}.
    \end{itemize}
\end{enumerate}
\vhhh

Define $\phi_3$ be a function on $K_3$ induced by \textsf{AFFIX-CONVERT}. By how \textsf{AFFIX-CONVERT} is defined, it is not difficult to verify that the image of $\phi_3$ is equal to the language $L'$ over $\Sigma'$ with following conditions.\begin{enumerate}
\item Prefix condition.\\
      A word $w$ must begin with $dd$, $db'd$, $da'd$, $dd_a'b_sd_a''$, $dd_a'b'd_a''$, $dd_a'a'd_a''$, $d\overline{b}d$, $dd_a'\overline{b}d_a''$, $dx'd$ or $dd_a'x'd_a''$. Therefore,\begin{enumerate}[$(a)$]
      \item $da'd$ is followed by a subword in $\{a'a'',a'a''a\}^*\{b_sa'',b'a'',\overline{b}a'',x'a''\}$.
      \item $dd_a'a'd_a''$ is followed by a subword in $\{a,\lambda\}\{a'a'',a'a''a\}^*\{b_sa'',b'a'',\overline{b}a'',x'a''\}$.
      \item $db'd$ is followed by a subword in $\{c,\lambda\}\{c',\lambda\}\{b''\}$ or the suffix $db''d_\l$, $d_c'b''d_c''d_\l$.
      \item $dd_a'b'd_a''$ is followed by a subword in $\{a,\lambda\}\{c,\lambda\}\{c',\lambda\}\{b''\}$ or the suffix $db''d_\l$, $d_c'b''d_c''d_\l$.
      \item $d\overline{b}d$ is followed by a subword in $\{c,\lambda\}\{\underline{d}d,\underline{\overline{d}}d\}$.
      \item $dd_a'\overline{b}d_a''$ is followed by a subword in $\{a,\lambda\}\{c,\lambda\}\{\underline{d}d,\underline{\overline{d}}d\}$.
      \item $dx'd$ is followed by a subword $v_1dv_2$ where $v_1$ is in $\{c,\lambda\}\{zx''\}$ and $v_2$ is in $\{y,\overline{y}\}$ or $\{y'\}\{a,\lambda\}\{c,\lambda\}\{y'',c'y'',zy''\}$.
      \item $dd_a'x'd_a''$ is followed by a subword $v_1dv_2$ where $v_1$ is in $\{a,\lambda\}\{c,\lambda\}\{zx''\}$ and $v_2$ is in $\{y,\overline{y}\}$ or
          $\{y'\}\{a,\lambda\}\{c,\lambda\}\{y'',c'y'',zy''\}$.
      \end{enumerate}
      Letters $d_a'$ and $d_a''$ are only permitted in the above listed prefixes.
\item Suffix condition.\\
      A word $w$ must end with $dd_\l$, $db''d_\l$, $dc''d_\l$, $d_c'b_sd_c''d_\l$, $d_c'b''d_c''d_\l$, $d_c'c''d_c''d_\l$, $d\underline{b}d_\l$, $d_c'\underline{b}d_c''d_\l$, $dy''d_\l$ or $d_c'y''d_c''d_\l$. Therefore,\begin{enumerate}[$(a)$]
      \item $dc''d_\l$ is preceded by a subword in $\{c'b_s,c'b'',c'\underline{b},c'y''\}\{c'c'',cc'c''\}^*$.
      \item $d_c'c''d_c''d_\l$ is preceded by a subword in $\{c'b_s,c'b'',c'\underline{b},c'y''\}\{c'c'',cc'c''\}^*\{c,\lambda\}$.
      \item $db''d_\l$ is preceded by a subword in $\{b'\}\{a'',\lambda\}\{a,\lambda\}$ or the prefix $db'd$, $dd_a'b'd_a''$.
      \item $d_c'b''d_c''d_\l$ is preceded by a subword in $\{b'\}\{a'',\lambda\}\{a,\lambda\}\{c,\lambda\}$ or the prefix $db'd$, $dd_a'b'd_a''$.
      \item $d\underline{b}d_\l$ is preceded by $\overline{d}d$ or $\underline{\overline{d}}d$.
      \item $d_c'\underline{b}d_c''d_\l$ is preceded by a subword in $\{\overline{d}d,\underline{\overline{d}}d\}\{c,\lambda\}$.
      \item $dy''d_\l$ is preceded by a subword $v_1dv_2$ where $v_1$ is in $\{zx,z\underline{x},zc'',zy''\}$ or $\{x'\}\{a'',\lambda\}\allowbreak\{a,\lambda\}\{c,\lambda\}\{zx''\}$ and $v_2$ is $y'$ or $y'a$.
      \item $d_c'y''d_c''d_\l$ is preceded by a subword $v_1dv_2$ where $v_1$ is in $\{zx,z\underline{x},zc'',zy''\}$ or $\{x'\}\{a'',\lambda\}\allowbreak\{a,\lambda\}\{c,\lambda\}\{zx''\}$ and $v_2$ is in $\{y'\}\{a,\lambda\}\{c,\lambda\}$.
      \end{enumerate}
      Letters $d_c'$, $d_c''$ and $d_\l$ are only permitted in the above listed suffixes. In particular, every $w$ must end with $d_\l$.
\item Conditions on $a'$ and $a''$.\\
      Every $a'$ and $a''$ in $w$ is a part of a prefix described in $(a)$ or $(b)$ in the prefix condition, or a part of a subword in $\{a',a'a\}\{a'a'',a'a''a\}^*\{b_sa'',b'a'',\overline{b}a'',x'a''\}$. Note the numbers of $a'$ and $a''$ are the same in each sequence, and thus, in the entire $w$.
\item Conditions on $c'$ and $c''$.\\
      Every $c'$ and $c''$ in $w$ is a part of a suffix described in $(a)$ or $(b)$ in the suffix condition, or a part of a subword in $\{c'b_s,c'b'',c'\underline{b},c'y''\}\{c'c'',cc'c''\}^*\{c'',cc'',zc'',czc''\}$. Note the numbers of $c'$ and $c''$ are the same in each sequence, and thus, in the entire $w$.
\item Conditions on $b_s$.\\
      The letter $b_s$ is only allowed in the prefix $dd_a'b_sd_a''$, in the suffix $d_c'b_sd_c''d_\l$, as a part of a subword of $a'$ and $a''$ listed in the third condition, or as a part of a subword of $c'$ and $c''$, listed in the fourth condition.
\item Conditions on $b'$ and $b''$.\\
      Every $b'$ and $b''$ in $w$ is a part of a prefix described in $(c)$ or $(d)$ in the prefix condition, a part of a suffix described in $(c)$ or $(d)$ in the suffix condition, or a part of a subword in $\{b'\}\{a'',\lambda\}\{a,\lambda\}\{c,\lambda\}\{c',\lambda\}\{b''\}$ with at least one of $a''$, $a$, $c$ or $c'$ being present. Note the numbers of $b'$ and $b''$ are the same in each sequence, and thus, in the entire $w$.
\item Conditions on letters with overlines and underlines.\\
      For every letter with an overline, there is a corresponding letter with an underline.\begin{itemize}
      \item Every $\overline{b}$ in $w$ is a part of a prefix described in $(e)$ or $(f)$ in the prefix condition, or a part of a subword in $\{\overline{b}\}\{a'',\lambda\}\{a,\lambda\}\{c,\lambda\}\{\underline{d}d,\underline{\overline{d}}d\}$.
      \item Every $\overline{y}$ in $w$ is a part of a subword in $\{\overline{y}\}\{a,\lambda\}\{c,\lambda\}\{\underline{d}d,\underline{\overline{d}}d\}$.
      \item Every $\overline{d}$ and $\underline{\overline{d}}$ in $w$ is a part of a subword in one of the following.\begin{itemize}
          \item $\{\overline{d}d,\underline{\overline{d}}d\}\{d\underline{b}d_\l\}$.
          \item $\{\overline{d}d,\underline{\overline{d}}d\}\{c,\lambda\}\{d_c'\underline{b}d_c''d_\l\}$.
          \item $\{\overline{d}d,\underline{\overline{d}}d\}\{c,\lambda\}\{c',\lambda\}\{\underline{b}\}$.
          \item $\{\overline{d}d,\underline{\overline{d}}d\}\{c,\lambda\}\{z\underline{x}\}$.
          \item $\{\overline{d}d,\underline{\overline{d}}d\}\{c,\lambda\}\{\underline{d}d,\underline{\overline{d}}d\}$.
          \end{itemize}
      \end{itemize}
\item Conditions on $x$, $y$, $z$ and other related letters.\\
      In addition to $(g)$, $(h)$ in the prefix condition and $(g)$, $(h)$ in the suffix condition, every $x$, $x'$, $x''$, $\underline{x}$, $y$, $y'$, $y''$, $\overline{y}$ and $z$ is a part of a subword $v_1dv_2$ where $v_1$ is in
          \begin{itemize}
          \item $\{zx,z\underline{x},zc'',zy''\}$ or
          \item $\{x'\}\{a'',\lambda\}\{a,\lambda\}\{c,\lambda\}\{zx''\}$, and
          \end{itemize}
      $v_2$ is in
          \begin{itemize}
          \item $\{y,\overline{y}\}$, or
          \item $\{y'\}\{a,\lambda\}\{c,\lambda\}\{y'',c'y'',zy''\}$.
          \end{itemize}
\item Other restrictions.\\
      $w$ must not contain $aa$, $bb$, $cc$, $da$ and $cdd_\l$.
\end{enumerate}
\vhhh

Finally, with every bijection we defined, let us define $\phi'$ on $H'$ as $\phi'=\phi_3\circ\phi_2\circ\phi_1$. Then, we have the following proposition. The proof is immediate with the lemmas we established.
\begin{proposition}\label{prop:6.1}
\textit{The encoding function $\phi'$ is a bijection between $H'$ and $L'$.}
\end{proposition}

To visualize the encoding of a permutation in $H'$ to $L'$, we refer to the Python code provided in Appendix \ref{app:B}.
\subsection{Defining the automaton $M'$}
Recall that we defined the language $\overline{L}$ associated with $L$ in Chapter \ref{chap:4}. While every word in $L$ must begin with $dd$, we defined $\overline{L}$ to be the set of words in $L$ without the prefix $dd$, and showed $\overline{L}$ is regular.\\

For $L'$, words have 10 distinct prefixes as described. We prepare associated 10 languages $\overline{L}_i$ ($1\leq i\leq 10$), one for each prefix, and construct an automaton for each of these languages to show that they are all regular languages. For each $i$ with $1\leq i\leq 10$, we define $\overline{L}_i$ as the following.\[
\begin{array}{lcl}
\overline{L}_1=\{w\in\Sigma'^*:ddw\in L'\}, & \hspace{0.4in}\textrm{ } &\overline{L}_2=\{w\in\Sigma'^*:dd_a'b_sd_a''w\in L'\},\\
\overline{L}_3=\{w\in\Sigma'^*:da'dw\in L'\}, & \hspace{0.4in} &\overline{L}_4=\{w\in\Sigma'^*:dd_a'a'd_a''w\in L'\},\\
\overline{L}_5=\{w\in\Sigma'^*:db'dw\in L'\}, & \hspace{0.4in} &\overline{L}_6=\{w\in\Sigma'^*:dd_a'b'd_a''w\in L'\},\\
\overline{L}_7=\{w\in\Sigma'^*:dx'dw\in L'\}, & \hspace{0.4in} &\overline{L}_8=\{w\in\Sigma'^*:dd_a'x'd_a''w\in L'\},\\
\overline{L}_9=\{w\in\Sigma'^*:d\overline{b}dw\in L'\}, & \hspace{0.4in} &\overline{L}_{10}=\{w\in\Sigma'^*:dd_a'\overline{b}d_a''w\in L'\}.
\end{array}\]
Each $\overline{L}_i$ is a language over $\overline{\Sigma}=\Sigma'\setminus\{d_a',d_a''\}$ that shares all conditions of $L'$, except for the prefix condition. Since an arbitrary $w$ in each $\overline{L}_i$ is a word in $L'$ without the certain prefix, the prefix condition for each $\overline{L}_i$ shall be stated with what the new prefix has to be. Thus, the prefix condition for $\overline{L}_i$ ($1\leq i\leq 10$) is as the following.
\begin{enumerate}[$L_1$:]
\item The first letter must be in $\{a',b,b',\overline{b},c,c',d,d_c',\overline{d},x',z\}$.
\item The first letter must be in $\{a,a',b,b',\overline{b},c,c',d,d_c',\overline{d},x',z\}$.
\item A word $w$ must begin with a subword in $\{a'a'',a'a''a\}^*\{b_sa'',b'a'',\overline{b}a'',x'a''\}$.
\item A word $w$ must begin with a subword in $\{a,\lambda\}\{a'a'',a'a''a\}^*\{b_sa'',b'a'',\overline{b}a'',x'a''\}$.
\item A word $w$ must begin with a subword in $\{c,\lambda\}\{c',\lambda\}\{b''\}$ or $w=db''d_\l$, $d_c'b''d_c''d_\l$.
\item A word $w$ must begin with a subword in $\{a,\lambda\}\{c,\lambda\}\{c',\lambda\}\{b''\}$ or $w=db''d_\l$, $d_c'b''d_c''d_\l$.
\item A word $w$ must begin with a subword $v_1dv_2$ where $v_1$ is in $\{c,\lambda\}\{zx''\}$ and $v_2$ is in $\{y,\overline{y}\}$ or $\{y'\}\{a,\lambda\}\{c,\lambda\}\{y'',c'y'',zy''\}$.
\item A word $w$ must begin with a subword $v_1dv_2$ where $v_1$ is in $\{a,\lambda\}\{c,\lambda\}\{zx''\}$ and $v_2$ is in $\{y,\overline{y}\}$ or $\{y'\}\{a,\lambda\}\{c,\lambda\}\{y'',c'y'',zy''\}$.
\item A word $w$ must begin with a subword in $\{c,\lambda\}\{\underline{d}d,\underline{\overline{d}}d\}$.
\item[$L_{10}$] A word $w$ must begin with a subword in $\{a,\lambda\}\{c,\lambda\}\{\underline{d}d,\underline{\overline{d}}d\}$.
\end{enumerate}
\vhhh

Next, we define 10 deterministic finite-state automatons $M'_i$ ($1\leq i\leq 10$), and show that $\mathcal{L}(M'_i)=\overline{L}_i$ for each $i$. Before we do so, we note that the only difference among each $M'_i$ is the initial state. Each automaton runs over $\Sigma'$ and shares the same set of states, transition function and accept state. The set of states $Q'$ for each $M'_i$ contains 83 distinct states. For this reason, it is unarguably not reasonable to present their transition function as a diagram. The description of the transition function is given in Table \ref{tab:A.1} in Appendix \ref{app:A}. We ask readers to refer to the same table to see each of 83 states in as well.\\

Let us now introduce $M'_i$ for each $i$ with $1\leq i\leq 10$. For every $i$ ($1\leq i\leq 10)$, we define an automaton $M'_i=(Q',\Sigma',\d',q_i,\{D_\l\})$ where $Q$ contains all states presented in Table \ref{tab:A.1}, $\d'$ is the transition function described in Table \ref{tab:A.1} in Appendix \ref{app:A}, and $q_i$ is the initial state such that\[
\begin{array}{lclclclcl}
q_1=A & \quad & q_2=A'' & \quad & q_3=A[A] & \quad & q_4=A''[A] & \quad & q_5=A[B]\\
q_6=A''[B] & \quad & q_7=A[X] & \quad & q_8=A''[X] & \quad & q_9=\overline{B}A & \quad & q_{10}=\overline{B}A''\\
\end{array}\]

Before we show $\mathcal{L}(M'_i)=\overline{L}_i$ for each $i$, we give a general description of each state in $Q$. Frankly, every state is denoted based on which letter in $\Sigma'$ is previously used to obtain the state.\\

A state in
\[
\{A,A'',B,B'',\underline{B},\overline{B},C,C'',D,D_c'',\underline{D},\overline{D},\underline{\overline{D}},D_\l,X,X'',\underline{X},Y,Y'',\overline{Y},Z\}
\]
is denoted by a single capital letter, which indicates the transition that was used to arrive at the state. For instance, $\underline{D}=\d(q,\underline{d})$ for some state $q\in Q$.\\

\begin{figure}[b!]
\[\begin{tikzpicture}[->,>=stealth',shorten >=1pt,auto,node distance=2.8cm,semithick,initial/.style={},scale=0.8]
  \tikzstyle{every state}=[circle,fill=black!25,minimum size=32pt,inner sep=0pt]
  \node[state] at (0,1)  (q1) {\tiny $A[B]$};
  \node[state] at (0,-1) (q2) {\tiny $B'[B]$};
  \node[state] at (4,3)  (q3) {\tiny $C[B]$};
  \node[state] at (4,1)  (q4) {\tiny $C'[B,C]$};
  \node[state] at (4,-1) (q5) {\tiny $D[B]$};
  \node[state] at (4,-3) (q6) {\tiny $D_c'[B,D]$};
  \node[state] at (8,3)  (q7) {\tiny $B''$};
  \node[state] at (8,1)  (q8) {\tiny $B''[C]$};
  \node[state] at (8,-1) (q9) {\tiny $DB''$};
  \node[state] at (8,-3) (q10){\tiny $B''[D]$};
  \path (q1) edge (q3)
             edge (q4)
             edge (q5)
             edge (q6)
             edge (q7)
        (q2) edge (q1)
             edge (q3)
             edge (q4)
             edge (q5)
             edge (q6)
        (q3) edge (q4)
             edge [bend left=40] (q6)
             edge (q7)
        (q4) edge (q8)
        (q5) edge (q9)
        (q6) edge (q10);
\end{tikzpicture}\]
\caption{Partial state diagram of $M'$.}
\label{fig:6.1}
\end{figure}
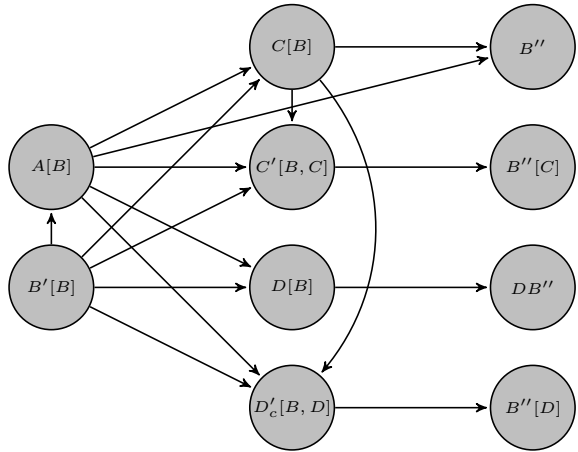

Many states come with square brackets with one or two capital letters inside, such as $A'[A]$ and $A'[A,A]$. The letter outside of the brackets indicates that the associated lower case letter in $\Sigma'$ that was just used for transition to arrive at the state. Each letter in the brackets, on the other hand, states that there must eventually be a transition using the associated lower case letter with a double prime. The appearance of a letter in square brackets is initially caused by a transition with a prime sign in $\Sigma'$ such as $a'$ and $b'$. For example, from a state with no brackets, say $A$, we have $\d'(A,b')=B'[B]$ as described in Table \ref{tab:A.1}. Until its paired $b''$ is used for transition, the letter $B$ is shown in brackets of the states following $B'[B]$. Indeed, from $B'[B]$, we eventually arrive at $B''$, $B''[C]$, $DB''$ or $B''[D]$ by $b''$ from a previous state as shown in Figure \ref{fig:6.1}. Note that labels of transitions are omitted in Figure \ref{fig:6.1}.
\\

Notice that, as we see in $\d'(B'[B],c')=C'[B,C]$, it is possible to have another transition with a prime sign until we reach a state with $B''$ involved. In this case, we also must have the transition $c''$ eventually, so we include $C$ in the square brackets as well.\\

Next, we give a description for states denoted by two consecutive upper case letters. Namely, the list of these states are given as follows.\[
\{\overline{B}A,\overline{B}A'',\overline{B}C,CD,DB'',D\underline{B},DC'',\underline{D}D,\overline{D}D,DY'',XD,ZC'',ZY''\}
\]
The transition to a state $q$ in the above set is the lower case of the second letter in the expression of $q$. For instance, to arrive at $ZC''$, we use $c''$ from some appropriate state. On the other hand, the first letter does not indicate anything relevant about a transition. One may think that the first letter informs which letter was used to obtain the state immediately before the current state, but this is only true for some cases. For example, the only way to arrive at $\overline{B}A''$ is from $\overline{B}[A]$ with the transition $a''$. From here, it is possible to go to $\overline{B}A$ with the transition $a$. As we can see, the previous transition $a''$ is not recorded in the notation $\overline{B}A$.\\

Finally, we explain the states with arrows involved, that is, the states in the set\[
\{C\ra\{\underline{B},\underline{D},\underline{\overline{D}},\underline{X}\},C'[C] \ra \underline{B}[C],D \ra\underline{B},D_c'[D]\ra\underline{B}[D],Z \ra \underline{X}\}.
\]
Each state in this set is denoted based on which transition with an underline is used in the future. In order to arrive at a state in the above set, we must first obtain $\overline{D}D$. Once the transitions $\overline{d}$ or $\underline{\overline{d}}$ occur, we obtain $\overline{D}$ or $\underline{\overline{D}}$ respectively, and as a result, we arrive at $\overline{D}D$. From here, we are required to have either $\underline{b}$, $\underline{d}$, $\underline{\overline{d}}$ or $\underline{x}$ in the near future. Hence, we arrive at $\underline{B}$, $\underline{B}[C]$, $\underline{B}[D]$, $\underline{D}$, $\underline{\overline{D}}$ or $\underline{X}$. These states with arrows are necessary to distinguish them from other states, which do not require future transitions to underlined sates.\\

In each expression of a state with an arrow, the letter at the tail of the arrow indicates the associated transition to reach there. To reach $C\ra\{\underline{B},\underline{D},\underline{\overline{D}},\underline{X}\}$, for example, we must have $c$ from $\overline{D}D$. The state, or the set of states at the head of the arrow shows which transition with an underline will be possibly used. This transition with an underline, however, does not have to occur immediately from the state with an arrow, as we can observe with $\d'(C\ra\{\underline{B},\underline{D},\underline{\overline{D}},\underline{X}\},c')=C'[C]\ra\underline{B}[C]$. Square brackets are used for the same purpose as before.\\

We are now ready to prove the following proposition.
\begin{proposition}\label{prop:6.2}
\textit{For every $i$ with $1\leq i\leq 10$, $\mathcal{L}(M'_i)=\overline{L}_i$.}
\end{proposition}
\textit{Proof.} We first show that $\mathcal{L}(M'_i)\subseteq\overline{L}_i$ for each $i$ with $1\leq i\leq 10$. We define the set of letters $\Gamma$ as\[
\Gamma=\{a,a',b,b',\overline{b},c,c',d,d_c',\overline{d},x',z\}.
\]
We also let $P'$ be a subset of $Q'$ defined as $P'=\{A,A'',B,B'',\underline{B},C,C'',CD,D,\underline{D}D,Y,Y''\}$.\\

Let $i$ be arbitrary with $1\leq i\leq 10$. In $M_i'$, we verify that every sequence of transitions from a state to a non-jail state constructs a subword obeying the conditions of $\overline{L}_i$. As we go on, we also make sure that the initial state for each $M_i'$ is the appropriate one.\\

Let $w$ be in $\mathcal{L}(M_i')$. Due to certain conditions of $\overline{L}_i$, for the letters in \[\overline{\Sigma}\setminus\Gamma=\{a'',b'',b_s,\underline{b},c'',d_c'',\underline{d},\underline{\overline{d}},d_\l,x,\underline{x},y,y',y'',\overline{y}\}\]
to appear in $w$, it requires some specific letters to be previously appearing in $w$. For instance, $b''$ is not allowed to appear in $w$ unless there is an unmatched $b'$ previously showing up in $w$. Hence, the letters in $\Gamma$ are the ones which do not require any particular preceding letters in $w$.\\

On the other hand, the letters in the set\[
\Delta=\{a,a',b,b',b_s,\overline{b},c,c',d,d_c',d_c'',\underline{d},\overline{d},\underline{\overline{d}},x,x',x'',\underline{x},y',\overline{y},z\}
\]
restricts what the next letter or the next few letters can be. For instance, $a$ cannot occur after $a$, a sequence of letters following $b'$ must be in $\{a'',\lambda\}\{a,\lambda\}\{c,\lambda\}\{c',\lambda\}\{b''\}$, and so forth. Therefore, when there is no restriction coming from other letters, we are allowed to have any letter in $\Gamma$ from a letter in $\Delta^c\setminus\{d_\l\}=\{a'',b'',\underline{b},c'',y,y''\}$. Achieving each letter in $\Delta^c\setminus\{d_\l\}$ when there is no restriction corresponds to the states $A''$, $B''$, $\underline{B}$, $C''$, $Y$ and $Y''$ respectively. In Table \ref{tab:A.1}, we can see that from each one of these states, we have every transition of letters in $\Gamma$ transiting to a non-jail states. For $\mathcal{L}(M_2')$, since $w$ in $\overline{L}_2$ must begin with a letter in $\Gamma$, it is appropriate to have $A''$ as its initial state.\\

In addition to letters in $\Delta^c\setminus\{d_\l\}$, the only restriction what a letter after $a$ can be is that it cannot be $a$ due to the ninth condition of $L'$. Indeed, every letter in $\Gamma\{a\}$ is a transition from $A$ to a non-jail state. This also verifies that the initial state of $\mathcal{L}(M_1')$ being $A$ is valid. Similarly, we can confirm transitions from $B$, $C$, $D$ and $\underline{D}D$ are all valid as well. From $D$, we can transit to $D_\l$ with $d_\l$, ending $w$ with $dd_\l$ which satisfies the suffix condition of $\overline{L}_i$. The state $CD=\d'(C,d)$ is to avoid having $d_\l$ after having $cd$, as $cdd_\l$ is forbidden in $L'$. Except that we cannot have the transition $d_\l$ from $CD$, other transitions are identical to the ones from $D$, verifying the row of $CD$.\\

From any of the states we have discussed so far, having $d_c'$ will take us to $D_c'[D]$. From here, the only path is to $B_s[D]$, $D_c''$, and then $D_\l$, which gives us the suffix $d_c'b_sd_c''d_\l$. This is one of the acceptable suffixes of $\overline{L}_i$, so the transitions from $D_c'[D]$ and $B_s[D]$ are confirmed.\\

Next, we examine transitions involving states having $[A]$ in their expressions. Except for the initial states of $\mathcal{L}(M_3')$ and $\mathcal{L}(M_4')$, in order to enter states with $[A]$, we always must go to $A'[A]$ first, since $[A]$ indicates that $a'$ has been mentioned and currently unmatched. On the other hand, to exit out of states with $[A]$, we must have $a''$ from either $B_s[A]$, $\overline{B}[A]$, $B'[A,B]$ or $X'[A,X]$ to transit to $B''$, $\overline{B}A''$, $A''[B]$ or $A''[X]$ respectively. Before we arrive at any of $B_s[A]$, $\overline{B}[A]$, $B'[A,B]$ and $X'[A,X]$, we can bypass through states $A[A]$, $A'[A,A]$, $A''[A]$ in this order finite number of times with $A[A]$ being optional. In this case, we can transit out to $B_s[A]$, $\overline{B}[A]$, $B'[A,B]$ or $X'[A,X]$ from $A[A]$ or $A''[A]$, but not $A'[A,A]$. Observing these cases, the set
of the sequence of letters allowed by these transitions is\[
\{a'\}\{aa'a'',a'a''\}^*\{a,\lambda\}\{b_sa'',b'a'',\overline{b}a'',x'a''\}=\{a',a'a\}\{a'a'',a'a''a\}^*\{b_sa'',b'a'',\overline{b}a'',x'a''\},
\]
which is the set of subwords of a word in $L'$ containing $a'$ and $a''$ as it is listed in Condition 3 in the definition of $L'$. We can also verify that the initial states of $\mathcal{L}(M_3')$ and $\mathcal{L}(M_4')$ are appropriate with the prefix conditions of $\overline{L}_3$ and $\overline{L}_4$.\\

We now look at transitions from states with $[B]$. Besides the initial states of $\mathcal{L}(M_5')$ and $\mathcal{L}(M_6')$, two ways to enter states with $[B]$ are either by transiting to $B'[A,B]$ from $A[A]$, $A'[A]$ or $A''[A]$, or by transiting to $B'[B]$ from any other states that allows the transition $b'$. As previously observed, from $B'[A,B]$, we only have one transition available, which is $a''$ to move to $A''[B]$. Now, whether from $A''[B]$ or $B'[B]$, we must arrive at $B''$, $B''[C]$, $B''[D]$ or $DB''$ to exit out of states with $[B]$. Until we arrive at $B''$, we can bypass through $A[B]$ and $C[B]$ in that order. From $A'[B]$, going through $A[B]$ and $C[B]$ are unnecessary. On the other hand, if we are coming from $B'[B]$, stopping at one of them is necessary, since we do not have the transition $b''$ from $B'[B]$. Similarly, we can observe the possibilities to arrive at each one of $B''[C]$, $B''[D]$ and $DB''$ from $A''[B]$ and $B'[B]$ to see all possible subwords induced by the transitions after $b'$ up to $b''$ are in one of the following sets.\[
\{a'',\lambda\}\{a,\lambda\}\{c,\lambda\}\{c',\lambda\}\{b''\},\quad\{a'',\lambda\}\{a,\lambda\}\{db''\},\quad\{a'',\lambda\}\{a,\lambda\}\{c,\lambda\}\{d_c'b''\}.
\]
Having the initial states of $\mathcal{L}(M_5')$ and $\mathcal{L}(M_6')$ as $A[B]$ and $A''[B]$ respectively, we can confirm that $w$ in $\mathcal{L}(M_5')$ satisfies the prefix condition of $\overline{L}_5$ and $w$ in $\mathcal{L}(M_6')$ satisfies the prefix condition of $\overline{L}_6$. Also, obtaining $db''$ takes us from $D[B]$ to $DB''$, and the only transition from there is $d_\l$ to go to $D_\l$, which meets one of the suffix condition of $\overline{L}_i$. Similarly, $d_c'b''$ corresponds to the transition from $D_c'[B,D]$ to $B''[D]$, and this continues with $D_c''$ and then $D_\l$. Hence, we have $d_c'b''d_c''d_\l$, an acceptable suffix of $w$ in $\overline{L}_i$.\\

We move onto the transitions with overlines and underlines. From any state in $P'$, we arrive at $\overline{B}$ with $\overline{b}$ from which, we transit to either $\underline{D}$ or $\underline{\overline{D}}$ with optional paths through $\overline{B}A$ and $\overline{B}C$ in that order. If we come from $A[A]$, $A'[A]$ or $A''[A]$, we transit to $\overline{B}[A]$, and then $\overline{B}A''$ with $a''$. Again, having optional states $\overline{B}A$ and $\overline{B}C$, we arrive at $\underline{D}$ or $\underline{\overline{D}}$. The only transition from $\underline{D}$ and $\underline{\overline{D}}$ is $d$ to $\underline{D}D$ and $\overline{D}D$ respectively. Thus, a sequence of letters constructed by these paths is in\[
\{\overline{b}\}\{a'',\lambda\}\{a,\lambda\}\{c,\lambda\}\{\underline{d}d,\underline{\overline{d}}d\}
\]
which satisfies the seventh condition of $\overline{L}_i$. For $\mathcal{L}(M_9')$ and $\mathcal{L}(M_{10}')$, the initial states are $\overline{B}A$ and $\overline{B}A''$ respectively. Thus, a word in $\mathcal{L}(M_9')$ has a prefix in $\{c,\lambda\}\{\underline{d}d,\underline{\overline{d}}d\}$, and a word in $\mathcal{L}(M_{10}')$ has a prefix in $\{a,\lambda\}\{c,\lambda\}\{\underline{d}d,\underline{\overline{d}}d\}$.\\

Instead of $\overline{b}$, we can also use $\overline{d}$, $\underline{\overline{d}}$ or $\overline{y}$ to enter states with overlines, and arrive at states with underlines later. For now, we examine the cases of $\overline{d}$ and $\underline{\overline{d}}$ to have $\overline{D}$ and $\underline{\overline{D}}$ respectively, and we will discuss the case of $\overline{y}$ later. When we have a transition $\overline{d}$ or $\underline{\overline{d}}$ from a certain state, we arrive at $\overline{D}$ or $\underline{\overline{D}}$ respectively, which are both followed by $\overline{D}D$. From $\overline{D}D$, there are several ways to achieve a state with an underline. First, going to $D\ra\underline{B}$, $D\underline{B}$ and then $D_\l$ gives us either $\overline{d}dd\underline{b}d_\l$ or $\underline{\overline{d}}dd\underline{b}d_\l$, which are appropriate suffixes of a word in $\overline{L}_i$. Similarly, after the optional $C\ra\{\underline{B},\underline{D},\underline{\overline{D}},\underline{X}\}$, we can transit through $D_c'[D]\ra\underline{B}[D]$, $\underline{B}[D]$, $D_c''$ and then $D_\l$ to have a suffix in $\{\overline{d}d,\underline{\overline{d}}d\}\{c,\lambda\}\{d_c'\underline{b}d_c''d_\l\}$, which is also a valid suffix in $\overline{L}_i$. Additionally, there are five more cases. In each case, $C\ra\{\underline{B},\underline{D},\underline{\overline{D}},\underline{X}\}$ is available after $\overline{D}D$ as an optional state. Thus, these five cases are:\begin{itemize}
\item $\underline{B}$ to obtain $\{\overline{d}d,\underline{\overline{d}}d\}\{c,\lambda\}\{\underline{b}\}$.
\item $C'[C]\ra\underline{B}[C]$ and $\underline{B}[C]$ to obtain $\{\overline{d}d,\underline{\overline{d}}d\}\{c,\lambda\}\{c'\underline{b}\}$.
\item $\underline{D}$ and $\underline{D}D$ to obtain $\{\overline{d}d,\underline{\overline{d}}d\}\{c,\lambda\}\{\underline{d}d\}$.
\item $\underline{\overline{D}}$ and $\overline{D}D$ to obtain $\{\overline{d}d,\underline{\overline{d}}d\}\{c,\lambda\}\{\underline{\overline{d}}d\}$.
\item $Z\ra\underline{X}$, $\underline{X}$ and $XD$ to obtain $\{\overline{d}d,\underline{\overline{d}}d\}\{c,\lambda\}\{z\underline{x}d\}$.
\end{itemize}
All of them are acceptable under the seventh condition of $\overline{L}_i$.\\

Let us now look at transitions from states with $[C]$. Entering these states can be done in multiple ways. Namely, from $A[B]$, $A''[B]$, $B'[B]$ and $C[B]$ to $C'[B,C]$ while $b'$ is unmatched, from $A[Y]$, $C[Y]$ and $Y'[Y]$ to $C'[Y,C]$ while $y'$ is unmatched, from $\overline{D}D$ and $C\ra\{\underline{B},\underline{D},\underline{\overline{D}},\underline{X}\}$ to $C'[C]\ra\underline{B}[C]$ while $\overline{d}$ is unmatched, and simply from a state in $P'$ to $C'[C]$. The states immediately after each case are $B''[C]$, $Y''[C]$, $\underline{B}[C]$ and $B_s[C]$ respectively, and they share the same transitions to the same states afterwards. Similar to the case of $[A]$, we can transit through $C[C]$, $C'[C,C]$ and $C''[C]$ in this order finite number of times with $C[C]$ being optional. We can exit out of these states from $C[C]$ by having $c''$ to $C''$, $d_c'$ to $D_c'[C,D]$ or $z$ to $Z[C]$, or from $B''[C]$, $Y''[C]$, $\underline{B}[C]$, $B_s[C]$ and $C''[C]$ by having $d$ to $D[C]$, $d_c'$ to $D_c'[C,D]$ or $z$ to $Z[C]$. Since $Z[C]$ is followed by $ZC''$, if we arrive at $C''$ or $Z[C]$, the sequence of letters constructed by these paths are in\[
\{c'b_s,c'b'',c'\underline{b},c'y''\}\{c'c'',cc'c''\}^*\{c'',cc'',zc'',czc''\},
\]
which obeys the fifth condition of $\overline{L}_i$. Similarly, if we arrive at $D[C]$ or $D_c'[C,D]$, the rest of the run of $M_i'$ are $DC''$ then $D_\l$ and $C''[D]$, $D_c''$ then $D_\l$ respectively. Suffixes constructed by paths are respectively in $\{c'b_s,c'b'',c'\underline{b},c'y''\}\{c'c'',cc'c''\}^*\{dc''d_\l\}$ and $\{c'b_s,c'b'',c'\underline{b},c'y''\}\{c'c'',cc'c''\}^*\allowbreak\{c,\lambda\}\{d_c'c''d_c''d_\l\}$, which satisfy the suffix condition of $\overline{L}_i$.\\

Finally, we observe transitions from the rest of the states, which are related to $X$, $Y$ and $Z$. We first describe four ways to arrive at $XD$ first, and any sequence of letters up to $XD$ obeys the eighth condition of $\overline{L}_i$. Frankly, these four cases are recognized by which state of $Z$, $Z[X]$, $Z\ra\underline{X}$ and $Z[C]$ we obtain by the transition $z$. There is one more state that we can arrive at with $z$, which is $Z[Y]$, but we will look at this case later.\\

The transition $x'$ can occur if and only if $b'$ can happen. From $A[A]$, $A'[A]$ or $A''[A]$, we arrive at $X'[A,X]$, and from any other states that allow the transition $x'$, we transit to $X'[X]$. From $X'[A,X]$, $a''$ is the only transition available to $A''[X]$. Whether from $A''[X]$ or $X'[X]$, we have the optional states $A[X]$ and $C[X]$ in this order. Afterwards, we must arrive at $Z[X]$ which is followed by $X''$, and then $XD$. Therefore, words constructed by these sequences of transitions starting from $x'$ are in $\{x'\}\{a'',\lambda\}\{a,\lambda\}\{c,\lambda\}\{zx''\}$, which is listed under the eighth condition of $\overline{L}_i$. For $\mathcal{L}(M_7')$ and $\mathcal{L}(M_8')$, a word begins with a prefix in $\{c,\lambda\}\{zx''\}$ and $\{a,\lambda\}\{c,\lambda\}\{zx''\}$ respectively, and it is followed by $d$ afterwards.\\

From any state in $P'$, we have the transition $z$ to $Z$ followed by $X$ and $XD$, so this path produces $zxd$. Also, as explained before, having $\overline{d}$ or $\underline{\overline{d}}$ previously allows us to transit to $Z\ra\underline{X}$, $\underline{X}$ and $XD$ to obtain $z\underline{x}d$. The last case, $Z[C]$ can be obtained whenever $C'[C,C]$ is available. Once we arrive at $Z[C]$, we transit to $ZC''$ with $c''$ instead of $x$, then $XD$, giving us $zc''d$. Hence, in all four cases of passing through $Z$, $Z[X]$, $Z\ra\underline{X}$ or $Z[C]$, a sequence of letters up to $XD$ obeys the eighth condition of $\overline{L}_i$.\\

We note that primarily, $x$ after $z$ is a necessary transition, but this $x$ is replaceable with $x''$, $\underline{x}$ or $c''$. In the case $x'$ or $c'$ is unmatched, we have $x''$ or $c''$ respectively instead of $x$. Similarly, if $\overline{d}$ or $\underline{\overline{d}}$ is unmatched with a transition with an underline, then we have $\underline{x}$.\\

From $XD$, there are three transitions available. Having $y$ takes us to $Y$, which completes the sequence of transitions related to $x$, $y$ and $z$. On the other hand, if we transit to $\overline{Y}$ with $\overline{y}$, we have optional $\overline{B}A$ and $\overline{B}C$ before we get to $\underline{D}$ or $\underline{\overline{D}}$, constructing a word in $\{\overline{y}\}\{a,\lambda\}\{c,\lambda\}\{\underline{d}d,\underline{\overline{d}}d\}$. Finally, if we go to $Y'[Y]$ with $y'$, there are a few different ways to achieve $y''$. With optional $A[Y]$, $C[Y]$ and $C'[Y,C]$ in this order, we can arrive at either $Y''$ or $Y''[C]$, giving us $\{y'\}\{a,\lambda\}\{c,\lambda\}\{y'',c'y''\}$. Now, after optional $A[Y]$ and $C[Y]$, it is also possible to have $Z[Y]$ followed by $ZY''$ to have $\{y'\}\{a,\lambda\}\{c,\lambda\}\{zy''\}$. In this case, $y''$ is taking the primary role of $x$, so it goes to $XD$ afterwards, creating another sequence of letters starting with $z$. The other two states after $Y'$ are $D[Y]$, $DY''$ and then $D_\l$, or $D_c'[Y,D]$, $Y''[D]$, $D_c''$ and then $D_\l$, resulting in suffixes $dy''d_\l$ or $d_c'y''d_c''d_\l$ respectively.\\

With any sequence of transitions up to $XD$ combined with a sequence of transitions from $XD$, the constructed subword satisfies the conditions of $\overline{L}_i$.\\

We have examined every sequence of transitions recursively occurring to form a word obeys all conditions of $\overline{L}_i$. Consequently, a word in $\mathcal{L}_(M_i')$ is in $\overline{L}_i$, completing the proof of $\mathcal{L}(M_i')\subseteq\overline{L}_i$.\\

To show that $\overline{L}_i\subseteq\mathcal{L}(M_i')$, we take the same approach as we did in Chapter \ref{chap:4} for $\overline{L}$ and $\mathcal{L}(M)$. That is, suppose $w$ is not in $\mathcal{L}(M_i')$. The only ways $w$ cannot be accepted by $M_i$ are either the run of $M_i$ on $w$ contains the jail state or the last state is not $D_\l$. The latter violates the suffix condition of $\overline{L}_i$, so $w$ is not in $\overline{L}_i$. For the case the run of $M$ on $w$ contains the jail state, we need to show that every transition to the jail state is due to a failure of $w$ to meet one of the conditions of $\overline{L}_i$.\\

For some cases such as $w$ containing the jail state due to $(A,a)$ or $(B,b)$, $w$ must contain a subword that is prohibited in the ninth condition of $\overline{L}_i$. In other cases, which are the majority of the transitions to the jail state, $w$ must disobey at least one of the second to the eighth conditions of $\overline{L}_i$. For instance, having $a''$ with any state $q$ without the notation of $[A]$ implies that there is no $a'$ previously. Thus, this violates the third condition of $\overline{L}_i$. As we can see, it is a straightforward exercise (but long and tedious one) to verify that every transition to the jail state violates some conditions of $\overline{L}_i.$\\

Consequently, this completes the proof of $\mathcal{L}(M_i')=\overline{L}_i$.\hfill$\blacksquare$\\

With Proposition \ref{prop:6.1} and \ref{prop:6.2}, we are finally at the place to derive the generating function for all simple permutations of length 4 or more in $\A'$.\\
\\
\textit{Proof of Theorem \ref{thm:6.1}.} Let us apply the transfer matrix method to each $\mathcal{L}(M_i')$ ($1\leq i\leq 10$). Since each $M_i'$ shares the same transitions with different initial states, we only need to provide one adjacency matrix. The adjacency matrix for each $M_i'$ is extremely large, so it is given in Appendix \ref{app:A}.\\

In this matrix, there are three weights, $x$, $\bar{F}$ and $\bar{G}$. These weights are used for inflation in the next section. In order to we replace all weights $\bar{F}$ and $\bar{G}$ with $x$. Let $P'$ be such adjacency matrix. We denote by $(I-P)^{-1}_{q_1,q_2}$ the $(q_1,q_2)$-entry of $(I-P)^{-1}$. For each $\mathcal{L}(M_i')$, we examine $(I-P)^{-1}_{q_i,D_\l}$ where $q_i$ is the initial state for $M_i'$. Each one gives the generating function for $n$-letter words in $\mathcal{L}(M_i')$, and hence, for $n$-letter words in $\overline{L}_i$. A word in each $\overline{L}_i$, is missing a certain prefix, so we need to multiply either $x^2$, $x^3$ or $x^4$ accordingly to count all associated words in $L'$. We then add all of these generating functions to obtain the generating function for $L'$, which enumerates all permutations in $H'$. Thus,\[\begin{array}{rcl}
f_{H'} & = & x^2\cdot(I-P)^{-1}_{A,D_\l}+x^4\cdot(I-P)^{-1}_{A'',D_\l}+x^3\cdot(I-P)^{-1}_{A[A],D_\l}+x^4\cdot(I-P)^{-1}_{A''[A],D_\l}\vhhh\\
& + & x^3\cdot(I-P)^{-1}_{A[B],D_\l}+x^4\cdot(I-P)^{-1}_{A''[B],D_\l}+x^3\cdot(I-P)^{-1}_{A[X],D_\l}+x^4\cdot(I-P)^{-1}_{A''[X],D_\l}\vhhh\\
& + & x^3\cdot(I-P)^{-1}_{\overline{B}A,D_\l}+x^4\cdot(I-P)^{-1}_{\overline{B}A'',D_\l}\vhhh\\
 & = & \ds\frac{x^4 \left(x^{10}+7 x^9+18 x^8+23 x^7+16 x^6+10 x^5+12 x^4+9 x^3+2 x^2+1\right)}{(x+1) \left(2 x^9+12 x^8+16 x^7+3 x^6-11 x^5-5 x^4-3 x^2-3
   x+1\right)}.\end{array}\]
By doubling this result, we obtain the desired generation function $f_{\textrm{Si}(\A')\setminus\mathcal{S}_2}$.\hfill$\blacksquare$
\section{Enumeration of the whole class $\A'$}
In order to finish the enumeration of the whole class, we show that every simple permutation in $\A'$ satisfies the hypothesis of Proposition \ref{prop:2.4}, just as we did in Chapter \ref{chap:4}.\\

Before we state and prove the statement, we make a few important notes. For $\pi\in N$, values in $\LRmax(\pi)$ have positions in $[1,\pi^{-1}(n)]$. In particular, except the value $n$, a value of $\pi$ is in $\LRmax(\pi)$ if and only if it is 1 of 1, 2 of 21 in 231-value chains, or 2 or 3 of 231 in 231-value chains in Equation \ref{eqn:5.2}. This implies that $\LRmax(\pi)\setminus\{n\}$ is the set of values that are encoded as $a$ or $a''$ by \textsf{N-ENCODE}. By the reverse complement symmetry, we also know that $\RLmin(\pi)\setminus\{1\}$ is the set of values that are encoded as $c$ or $c'$ by \textsf{N-ENCODE}. Hence, values $\pi(i)$ that are covered in the fourth condition of the latter four must be playing the roles of 1 of 21 or 231 in 231-value chains in Equation \ref{eqn:5.2}, 1 of 21 or 312 in 312-value chains in Equation \ref{eqn:5.3}, 1 of 1 in Equation \ref{eqn:5.1} which is a scissor of a value chain, or 1 or 2 of 12 in Equation \ref{eqn:5.1}. Respectively, these are the values encoded as $a'$, $c''$, $b_s$, $b'$ and $b''$ by \textsf{N-ENCODE}.\\

We are now ready to prove the following lemma.
\begin{lemma}\label{lem:6.4}
\textit{Let $\pi$ be a simple permutation of extreme pattern 2413 or 3142 whose length is $n$. The condition $\a=\s[\s_1,\ldots,\s_n]\in\A'$ is equivalent to the condition stating that for all $i$ with $1\leq i\leq n$,}\begin{itemize}
\item \textit{if $\pi(i)\in\LRmax(\pi)$, then $\s_i\in\Av(4123,4213,4132)$,}
\item \textit{if $\pi(i)\in\RLmin(\pi)$, then $\s_i\in\Av(2341,3241,2431)$,}
\item \textit{if $\pi(i)$ is 1 of 1 in Equation \ref{eqn:5.1} and it is not a scissor of either value chain, then $\s_i\in\Av(123,213,132)$, and}
\item \textit{otherwise, $\s_i\in\Av(12,21)$.}
\end{itemize}
\end{lemma}
\textit{Proof.} We only consider the case of $\pi$ is of extreme pattern 2413, since we can apply the inverse symmetry to prove the case of extreme pattern 3142.\\

First, we show that the first condition implies the latter four conditions. Suppose the latter condition is false. That is, at least one of the above four conditions is not met. Assume it is the first one. Then there exists $i$ ($1\leq i\leq n$) such that $\pi(i)\in\LRmax(\pi)$ and $\s_i$ contains at least one of 4123, 4213 and 4132. Now, whichever point $(i,\pi(i))$ is, the point with the value 1 is located to the right. Thus, the value 1 of $\a$ together with the subsequence corresponding to 4123, 4213 or 4132 of $\a$ within the inflation $\s_i$, $\a$ contains 52341, 53241 or 52431 respectively, so $\a\notin\A'$. If the second condition is not met, we can apply the reverse complement argument of the previous one to show $\a$ contains one of 52341, 53241 and 52431. If the third condition is false, say at least one of 123, 213 and 132 is contained in $\s_i$ for some $i$, then with the values $|\a|$ and $1$, $\a$ contains 52341, 53241 or 52431, $\a\notin\A'$. Finally, for the last condition, each of these values is either 2 or 3 of 4231 pattern in $\pi$. Hence, inflating with $\s_i$ which contains 12 or 21 will cause $\a$ to contain 52341, 53241 or 52431, so again, $\a\notin\A'$. Consequently, $\a=\pi[\s_1,\ldots,\s_n]\in\A'$ implies the latter condition.\\

Next, assume a permutation $\a=\pi[\s_1,\ldots,\s_n]$ where $\pi\in H'$ is not in $\A'$. Thus, $\a$ contains at least one permutation $\b$ in the basis. Since $\pi$ avoids every permutation in the basis, it means there exist $\s_{i_1},\ldots,\s_{i_k}$ (for all $j\in\{1,\ldots,k\}$, $1\leq i_j\leq n$) such that $\a$ contains $\b$ within the union of subintervals corresponding to $\s_{i_1},\ldots,\s_{i_k}$.\\

Because every permutation in $\{35142,42513,351624\}$ is simple, so we only need to consider $\b\in\{52341,53241,52431\}$. For now, suppose $\b=52341$. If there exist $\s_{i_1},\s_{i_2},\s_{i_3},\s_{i_4},\s_{i_5}$ such that each point of $\b$ is contained in intervals corresponding to $\s_{i_1},\s_{i_2},\s_{i_3},\s_{i_4},\s_{i_5}$ respectively, then $\b\cont\pi$ which cannot be true. On the other hand, if there exists a single $\s_i$ such that $\b\cont\s_i$, then one of the latter conditions is false, since $\b$ contains whatever $\s_i$ cannot contain. Hence, there exist two, three or four subintervals of $\a$ such that the containment of $\b$ is involved in.\\

We first assume that a containment of $\b$ occurs among four subintervals of $\a$. Let $\s_{i_1},\s_{i_2},\s_{i_3},\s_{i_4}$ with $i_1<i_2<i_3<i_4$ be the permutations from the expression $\a=\pi[\s_1,\ldots,\s_n]$ which correspond to each subinterval. By the pigeonhole principle, there is one subinterval containing two values of $\b$ which have consecutive positions. Due to inflation, these two values cannot be 52 and 41 of $\b$. Hence, it is either 23 or 34, and each case implies $12\cont\s_{i_2}$ or $12\cont\s_{i_3}$ respectively. Suppose the case of 23 contained in one subinterval, so we have $12\cont\s_{i_2}$. Note that $\pi(i_2)$ cannot be in $\LRmax(\pi)$ because then we do not have any value for $\pi(i_1)$ to the left of $\pi(i_2)$. Similarly, $\pi(i_2)\notin\RLmin(\pi)$. Also, if $\pi(i_2)$ plays the role of 1 of 21 or 231 in 231-value chains in Equation \ref{eqn:5.2} or 1 of 21 or 312 in 312-value chains in Equation \ref{eqn:5.3}, then the fourth condition is not met. Therefore, suppose $i_2\in(\pi^{-1}(n),\pi^{-1}(1))$. If $\pi(i_2)$ plays the role of 1 or 2 of 12 in Equation \ref{eqn:5.1} or 1 of 1 in Equation \ref{eqn:5.1} that is a scissor of some value chain, then again, the fourth condition is not satisfied. Thus, $\pi(i_2)$ plays the roles of 1 of 1 in Equation \ref{eqn:5.1} that is not a scissor of a value chain. Then $i_3,i_4\in[\pi^{-1}(1),n]$ with $\pi(i_4)<\pi(i_2)<\pi(i_3)$, but this means $\pi(i_2)$ is a scissor of a 312-value chain, so we have a contradiction. We achieve same results by assuming 34 is contained in one subinterval. Consequently, if a containment of $\b$ occurs among four subintervals of $\a$, then the latter conditions are false.\\

Next, suppose a containment of $\b$ involves three subintervals of $\a$. Let $\s_{i_1},\s_{i_2},\s_{i_3}$ with $i_1<i_2<i_3$ be the permutations from the expression $\a=\pi[\s_1,\ldots,\s_n]$ which correspond to each subinterval. In this case, the only possibility is that the middle interval corresponding to $\s_{i_2}$ contains the positions of 234, and each of the other two interval contains the positions of 5 and 1. Hence, we have $123\cont\s_{i_2}$. As before, $\pi(i_2)$ cannot be in $\LRmax(\pi)$ and $\RLmin(\pi)$. Now, if $\pi(i_2)$ corresponds to 1 of 1 in Equation \ref{eqn:5.1} that is not a scissor of a value chain, then the third condition is disobeyed. Otherwise, the fourth condition is false, since $12\cont 123\cont\s_{i_2}$.\\

Finally, suppose the containment of $\b$ is shared by two subintervals of $\a$, say the ones of $\s_{i_1}$ and $\s_{i_2}$ with $i_1<i_2$. Then this must imply $4123\in\s_{i_1}$ or $2341\in\s_{i_2}$. Suppose $4123\in\s_{i_1}$. Then $\pi(i_1)\notin\RLmin(\pi)$ since there is no value for $\pi(i_2)$ to the right of $\pi(i_1)$. However, because 12 and 123 are contained in 4123, this implies $\s_1$ disobeys one of the latter conditions. We obtain the same result for $2341\in\s_{i_2}$. So once again, the latter conditions are false.\\

Proving for the cases of $\b=53241$ and $\b=52431$ are very similar. With every observation we made, the latter conditions imply $\a=\pi[\s_1,\ldots,\s_n]\in\A'$. Consequently, those two statements are equivalent.\hfill$\blacksquare$\\

By using Lemma \ref{lem:6.4}, we claim that every simple permutation in $H'$ also satisfies the hypothesis of Proposition \ref{prop:2.4}. Let $\s$, $\tau$ be simple permutations of extreme pattern 2413 and 3142 respectively where $|\s|=m$ and $|\tau|=n$. We examine $\s\nwsum_1^1$ as an example. Let $i$ be the position of $m$ in $\s$. In $\s\nwsum_1^0\tau$, for every position $k$ except for $i$, $\s'(k)=\s(k)$. Because all values coming from $\tau$ are greater than $\s'(k)$, subclasses of inflation for each $\s'(k)$ are exactly the same for $\s(k)$ of $\s$. Similarly, for $\tau'(k)$ with $3\leq k\leq n$, as $\s'(i)$ taking the place of $\tau(1)$, subclasses of inflation for each $\tau'(k)$ are exactly the same for $\tau(k)$ of $\tau$. For $\s'(i)$, since it is now in $\LRmax(\s\nwsum_1^0\tau)$, its subclass of inflation is $\Av(4123,4213,4132)$, which is the same as the one for $\s(i)$ in $\s$. Hence, for every value of $\s\nwsum_1^0\tau$, its subclass of inflation remains the same. We obtain the same result for $\s\nwsum_1^0\tau$.\\

For $(x,y)\in\{(2,0),(2,1),(3,0),(4,0)\}$, there is another value of $\s$ that is shifted upward in $\s\nwsum_x^y\tau$. In particular, for $(x,y)=(2,0)$ and $(x,y)=(2,1)$, it is $\s(i+2)$, for $(x,y)=(3,0)$, it is $\s(m-2)$, and for $(x,y)=(4,0)$, it is $\s(s)$ where $s$ is the position of the value $m-2$. All of these values correspond to either 2 or 3 of 4231 pattern in $\s$, so their subclass of inflation is $\Av(12,21)$. After they are shifted by $\nwsum_x^y$, they are still playing the roles of 2 or 3, so their subclass of inflation do not change. In addition, just like the case of $\nwsum_1^1$, other values of $\s$ and values of $\tau$ keep their subclasses of inflation when $\s$ and $\tau$ are operated by $\nwsum_x^y$.\\

We can establish the same result for SE glue sums. Thus looking at how \textsf{W-COMBINE} and \textsf{AFFIX-CONVERT} replaces letters, we have the following proposition as a consequence of Lemma \ref{lem:6.4}.
\begin{proposition}
\textit{Let $\pi$ be a simple permutation in $H'$. The condition $\a=\pi[\s_1,\ldots,\s_n]\in\A'$ is equivalent to the condition stating that for all $i$ with $1\leq i\leq n$,}\begin{itemize}
\item \textit{if $\pi(i)$ is encoded by $\phi'$ as a letter in $\{a,a'',c,c',d,d_a'',d_c',d_\l,z\}$, then $\s_i\in\Av(4123,4213,4132)$ or $\s_i\in\Av(2341,3241,2431)$ (depending on the specific letter and whether $\pi(i)$ originally belonged to $N$ or $S$),}
\item \textit{if $\pi(i)$ is encoded as $b$, then $\s_i\in\Av(123,213,132)$, and}
\item \textit{if $\pi(i)$ is encoded as a letter in $\{a',b_s,b',b'',\underline{b},\overline{b},c'',\underline{d},\overline{d}, \underline{\overline{d}},x,x',x'',\underline{x},y,y',y'',\overline{y}\}$, then $\s_i\in\Av(12,21)$.}
\end{itemize}
\end{proposition}
\vhhh

The generating functions for $\Av(4123,4213,4132)$ and $\s_i\in\Av(2341,3241,2431)$ are exactly the same by reverse complement, as stated in Lemma \ref{lem:3.2}. For their inflation, we need the generating function $\bar{G}$. The generating function for $\Av(123,213,132)$ was also previously explained in Chapter \ref{chap:3} as Lemma \ref{lem:3.1}. We denote by $F$ the generating function for this class, and use $\bar{F}$ for inflation. Finally, $\Av(12,21)=\{\e,1\}$, so $f_{\Av(12,21)}=1+x$, implying $\bar{f}_{\Av(12,21)}=x$.\\

We now revisit the adjacency matrix $P'$. By replacing weights of transitions in $\{a,a'',c,c',\allowbreak d,d_a'',d_c',d_\l,z\}$ with $\bar{G}$ and the transition $b$ with $\bar{F}$, we obtain the adjacency matrix $\hat{P}'$ shown in Table \ref{tab:A.2}. The generating function for permutations in $\A'$ that are obtained by inflation of simple permutations in $H'$ is
\begin{eqnarray}
f_{\textrm{ifl}(H')} & = & \bar{G}^2\cdot(I-\hat{P})^{-1}_{A,D_\l}+x^2\cdot\bar{G}^2\cdot(I-\hat{P})^{-1}_{A'',D_\l}+x\cdot\bar{G}^2\cdot(I-\hat{P})^{-1}_{A[A],D_\l}\nn\vhhh\\
& + &x^2\cdot\bar{G}^2\cdot(I-\hat{P})^{-1}_{A''[A],D_\l}+x\cdot\bar{G}^2\cdot(I-\hat{P})^{-1}_{A[B],D_\l}+x^2\cdot\bar{G}^2\cdot(I-\hat{P})^{-1}_{A''[B],D_\l}\nn\vhhh\\
& + & x\cdot\bar{G}^2\cdot(I-\hat{P})^{-1}_{A[X],D_\l}+x^2\cdot\bar{G}^2\cdot(I-\hat{P})^{-1}_{A''[X],D_\l}+x\cdot\bar{G}^2\cdot(I-\hat{P})^{-1}_{\overline{B}A,D_\l}\nn\vhhh\\
& + & x^2\cdot\bar{G}^2\cdot(I-\hat{P})^{-1}_{\overline{B}A'',D_\l}.\label{eqn:6.1}
\end{eqnarray}
As usual, we multiply 2 to include the inverse case, so we have $f_{\textrm{ifl}(\textrm{Si}(\A')\setminus\mathcal{S}_2)}=2f_{\textrm{ifl(H')}}$.\\

For the inflation of $\pi=21$, we need to inflate $\pi(1)=2$ with a skew-indecomposable permutation $\s_1$. We claim that $\a=21[\s_1,\s_2]\in\A'$ with skew-indecomposable $\s_1$ is equivalent to the condition that $\s_1\in\Av(4123,4213,4132)$ and $\s_1$ is skew-indecomposable, and $\s_2\in\Av(2341,3241,2431)$. The condition of $\s_1$ being skew-indecomposable cannot be dropped to enforce the uniqueness of inflation. It is clear that either of $\b_1\cont\s_1$ for any $\b_1\{4123,4213,4132\}$ or $\b_2\cont\s_2$ for any $\b_2\{2341,3241,2431\}$ implies $\a\notin\A'$. Suppose that $\a\notin\A'$. Then $\b\in\a$ for some $\b\in\{52341,53241,52431,35142,42513,351624\}$. Since 35142, 42513 and 351624 are simple, if $\b\in\{35142,42513,351624\}$, then $\b\cont\s_1$ or $\b\cont\s_2$. In either case, either one of 4123, 4213, 4132 is contained in $\s_1$ or either one of 2341, 3241, 2431 is contained in $\s_2$, so the second condition is not met. If $\b=52341$, then it is immediate that $4123\cont\s_1$ or $2341\cont\s_2$. We obtain the similar results for $\b=53241$ and $\b=52431$. Thus, the condition $\a=21[\s_1,\s_2]\in\A'$ with skew decomposable $\s_1$ and the condition $\s_1\in\Av(4123,4213,4132)$ where $\s_1$ is skew-indecomposable are equivalent.\\

By Proposition \ref{prop:2.4},
\[\bar{f}_{\textrm{ifl}(21)}=\bar{f}_{\Av(4123,4213,4132)}^\ominus\cdot \bar{f}_{\Av(2341,3241,2431)}.
\]
By Lemma \ref{lem:3.3}, we know $\bar{f}_{\Av(4123,4213,4132)}^\ominus = (1-x-x^2)\bar{G}$. With $\bar{f}_{\Av(2341,3241,2431)}=\bar{G}$, we obtain\[
f_{\textrm{ifl}(21)}=(1-x-x^2)\bar{G}^2.
\]\vspace{0.05in}

The last case we need to consider is $\pi=12$. We can inflate both 1 and 2 by any permutations $\s_1$ and $\s_2$ of $\A'$ itself, provided that $\s_1$ is a sum-indecomposable permutation in $\A'$. Three cases for $\s_1$ being sum-indecomposable are $\s_1=1$, $\s_1$ is skew-decomposable, or $\s_1$ is an inflated permutation of $\pi\in\textrm{Si}(\A')\setminus\mathcal{S}_2$. Generating functions for each case are $x$, $f_{\textrm{ifl}(21)}$ and $f_{\textrm{ifl}(\textrm{Si}(\A')\setminus\mathcal{S}_2)}$ respectively. Thus,
\[
f_{\textrm{ifl}(12)}=(x+f_{\textrm{ifl}(21)}+f_{\textrm{ifl}(\textrm{Si}(\A)\setminus\mathcal{S}_2)})\cdot\bar{f}_{\A'}.
\]

Consequently, the generating function for $\A'$ satisfies the functional equation\\
\[\begin{array}{rcl}
f_{\A'} & = & 1+x+f_{\textrm{ifl}(12)}+f_{\textrm{ifl}(21)}+f_{\textrm{ifl}(\textrm{Si}(\A')\setminus\mathcal{S}_2)}\\
 & = & 1+x+(x+f_{\textrm{ifl}(21)}+f_{\textrm{ifl}(\textrm{Si}(\A')\setminus\mathcal{S}_2)})\cdot \bar{f}_\A'+f_{\textrm{ifl}(21)}+f_{\textrm{ifl}(\textrm{Si}(\A')\setminus\mathcal{S}_2)}.
 \end{array}
\]
\\
Finally, with $\bar{f}_{\A'}=f_{\A'}-1$, we solve for $f_{\A'}$. Then, we obtain\\
\[\begin{array}{rcl}
f_{\A'} & = & \ds\frac{1}{1-x-f_{\textrm{ifl}(21)}-f_{\textrm{ifl}(\textrm{Si}(\A')\setminus\mathcal{S}_2)}}.
\end{array}
\]
We substitute $f_{\textrm{ifl}(21)}$ with $(1-x-x^2)\bar{G}^2$ and $f_{\textrm{ifl}(\textrm{Si}(\A')\setminus\mathcal{S}_2)}=2f_{\textrm{inf}(H')}$ with Equation \ref{eqn:6.1}. After we simplify the whole expression of $f_{\A'}$ as well as Equation \ref{eqn:6.1}, we obtain the final desired result as the following.
\begin{theorem}
\label{thm:6.2}
\textit{The generating function for the class $\A'$ is defined by}\\
\[
f_{\A'}=\frac{\sum_{i=0}^{5} a_i \bar{G}^i}{\sum_{i=0}^{6} b_i \bar{G}^i}
\]
\\
\textit{where $\bar{G}=G-1$ and $G$ is the generating function satisfying the equation}\\
\[
G=1+\frac{xG}{1-xG^2},
\]
\textit{and}\\
\[
\begin{array}{l}
a_0=-1 + 14 x - 39 x^2 + 28 x^3 + 9 x^4 - 11 x^5 + x^6,\\
a_1=-12 + 81 x - 100 x^2 + 15 x^3 + 46 x^4 - 19 x^5,\\
a_2=-8 + 35 x - 20 x^2 - 25 x^3 + 31 x^4 - 6 x^5 - x^6,\\
a_3=7,\qquad a_4=1,\qquad a_5=-2.
\end{array}\]
\end{theorem}
\[
\begin{array}{l}
b_0=-1 + 57 x - 125 x^2 + 143 x^3 - 48 x^4 - 64 x^5 + 51 x^6 - 2 x^8,\\
b_1=-54 + 260 x - 386 x^2 + 250 x^3 + 81 x^4 - 226 x^5 + 74 x^6 + 15 x^7 - 3 x^8,\\
b_2=-18 + 114 x - 104 x^2 - 22 x^3 + 148 x^4 - 123 x^5 + 11 x^6 + 14 x^7 - x^8,\\
b_3=24,\qquad b_4=-2,\qquad b_5=-5,\qquad b_6=1.
\end{array}
\]
The first several terms in power series expression are\[
f_{\A'}=1+x+2x^2+6x^3+24x^4+115x^5+607x^6+3370x^7+19235x^8+111571x^9+\cdots
\]
\section{Conclusions}
As authors noted in \cite{AB2013}, the technique they used to enumerate $\A$ was indeed applicable to enumerate the permutations indexing local complete intersection Schubert varieties. Although the whole process turned out to be extremely complicated, as of now, examining the geometric structures and constructing an encoding to apply the transfer matrix method is the only known way to enumerate this class. There may be a much more efficient method by closely learning the class $\A'$ and observing Rothe diagrams of permutations described in \cite{UW2013}.\\

Lastly, this dissertation owes huge thanks to PermLab \cite{PermLab1.0} and its author, Michael Albert. Without this extremely efficient program, viewing structures of permutations in $\A'$ would have been impossible.
\newpage
\phantomsection
\addcontentsline{toc}{chapter}{Appendices}
\chapter*{Appendices}
\setcounter{section}{0}
\renewcommand{\thesection}{\Alph{section}}
\renewcommand{\theHsection}{appendixsection.\Alph{section}}
\renewcommand{\thetable}{\Alph{section}.\arabic{table}}
\renewcommand{\theHtable}{appendixtable.\Alph{section}.\arabic{table}}
\section{Transitions of $M'_i$ ($1\leq i\leq 10$) and adjacency matrix associated with $M'_i$}\label{app:A}
Two large tables which are referred in Chapter \ref{chap:6} are placed in this section. Table \ref{tab:A.1} shows the description of $\d'$, the transition function for the automaton $M'_i$ for each $i$ with $1\leq i\leq 10$. Note that Jail states and transitions to them are omitted.\\

Table \ref{tab:A.2} shows the associated adjacency matrix $\hat{P}'$ with weights $x$, $\bar{F}$ and $\bar{G}$. Since the matrix is extremely large, it is presented by dividing into eight sub-matrices $A$ through $H$ in alphabetical order, where\[
\hat{P}'=\left[\begin{array}{c|c|c|c}
A & B & C & D\\
\hline
E & F & G & H
\end{array}\right].
\]\\

As it was described in Chapter \ref{chap:6}, in order to obtain $P'$, the adjacency matrix for each $M_i'$, we need replace all weights $\bar{F}$ and $\bar{G}$ with $x$. Let $P'$ be such adjacency matrix.\\
\begin{landscape}
\begin{table}
\begin{footnotesize}
\[
\begin{array}{C{\mycolwd}|C{\mycolwd}C{\mycolwd}C{\mycolwd}C{\mycolwd}C{\mycolwd}C{\mycolwd}C{\mycolwd}}
 & a & a' & a'' & b & b' & b'' & b_s \\
\hline
A & & A'[A] & & B & B'[B]\\
A'' & A & A'[A] & & B & B'[B]\\
A[A] & & A'[A,A] & & & B'[A,B] & & B_s[A]\\
A'[A] & A[A] & A'[A,A] & & & B'[A,B] & & B_s[A]\\
A''[A] & A[A] & A'[A,A] & & & B'[A,B] & & B_s[A]\\
\hline
\end{array}
\]
\\
\[
\begin{array}{C{\mycolwd}|C{\mycolwd}C{\mycolwd}C{\mycolwd}C{\mycolwd}C{\mycolwd}C{\mycolwd}C{\mycolwd}}
 & \underline{b} & \overline{b}& c & c' & c'' & d & d_c' \\
\hline
A & & \overline{B} & C & C'[C] & & D & D_c'[D]\\
A'' & & \overline{B} & C & C'[C] & & D & D_c'[D]\\
A[A] & & \overline{B}[A]\\
A'[A] & & \overline{B}[A]\\
A''[A] & & \overline{B}[A]\\
\hline
\end{array}
\]
\\
\[
\begin{array}{C{\mycolwd}|C{\mycolwd}C{\mycolwd}C{\mycolwd}C{\mycolwd}C{\mycolwd}C{\mycolwd}C{\mycolwd}}
 & d_c'' & \underline{d} & \overline{d} & \underline{\overline{d}} & d_\l & x & x' \\
\hline
A & & & \overline{D} & & & & X'[X]\\
A'' & & & \overline{D} & & & & X'[X]\\
A[A] & & & & & & & X'[A,X]\\
A'[A] & & & & & & & X'[A,X]\\
A''[A] & & & & & & & X'[A,X]\\
\hline
\end{array}
\]
\\
\[
\begin{array}{C{\mycolwd}|C{\mycolwd}C{\mycolwd}C{\mycolwd}C{\mycolwd}C{\mycolwd}C{\mycolwd}C{\mycolwd}}
 & x'' & \underline{x} & y & y' & y'' & \overline{y} & z \\
\hline
A & & & & & & & Z\\
A'' & & & & & & & Z\\
A[A]\\
A'[A]\\
A''[A]\\
\hline
\end{array}
\]
\caption{Partial state diagram of $M'$.}
\label{tab:A.1}
\end{footnotesize}
\end{table}
\end{landscape}
\newpage
\begin{landscape}
\begin{footnotesize}
\[
\begin{array}{C{\mycolwd}|C{\mycolwd}C{\mycolwd}C{\mycolwd}C{\mycolwd}C{\mycolwd}C{\mycolwd}C{\mycolwd}}
 & a & a' & a'' & b & b' & b'' & b_s \\
\hline
A[B] & & & & & & B'' & \\
A''[B] & A[B] & & & & & B'' & \\
A[X]\\
A''[X]\\
A[Y]\\
\hline
\end{array}
\]
\\
\[
\begin{array}{C{\mycolwd}|C{\mycolwd}C{\mycolwd}C{\mycolwd}C{\mycolwd}C{\mycolwd}C{\mycolwd}C{\mycolwd}}
 & \underline{b} & \overline{b}& c & c' & c'' & d & d_c' \\
\hline
A[B] & & & C[B] & C'[B,C] & & D[B] & D_c'[B,D] \\
A''[B] & & & C[B] & C'[B,C] & & D[B] & D_c'[B,D] \\
A[X] & & & C[X] & & & & \\
A''[X] & & & C[X] & & & & \\
A[Y] & & & C[Y] & C'[Y,C] & & D[Y] & D_c'[Y,D] \\
\hline
\end{array}
\]
\\
\[
\begin{array}{C{\mycolwd}|C{\mycolwd}C{\mycolwd}C{\mycolwd}C{\mycolwd}C{\mycolwd}C{\mycolwd}C{\mycolwd}}
 & d_c'' & \underline{d} & \overline{d} & \underline{\overline{d}} & d_\l & x & x' \\
\hline
A[B]\\
A''[B]\\
A[X]\\
A''[X]\\
A[Y]\\
\hline
\end{array}
\]
\\
\[
\begin{array}{C{\mycolwd}|C{\mycolwd}C{\mycolwd}C{\mycolwd}C{\mycolwd}C{\mycolwd}C{\mycolwd}C{\mycolwd}}
 & x'' & \underline{x} & y & y' & y'' & \overline{y} & z \\
\hline
A[B]\\
A''[B]\\
A[X] & & & & & & & Z[X]\\
A''[X] & & & & & & & Z[X]\\
A[Y] & & & & & Y'' & & Z[Y]\\
\hline
\end{array}
\]
\newpage
\[
\begin{array}{C{\mycolwd}|C{\mycolwd}C{\mycolwd}C{\mycolwd}C{\mycolwd}C{\mycolwd}C{\mycolwd}C{\mycolwd}}
 & a & a' & a'' & b & b' & b'' & b_s \\
\hline
A'[A,A] & & & A''[A] & & & &\\
B & A & A'[A] & & & B'[B]\\
B'' & A & A'[A] & & B & B'[B]\\
\underline{B} & A & A'[A] & & B & B'[B]\\
\overline{B} & \overline{B}A\\
\hline
\end{array}
\]
\\
\[
\begin{array}{C{\mycolwd}|C{\mycolwd}C{\mycolwd}C{\mycolwd}C{\mycolwd}C{\mycolwd}C{\mycolwd}C{\mycolwd}}
 & \underline{b} & \overline{b}& c & c' & c'' & d & d_c' \\
\hline
A'[A,A]\\
B & & \overline{B} & C & C'[C] & & D & D_c'[D]\\
B'' & & \overline{B} & C & C'[C] & & D & D_c'[D]\\
\underline{B} & & \overline{B} & C & C'[C] & & D & D_c'[D]\\
\overline{B} & & & \overline{B}C\\
\hline
\end{array}
\]
\\
\[
\begin{array}{C{\mycolwd}|C{\mycolwd}C{\mycolwd}C{\mycolwd}C{\mycolwd}C{\mycolwd}C{\mycolwd}C{\mycolwd}}
 & d_c'' & \underline{d} & \overline{d} & \underline{\overline{d}} & d_\l & x & x' \\
\hline
A'[A,A]\\
B & & & \overline{D} & & & & X'[X]\\
B'' & & & \overline{D} & & & & X'[X]\\
\underline{B} & & & \overline{D} & & & & X'[X]\\
\overline{B} & & \underline{D} & & \underline{\overline{D}}\\
\hline
\end{array}
\]
\\
\[
\begin{array}{C{\mycolwd}|C{\mycolwd}C{\mycolwd}C{\mycolwd}C{\mycolwd}C{\mycolwd}C{\mycolwd}C{\mycolwd}}
 & x'' & \underline{x} & y & y' & y'' & \overline{y} & z \\
\hline
A'[A,A]\\
B & & & & & & & Z\\
B'' & & & & & & & Z\\
\underline{B} & & & & & & & Z\\
\overline{B}\\
\hline
\end{array}
\]
\newpage
\[
\begin{array}{C{\mycolwd}|C{\mycolwd}C{\mycolwd}C{\mycolwd}C{\mycolwd}C{\mycolwd}C{\mycolwd}C{\mycolwd}}
 & a & a' & a'' & b & b' & b'' & b_s \\
\hline
\overline{B}A\\
\overline{B}A'' & \overline{B}A\\
\overline{B}C\\
B_s[A] & & & A''\\
\overline{B}[A] & & & \overline{B}A''\\
\hline
\end{array}
\]
\\
\[
\begin{array}{C{\mycolwd}|C{\mycolwd}C{\mycolwd}C{\mycolwd}C{\mycolwd}C{\mycolwd}C{\mycolwd}C{\mycolwd}}
 & \underline{b} & \overline{b}& c & c' & c'' & d & d_c' \\
\hline
\overline{B}A & & & \overline{B}C\\
\overline{B}A'' & & & \overline{B}C\\
\overline{B}C\\
B_s[A]\\
\overline{B}[A]\\
\hline
\end{array}
\]
\\
\[
\begin{array}{C{\mycolwd}|C{\mycolwd}C{\mycolwd}C{\mycolwd}C{\mycolwd}C{\mycolwd}C{\mycolwd}C{\mycolwd}}
 & d_c'' & \underline{d} & \overline{d} & \underline{\overline{d}} & d_\l & x & x' \\
\hline
\overline{B}A & & \underline{D} & & \underline{\overline{D}}\\
\overline{B}A'' & & \underline{D} & & \underline{\overline{D}}\\
\overline{B}C  & & \underline{D} & & \underline{\overline{D}}\\
B_s[A]\\
\overline{B}[A]\\
\hline
\end{array}
\]
\\
\[
\begin{array}{C{\mycolwd}|C{\mycolwd}C{\mycolwd}C{\mycolwd}C{\mycolwd}C{\mycolwd}C{\mycolwd}C{\mycolwd}}
 & x'' & \underline{x} & y & y' & y'' & \overline{y} & z \\
\hline
\overline{B}A\\
\overline{B}A''\\
\overline{B}C\\
B_s[A]\\
\overline{B}[A]\\
\hline
\end{array}
\]
\newpage
\[
\begin{array}{C{\mycolwd}|C{\mycolwd}C{\mycolwd}C{\mycolwd}C{\mycolwd}C{\mycolwd}C{\mycolwd}C{\mycolwd}}
 & a & a' & a'' & b & b' & b'' & b_s \\
\hline
B'[B] & A[B]\\
B_s[C]\\
B''[C]\\
\underline{B}[C]\\
B_s[D]\\
\hline
\end{array}
\]
\\
\[
\begin{array}{C{\mycolwd}|C{\mycolwd}C{\mycolwd}C{\mycolwd}C{\mycolwd}C{\mycolwd}C{\mycolwd}C{\mycolwd}}
 & \underline{b} & \overline{b}& c & c' & c'' & d & d_c' \\
\hline
B'[B] & & & C[B] & C'[B,C] & & D[B] & D_c'[B,D] \\
B_s[C] & & & C[C] & C'[C,C] & C'' & D[C] & D_c'[C,D]\\
B''[C] & & & C[C] & C'[C,C] & C'' & D[C] & D_c'[C,D]\\
\underline{B}[C] & & & C[C] & C'[C,C] & C'' & D[C] & D_c'[C,D]\\
B_s[D]\\
\hline
\end{array}
\]
\\
\[
\begin{array}{C{\mycolwd}|C{\mycolwd}C{\mycolwd}C{\mycolwd}C{\mycolwd}C{\mycolwd}C{\mycolwd}C{\mycolwd}}
 & d_c'' & \underline{d} & \overline{d} & \underline{\overline{d}} & d_\l & x & x' \\
\hline
B'[B]\\
B_s[C]\\
B''[C]\\
\underline{B}[C]\\
B_s[D] & D_c''\\
\hline
\end{array}
\]
\\
\[
\begin{array}{C{\mycolwd}|C{\mycolwd}C{\mycolwd}C{\mycolwd}C{\mycolwd}C{\mycolwd}C{\mycolwd}C{\mycolwd}}
 & x'' & \underline{x} & y & y' & y'' & \overline{y} & z \\
\hline
B'[B]\\
B_s[C] & & & & & & & Z[C]\\
B''[C] & & & & & & & Z[C]\\
\underline{B}[C] & & & & & & & Z[C]\\
B_s[D]\\
\hline
\end{array}
\]
\newpage
\[
\begin{array}{C{\mycolwd}|C{\mycolwd}C{\mycolwd}C{\mycolwd}C{\mycolwd}C{\mycolwd}C{\mycolwd}C{\mycolwd}}
 & a & a' & a'' & b & b' & b'' & b_s \\
\hline
B''[D]\\
\underline{B}[D]\\
B'[A,B] & & & A''[B]\\
C & A & A'[A] & & B & B'[B]\\
C'' & A & A'[A] & & B & B'[B]\\
\hline
\end{array}
\]
\\
\[
\begin{array}{C{\mycolwd}|C{\mycolwd}C{\mycolwd}C{\mycolwd}C{\mycolwd}C{\mycolwd}C{\mycolwd}C{\mycolwd}}
 & \underline{b} & \overline{b}& c & c' & c'' & d & d_c' \\
\hline
B''[D]\\
\underline{B}[D]\\
B'[A,B]\\
C & & \overline{B} & & C'[C] & & CD & D_c'[D]\\
C'' & & \overline{B} & C & C'[C] & & D & D_c'[D]\\
\hline
\end{array}
\]
\\
\[
\begin{array}{C{\mycolwd}|C{\mycolwd}C{\mycolwd}C{\mycolwd}C{\mycolwd}C{\mycolwd}C{\mycolwd}C{\mycolwd}}
 & d_c'' & \underline{d} & \overline{d} & \underline{\overline{d}} & d_\l & x & x' \\
\hline
B''[D] & D_c''\\
\underline{B}[D] & D_c''\\
B'[A,B]\\
C & & & \overline{D} & & & & X'[X]\\
C'' & & & \overline{D} & & & & X'[X]\\
\hline
\end{array}
\]
\\
\[
\begin{array}{C{\mycolwd}|C{\mycolwd}C{\mycolwd}C{\mycolwd}C{\mycolwd}C{\mycolwd}C{\mycolwd}C{\mycolwd}}
 & x'' & \underline{x} & y & y' & y'' & \overline{y} & z \\
\hline
B''[D]\\
\underline{B}[D]\\
B'[A,B]\\
C & & & & & & & Z\\
C'' & & & & & & & Z\\
\hline
\end{array}
\]
\newpage
\[
\begin{array}{C{\mycolwd}|C{\mycolwd}C{\mycolwd}C{\mycolwd}C{\mycolwd}C{\mycolwd}C{\mycolwd}C{\mycolwd}}
 & a & a' & a'' & b & b' & b'' & b_s \\
\hline
C\!\ra\!\{\underline{B},\underline{D},\underline{\overline{D}},\underline{X}\}\\
CD & & A'[A] & & B & B'[B]\\
C[B] & & & & & & B''\\
C[C]\\
C'[C] & & & & & & & B_s[C]\\
\hline
\end{array}
\]
\\
\[
\begin{array}{C{\mycolwd}|C{\mycolwd}C{\mycolwd}C{\mycolwd}C{\mycolwd}C{\mycolwd}C{\mycolwd}C{\mycolwd}}
 & \underline{b} & \overline{b} & c & c' & c'' & d & d_c' \\
\hline
C\!\ra\!\{\underline{B},\underline{D},\underline{\overline{D}},\underline{X}\} & \underline{B} & & & C'[C]\ra\underline{B}[C] & & & D_c'[D]\ra\underline{B}[D]\\
CD & & \overline{B} & C & C'[C] & & D & D_c'[D]\\
C[B] & & & & C'[B,C] & & & D_c'[C,D]\\
C[C] & & & & C'[C,C] & C'' & & D_c'[C,D]\\
C'[C]\\
\hline
\end{array}
\]
\\
\[
\begin{array}{C{\mycolwd}|C{\mycolwd}C{\mycolwd}C{\mycolwd}C{\mycolwd}C{\mycolwd}C{\mycolwd}C{\mycolwd}}
 & d_c'' & \underline{d} & \overline{d} & \underline{\overline{d}} & d_\l & x & x' \\
\hline
C\!\ra\!\{\underline{B},\underline{D},\underline{\overline{D}},\underline{X}\} & & \underline{D} & & \underline{\overline{D}}\\
CD & & & \overline{D} & & & & X'[X]\\
C[B]\\
C[C]\\
C'[C]\\
\hline
\end{array}
\]
\\
\[
\begin{array}{C{\mycolwd}|C{\mycolwd}C{\mycolwd}C{\mycolwd}C{\mycolwd}C{\mycolwd}C{\mycolwd}C{\mycolwd}}
 & x'' & \underline{x} & y & y' & y'' & \overline{y} & z \\
\hline
C\!\ra\!\{\underline{B},\underline{D},\underline{\overline{D}},\underline{X}\} & & & & & & & Z\ra\underline{X}\\
CD  & & & & & & & Z\\
C[B]\\
C[C] & & & & & & & Z[C]\\
C'[C]\\
\hline
\end{array}
\]
\newpage
\[
\begin{array}{C{\mycolwd}|C{\mycolwd}C{\mycolwd}C{\mycolwd}C{\mycolwd}C{\mycolwd}C{\mycolwd}C{\mycolwd}}
 & a & a' & a'' & b & b' & b'' & b_s \\
\hline
C'[C]\ra\underline{B}[C]\\
C''[C]\\
C''[D]\\
C[X]\\
C[Y]\\
\hline
\end{array}
\]
\\
\[
\begin{array}{C{\mycolwd}|C{\mycolwd}C{\mycolwd}C{\mycolwd}C{\mycolwd}C{\mycolwd}C{\mycolwd}C{\mycolwd}}
 & \underline{b} & \overline{b}& c & c' & c'' & d & d_c' \\
\hline
C'[C]\ra\underline{B}[C] & \underline{B}[C]\\
C''[C] & & & C[C] & C'[C,C] & C'' & D[C] & D_c'[C,D]\\
C''[D]\\
C[X]\\
C[Y] & & & & C'[Y,C] & & & D_c'[Y,D]\\
\hline
\end{array}
\]
\\
\[
\begin{array}{C{\mycolwd}|C{\mycolwd}C{\mycolwd}C{\mycolwd}C{\mycolwd}C{\mycolwd}C{\mycolwd}C{\mycolwd}}
 & d_c'' & \underline{d} & \overline{d} & \underline{\overline{d}} & d_\l & x & x' \\
\hline
C'[C]\ra\underline{B}[C]\\
C''[C]\\
C''[D] & D_c''\\
C[X]\\
C[Y]\\
\hline
\end{array}
\]
\\
\[
\begin{array}{C{\mycolwd}|C{\mycolwd}C{\mycolwd}C{\mycolwd}C{\mycolwd}C{\mycolwd}C{\mycolwd}C{\mycolwd}}
 & x'' & \underline{x} & y & y' & y'' & \overline{y} & z \\
\hline
C'[C]\ra\underline{B}[C]\\
C''[C] & & & & & & & Z[C]\\
C''[D]\\
C[X] & & & & & & & Z[X]\\
C[Y] & & & & & Y'' & & Z[Y]\\
\hline
\end{array}
\]
\newpage
\[
\begin{array}{C{\mycolwd}|C{\mycolwd}C{\mycolwd}C{\mycolwd}C{\mycolwd}C{\mycolwd}C{\mycolwd}C{\mycolwd}}
 & a & a' & a'' & b & b' & b'' & b_s \\
\hline
C'[B,C] & & & & & & B''[C]\\
C'[C,C]\\
C'[Y,C]\\
D & & A'[A] & & B & B'[B]\\
D_c''\\
\hline
\end{array}
\]
\\
\[
\begin{array}{C{\mycolwd}|C{\mycolwd}C{\mycolwd}C{\mycolwd}C{\mycolwd}C{\mycolwd}C{\mycolwd}C{\mycolwd}}
 & \underline{b} & \overline{b}& c & c' & c'' & d & d_c' \\
\hline
C'[B,C]\\
C'[C,C] & & & & & C''[C]\\
C'[Y,C]\\
D & & \overline{B} & C & C'[C] & & D & D_c'[D]\\
D_c''\\
\hline
\end{array}
\]
\\
\[
\begin{array}{C{\mycolwd}|C{\mycolwd}C{\mycolwd}C{\mycolwd}C{\mycolwd}C{\mycolwd}C{\mycolwd}C{\mycolwd}}
 & d_c'' & \underline{d} & \overline{d} & \underline{\overline{d}} & d_\l & x & x' \\
\hline
C'[B,C]\\
C'[C,C]\\
C'[Y,C]\\
D & & & \overline{D} & & D_\l & & X'[X]\\
D_c'' & & & & & D_\l\\
\hline
\end{array}
\]
\\
\[
\begin{array}{C{\mycolwd}|C{\mycolwd}C{\mycolwd}C{\mycolwd}C{\mycolwd}C{\mycolwd}C{\mycolwd}C{\mycolwd}}
 & x'' & \underline{x} & y & y' & y'' & \overline{y} & z \\
\hline
C'[B,C]\\
C'[C,C]\\
C'[Y,C] & & & & & Y''[C]\\
D & & & & & & & Z\\
D_c''\\
\hline
\end{array}
\]
\newpage
\[
\begin{array}{C{\mycolwd}|C{\mycolwd}C{\mycolwd}C{\mycolwd}C{\mycolwd}C{\mycolwd}C{\mycolwd}C{\mycolwd}}
 & a & a' & a'' & b & b' & b'' & b_s \\
\hline
\underline{D}\\
\overline{D}\\
\underline{\overline{D}}\\
D_\l\\
DB''\\
\hline
\end{array}
\]
\\
\[
\begin{array}{C{\mycolwd}|C{\mycolwd}C{\mycolwd}C{\mycolwd}C{\mycolwd}C{\mycolwd}C{\mycolwd}C{\mycolwd}}
 & \underline{b} & \overline{b}& c & c' & c'' & d & d_c' \\
\hline
\underline{D} & & & & & & \underline{D}D\\
\overline{D} & & & & & & \overline{D}D\\
\underline{\overline{D}} & & & & & & \underline{D}D\\
D_\l\\
DB''\\
\hline
\end{array}
\]
\\
\[
\begin{array}{C{\mycolwd}|C{\mycolwd}C{\mycolwd}C{\mycolwd}C{\mycolwd}C{\mycolwd}C{\mycolwd}C{\mycolwd}}
 & d_c'' & \underline{d} & \overline{d} & \underline{\overline{d}} & d_\l & x & x' \\
\hline
\underline{D}\\
\overline{D}\\
\underline{\overline{D}}\\
D_\l\\
DB'' & & & & & D_\l\\
\hline
\end{array}
\]
\\
\[
\begin{array}{C{\mycolwd}|C{\mycolwd}C{\mycolwd}C{\mycolwd}C{\mycolwd}C{\mycolwd}C{\mycolwd}C{\mycolwd}}
 & x'' & \underline{x} & y & y' & y'' & \overline{y} & z \\
\hline
\underline{D}\\
\overline{D}\\
\underline{\overline{D}}\\
D_\l\\
DB''\\
\hline
\end{array}
\]
\newpage
\[
\begin{array}{C{\mycolwd}|C{\mycolwd}C{\mycolwd}C{\mycolwd}C{\mycolwd}C{\mycolwd}C{\mycolwd}C{\mycolwd}}
 & a & a' & a'' & b & b' & b'' & b_s \\
\hline
D\underline{B}\\
D\ra\underline{B}\\
DC''\\
\underline{D}D & & A'[A] & & B & B'[B]\\
\overline{D}D\\
\hline
\end{array}
\]
\\
\[
\begin{array}{C{\mycolwd}|C{\mycolwd}C{\mycolwd}C{\mycolwd}C{\mycolwd}C{\mycolwd}C{\mycolwd}C{\mycolwd}}
 & \underline{b} & \overline{b}& c & c' & c'' & d & d_c' \\
\hline
D\underline{B}\\
D\ra\underline{B} & D\underline{B}\\
DC''\\
\underline{D}D & & \overline{B} & C & C'[C] & & D & D_c'[D]\\
\overline{D}D & \underline{B} & & C\!\ra\!\{\underline{B},\underline{D},\underline{\overline{D}},\underline{X}\} & C'[C]\ra\underline{B}[C] & & D\ra\underline{B} & D_c'[D]\ra\underline{B}[D]\\
\hline
\end{array}
\]
\\
\[
\begin{array}{C{\mycolwd}|C{\mycolwd}C{\mycolwd}C{\mycolwd}C{\mycolwd}C{\mycolwd}C{\mycolwd}C{\mycolwd}}
 & d_c'' & \underline{d} & \overline{d} & \underline{\overline{d}} & d_\l & x & x' \\
\hline
D\underline{B} & & & & &D_\l\\
D\ra\underline{B}\\
DC'' & & & & & D_\l\\
\underline{D}D & & & \overline{D} & & & & X'[X]\\
\overline{D}D & & \underline{D} & & \underline{\overline{D}}\\
\hline
\end{array}
\]
\\
\[
\begin{array}{C{\mycolwd}|C{\mycolwd}C{\mycolwd}C{\mycolwd}C{\mycolwd}C{\mycolwd}C{\mycolwd}C{\mycolwd}}
 & x'' & \underline{x} & y & y' & y'' & \overline{y} & z \\
\hline
D\underline{B}\\
D\ra\underline{B}\\
DC''\\
\underline{D}D & & & & & & & Z\\
\overline{D}D & & & & & & & Z\ra\underline{X}\\
\hline
\end{array}
\]
\newpage
\[
\begin{array}{C{\mycolwd}|C{\mycolwd}C{\mycolwd}C{\mycolwd}C{\mycolwd}C{\mycolwd}C{\mycolwd}C{\mycolwd}}
 & a & a' & a'' & b & b' & b'' & b_s \\
\hline
DY''\\
D[B] & & & & & & DB''\\
D[C]\\
D_c'[D] & & & & & & & B_s[D]\\
D_c'[D]\ra\underline{B}[D]\\
\hline
\end{array}
\]
\\
\[
\begin{array}{C{\mycolwd}|C{\mycolwd}C{\mycolwd}C{\mycolwd}C{\mycolwd}C{\mycolwd}C{\mycolwd}C{\mycolwd}}
 & \underline{b} & \overline{b}& c & c' & c'' & d & d_c' \\
\hline
DY''\\
D[B]\\
D[C] & & & & & DC''\\
D_c'[D]\\
D_c'[D]\ra\underline{B}[D] & \underline{B}[D]\\
\hline
\end{array}
\]
\\
\[
\begin{array}{C{\mycolwd}|C{\mycolwd}C{\mycolwd}C{\mycolwd}C{\mycolwd}C{\mycolwd}C{\mycolwd}C{\mycolwd}}
 & d_c'' & \underline{d} & \overline{d} & \underline{\overline{d}} & d_\l & x & x' \\
\hline
DY'' & & & & & D_\l\\
D[B]\\
D[C]\\
D_c'[D]\\
D_c'[D]\ra\underline{B}[D]\\
\hline
\end{array}
\]
\\
\[
\begin{array}{C{\mycolwd}|C{\mycolwd}C{\mycolwd}C{\mycolwd}C{\mycolwd}C{\mycolwd}C{\mycolwd}C{\mycolwd}}
 & x'' & \underline{x} & y & y' & y'' & \overline{y} & z \\
\hline
DY''\\
D[B]\\
D[C]\\
D_c'[D]\\
D_c'[D]\ra\underline{B}[D]\\
\hline
\end{array}
\]
\newpage
\[
\begin{array}{C{\mycolwd}|C{\mycolwd}C{\mycolwd}C{\mycolwd}C{\mycolwd}C{\mycolwd}C{\mycolwd}C{\mycolwd}}
 & a & a' & a'' & b & b' & b'' & b_s \\
\hline
D[Y]\\
D_c'[B,D] & & & & & & B''[D]\\
D_c'[C,D]\\
D_c'[Y,D]\\
X\\
\hline
\end{array}
\]
\\
\[
\begin{array}{C{\mycolwd}|C{\mycolwd}C{\mycolwd}C{\mycolwd}C{\mycolwd}C{\mycolwd}C{\mycolwd}C{\mycolwd}}
 & \underline{b} & \overline{b}& c & c' & c'' & d & d_c' \\
\hline
D[Y]\\
D_c'[B,D]\\
D_c'[C,D] & & & & & C''[D]\\
D_c'[Y,D]\\
X & & & & & & XD\\
\hline
\end{array}
\]
\\
\[
\begin{array}{C{\mycolwd}|C{\mycolwd}C{\mycolwd}C{\mycolwd}C{\mycolwd}C{\mycolwd}C{\mycolwd}C{\mycolwd}}
 & d_c'' & \underline{d} & \overline{d} & \underline{\overline{d}} & d_\l & x & x' \\
\hline
D[Y]\\
D_c'[B,D]\\
D_c'[C,D]\\
D_c'[Y,D]\\
X\\
\hline
\end{array}
\]
\\
\[
\begin{array}{C{\mycolwd}|C{\mycolwd}C{\mycolwd}C{\mycolwd}C{\mycolwd}C{\mycolwd}C{\mycolwd}C{\mycolwd}}
 & x'' & \underline{x} & y & y' & y'' & \overline{y} & z \\
\hline
D[Y] & & & & & DY''\\
D_c'[B,D]\\
D_c'[C,D]\\
D_c'[Y,D] & & & & & Y''[D]\\
X\\
\hline
\end{array}
\]
\newpage
\[
\begin{array}{C{\mycolwd}|C{\mycolwd}C{\mycolwd}C{\mycolwd}C{\mycolwd}C{\mycolwd}C{\mycolwd}C{\mycolwd}}
 & a & a' & a'' & b & b' & b'' & b_s \\
\hline
X''\\
\underline{X}\\
XD\\
X'[X] & A[X]\\
X'[A,X] & & & A''[X]\\
\hline
\end{array}
\]
\\
\[
\begin{array}{C{\mycolwd}|C{\mycolwd}C{\mycolwd}C{\mycolwd}C{\mycolwd}C{\mycolwd}C{\mycolwd}C{\mycolwd}}
 & \underline{b} & \overline{b}& c & c' & c'' & d & d_c' \\
\hline
X'' & & & & & & XD\\
\underline{X} & & & & & & XD\\
XD\\
X'[X] & & & C[X]\\
X'[A,X]\\
\hline
\end{array}
\]
\\
\[
\begin{array}{C{\mycolwd}|C{\mycolwd}C{\mycolwd}C{\mycolwd}C{\mycolwd}C{\mycolwd}C{\mycolwd}C{\mycolwd}}
 & d_c'' & \underline{d} & \overline{d} & \underline{\overline{d}} & d_\l & x & x' \\
\hline
X''\\
\underline{X}\\
XD\\
X'[X]\\
X'[A,X]\\
\hline
\end{array}
\]
\\
\[
\begin{array}{C{\mycolwd}|C{\mycolwd}C{\mycolwd}C{\mycolwd}C{\mycolwd}C{\mycolwd}C{\mycolwd}C{\mycolwd}}
 & x'' & \underline{x} & y & y' & y'' & \overline{y} & z \\
\hline
X''\\
\underline{X}\\
XD & & & Y & Y'[Y] & & \overline{Y}\\
X'[X] & & & & & & & Z[X]\\
X'[A,X]\\
\hline
\end{array}
\]
\newpage
\[
\begin{array}{C{\mycolwd}|C{\mycolwd}C{\mycolwd}C{\mycolwd}C{\mycolwd}C{\mycolwd}C{\mycolwd}C{\mycolwd}}
 & a & a' & a'' & b & b' & b'' & b_s \\
\hline
Y & A & A'[A] & & B & B'[B]\\
Y'' & A & A'[A] & & B & B'[B]\\
\overline{Y} & \overline{B}A\\
Y''[C]\\
Y''[D]\\
\hline
\end{array}
\]
\\
\[
\begin{array}{C{\mycolwd}|C{\mycolwd}C{\mycolwd}C{\mycolwd}C{\mycolwd}C{\mycolwd}C{\mycolwd}C{\mycolwd}}
 & \underline{b} & \overline{b}& c & c' & c'' & d & d_c' \\
\hline
Y & & \overline{B} & C & C'[C] & & D & D_c'[D]\\
Y'' & & \overline{B} & C & C'[C] & & D & D_c'[D]\\
\overline{Y} & & & \overline{B}{C}\\
Y''[C] & & & C[C] & C'[C,C] & C'' & D[C] & D_c'[C,D]\\
Y''[D]\\
\hline
\end{array}
\]
\\
\[
\begin{array}{C{\mycolwd}|C{\mycolwd}C{\mycolwd}C{\mycolwd}C{\mycolwd}C{\mycolwd}C{\mycolwd}C{\mycolwd}}
 & d_c'' & \underline{d} & \overline{d} & \underline{\overline{d}} & d_\l & x & x' \\
\hline
Y & & & \overline{D} & & & & X'[X]\\
Y'' & & & \overline{D} & & & & X'[X]\\
\overline{Y} & & \underline{D} & & \underline{\overline{D}}\\
Y''[C]\\
Y''[D] & D_c''\\
\hline
\end{array}
\]
\\
\[
\begin{array}{C{\mycolwd}|C{\mycolwd}C{\mycolwd}C{\mycolwd}C{\mycolwd}C{\mycolwd}C{\mycolwd}C{\mycolwd}}
 & x'' & \underline{x} & y & y' & y'' & \overline{y} & z \\
\hline
Y & & & & & & & Z\\
Y'' & & & & & & & Z\\
\overline{Y}\\
Y''[C] & & & & & & & Z[C]\\
Y''[D]\\
\hline
\end{array}
\]
\newpage\text{ }\\
\\
\[
\begin{array}{C{\mycolwd}|C{\mycolwd}C{\mycolwd}C{\mycolwd}C{\mycolwd}C{\mycolwd}C{\mycolwd}C{\mycolwd}}
 & a & a' & a'' & b & b' & b'' & b_s \\
\hline
Y'[Y] & A[Y]\\
Z\\
ZC''\\
Z\ra\underline{X}\\
\hline
\end{array}
\]
\\
\[
\begin{array}{C{\mycolwd}|C{\mycolwd}C{\mycolwd}C{\mycolwd}C{\mycolwd}C{\mycolwd}C{\mycolwd}C{\mycolwd}}
 & \underline{b} & \overline{b}& c & c' & c'' & d & d_c' \\
\hline
Y'[Y] & & & C[Y] & C'[Y,C] & & D[Y] & D_c'[Y,D]\\
Z\\
ZC'' & & & & & & XD\\
Z\ra\underline{X}\\
\hline
\end{array}
\]
\\
\[
\begin{array}{C{\mycolwd}|C{\mycolwd}C{\mycolwd}C{\mycolwd}C{\mycolwd}C{\mycolwd}C{\mycolwd}C{\mycolwd}}
 & d_c'' & \underline{d} & \overline{d} & \underline{\overline{d}} & d_\l & x & x' \\
\hline
Y'[Y]\\
Z & & & & & & X\\
ZC''\\
Z\ra\underline{X}\\
\hline
\end{array}
\]
\\
\[
\begin{array}{C{\mycolwd}|C{\mycolwd}C{\mycolwd}C{\mycolwd}C{\mycolwd}C{\mycolwd}C{\mycolwd}C{\mycolwd}}
 & x'' & \underline{x} & y & y' & y'' & \overline{y} & z \\
\hline
Y'[Y] & & & & & Y'' & & Z[Y]\\
Z\\
ZC''\\
Z\ra\underline{X} & & \underline{X}\\
\hline
\end{array}
\]
\newpage\text{ }\\
\\
\[
\begin{array}{C{\mycolwd}|C{\mycolwd}C{\mycolwd}C{\mycolwd}C{\mycolwd}C{\mycolwd}C{\mycolwd}C{\mycolwd}}
 & a & a' & a'' & b & b' & b'' & b_s \\
\hline
ZY''\\
Z[C]\\
Z[X]\\
Z[Y]\\
\hline
\end{array}
\]
\\
\[
\begin{array}{C{\mycolwd}|C{\mycolwd}C{\mycolwd}C{\mycolwd}C{\mycolwd}C{\mycolwd}C{\mycolwd}C{\mycolwd}}
 & \underline{b} & \overline{b}& c & c' & c'' & d & d_c' \\
\hline
ZY'' & & & & & & XD\\
Z[C] & & & & & ZC''\\
Z[X]\\
Z[Y]\\
\hline
\end{array}
\]
\\
\[
\begin{array}{C{\mycolwd}|C{\mycolwd}C{\mycolwd}C{\mycolwd}C{\mycolwd}C{\mycolwd}C{\mycolwd}C{\mycolwd}}
 & d_c'' & \underline{d} & \overline{d} & \underline{\overline{d}} & d_\l & x & x' \\
\hline
ZY''\\
Z[C]\\
Z[X]\\
Z[Y]\\
\hline
\end{array}
\]
\\
\[
\begin{array}{C{\mycolwd}|C{\mycolwd}C{\mycolwd}C{\mycolwd}C{\mycolwd}C{\mycolwd}C{\mycolwd}C{\mycolwd}}
 & x'' & \underline{x} & y & y' & y'' & \overline{y} & z \\
\hline
ZY''\\
Z[C]\\
Z[X] & X''\\
Z[Y] & & & & & ZY''\\
\hline
\end{array}
\]
\end{footnotesize}
\end{landscape}
\begin{landscape}
\begin{table}
\begin{tiny}
\[\kbordermatrix{ & A & A" & A[A] & A'[A] & A"[A] & A[B] & A"[B] & A[X] & A"[X] & A[Y] & A'[A,A] & B & B" & \underline{B} & \overline{B} & \overline{B}A & \overline{B}A" & \overline{B}C & B_s[A] & \overline{B}[A] & B'[B]\\
A                        & 0 & 0 & 0 & x & 0 & 0 & 0 & 0 & 0 & 0 & 0 & \bar{F} & 0 & 0 & x & 0 & 0 & 0 & 0 & 0 & x \\
A"                       & \bar{G} & 0 & 0 & x & 0 & 0 & 0 & 0 & 0 & 0 & 0 & \bar{F} & 0 & 0 & x & 0 & 0 & 0 & 0 & 0 & x \\
A[A]                     & 0 & 0 & 0 & 0 & 0 & 0 & 0 & 0 & 0 & 0 & x & 0 & 0 & 0 & 0 & 0 & 0 & 0 & x & x & 0 \\
A'[A]                    & 0 & 0 & \bar{G} & 0 & 0 & 0 & 0 & 0 & 0 & 0 & x & 0 & 0 & 0 & 0 & 0 & 0 & 0 & x & x & 0 \\
A"[A]                    & 0 & 0 & \bar{G} & 0 & 0 & 0 & 0 & 0 & 0 & 0 & x & 0 & 0 & 0 & 0 & 0 & 0 & 0 & x & x & 0 \\
A[B]                     & 0 & 0 & 0 & 0 & 0 & 0 & 0 & 0 & 0 & 0 & 0 & 0 & x & 0 & 0 & 0 & 0 & 0 & 0 & 0 & 0 \\
A"[B]                    & 0 & 0 & 0 & 0 & 0 & \bar{G} & 0 & 0 & 0 & 0 & 0 & 0 & x & 0 & 0 & 0 & 0 & 0 & 0 & 0 & 0 \\
A[X]                     & 0 & 0 & 0 & 0 & 0 & 0 & 0 & 0 & 0 & 0 & 0 & 0 & 0 & 0 & 0 & 0 & 0 & 0 & 0 & 0 & 0 \\
A"[X]                    & 0 & 0 & 0 & 0 & 0 & 0 & 0 & \bar{G} & 0 & 0 & 0 & 0 & 0 & 0 & 0 & 0 & 0 & 0 & 0 & 0 & 0 \\
A[Y]                     & 0 & 0 & 0 & 0 & 0 & 0 & 0 & 0 & 0 & 0 & 0 & 0 & 0 & 0 & 0 & 0 & 0 & 0 & 0 & 0 & 0 \\
A'[A,A]                  & 0 & 0 & 0 & 0 & x & 0 & 0 & 0 & 0 & 0 & 0 & 0 & 0 & 0 & 0 & 0 & 0 & 0 & 0 & 0 & 0 \\
B                        & \bar{G} & 0 & 0 & x & 0 & 0 & 0 & 0 & 0 & 0 & 0 & 0 & 0 & 0 & x & 0 & 0 & 0 & 0 & 0 & x \\
B"                       & \bar{G} & 0 & 0 & x & 0 & 0 & 0 & 0 & 0 & 0 & 0 & \bar{F} & 0 & 0 & x & 0 & 0 & 0 & 0 & 0 & x \\
\underline{B}            & \bar{G} & 0 & 0 & x & 0 & 0 & 0 & 0 & 0 & 0 & 0 & \bar{F} & 0 & 0 & x & 0 & 0 & 0 & 0 & 0 & x \\
\overline{B}             & 0 & 0 & 0 & 0 & 0 & 0 & 0 & 0 & 0 & 0 & 0 & 0 & 0 & 0 & 0 & \bar{G} & 0 & \bar{G} & 0 & 0 & 0 \\
\overline{B}A            & 0 & 0 & 0 & 0 & 0 & 0 & 0 & 0 & 0 & 0 & 0 & 0 & 0 & 0 & 0 & 0 & 0 & \bar{G} & 0 & 0 & 0 \\
\overline{B}A"           & 0 & 0 & 0 & 0 & 0 & 0 & 0 & 0 & 0 & 0 & 0 & 0 & 0 & 0 & 0 & \bar{G} & 0 & \bar{G} & 0 & 0 & 0 \\
\overline{B}C            & 0 & 0 & 0 & 0 & 0 & 0 & 0 & 0 & 0 & 0 & 0 & 0 & 0 & 0 & 0 & 0 & 0 & 0 & 0 & 0 & 0 \\
B_s[A]                   & 0 & \bar{G} & 0 & 0 & 0 & 0 & 0 & 0 & 0 & 0 & 0 & 0 & 0 & 0 & 0 & 0 & 0 & 0 & 0 & 0 & 0 \\
\overline{B}[A]          & 0 & 0 & 0 & 0 & 0 & 0 & 0 & 0 & 0 & 0 & 0 & 0 & 0 & 0 & 0 & 0 & \bar{G} & 0 & 0 & 0 & 0 \\
B'[B]                    & 0 & 0 & 0 & 0 & 0 & \bar{G} & 0 & 0 & 0 & 0 & 0 & 0 & 0 & 0 & 0 & 0 & 0 & 0 & 0 & 0 & 0 \\
B_s[C]                   & 0 & 0 & 0 & 0 & 0 & 0 & 0 & 0 & 0 & 0 & 0 & 0 & 0 & 0 & 0 & 0 & 0 & 0 & 0 & 0 & 0 \\
B"[C]                    & 0 & 0 & 0 & 0 & 0 & 0 & 0 & 0 & 0 & 0 & 0 & 0 & 0 & 0 & 0 & 0 & 0 & 0 & 0 & 0 & 0 \\
\underline{B}[C]         & 0 & 0 & 0 & 0 & 0 & 0 & 0 & 0 & 0 & 0 & 0 & 0 & 0 & 0 & 0 & 0 & 0 & 0 & 0 & 0 & 0 \\
B_s[D]                   & 0 & 0 & 0 & 0 & 0 & 0 & 0 & 0 & 0 & 0 & 0 & 0 & 0 & 0 & 0 & 0 & 0 & 0 & 0 & 0 & 0 \\
B"[D]                    & 0 & 0 & 0 & 0 & 0 & 0 & 0 & 0 & 0 & 0 & 0 & 0 & 0 & 0 & 0 & 0 & 0 & 0 & 0 & 0 & 0 \\
\underline{B}[D]         & 0 & 0 & 0 & 0 & 0 & 0 & 0 & 0 & 0 & 0 & 0 & 0 & 0 & 0 & 0 & 0 & 0 & 0 & 0 & 0 & 0 \\
B'[A,B]                  & 0 & 0 & 0 & 0 & 0 & 0 & \bar{G} & 0 & 0 & 0 & 0 & 0 & 0 & 0 & 0 & 0 & 0 & 0 & 0 & 0 & 0 \\
C                        & \bar{G} & 0 & 0 & x & 0 & 0 & 0 & 0 & 0 & 0 & 0 & \bar{F} & 0 & 0 & x & 0 & 0 & 0 & 0 & 0 & x \\
C"                       & \bar{G} & 0 & 0 & x & 0 & 0 & 0 & 0 & 0 & 0 & 0 & \bar{F} & 0 & 0 & x & 0 & 0 & 0 & 0 & 0 & x \\
C\ra\{\underline{B},\underline{D},\underline{\overline{D}},\underline{X}\}        & 0 & 0 & 0 & 0 & 0 & 0 & 0 & 0 & 0 & 0 & 0 & 0 & 0 & x & 0 & 0 & 0 & 0 & 0 & 0 & 0 \\
CD                       & 0 & 0 & 0 & x & 0 & 0 & 0 & 0 & 0 & 0 & 0 & \bar{F} & 0 & 0 & x & 0 & 0 & 0 & 0 & 0 & x \\
C[B]                     & 0 & 0 & 0 & 0 & 0 & 0 & 0 & 0 & 0 & 0 & 0 & 0 & x & 0 & 0 & 0 & 0 & 0 & 0 & 0 & 0 \\
C[C]                     & 0 & 0 & 0 & 0 & 0 & 0 & 0 & 0 & 0 & 0 & 0 & 0 & 0 & 0 & 0 & 0 & 0 & 0 & 0 & 0 & 0 \\
C'[C]                    & 0 & 0 & 0 & 0 & 0 & 0 & 0 & 0 & 0 & 0 & 0 & 0 & 0 & 0 & 0 & 0 & 0 & 0 & 0 & 0 & 0 \\
C'[C]\ra\underline{B}[C] & 0 & 0 & 0 & 0 & 0 & 0 & 0 & 0 & 0 & 0 & 0 & 0 & 0 & 0 & 0 & 0 & 0 & 0 & 0 & 0 & 0 \\
C"[C]                    & 0 & 0 & 0 & 0 & 0 & 0 & 0 & 0 & 0 & 0 & 0 & 0 & 0 & 0 & 0 & 0 & 0 & 0 & 0 & 0 & 0 \\
C"[D]                    & 0 & 0 & 0 & 0 & 0 & 0 & 0 & 0 & 0 & 0 & 0 & 0 & 0 & 0 & 0 & 0 & 0 & 0 & 0 & 0 & 0 \\
C[X]                     & 0 & 0 & 0 & 0 & 0 & 0 & 0 & 0 & 0 & 0 & 0 & 0 & 0 & 0 & 0 & 0 & 0 & 0 & 0 & 0 & 0 \\
C[Y]                     & 0 & 0 & 0 & 0 & 0 & 0 & 0 & 0 & 0 & 0 & 0 & 0 & 0 & 0 & 0 & 0 & 0 & 0 & 0 & 0 & 0 \\
C'[B,C]                  & 0 & 0 & 0 & 0 & 0 & 0 & 0 & 0 & 0 & 0 & 0 & 0 & 0 & 0 & 0 & 0 & 0 & 0 & 0 & 0 & 0 \\
C'[C,C]                  & 0 & 0 & 0 & 0 & 0 & 0 & 0 & 0 & 0 & 0 & 0 & 0 & 0 & 0 & 0 & 0 & 0 & 0 & 0 & 0 & 0 \\
}\]\caption{Adjacency matrix associated with $M'_i$.}
\label{tab:A.2}
\end{tiny}
\end{table}
\end{landscape}
\newpage
\begin{landscape}
\begin{tiny}
\[\kbordermatrix{ & B_s[C] & B"[C] & \underline{B}[C] & B_s[D] & B"[D] & \underline{B}[D] & B'[A,B] & C & C" & C\ra\{\underline{B},\underline{D},\underline{\overline{D}},\underline{X}\}& CD & C[B] & C[C] & C'[C] & C'[C]\ra\underline{B}[C] & C"[C] & C"[D] & C[X] & C[Y]\\
& 0 & 0 & 0 & 0 & 0 & 0 & 0 & \bar{G} & 0 & 0 & 0 & 0 & 0 & \bar{G} & 0 & 0 & 0 & 0 & 0 \\
& 0 & 0 & 0 & 0 & 0 & 0 & 0 & \bar{G} & 0 & 0 & 0 & 0 & 0 & \bar{G} & 0 & 0 & 0 & 0 & 0 \\
& 0 & 0 & 0 & 0 & 0 & 0 & x & 0 & 0 & 0 & 0 & 0 & 0 & 0 & 0 & 0 & 0 & 0 & 0 \\
& 0 & 0 & 0 & 0 & 0 & 0 & x & 0 & 0 & 0 & 0 & 0 & 0 & 0 & 0 & 0 & 0 & 0 & 0 \\
& 0 & 0 & 0 & 0 & 0 & 0 & x & 0 & 0 & 0 & 0 & 0 & 0 & 0 & 0 & 0 & 0 & 0 & 0 \\
& 0 & 0 & 0 & 0 & 0 & 0 & 0 & 0 & 0 & 0 & 0 & \bar{G} & 0 & 0 & 0 & 0 & 0 & 0 & 0 \\
& 0 & 0 & 0 & 0 & 0 & 0 & 0 & 0 & 0 & 0 & 0 & \bar{G} & 0 & 0 & 0 & 0 & 0 & 0 & 0 \\
& 0 & 0 & 0 & 0 & 0 & 0 & 0 & 0 & 0 & 0 & 0 & 0 & 0 & 0 & 0 & 0 & 0 & \bar{G} & 0 \\
& 0 & 0 & 0 & 0 & 0 & 0 & 0 & 0 & 0 & 0 & 0 & 0 & 0 & 0 & 0 & 0 & 0 & \bar{G} & 0 \\
& 0 & 0 & 0 & 0 & 0 & 0 & 0 & 0 & 0 & 0 & 0 & 0 & 0 & 0 & 0 & 0 & 0 & 0 & \bar{G} \\
& 0 & 0 & 0 & 0 & 0 & 0 & 0 & 0 & 0 & 0 & 0 & 0 & 0 & 0 & 0 & 0 & 0 & 0 & 0 \\
& 0 & 0 & 0 & 0 & 0 & 0 & 0 & \bar{G} & 0 & 0 & 0 & 0 & 0 & \bar{G} & 0 & 0 & 0 & 0 & 0 \\
& 0 & 0 & 0 & 0 & 0 & 0 & 0 & \bar{G} & 0 & 0 & 0 & 0 & 0 & \bar{G} & 0 & 0 & 0 & 0 & 0 \\
& 0 & 0 & 0 & 0 & 0 & 0 & 0 & \bar{G} & 0 & 0 & 0 & 0 & 0 & \bar{G} & 0 & 0 & 0 & 0 & 0 \\
& 0 & 0 & 0 & 0 & 0 & 0 & 0 & 0 & 0 & 0 & 0 & 0 & 0 & 0 & 0 & 0 & 0 & 0 & 0 \\
& 0 & 0 & 0 & 0 & 0 & 0 & 0 & 0 & 0 & 0 & 0 & 0 & 0 & 0 & 0 & 0 & 0 & 0 & 0 \\
& 0 & 0 & 0 & 0 & 0 & 0 & 0 & 0 & 0 & 0 & 0 & 0 & 0 & 0 & 0 & 0 & 0 & 0 & 0 \\
& 0 & 0 & 0 & 0 & 0 & 0 & 0 & 0 & 0 & 0 & 0 & 0 & 0 & 0 & 0 & 0 & 0 & 0 & 0 \\
& 0 & 0 & 0 & 0 & 0 & 0 & 0 & 0 & 0 & 0 & 0 & 0 & 0 & 0 & 0 & 0 & 0 & 0 & 0 \\
& 0 & 0 & 0 & 0 & 0 & 0 & 0 & 0 & 0 & 0 & 0 & 0 & 0 & 0 & 0 & 0 & 0 & 0 & 0 \\
& 0 & 0 & 0 & 0 & 0 & 0 & 0 & 0 & 0 & 0 & 0 & \bar{G} & 0 & 0 & 0 & 0 & 0 & 0 & 0 \\
& 0 & 0 & 0 & 0 & 0 & 0 & 0 & 0 & x & 0 & 0 & 0 & \bar{G} & 0 & 0 & 0 & 0 & 0 & 0 \\
& 0 & 0 & 0 & 0 & 0 & 0 & 0 & 0 & x & 0 & 0 & 0 & \bar{G} & 0 & 0 & 0 & 0 & 0 & 0 \\
& 0 & 0 & 0 & 0 & 0 & 0 & 0 & 0 & x & 0 & 0 & 0 & \bar{G} & 0 & 0 & 0 & 0 & 0 & 0 \\
& 0 & 0 & 0 & 0 & 0 & 0 & 0 & 0 & 0 & 0 & 0 & 0 & 0 & 0 & 0 & 0 & 0 & 0 & 0 \\
& 0 & 0 & 0 & 0 & 0 & 0 & 0 & 0 & 0 & 0 & 0 & 0 & 0 & 0 & 0 & 0 & 0 & 0 & 0 \\
& 0 & 0 & 0 & 0 & 0 & 0 & 0 & 0 & 0 & 0 & 0 & 0 & 0 & 0 & 0 & 0 & 0 & 0 & 0 \\
& 0 & 0 & 0 & 0 & 0 & 0 & 0 & 0 & 0 & 0 & 0 & 0 & 0 & 0 & 0 & 0 & 0 & 0 & 0 \\
& 0 & 0 & 0 & 0 & 0 & 0 & 0 & 0 & 0 & 0 & \bar{G} & 0 & 0 & \bar{G} & 0 & 0 & 0 & 0 & 0 \\
& 0 & 0 & 0 & 0 & 0 & 0 & 0 & \bar{G} & 0 & 0 & 0 & 0 & 0 & \bar{G} & 0 & 0 & 0 & 0 & 0 \\
& 0 & 0 & 0 & 0 & 0 & 0 & 0 & 0 & 0 & 0 & 0 & 0 & 0 & 0 & \bar{G} & 0 & 0 & 0 & 0 \\
& 0 & 0 & 0 & 0 & 0 & 0 & 0 & \bar{G} & 0 & 0 & 0 & 0 & 0 & \bar{G} & 0 & 0 & 0 & 0 & 0 \\
& 0 & 0 & 0 & 0 & 0 & 0 & 0 & 0 & 0 & 0 & 0 & 0 & 0 & 0 & 0 & 0 & 0 & 0 & 0 \\
& 0 & 0 & 0 & 0 & 0 & 0 & 0 & 0 & x & 0 & 0 & 0 & 0 & 0 & 0 & 0 & 0 & 0 & 0 \\
& x & 0 & 0 & 0 & 0 & 0 & 0 & 0 & 0 & 0 & 0 & 0 & 0 & 0 & 0 & 0 & 0 & 0 & 0 \\
& 0 & 0 & x & 0 & 0 & 0 & 0 & 0 & 0 & 0 & 0 & 0 & 0 & 0 & 0 & 0 & 0 & 0 & 0 \\
& 0 & 0 & 0 & 0 & 0 & 0 & 0 & 0 & x & 0 & 0 & 0 & \bar{G} & 0 & 0 & 0 & 0 & 0 & 0 \\
& 0 & 0 & 0 & 0 & 0 & 0 & 0 & 0 & 0 & 0 & 0 & 0 & 0 & 0 & 0 & 0 & 0 & 0 & 0 \\
& 0 & 0 & 0 & 0 & 0 & 0 & 0 & 0 & 0 & 0 & 0 & 0 & 0 & 0 & 0 & 0 & 0 & 0 & 0 \\
& 0 & 0 & 0 & 0 & 0 & 0 & 0 & 0 & 0 & 0 & 0 & 0 & 0 & 0 & 0 & 0 & 0 & 0 & 0 \\
& 0 & x & 0 & 0 & 0 & 0 & 0 & 0 & 0 & 0 & 0 & 0 & 0 & 0 & 0 & 0 & 0 & 0 & 0 \\
& 0 & 0 & 0 & 0 & 0 & 0 & 0 & 0 & 0 & 0 & 0 & 0 & 0 & 0 & 0 & x & 0 & 0 & 0
}\]
\newpage
\text{ }\\
\[\kbordermatrix{ & C'[B,C] & C'[C,C] & C'[Y,C] & D & D" & \underline{D} & \overline{D} & \underline{\overline{D}} & D_\l & DB" & D\underline{B} & D\ra\underline{B} & DC" & \underline{D}D & \overline{D}D & DY" & D[B] & D[C] & D_c'[D] & D_c'[D]\ra\underline{B}[D] & D[Y]\\
& 0 & 0 & 0 & \bar{G} & 0 & 0 & x & 0 & 0 & 0 & 0 & 0 & 0 & 0 & 0 & 0 & 0 & 0 & \bar{G} & 0 & 0 \\
& 0 & 0 & 0 & \bar{G} & 0 & 0 & x & 0 & 0 & 0 & 0 & 0 & 0 & 0 & 0 & 0 & 0 & 0 & \bar{G} & 0 & 0 \\
& 0 & 0 & 0 & 0 & 0 & 0 & 0 & 0 & 0 & 0 & 0 & 0 & 0 & 0 & 0 & 0 & 0 & 0 & 0 & 0 & 0 \\
& 0 & 0 & 0 & 0 & 0 & 0 & 0 & 0 & 0 & 0 & 0 & 0 & 0 & 0 & 0 & 0 & 0 & 0 & 0 & 0 & 0 \\
& 0 & 0 & 0 & 0 & 0 & 0 & 0 & 0 & 0 & 0 & 0 & 0 & 0 & 0 & 0 & 0 & 0 & 0 & 0 & 0 & 0 \\
& \bar{G} & 0 & 0 & 0 & 0 & 0 & 0 & 0 & 0 & 0 & 0 & 0 & 0 & 0 & 0 & 0 & \bar{G} & 0 & 0 & 0 & 0 \\
& \bar{G} & 0 & 0 & 0 & 0 & 0 & 0 & 0 & 0 & 0 & 0 & 0 & 0 & 0 & 0 & 0 & \bar{G} & 0 & 0 & 0 & 0 \\
& 0 & 0 & 0 & 0 & 0 & 0 & 0 & 0 & 0 & 0 & 0 & 0 & 0 & 0 & 0 & 0 & 0 & 0 & 0 & 0 & 0 \\
& 0 & 0 & 0 & 0 & 0 & 0 & 0 & 0 & 0 & 0 & 0 & 0 & 0 & 0 & 0 & 0 & 0 & 0 & 0 & 0 & 0 \\
& 0 & 0 & \bar{G} & 0 & 0 & 0 & 0 & 0 & 0 & 0 & 0 & 0 & 0 & 0 & 0 & 0 & 0 & 0 & 0 & 0 & \bar{G} \\
& 0 & 0 & 0 & 0 & 0 & 0 & 0 & 0 & 0 & 0 & 0 & 0 & 0 & 0 & 0 & 0 & 0 & 0 & 0 & 0 & 0 \\
& 0 & 0 & 0 & \bar{G} & 0 & 0 & x & 0 & 0 & 0 & 0 & 0 & 0 & 0 & 0 & 0 & 0 & 0 & \bar{G} & 0 & 0 \\
& 0 & 0 & 0 & \bar{G} & 0 & 0 & x & 0 & 0 & 0 & 0 & 0 & 0 & 0 & 0 & 0 & 0 & 0 & \bar{G} & 0 & 0 \\
& 0 & 0 & 0 & \bar{G} & 0 & 0 & x & 0 & 0 & 0 & 0 & 0 & 0 & 0 & 0 & 0 & 0 & 0 & \bar{G} & 0 & 0 \\
& 0 & 0 & 0 & 0 & 0 & x & 0 & x & 0 & 0 & 0 & 0 & 0 & 0 & 0 & 0 & 0 & 0 & 0 & 0 & 0 \\
& 0 & 0 & 0 & 0 & 0 & x & 0 & x & 0 & 0 & 0 & 0 & 0 & 0 & 0 & 0 & 0 & 0 & 0 & 0 & 0 \\
& 0 & 0 & 0 & 0 & 0 & x & 0 & x & 0 & 0 & 0 & 0 & 0 & 0 & 0 & 0 & 0 & 0 & 0 & 0 & 0 \\
& 0 & 0 & 0 & 0 & 0 & x & 0 & x & 0 & 0 & 0 & 0 & 0 & 0 & 0 & 0 & 0 & 0 & 0 & 0 & 0 \\
& 0 & 0 & 0 & 0 & 0 & 0 & 0 & 0 & 0 & 0 & 0 & 0 & 0 & 0 & 0 & 0 & 0 & 0 & 0 & 0 & 0 \\
& 0 & 0 & 0 & 0 & 0 & 0 & 0 & 0 & 0 & 0 & 0 & 0 & 0 & 0 & 0 & 0 & 0 & 0 & 0 & 0 & 0 \\
& \bar{G} & 0 & 0 & 0 & 0 & 0 & 0 & 0 & 0 & 0 & 0 & 0 & 0 & 0 & 0 & 0 & \bar{G} & 0 & 0 & 0 & 0 \\
& 0 & \bar{G} & 0 & 0 & 0 & 0 & 0 & 0 & 0 & 0 & 0 & 0 & 0 & 0 & 0 & 0 & 0 & \bar{G} & 0 & 0 & 0 \\
& 0 & \bar{G} & 0 & 0 & 0 & 0 & 0 & 0 & 0 & 0 & 0 & 0 & 0 & 0 & 0 & 0 & 0 & \bar{G} & 0 & 0 & 0 \\
& 0 & \bar{G} & 0 & 0 & 0 & 0 & 0 & 0 & 0 & 0 & 0 & 0 & 0 & 0 & 0 & 0 & 0 & \bar{G} & 0 & 0 & 0 \\
& 0 & 0 & 0 & 0 & x & 0 & 0 & 0 & 0 & 0 & 0 & 0 & 0 & 0 & 0 & 0 & 0 & 0 & 0 & 0 & 0 \\
& 0 & 0 & 0 & 0 & x & 0 & 0 & 0 & 0 & 0 & 0 & 0 & 0 & 0 & 0 & 0 & 0 & 0 & 0 & 0 & 0 \\
& 0 & 0 & 0 & 0 & x & 0 & 0 & 0 & 0 & 0 & 0 & 0 & 0 & 0 & 0 & 0 & 0 & 0 & 0 & 0 & 0 \\
& 0 & 0 & 0 & 0 & 0 & 0 & 0 & 0 & 0 & 0 & 0 & 0 & 0 & 0 & 0 & 0 & 0 & 0 & 0 & 0 & 0 \\
& 0 & 0 & 0 & 0 & 0 & 0 & x & 0 & 0 & 0 & 0 & 0 & 0 & 0 & 0 & 0 & 0 & 0 & \bar{G} & 0 & 0 \\
& 0 & 0 & 0 & \bar{G} & 0 & 0 & x & 0 & 0 & 0 & 0 & 0 & 0 & 0 & 0 & 0 & 0 & 0 & \bar{G} & 0 & 0 \\
& 0 & 0 & 0 & 0 & 0 & x & 0 & x & 0 & 0 & 0 & 0 & 0 & 0 & 0 & 0 & 0 & 0 & 0 & \bar{G} & 0 \\
& 0 & 0 & 0 & \bar{G} & 0 & 0 & x & 0 & 0 & 0 & 0 & 0 & 0 & 0 & 0 & 0 & 0 & 0 & \bar{G} & 0 & 0 \\
& \bar{G} & 0 & 0 & 0 & 0 & 0 & 0 & 0 & 0 & 0 & 0 & 0 & 0 & 0 & 0 & 0 & 0 & 0 & 0 & 0 & 0 \\
& 0 & \bar{G} & 0 & 0 & 0 & 0 & 0 & 0 & 0 & 0 & 0 & 0 & 0 & 0 & 0 & 0 & 0 & 0 & 0 & 0 & 0 \\
& 0 & 0 & 0 & 0 & 0 & 0 & 0 & 0 & 0 & 0 & 0 & 0 & 0 & 0 & 0 & 0 & 0 & 0 & 0 & 0 & 0 \\
& 0 & 0 & 0 & 0 & 0 & 0 & 0 & 0 & 0 & 0 & 0 & 0 & 0 & 0 & 0 & 0 & 0 & 0 & 0 & 0 & 0 \\
& 0 & \bar{G} & 0 & 0 & 0 & 0 & 0 & 0 & 0 & 0 & 0 & 0 & 0 & 0 & 0 & 0 & 0 & \bar{G} & 0 & 0 & 0 \\
& 0 & 0 & 0 & 0 & x & 0 & 0 & 0 & 0 & 0 & 0 & 0 & 0 & 0 & 0 & 0 & 0 & 0 & 0 & 0 & 0 \\
& 0 & 0 & 0 & 0 & 0 & 0 & 0 & 0 & 0 & 0 & 0 & 0 & 0 & 0 & 0 & 0 & 0 & 0 & 0 & 0 & 0 \\
& 0 & 0 & \bar{G} & 0 & 0 & 0 & 0 & 0 & 0 & 0 & 0 & 0 & 0 & 0 & 0 & 0 & 0 & 0 & 0 & 0 & 0 \\
& 0 & 0 & 0 & 0 & 0 & 0 & 0 & 0 & 0 & 0 & 0 & 0 & 0 & 0 & 0 & 0 & 0 & 0 & 0 & 0 & 0 \\
& 0 & 0 & 0 & 0 & 0 & 0 & 0 & 0 & 0 & 0 & 0 & 0 & 0 & 0 & 0 & 0 & 0 & 0 & 0 & 0 & 0
}\]
\newpage
\[\kbordermatrix{ & D_c'[B,D] & D_c'[C,D] & D_c'[Y,D] & X & X" & \underline{X} & XD & X'[X] & X'[A,X] & Y & Y" & \overline{Y} & Y"[C] & Y"[D] & Y'[Y] & Z & ZC" & Z\ra\underline{X} & ZY" & Z[C] & Z[X] & Z[Y]\\
& 0 & 0 & 0 & 0 & 0 & 0 & 0 & x & 0 & 0 & 0 & 0 & 0 & 0 & 0 & \bar{G} & 0 & 0 & 0 & 0 & 0 & 0 \\
& 0 & 0 & 0 & 0 & 0 & 0 & 0 & x & 0 & 0 & 0 & 0 & 0 & 0 & 0 & \bar{G} & 0 & 0 & 0 & 0 & 0 & 0 \\
& 0 & 0 & 0 & 0 & 0 & 0 & 0 & 0 & x & 0 & 0 & 0 & 0 & 0 & 0 & 0 & 0 & 0 & 0 & 0 & 0 & 0 \\
& 0 & 0 & 0 & 0 & 0 & 0 & 0 & 0 & x & 0 & 0 & 0 & 0 & 0 & 0 & 0 & 0 & 0 & 0 & 0 & 0 & 0 \\
& 0 & 0 & 0 & 0 & 0 & 0 & 0 & 0 & x & 0 & 0 & 0 & 0 & 0 & 0 & 0 & 0 & 0 & 0 & 0 & 0 & 0 \\
& \bar{G} & 0 & 0 & 0 & 0 & 0 & 0 & 0 & 0 & 0 & 0 & 0 & 0 & 0 & 0 & 0 & 0 & 0 & 0 & 0 & 0 & 0 \\
& \bar{G} & 0 & 0 & 0 & 0 & 0 & 0 & 0 & 0 & 0 & 0 & 0 & 0 & 0 & 0 & 0 & 0 & 0 & 0 & 0 & 0 & 0 \\
& 0 & 0 & 0 & 0 & 0 & 0 & 0 & 0 & 0 & 0 & 0 & 0 & 0 & 0 & 0 & 0 & 0 & 0 & 0 & 0 & \bar{G} & 0 \\
& 0 & 0 & 0 & 0 & 0 & 0 & 0 & 0 & 0 & 0 & 0 & 0 & 0 & 0 & 0 & 0 & 0 & 0 & 0 & 0 & \bar{G} & 0 \\
& 0 & 0 & \bar{G} & 0 & 0 & 0 & 0 & 0 & 0 & 0 & x & 0 & 0 & 0 & 0 & 0 & 0 & 0 & 0 & 0 & 0 & \bar{G} \\
& 0 & 0 & 0 & 0 & 0 & 0 & 0 & 0 & 0 & 0 & 0 & 0 & 0 & 0 & 0 & 0 & 0 & 0 & 0 & 0 & 0 & 0 \\
& 0 & 0 & 0 & 0 & 0 & 0 & 0 & x & 0 & 0 & 0 & 0 & 0 & 0 & 0 & \bar{G} & 0 & 0 & 0 & 0 & 0 & 0 \\
& 0 & 0 & 0 & 0 & 0 & 0 & 0 & x & 0 & 0 & 0 & 0 & 0 & 0 & 0 & \bar{G} & 0 & 0 & 0 & 0 & 0 & 0 \\
& 0 & 0 & 0 & 0 & 0 & 0 & 0 & x & 0 & 0 & 0 & 0 & 0 & 0 & 0 & \bar{G} & 0 & 0 & 0 & 0 & 0 & 0 \\
& 0 & 0 & 0 & 0 & 0 & 0 & 0 & 0 & 0 & 0 & 0 & 0 & 0 & 0 & 0 & 0 & 0 & 0 & 0 & 0 & 0 & 0 \\
& 0 & 0 & 0 & 0 & 0 & 0 & 0 & 0 & 0 & 0 & 0 & 0 & 0 & 0 & 0 & 0 & 0 & 0 & 0 & 0 & 0 & 0 \\
& 0 & 0 & 0 & 0 & 0 & 0 & 0 & 0 & 0 & 0 & 0 & 0 & 0 & 0 & 0 & 0 & 0 & 0 & 0 & 0 & 0 & 0 \\
& 0 & 0 & 0 & 0 & 0 & 0 & 0 & 0 & 0 & 0 & 0 & 0 & 0 & 0 & 0 & 0 & 0 & 0 & 0 & 0 & 0 & 0 \\
& 0 & 0 & 0 & 0 & 0 & 0 & 0 & 0 & 0 & 0 & 0 & 0 & 0 & 0 & 0 & 0 & 0 & 0 & 0 & 0 & 0 & 0 \\
& 0 & 0 & 0 & 0 & 0 & 0 & 0 & 0 & 0 & 0 & 0 & 0 & 0 & 0 & 0 & 0 & 0 & 0 & 0 & 0 & 0 & 0 \\
& \bar{G} & 0 & 0 & 0 & 0 & 0 & 0 & 0 & 0 & 0 & 0 & 0 & 0 & 0 & 0 & 0 & 0 & 0 & 0 & 0 & 0 & 0 \\
& 0 & \bar{G} & 0 & 0 & 0 & 0 & 0 & 0 & 0 & 0 & 0 & 0 & 0 & 0 & 0 & 0 & 0 & 0 & 0 & \bar{G} & 0 & 0 \\
& 0 & \bar{G} & 0 & 0 & 0 & 0 & 0 & 0 & 0 & 0 & 0 & 0 & 0 & 0 & 0 & 0 & 0 & 0 & 0 & \bar{G} & 0 & 0 \\
& 0 & \bar{G} & 0 & 0 & 0 & 0 & 0 & 0 & 0 & 0 & 0 & 0 & 0 & 0 & 0 & 0 & 0 & 0 & 0 & \bar{G} & 0 & 0 \\
& 0 & 0 & 0 & 0 & 0 & 0 & 0 & 0 & 0 & 0 & 0 & 0 & 0 & 0 & 0 & 0 & 0 & 0 & 0 & 0 & 0 & 0 \\
& 0 & 0 & 0 & 0 & 0 & 0 & 0 & 0 & 0 & 0 & 0 & 0 & 0 & 0 & 0 & 0 & 0 & 0 & 0 & 0 & 0 & 0 \\
& 0 & 0 & 0 & 0 & 0 & 0 & 0 & 0 & 0 & 0 & 0 & 0 & 0 & 0 & 0 & 0 & 0 & 0 & 0 & 0 & 0 & 0 \\
& 0 & 0 & 0 & 0 & 0 & 0 & 0 & 0 & 0 & 0 & 0 & 0 & 0 & 0 & 0 & 0 & 0 & 0 & 0 & 0 & 0 & 0 \\
& 0 & 0 & 0 & 0 & 0 & 0 & 0 & x & 0 & 0 & 0 & 0 & 0 & 0 & 0 & \bar{G} & 0 & 0 & 0 & 0 & 0 & 0 \\
& 0 & 0 & 0 & 0 & 0 & 0 & 0 & x & 0 & 0 & 0 & 0 & 0 & 0 & 0 & \bar{G} & 0 & 0 & 0 & 0 & 0 & 0 \\
& 0 & 0 & 0 & 0 & 0 & 0 & 0 & 0 & 0 & 0 & 0 & 0 & 0 & 0 & 0 & 0 & 0 & \bar{G} & 0 & 0 & 0 & 0 \\
& 0 & 0 & 0 & 0 & 0 & 0 & 0 & x & 0 & 0 & 0 & 0 & 0 & 0 & 0 & \bar{G} & 0 & 0 & 0 & 0 & 0 & 0 \\
& \bar{G} & 0 & 0 & 0 & 0 & 0 & 0 & 0 & 0 & 0 & 0 & 0 & 0 & 0 & 0 & 0 & 0 & 0 & 0 & 0 & 0 & 0 \\
& 0 & \bar{G} & 0 & 0 & 0 & 0 & 0 & 0 & 0 & 0 & 0 & 0 & 0 & 0 & 0 & 0 & 0 & 0 & 0 & \bar{G} & 0 & 0 \\
& 0 & 0 & 0 & 0 & 0 & 0 & 0 & 0 & 0 & 0 & 0 & 0 & 0 & 0 & 0 & 0 & 0 & 0 & 0 & 0 & 0 & 0 \\
& 0 & 0 & 0 & 0 & 0 & 0 & 0 & 0 & 0 & 0 & 0 & 0 & 0 & 0 & 0 & 0 & 0 & 0 & 0 & 0 & 0 & 0 \\
& 0 & \bar{G} & 0 & 0 & 0 & 0 & 0 & 0 & 0 & 0 & 0 & 0 & 0 & 0 & 0 & 0 & 0 & 0 & 0 & \bar{G} & 0 & 0 \\
& 0 & 0 & 0 & 0 & 0 & 0 & 0 & 0 & 0 & 0 & 0 & 0 & 0 & 0 & 0 & 0 & 0 & 0 & 0 & 0 & 0 & 0 \\
& 0 & 0 & 0 & 0 & 0 & 0 & 0 & 0 & 0 & 0 & 0 & 0 & 0 & 0 & 0 & 0 & 0 & 0 & 0 & 0 & \bar{G} & 0 \\
& 0 & 0 & \bar{G} & 0 & 0 & 0 & 0 & 0 & 0 & 0 & x & 0 & 0 & 0 & 0 & 0 & 0 & 0 & 0 & 0 & 0 & \bar{G} \\
& 0 & 0 & 0 & 0 & 0 & 0 & 0 & 0 & 0 & 0 & 0 & 0 & 0 & 0 & 0 & 0 & 0 & 0 & 0 & 0 & 0 & 0 \\
& 0 & 0 & 0 & 0 & 0 & 0 & 0 & 0 & 0 & 0 & 0 & 0 & 0 & 0 & 0 & 0 & 0 & 0 & 0 & 0 & 0 & 0
}\]
\newpage\text{ }\\
\[\kbordermatrix{ & A & A" & A[A] & A'[A] & A"[A] & A[B] & A"[B] & A[X] & A"[X] & A[Y] & A'[A,A] & B & B" & \underline{B} & \overline{B} & \overline{B}A & \overline{B}A" & \overline{B}C & B_s[A] & \overline{B}[A] & B'[B]\\
C'[Y,C]                   & 0 & 0 & 0 & 0 & 0 & 0 & 0 & 0 & 0 & 0 & 0 & 0 & 0 & 0 & 0 & 0 & 0 & 0 & 0 & 0 & 0 \\
D                         & 0 & 0 & 0 & x & 0 & 0 & 0 & 0 & 0 & 0 & 0 & \bar{F} & 0 & 0 & x & 0 & 0 & 0 & 0 & 0 & x \\
D"                        & 0 & 0 & 0 & 0 & 0 & 0 & 0 & 0 & 0 & 0 & 0 & 0 & 0 & 0 & 0 & 0 & 0 & 0 & 0 & 0 & 0 \\
\underline{D}             & 0 & 0 & 0 & 0 & 0 & 0 & 0 & 0 & 0 & 0 & 0 & 0 & 0 & 0 & 0 & 0 & 0 & 0 & 0 & 0 & 0 \\
\overline{D}              & 0 & 0 & 0 & 0 & 0 & 0 & 0 & 0 & 0 & 0 & 0 & 0 & 0 & 0 & 0 & 0 & 0 & 0 & 0 & 0 & 0 \\
\underline{\overline{D}}  & 0 & 0 & 0 & 0 & 0 & 0 & 0 & 0 & 0 & 0 & 0 & 0 & 0 & 0 & 0 & 0 & 0 & 0 & 0 & 0 & 0 \\
D_\l                      & 0 & 0 & 0 & 0 & 0 & 0 & 0 & 0 & 0 & 0 & 0 & 0 & 0 & 0 & 0 & 0 & 0 & 0 & 0 & 0 & 0 \\
DB"                       & 0 & 0 & 0 & 0 & 0 & 0 & 0 & 0 & 0 & 0 & 0 & 0 & 0 & 0 & 0 & 0 & 0 & 0 & 0 & 0 & 0 \\
D\underline{B}            & 0 & 0 & 0 & 0 & 0 & 0 & 0 & 0 & 0 & 0 & 0 & 0 & 0 & 0 & 0 & 0 & 0 & 0 & 0 & 0 & 0 \\
D\ra\underline{B}         & 0 & 0 & 0 & 0 & 0 & 0 & 0 & 0 & 0 & 0 & 0 & 0 & 0 & 0 & 0 & 0 & 0 & 0 & 0 & 0 & 0 \\
DC"                       & 0 & 0 & 0 & 0 & 0 & 0 & 0 & 0 & 0 & 0 & 0 & 0 & 0 & 0 & 0 & 0 & 0 & 0 & 0 & 0 & 0 \\
\underline{D}D            & 0 & 0 & 0 & x & 0 & 0 & 0 & 0 & 0 & 0 & 0 & \bar{F} & 0 & 0 & x & 0 & 0 & 0 & 0 & 0 & x \\
\overline{D}D             & 0 & 0 & 0 & 0 & 0 & 0 & 0 & 0 & 0 & 0 & 0 & 0 & 0 & x & 0 & 0 & 0 & 0 & 0 & 0 & 0 \\
DY"                       & 0 & 0 & 0 & 0 & 0 & 0 & 0 & 0 & 0 & 0 & 0 & 0 & 0 & 0 & 0 & 0 & 0 & 0 & 0 & 0 & 0 \\
D[B]                      & 0 & 0 & 0 & 0 & 0 & 0 & 0 & 0 & 0 & 0 & 0 & 0 & 0 & 0 & 0 & 0 & 0 & 0 & 0 & 0 & 0 \\
D[C]                      & 0 & 0 & 0 & 0 & 0 & 0 & 0 & 0 & 0 & 0 & 0 & 0 & 0 & 0 & 0 & 0 & 0 & 0 & 0 & 0 & 0 \\
D_c'[D]                     & 0 & 0 & 0 & 0 & 0 & 0 & 0 & 0 & 0 & 0 & 0 & 0 & 0 & 0 & 0 & 0 & 0 & 0 & 0 & 0 & 0 \\
D_c'[D]\ra\underline{B}[D]  & 0 & 0 & 0 & 0 & 0 & 0 & 0 & 0 & 0 & 0 & 0 & 0 & 0 & 0 & 0 & 0 & 0 & 0 & 0 & 0 & 0 \\
D[Y]                      & 0 & 0 & 0 & 0 & 0 & 0 & 0 & 0 & 0 & 0 & 0 & 0 & 0 & 0 & 0 & 0 & 0 & 0 & 0 & 0 & 0 \\
D_c'[B,D]                   & 0 & 0 & 0 & 0 & 0 & 0 & 0 & 0 & 0 & 0 & 0 & 0 & 0 & 0 & 0 & 0 & 0 & 0 & 0 & 0 & 0 \\
D_c'[C,D]                   & 0 & 0 & 0 & 0 & 0 & 0 & 0 & 0 & 0 & 0 & 0 & 0 & 0 & 0 & 0 & 0 & 0 & 0 & 0 & 0 & 0 \\
D_c'[Y,D]                   & 0 & 0 & 0 & 0 & 0 & 0 & 0 & 0 & 0 & 0 & 0 & 0 & 0 & 0 & 0 & 0 & 0 & 0 & 0 & 0 & 0 \\
X                         & 0 & 0 & 0 & 0 & 0 & 0 & 0 & 0 & 0 & 0 & 0 & 0 & 0 & 0 & 0 & 0 & 0 & 0 & 0 & 0 & 0 \\
X"                        & 0 & 0 & 0 & 0 & 0 & 0 & 0 & 0 & 0 & 0 & 0 & 0 & 0 & 0 & 0 & 0 & 0 & 0 & 0 & 0 & 0 \\
\underline{X}             & 0 & 0 & 0 & 0 & 0 & 0 & 0 & 0 & 0 & 0 & 0 & 0 & 0 & 0 & 0 & 0 & 0 & 0 & 0 & 0 & 0 \\
XD                        & 0 & 0 & 0 & 0 & 0 & 0 & 0 & 0 & 0 & 0 & 0 & 0 & 0 & 0 & 0 & 0 & 0 & 0 & 0 & 0 & 0 \\
X'[X]                     & 0 & 0 & 0 & 0 & 0 & 0 & 0 & \bar{G} & 0 & 0 & 0 & 0 & 0 & 0 & 0 & 0 & 0 & 0 & 0 & 0 & 0 \\
X'[A,X]                   & 0 & 0 & 0 & 0 & 0 & 0 & 0 & 0 & \bar{G} & 0 & 0 & 0 & 0 & 0 & 0 & 0 & 0 & 0 & 0 & 0 & 0 \\
Y                         & \bar{G} & 0 & 0 & x & 0 & 0 & 0 & 0 & 0 & 0 & 0 & \bar{F} & 0 & 0 & x & 0 & 0 & 0 & 0 & 0 & x \\
Y"                        & \bar{G} & 0 & 0 & x & 0 & 0 & 0 & 0 & 0 & 0 & 0 & \bar{F} & 0 & 0 & x & 0 & 0 & 0 & 0 & 0 & x \\
\overline{Y}              & 0 & 0 & 0 & 0 & 0 & 0 & 0 & 0 & 0 & 0 & 0 & 0 & 0 & 0 & 0 & \bar{G} & 0 & \bar{G} & 0 & 0 & 0 \\
Y"[C]                     & 0 & 0 & 0 & 0 & 0 & 0 & 0 & 0 & 0 & 0 & 0 & 0 & 0 & 0 & 0 & 0 & 0 & 0 & 0 & 0 & 0 \\
Y"[D]                     & 0 & 0 & 0 & 0 & 0 & 0 & 0 & 0 & 0 & 0 & 0 & 0 & 0 & 0 & 0 & 0 & 0 & 0 & 0 & 0 & 0 \\
Y'[Y]                     & 0 & 0 & 0 & 0 & 0 & 0 & 0 & 0 & 0 & \bar{G} & 0 & 0 & 0 & 0 & 0 & 0 & 0 & 0 & 0 & 0 & 0 \\
Z                         & 0 & 0 & 0 & 0 & 0 & 0 & 0 & 0 & 0 & 0 & 0 & 0 & 0 & 0 & 0 & 0 & 0 & 0 & 0 & 0 & 0 \\
ZC"                       & 0 & 0 & 0 & 0 & 0 & 0 & 0 & 0 & 0 & 0 & 0 & 0 & 0 & 0 & 0 & 0 & 0 & 0 & 0 & 0 & 0 \\
Z\ra\underline{X}         & 0 & 0 & 0 & 0 & 0 & 0 & 0 & 0 & 0 & 0 & 0 & 0 & 0 & 0 & 0 & 0 & 0 & 0 & 0 & 0 & 0 \\
ZY"                       & 0 & 0 & 0 & 0 & 0 & 0 & 0 & 0 & 0 & 0 & 0 & 0 & 0 & 0 & 0 & 0 & 0 & 0 & 0 & 0 & 0 \\
Z[C]                      & 0 & 0 & 0 & 0 & 0 & 0 & 0 & 0 & 0 & 0 & 0 & 0 & 0 & 0 & 0 & 0 & 0 & 0 & 0 & 0 & 0 \\
Z[X]                      & 0 & 0 & 0 & 0 & 0 & 0 & 0 & 0 & 0 & 0 & 0 & 0 & 0 & 0 & 0 & 0 & 0 & 0 & 0 & 0 & 0 \\
Z[Y]                      & 0 & 0 & 0 & 0 & 0 & 0 & 0 & 0 & 0 & 0 & 0 & 0 & 0 & 0 & 0 & 0 & 0 & 0 & 0 & 0 & 0
}\]
\newpage
\[\kbordermatrix{ & B_s[C] & B"[C] & \underline{B}[C] & B_s[D] & B"[D] & \underline{B}[D] & B'[A,B] & C & C" & C\ra\{\underline{B},\underline{D},\underline{\overline{D}},\underline{X}\}& CD & C[B] & C[C] & C'[C] & C'[C]\ra\underline{B}[C] & C"[C] & C"[D] & C[X] & C[Y]\\
& 0 & 0 & 0 & 0 & 0 & 0 & 0 & 0 & 0 & 0 & 0 & 0 & 0 & 0 & 0 & 0 & 0 & 0 & 0\\
& 0 & 0 & 0 & 0 & 0 & 0 & 0 & \bar{G} & 0 & 0 & 0 & 0 & 0 & \bar{G} & 0 & 0 & 0 & 0 & 0\\
& 0 & 0 & 0 & 0 & 0 & 0 & 0 & 0 & 0 & 0 & 0 & 0 & 0 & 0 & 0 & 0 & 0 & 0 & 0\\
& 0 & 0 & 0 & 0 & 0 & 0 & 0 & 0 & 0 & 0 & 0 & 0 & 0 & 0 & 0 & 0 & 0 & 0 & 0\\
& 0 & 0 & 0 & 0 & 0 & 0 & 0 & 0 & 0 & 0 & 0 & 0 & 0 & 0 & 0 & 0 & 0 & 0 & 0\\
& 0 & 0 & 0 & 0 & 0 & 0 & 0 & 0 & 0 & 0 & 0 & 0 & 0 & 0 & 0 & 0 & 0 & 0 & 0\\
& 0 & 0 & 0 & 0 & 0 & 0 & 0 & 0 & 0 & 0 & 0 & 0 & 0 & 0 & 0 & 0 & 0 & 0 & 0\\
& 0 & 0 & 0 & 0 & 0 & 0 & 0 & 0 & 0 & 0 & 0 & 0 & 0 & 0 & 0 & 0 & 0 & 0 & 0\\
& 0 & 0 & 0 & 0 & 0 & 0 & 0 & 0 & 0 & 0 & 0 & 0 & 0 & 0 & 0 & 0 & 0 & 0 & 0\\
& 0 & 0 & 0 & 0 & 0 & 0 & 0 & 0 & 0 & 0 & 0 & 0 & 0 & 0 & 0 & 0 & 0 & 0 & 0\\
& 0 & 0 & 0 & 0 & 0 & 0 & 0 & 0 & 0 & 0 & 0 & 0 & 0 & 0 & 0 & 0 & 0 & 0 & 0\\
& 0 & 0 & 0 & 0 & 0 & 0 & 0 & \bar{G} & 0 & 0 & 0 & 0 & 0 & \bar{G} & 0 & 0 & 0 & 0 & 0\\
& 0 & 0 & 0 & 0 & 0 & 0 & 0 & 0 & 0 & \bar{G} & 0 & 0 & 0 & 0 & \bar{G} & 0 & 0 & 0 & 0\\
& 0 & 0 & 0 & 0 & 0 & 0 & 0 & 0 & 0 & 0 & 0 & 0 & 0 & 0 & 0 & 0 & 0 & 0 & 0\\
& 0 & 0 & 0 & 0 & 0 & 0 & 0 & 0 & 0 & 0 & 0 & 0 & 0 & 0 & 0 & 0 & 0 & 0 & 0\\
& 0 & 0 & 0 & 0 & 0 & 0 & 0 & 0 & 0 & 0 & 0 & 0 & 0 & 0 & 0 & 0 & 0 & 0 & 0\\
& 0 & 0 & 0 & x & 0 & 0 & 0 & 0 & 0 & 0 & 0 & 0 & 0 & 0 & 0 & 0 & 0 & 0 & 0\\
& 0 & 0 & 0 & 0 & 0 & x & 0 & 0 & 0 & 0 & 0 & 0 & 0 & 0 & 0 & 0 & 0 & 0 & 0\\
& 0 & 0 & 0 & 0 & 0 & 0 & 0 & 0 & 0 & 0 & 0 & 0 & 0 & 0 & 0 & 0 & 0 & 0 & 0\\
& 0 & 0 & 0 & 0 & x & 0 & 0 & 0 & 0 & 0 & 0 & 0 & 0 & 0 & 0 & 0 & 0 & 0 & 0\\
& 0 & 0 & 0 & 0 & 0 & 0 & 0 & 0 & 0 & 0 & 0 & 0 & 0 & 0 & 0 & 0 & x & 0 & 0\\
& 0 & 0 & 0 & 0 & 0 & 0 & 0 & 0 & 0 & 0 & 0 & 0 & 0 & 0 & 0 & 0 & 0 & 0 & 0\\
& 0 & 0 & 0 & 0 & 0 & 0 & 0 & 0 & 0 & 0 & 0 & 0 & 0 & 0 & 0 & 0 & 0 & 0 & 0\\
& 0 & 0 & 0 & 0 & 0 & 0 & 0 & 0 & 0 & 0 & 0 & 0 & 0 & 0 & 0 & 0 & 0 & 0 & 0\\
& 0 & 0 & 0 & 0 & 0 & 0 & 0 & 0 & 0 & 0 & 0 & 0 & 0 & 0 & 0 & 0 & 0 & 0 & 0\\
& 0 & 0 & 0 & 0 & 0 & 0 & 0 & 0 & 0 & 0 & 0 & 0 & 0 & 0 & 0 & 0 & 0 & 0 & 0\\
& 0 & 0 & 0 & 0 & 0 & 0 & 0 & 0 & 0 & 0 & 0 & 0 & 0 & 0 & 0 & 0 & 0 & \bar{G} & 0\\
& 0 & 0 & 0 & 0 & 0 & 0 & 0 & 0 & 0 & 0 & 0 & 0 & 0 & 0 & 0 & 0 & 0 & 0 & 0\\
& 0 & 0 & 0 & 0 & 0 & 0 & 0 & \bar{G} & 0 & 0 & 0 & 0 & 0 & \bar{G} & 0 & 0 & 0 & 0 & 0\\
& 0 & 0 & 0 & 0 & 0 & 0 & 0 & \bar{G} & 0 & 0 & 0 & 0 & 0 & \bar{G} & 0 & 0 & 0 & 0 & 0\\
& 0 & 0 & 0 & 0 & 0 & 0 & 0 & 0 & 0 & 0 & 0 & 0 & 0 & 0 & 0 & 0 & 0 & 0 & 0\\
& 0 & 0 & 0 & 0 & 0 & 0 & 0 & 0 & x & 0 & 0 & 0 & \bar{G} & 0 & 0 & 0 & 0 & 0 & 0\\
& 0 & 0 & 0 & 0 & 0 & 0 & 0 & 0 & 0 & 0 & 0 & 0 & 0 & 0 & 0 & 0 & 0 & 0 & 0\\
& 0 & 0 & 0 & 0 & 0 & 0 & 0 & 0 & 0 & 0 & 0 & 0 & 0 & 0 & 0 & 0 & 0 & 0 & \bar{G}\\
& 0 & 0 & 0 & 0 & 0 & 0 & 0 & 0 & 0 & 0 & 0 & 0 & 0 & 0 & 0 & 0 & 0 & 0 & 0\\
& 0 & 0 & 0 & 0 & 0 & 0 & 0 & 0 & 0 & 0 & 0 & 0 & 0 & 0 & 0 & 0 & 0 & 0 & 0\\
& 0 & 0 & 0 & 0 & 0 & 0 & 0 & 0 & 0 & 0 & 0 & 0 & 0 & 0 & 0 & 0 & 0 & 0 & 0\\
& 0 & 0 & 0 & 0 & 0 & 0 & 0 & 0 & 0 & 0 & 0 & 0 & 0 & 0 & 0 & 0 & 0 & 0 & 0\\
& 0 & 0 & 0 & 0 & 0 & 0 & 0 & 0 & 0 & 0 & 0 & 0 & 0 & 0 & 0 & 0 & 0 & 0 & 0\\
& 0 & 0 & 0 & 0 & 0 & 0 & 0 & 0 & 0 & 0 & 0 & 0 & 0 & 0 & 0 & 0 & 0 & 0 & 0\\
& 0 & 0 & 0 & 0 & 0 & 0 & 0 & 0 & 0 & 0 & 0 & 0 & 0 & 0 & 0 & 0 & 0 & 0 & 0
}\]
\newpage\text{ }\\
\[\kbordermatrix{ & C'[B,C] & C'[C,C] & C'[Y,C] & D & D" & \underline{D} & \overline{D} & \underline{\overline{D}} & D_\l & DB" & D\underline{B} & D\ra\underline{B} & DC" & \underline{D}D & \overline{D}D & DY" & D[B] & D[C] & D_c'[D] & D_c'[D]\ra\underline{B}[D] & D[Y]\\
& 0 & 0 & 0 & 0 & 0 & 0 & 0 & 0 & 0 & 0 & 0 & 0 & 0 & 0 & 0 & 0 & 0 & 0 & 0 & 0 & 0 \\
& 0 & 0 & 0 & \bar{G} & 0 & 0 & x & 0 & \bar{G} & 0 & 0 & 0 & 0 & 0 & 0 & 0 & 0 & 0 & \bar{G} & 0 & 0 \\
& 0 & 0 & 0 & 0 & 0 & 0 & 0 & 0 & \bar{G} & 0 & 0 & 0 & 0 & 0 & 0 & 0 & 0 & 0 & 0 & 0 & 0 \\
& 0 & 0 & 0 & 0 & 0 & 0 & 0 & 0 & 0 & 0 & 0 & 0 & 0 & \bar{G} & 0 & 0 & 0 & 0 & 0 & 0 & 0 \\
& 0 & 0 & 0 & 0 & 0 & 0 & 0 & 0 & 0 & 0 & 0 & 0 & 0 & 0 & \bar{G} & 0 & 0 & 0 & 0 & 0 & 0 \\
& 0 & 0 & 0 & 0 & 0 & 0 & 0 & 0 & 0 & 0 & 0 & 0 & 0 & 0 & \bar{G} & 0 & 0 & 0 & 0 & 0 & 0 \\
& 0 & 0 & 0 & 0 & 0 & 0 & 0 & 0 & 0 & 0 & 0 & 0 & 0 & 0 & 0 & 0 & 0 & 0 & 0 & 0 & 0 \\
& 0 & 0 & 0 & 0 & 0 & 0 & 0 & 0 & \bar{G} & 0 & 0 & 0 & 0 & 0 & 0 & 0 & 0 & 0 & 0 & 0 & 0 \\
& 0 & 0 & 0 & 0 & 0 & 0 & 0 & 0 & \bar{G} & 0 & 0 & 0 & 0 & 0 & 0 & 0 & 0 & 0 & 0 & 0 & 0 \\
& 0 & 0 & 0 & 0 & 0 & 0 & 0 & 0 & 0 & 0 & 0 & 0 & 0 & 0 & 0 & 0 & 0 & 0 & 0 & 0 & 0 \\
& 0 & 0 & 0 & 0 & 0 & 0 & 0 & 0 & \bar{G} & 0 & 0 & 0 & 0 & 0 & 0 & 0 & 0 & 0 & 0 & 0 & 0 \\
& 0 & 0 & 0 & \bar{G} & 0 & 0 & x & 0 & 0 & 0 & 0 & 0 & 0 & 0 & 0 & 0 & 0 & 0 & \bar{G} & 0 & 0 \\
& 0 & 0 & 0 & 0 & 0 & x & 0 & x & 0 & 0 & 0 & \bar{G} & 0 & 0 & 0 & 0 & 0 & 0 & 0 & \bar{G} & 0 \\
& 0 & 0 & 0 & 0 & 0 & 0 & 0 & 0 & \bar{G} & 0 & 0 & 0 & 0 & 0 & 0 & 0 & 0 & 0 & 0 & 0 & 0 \\
& 0 & 0 & 0 & 0 & 0 & 0 & 0 & 0 & 0 & 0 & x & 0 & 0 & 0 & 0 & 0 & 0 & 0 & 0 & 0 & 0 \\
& 0 & 0 & 0 & 0 & 0 & 0 & 0 & 0 & 0 & 0 & 0 & 0 & x & 0 & 0 & 0 & 0 & 0 & 0 & 0 & 0 \\
& 0 & 0 & 0 & 0 & 0 & 0 & 0 & 0 & 0 & 0 & 0 & 0 & 0 & 0 & 0 & 0 & 0 & 0 & 0 & 0 & 0 \\
& 0 & 0 & 0 & 0 & 0 & 0 & 0 & 0 & 0 & 0 & 0 & 0 & 0 & 0 & 0 & 0 & 0 & 0 & 0 & 0 & 0 \\
& 0 & 0 & 0 & 0 & 0 & 0 & 0 & 0 & 0 & 0 & 0 & 0 & 0 & 0 & 0 & x & 0 & 0 & 0 & 0 & 0 \\
& 0 & 0 & 0 & 0 & 0 & 0 & 0 & 0 & 0 & 0 & 0 & 0 & 0 & 0 & 0 & 0 & 0 & 0 & 0 & 0 & 0 \\
& 0 & 0 & 0 & 0 & 0 & 0 & 0 & 0 & 0 & 0 & 0 & 0 & 0 & 0 & 0 & 0 & 0 & 0 & 0 & 0 & 0 \\
& 0 & 0 & 0 & 0 & 0 & 0 & 0 & 0 & 0 & 0 & 0 & 0 & 0 & 0 & 0 & 0 & 0 & 0 & 0 & 0 & 0 \\
& 0 & 0 & 0 & 0 & 0 & 0 & 0 & 0 & 0 & 0 & 0 & 0 & 0 & 0 & 0 & 0 & 0 & 0 & 0 & 0 & 0 \\
& 0 & 0 & 0 & 0 & 0 & 0 & 0 & 0 & 0 & 0 & 0 & 0 & 0 & 0 & 0 & 0 & 0 & 0 & 0 & 0 & 0 \\
& 0 & 0 & 0 & 0 & 0 & 0 & 0 & 0 & 0 & 0 & 0 & 0 & 0 & 0 & 0 & 0 & 0 & 0 & 0 & 0 & 0 \\
& 0 & 0 & 0 & 0 & 0 & 0 & 0 & 0 & 0 & 0 & 0 & 0 & 0 & 0 & 0 & 0 & 0 & 0 & 0 & 0 & 0 \\
& 0 & 0 & 0 & 0 & 0 & 0 & 0 & 0 & 0 & 0 & 0 & 0 & 0 & 0 & 0 & 0 & 0 & 0 & 0 & 0 & 0 \\
& 0 & 0 & 0 & 0 & 0 & 0 & 0 & 0 & 0 & 0 & 0 & 0 & 0 & 0 & 0 & 0 & 0 & 0 & 0 & 0 & 0 \\
& 0 & 0 & 0 & \bar{G} & 0 & 0 & x & 0 & 0 & 0 & 0 & 0 & 0 & 0 & 0 & 0 & 0 & 0 & \bar{G} & 0 & 0 \\
& 0 & 0 & 0 & \bar{G} & 0 & 0 & x & 0 & 0 & 0 & 0 & 0 & 0 & 0 & 0 & 0 & 0 & 0 & \bar{G} & 0 & 0 \\
& 0 & 0 & 0 & 0 & 0 & x & 0 & x & 0 & 0 & 0 & 0 & 0 & 0 & 0 & 0 & 0 & 0 & 0 & 0 & 0 \\
& 0 & \bar{G} & 0 & 0 & 0 & 0 & 0 & 0 & 0 & 0 & 0 & 0 & 0 & 0 & 0 & 0 & 0 & \bar{G} & 0 & 0 & 0 \\
& 0 & 0 & 0 & 0 & x & 0 & 0 & 0 & 0 & 0 & 0 & 0 & 0 & 0 & 0 & 0 & 0 & 0 & 0 & 0 & 0 \\
& 0 & 0 & \bar{G} & 0 & 0 & 0 & 0 & 0 & 0 & 0 & 0 & 0 & 0 & 0 & 0 & 0 & 0 & 0 & 0 & 0 & \bar{G} \\
& 0 & 0 & 0 & 0 & 0 & 0 & 0 & 0 & 0 & 0 & 0 & 0 & 0 & 0 & 0 & 0 & 0 & 0 & 0 & 0 & 0 \\
& 0 & 0 & 0 & 0 & 0 & 0 & 0 & 0 & 0 & 0 & 0 & 0 & 0 & 0 & 0 & 0 & 0 & 0 & 0 & 0 & 0 \\
& 0 & 0 & 0 & 0 & 0 & 0 & 0 & 0 & 0 & 0 & 0 & 0 & 0 & 0 & 0 & 0 & 0 & 0 & 0 & 0 & 0 \\
& 0 & 0 & 0 & 0 & 0 & 0 & 0 & 0 & 0 & 0 & 0 & 0 & 0 & 0 & 0 & 0 & 0 & 0 & 0 & 0 & 0 \\
& 0 & 0 & 0 & 0 & 0 & 0 & 0 & 0 & 0 & 0 & 0 & 0 & 0 & 0 & 0 & 0 & 0 & 0 & 0 & 0 & 0 \\
& 0 & 0 & 0 & 0 & 0 & 0 & 0 & 0 & 0 & 0 & 0 & 0 & 0 & 0 & 0 & 0 & 0 & 0 & 0 & 0 & 0 \\
& 0 & 0 & 0 & 0 & 0 & 0 & 0 & 0 & 0 & 0 & 0 & 0 & 0 & 0 & 0 & 0 & 0 & 0 & 0 & 0 & 0
}\]
\newpage
\[\kbordermatrix{ & D_c'[B,D] & D_c'[C,D] & D_c'[Y,D] & X & X" & \underline{X} & XD & X'[X] & X'[A,X] & Y & Y" & \overline{Y} & Y"[C] & Y"[D] & Y'[Y] & Z & ZC" & Z\ra\underline{X} & ZY" & Z[C] & Z[X] & Z[Y]\\
& 0 & 0 & 0 & 0 & 0 & 0 & 0 & 0 & 0 & 0 & 0 & 0 & x & 0 & 0 & 0 & 0 & 0 & 0 & 0 & 0 & 0 \\
& 0 & 0 & 0 & 0 & 0 & 0 & 0 & x & 0 & 0 & 0 & 0 & 0 & 0 & 0 & \bar{G} & 0 & 0 & 0 & 0 & 0 & 0 \\
& 0 & 0 & 0 & 0 & 0 & 0 & 0 & 0 & 0 & 0 & 0 & 0 & 0 & 0 & 0 & 0 & 0 & 0 & 0 & 0 & 0 & 0 \\
& 0 & 0 & 0 & 0 & 0 & 0 & 0 & 0 & 0 & 0 & 0 & 0 & 0 & 0 & 0 & 0 & 0 & 0 & 0 & 0 & 0 & 0 \\
& 0 & 0 & 0 & 0 & 0 & 0 & 0 & 0 & 0 & 0 & 0 & 0 & 0 & 0 & 0 & 0 & 0 & 0 & 0 & 0 & 0 & 0 \\
& 0 & 0 & 0 & 0 & 0 & 0 & 0 & 0 & 0 & 0 & 0 & 0 & 0 & 0 & 0 & 0 & 0 & 0 & 0 & 0 & 0 & 0 \\
& 0 & 0 & 0 & 0 & 0 & 0 & 0 & 0 & 0 & 0 & 0 & 0 & 0 & 0 & 0 & 0 & 0 & 0 & 0 & 0 & 0 & 0 \\
& 0 & 0 & 0 & 0 & 0 & 0 & 0 & 0 & 0 & 0 & 0 & 0 & 0 & 0 & 0 & 0 & 0 & 0 & 0 & 0 & 0 & 0 \\
& 0 & 0 & 0 & 0 & 0 & 0 & 0 & 0 & 0 & 0 & 0 & 0 & 0 & 0 & 0 & 0 & 0 & 0 & 0 & 0 & 0 & 0 \\
& 0 & 0 & 0 & 0 & 0 & 0 & 0 & 0 & 0 & 0 & 0 & 0 & 0 & 0 & 0 & 0 & 0 & 0 & 0 & 0 & 0 & 0 \\
& 0 & 0 & 0 & 0 & 0 & 0 & 0 & 0 & 0 & 0 & 0 & 0 & 0 & 0 & 0 & 0 & 0 & 0 & 0 & 0 & 0 & 0 \\
& 0 & 0 & 0 & 0 & 0 & 0 & 0 & x & 0 & 0 & 0 & 0 & 0 & 0 & 0 & \bar{G} & 0 & 0 & 0 & 0 & 0 & 0 \\
& 0 & 0 & 0 & 0 & 0 & 0 & 0 & 0 & 0 & 0 & 0 & 0 & 0 & 0 & 0 & 0 & 0 & \bar{G} & 0 & 0 & 0 & 0 \\
& 0 & 0 & 0 & 0 & 0 & 0 & 0 & 0 & 0 & 0 & 0 & 0 & 0 & 0 & 0 & 0 & 0 & 0 & 0 & 0 & 0 & 0 \\
& 0 & 0 & 0 & 0 & 0 & 0 & 0 & 0 & 0 & 0 & 0 & 0 & 0 & 0 & 0 & 0 & 0 & 0 & 0 & 0 & 0 & 0 \\
& 0 & 0 & 0 & 0 & 0 & 0 & 0 & 0 & 0 & 0 & 0 & 0 & 0 & 0 & 0 & 0 & 0 & 0 & 0 & 0 & 0 & 0 \\
& 0 & 0 & 0 & 0 & 0 & 0 & 0 & 0 & 0 & 0 & 0 & 0 & 0 & 0 & 0 & 0 & 0 & 0 & 0 & 0 & 0 & 0 \\
& 0 & 0 & 0 & 0 & 0 & 0 & 0 & 0 & 0 & 0 & 0 & 0 & 0 & 0 & 0 & 0 & 0 & 0 & 0 & 0 & 0 & 0 \\
& 0 & 0 & 0 & 0 & 0 & 0 & 0 & 0 & 0 & 0 & 0 & 0 & 0 & 0 & 0 & 0 & 0 & 0 & 0 & 0 & 0 & 0 \\
& 0 & 0 & 0 & 0 & 0 & 0 & 0 & 0 & 0 & 0 & 0 & 0 & 0 & 0 & 0 & 0 & 0 & 0 & 0 & 0 & 0 & 0 \\
& 0 & 0 & 0 & 0 & 0 & 0 & 0 & 0 & 0 & 0 & 0 & 0 & 0 & 0 & 0 & 0 & 0 & 0 & 0 & 0 & 0 & 0 \\
& 0 & 0 & 0 & 0 & 0 & 0 & 0 & 0 & 0 & 0 & 0 & 0 & 0 & x & 0 & 0 & 0 & 0 & 0 & 0 & 0 & 0 \\
& 0 & 0 & 0 & 0 & 0 & 0 & \bar{G} & 0 & 0 & 0 & 0 & 0 & 0 & 0 & 0 & 0 & 0 & 0 & 0 & 0 & 0 & 0 \\
& 0 & 0 & 0 & 0 & 0 & 0 & \bar{G} & 0 & 0 & 0 & 0 & 0 & 0 & 0 & 0 & 0 & 0 & 0 & 0 & 0 & 0 & 0 \\
& 0 & 0 & 0 & 0 & 0 & 0 & \bar{G} & 0 & 0 & 0 & 0 & 0 & 0 & 0 & 0 & 0 & 0 & 0 & 0 & 0 & 0 & 0 \\
& 0 & 0 & 0 & 0 & 0 & 0 & 0 & 0 & 0 & x & 0 & x & 0 & 0 & x & 0 & 0 & 0 & 0 & 0 & 0 & 0 \\
& 0 & 0 & 0 & 0 & 0 & 0 & 0 & 0 & 0 & 0 & 0 & 0 & 0 & 0 & 0 & 0 & 0 & 0 & 0 & 0 & \bar{G} & 0 \\
& 0 & 0 & 0 & 0 & 0 & 0 & 0 & 0 & 0 & 0 & 0 & 0 & 0 & 0 & 0 & 0 & 0 & 0 & 0 & 0 & 0 & 0 \\
& 0 & 0 & 0 & 0 & 0 & 0 & 0 & x & 0 & 0 & 0 & 0 & 0 & 0 & 0 & \bar{G} & 0 & 0 & 0 & 0 & 0 & 0 \\
& 0 & 0 & 0 & 0 & 0 & 0 & 0 & x & 0 & 0 & 0 & 0 & 0 & 0 & 0 & \bar{G} & 0 & 0 & 0 & 0 & 0 & 0 \\
& 0 & 0 & 0 & 0 & 0 & 0 & 0 & 0 & 0 & 0 & 0 & 0 & 0 & 0 & 0 & 0 & 0 & 0 & 0 & 0 & 0 & 0 \\
& 0 & \bar{G} & 0 & 0 & 0 & 0 & 0 & 0 & 0 & 0 & 0 & 0 & 0 & 0 & 0 & 0 & 0 & 0 & 0 & \bar{G} & 0 & 0 \\
& 0 & 0 & 0 & 0 & 0 & 0 & 0 & 0 & 0 & 0 & 0 & 0 & 0 & 0 & 0 & 0 & 0 & 0 & 0 & 0 & 0 & 0 \\
& 0 & 0 & \bar{G} & 0 & 0 & 0 & 0 & 0 & 0 & 0 & x & 0 & 0 & 0 & 0 & 0 & 0 & 0 & 0 & 0 & 0 & \bar{G} \\
& 0 & 0 & 0 & x & 0 & 0 & 0 & 0 & 0 & 0 & 0 & 0 & 0 & 0 & 0 & 0 & 0 & 0 & 0 & 0 & 0 & 0 \\
& 0 & 0 & 0 & 0 & 0 & 0 & \bar{G} & 0 & 0 & 0 & 0 & 0 & 0 & 0 & 0 & 0 & 0 & 0 & 0 & 0 & 0 & 0 \\
& 0 & 0 & 0 & 0 & 0 & x & 0 & 0 & 0 & 0 & 0 & 0 & 0 & 0 & 0 & 0 & 0 & 0 & 0 & 0 & 0 & 0 \\
& 0 & 0 & 0 & 0 & 0 & 0 & \bar{G} & 0 & 0 & 0 & 0 & 0 & 0 & 0 & 0 & 0 & 0 & 0 & 0 & 0 & 0 & 0 \\
& 0 & 0 & 0 & 0 & 0 & 0 & 0 & 0 & 0 & 0 & 0 & 0 & 0 & 0 & 0 & 0 & \bar{G} & 0 & 0 & 0 & 0 & 0 \\
& 0 & 0 & 0 & 0 & x & 0 & 0 & 0 & 0 & 0 & 0 & 0 & 0 & 0 & 0 & 0 & 0 & 0 & 0 & 0 & 0 & 0 \\
& 0 & 0 & 0 & 0 & 0 & 0 & 0 & 0 & 0 & 0 & 0 & 0 & 0 & 0 & 0 & 0 & 0 & 0 & x & 0 & 0 & 0
}\]
\end{tiny}
\end{landscape}
\newpage
\section{Python code}\label{app:B}
In this section, we provide Python programming code which demonstrate the encoding of a simple permutation in $H'$ of length $n$ into a $n$-letter word word in $L'$. Together with PermLab \cite{PermLab1.0} and all Python files provided in this section, we are able to construct a list of all $n$-letter words in $L'$ as a text file where $n$ is a reasonable length.\\

Note that the code was written first, and was used greatly to discover the encoding rules stated in Chapter \ref{chap:6}. What we provide here is kept as original, and was not modified even after proofs in Chapter \ref{chap:6} were established. For this reason, there may be a slight difference in the process of encoding between Python code and the sequence of encoding defined in Chapter \ref{chap:6}, but the obtained result is the same. For instance, \verb+prefixConverter.py+ and \verb+suffixConverter.py+ take the places of \textsf{W-COMBINE} and \textsf{AFFIX-CONVERT}.\\
\\
\verb+makeList.py+
\begin{lstlisting}[language=python]
def readFile(fname):
    """
    reads the list of permutations from txt file line by line.
    """
    with open (fname) as f:
        permList = f.read().splitlines()
    return permList

def convList1(permList):
    """
    converts the list of string type permutations in permList
    to int type, specifically length <10.
    """
    listofperm=[]
    for i in range(len(permList)):
        temp=[]
        for k in range(len(permList[i])):
            temp.append(int(permList[i][k]))
        listofperm.append(temp)
    return listofperm

def convList2(permList):
    """
    converts the list of string type permutations in permList
    to int type, specifically length >=10.
    """
    listofperm=[]
    for i in range(len(permList)):
        temp = permList[i].split()
        for k in range(len(temp)):
            temp[k]=int(temp[k])
        listofperm.append(temp)
    return listofperm
\end{lstlisting}
\vhh
\verb+list2str.py+
\begin{lstlisting}[language=python]
"""
The following function converts a list to a string.
"""

def list2Str(word):
    string = ' '.join(word)
    return string

def test(blah):
    word = []
    while len(blah) > 0:
        i = 0
        while blah[i] != ' ':
            i += 1
        temp = ''.join(blah[0:i])
        word.append(temp)
        temp2 = temp2[i+1:len(temp2)]
    return word
\end{lstlisting}
\vhh
\verb+invert.py+
\begin{lstlisting}[language=python]
"""
The following function simply inverts a permutation.
"""

def invPerm(perm):
    inverse = []
    for i in range(1,max(perm)+1):
        if i in perm:
            inverse.append(perm.index(i)+1)
        else:
            pass
    return inverse
\end{lstlisting}
\vhh
\verb+ep.py+
\begin{lstlisting}[language=python]
def ep2413(perm):
    """
    keeps permutations of extreme pattern 2413, and discard others.
    """
    if perm[0]<perm[len(perm)-1] and invPerm(perm)[0]>
    invPerm(perm)[len(invPerm(perm))-1]:
        return perm
    else:
        pass

def ep3142(perm):
    """
    keeps permutations of extreme pattern 3142, and discard others.
    """
    if perm[0]>perm[len(perm)-1] and invPerm(perm)[0]<
    invPerm(perm)[len(invPerm(perm))-1]:
        return perm
    else:
        pass

def Nstart(perm):
    """
    This function eliminates a simple adbi permutation which starts
    with S set.
    """
    if (perm[0] == 2 or perm[0] == 3 or perm[0] == 4) and
    perm[1] != 1:
        return perm
    else:
        pass
\end{lstlisting}
\vhh
\verb+createFile.py+
\begin{lstlisting}[language=python]
"""
The following function creates a text file listing all the desired
words.
"""

from os import path

def createFile(dest,words):
    if not (path.isfile(dest)):
        f = open(dest,'w')
        for i in words:
            f.write(i+'\n')
        f.close
\end{lstlisting}
\vhh
\verb+counter.py+
\begin{lstlisting}[language=python]
"""
The following function will count how many N's and S's are embedded
in a permutation.
"""

def counter(perm):
    d1 = perm[0] #The first number. It is either 2, 3 or 4.#
    d2 = 1 #This is always 1.#
    temp1 = perm[perm.index(d1):perm.index(d2)+1]
    #Get the interval from d1 to d2"
    d3 = max(temp1) #The maximum number of temp 1 is d3.#
    temp2 = perm[0:len(perm)] #Get a copy of a permutation.#
    for x in range(d3+1,max(temp2)+1):
        temp2.remove(x) #Erase everything > d3.#
    d4 = temp2[-1] #The last number of temp 2 is d4.#
    count = 1
    while d2 != d4:
        d1 = d3
        d2 = d4 #Move up d1 and d2.#
        temp1 = perm[perm.index(d1):perm.index(d2)+1]
        d3 = max(temp1) #Determine the new d3.#
        count += 1 #Count goes up.#
        if d1 != d3: #Checking if d1 and d3 are the same.#
        #If not, determine the new d4.#
            temp2 = perm[0:len(perm)]
            for x in range(d3+1,max(temp2)+1):
                temp2.remove(x)
            d4 = temp2[-1]
            count += 1
        else: #If d1 and d3 become the same, we end it here.#
            break
    return count

"""
The number n = 'count' tells us how many N's and S's we have.
If n is odd, then there are (n-1)/2 N's and (n-1)/2 S's in a
permutation. If n is even, then there are n/2 N's and (n/2)-1
S's in a permutation.
"""
\end{lstlisting}
\vhh
\verb+Nsplit.py+
\begin{lstlisting}[language=python]
"""
We define two functions. The first one determines where the first
N shape ends. It returns the split N segment and some extras. The
second one takes the N set away from the whole permutation. It
returns the rest.
"""

def Nsplit(perm):
    """
    First, we detect d1, the greatest element between the first
    and the least numbers.
    """
    temp1 = perm[0:perm.index(1)+1]
    d1 = max(temp1)
    """
    For convenience, we determine d2 as well. Remove everything
    greater than d1. Let the last number be d2.
    """
    temp2 = perm[0:len(perm)]
    for x in range(d1+1,max(temp2)+1):
        temp2.remove(x)
    d2 = temp2[-1]
    """
    Finally, we identify the splitter number. This is necessary
    only when d1 is scissored vertically by a number that is not
    equal to d2.
    """
    splitter = d1-1
    """
    We determine if d1 is the max of the whole permutation. If so,
    d2 is the last number in N-C. If not, we see which one of d1-2,
    d1-1 and d1+1 appears first after d1.
    """
    if d1 == max(perm):
        CLast = d2
        """
        If d1-2 appears first, one of the following three is
        happening.
        Case 1: The following S set is either empty, or just one
        number in C segment.
        Case 2: d1 is scissored vertically by d1-1.
        Case 3: d1-1 and d1-2 are a scission pair, or d1-1 and d2
        are a scission pair with d1-2 being a scissor.
        """
    elif perm.index(d1-2) < perm.index(d1+1) and perm.index(d1-2)
    < perm.index(d1-1):
        temp3 = temp1[0:len(perm)]
        temp3.remove(d1)
        #Case 1a: d1-2 comes from A or B segment.#
        if max(temp3) == d1-2:
        #Detecting if d1-2 belongs to either A or B segment.#
        #In order to distinguish from Case 2, we need to find#
        #what separates d1 and d1-2.#
        #If d1-1 is not d2, then it is Case 2. Otherwise, it is#
        #Case 1a.#
            if d2 != d1-1: #Case 2#
                if perm[perm.index(d1-1)+1] == d2 or
                perm[perm.index(d1-1)+2] == d2:
                #One of scission pair of d2 is being a scissor.#
                    temp4 = perm[perm.index(1):perm.index(d2)+1]
                    temp5 = []
                    for x in temp4:
                        if x < d2:
                            temp5.append(x)
                        else:
                            pass
                    CLast = temp5[-1]
                elif perm[perm.index(d1-1)-2]>d1-1:
                #Case of d d^ c c' b_ ...#
                    CLast = perm[perm.index(d1-1)-3]
                elif perm[perm.index(d1-1)-1]>d1-1:
                #Case of d d^ c' b_ ... or d d^ c b_ ...#
                    CLast = perm[perm.index(d1-1)-2]
                else:
                    CLast = perm[perm.index(d1-1)-1]
            else:
                temp4 = [] #REVISED (3/4/15)#
                for i in perm:
                    if i < d1-2:
                        temp4.append(i)
                    else:
                        pass
                CLast = temp4[-1]
        #Case 1b: d1-2 comes from C segment.#
        elif d1-2 < d2:
            if perm[perm.index(d1-2)+1]<d1-2:
            #Detecting if d1-2 is a part of scission.#
                if perm.index(d1-2)<perm.index(d1-3):
                #Detecting if d1-2 scission pairs are not the last#
                #numbers in C.#
                    CLast = d1-3
                else:
                #d1-2 scission pairs are the last numbers in C.#
                    CLast = perm[perm.index(d1-2)+1]
            else: #d1-2 is not a part of scission.#
                CLast = d1-2
        #The last one is Case 3: d1-1 and d1-2 are a scission pair#
        #in B segment of S shape, or d1-1 and d2 are scission pairs#
        #with d1-2 being a scissor.#
        #Note that for the former case, d1-1 and d1-2 can be y' and#
        #y'', and for the latter case, d1-2 can be y^.#
        else:
            if perm[perm.index(d1-2)+1] > d1-2:
            #No scissor in A segment.#
                CLast = perm[perm.index(d1-2)-1]
            else:
                s1 = perm[perm.index(d1-2)+1]
                s2 = s1+1
                temp6 = perm[perm.index(1):perm.index(s1)+1]
                while s2 in temp6 and perm[perm.index(s2)+1] < s2:
                    s1 = perm[perm.index(s2)+1]
                    s2 = s1+1
                if perm.index(s2) < perm.index(s1) and s2 in temp6:
                #Case of S-A segment self-scission chain. 2nd one#
                #is a scission.#
                    CLast = perm[perm.index(s2)-1]
                elif perm.index(s2) > perm.index(s1):
                #Case of S-A segment self-scission chain. 2nd one#
                #is not a scission.#
                    CLast = perm[perm.index(s1)-2]
                elif perm.index(s2) < perm.index(s1) and s2
                not in temp6:
                #Case of N-C segment self-scission chain with the#
                #last one being z.#
                    CLast = perm[perm.index(d1-2)+1]
        """
        In the case d1+1 appears first, the number immediately
        left of d1+1 is the last number in C segment.
        """
    elif perm.index(d1+1) < perm.index(d1-1) and perm.index(d1+1)
    < perm.index(d1-2):
        CLast = perm[perm.index(d1+1)-1]
        """
        Finally, in the case d1-1 shows up first, there are
        following six cases.
        Case 1: d1 and d1-1 are a scission pair.
        Case 2: d1-1 is encoded as y.
        Case 3: d1-1 scissors a chain of S-A self-scissions or
        one of N-C self-scissions as d1-1 is being y.
        Case 4: d1-1 and d2 are a scission pair.
        Case 5: d1-1 scissors c' and c'' in S set.
        Case 6: Nothing special. Encoding goes d b, or d b'.
        """
    elif perm.index(d1-1) < perm.index(d1+1) and perm.index(d1-1)
    < perm.index(d1-2):
        if d1-1 == perm[perm.index(d1)+2] and perm[perm.index(d1)+1]
        != 1: #Case 1, d1 and d1-1 are a scission pair.#
            temp7 = perm[0:len(perm)]
            for x in range(1,temp7[temp7.index(d1)+1]):
                temp7.remove(x)
            if temp7[temp7.index(d1-1)+1]
            == temp7[temp7.index(d1)+1]+1:
                CLast = temp7[temp7.index(d1)+1]+1
            elif temp7[temp7.index(d1-1)+1]
            == temp7[temp7.index(d1)+1]+2 and
            temp7.index(temp7[temp7.index(d1)+1]+1)
            < temp7.index(temp7[temp7.index(d1)+1]+2):
                CLast = temp7[temp7.index(d1)+1]+2
            else:
                temp7 = perm[0:len(perm)]
                for x in
                range(temp7[temp7.index(d1)+1],max(temp7)+1):
                    temp7.remove(x)
                CLast = temp7[-1]

        elif perm[perm.index(d1-1)+1] < d1-1 and
        perm[perm.index(d1-1)+1]+1 in
        perm[0:perm.index(perm[perm.index(d1-1)+1])+1]:
        #Case 2 and 3.#
            s1 = perm[perm.index(d1-1)+1]
            s2 = s1+1
            temp8 = perm[perm.index(1):perm.index(s1)+1]
            while s2 in temp8 and perm[perm.index(s2)+1] < s2:
                s1 = perm[perm.index(s2)+1]
                s2 = s1+1
            if perm.index(s2) < perm.index(s1) and s2 in temp8:
            #Case 3, S-A segment self-scission chain. 2nd one is#
            #a scission.#
                CLast = perm[perm.index(s2)-1]
            elif perm.index(s2) > perm.index(s1):
            #Case 3, S-A segment self-scission chain. 2nd one is#
            #not a scission.#
                CLast = perm[perm.index(s1)-2]
            elif perm.index(s2) < perm.index(s1) and s2 not in
            temp8: #Case 2, d1-1 is encoded as y. There could be#
            #an N-C segment self-scission chain.#
                CLast = perm[perm.index(d1-1)+1]
        elif d1-2 == d2 and perm[perm.index(d1-1)-1] > d1 and
        perm.index(d2)-2 == perm.index(d1-1): #Case 4, d1-1 and d2#
        #are a scission pair.# #REVISED (3/5/15)#
            temp9 = perm[0:len(perm)]
            for x in range(d2,max(temp9)+1):
                temp9.remove(x)
            CLast = temp9[-1]
        elif perm[perm.index(d1-1)-1] == d1+2: #Case 5, d1-1#
        #scissors c' and c'' in S set.#
            CLast = perm[perm.index(d1-1)-2]
        else: #Case 6, Nothing special is present.#
            CLast = perm[perm.index(d1-1)-1]
    subperm = [] #This is the split subpermutation.#
    for i in range(0,perm.index(CLast)+1):
            subperm.append(perm[i]) #subperm now is the#
            #subpermutation up to the last number in C segment that#
            #we identified.#
            #*** Refer this list as (1) in the latter function.#
    if perm[perm.index(d1)+1] == d1-2 or perm[perm.index(d1)+2]
    == d1-2: #Determine if d1 is scissored.#
        if perm.index(d1-2) > perm.index(1): #REVISED (3/4/15)#
            pass
        else:
            if splitter not in subperm and splitter != d2:
                subperm.append(splitter) #We add the splitter number#
                #in case it belongs to S-B or the scission of d2.#
            else:
                pass
    else:
        pass
    if CLast != d2:
        subperm.append(d2)
        #Finally, add d2 to the subperm list, and we are done.#
    return subperm, CLast

def setMinus(perm,CLast,subword):
    """
    perm is the whole permutation, CLast is the last number in N-C
     obtained by Nsplit function. subword is the encoded
    word after the modification of both prefix/suffix converters.
    What we do with this function are the following:
    1. We first take away everything from perm up to CLast number
    (including itself).
    2. By looking at the suffix of subword, we add back the necessary
    numbers.
    """
    subperm = [] #This is the split subpermutation.#
    for i in perm[0:perm.index(CLast)+1]:
        subperm.append(i) #*** Same list as (1)#
        #from the above function.#
    for i in perm[0:perm.index(CLast)+1]:
        perm.remove(i)
    orig = subperm + perm #We need to refer back in some cases.#
    d0 = min(subperm) #i.e., 1. We always need to add this one back.#
    d1 = max(subperm) #Similarly, we need to add this back.#
    if subword[-2:len(subword)] == ['d^', 'd']: #In this case,#
    #d = d1. We need to add back d^ as well. d^ is simply d1-2.#
        recover = [d1, d1-2, d0]
    elif subword[-2:len(subword)] == ['d_^', 'd']: #Similar to the#
    #above case. d_^ is d1-2.#
        recover = [d1, d1-2, d0]
    elif subword[-3:len(subword)] == ['z', 'x', 'd'] or
    subword[-3:len(subword)] == ['z', 'x"', 'd'] or
    subword[-3:len(subword)] == ['z', 'x_', 'd'] or
    subword[-3:len(subword)] == ['z', 'c"', 'd'] or
    subword[-3:len(subword)] == ['z', 'y"', 'd']:
        for i in subperm[0:subperm.index(d1)+1]:
            subperm.remove(i) #Get rid of A segment, including d1.#
        for i in
        subperm[subperm.index(d0)+1:subperm.index(subperm[-1])+1]:
            subperm.remove(i) #Get rid of C segment,#
            #NOT including d0.#
        x = max(subperm) #This is the max of B segment, i.e., x, x",#
        #x_, y" or b (b") which is#
        #associated with z _ d y encoding. We need to add this#
        #back to the next perm.#
        if x-1 == CLast: #All the cases except x being b or b",#
        #causing the C-chain scissions.(Need revision here???)#
        #In these cases, simply CLast is z. All we need to add back#
        #is d1, x, d0, y and z.#
            y = orig[orig.index(CLast)-1]
            recover = [d1, x, d0, y, CLast]
        else: #In this case, we must add back the whole C chain#
        #including x which initiates the chain.#
            firstcPrime = x-1
            firstcDPrime = orig[orig.index(firstcPrime)-1]
            CChain =
            orig[orig.index(firstcDPrime):orig.index(CLast)+1]
            recover = [d1, x, d0] + CChain
    else:
        recover = [d1, d0]
    perm = recover + perm
    return perm
\end{lstlisting}
\vhh
\verb+encodeNS.py+
\begin{lstlisting}[language=python]
import counter
import ep
from os import path
import invert
import Nsplit

"""
Looking at a permutation, or subsequence embedded in a permutation
which has the extreme pattern 2413 (N shape), the following function
encodes it into a word under certain rules. The Alphabet is the set
containing following 30 letters.

# a, a', a", b, b', b", bs, b_, b^, c, c', c", d, da', da",  #
# dc', dc", d_, d^, d_^, dl, x, x', x", x_, y, y', y", y^, z #

'perm' is an arbitrary permutation of extreme pattern 2413.
'locator' tells us where the subsequence is in a main permutation.
'transit1' is the information about the previous S shape.
'transit2' is the information about the next S shape.
"""

def encodeS(perm):
    """
    encodes a permutation of extreme pattern 3142 into a word.
    """
    if perm[0]>perm[len(perm)-1] and perm.index(len(perm)) >
    perm.index(1):
        word = []
        for i in range(len(perm)):
            if i == 0 or i == len(perm)-1:
                word.append('d') #First and last letters must be d.#
            ##=====Prefix coding=====##
            elif i == 1 and (perm[i] == 1 or perm[i] == 2):
                if perm[i] == 1:
                    word.append('d')
                else:
                    word.append('da\'')
            elif i == 2 and perm[i] == 1:
                word.append('d')
            elif i == 3 and perm[i] == 1:
                word.append('da\"')
            ##=====Surffix coding=====##
            elif i == len(perm)-2 and (perm[i] == len(perm) or
            perm[i] == len(perm)-1):
                if perm[i] == len(perm):
                    word.append('d')
                else:
                    word.append('dc\"')
            elif i == len(perm)-3 and perm[i] == len(perm):
                word.append('d')
            elif i == len(perm)-4 and perm[i] == len(perm):
                word.append('dc\'')
            ##=====Infix coding=====##
            elif perm[i] > perm[0]:
                if perm.index(perm[i]-1) > i:
                    word.append('c\'')
                elif perm.index(perm[i]+1) < i:
                    word.append('c\"')
                else:
                    word.append('c')
            elif perm[i] < perm[0] and perm[i] > perm[len(perm)-1]:
                if perm.index(perm[i]+1) > i:
                    word.append('b\'')
                elif perm.index(perm[i]-1) < i:
                    word.append('b\"')
                else:
                    word.append('b')
            else:
                if perm.index(perm[i]-1) > i:
                    word.append('a\'')
                elif perm.index(perm[i]+1) < i:
                    word.append('a\"')
                else:
                    word.append('a')
        for i in range(len(word)):
            if word[i] == 'b':
                if word[i-1] == 'c\'' or word[i+1] == 'a\"' or
                word[i-1] == 'dc\'' or word[i+1] == 'da\"':
                    word[i] = 'bs'
                else:
                    pass
            else:
                pass
        return word
    else:
        print("Error. The entered permutation does not have extreme
        pattern 3142.")


def encodeN(perm):
    """
    encodes a permutation of extreme pattern 2413 into a word.
    """
    inverse = invert.invPerm(perm)
    word = encodeS(inverse)
    return word
\end{lstlisting}
\vhh
\verb+prefixConverter.py+
\begin{lstlisting}[language=python]
"""
Taking an arbitrary N or S set out of a simple adbi permutation of
extreme pattern 2143, the prefix & suffix of it must be modified
according to:
1. where it comes from, and
2. how it connects to the previous/next N & S sets.
The following functions will modify the prefix of the N or S set that
we are currently encoding.
"""

import Nsplit
import encodeNS

def initialPrefix(word):
    """
    This function is particularly applied if the N set we are
    encoding is the initial set of a whole permutation, but not the
    final set. In addition to the existing 6 ways to start the
    encoded word, there are 4 more possibilities associated with
    x' and b^.
    """
    codedWord = word[0:len(word)]
    if codedWord.count("b'") == 1:
    #This line ensures that the first b' showing up needs to be#
    #modified to x' or b^.#
        if codedWord[0:3] == ['d', "b'", 'd'] and
        codedWord[-5:len(codedWord)] == ["c'", 'b"', 'd', 'c"', 'd']:
            codedWord[1] = "x'"
            #This is the modification for Case 9: d x' d...#
        elif codedWord[0:4] == ['d', "da'", "b'", 'da"'] and
        codedWord[-5:len(codedWord)] == ["c'", 'b"', 'd', 'c"', 'd']:
            codedWord[2] = "x'"
            #This is the modification for Case 10: d da' x' da"...#
        elif codedWord[0:3] == ['d', "b'", 'd'] and
        (codedWord[-3:len(codedWord)] == ['d', 'b"', 'd'] or
        codedWord[-4:len(codedWord)] == ["dc'", 'b"', 'dc"', 'd']):
            codedWord[1] = 'b^'
            #This is the modification for Case 18: d b^ d...#
        elif codedWord[0:4] == ['d', "da'", "b'", 'da"'] and
        (codedWord[-3:len(codedWord)] == ['d', 'b"', 'd'] or
        codedWord[-4:len(codedWord)] == ["dc'", 'b"', 'dc"', 'd']):
            codedWord[2] = 'b^'
            #This is the modification for Case 19: d da' b^ da"...#
        else:
            pass
    else:
        pass
    return codedWord

def noninitialPrefix(word,m,n):
    """
    This function is applied if the N set we are encoding is not the
    initial set of a whole permutation. word is the encoded word by
    encodeN(word) function, m is the current set number
    (3 leq m leq n), and n is the final set (found by counter(perm)
    function).
    """
    codedWord = word[0:len(word)]
    if codedWord[0:2] == ['d', 'd']:
    #This may happen if the previous set ended with:#
    #1. Simply with d d.#
    #2. d_ d. We do not carry d_ to the following set.#
    #We simply erase d d, since they are previously encoded.#
        codedWord = codedWord[2:]
    elif codedWord[0:3] == ['d', "a'", 'd']:
    #This is associated with z x d y interaction.#
    #a' could be various things. The associated a" is z, or . Those#
    #are previously encoded,#
    #so we make sure to erase d a' d and a".#
        codedWord = codedWord[3:]
    #There are three possiblilities immediately after d a' d, namely#
    #bs, b' or a'.#
    #If it is a', then the suffix of the previous set was encoded by#
    #z c" d.#
    #We first get rid of this case.#
        while codedWord[0] != 'bs' and codedWord[0] != "b'":
            codedWord = codedWord[1:]
        codedWord.remove('a"')
    #Now, it's down to two possibilities; bs or b'.#
        if codedWord[0] == 'bs':
        #If it is followed by bs, then all we have to do is to#
        #replace this with y.#
            codedWord[0] = 'y'
        elif codedWord[0] == "b'":
        #If it is followed by b', this should be y' or y^.#
        #Now, if the number of b' & b" pair is 1, AND b" happens#
        #in the suffix,#
        #then this case needs a special attention.#
            if codedWord.count("b'") == 1:
                if codedWord[-3:len(codedWord)] == ['d', 'b"', 'd']:
                #We do not modify 'b"' because it is#
                #included in the suffixConverter.#
                    if m != n:
                    #The current set we are encoding is not the
                    #final set. Then b' -> y^.#
                        codedWord[0] = 'y^'
                    else:
                    #I.e., m == n. In this case, b' -> y'.#
                    #Then d y" d at the end.#
                        codedWord[0] = "y'"
                elif codedWord[-4:len(codedWord)] ==
                ["dc'", 'b"', 'dc"', 'd']:
                #Similar to the above case.#
                #The difference is whether the following b"#
                #becomes d_ (previous case) or d_^ (this case).#
                #Again, We do not modify 'b"' because it is
                #included in the suffixConverter.#
                    if m != n:
                    #The current set we are encoding is not the#
                    #final set. Then b' -> y^.#
                        codedWord[0] = 'y^'
                    else: #I.e., m == n. In this case, b' -> y'.#
                    #Then dc' y" dc" d at the end.#
                        codedWord[0] = "y'"
                elif codedWord[-5:len(codedWord)] ==
                ["c'", 'b"', 'd', 'c"', 'd']:
                #Very special case.#
                #This is the structure of#
                #d a' d b' a" (a) (c) c' b" d c" d, resulting#
                #y' z y" d y encoding. We do not convert b"#
                #since it is a part of the suffix.#
                    codedWord[0] = "y'"
                else:
                #The case b'& b" pair happens once, but b" is not#
                #in the suffix.#
                #In this case, b" -> y" as well as b' -> y'.#
                #Since b" isn't included in#
                #the suffix, we must convert it here.#
                    codedWord[0] = "y'"
                    codedWord[codedWord.index('b"')] = 'y"'
            else:
            #The case b' & b" pair happens more than once.#
            #Same as before, b" -> y" as well as#
            #b' -> y' because b" isn't included in the suffix.#
                codedWord[0] = "y'"
                codedWord[codedWord.index('b"')] = 'y"'
    elif codedWord[0:3] == ['d', "b'", 'd']:
    #This is associated with d^ b_ interaction.#
    #b' is either d^ or d_^. This is previously encoded, so we#
    #simply erase d b' d.#
    #The associated b" can be many different things, such as#
    #b_, d_, d_^, and x_.#
    #If b" is a part of the suffix, then we postpone the#
    #modification till suffixConverter.#
    #If this isn't the case, then we modify b" here.#
        codedWord = codedWord[3:]
        if codedWord.count('b"') == 1:
        #Since we have erased b', we check if b" only shows up once.#
            if codedWord[-3:len(codedWord)] == ['d', 'b"', 'd']:
            #We do not modify b" because it is#
            #included in the suffixConverter.#
                pass
            elif codedWord[-4:len(codedWord)] ==
            ["dc'", 'b"', 'dc"', 'd']:
            #We do not modify b" because it is#
            #included in the suffixConverter.#
                pass
            elif codedWord[-5:len(codedWord)] ==
            ["c'", 'b"', 'd', 'c"', 'd']:
                if m != n:
                #If the set isn't final, we leave b" for the#
                #suffixConverter.#
                    pass
                else: #Otherwise, b" needs to be modified to b_,#
                #since finalSuffix function don't#
                #deal with the case of ending with d c" d.#
                    codedWord[codedWord.index('b"')] = 'b_'
            else: #The case b'& b" pair happens once, but b" is not#
            #in the suffix.#
            #In this case, we must convert b" -> b_.#
                codedWord[codedWord.index('b"')] = 'b_'
        else: #The case b' & b" pair happens more than once. Same as#
        #before, b" -> b_#
        #because b" isn't included in the suffix.#
            codedWord[codedWord.index('b"')] = 'b_'
    else: #Any other prefix should never appear due to the way we#
    #split each set.#
        codedWord = print("Error. Something went wrong...")
    return codedWord
\end{lstlisting}
\vhh
\verb+suffixConverter.py+
\begin{lstlisting}[language=python]
"""
Taking an arbitrary N or S set out of a simple adbi permutation of
extreme pattern 2143, the prefix & suffix of it must be modified
according to:
1. where it comes from, and
2. how it connects to the previous/next N & S sets.
The following functions will modify the suffix of the N or S set
that we are currently encoding.
There is a slight difference from prefixConverter functions. word is
the original, and codedWord is the encoded word after prefixConverter
is applied.
"""

import Nsplit
import encodeNS
import prefixConverter

def finalSuffix(word,codedWord):
    """
    This function is particularly applied if the N set we are
    encoding is the final set of a whole permutation, but not the
    initial set. In addition to the existing 6 ways to end the
    encoded word, there are 4 more possibilities associated with
    y" and b_.
    """
    if word.count("b'") == 1: # This line ensures that the first b'#
    #showing up needs to be modified to x' or b^.#
        if word[-3:len(word)] == ['d', 'b"', 'd'] and word[0:5] ==
        ['d', "a'", 'd', "b'", 'a"']:
            codedWord[-2] = 'y"' #This is the modification for the#
            #reverse of Case 9: ...d y" d ending case.#
        elif word[-4:len(word)] == ["dc'", 'b"', 'dc"', 'd'] and
        word[0:5] == ['d', "a'", 'd', "b'", 'a"']:
            codedWord[-3] = 'y"' #This is the modification for the#
            #reverse of Case 10: ...dc' y" dc" d#
        elif word[-3:len(word)] == ['d', 'b"', 'd'] and (word[0:3] ==
        ['d', "b'", 'd'] or word[0:4] == ['d', "da'", "b'", 'da"']):
            codedWord[-2] = 'b_' #This is the modification for the#
            #reverse of Case 18: ...d b_ d#
        elif word[-4:len(word)] == ["dc'", 'b"', 'dc"', 'd'] and
        (word[0:3] == ['d', "b'", 'd'] or word[0:4] ==
        ['d', "da'", "b'", 'da"']):
            codedWord[-3] = 'b_' #This is the modification for the#
            #reverse of Case 19: ...dc' b_ dc" d#
        else:
            pass
    else:
        pass
    return codedWord

def nonfinalSuffix(word,codedWord,m):
    """
    This function is applied if the N set we are encoding is not the
    final set of a whole permutation. word is the original encoded
    word which is the result of encodeN(word) function, codedWord is
    the encoded word by either initialPrefix(word) or
    noninitialPrefix(word) function, m is the current set number
    (2 leq m leq (# of final set)-1).
    """
    last1 = len(word)
    last2 = len(codedWord) #Just to make things a little easier...#
    if word[-2:last1] == ['d', 'd']: #We erase the second d, since#
    #this will be encoded in the next set.#
        codedWord = codedWord[0:last2-1]
    elif word[-4:last1] == ["dc'", 'bs', 'dc"', 'd']:
    #The case which didn't rise in the prefixConverter.#
    #We need to erase dc' and dc", and convert b -> d^. dc" will be#
    #converted in the next set.#
        codedWord = codedWord[0:last2-4]
        #Simply erase all four of them first.#
        codedWord.append('d^') #Then add d^, and#
        codedWord.append('d') #add back d.#
    elif word[-3:last1] == ['d', 'c"', 'd']:
    #This is associated with z x d y interaction.#
    #c' is z. The associated c" is y. The scissoring element for#
    #these c' & c" could be various things.#
    #We first erase the suffix d c" d, modify whatever needs to be#
    #modified, then add back necessary stuff.#
    #c" and the first d are encoded in the future, so for now, we#
    #will not add these back.#
        word = word[0:last1-3]
        codedWord = codedWord[0:last2-3]
        codedWord[-2] = 'z' #Fix c' to be z here.#
    #There are three possiblilities immediately before d c" d#
    #namely bs, b" or c".#
    #If it is c", then we keep it as c", so we have z c" d.#
        if word[-1] == 'c"':
            codedWord.append('d')
        elif word[-1] == 'bs': #If it is bs, then all we have to do#
        #is to replace this with x.#
            codedWord[-1] = 'x'
            codedWord.append('d')
        elif word[-1] == 'b"':
        #If it is b", this should be x", x_ or y".#
        #Now, if the number of b' & b" pair is 1, AND b' happens in#
        #the prefix,#
        #then this case needs a special attention.#
            if word.count("b'") == 1:
                if word[0:3] == ['d', "b'", 'd']:
                #We do not modify "b'" because it is#
                #included in the prefixConverter.#
                    if m != 2:
                    #The current set we are encoding is not the#
                    #initial set. Then b" -> x_.#
                        codedWord[-1] = 'x_'
                        codedWord.append('d')
                    else: #I.e., m == 2. In this case, b" -> x".#
                    #d x" d should be the actual suffix.#
                        codedWord[-1] = 'x"'
                        codedWord.append('d')
                elif word[0:4] == ['d', "da'", "b'", 'da"']:
                #We do not modify b' because it is#
                #included in the prefixConverter.#
                    codedWord[-1] = 'x"'
                    #Now, Since this particular prefix only#
                    #happens when the set is#
                    #initial, we may assume m == 2 here.#
                    #Thus, this can never be x_.#
                    codedWord.append('d')
                elif word[0:5] == ['d', "a'", 'd', "b'", 'a"']:
                #Very special case.#
                #This is the structure of#
                #d a' d b' a" (a) (c) c' b" d c" d, resulting#
                #y' z y" d y encoding. We do not convert b'#
                #since it is a part of the prefix.#
                    if m != 2: #The current set we are encoding is#
                    #not the initial set. Then b" -> y".#
                        codedWord[-1] = 'y"'
                        codedWord.append('d')
                    else: #If it is the initial set, then b" becomes#
                    #x" instead. In addition, we must change#
                    #the associated b' to x'.#
                        codedWord[-1] = 'x"'
                        codedWord[codedWord.index("b'")] = "x'"
                        codedWord.append('d')
                else: #The case b'& b" pair happens once, but b' is#
                #not in the prefix.#
                #In this case, b' -> x' as well as b" -> x".#
                #Since b' isn't included in#
                #the prefix, we must convert it here.#
                    codedWord[-1] = 'x"'
                    codedWord[codedWord.index("b'")] = "x'"
                    codedWord.append('d')
            else: #The case b' & b" pair happens more than once.#
            #Same as before, b' -> x' as well as#
            #b" -> x" because b' isn't included in the prefix.#
                codedWord[-1] = 'x"'
                bPrimeIndexList =
                [i for i, j in enumerate(codedWord) if j == "b'"]
                #Make a list with indices of all b' in#
                #codedWord list.#
                codedWord[bPrimeIndexList[-1]] = "x'"
                codedWord.append('d')
    elif word[-3:len(word)] == ['d', 'b"', 'd']:
    #This is associated with b^ d_ interaction.#
    #b" is d_. As usual, after erasing d b" d, we add back#
    #whatever we need. The associated b' can be many different#
    #things, such as b^, d^, d_^, and y^. If the associated b' is a#
    part of the prefix, then we must have modified it previously.#
    #If this isn't the case, then we modify b' here.#
        word = word[0:last1-3]
        codedWord = codedWord[0:last2-3]
        if word.count("b'") == 1:
        #As usual, count the number of b' appearing.#
            if word[0:3] == ['d', "b'", 'd']:
            #We do not modify b' because it is#
            #included in the prefixConverter.#
                pass
            elif word[0:4] == ['d', "da'", "b'", 'da"']:
            #We do not modify b' because it is#
            #included in the prefixConverter.#
                pass
            elif word[0:5] == ['d', "a'", 'd', "b'", 'a"']:
                if m != 2:
                #If the set isn't initial, we modify b' in the#
                #prefixConverter.#
                    pass
                else: #Otherwise, b' needs to be modified to b^,#
                #since initialPrefix function don't#
                #deal with the case of beginning with d a' d.#
                    codedWord[codedWord.index("b'")] = 'b^'
            else: #The case b'& b" pair happens once, but b' is#
            #not in the prefix. In this case, we must convert#
            #b' -> b^.#
                if codedWord.count("b'") == 1:
                #This double checks if it really isn't fixed in the#
                #prefixConverter. In the case of N-C chain, the#
                #above three cannot cover it.#
                    codedWord[codedWord.index("b'")] = 'b^'
                else:
                    pass
        else: #The case b' & b" pair happens more than once.#
        #Same as before, b' -> b^#
        #because b' isn't included in the prefix.#
            bPrimeIndexList =
            [i for i, j in enumerate(codedWord) if j == "b'"]
            #Make a list with indices of all b' in codedWord list.#
            codedWord[bPrimeIndexList[-1]] = 'b^'
        codedWord.append('d_')
        codedWord.append('d')
    elif word[-4:len(word)] == ["dc'", 'b"', 'dc"', 'd']:
    #This is associated with b^ d_^ interaction. The case which#
    #didn't rise in the prefixConverter.
    #The modification is identical to the case of b^ d_.#
        word = word[0:last1-4]
        codedWord = codedWord[0:last2-4]
        if word.count("b'") == 1:
        #As usual, count the number of b' appearing.#
            if word[0:3] == ['d', "b'", 'd']:
            #We do not modify b' because it is#
            #included in the prefixConverter.#
                pass
            elif word[0:4] == ['d', "da'", "b'", 'da"']:
            #We do not modify b' because it is#
            #included in the prefixConverter.#
                pass
            elif word[0:5] == ['d', "a'", 'd', "b'", 'a"']:
                if m != 2: #If the set isn't initial, we modify b'#
                #in the prefixConverter.#
                    pass
                else: #Otherwise, b' needs to be modified to b^,#
                #since initialPrefix function don't#
                #deal with the case of ending with d a' d.#
                    codedWord[codedWord.index("b'")] = 'b^'
            else: #The case b'& b" pair happens once, but b' is not#
            #in the prefix. In this case, we must convert b' -> b^.#
                if codedWord.count("b'") == 1:
                #This double checks if it really isn't fixed in the#
                #prefixConverter. In the case of N-C chain, the#
                #above three cannot cover it.#
                    codedWord[codedWord.index("b'")] = 'b^'
                else:
                    pass
        else: #The case b' & b" pair happens more than once.#
        #Same as before, b' -> b^#
        #because b' isn't included in the prefix.#
            bPrimeIndexList =
            #[i for i, j in enumerate(codedWord) if j == "b'"]
            #Make a list with indices of all b' in codedWord list.#
            codedWord[bPrimeIndexList[-1]] = 'b^'
        codedWord.append('d_^')
        codedWord.append('d')
    else: #Any other prefix should never appear due to the way we#
    #split each set.#
        print("Error.")
        codedWord += ['<- Something went wrong here.']
    return codedWord
\end{lstlisting}
\vhh
\verb+encode.py+
\begin{lstlisting}[language=python]
"""
The following function takes a permutation having extreme pattern of
either 2413, or 2143 starting with N set.
"""

import counter
import encodeNS
import invert
import Nsplit
import prefixConverter
import suffixConverter

def encode(perm):
    word = [] #We will return this list as the final word.#
    m = 2 #Initial number.#
    n = counter.counter(perm) #Final number.#
    if n == 2: #The permutation is of extreme pattern 2413.#
        word = encodeNS.encodeN(perm)
    else:
        while m != n+1:
            subperm = Nsplit.Nsplit(perm)[0] #Split it.#
            CLast = Nsplit.Nsplit(perm)[1]
            subword1 = encodeNS.encodeN(subperm)
            #Encode the split subpermutation.#
            #Perform the modifications of prefix/suffix Converters.#
            if m == 2: #Initial set.#
                subword2 = prefixConverter.initialPrefix(subword1)
                #Result of the prefixConverter.#
                subword2 =
                suffixConverter.nonfinalSuffix(subword1,subword2,m)
                #Result of the suffixConverter.#
            elif m == n: #Final set.#
                subword2 =
                prefixConverter.noninitialPrefix(subword1,m,n)
                #Result of the prefixConverter.#
                subword2 =
                suffixConverter.finalSuffix(subword1,subword2)
                #Result of the suffixConverter.#
            else: #Neither initial nor final.#
                subword2 =
                prefixConverter.noninitialPrefix(subword1,m,n)
                #Result of the prefixConverter.#
                subword2 =
                suffixConverter.nonfinalSuffix(subword1,subword2,m)
                #Result of the suffixConverter.#
            word += subword2
            perm = Nsplit.setMinus(perm,CLast,subword2)
            perm = invert.invPerm(perm)
            m += 1
    word[-1] = 'dl'
    return word
\end{lstlisting}
\vhh
\verb+main.py+
\begin{lstlisting}[language=python]
import makeList
import createFile
import ep
import encode
import list2Str
import Nsplit

input("We first read a text file of permutations, and convert the
string type to the integer type.")
print("Enter a text file directory.")

while(True):
    try:
        permList = input()
        permList = makeList.readFile(permList)
        break
    except:
        print("There is no such a file directory.")
if len(permList[0])<10:
    permList = makeList.convList1(permList)
else:
    permList = makeList.convList2(permList)

print("Done. Permutations in the txt file are stored in permList.\n")
operation1 = input("Now, we delete all permutations which satisfy at
least one of the following conditions.\n 1.The first number is not 2,
3, or 4. \n 2.The second number is 1.")
halfList=[]
for i in permList:
        halfList.append(ep.Nstart(i))
while(True):
    try:
        halfList.remove(None)
    except:
        break
print("Done. Permutations starting with 2, 3, or 4 AND the second
number is not 1 are stored in halfList.")
operation2 = input("Type encode to create a text file of coded words.
Press enter to exist.\n")
words=[]
while operation2 != 'encode' and operation2 != '':
    print("Type encode or press enter to exit.")
    operation2 = input()
if operation2 == 'encode':
    destination =
    'C:\\Users\\Ikeda\\Desktop\\encode(%d).txt'%len(permList[0])
    for i in halfList:
        words.append(list2Str.list2Str(encode.encode(i)))
    createFile.createFile(destination,words)
    print("Done. The text file is now created.")
else:
    pass

def uniqueCheck(words):
    for i in range(len(words)-1):
        for j in words[i+1:]:
            if words[i] != j:
                pass
            else:
                print(j)

def createWordList(halfList):
    wordList = []
    for i in halfList:
        wordList.append(encode.encode(i))
    return wordList

def lengthCheck(wordList):
    n = len(wordList[0])
    for i in wordList:
        if len(i) != n:
            print(i)
        else:
            pass
\end{lstlisting}
\newpage
\phantomsection
\addcontentsline{toc}{chapter}{References}
\bibliography{thesis.bbl}
\end{document}